\documentclass[envcountsame,envcountchap]{svmono}

\usepackage{amsmath,amssymb}
\usepackage{appendix}
\usepackage[all]{xypic}
\usepackage{makeidx}         
\usepackage{graphicx}        
\usepackage{multicol}        
\usepackage[bottom]{footmisc}

\usepackage{color}   
\usepackage{hyperref}

\makeindex             

\smartqed

\newtheorem{terminology}[theorem]{Terminology}

\font\Ms=rsfs10
\font\ms=rsfs7
\def\map{\mathop{\rm map}\nolimits}
\def\imm{\looparrowright}
\def\emb{\hookrightarrow}\def\iso{\xymatrix{\ar[r]^-{\di{\cong}}&}}
\def\sto{\mathrel{\kern1.2mm\mapstochar\kern-1.2mm\rightarrow}}
\def\AA{{\cal A}}

\def\CC{{\cal C}}
\def\DD{{\cal D}}
\def\Cc{\hbox{\Ms C}\,}
\def\Mm{\hbox{\Ms M}\,}
\def\cc{\hbox{\ms C}\,}
\def\di{\displaystyle}
\def\LL{\mathbb L}
\def\N{\mathbb N}
\def\P{\mathbb P}

\def\R{\mathbb R}
\def\ZZ{\mathbb Z}
\def\fN{\frak N}
\def\st{\,\vert\,}

\def\Qq{{\mathbb Q}}
\def\Nn{{\mathbb N}}
\def\Rr{{\mathbb R}}
\def\Tt{{\mathbb T}}
\def\into{\hookrightarrow}
\def\st{\mid}
\def\map{{\rm map}}

\def\O{{\rm O}}

\def\O{{\rm O}}

\def\NN{{\mathcal N}}

\def\res{{\rm res}}
\def\free{{\rm free}}
\def\PT{\varphi}
\def\KO{K{\rm O}}
\def\dirlim{\varinjlim}
\def\OO{\mathcal{O}}

\def\SS{\mathcal{S}}

\def\im{{\rm im}\,}

\def\MM{\mathfrak{M}}
\def\PP{\mathfrak{P}}

\def\Rarr#1#2{\xrightarrow[#2]{#1}}

\def\Darr#1#2{{\scriptstyle #1}\downarrow{\scriptstyle #2}}
\def\Uarr#1#2{{\scriptstyle #1}\uparrow{\scriptstyle #2}}
\long\def\Ref#1#2#3#4#5#6{
\bibitem{#1}
{\rm #2,}
\textit{#3.}
{\rm #4}
\textbf{#5}
{\rm #6.}
}
\long\def\Refb#1#2#3#4{
\bibitem{#1}
{\rm #2,}
\textit{#3.}
#4.
}

\def\Rarr#1#2{\xrightarrow[#2]{#1}}

\def\Darr#1#2{{\scriptstyle #1}\downarrow{\scriptstyle #2}}
\def\Uarr#1#2{{\scriptstyle #1}\uparrow{\scriptstyle #2}}
\long\def\Ref#1#2#3#4#5#6{
\bibitem{#1}
{\rm #2,}
\textit{#3.}
{\rm #4}
\textbf{#5}
{\rm #6.}
}
\long\def\Refb#1#2#3#4{
\bibitem{#1}
{\rm #2,}
\textit{#3.}
#4.
}
\def\Qq{{\mathbb Q}}
\def\Zz{{\mathbb Z}}
\def\Nn{{\mathbb N}}
\def\Rr{{\mathbb R}}
\def\into{\hookrightarrow}
\def\iso{\cong}
\def\st{\mid}

\def\O{{\rm O}}
\def\NN{{\mathcal N}}
\def\res{{\rm res}}
\def\free{{\rm free}}
\def\ori{{\rm or}}
\def\PT{\varphi}
\def\KO{K{\rm O}}
\def\dirlim{\varinjlim}
\def\OO{\mathcal{O}}

\def\RO{R{\rm O}}
\def\im{{\rm im}\,}

\def\AA{\mathcal{A}}
\def\comp{\circ}

\def\nn{\mathfrak{n}}
\def\Thom{T}
\def\cpt{\infty}
\def\bpt{^+}

\begin{document}

\author{Michael Crabb and Andrew Ranicki}
\title{The geometric Hopf invariant and surgery theory}

\maketitle


\include{dedic}

\begin{center}
~\\[10ex]
\includegraphics[width=12cm]{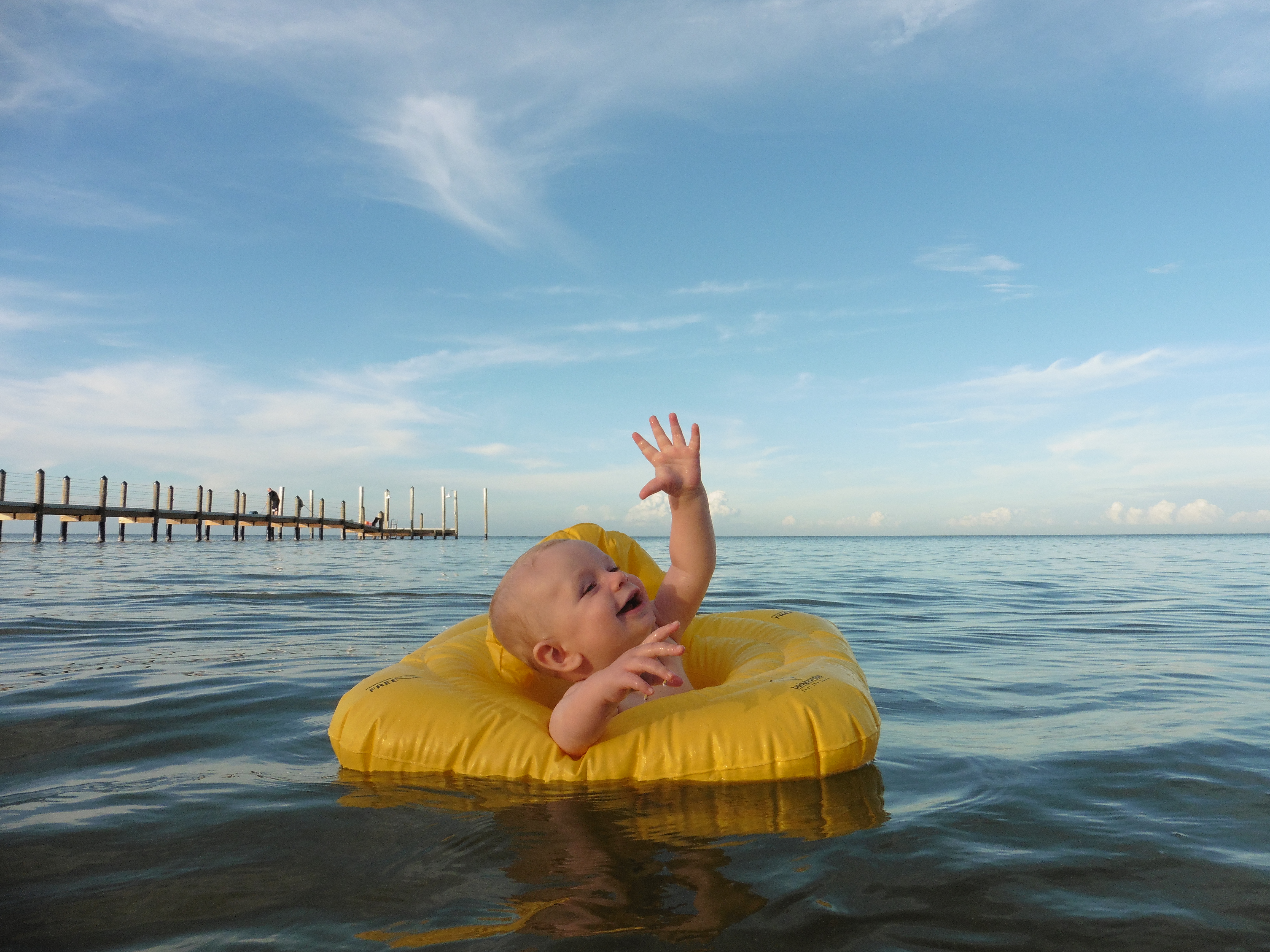}
~\\[10ex]
{\Large For Nico Marcel Vallauri}
\end{center}

\pagenumbering{roman}

\begin{preface}

This joint project is an outcome of two separate projects started by us in the 1970's,
which independently dealt with double points of maps of manifolds in geometry and quadratic structures in homotopy theory and algebra, generalizing the Hopf invariant.
However, the first-named author was more concerned with $\ZZ_2$-equivariant homotopy theory and simply-connected manifolds, while the second-named author was more concerned with chain complexes and the surgery theory of non-simply-connected manifolds.

The {\it geometric Hopf invariant}\index{geometric Hopf invariant!$h(F)$} of a stable map $F:\Sigma^{\infty}X \to
\Sigma^{\infty}Y$ for pointed spaces $X,Y$ is the stable $\ZZ_2$-equivariant map
$$h(F)~=~(F\wedge F)\Delta_X-\Delta_YF~:~\Sigma^\infty X\to
\Sigma^\infty(Y\wedge Y)$$
measuring the failure of $F$ to preserve the diagonal maps of $X$ and $Y$.
Here $\ZZ_2$ acts on $Y \wedge Y$ by transposition
$$T~:~ Y \wedge Y \to Y \wedge Y~;~(y_1,y_2) \mapsto (y_2,y_1)~.$$
The stable $\ZZ_2$-equivariant homotopy class of $h(F)$ is the primary obstruction to desuspending $F$, i.e. to $F$ being stable homotopic to an unstable map.

The original Hopf invariant of a map $F:S^3 \to S^2$ (\cite{hopf2}) was defined geometrically to be the linking number
$$H(F)~=~{\rm Lk}(F^{-1}(x),F^{-1}(y)) \in \ZZ$$
for generic $x \neq y  \in S^2$ with $F^{-1}(x),F^{-1}(y) \subset S^3$
unions of  disjoint circles. In the 1930's the Hopf invariant map was proved to be an isomorphism
$$H~:~\pi_3(S^2) \to \ZZ~;~F \mapsto H(F)~,$$
with the projection of the Hopf bundle $$\xymatrix{S^1 \ar[r] & S^3 \ar[r]^-{\displaystyle{\eta}} & S^2}$$
such that the fibres $\eta^{-1}(x),\eta^{-1}(y) \subset S^3$ are
disjoint circles with linking number $H(\eta)=1$. Samelson \cite{samelson} described the historical antecedents (Clifford, Klein) of the Hopf bundle\footnote{Samelson's paper includes the letter sent by Hopf to Freudenthal giving a first-hand account of the Hopf invariant. The letter was sent on 17th August 1928 from 30 Murray Place, Princeton. By a coincidence, this was the house occupied 50 years later by one of the authors (A.R.)}.
The original definition of the Hopf invariant $H(F) \in \ZZ$ applied to any map $F:S^{2n-1} \to S^n$ ($n \geqslant 2$), with $H:\pi_{2n-1}(S^n) \to \ZZ$ an isomorphism
for $n=2,4,8$.

In the 1940's Steenrod \cite{steenrod} used algebraic topology to interpret the Hopf invariant of a map $F:S^{2n-1} \to S^n$ ($n \geqslant 2$) as the evaluation of the cup product in the mapping cone $\Cc(F)=S^n\cup_FD^{2n}$
$$\begin{array}{l}
H^n(\Cc(F)) \times H^n(\CC(f))~=~\ZZ \times \ZZ \to H^{2n}(\Cc(F))~=~\ZZ~;\\[1ex]
\hskip100pt
(1,1) \mapsto 1 \cup 1~=~H(F)~(=~0~{\rm for~odd}~n)~.
\end{array}$$
 The mod 2 Hopf invariant $H_2(F) \in \ZZ_2$ of a map $F:S^j \to S^k$ was then defined  using the Steenrod square in
$\Cc(F)=S^k \cup_F D^{j+1}$
$$Sq^{j-k+1}~:~H^k(\Cc(F);\ZZ_2)~=~\ZZ_2 \to H^{j+1}(\Cc(F);\ZZ_2)~=~\ZZ_2~;~1 \mapsto H_2(F)$$
or equivalently the  functional Steenrod square $Sq^{j-k+1}_F$.

For a pointed space $X$ write the loop space and the suspension as usual
$$\Omega X~=~{\rm map}_*(S^1,X) ~,~ \Sigma X~=~S^1 \wedge X~.$$
\indent In the 1950's and 1960's it was realized that the Hopf invariant was the first example of the primary obstruction to desuspending a stable map $F:\Sigma^k X \to \Sigma^k Y$ for any pointed spaces $X,Y$ (not necessarily spheres),
with $1 \leqslant k \leqslant \infty$.  A $k$-stable map $F:\Sigma^kX \to \Sigma^kY$
can be desuspended, i.e. is homotopic to $\Sigma^kF_0$ for a map $F_0:X \to Y$,
 if and only if the adjoint map
$${\rm adj}(F)~:~X \to \Omega^k\Sigma^kY~;~x \mapsto (s \mapsto F(s,x))~(s \in S^k)$$
can be compressed into
$$Y \subset \Omega^k\Sigma^kY~;~y \mapsto (s \mapsto (s,y))~.$$
\indent For any space $Y$ the homotopy groups of $\Omega\Sigma Y$ are
$$\pi_n(\Omega \Sigma Y)~=~\pi_{n+1}(\Sigma Y)~.$$
James \cite{imj2} showed that for a connected $Y$ there is a homology equivalence
$$\Omega \Sigma Y~\simeq~ \bigvee\limits^{\infty}_{j=1} \bigwedge \limits_j Y~.$$
The adjoint of a 1-stable map $F:\Sigma X \to \Sigma Y$ is
a map ${\rm adj}(F):X \to \Omega \Sigma Y$. The primary
obstruction to desuspending $F$ is the composite
$$H_*(X) \xymatrix{\ar[r]^{{\rm adj}(F)_*}&}
H_*(\Omega\Sigma Y)~=~
\bigoplus\limits^{\infty}_{j=1} H_*(\bigwedge \limits_j Y)
\to H_*(Y \wedge Y)~.$$
For $F:S^3=\Sigma(S^2)\to S^2=\Sigma(S^1)$ this
is just the Hopf invariant
$$H_2(S^2)=\ZZ \xymatrix{\ar[r]^{{\rm adj}(F)_*}&}
H_2(\Omega S^2)~=~
\bigoplus\limits^{\infty}_{j=1} H_2(S^j) = H_2(S^2)=\ZZ~;~1 \mapsto H(F)~.$$

The homotopy groups of the space
$$\Omega^{\infty}\Sigma^{\infty}Y~=~\varinjlim\limits_k\Omega^k\Sigma^kY$$
are the stable homotopy groups of $Y$
$$\pi_n(\Omega^\infty\Sigma^\infty Y)~=~\varinjlim\limits_k\pi_{n+k}(\Sigma^kY)~.$$
The desuspension of $F$ is equivalent to compressing the adjoint map
$$X \to \Omega^{\infty}\Sigma^{\infty}Y~;~x \mapsto (s \mapsto F(s,x))~(s \in S^{\infty})$$
into
$$Y \subset \Omega^{\infty}\Sigma^{\infty}Y~;~y \mapsto (s \mapsto (s,y))~.$$
For connected $Y$ the space $\Omega^{\infty}\Sigma^{\infty}Y$ has a filtration in stable homotopy theory with successive quotients
$$(E\Sigma_k)^+\wedge_{\Sigma_k} \bigwedge\limits_k Y~(k \geqslant 1)~.$$
Here, $\Sigma_k$ is the permutation group on $k$ letters acting by
$$\Sigma_k \times \bigwedge\limits_k Y \to \bigwedge\limits_k Y~;~
(\sigma,(y_1,y_2,\dots,y_k)) \mapsto (y_{\sigma(1)},y_{\sigma(2)},\dots,y_{\sigma(k)})$$
and $E\Sigma_k$ a contractible space with a free $\Sigma_k$-action. (There is also a version for disconnected $Y$). The $k$th filtration quotient correspond to the $k$-tuple points of maps of manifolds.
We shall only be concerned with the first filtration quotient, for $k=2$, as we are
particularly interested in double points. We write
$$E\Sigma_2~=~S(\infty)~,~\Sigma_2~=~\ZZ_2~.$$
\indent  For any pointed space $Y$ define the {\it quadratic construction}\index{quadratic construction!on space $Y$, $Q_\bullet(Y)$} pointed space
$$Q_\bullet(Y)~=~S(\infty)^+ \wedge_{\ZZ_2}(Y \wedge Y)~.$$
Such quadratic spaces were introduced in homotopy theory
by Toda, in the 1950's. For  any pointed spaces $X,Y$ the stable
homotopy groups $\{X,Y\}$, $\{X,Q_\bullet(Y)\}$ and the stable $\ZZ_2$-equivariant homotopy group $\{X,Y \wedge Y\}_{\ZZ_2}$
fit into a split short exact sequence
$$\xymatrix{0 \ar[r] & \{X,Q_\bullet(Y)\} \ar[r] & \{X,Y \wedge Y\}_{\ZZ_2} \ar[r]^-{\rho} & \{X,Y\} \ar[r]&0}$$
with $\rho$ the fixed point map. Following Crabb \cite{crabb} and Crabb and James \cite{crabbjames} the primary obstruction to desuspending a stable map $F:\Sigma^{\infty}X \to \Sigma^{\infty}Y$ is the stable $\ZZ_2$-equivariant homotopy class of the geometric Hopf invariant $h(F)$, which can also be
viewed as the stable homotopy class of a stable map
$$h'(F)~:~\Sigma^{\infty}X \to \Sigma^{\infty}Q_\bullet(Y)~,$$
that is
$$h(F)~=~h'(F) \in {\rm ker}(\rho)~=~\{X,Q_\bullet(Y)\}~.$$
\indent The algebraic theory of surgery of Ranicki \cite{ranicki1,ranicki2} was based on the {\it quadratic construction}\index{quadratic construction!on chain complex $C$, $Q_\bullet(C)$} on a chain complex $C$, the chain complex defined by
$$Q_\bullet(C)~=~W\otimes_{\ZZ[\ZZ_2]}(C \otimes_\ZZ C)$$
with $W=C(S(\infty))$ a free $\ZZ[\ZZ_2]$-module resolution of $\ZZ$.
A stable map $F:\Sigma^{\infty}X \to \Sigma^{\infty}Y$ was shown to induce a natural chain homotopy class of chain maps
$$\psi_F~:~\dot C(X) \to Q_\bullet(\dot C(Y))~=~\dot C(Q_\bullet(Y))$$
called the {\it quadratic construction}\index{quadratic construction!on stable map $F$, $\psi_F$} on $F$.

In May 1999 it became clear to us that the quadratic construction $\psi_F$ of
\cite{ranicki1,ranicki2} was induced by the geometric Hopf invariant $h(F)$ of \cite{crabb,crabbjames}, in the sense of
$$\psi_F~=~h'(F)~:~\dot C(X) \to \dot C(Q_\bullet(Y))$$
up to chain homotopy.
In this book we unite the two approaches into a single theory, together with the applications to the double points of immersions and to various constructions in the algebraic theory of surgery. In the applications it is necessary
to consider $\pi$-equivariant, fibrewise and local versions of the quadratic construction/geometric Hopf invariant.

Throughout the book ``manifold'' will mean ``compact differentiable manifold'',
except that in dealing with the total surgery obstruction in \S\ref{total} we shall be considering compact topological manifolds.

The double point set of a map $f:M \to N$ is the free $\ZZ_2$-space
$$D_2(f)~=~\{(x,y) \in M \times M \,\vert\, x \neq y \in M,\,f(x)=f(y) \in N\}$$
with $T \in \ZZ_2$ acting by
$$T~:~D_2(f) \to D_2(f)~;~(x,y) \mapsto (y,x)~.$$
An immersion of manifolds $f:M^m \imm N^n$ has a normal $(n-m)$-plane bundle $\nu_f$, with  a Thom space $T(\nu_f)$. There exists a map $e:M \to \R^k$ ($k$ large)
such that
$$M \emb \R^k \times N~;~x \mapsto (e(x),f(x))$$
is an embedding, and the Pontrjagin-Thom construction
gives a stable Umkehr map
$$F~:~\Sigma^{\infty}N^+\to \Sigma^{\infty}T(\nu_f)~.$$
$f$ is an embedding if and only if $D_2(f)=\emptyset$, in which case $F$ is unstable.
Our main result is the Double Point Theorem \ref{doubleHopf}, which stably factors
the geometric Hopf invariant $h(F)$ for a generic immersion $f$ through
the double point set
$$h(F)~:~N^+ \to D_2(f)^+ \to Q_\bullet(T(\nu_f))~.$$
The stable homotopy class of $h(F)$ depends only on the
regular homotopy class of $f$.  If $f$ is regular homotopic to an
embedding then $h(F)$ is stably null-homotopic.

We use the geometric Hopf invariant to provide the homotopy theory underpinnings of the following constructions of surgery theory:

\noindent 1. The double point invariant of Wall \cite{wall2} for an immersion $f:S^m \imm N^{2m}$
$$\mu(f) \in \ZZ[\pi_1(N)]/\{ x-(-)^m\overline{x}\,\vert\, x \in \ZZ[\pi_1(N)]\}~(\overline{g}=g^{-1} \in \pi_1(N))
$$
such that $\mu(f)=0$ if (and, for $m > 2$, only if)  $f$ is regular
homotopic to an embedding. The expression of $\mu(f)$  in terms of the geometric Hopf invariant is in \S\ref{self}.

\noindent 2. The surgery obstruction $\sigma_*(f,b) \in L_n(\ZZ[\pi_1(X)])$ (Wall \cite{wall2}) of an $n$-dimensional normal map $(f,b):M \to X$
is such that $\sigma_*(f,b)=0$  if (and, for $n > 4$, only if)  $(f,b)$ is normal bordant to a homotopy equivalence. The invariant was expressed in Ranicki \cite{wall2} as the cobordism class of an $n$-dimensional
$\ZZ[\pi_1(X)]$-module chain complex $C$ with a quadratic Poincar\'e duality determined
via  the quadratic construction $\psi_F$ of a $\pi_1(X)$-equivariant stable
Umkehr map $F:\Sigma^{\infty}\widetilde{X}^+ \to \Sigma^{\infty}\widetilde{M}^+$.
Here $\widetilde{X}$ is the universal cover of $X$, $\widetilde{X}^+=\widetilde{X} \cup \{{\rm pt.}\}$, $\widetilde{M}=f^*\widetilde{X}$ is
the pullback cover of $M$, and
$$H_*(C)~=~{\rm ker}(\widetilde{f}_*:H_*(\widetilde{M}) \to H_*(\widetilde{X}))~.$$
The expression of $\psi_F$ in terms of the geometric Hopf invariant is in \S\ref{quad}.

\noindent 3. The total surgery obstruction $s(X) \in \SS_n(X)$ (Ranicki \cite{ranicki0,ranicki3})
of an $n$-dimensional geometric Poincar\'e complex $X$ is such that $s(X)=0$  if (and, for $n > 4$, only if)  $X$ is homotopy equivalent to an $n$-dimensional topological manifold.  The invariant  is
the cobordism class of the $\ZZ[\pi_1(X)]$-contractible $\ZZ$-module chain complex
$$C={\cal C}([X] \cap-:C(X)^{n-*} \to C(X))_{*+1}$$
with an $X$-local $(n-1)$-dimensional quadratic Poincar\'e duality determined
via  the spectral quadratic construction $s\psi_F$ of a $\pi_1(X)$-equivariant semistable homotopy equivalence
$F:\Sigma^{\infty}T(\nu_{\widetilde{X}})^* \to \Sigma^{\infty}\widetilde{X}^+$ inducing the chain equivalence $[X] \cap -$.
The expression of $s(X)$ in terms of the geometric Hopf invariant is in \S\ref{total}.

Chapters 1-7 develop the relationship between the geometric Hopf invariant and double points in such a natural way that the $\pi$-equivariant
constructions required for non-simply-connected surgery theory fall out without extra work. Along the way,
we point out the many particular instances of this relationship in the literature, and explain how they fit into our development. The applications to surgery theory  are then considered in Chapter 8.

In principle, it is possible to develop surgery theory in the context
of the fibrewise homotopy theory of Crabb and James \cite{crabbjames},
relating it to the controlled and bounded surgery theories. However,
we shall not do this here, restricting ourselves to the application
of $\ZZ_2$-equivariant fibrewise homotopy theory in Appendix \ref{appendix1} to obtain the $\pi$-equivariant quadratic construction
used in the conventional non-simply-connected surgery theory.

Appendix \ref{appendix2} provides a treatment of $\ZZ_2$-equivariant bordism theory.

In 2010 we published a joint paper {\it The geometric Hopf invariant and double points} (\cite{crabbranicki}) analyzing double points using a combination of stable fibrewise homotopy and immersion theories, and obtaining the double point theorem of Chapter 6 from a somewhat different viewpoint. The paper interprets the Smale-Hirsch-Haefliger regular homotopy classification of immersions $f$ in the metastable dimension range $3m<2n-1$ (when a  generic $f$ has no triple points) in terms of the geometric Hopf invariant $h(F)$. The paper is reprinted here in Appendix \ref{appendix3}, and may be regarded as an extended introduction to this book.

The photo of the dedicatee on page ii is by Carla Ranicki.

We are grateful to the referees for their valuable comments.

\noindent Michael Crabb, University of Aberdeen, m.crabb@abdn.ac.uk\\[0.5ex]
Andrew Ranicki, University of Edinburgh, a.ranicki@ed.ac.uk

\hfill October 2017

\end{preface}

\tableofcontents

\pagenumbering{arabic}

\mainmatter

\chapter{The difference construction} \label{difference}

Chapter \ref{difference} describes the difference construction in both homotopy and chain homotopy theory. The homotopy version will be used in Chapter \ref{umkehr} to define the subtraction of maps.

We shall be dealing with both pointed and unpointed spaces, and
likewise for maps and homotopies.

\section{Sum and difference}

For a pointed space $X$ we shall denote the base point by $* \in X$,
and if $Y$ is also a pointed space then $*:X \to Y$ is the constant pointed map. The wedge of pointed spaces $X,Y$ is the space
$$X \vee Y~=~ X \times \{*_Y\} \cup \{*_X\} \times Y \subseteq X \times Y$$
with base point $(*_X,*_Y) \in X \vee Y$.

\begin{definition} \label{sumdifference}
Let $X$ be a pointed space.\\
(i) A {\it sum map} $\nabla:X \to X \vee X$ is a pointed map  such that the composites
$$\begin{array}{l}
\xymatrix{(1 \vee *)\nabla~:~X \ar[r]^-{\di{\nabla}} & X \vee X
\ar[r]^-{\di{1 \vee *}} &X}~,\\[1ex]
\xymatrix{(* \vee 1)\nabla~:~X \ar[r]^-{\di{\nabla}} & X \vee X
\ar[r]^-{\di{* \vee 1}} &X}
\end{array}$$
are  pointed homotopic to the identity.\\
(ii)  For $\nabla$ as in (i) the {\it sum} of pointed maps $f_1,f_2:X \to Y$ is defined by
$$f_1+f_2~=~(f_2 \vee f_1)\nabla~:~X \to Y~.$$
(iii) A {\it difference map} $\overline{\nabla}:X \to X \vee X$ is a pointed map  such that the composite
$$\xymatrix{(* \vee 1)\overline{\nabla}~:~X \ar[r]^-{\di{\overline{\nabla}}} & X \vee X
\ar[r]^-{\di{* \vee 1}} &X}$$
is pointed homotopic to the identity, and the composite
$$\xymatrix{(1 \vee 1)\overline{\nabla}~:~X \ar[r]^-{\di{\overline{\nabla}}} & X \vee X
\ar[r]^-{\di{1 \vee 1}} &X}$$
is pointed null-homotopic.\\
(iv)  For $\overline{\nabla}$ as in (iii) the {\it difference} of pointed maps $f_1,f_2:X \to Y$ is defined by
$$f_1-f_2~=~(f_2 \vee f_1)\overline{\nabla}~:~X \to Y~.$$
\hfill\qed
\end{definition}

\begin{example} \label{circle}
{\rm The standard sum and difference maps on $S^1$
$$\begin{array}{l}
\nabla_{S^1}~:~S^1 \to S^1 \vee S^1~;~
e^{2\pi i t} \mapsto
\begin{cases} (e^{4\pi i t})_1&{\rm if}~0 \leqslant t \leqslant 1/2 \\
(e^{2\pi i (2t-1)})_2&{\rm if}~1/2 \leqslant t \leqslant 1~,
\end{cases}\\[3ex]
\overline{\nabla}_{S^1}~:~S^1 \to S^1 \vee S^1~;~
e^{2\pi i t} \mapsto
\begin{cases} (e^{2\pi i (1-2t)})_1&{\rm if}~0 \leqslant t \leqslant 1/2 \\
(e^{2\pi i (2t-1)})_2&{\rm if}~1/2 \leqslant t \leqslant 1
\end{cases}
\end{array}
$$
satisfy the conditions of Definition \ref{sumdifference}, and
induce the usual sum and difference maps on the fundamental
group $\pi_1(X)$ of a pointed space $X$.
}
\hfill\qed
\end{example}

\section{Joins, cones and suspensions}

The {\it one-point adjunction} of an unpointed space $X$ is the pointed space\index{one-point!adjunction, $X^+$}
$$X^+~=~X \sqcup \{*\}$$
obtained by adjoining a base point $*$, with open subsets $U, U\sqcup
\{*\} \subseteq X^+$ for open subsets $U \subseteq X$.

The {\it smash product} of pointed spaces $X,Y$ is the pointed space\index{smash product, $X\wedge Y$}
$$X\wedge Y~=~X \times Y/(X \times \{*\} \cup\{*\} \times Y)~.$$
In particular, for unpointed spaces $A,B$
$$A^+\wedge B^+~=~(A \times B)^+~.$$
\indent
The {\it join}  of spaces $X,Y$ is the space defined as usual by\index{join, $X*Y$}
$$X*Y~=~I \times X \times Y/\{(0,x,y_1) \sim (0,x,y_2)\,,\,(1,x_1,y)\sim (1,x_2,y)\}~,$$
so that
$$X~=~[0 \times X \times Y]~,~Y~=~[1 \times X \times Y] \subset X*Y~.$$
\indent
The {\it cone} on a space $X$ is the space\index{cone, $cX$}
$$cX~=~X*\{{\rm pt.}\}~=~I  \times X/\{1\} \times X~=~I \wedge X^+$$
with $I$ pointed at 1.
Identify $X=\{0\}\times X \subset cX$.
A map $cX \to Y$ is a map $X \to Y$ with a homotopy to a constant map.
The {\it reduced cone} on a pointed space $X$ is the pointed space\index{reduced!cone, $CX$}
$$CX~=~cX/c\{*\}~=~I \times X/(\{1\} \times X \cup I \times \{*\})~=~I \wedge X~.$$
\indent
A pointed homotopy $h:f \simeq g:X \to Y$ is a pointed map
$h:I^+ \wedge X \to Y$ with
$$h(0,x)~=~f(x)~,~h(1,x)~=~g(x)\in Y~~(x \in X)~.$$
A {\it null-homotopy} of a pointed map $f:X \to Y$ is a pointed
homotopy $f \simeq *:X \to Y$, which is the same as an extension of $f$
to a pointed map $\delta f:CX \to Y$.\index{null-homotopy}
Let $[X,Y]$ be the pointed set of pointed homotopy classes of pointed maps
$f:X \to Y$.

The {\it suspension} of a space $X$ is the space\index{suspension!space, $sX$}
$$sX~=~S^0*X~=~I\times X/\{(i,x) \sim (i,x')\,\vert\, i=0,1,~ x,x' \in X\}~.$$
A map $sX \to Y$ is a map $X \to Y$ with two homotopies to constant maps.
The {\it reduced suspension} of a pointed space $X$ is the pointed space
\index{reduced!suspension, $\Sigma X$}
$$\Sigma X~=~sX/s\{*\}~=~I \times X/\{0,1\} \times X \cup I \times \{*\})~=~S^1 \wedge X~.$$
A pointed map $\Sigma X \to Y$ is a pointed map $X \to Y$ with two
null-homotopies.

The {\it mapping cylinder} of a map $F:X \to Y$ is the
space\index{mapping!cylinder, $\Mm(F)$}
$$m(F)~=~I  \times X \cup_F Y/\{(0,x) \sim F(x)\,\vert\, x \in X\}~.$$
The {\it reduced mapping cylinder} of a pointed map $F:X \to Y$ is the pointed
space\index{mapping!cylinder, reduced $\Mm(F)$}
$$\Mm(F)~=~m(F)/I \times\{*\}~=~I^+ \wedge X \cup_F Y~.$$
\indent
The {\it mapping cone} of a map $F:X \to Y$ is the space\index{mapping!cone, $\Cc(F)$}
$$c(F)~=~m(F)/\{1\} \times X~=~cX \cup_F Y~.$$
The {\it reduced mapping cone} of a pointed map $F$ is the pointed
space\index{mapping!cone, reduced $\Cc(F)$}
$$\Cc(F)~=~c(F)/\{1\} \times X~=~CX \cup_F Y~,$$
Let $G:Y \to \Cc(F)$ be the inclusion.  The sequence of pointed spaces
and pointed maps
$$\xymatrix@C+10pt{X \ar[r]^-{\di{F}} &
Y \ar[r]^-{\di{G}} & \Cc(F)}$$
is a cofibration, such that for any pointed space $W$ there is induced an exact sequence of pointed homotopy sets
$$\xymatrix@C+10pt{[\Cc(F),W] \ar[r]^-{\di{G^*}} &
[Y,W] \ar[r]^-{\di{F^*}} & [X,W]~.}$$
Let $H:\Cc(F) \to \Sigma X$ be the composite
$$\xymatrix@C+20pt{H~:~\Cc(F)\ar[r]^-{inclusion} & \Cc(G)
\ar[r]^-{projection} & \Sigma X~.}$$
The projection $\Cc(G) \to \Sigma X$ is a homotopy equivalence.  The sequence of pointed spaces and pointed maps
$$\xymatrix@C+10pt{X \ar[r]^-{\di{F}} &
Y \ar[r]^-{\di{G}} & \Cc(F) \ar[r]^-{\di{H}} &  \Sigma X \ar[r]^-{\di{\Sigma F}} &\Sigma Y \ar[r] &\dots}$$
is a  {\it homotopy cofibration}.

\begin{proposition}~ \label{cone}\label{BP1}
{\rm (i)} For any map $F:X \to Y$ the maps
$$\begin{array}{l}
i~:~X \to m(F)~;~x \mapsto [1,x]~,\\[1ex]
j~:~Y \to m(F)~;~y \mapsto [y]~,\\[1ex]
k~:~I \times X \to m(F)~;~(t,x) \mapsto [t,x]
\end{array}$$
are such that $j$ is a homotopy equivalence, and $k:jF\simeq i:X\to m(F)$.
Similarly for a pointed map $F$ and the pointed mapping cylinder $\Mm(F)$.\\
{\rm (ii)} For any pointed map $F:X \to Y$ the reduced
mapping cone $Z=\Cc(F)$ fits into the homotopy cofibration sequence
$$\xymatrix@C+10pt{X \ar[r]^-{\di{F}} &
Y \ar[r]^-{\di{G}} &
Z\ar[r]^-{\di{H}}& \Sigma X \ar[r]^-{\di{\Sigma F}} & \Sigma Y \ar[r] &\dots}$$
with the inclusion $CX \to Z$ defining a null-homotopy $GF\simeq *:X \to Z$. For any pointed space $W$ there is induced a
Barratt-Puppe exact sequence of pointed homotopy sets
$$\xymatrix{
\dots \ar[r] & [\Sigma Y,W] \ar[r]^-{\di{\Sigma F^*}} &[\Sigma X,W] \ar[r]^-{\di{H^*}} &[\Cc(F),W] \ar[r]^-{\di{G^*}} & [Y,W] \ar[r]^-{\di{F^*}} & [X,W]~.}$$
\hfill\qed
\end{proposition}

\begin{example}{\rm
For an unpointed space $X$ define the pointed map
$$F~:~X^+ \to S^0~;~x \mapsto -1~,~* \mapsto +1~.$$
The reduced mapping cylinder and cone of $F$ are the cone and suspension of $X$
$$\Mm(F)~=~cX~,~\Cc(F)~=~sX~.$$
By Proposition \ref{cone}
there is defined a homotopy cofibration sequence of pointed spaces
$$\xymatrix@C+10pt{X^+ \ar[r]^-{\di{F}} &
S^0 \ar[r]^-{\di{G}} &
sX\ar[r]^-{\di{H}}& \Sigma X^+ \ar[r]^-{\di{\Sigma F}} & S^1 \ar[r] &\dots}$$
with
$$\begin{array}{l}
G~:~S^0 \to sX~;~-1 \mapsto (0,X)~,~1 \mapsto (1,X)~,\\[1ex]
H~:~sX \to sX/S^0~=~\Sigma X^+~;~(t,x) \mapsto (t,x)~,\\[1ex]
\Sigma F~:~\Sigma X^+ \to S^1~=~I/(0=1)~;~(t,x) \mapsto t~.
\end{array}$$
The reduced suspension $\Sigma X^+$ is obtained from $sX$ by
identifying the two suspension points $S^0 \subset sX$, with
$G$ the inclusion and $H$ the projection.
If $X$ is non-empty then for any $x \in X$ the pointed maps
$$\begin{array}{l}
E~:~S^0 \to X^+~;~1 \mapsto *~,~-1 \mapsto x~,\\[1ex]
\Sigma E~:~\Sigma (S^0)~=~S^1 \to \Sigma X^+~;~t \mapsto (t,x)
\end{array}$$
are such that $FE=1:S^0 \to S^0$ and
$$H \vee \Sigma E~:~sX \vee S^1 \to \Sigma X^+~.$$
is a homotopy equivalence.
\hfill\qed}
\end{example}

\section{The difference construction $\delta(f,g)$}\label{difcon}

The {\it sum} of  homotopies\index{sum of homotopies, $f+g$}
$$f~:~c~\simeq~d~,~g~:~d~\simeq~e~:~A \to B$$
is the concatenation homotopy
$$f+g~:~c~\simeq~e~:~A \to B$$
defined by
$$f+g~:~A \times I \to B~;~(a,t) \mapsto
\begin{cases}
f(a,2t)&\hbox{\rm if $0 \leqslant t \leqslant 1/2$}\\
g(a,2t-1)&\hbox{\rm if $1/2 \leqslant t \leqslant 1$}~.
\end{cases}$$
The {\it reverse} of a homotopy $f:c \simeq d:A \to B$ is the
homotopy \index{reverse of a homotopy, $-f$}
$$-f~:~d~\simeq~c~:~A \to B$$
defined by
$$-f~:~A \times I \to B~;~(a,t) \mapsto f(a,1-t)~.$$
\indent
As usual, given pointed spaces $A,B$ let $[A,B]$ be the set of
pointed homotopy classes of pointed maps $A \to B$.

A pointed map $d:\Sigma A \to B$ is essentially the same as a pointed homotopy
$$d~:~\{*\}~\simeq~\{*\}~:~A \to B$$
and $[\Sigma A,B]$ is a group, with sum and differences given by
$$\begin{array}{l}
[\Sigma A,B] \times [\Sigma A,B] \to [\Sigma A,B]~;~
(f,g) \mapsto \begin{cases}
f+g=(g \vee f)\nabla_{\Sigma A}\\[1ex]
f-g=(g \vee f)\overline{\nabla}_{\Sigma A}
\end{cases}
\end{array}$$
using the sum and difference maps (\ref{sumdifference}) given by
Example \ref{circle}
$$ \begin{cases}
\nabla_{\Sigma A}~=~\nabla_{S^1}\wedge 1_A\\[1ex]
\overline{\nabla}_{\Sigma A}~=~\overline{\nabla}_{S^1}\wedge 1_A
\end{cases}
~:~
\Sigma A ~=~S^1 \wedge A\to \Sigma A \vee \Sigma A~=~(S^1 \vee S^1) \wedge A~.$$

\begin{definition}~\label{reldif1}
{\rm
Let $e:A \to B$ be a pointed map, and let
$$f~:~e~\simeq~*~,~g~:~e~ \simeq~ *~:~A \to B$$
be two null-homotopies,
corresponding to two extensions $f,g:CA \to B$ of $e$.
The {\it rel $A$ difference} of $f,g$
is the pointed map\index{relative difference!maps, $\delta(f,g)$}
$$\delta(f,g)~=~f\cup -g~:~\Sigma A \to B~;~(t,a) \mapsto
\begin{cases}
g(1-2t,a)&\hbox{if $0 \leqslant  t \leqslant  1/2$}\\[1ex]
f(2t-1,a)&\hbox{if $1/2 \leqslant  t \leqslant 1$~.}
\end{cases}$$
Equivalently
$$\xymatrix@C+25pt{
\delta(f,g)~:~ \Sigma A\ar[r]^-{\di{-1\cup 1}} & CA \cup_A CA
\ar[r]^-{\di{g \cup f}} & B}$$
with $-1\cup 1$ the homeomorphism
$$-1\cup 1~:~\Sigma A \to CA \cup_A CA~;~(t,a) \mapsto
\begin{cases}(1-2t,a)_1&\hbox{if $0 \leqslant  t \leqslant  1/2$}\\
(2t-1,a)_2&\hbox{if $1/2 \leqslant  t \leqslant  1$.}\end{cases}$$
}
\hfill\qed
\end{definition}

\begin{proposition}~Every pointed map $d:\Sigma A \to B$ is the rel $A$
difference $\delta(f,g)$ of null-homotopies $f:e \simeq *$, $g:e \simeq *$
of a map $e:A \to B$.
\end{proposition}
\begin{proof} The pointed maps
$$\begin{array}{l}
e~:~A \to B~;~a \mapsto d(1/2,a)~,\\[1ex]
f~:~CA \to B~;~(t,a) \mapsto d((1+t)/2,a)~,\\[1ex]
g~:~CA \to B~;~(t,a) \mapsto d((1-t)/2,a)
\end{array}$$
are such that
$$e(a)~=~f(0,a)~=~g(0,a) \in B~(a \in A)~,~\delta(f,g)~=~d~.$$
\hfill\qed\end{proof}

\begin{remark}{\rm
(i) In the special case $e=*:A \to B$
the null-homotopies $f,g:e\simeq \{*\}$ in  Proposition \ref{reldif1} are just
maps $f,g:\Sigma A \to B$. The sum is
$$f+g~=~(g \vee f)\nabla_{\Sigma A}~=~f \cup g~:~\Sigma A~=~CA\cup_A CA \to B$$
and the rel $A$ difference is just
$$\delta(f,g)~=~f-g~=~(g \vee f)\overline{\nabla}_{\Sigma A}~=
~f \cup -g~:~\Sigma A~=~CA\cup_A CA \to B~.$$
(ii) In the special case $A=S^n$ the rel $A$ difference
of maps $f,g:CA=D^{n+1} \to B$ such that $f\vert = g\vert:S^n \to B$
is the {\it separation element} $\delta(f,g):\Sigma A=S^{n+1} \to B$
of James \cite[Appendix]{imj3}.\index{separation element, $\delta(f,g)$}}\\
\qed
\end{remark}

\begin{proposition}~ \label{main}
{\rm ~(i)} A homotopy $h:f \simeq g:CA \to B$ which is constant on $A$ determines a null-homotopy
$$\delta(h)~:~\delta(f,g)~ \simeq~ *~:~\Sigma A\to B~.$$
{\rm (ii)} For any map $f:CA \to B$ there is a canonical null-homotopy
$$\delta(f,f)~ \simeq~ *~ :~\Sigma A \to B~.$$
{\rm (iii)} Let $i:X \to A$ be the inclusion of a subspace $X \subseteq A$
such that the projection $\Cc(i) \to A/X$ is a homotopy equivalence,
and let $j:A \to A/X$ be the projection.
If $f,g:CA \to B$ are maps which agree on $A \cup CX \subseteq CA$
the composite
$$\xymatrix@C+20pt{\delta(f,g)(\Sigma i)~:~
\Sigma X \ar[r]^-{\di{\Sigma i}} & \Sigma A \ar[r]^-{\di{\delta(f,g)}} & B}$$
is equipped with a canonical null-homotopy
$h:\delta(f,g)(\Sigma i) \simeq *$, giving rise to a rel $A \cup CX$ difference map
$$\xymatrix@C+25pt{
\delta(f,g)~:~\Sigma(A/X) ~=~CA/(A \cup CX)~
\simeq~ \Cc(\Sigma i) \ar[r]^-{\di{\delta(f,g)\cup h}} &B}$$
such that the rel $A$ difference map factors up to homotopy as
$$\xymatrix@C+20pt{
\delta(f,g)~:~\Sigma A \ar[r]^-{\di{\Sigma j}} & \Sigma (A/X)
\ar[r]^-{\di{\delta(f,g)}} &B~.}$$
\end{proposition}
\begin{proof}
(i) The map defined by
$$\begin{array}{l}
\delta(h)~:~C(\Sigma A) \to B~;\\[1ex]
(s,t,a) \mapsto \begin{cases}
h(2s,1-2t,a)&{\rm if}~0 \leqslant s,t \leqslant 1/2\\
f(2t-1,a)& {\rm if}~0 \leqslant s \leqslant 1/2 \leqslant t
\leqslant 1 \\
f(1+4(s-1)t,a) &{\rm if}~0 \leqslant t \leqslant
1/2\leqslant s \leqslant 1\\
f(4s-3+4(1-s)t,a)&{\rm if}~1/2 \leqslant s,t \leqslant 1
\end{cases}
\end{array}$$
is a null-homotopy of $\delta(f,g)$.\\
(ii) This is the special case $f=g$ of (i), with $h:f\simeq g$ the
constant homotopy.\\
(iii) Let
$$e~=~f\vert~=~g\vert~:~A \cup CX \to B~,$$
so that
$$\delta(f,g)(\Sigma i)~:~\Sigma X \to B~;~
(t,x) \mapsto \begin{cases}
e(1-2t,x)&\hbox{if $0 \leqslant  t \leqslant  1/2$}\\[1ex]
e(2t-1,x)&\hbox{if $1/2 \leqslant  t \leqslant 1$}
\end{cases}$$
with a canonical null-homotopy $h:\delta(f,g)(\Sigma i)\simeq *$.\\
\hfill\qed\end{proof}

Thus if $f,g:CA \to B$ agree on $A\cup CX\subseteq CA$
the homotopy class of $\delta(f,g):\Sigma (A/X) \to B$ is an
obstruction to the existence of a homotopy $f \simeq g:CA \to B$
which is constant on $A\cup CX$.

If $f,g$ agree on a neighbourhood of $A \subset CA$ and
$g=*$ on the frontier and outside the neighbourhood, then
$\delta(f,g)$ is homotopic to the restriction of $ff$ to the
complement of the neighbourhood~:

\begin{proposition}~\label{near}
{\rm ~(i)} If $f,g:CA \to B$ are  maps such that
for some neighbourhood $U \subseteq CA$ of $A \subseteq CA$
$$g(t,a)~=~\begin{cases}
f(t,a)&{\it if}~ (t,a) \in U\\
*&{\it if}~ (t,a) \in \overline{CA \backslash U}
\end{cases}$$
then there exists a homotopy\index{$\gamma(f,g)$}
$$\gamma(f,g)~:~\delta(f,g)~\simeq~f'~:~\Sigma A \to B$$
with $f'$ defined by
$$f'~:~\Sigma A=CA/A \to B~;~(t,a) \mapsto
\begin{cases}
*&{\it if}~ (t,a) \in U\\
f(t,a)&{\it if}~ (t,a) \in \overline{CA \backslash U}
\end{cases}~.$$
{\rm (ii)} The homotopy $\gamma(f,g)$ in {\rm ~(i)} can be chosen to be
natural, meaning that given commutative squares of maps
$$\xymatrix{
CA_1 \ar[r]^-{\di{f_1}} \ar[d]_-{\di{h}} &  B_1\ar[d]^-{\di{k}} \\
CA_2 \ar[r]^-{\di{f_2}} &  B_2}~\lower20pt\hbox{,}~\xymatrix{
CA_1 \ar[r]^-{\di{g_1}} \ar[d]_-{\di{h}} &  B_1\ar[d]^-{\di{k}} \\
CA_2 \ar[r]^-{\di{g_2}} &  B_2}$$
such that for neighbourhoods
$U_i \subseteq CA_i$ of $A_i \subseteq CA_i$ with $U_2=(1\times h)(U_1)$
$$g_i(t,a)~=~\begin{cases}
f_i(t,a)&{\it if}~ (t,a) \in U_i\\
*&{\it if}~ (t,a) \in \overline{CA_i \backslash U_i}
\end{cases}~~(i=1,2)$$
there is defined a commutative square
$$\xymatrix@C+15pt@R+15pt{
I \times \Sigma A_1 \ar[r]^-{\di{\gamma(f_1,g_1)}}
\ar[d]_-{\di{1\times\Sigma h}} &
B_1\ar[d]^-{\di{k}} \\
I \times \Sigma A_2 \ar[r]^-{\di{\gamma(f_2,g_2)}} &  B_2}$$
\end{proposition}
\begin{proof} (i) We have
$$\begin{array}{ll}
\delta(f,g)(t,a)
&=~\begin{cases}
g(1-2t,a)&{\rm if}~0 \leqslant t \leqslant 1/2\\
f(2t-1,a)&{\rm if}~1/2 \leqslant t \leqslant 1
\end{cases}\\[3ex]
&=~\begin{cases}
*&{\rm if}~0 \leqslant t \leqslant 1/2~{\rm and}~(1-2t,a) \in
\overline{CA \backslash U}\\
g(1-2t,a)&{\rm if}~0 \leqslant t \leqslant 1/2
~{\rm and}~(1-2t,a) \in U\\
g(2t-1,a)&{\rm if}~1/2 \leqslant t \leqslant 1~{\rm and}~(2t-1,a) \in U\\
f'(2t-1,a)&{\rm if}~1/2 \leqslant t \leqslant 1~{\rm and}~(2t-1,a) \in
\overline{CA \backslash U}~.
\end{cases}
\end{array}$$
Define a homotopy $\gamma(f,g):\delta(f,g)\simeq f'$ by
$$\begin{array}{l}
\gamma(f,g)~:~I \times \Sigma A \to X~;~
(s,t,a) \mapsto\\[2ex]
\begin{cases}
g((1-2s)(1-2t)+2s,a)&
\hbox{if}~0 \leqslant s,t \leqslant 1/2~\hbox{and}~(1-2t,a) \in U\\[1ex]
g((1-2s)(2t-1)+2s,a)&\hbox{if}~0 \leqslant s \leqslant 1/2 \leqslant
t \leqslant 1~
\hbox{and}~(2t-1,a) \in U\\[1ex]
f'(2t-1,a)&\hbox{if}~0 \leqslant s \leqslant 1/2~\hbox{and}~(2t-1,a) \in
\overline{CA \backslash U}\\[1ex]
f'((2-2s)(2t-1)+(2s-1)t,a)&\hbox{if}~1/2 \leqslant s,t \leqslant 1\\[1ex]
*&\hbox{otherwise~.}
\end{cases}
\end{array}$$
(ii) The homotopy $\gamma(f,g)$ constructed in (i) is natural.\\
\hfill\qed\end{proof}

\begin{example} {\rm Here are two special cases of Proposition \ref{near}.\\
(i) For any map $f:CA \to B$ take $f=g$, $U=CA$ to obtain a homotopy
$$\gamma(f,f)~:~\delta(f,f)~\simeq~*~:~\Sigma A \to B~.$$
(ii) For any map $g:CA \to B$ such that
$$g(0,a)~=~* \in B~~(a \in A)$$
take $f=*$, $U=A$ to obtain a homotopy
$$\gamma(*,g)~:~\delta(*,g)~\simeq~[g]~:~\Sigma A \to B~.$$}
\qed
\end{example}

As in Proposition \ref{BP1} for any pointed map $F:A \to B$ the homotopy cofibration
$$\xymatrix{A \ar[r]^-{\di{F}} & B  \ar[r]^-{\di{G}} & \Cc(F)  \ar[r]^-{\di{H}} & \Sigma A  \ar[r]^-{\di{\Sigma F}} & \Sigma B  \ar[r]&\dots}$$
induces a Barratt-Puppe exact sequence of  pointed sets
for any pointed space $X$
$$\xymatrix{\dots \ar[r]&
[\Sigma B,X] \ar[r]^-{\di{\Sigma F}} & [\Sigma A,X] \ar[r]^-{\di{H}} &
[\Cc(F),X] \ar[r]^-{\di{G}} & [B,X] \ar[r]^-{\di{F}} & [A,X]~.}$$
A map $q \cup r:\Cc(F) \to X$ is the same as a map $r:B \to X$ together
with a null-homotopy $q:rF \simeq *:A \to X$. Use the projection
$$\begin{array}{l}
\nabla~:~\Cc(F) \to \Cc(F) \vee \Sigma A~;\\[1ex]
\hskip25pt (t,x) \mapsto \begin{cases}
(2t,x) \in \Cc(F)&{\rm if}~ 0 \leqslant t \leqslant 1/2 \\
(2t-1,x) \in \Sigma A&{\rm if}~ 1/2 \leqslant t \leqslant 1~,
\end{cases}~~(x \in A)\\[1ex]
\hskip25pt y \mapsto y~~(y \in B)
\end{array}$$
to define a pairing
$$ [\Cc(F),X]\times [\Sigma A,X] \to [\Cc(F),X]~;~
(q \cup r,s) \mapsto ((q \cup r) \vee s)\nabla~.$$
\begin{proposition}~
Suppose that $F:A \to B$, $r:B \to X$ are maps and there are given two null-homotopies
$p,q:CA \to X$ of $rF:A \to X$
$$p,q~:~rF \simeq *~:~A \to X~.$$
The maps $p \cup r,q \cup r:\Cc(F) \to X$ are both extensions
of $r$, with
$$(p \cup r)G~=~(q \cup r)G~=~r \in [B,X]~,$$
and the rel $A$ difference
$$\delta(p,q)~:~\Sigma A \to X~;~(t,x) \mapsto
\begin{cases}
q(1-2t,x)&{\it if}~ 0 \leqslant t \leqslant 1/2 \\
p(2t-1,x)&{\it if}~ 1/2 \leqslant t \leqslant 1
\end{cases}$$
is such that there exists a homotopy
$$h~:~p \cup r~\simeq~ ((q \cup r) \vee \delta(p,q))\nabla~:~\Cc(F) \to X~.$$
Thus if $[\Cc(F),X]$ has a group structure (e.g. if $\Cc(F)$ is a suspension
or $X$ is an $H$-space)
$$\begin{array}{l}
p \cup r - q \cup r ~=~H(\delta(p,q))\\[1ex]
\in {\rm ker}(G:[\Cc(F),X] \to [B,X])~=~
{\rm im}(H:[\Sigma A,X] \to [\Cc(F),X])~.
\end{array}$$
\end{proposition}
\begin{proof} By construction
$$\begin{array}{l}
((q \cup r) \vee \delta(p,q))\nabla~:~\Cc(F) \to X~;\\[1ex]
\hskip25pt (t,x) \mapsto \begin{cases}
q(2t,x)&{\rm if}~0 \leqslant t \leqslant 1/2\\
q(3-4t,x)&{\rm if}~1/2 \leqslant t \leqslant 3/4\\
p(4t-3,x)&{\rm if}~3/4 \leqslant t \leqslant 1~,
\end{cases}~(x \in A)\\[1ex]
\hskip25pt y \mapsto y~(y \in B)~.
\end{array}$$
Define a homotopy  $h:p \cup r \simeq ((q \cup r) \vee \delta(p,q))\nabla$
by
$$\begin{array}{l}
h~:~I \wedge  \Cc(F) \to X~;\\[1ex]
\hskip25pt (s,(t,x)) \mapsto
\begin{cases} q(2st,x)&{\rm if}~0 \leqslant t \leqslant s/2\\[1ex]
q(3s-4t,x)&{\rm if}~s/2 \leqslant t \leqslant 3s/4\\[1ex]
p((4t-3s)/(4-3s),t)&{\rm if}~ 3s/4 \leqslant t \leqslant 1~,\end{cases}~(x \in A)\\[1ex]
\hskip25pt y \mapsto y~(y \in B)~.
\end{array}$$
\hfill\qed\end{proof}

\begin{proposition}~  \label{reldif2}
{\rm ~(i)} A commutative diagram of  spaces and maps
$$\xymatrix@C+20pt@R+20pt{
A \ar[r]^-{\di{F}} \ar[d]_-{\di{a}}
&B\ar[d]^-{\di{b}} \\
A' \ar[r]^-{\di{F'}} &B'}$$
induces a natural transformation of homotopy cofibration sequences
$$\xymatrix@C+20pt@R+20pt{
A \ar[r]^-{\di{F}} \ar[d]^-{\di{a}}& B \ar[r]^-{\di{G}} \ar[d]^-{\di{b}}&
\Cc(F) \ar[r]^-{\di{H}} \ar[d]^-{\di{(a,b)}}&
\Sigma A \ar[d]^-{\di{\Sigma a}} \ar[r]^-{\di{\Sigma F}} & \dots \\
A' \ar[r]^-{\di{F'}} & B' \ar[r]^-{\di{G'}} &
\Cc(F') \ar[r]^-{\di{H'}}& \Sigma A'\ar[r]^-{\di{\Sigma F'}}  & \dots}$$
{\rm (ii)} Suppose given a commutative diagram as in {\rm ~(i)}, with
the maps $a,b$  null-homotopic.  For any null-homotopies $\delta a:CA
\to A'$, $\delta b:CB \to B'$ of $a,b$ the map $(a,b)$ is determined up
to homotopy by the rel $A$ difference of null-homotopies of $F'a=bF:A
\to B'$,
$$f~=~(\delta b) C F~,~g~=~F'(\delta a)~:~CA \to B'$$
with
$$(a,b)~\simeq~G' \delta(f,g) H~:~\Cc(F) \to \Cc(F')~.$$
{\rm (iii)} Suppose given a commutative diagram as in {\rm (i)}, with
$b$  null-homotopic.  For any  null-homotopy $\delta b:CB \to B'$ of $b$
the map $(a,b)$ is homotopic to the rel $A$ difference of the
inclusion $i:CA \to \Cc(F')$ and the composite
$$j~=~G'(\delta b) C F~:~C A \to CB \to B' \to \Cc(F')~,$$
that is
$$(a,b) ~\simeq~\delta(j,i)\circ H~:~\Cc(F)\to \Cc(F')~.\eqno{\hbox{\qed}}$$
\end{proposition}

\section{Chain complexes}

The difference construction for maps of spaces induces
the difference construction for chain maps of chain complexes.

Let $R$ be a ring.

\begin{definition}~ \label{chaincone} {\rm
(i) The {\it suspension} of an $R$-module chain complex $C$ is the $R$-module
chain complex $SC$ with\index{suspension!chain complex, $SC$}
$$d_{SC}~=~d_C~:~(SC)_r~=~C_{r-1} \to (SC)_{r-1}~=~C_{r-2}~.$$
(ii) Let $f:C \to D$ be an $R$-module chain map.
The {\it relative difference} of chain homotopies\index{relative difference!chain
homotopies, $\delta(p,q)$}
$$p~:~f \simeq~0~:~C \to D~,~q~:~f \simeq~0~:~C \to D$$
is the $R$-module chain map
$$\delta(p,q)~=~q-p~:~SC \to D~.$$
(iii) The {\it algebraic mapping cone} of an
$R$-module chain map $f:C \to D$ is the $R$-module chain complex $\Cc(f)$ with
\index{algebraic mapping cone, $\Cc(f)$}
$$d_{\cc(f)}~=~\begin{pmatrix} d_D & (-)^rf \\ 0 & d_C \end{pmatrix}~:~
\Cc(f)_r~=~D_r \oplus C_{r-1} \to \Cc(f)_{r-1}~=~D_{r-1} \oplus C_{r-2}~.$$
The algebraic mapping cone fits into a short exact sequence
$$\xymatrix{0 \ar[r] &D \ar[r]^-{\di{g}} & \Cc(f) \ar[r]^-{\di{h}} & SC \ar[r] & 0}$$
with
$$\begin{array}{l}
g~=~\begin{pmatrix} 1 \\ 0  \end{pmatrix}~:~
D_r  \to \Cc(f)_r~=~D_r \oplus C_{r-1}~,\\[2ex]
h~=~\begin{pmatrix} 0 & 1  \end{pmatrix}~:~
\Cc(f)_r~=~D_r \oplus C_{r-1} \to SC_r~=~C_{r-1}~.
\end{array}$$
\hfill\qed}
\end{definition}

Given a space $X$ let $C(X)$ denote the singular $\ZZ$-module chain complex,
and for a pointed space let
$$\dot C(X)~=~C(X,\{*\})$$
be the reduced singular chain complex. The Eilenberg-Zilber theorem
gives a natural chain equivalence
$$S \dot C(X)~\simeq~\dot C(\Sigma X)~.$$
\begin{proposition}~ \label{chainreldif}
{\rm (i)} A map of spaces $F:X \to Y$ induces a chain map $f:C(X) \to C(Y)$.\\
{\rm (ii)} A homotopy $H:F \simeq G:X \to Y$ induces a chain homotopy
$$h~:~f~\simeq~g~:~C(X) \to C(Y)~.$$
{\rm (iii)} The reverse homotopy $-H:G\simeq F:X \to Y$
induces the reverse chain homotopy $-h:g \simeq f:C(X) \to C(Y)$
(up to a natural higher chain homotopy).\\
{\rm (iv)} The sum $H_1+H_2:F_1 \simeq F_3$ of homotopies
$$H_1~:~F_1~\simeq~F_2~,~H_2~:~F_2~\simeq~F_3~:~X \to Y$$
induces the sum chain homotopy
$$h_1+h_2~:~f_1~\simeq~f_3~:~C(X) \to C(Y)$$
(up to a natural higher chain homotopy).\\
{\rm (v)}  A pointed map $F:X \to Y$ induces a chain map $f:\dot C(X) \to \dot C(Y)$,
and similarly for chain homotopies.
The reduced singular chain complex of the mapping cone $\Cc(F)$
is chain equivalent to the algebraic mapping cone of $f$
$$\dot C(\Cc(F))~\simeq~ \Cc(f)~,$$
and the homotopy cofibration sequence of pointed spaces (\ref{cone})
$$\xymatrix@C+10pt{X \ar[r]^-{\di{F}} &
Y \ar[r]^-{\di{G}} &
\Cc(F)\ar[r]^-{\di{H}}& \Sigma X \ar[r]^-{\di{\Sigma F}} & \Sigma Y \ar[r] &\dots}$$
induces the sequence of chain complexes
$$\xymatrix@C+10pt{\dot C(X) \ar[r]^-{\di{f}} &
\dot C(Y) \ar[r]^-{\di{g}} &
\Cc(f)\ar[r]^-{\di{h}}& S \dot C(X) \ar[r]^-{\di{Sf}} & S \dot C(Y) \ar[r] &\dots}$$
with $g,h$ as in \ref{chaincone}.\\
{\rm (vi)} The relative difference $\delta(P,Q):\Sigma X \to Y$
of null-homotopies
$$P~:~F ~\simeq~ \{*\}~,~Q~:~F ~\simeq~ \{*\}~:~X \to Y$$
of a pointed map $F:X \to Y$ induces the
relative difference of the induced null-chain homotopies
$$p~:~f ~\simeq~ 0~, ~q~:~f ~\simeq~ 0~:~\dot C(X) \to \dot C(Y)$$
i.e. there is defined a commutative diagram
$$\xymatrix{S\dot C(X)\ar[dr]^-{\di{\simeq}}
\ar[rr]^-{\di{\delta(p,q)=q-p}} && \dot C(Y) \\
& \dot C(\Sigma X) \ar[ur]_-{\di{\delta(P,Q)}} &
}$$
\hfill\qed
\end{proposition}

\chapter{Umkehr maps and inner product spaces}\label{umkehr}

This Chapter develops a coordinate-free approach to stable homotopy theory,
using inner product spaces.

\section{Adjunction and compactification}

\begin{definition} {\rm
The {\it one-point compactification} $X^{\infty}$ of a locally compact Hausdorff space $X$ is the
space $X \cup \{\infty\}$ with open sets $U,
(X \backslash K)\cup\{\infty\} \subseteq X^{\infty}$
for open subsets $U \subseteq X$ and compact subsets $K \subseteq X$.
\index{one-point!compactification, $X^{\infty}$}
\\
\hfill\qed}
\end{definition}

The canonical pointed map $X^+ \to X^{\infty}$ from the one-point adjunction to the one-point compactification is a continuous bijection
which is a homeomorphism if and only if $X$ is compact. For any locally compact Hausdorff spaces $X,Y$
$$(X \times Y)^{\infty}~=~X^{\infty} \wedge Y^{\infty}~.$$
\indent For any pair of spaces $(A,B \subseteq A)$ the quotient $A/B$ is a
pointed space with base point $*=[B] \in A/B$.  The quotient space
$A/B$ is understood to be $A^+$ if $B=\emptyset$.  The projection
$\pi^+:A^+\to A/B$ is a map such that a subspace $U \subseteq
A/B$ is open if and only if $(\pi^+)^{-1}(U) \subseteq A^+$ is open.

If $A$ is locally compact Hausdorff, $B \subseteq A$ is closed and
 $A\backslash B \subseteq K \subseteq A$ for a compact subspace $K \subseteq A$ there is defined a commutative diagram
of pointed spaces and pointed maps
$$\xymatrix{A^\infty \ar[rr] \ar[dr]
&& (A\backslash B)^{\infty}\\
& A/B \ar[ur] &}$$
with
$$A^{\infty} \to (A \backslash B)^{\infty}~;~
a \mapsto \begin{cases} a&{\rm if}~a \in A\backslash B\\
\infty&{\rm if}~a=\infty~{\rm or}~ a \in B
\end{cases}$$
and $A/B \to (A\backslash B)^{\infty}$
a homeomorphism.

.

\begin{definition}~  \label{umkehradjunct}\label{umkehrcompact}
{\rm
(i) The {\it adjunction Umkehr map} of an injective map
$f:X \emb Y$ with respect to a subspace $Y_0 \subseteq Y\backslash f(X)$
is the projection \index{Umkehr!adjunction}
$$F^+~:~Y/Y_0 \to Y/(Y\backslash f(X))~.$$
(ii) An {\it open embedding} $f:X \emb Y$ is an injective map
such that a subset $U \subseteq X$ is open if and only if $f(U)
\subseteq Y$ is open. It follows that $f\vert:K \to f(K)$ is a
homeomorphism for any subspace $K \subseteq X$, giving $f(K)\subseteq Y$ the
subspace topology.\index{open!embedding}\\
(iii) The {\it compactification Umkehr} map of an open embedding $f$ as in (ii) with $X,Y$ locally compact Hausdorff is the surjective function\index{Umkehr!compactification}
$$F^{\infty}~:~Y^{\infty} \to X^{\infty}~;~y \mapsto
\begin{cases} x&{\rm if}~y=f(x) \in f(X) \\
\infty&{\rm otherwise}
\end{cases}$$
which preserves the base points and is continuous:
for any open subset $U\subseteq X$ and any closed compact subset $K \subseteq X$
$$\begin{array}{l}
(F^{\infty})^{-1}(U)~=~f(U) \subset Y^{\infty}~,\\[1ex]
(F^{\infty})^{-1}((X \backslash K)\cup\{\infty\})~=~
(Y \backslash f(K))\cup\{\infty\} \subseteq Y^{\infty}
\end{array}$$
with $f(U) \subseteq Y$ open and $f(K)\subseteq Y$ closed compact.\hfill\qed}
\end{definition}

\begin{example} {\rm Let $(X,\partial X)$, $(Y,\partial Y)$ be
$n$-dimensional
manifolds with boundary, and let $f:X \emb Y$ be an injective map such that
$$f(\partial X) \cap \partial Y~=~\emptyset~.$$
(i) The map
$$X/\partial X \to Y/(Y\backslash f(X\backslash \partial X))~;~[x] \mapsto [f(x)]$$
is a homeomorphism. The adjunction Umkehr map of $f$ with
respect to $Y_0=\partial Y \subset Y$ is the Pontrjagin-Thom map
$$F^+~:~Y/\partial Y \to Y/(Y\backslash f(X\backslash \partial X)))~=~X/\partial X~.$$
(ii) If $X,Y$ are compact the compactification Umkehr of the open embedding defined
by the restriction to the interiors
$$g~=~f\vert~:~X\backslash \partial X \emb Y \backslash \partial Y~,$$
is the adjunction Umkehr $F^+$ of (i)
$$G^{\infty}~=~F^+~:~(Y \backslash \partial Y)^{\infty}~=~Y/\partial Y
\to (X \backslash \partial X)^{\infty}~=~X/\partial X~.$$
\hfill\qed}
\end{example}

From now on we shall only be concerned with locally compact Hausdorff spaces.

\section{Inner product spaces}

Inner product spaces are given the metric topology.

\begin{terminology}~ {\rm
Given an inner product space $V$ write the unit disc,
unit sphere and projective space of $V$ as
\index{unit!disc, $D(V)$}\index{unit!sphere, $S(V)$}\index{projective space $P(V)$}
$$\begin{array}{l}
D(V)~=~\{v \in V\st \Vert v \Vert\leqslant  1\}~,\\[1ex]
S(V)~=~\{v \in V\st \Vert v \Vert=1\}~,\\[1ex]
P(V)~=~S(V)/\{v \sim -v\}~.
\end{array}$$
For finite-dimensional $V$ the spaces $D(V),S(V),P(V)$ are compact.}\\
\hfill\qed
\end{terminology}

\begin{proposition}~\label{sphere}
Let $V$ be an inner product space.\\
{\rm (i)} The maps
$$\begin{cases}
D(V)\backslash S(V) \to V~;~v \mapsto \di{v \over 1-\Vert v \Vert}\\[4ex]
V \to D(V)\backslash S(V)~;~v \mapsto \di{v \over 1 + \Vert v \Vert}
\end{cases}$$
are inverse homeomorphisms.\\
{\rm (ii)} The maps
$$\begin{cases}
D(V)/S(V) \to V^{\infty}~;~v \mapsto \di{v \over 1-\Vert v \Vert}\\[4ex]
V^{\infty} \to D(V)/S(V)~;~
v \mapsto \begin{cases} \di{v \over 1 + \Vert v \Vert}&
{for}~v \in V\\[1ex]
[S(V)]&{for}~v=\infty
\end{cases}
\end{cases}$$
are inverse bijections which are inverse homeomorphisms for finite-dimensional $V$.\\
{\rm (iii)} The maps
$$\begin{cases}
V^{\infty} \to sS(V)~;~v \mapsto
\begin{cases}
(0,S(V))&{for}~v=0\in V\\[1ex]
\big(\di{\Vert v\Vert \over 1 + \Vert v \Vert},{v \over \Vert v \Vert}\big)&
{for}~v \in V\backslash \{0\} \\[1ex]
(1,S(V))&{for}~v=\infty \in V^{\infty}
\end{cases}\\[4ex]
sS(V) \to V^{\infty}~;~(t,u)\mapsto [t,u]~=~\di{tu \over 1-t}
\end{cases}$$
are inverse bijections which are inverse homeomorphisms for finite-dimensional $V$.\\
{\rm (iv)} The map
$$S(V) \times \R \to V\backslash \{0\}~;~(v,x) \mapsto ve^x$$
is a homeomorphism, inducing a homeomorphism
$$(S(V) \times \R)^{\infty}~=~\Sigma S(V)^{\infty} \to (V\backslash \{0\})^{\infty}~.$$
Thus for finite-dimensional $V$
$$(V \backslash \{0\})^{\infty}~=~\Sigma S(V)^+~.$$
\hfill\qed
\end{proposition}

For $V=\R^k$ there are identifications
$$\begin{array}{l}
D(\R^k)~=~D^k~,~S(\R^k)~=~S^{k-1}~,\\[1ex]
(\R^k)^{\infty}~=~D^k/S^{k-1}~=~S^k~,~P(\R^k)~=~\R\P^{k-1}
\end{array}$$
and for any space $X$
$$\begin{array}{l}
(\R^k \times X)^{\infty}~=~S^k \wedge X^{\infty}~=~\Sigma^kX^{\infty}~,\\[1ex]
(D^k \times X)/(S^{k-1} \times X)~=~(D^k/S^{k-1}) \wedge X^+~=~\Sigma^kX^+~.
\end{array}$$

\begin{proposition}~\label{pushout1}
Let $V$ be an inner product space. \\
{\rm (i)} The map
$$S(V) \times  [0,1)\to V~;~(t,u) \mapsto [t,u]~=~\di{tu \over 1-t}$$
is a surjection, which induces a homeomorphism
$$(S(V) \times  [0,1))/(S(V) \times \{0\})~=~S(V)^+ \wedge [0,1) \to V$$
with $[0,1)$ based at $0$.\\
{\rm (ii)} For finite-dimensional $V$ there is defined a pushout diagram
$$\xymatrix@C+10pt@R+10pt{
S(V)^+ \ar[r]^-{\di{s_V}} \ar[d] &0^+\ar[d]^-{\di{0_V}} \\
C(S(V)^+) \ar[r] &V^{\infty}}$$
with
$$\begin{array}{l}
s_V~=~projection~:~S(V)^+ \to S(V)^+/S(V)~=~S^0~,\\[1ex]
C(S(V)^+)~=~S(V)^+ \wedge [0,1]~~(\hbox{$[0,1]$ based at $1$})~,\\[1ex]
S(V)^+ \to C(S(V)^+)~;~u \mapsto (0,u)~,\\[1ex]
C(S(V)^+) \to V^{\infty}~;~(t,u) \mapsto [t,u]~.
\end{array}$$
\end{proposition}
\begin{proof} By construction.
\hfill\qed\end{proof}

In view of the pushout square of Proposition \ref{pushout1} (ii) for
finite-dimensional $V$ a pointed map
$p:V^{\infty} \wedge X \to Y$ can be viewed as a map $r:X \to Y$ together
with a null-homotopy
$$q~:~r(s_V \wedge 1_X) ~\simeq~ \{*\}~:~S(V)^+\wedge X \to Y~.$$
More precisely:

\begin{proposition}~ \label{difcon6}{\it
For any pointed spaces $X,Y$ and a finite-dimensional inner product space $V$
there are natural one-one correspondences between
\begin{itemize}
\item[\rm (a)]~ maps $p:V^{\infty} \wedge X \to Y$,
\item[\rm (b)]~ maps $q:CS(V)^+\wedge X \to Y$ such that
$$q(0,u,x)~=~q(0,v,x) \in Y~(u,v \in S(V),x \in X)~,$$
\item[\rm (c)]~ maps $r:X \to Y$ together with a null-homotopy
$$\delta r ~:~CS(V)^+\wedge X \to Y$$
of
$$r(s_V \wedge 1_X)~:~S(V)^+ \wedge X \to Y~;~(u,x) \mapsto r(x)~.$$
\end{itemize}}
\end{proposition}
\begin{proof}
(a) $\Longrightarrow$ (b) Given a map $p:V^{\infty} \wedge X \to Y$
define a null-homotopy
$$\delta p ~:~CS(V)^+\wedge X \to Y~;~(t,u,x) \mapsto p([t,u],x)$$
of
$$(p\vert )(s_V \wedge 1_X)~:~S(V)^+ \wedge X \to Y~;~(u,x) \mapsto p(0,x)~.$$
(b) $\Longrightarrow$ (a) Given a map $q:CS(V)^+\wedge X \to Y$ such that
$$q(0,u,x)~=~q(0,v,x) \in Y~(u,v \in S(V),x \in X)$$
define a map
$$p~:~V^{\infty} \wedge X \to Y~;~([t,u],x) \mapsto q([t,u],x)$$
with $q=\delta p$.\\
(b) $\Longrightarrow$ (c) Given $q$ define
$$r~:~X \to Y~;~x \mapsto q([ 0,u],x) ~(\hbox{\rm for any}~u \in S(V))~.$$
(c) $\Longrightarrow$ (b) Given $r,\delta r$ define $q = \delta r$.\\
\hfill\qed\end{proof}

\begin{proposition}~ \label{pushout2}
{\rm (i)} For any inner product spaces $U,V$
$$\begin{array}{l}
D(U \oplus V)~=~D(U) \times D(V)~,\\[1ex]
S(U \oplus V)~=~S(U) \times D(V) \cup D(U) \times S(V)
\end{array}$$
and there is defined a homeomorphism
$$\lambda_{U,V}~:~S(U)*S(V) \to S(U \oplus V)~;~
(t,u,v) \mapsto \big(u\cos(\pi t/2),v\sin(\pi t/2)\big)$$
with inverse
$$\lambda_{U,V}^{-1}~:~S(U \oplus V)\to S(U)*S(V) ~;~
(x,y) \mapsto
\bigg(\dfrac{2}{\pi}\tan^{-1}\big(\dfrac{\Vert y \Vert}{\Vert x \Vert}\big),
\dfrac{x}{\Vert x \Vert},
\dfrac{y}{\Vert y \Vert}\bigg)~.$$
{\rm (ii)} For finite-dimensional $U,V$ there is
a homotopy cofibration sequence of pointed spaces
$$\xymatrix{
S(V)^+ \ar[r]^-{\di{i_{U,V}}}&
S(U\oplus V)^+ \ar[r]^-{\di{j_{U,V}}} & S(U)^+\wedge V^{\infty}
\ar[r]^-{\di{k_{U,V}}} & \Sigma S(V)^+}$$
and a pushout diagram of pointed spaces
$$\xymatrix@C+10pt@R+10pt{
S(U \oplus V)^+ \ar[r]^-{\di{j_{U,V}}} \ar[d] &S(U)^+ \wedge V^{\infty}
\ar[d]^-{\di{k_{U,V}}} \\
CS(U \oplus V)^+ \ar[r]^-{\di{\delta j_{U,V}}} &\Sigma S(V)^+}$$
with
$$\begin{array}{l}
i_{U,V}~=~{\rm inclusion}~:~S(V)^+ \to S(U \oplus V)^+~;~v \mapsto (0,v)~,\\[1ex]
j_{U,V}~=~{\rm projection}~:~S(U \oplus V)^+~=~(S(U) * S(V))^+\\[1ex]
\hskip25pt \to
S(U \oplus V)/(S(V) \times D(U))~=~(D(V) \times S(U)/S(V) \times S(U))\\[1ex]
\hskip100pt =~S(U)^+ \wedge V^{\infty}~;~(t,u,v) \mapsto (u,[t,v])~,\\[1ex]
\delta j_{U,V}~:~CS(U \oplus V)^+ \to \Sigma S(V)^+~;~
(s,(t,u,v)) \mapsto (s+(1-s)t,v)~,\\[1ex]
k_{U,V}~:~S(U)^+\wedge V^{\infty} \to \Sigma S(V)^+~;~(u,[t,v]) \mapsto (t,v)~.
\end{array}$$
{\rm (iii)} For finite-dimensional $V$
there is a natural one-one correspondence between the homotopy classes of
pointed maps $f:\Sigma S(V)^+ \wedge X \to Y$ and the homotopy classes of
pairs $(g,h)$ defined by a pointed map $g:S(U)^+ \wedge V^{\infty} \wedge X \to Y$ and a null-homotopy
$$h~:~g(j_{U,V} \wedge 1_X)~\simeq~\{*\}~:~S(U \oplus V)^+ \wedge X \to Y~.$$
\end{proposition}
\begin{proof} (i)+(ii) By construction.\\
(iii) By (ii) there is defined a pushout square
$$\xymatrix@C+10pt@R+10pt{
S(U \oplus V)^+\wedge X \ar[r]^-{\di{j_{U,V}\wedge 1_X}} \ar[d] &S(U)^+ \wedge V^{\infty}\wedge X
\ar[d]^-{\di{k_{U,V}\wedge 1_X}} \\
CS(U \oplus V)^+ \wedge X\ar[r]^-{\di{\delta j_{U,V}\wedge 1_X}} &\Sigma S(V)^+\wedge X}$$
\hfill\qed\end{proof}

\section{The addition and subtraction of maps}

The compactification Umkehr construction is now used to define the
addition and subtraction of maps $F:V^{\infty} \wedge X \to Y$
for a non-zero inner product space $V$.

In dealing with open embeddings $e:V \emb V$ we shall only be concerned
with ones (e.g.  smooth) such that the compactification Umkehr map is a
homotopy equivalence $F:V^{\infty} \to V^{\infty}$.

\begin{definition}~ \label{sum-map}{\rm
(i) A homotopy equivalence $F:V^{\infty} \to V^{\infty}$ is\\
$\begin{cases} \hbox{\it orientation-preserving} \\
\hbox{\it orientation-reversing}
\end{cases}$ if\index{orientation-preserving!homotopy equivalence}
\index{orientation-reversing!homotopy equivalence}\\
$${\rm degree}(F)~=~\begin{cases} +1 \\ -1~.
\end{cases}$$
(ii) An embedding $e:V \emb V$ is
$\begin{cases} \hbox{\it orientation-preserving} \\
\hbox{\it orientation-reversing}
\end{cases}$
if the compactification Umkehr map $F^{\infty}:V^{\infty} \to V^{\infty}$
is an $\begin{cases} \hbox{orientation-preserving} \\
\hbox{orientation-reversing}\end{cases}$ homotopy equivalence.
\index{orientation-preserving!embedding}
\index{orientation-reversing!embedding}\\
(iii) A {\it sum map} for $V$ is any map\index{sum map, $\nabla_V$}
$$\nabla_V~:~V^{\infty} \to V^{\infty} \vee V^{\infty}$$
in the homotopy class of the compactification Umkehr map
of an open embedding $e:V \times \{1,2\} \emb V$
such that the restrictions $e\vert: V\times \{i\} \emb V$ ($i=1,2$)
are orientation-preserving. Then $\nabla_V$ is a sum map for the pointed space $V^{\infty}$ in the sense of Definition \ref{sumdifference} (i).\\
(iv) A {\it difference map} for $V$ is a map
\index{difference map!$\overline{\nabla}_V$}
$$\overline{\nabla}_V~:~V^{\infty} \to V^{\infty} \vee V^{\infty}$$
in the homotopy class of the compactification Umkehr map of
an open embedding $\overline{e}:V \times \{1,2\} \emb V$
such that the restriction $\overline{e}\vert: V \times \{i\} \emb V$
is orientation-reversing for $i=1$ and orientation-preserving for $i=2$.
Then $\overline{\nabla}_V$ is a difference map for the pointed space $V^{\infty}$ in the sense of
Definition \ref{sumdifference} (iii).
\hfill\qed}
\end{definition}

\begin{remark} {\rm For $i=1,2$ let
$$\pi_i~:~V^{\infty} \vee V^{\infty} \to V^{\infty}~;~
v_j \mapsto \begin{cases} v_j&{\rm if}~i=j \\
\infty&{\rm if}~i \neq j
\end{cases}~(j=1,2)~.$$
(i) If ${\rm dim}(V) \geqslant 2$ a map $\nabla_V:V^{\infty} \to
V^{\infty} \vee V^{\infty}$ is a sum map if and only if the
maps $\pi_1\nabla_V,\pi_2\nabla_V:V^{\infty} \to V^{\infty}$ are
orientation-preserving homotopy equivalences.\\
(ii) If ${\rm dim}(V) \geqslant 2$ a map $\overline{\nabla}_V:V^{\infty} \to
V^{\infty} \vee V^{\infty}$ is a difference map if and only if the
maps $\pi_1\overline{\nabla}_V,\pi_2\overline{\nabla}_V:V^{\infty} \to V^{\infty}$ are
homotopy equivalences with $\pi_1\nabla_V$ orientation-reversing and
$\pi_2\nabla_V$ orientation-preserving.
\hfill\qed}
\end{remark}

\begin{example} \label{sumdif} {\rm
(i) For any $x \in V$ there is defined an orientation-preserving open embedding
$$e_x~:~V \to V~;~v \mapsto x+ \di{v \over 1+\Vert v \Vert }$$
such that
$$e_x(0)~=~x \in V~,~e_x(V)~=~\{v \in V\,\vert\, \Vert v-x \Vert < 1\}~.$$
(ii) For any $y \in S(V)$ let
$$L_y~=~\{t y\,\vert\, t \in \R\}~,~
H_y~=~\{z \in V\,\vert\, \langle y,z \rangle=0 \in \R\} \subset V~,$$
so that $H_y$ the hyperplane  orthogonal to the line $L_y \subset V$
containing $y$.
Reflection in $H_y$ defines an orientation-reversing isomorphism
$$\begin{array}{l}
R_y~=~-1_{L_y} \oplus 1_{H_y}~:~V~=~L_y \oplus H_y \to
V~=~L_y \oplus H_y~;\\[2ex]
\hspace*{100pt} v=(ty,z) \mapsto v - 2 \langle v,y \rangle y=(-ty,z)~.
\end{array}$$
For any $x\in V$ the composite $\overline{e}_{x,y}=e_xR_y:V \to V$ is an
orientation-reversing open embedding such that
$$\overline{e}_{x,y}(0)~=~x \in V~,~\overline{e}_{x,y}(V)~=~
\{v \in V\,\vert\, \Vert v-x \Vert < 1\}~.$$
(iii) For any $v \in V$ with $\Vert v \Vert =1$ define an open embedding
$$f_v~:~V \times \{v,-v\} \emb V~;~
(u,\pm v) \mapsto \pm v+ \dfrac{u}{1+\Vert u \Vert}$$
with compactification Umkehr map a sum map
$$\begin{array}{l}
F_v~=~\nabla_v~:~V^{\infty} \to (V \times
\{v,-v\})^{\infty}~=~V^{\infty} \vee V^{\infty}~;\\[1ex]
\hspace*{150pt}
u \mapsto
\begin{cases}
\bigg(\di{u-v \over 1- \Vert u-v \Vert}\bigg)_1&
{\rm if}~ \Vert u-v \Vert <1\\[2ex]
\bigg(\di{u+v \over 1- \Vert u+v \Vert}\bigg)_2&
{\rm if}~ \Vert u+v \Vert <1\\[2ex]
*&{\rm otherwise}
\end{cases}
\end{array}$$
such that $F_{-v}=TF_v$ with $T:V^{\infty} \vee V^{\infty} \to
V^{\infty} \vee V^{\infty}$ the transposition involution.\\
(iv) If $\nabla_V:V^{\infty} \to V^{\infty} \vee V^{\infty}$ is a sum map
then
$$\overline{\nabla}_V~=~(R_y \vee 1)\nabla_V~:~V^{\infty} \to
V^{\infty} \vee V^{\infty}$$
is a difference map, with $R_y$ as in (ii). \\
(v) For a decomposition $V=\R \oplus W$ as in (ii) use the homeomorphism
$$V^{\infty}\xymatrix{\ar[r]^-{\di{\cong}}&}  \Sigma W^{\infty} ~;~(t,w) \mapsto
(\di{ \frac{e^t}{e^t+1}},w)$$
to identify $V^{\infty}=\Sigma W^{\infty}$ and define sum and difference maps
$$\begin{array}{l}
\nabla_V~:~V^{\infty}\to V^{\infty}\vee V^{\infty}~;~
(t,w) \mapsto \begin{cases}
(2t,w)_1&\hbox{if $0 \leqslant  t \leqslant  1/2$}\\[1ex]
(2t-1,w)_2&\hbox{if $1/2 \leqslant  t \leqslant  1$~,}
\end{cases}\\[3ex]
\overline{\nabla}_V~:~V^{\infty}\to V^{\infty}\vee V^{\infty}~;~
(t,w) \mapsto \begin{cases}
(1-2t,w)_1&\hbox{if $0 \leqslant  t \leqslant  1/2$}\\[1ex]
(2t-1,w)_2&\hbox{if $1/2 \leqslant  t \leqslant  1$~.}
\end{cases}
\end{array}$$
\hfill\qed}
\end{example}

\begin{proposition}~ \label{cofib1}
For any pointed spaces $X,Y$ and an inner product space $V$ the set
$[V^{\infty} \wedge X,Y]$ of homotopy classes of maps $V^{\infty}
\wedge X \to Y$ has a group structure for ${\rm dim}(V) \geqslant 1$
{\rm (}abelian for ${\rm dim}(V) \geqslant 2${\rm )}, with the sum
and difference of $f,g:V^{\infty}\wedge X \to Y$ given by
$$\begin{array}{l}
f+g~=~(g \vee f)(\nabla_V\wedge 1_X)~:~V^{\infty}\wedge X \to Y~,\\[1ex]
f-g~=~(g \vee f)(\overline{\nabla}_V\wedge 1_X)~:~V^{\infty}\wedge X \to Y~.
\end{array}$$
\hfill\qed
\end{proposition}

\begin{example} {\rm
(i) The map
$$\begin{array}{l}
\nabla_V~:~V^{\infty} \to V^{\infty} \vee V^{\infty}~;\\[1ex]
[t,u]=v \mapsto  \begin{cases}
[2t,u]_1=\bigg(\di{ 2v \over 1- \Vert v \Vert }\bigg)_1&
{\rm if}~0 \leqslant t \leqslant 1/2\\[2ex]
[2t-1,u]_2=\bigg(
\di{(\Vert v \Vert  -1)v \over 2 \Vert v \Vert }\bigg)_2
&{\rm if}~1/2 \leqslant t \leqslant 1
\end{cases}
\end{array}$$
is a sum map such that $\nabla_V(0)=0_1$, $\nabla_V(\infty)=\infty_2$.\\
(ii) The composite of the projection $V^{\infty} \to V^{\infty}/S(V)$
and the homeomorphism
$$\begin{array}{l}
V^{\infty}/S(V) \to V^{\infty} \vee V^{\infty}~;\\[1ex]
[t,u]=v \mapsto  \begin{cases}
[1-2t,u]_1=\bigg(\di{ (1- \Vert v \Vert)v \over 2 \Vert v \Vert }\bigg)_1&
{\rm if}~0 \leqslant t \leqslant 1/2\\[2ex]
[2t-1,u]_2=\bigg(
\di{(\Vert v \Vert  -1)v \over 2 \Vert v \Vert }\bigg)_2
&{\rm if}~1/2 \leqslant t \leqslant 1
\end{cases}
\end{array}$$
is a difference map
$$\begin{array}{l}
\overline{\nabla}_V~:~V^{\infty}\to V^{\infty}/S(V)~\cong~V^{\infty}\vee V^{\infty}~;~
[t,u] \mapsto \begin{cases}
[1-2t,u]_1&\hbox{if $0 \leqslant  t \leqslant  1/2$}\\[1ex]
[2t-1,u]_2&\hbox{if $1/2 \leqslant  t \leqslant  1$}
\end{cases}
\end{array}$$
such that $\overline{\nabla}_V(0)=\overline{\nabla}_V(\infty)=\infty$.\\
\hfill\qed}
\end{example}

\begin{proposition}~ \label{cofib2}
Let $V$ be a non-zero finite-dimensional inner product space.\\
{\rm (i)} There is defined a commutative diagram of maps and spaces
$$\xymatrix{
\R \times S(V) \ar[r] \ar[d]^-{\di{\cong}}
& S(\R \oplus V) \ar[r] \ar[d]^-{\di{\cong}}
& \Sigma S(V)^+ \ar[r] \ar[d]^-{\di{\cong}} & \Sigma S(V)^+/S(V)
\ar[d]^-{\di{\cong}} \\
V\backslash \{0\} \ar[r] & V^{\infty} \ar[r]& V^{\infty}/0^+ \ar[r] &
V^{\infty} \vee_0 V^{\infty}}$$
with the vertical maps homeomorphisms, and
$$V^{\infty} \vee_0 V^{\infty}~=~(V^{\infty} \vee V^{\infty})
/\{0_1\sim 0_2\}~.$$
{\rm (ii)} There is defined a commutative square of homeomorphisms
$$\xymatrix{(\R \times S(V))^{\infty}
\ar[r]^-{\di{\cong}} \ar[d]_-{\di{\cong}} &
\Sigma S(V)^+ \ar[d]^-{\di{\cong}} \\
(V \backslash \{0\})^{\infty} \ar[r]^-{\di{\cong}}
& V^{\infty}/0^+}$$
{\rm (iii)} The maps $\alpha_V,\beta_V,\gamma_V,s_V,0_V$ defined by
\index{$\alpha_V$} \index{$\beta_V$} \index{$\gamma_V$}
\index{$s_V$}\index{$z_V$}
$$\begin{array}{l}
\alpha_V~:~V^{\infty} \to \Sigma S(V)^+~;~[t,u] \mapsto (t,u)~,\\[1ex]
\beta_V~:~V^{\infty} \to V^{\infty} \vee_0 V^{\infty}~;\\[1ex]
\hphantom{\beta_V~:~}
[t,u]=v \mapsto  \begin{cases}
[ 1-2t,u]_1=\bigg(\di{1- \Vert v \Vert \over 2 \Vert v \Vert}v\bigg)_1
&\hbox{if $0 \leqslant t \leqslant 1/2$} \\[2ex]
[2t-1,u]_2=\bigg(\di{\Vert v \Vert -1 \over 2 \Vert v \Vert}v\bigg)_2
&\hbox{if $1/2 \leqslant t \leqslant 1$~,}
\end{cases}\\[7ex]
\gamma_V~:~\Sigma S(V)^+ \to \Sigma S(V)^+/S(V)~\cong~
V^{\infty} \vee_0 V^{\infty}~;\\[1ex]
\hphantom{\beta_V~:~} (t,u) \mapsto
\begin{cases}
[1-2t,u]_1&\hbox{if $0 \leqslant t \leqslant 1/2$}\\
[2t-1,u]_2&\hbox{if $1/2 \leqslant t \leqslant 1$~,}\end{cases}\\[3ex]
s_V~=~{\rm projection}~:~
S(V)^+ \to S(V)^+/S(V)~=~S^0~=~\{0,\infty\}~;\\[1ex]
\hskip150pt u \mapsto 0~(u \in S(V))~,\\[1ex]
0_V~:~S^0 \to V^{\infty}~;~ \infty \mapsto \infty~,~0 \mapsto 0
\end{array}$$
fit into a commutative braid of homotopy cofibrations of unpointed spaces
$$\xymatrix@C-20pt{
&S^0\ar[dr] \ar@/^2pc/[rr]^-{\di{0_V}} && V^{\infty}
\ar[dr]^-{\di{\alpha_V}}
\ar@/^2pc/[rr]^-{\di{\beta_V}} &&V^{\infty} \vee_0 V^{\infty}\\
S(V)^+ \ar[ur]^-{\di{s_V}} \ar[dr]&&S(V) \sqcup \{0,\infty\}\ar[ur] \ar[dr] &&
\Sigma S(V)^+\ar[ur]^-{\di{\gamma_V}} \ar[dr]&&\\
&(S(V) \vee  S(V))^+ \ar[ur] \ar@/_2pc/[rr]
&&S(V)\ar@/_2pc/[rr]\ar[ur]^-{\di{*}} &&S^1}$$

\medskip

\noindent {\rm (}which extends to a braid of homotopy cofibrations of
pointed spaces on the right{\rm )}. For any $u \in S(V)$
$$\begin{array}{l}
\alpha_V(u)~=~\alpha_V[1/2,u]~=~(1/2,u) \in \Sigma S(V)^+~,\\[1ex]
\beta_V(u)~=~\beta_V[1/2,u]~=~\gamma_V(1/2,u)~=~0_1~=~0_2 \in V^{\infty} \vee_0 V^{\infty}~.
\end{array}$$
The map $\alpha_V:V^{\infty} \to (\R\times
S(V))^{\infty}~\cong~\Sigma S(V)^+$ is the compactification Umkehr
map of the open embedding
$$\R\times S(V) \emb V~;~(x,u) \mapsto e^x u~,$$
and $\gamma_V$ pinches $S(V) \subset \Sigma S(V)^+$ to
$0 \in V^{\infty} \vee_0 V^{\infty}$.\\
{\rm (iii)} There exists a homotopy
$$\beta_V~\simeq~i\overline{\nabla}_V~:~V^{\infty} \to
V^{\infty}\vee_0V^{\infty}$$
with $\overline{\nabla}_V:V^{\infty} \to V^{\infty} \vee V^{\infty}$
a difference map for $V$, and $i:V^{\infty} \vee
V^{\infty} \to V^{\infty} \vee_0 V^{\infty}$ the projection.\\
{\rm (iv)} For finite-dimensional inner product spaces $U,V$
there is defined a commutative braid of homotopy cofibrations of pointed spaces
$$\xymatrix@C-30pt{
S(V)^+\ar[dr]^-{\di{i_{U,V}}} \ar@/^2pc/[rr]^-{\di{s_V}} &&
S^0 \ar[dr]^-{\di{0_V}}
\ar@/^2pc/[rr]^-{\di{0_{U \oplus V}}} &&
(U\oplus V)^{\infty}=U^{\infty} \wedge V^{\infty}
\ar[dr]^-{\di{\alpha_{U \oplus V}}} & \\
&S(U\oplus V)^+\ar[ur]^-{\di{s_{U\oplus V}}} \ar[dr]^-{\di{j_{U,V}}} &&
~V^{\infty}~
\ar[ur]\ar[dr]^-{\di{\alpha_V}}&&\Sigma S(U \oplus V)^+&\\
&& S(U)^+\wedge V^{\infty}\ar@/_2pc/[rr]^-{\di{k_{U,V}}} \ar[ur] &&\Sigma S(V)^+
\ar[ur] &}$$

\noindent
\end{proposition}
\begin{proof} (i) The projection
$$V^{\infty} \to (V\backslash \{0\})^{\infty}~=~V^{\infty}/0^+$$
is the compactification
Umkehr map of the open embedding $V\backslash \{0\} \emb V$, and
$$\begin{array}{l}
\R \times S(V) \iso V \backslash \{0\}~;~
(s,u) \mapsto e^su~,\\[1ex]
\R \times S(V) \xymatrix{\ar[r]&} S(\R \oplus V)~;~
(s,u) \mapsto \bigg(\di{e^{-s}- e^s \over e^{-s}+e^s},
\di{2u \over e^{-s}+e^s}\bigg)~,
\end{array}$$
and the other maps defined as follows.
Stereographic projection defines inverse homeomorphisms
$$\begin{array}{l}
S(\R \oplus V)\iso  V^{\infty} ~;~
(t,v) \mapsto \di{v \over 1-t}~,~(1,0) \mapsto \infty~,\\[2ex]
V^{\infty}\iso S(\R \oplus V)~;~
w \mapsto \di{\bigg({\Vert w
\Vert^2 -1 \over \Vert w \Vert^2 +1},{2 w \over \Vert w\Vert^2+1
}\bigg)}~,~ \infty \mapsto (1,0)~.
\end{array}$$
For any $t \in I$, $u \in S(V)^+$ define $[t,u] \in V^{\infty}$ by
$$[t,u]~=~\begin{cases}\di{tu \over 1-t}&\hbox{if $t <1$ and $u \neq \infty$}
\\[1ex]
\infty&\hbox{if $t=1$ or $u=\infty$~.}\end{cases}$$
Every $v\neq 0 \in V$ has a unique expression as $v=[t,u]$ with
$$t~=~\di{\Vert v \Vert \over 1+ \Vert v \Vert} \in (0,1)~,~
u~=~\di{v \over \Vert v \Vert} \in S(V)~.$$
The projection
$$I \wedge S(V)^+\to V^{\infty}~;~(t,u) \mapsto [t,u]$$
pinches $\{0\} \times S(V)$ to $0 \subset V^{\infty}$, inducing
the homeomorphism
$$\Sigma S(V)^+ \iso V^{\infty}/0^+~;~(t,u) \mapsto [t,u]$$
with inverse
$$V^{\infty}/0^+ \iso \Sigma S(V)^+~;~v \mapsto
\begin{cases}
\bigg(\di{\Vert v \Vert \over 1+\Vert v \Vert},\di{v \over \Vert v
\Vert}\bigg)&{\rm if}~v \neq 0 \\
\infty&{\rm if}~v=0~. \end{cases}$$
The elements of
$V^{\infty}$ can thus be written as $[t,u]$ with $0 \leqslant  t
\leqslant  1$, $u \in S(V)$, identifying
$$0~=~[ 0,u]~,~u~=~[ 1/2,u]~,~\infty~=~[ 1,u] \in V^{\infty}~,$$
and
$$S(V)~=~\{1/2\} \times S(V) \subset \Sigma S(V)^+~.$$
The homeomorphism $\Sigma S(V)^+/S(V)\cong V^{\infty}\vee_0V^{\infty}$
is defined by
$$\Sigma S(V)^+/S(V) \iso V^{\infty}\vee_0V^{\infty}~;~
(t,u) \mapsto \begin{cases}
[1-2t,u]_1&\hbox{if}~0 \leqslant t \leqslant 1/2\\
[2t-1,u]_2&\hbox{if}~1/2 \leqslant t \leqslant 1~.
\end{cases}$$
{\rm (ii)} Define
$$\begin{array}{l}
(\R \times S(V))^{\infty} \iso
(V \backslash \{0\})^{\infty}~;~(x,u) \mapsto e^x u~, \\[1ex]
(\R \times S(V))^{\infty} \iso \Sigma S(V)^+ ~;~(x,u) \mapsto
\bigg(\di{e^x \over 1+e^x},u\bigg)~.
\end{array}$$
(iii) Let $e:V \times \{1,2\} \emb V$ be an open embedding such that
\begin{enumerate}
\item $0 < \Vert e(0,1) \Vert < 1 < \Vert e(0,2) \Vert$,
\item $e(V \times \{1,2\}) \cap (\{0\}\cup S(V))=\emptyset$,
\item $e\vert:V \times \{1\} \emb V$ is orientation-reversing,
\item $e\vert:V \times \{2\} \emb V$ is orientation-preserving.
\end{enumerate}
The compactification Umkehr of $e$ is a difference map for $V$ (\ref{sum-map} (iv))
$$F~=~\overline{\nabla}_V~:~V^{\infty} \to V^{\infty}\vee V^{\infty}$$
such that
\begin{enumerate}
\item $F(0)=F(\infty)=\infty$,
\item $F(e(0,i))=0_i$ ($i=1,2$),
\item the induced map
$$[F]~:~V^{\infty}/S(V) \to V^{\infty}\vee V^{\infty}$$
is a homotopy equivalence,
\item the induced map
$$[F]~:~(V^{\infty}/S(V))/\{e(0,1)\sim e(0,2)\} \to V^{\infty} \vee_0
  V^{\infty}$$
is a homotopy equivalence.
\end{enumerate}
The composite
$$i\overline{\nabla}_V~:~V^{\infty} \to
(V^{\infty}/S(V))/\{e(0,1)\sim e(0,2)\}
\xymatrix{\ar[r]^-{\di{[F]}}&} V^{\infty} \vee_0
  V^{\infty}$$
is homotopic to
$$V^{\infty} \to (V^{\infty}/S(V))/\{0 \sim
\infty\}~=~(V^{\infty}/0^+)/S(V) \iso V^{\infty} \vee_0
    V^{\infty}$$
which is just $\beta_V$.\\
(iv) By construction.
\hfill\qed\end{proof}

\begin{definition}~ \label{difcon1}
{\rm  Let $p,q:V^{\infty} \wedge X \to Y$ be maps which agree on
$X=0^+ \wedge X\subset V^{\infty} \wedge X$, i.e. such that
$$p(0,x)~=~q(0,x) \in Y~(x \in X)~.$$
(i) The {\it difference map} of $p$ and $q$ is the map
\index{difference map of stable maps!$p-q$}
$$p-q~=~(q \vee p)(\overline{\nabla}_V\wedge 1_X)~:~
V^{\infty} \wedge X  \to Y~,$$
that is
$$p-q~:~V^{\infty} \wedge X \to Y~;~
([t,u],x) \mapsto \begin{cases}q([ 1-2t,u],x)&\hbox{if $0 \leqslant  t
\leqslant  1/2$}\\[1ex]
p([ 2t-1,u],x)&\hbox{if $1/2 \leqslant  t \leqslant 1$}~.
\end{cases}$$
(ii) The {\it relative difference of stable maps, $\delta(p,q)$} of $p$ and $q$ is
the rel $X$ difference
\index{relative difference!stable maps, $\delta(p,q)$}
$$\delta(p,q)~:~\Sigma S(V)^+ \wedge X \to Y~;~
(t,u,x) \mapsto \begin{cases}
q([ 1-2t,u],x)&\hbox{if $0 \leqslant  t \leqslant 1/2$}\\[1ex]
p([ 2t-1,u],x)&\hbox{if $1/2 \leqslant  t \leqslant 1$~,}
\end{cases}$$
that is
$$\delta(p,q)~=~(q \vee_0 p)(\gamma_V\wedge 1_X)~:~
\Sigma S(V)^+ \wedge X \to (V^{\infty} \vee_0 V^{\infty}) \wedge X \to Y~.$$
\hfill\qed
}\end{definition}

\begin{example} {\rm The relative difference of maps
$p,q:V^{\infty} \wedge X \to Y$ such that
$$p(0,x)~=~q(0,x)~=~* \in Y~(x\in X)$$
is just the difference of the induced maps
$$[p],[q]~:~\Sigma S(V)^+\wedge X~=~(V^{\infty}/0^+) \wedge X \to Y~,$$
that is
$$\delta(p,q)~=~[p]-[q]~:~\Sigma S(V)^+ \wedge X \to Y~,$$
with $[p]-[q]=([q] \vee [p])\overline{\nabla}$ defined using the difference map
$$\overline{\nabla}~:~\Sigma S(V)^+ \to \Sigma S(V)^+\vee\Sigma S(V)^+~;~
(t,u) \mapsto
\begin{cases}
(1-2t,u)&{\rm if}~0 \leqslant t \leqslant 1/2 \\
(2t-1,u)&{\rm if}~1/2 \leqslant t \leqslant 1~.
\end{cases}$$
\hfill\qed}
\end{example}

The difference map $p-q$ of Definition \ref{difcon1} (i) is
just an explicit representative of the difference of the homotopy classes
$p,q \in [V^{\infty}\wedge X,Y]$.
Writing $[t,u]=v \in V^{\infty}$ the difference of
$p,q:V^{\infty}\wedge X \to Y$ can be expressed as
$$p-q~:~V^{\infty} \wedge X \to Y~;~(v,x) \mapsto \begin{cases}
q(\di{1- \Vert v \Vert \over 2 \Vert v \Vert^2}v,x)
&{\rm if}~0 < \Vert v \Vert < 1 \\[2ex]
p(0,x)=q(0,x)&{\rm if}~\Vert v \Vert = 1 \\[2ex]
p(\di{\Vert v \Vert -1 \over 2 \Vert v \Vert}v,x)
&{\rm if}~\Vert v \Vert > 1\\
\infty&{\rm if}~v=0~{\rm or}~*~.
\end{cases}$$
The difference is the composite
$$p-q~:~V^{\infty} \wedge X \to (V^{\infty}/0^+) \wedge X~
\cong~ \Sigma S(V)^+ \wedge X
\xymatrix@C+10pt{\ar[r]^-{\di{\delta(p,q)}}&} Y~,$$
and the same formulae apply to the relative difference~:
 $$\delta(p,q)~:~(V^{\infty}/0^+) \wedge X
 \to Y~;~(v,x) \mapsto \begin{cases}
q(\di{1- \Vert v \Vert \over 2 \Vert v \Vert^2}v,x)
&{\rm if}~0 < \Vert v \Vert \leqslant 1 \\[2ex]
p(\di{\Vert v \Vert -1 \over 2 \Vert v \Vert}v,x)
&{\rm if}~\Vert v \Vert \geqslant 1\\
\infty&{\rm if}~v=0~{\rm or}~*~.
\end{cases}
$$

The relative difference of maps $p,q:V^{\infty} \wedge X \to Y$ which
agree on $0^{\infty} \wedge X$ can be interpreted as the difference
of two null-homotopies:

\begin{proposition}~ \label{difcon7}
Let $p,q:V^{\infty} \wedge X \to Y$ be pointed maps such that
$$p(0,x)~=~q(0,x) \in Y~(x \in X)~.$$
{\rm (i)} The pointed map defined by
$$f~:~ S(V)^+ \wedge X \to Y~;~(u,x) \mapsto p(0,x)=q(0,x)$$
has two null-homotopies
$$\begin{array}{l}
\delta p~:~C S(V)^+ \wedge X \to Y~;~(t,v,x) \mapsto p([t,v],x)~,\\[1ex]
\delta q~:~C S(V)^+ \wedge X \to Y~;~(t,v,x) \mapsto q([t,v],x)
\end{array}$$
and the relative difference of $p,q$ is given by
$$\begin{array}{l}
\delta(p,q)~=~-\delta q \cup \delta p~:\\[1ex]
\Sigma S(V)^+ \wedge X~=~ C S(V)^+\wedge X \cup_{S(V)^+\wedge X} CS(V)^+ \wedge X \to Y~.
\end{array}$$
{\rm (ii)} Under the one-one correspondence of Proposition \ref{pushout2} (ii)
between pointed maps $f:\Sigma S(V)^+ \wedge X \to Y$ and pairs $(g,h)$ defined by a
pointed map $g:S(U)^+ \wedge V^{\infty} \wedge X \to Y$ and a null-homotopy
$$h~:~g(j_{U,V} \wedge 1_X)~\simeq~\{*\}~:~S(U \oplus V)^+ \wedge X \to Y$$
the relative difference $f=\delta(p,q)$ corresponds to
$$\begin{array}{l}
g~=~\delta(p,q)(k_{U,V} \wedge 1_X)~=~(p-q)(s_U \wedge 1_{V^{\infty} \wedge X})~:~
S(U)^+ \wedge V^{\infty} \wedge X \to Y~,\\[1ex]
h~=~\delta(p,q)(\delta j_{U,V} \wedge 1_X)~:~
C(S(U \oplus V))^+ \wedge X \to Y~.
\end{array}$$
\end{proposition}
\begin{proof} (i) Immediate from Proposition \ref{difcon6} and Definition
\ref{reldif1}.\\
(ii) By construction.\\
\hfill\qed\end{proof}

A linear map $f:V \to W$ of inner product spaces is proper
if and only if $f$ is injective, in which case there is
defined an injective pointed map $f^{\infty}:V^{\infty} \to W^{\infty}$
of the one-point compactifications.

\begin{example} {\rm
(i) Let $f,g:V \to W$ be injective linear maps
of inner product spaces $V,W$, and let $X=S^0$, $Y=W^{\infty}$
in \ref{difcon7}.
The difference and relative difference maps of the pointed maps
$$f^{\infty}~,~g^{\infty}~:~V^{\infty}\wedge X~=~V^{\infty}
\to Y~=~W^{\infty}$$
are
$$\begin{array}{l}
f^{\infty}-g^{\infty}~=~
(g^{\infty} \vee f^{\infty})\overline{\nabla}_V~:~V^{\infty} \to W^{\infty}~,\\[1ex]
\delta(f^{\infty},g^{\infty})~=~
(g^{\infty} \vee f^{\infty})\gamma_V~:~\Sigma S(V)^+ \to W^{\infty}~.
\end{array}$$
(ii) For $f=-1$, $g=1:V \to W=V$ (i) gives the difference and relative
difference maps of the maps
$$p~=~1^{\infty}~=~1~,~q~=~(-1)^{\infty}~=~-1~:~V^{\infty}\wedge X~
=~V^{\infty} \to Y~=~V^{\infty}$$
to be
$$\begin{array}{l}
p-q~=~((-1)^{\infty} \vee 1^{\infty})\overline{\nabla}_V~\simeq~2~:~
V^{\infty} \to V^{\infty}~,\\[1ex]
\delta(p,q)~=~
((-1)^{\infty} \vee 1^{\infty})\gamma_V~:~\Sigma S(V)^+ \to V^{\infty}~.
\end{array}$$
\hfill\qed}
\end{example}

\begin{example} {\rm Let $X=S^0$.
The difference and relative difference maps of the maps
$$p,q~:~V^{\infty} \wedge X~=~V^{\infty} \to Y~=~V^{\infty}\vee_0V^{\infty}$$
defined by
$$p(v)~=~v_2~,~q(v)~=~v_1~(v \in V)~.$$
are
$$\begin{array}{l}
p-q~=~\beta_V~:~V^{\infty} \to V^{\infty}\vee_0V^{\infty}~,\\[1ex]
\delta(p,q)~=~\gamma_V~:~\Sigma S(V)^+ \to V^{\infty}\vee_0V^{\infty}~.
\end{array}$$
\hfill\qed}
\end{example}

\begin{proposition}~ \label{difcon2}
{\rm (i)} For any pointed  spaces $X,Y$  and inner product space $V$
there is defined an exact sequence of pointed sets
$$\begin{array}{l}
\xymatrix{\dots \ar[r] & [\Sigma X,Y] \ar[r]^-{\di{\Sigma s^*_V}} &
[\Sigma S(V)^+ \wedge X,Y] \ar[r]^-{\di{\alpha^*_V}} &
[V^{\infty} \wedge X,Y]}\\
\hskip150pt \xymatrix{\ar[r]^-{\di{0^*_V}}& [X,Y] \ar[r]^-{\di{s^*_V}} &
[S(V)^+ \wedge X,Y]}
\end{array}$$
with
$$0^*_V~:~[V^{\infty} \wedge X,Y]\to
[0^+ \wedge X,Y]~=~[X,Y] ~;~F \mapsto F\vert_{0^+ \wedge X}~.$$
For any maps $p,q:V^{\infty} \wedge X \to Y$ which agree on
$0^+ \wedge X$ the difference $p-q \in [V^{\infty} \wedge X,Y]$ has
image $* \in [X,Y]$, and $p-q$ is the image of $\delta(p,q)
\in [\Sigma S(V)^+ \wedge X,Y]$, with a homotopy commutative diagram
$$\xymatrix@C+30pt@R+10pt{
V^{\infty} \wedge X \ar[r]^-{\di{\alpha_V \wedge 1_X}}
\ar[dr]_{\di{p-q}} &
\Sigma S(V)^+ \wedge X \ar[d]^{\di{\delta(p,q)}}\\
    & Y}$$
that is
$$\begin{array}{l}
p-q~=~\alpha_V^*(\delta(p,q))
\in {\rm ker}(0^*_V:[V^{\infty} \wedge X,Y] \to [X,Y])\\[1ex]
\hskip110pt
=~{\rm im}(\alpha_V^*:[\Sigma S(V)^+ \wedge X,Y]\to [V^{\infty} \wedge X,Y])~.
\end{array}$$
{\rm (ii)} If the maps $p,q:V^{\infty}\wedge X \to Y$ in {\rm (i)}
are related by a homotopy
$$r~:~p ~\simeq~ q~:~V^{\infty} \wedge X \to Y$$
then
$$p-q~\simeq~*~:~V^{\infty} \wedge X \to Y$$
and  the relative difference $\delta(p,q) \in [\Sigma S(V)^+\wedge X,Y]$
is the image of the map
$$r_0~:~\Sigma X \to Y~;~(t,x) \mapsto r(t,0,x)~,$$
with a homotopy commutative diagram
$$\xymatrix@C+30pt@R+10pt{
\Sigma S(V)^+ \wedge X \ar[r]^-{\di{\Sigma s_V \wedge 1_X}}
\ar[dr]_{\di{\delta(p,q)}} &
S^1 \wedge X = \Sigma X\ar[d]^{\di{r_0}}\\
    & Y}$$
i.e. if $p=q \in [V^{\infty} \wedge X,Y]$ then
$$\begin{array}{l}
\alpha_V^*(\delta(p,q))~=~p-q~=~* \in [V^{\infty} \wedge X,Y]~,\\[1ex]
\delta(p,q)~=~\Sigma s^*_V(r_0) \in{\rm ker}(\alpha^*_V:
[\Sigma S(V)^+ \wedge X,Y] \to [V^{\infty} \wedge X,Y])\\[1ex]
\hskip100pt
=~{\rm im}(\Sigma s^*_V:[\Sigma X,Y] \to [\Sigma S(V)^+ \wedge X,Y])~.
\end{array}
$$
In particular, if
$$r(t,0,x)~=~p(0,x)~=~q(0,x) \in Y~(t \in I,x \in X)$$
then $\delta(p,q)\simeq *:\Sigma S(V)^+ \wedge X\to  Y$,
i.e. if $p,q$ are related by a homotopy which is constant on
$\{0\} \times X$ then
$$\delta(p,q)~=~* \in [\Sigma S(V)^+ \wedge X,Y]~.$$
{\rm (iii)} If $p,p',q,q':V^{\infty} \wedge X \to Y$ are
maps which agree on $0^+ \wedge X$ and
$$f~:~p~\simeq~p'~,~g~:~q~\simeq~q'~:~V^{\infty} \wedge X \to Y$$
are homotopies which are constant on $0^+ \wedge X$
there is induced a homotopy of the relative differences
$$\delta(f,g)~:~\delta(p,q)~\simeq~\delta(p',q')~:~
    \Sigma S(V)^+ \wedge X \to Y~.$$
{\rm (iv)} If $p,q,r:V^{\infty}
\wedge X \to Y$ are maps which agree on $0^+ \wedge
X$ the relative differences are related by a homotopy
$$h~:~\delta(p,q)+\delta(q,r)~\simeq~\delta(p,r)~:~
    \Sigma S(V)^+ \wedge X \to Y~.$$
{\rm (v)} If $p,q:V^{\infty}\wedge X \to Y$ are maps such that
$$p(t,x)~=~\begin{cases}
q(t,x)&{\it if}~ (t,x) \in U\\
*&{\it if}~ (t,x) \in \overline{(V^{\infty} \wedge X) \backslash U}
\end{cases}$$
for some neighbourhood $U \subseteq V^{\infty} \wedge X$ of
$0^+\wedge X\subset V^{\infty} \wedge X$
there exists a homotopy
$$\gamma(p,q)~:~\delta(p,q)~\simeq~q'~:~\Sigma S(V)^+ \wedge X \to Y$$
with $q'$ defined by
$$\begin{array}{l}
q'~:~\Sigma  S(V)^+\wedge X=(V^{\infty}/0^+)\wedge X\to Y~;\\[1ex]
\hphantom{q'~:~}(t,x) \mapsto \begin{cases}
*&{\it if}~ (t,x) \in U\\
q(t,x)&{\it if}~ (t,x) \in \overline{V^{\infty}\wedge X \backslash U}~.
\end{cases}
    \end{array}$$
{\rm (vi)} The homotopy $\gamma(p,q)$ in {\rm (v)} can be chosen
to be natural, meaning that given commutative squares of maps
$$\xymatrix{
V^{\infty} \wedge X_1 \ar[r]^-{\di{p_1}} \ar[d]_-{\di{1 \wedge h}} &  Y_1\ar[d]^-{\di{k}} \\
V^{\infty} \wedge X_2 \ar[r]^-{\di{p_2}} &  Y_2}~\lower19pt\hbox{,}~
\xymatrix{
V^{\infty} \wedge X_1 \ar[r]^-{\di{q_1}} \ar[d]_-{\di{1 \wedge h}} &  Y_1\ar[d]^-{\di{k}} \\
V^{\infty} \wedge X_2 \ar[r]^-{\di{q_2}} &  Y_2}$$
and neighbourhoods
$U_i \subseteq V^{\infty} \wedge X_i$ of
$0^+\wedge X_i\subset V^{\infty} \wedge X_i$
with $U_2=(1\wedge h)(U_1)$ such that
$$p_i(t,x)~=~\begin{cases}
q_i(t,x)&{\it if}~ (t,x) \in U_i\\
*&{\it if}~ (t,x) \in \overline{V^{\infty} \wedge X_i \backslash U_i}
\end{cases}$$
there is defined a commutative square
$$\xymatrix@C+15pt@R+15pt{
I \times (\Sigma S(V)^+ \wedge X_1) \ar[r]^-{\di{\gamma(p_1,q_1)}}
\ar[d]_-{\di{1\wedge h}} &
Y_1\ar[d]^-{\di{k}} \\
I \times (\Sigma S(V)^+ \wedge X_2) \ar[r]^-{\di{\gamma(p_2,q_2)}}
&  Y_2}$$
\end{proposition}
\begin{proof}
(i) Immediate from the homotopy cofibration sequence of pointed spaces
$$S(V)^+\xymatrix{\ar[r]^-{\di{s_V}}&}
S^0~=~0^+ \to V^{\infty} \xymatrix{\ar[r]^-{\di{\alpha_V}}&} \Sigma S(V)^+
\xymatrix{\ar[r]^-{\di{\Sigma s_V}}&}  S^1 \to \dots$$
given by Proposition \ref{cofib2} (i).\\
(ii) Immediate from the homotopy $\beta_V \simeq i\overline{\nabla}_V$
given by Proposition \ref{cofib2} (iii).\\
(iii) The homotopies
$$f,g~:~ I \times (V^{\infty} \wedge X) \to Y$$
which agree on $I \times 0^+ \wedge X$. The relative difference
$$\delta(f,g)~:~I \times (\Sigma S(V)^+ \wedge X)  \to Y~.$$
defines a homotopy $\delta(f,g):\delta(p,q)\simeq \delta(p',q')$.\\
(iv) The sum of $\delta(p,q)$ and $\delta(q,r)$ is defined by
$$\begin{array}{l}
\delta(p,q)+\delta(q,r)~:~\Sigma S(V)^+ \wedge X \to Y~;\\[1ex]
(t,u,x) \mapsto \begin{cases}\delta(q,r)(2t,u,x)&\hbox{if $0 \leqslant t
\leqslant  1/2$}\\[1ex]
\delta(p,q)(2t-1,u,x)&\hbox{if $1/2 \leqslant  t
\leqslant  1$~.}
\end{cases}
\end{array}$$
Define a homotopy
$h:\delta(p,q)+\delta(q,r)\simeq\delta(p,r)$ by
$$h(s,t,u,x)~=~
\begin{cases}r(\big[\di{1+s-4t \over 1+s},u\big],x)&\hbox{if $ 0 \leqslant
t \leqslant  (1+s)/4$}\\[1ex]
q(\big[
4t-1-s,u\big],x)&\hbox{if $(1+s)/4 \leqslant  t \leqslant
1/2$}\\[1ex]
q(\big[ 3-s-4t,u\big],x)&\hbox{if $1/2
\leqslant  t \leqslant  (3-s)/4$}\\[1ex]
p(\big[\di{s+4t-3 \over 1+s},u\big],x)&
\hbox{if $(3-s)/4 \leqslant t\leqslant  1$~.}
\end{cases}$$
(v)+(vi) These are special cases of Proposition \ref{near}.\\
\hfill\qed\end{proof}

\begin{example} \label{inj} {\rm
(i) For any inner product spaces $V,W$ and injective linear
maps $f,g:V \to W$ the maps
$f^{\infty},g^{\infty}:V^{\infty} \to W^{\infty}$ are such that
$$f^{\infty}(0_V)~=~g^{\infty}(0_V)~=~0_W \in W^{\infty}~,$$
so that there is defined a relative difference map
$$\delta(f^{\infty},g^{\infty})~:~\Sigma S(V)^+ \to W^{\infty}~;~
[t,u] \mapsto
\begin{cases}
g([1-2t,u])&{\rm if}~0 \leqslant t \leqslant 1/2\\
f([2t-1,u])&{\rm if}~1/2 \leqslant t \leqslant 1~.
\end{cases}$$
By Proposition \ref{difcon2} (i) there exists a homotopy
$$\delta(f^{\infty},g^{\infty}) \alpha_V~\simeq~f^{\infty}-g^{\infty}~:~
V^{\infty} \to W^{\infty}~.$$
(ii) Let $GL(V)$ be the space of linear automorphisms
$a:V \to V$, with base point $1:V \to V$.
For any map $c:X \to GL(V)$ the maps defined by
$$\begin{array}{l}
p~=~{\rm adj}(1)~:~V^{\infty} \wedge X \to V^{\infty}~;~(v,x) \mapsto v~,\\[1ex]
q~=~{\rm adj}(c)~:~V^{\infty} \wedge X \to V^{\infty}~;~(v,x) \mapsto c(x)(v)
\end{array}$$
are such that $p(0,x)=q(0,x)=0 \in V^{\infty}$ ($x \in X$).
Use the relative difference construction of (i) to define a function
$$\delta_V~:~[X,GL(V)] \to [\Sigma S(V)^+ \wedge X,V^{\infty}]~;~c
\mapsto \delta_V(c)~=~\delta(p,q)$$
with
$$\delta_V(c)([t,u],x)~=~\begin{cases}
c(x)^{\infty}[1-2t,u]&{\rm if}~0 \leqslant t \leqslant 1/2\\
[2t-1,u]&{\rm if}~1/2 \leqslant t \leqslant 1~.
\end{cases}$$
The composite
$$\xymatrix{[X,GL(V)] \ar[r]^-{\di{\delta_V}}& [\Sigma S(V)^+\wedge X,V^\infty]
\ar[r]^-{\di{\alpha^*_V}} & [V^{\infty} \wedge X,V^{\infty}]}$$
is given by
$$\alpha^*_V \delta_V~:~
[X,GL(V)] \to [V^{\infty} \wedge X,V^{\infty}]~;~c \mapsto p-q~.$$
}
\hfill\qed
\end{example}

\begin{example} {\rm
Here are two special cases of Proposition \ref{difcon2} (v).\\
(i) For any map $p:V^{\infty} \wedge X \to Y$ take
$q=p$, $U=V^{\infty} \wedge X$ to obtain a homotopy
$$\gamma(p,p)~:~\delta(p,p)~\simeq~*~:~\Sigma S(V)^+\wedge X \to Y~.$$
(ii) For any map $q:V^{\infty} \wedge X \to Y$ such that
$$q(0,x)~=~* \in Y~(x \in X)$$
take $q=*$, $U=0^+ \wedge X$ to obtain a homotopy
$$\gamma(p,*)~:~\delta(p,*)~\simeq~[p]~:~\Sigma S(V)^+ \wedge X \to Y~.$$
with $[p]$ induced from $p$ using the homeomorphism $\Sigma S(V)^+ \cong
V^{\infty}/0^+$.\hfill\qed}
\end{example}

\begin{proposition}~ {\it
The relative difference (\ref{difcon1} (ii)) of maps $p,q:V^{\infty} \wedge X
\to Y$ such that $p \vert = q\vert :X \to Y$
is the relative difference (\ref{reldif1}) of the null-homotopies $\delta p,\delta q$ of
$$(p\vert )(s_V \wedge 1_X)~=~(q\vert )(s_V \wedge 1_X)~:~S(V)^+ \wedge X \to Y$$
given by \ref{difcon6}, that is
$$\begin{array}{l}
\delta(p,q)~:~\xymatrix@C+40pt{\Sigma S(V)^+ \wedge X
\ar[r]^-{\di{(-1\cup 1 )\wedge 1_X}}&}\\[1ex]
\hskip40pt
\xymatrix@C+20pt{CS(V)^+ \wedge X\cup_{S(V)^+ \wedge X}
CS(V)^+ \wedge X
\ar[r]^-{\di{\delta q\cup \delta p}}&Y}
\end{array}$$
with}
$$\delta p(t,u,x)~=~p([t,u],x)~,~\delta q(t,u,x)~=~q([t,u],x)~.$$
\end{proposition}
\begin{proof} For all $t \in I$, $u \in S(V)$, $x \in X$
$$\delta(p,q)(t,u,x)~=~\begin{cases}
q(1-2t,u,x)&\hbox{if $0 \leqslant  t \leqslant  1/2$}\\[1ex]
p(2t-1,u,x)&\hbox{if $1/2 \leqslant  t \leqslant  1$~.}
\end{cases}$$
\hfill\qed
\end{proof}

\chapter{Stable homotopy theory}\label{stablehomotopy}

This Chapter develops stable homotopy theory and bordism theory using inner product spaces.

\section{Stable maps}

A {\it stable map} $F:X \to Y$ between pointed spaces $X,Y$ is a
map of the type
$$F~:~V^{\infty} \wedge X \to V^{\infty} \wedge Y$$
for some (finite-dimensional) inner product space $V$.

The {\it stable homotopy group} is the abelian group
$$\{X;Y\}~=~\mathop{\varinjlim}_V\, [V^{\infty}\wedge X,V^{\infty}\wedge Y]$$
with the direct limit is taken over all (finite-dimensional) inner product
spaces $V$, and addition and subtraction maps induced by the sum and difference maps
$\nabla_V,\overline{\nabla}_V:V^{\infty} \to V^{\infty}\vee V^{\infty}$.
By definition, an element $F \in \{X;Y\}$ is an equivalence class of stable maps $F:X \to Y$.

\begin{example} {\rm For an open embedding
$e:V \times M \emb V \times N$ the compactification Umkehr
$F:V^{\infty} \wedge N^\infty \to V^{\infty} \wedge M^+$
defines a stable map $F:N^\infty \to M^+$.
\hfill\qed}
\end{example}

Given a map $F:A \to B$ let
$G:B \to \Cc(F)$ be the inclusion in the mapping cone, and let
$H:\Cc(F) \to \Sigma A$ be the projection, so that
$$\xymatrix@C+10pt{A \ar[r]^-{\di{F}} &
B \ar[r]^-{\di{G}} &
\Cc(F)\ar[r]^-{\di{H}}& \Sigma A \ar[r]^-{\di{\Sigma F}} & \Sigma B \ar[r] &\dots}$$
is a homotopy cofibration sequence, as before. By analogy with the
Barratt-Puppe cohomotopy exact sequence (\ref{BP1}) there is a stable homotopy
exact sequence~:

\begin{proposition}~\label{BP2}
For any pointed space $X$ there is induced an exact sequence of stable homotopy groups
$$\begin{array}{l}
\xymatrix{\dots \ar[r]& \{\Sigma X;\Cc(F)\} \ar[r]^-{\di{H}} &
\{X;A\}\ar[r]^-{\di{F}} & \{X;B\} }\\
\hskip100pt
\xymatrix{\ar[r]^-{\di{G}} &\{X;\Cc(F)\} \ar[r]^-{\di{H}} & \{X;\Sigma A\} \ar[r] & \dots~.}
\end{array}$$
\hfill\qed
\end{proposition}

The homotopy cofibration sequence
$$\xymatrix{S(V)^+ \ar[r]^-{\di{s_V}} & S^0 \ar[r]^-{\di{0_V}} &
V^{\infty} \ar[r]^-{\di{\alpha_V}} & \Sigma S(V)^+ \ar[r]&\dots}$$
determines the following braid~:

\begin{proposition}~ \label{stablesequence0}
For any inner product spaces $U,V$ and pointed spaces $X,Y$
there is defined a commutative braid of exact sequences of stable
homotopy groups
$$\xymatrix@C-30pt{
A_1 \ar[dr] \ar@/^2pc/[rr]^-{\di{\alpha_V}} && \{X;S(V)^+ \wedge Y\}
\ar[dr]^-{\di{s_{V}}}\ar@/^2pc/[rr] &&
\{S(U)^+ \wedge X;Y\}\ar[dr]&\\
&~~~A_2~~~ \ar[ur] \ar[dr] &&
\{X;Y\} \ar[ur]^-{\di{s^*_{U}}} \ar[dr]^-{\di{0_{V}}}&&A_3\\
A_4\ar[ur] \ar@/_2pc/[rr]^-{\di{\alpha^*_U}}
&&\{U^{\infty} \wedge X;Y\}
\ar@/_2pc/[rr]\ar[ur]^-{\di{0^*_{U}}} &&
\{X;V^{\infty}\wedge Y\}\ar[ur]&}$$

\bigskip

\noindent with
$$\begin{array}{l}
A_1~=~\{\Sigma X;V^{\infty} \wedge Y\}~,~
A_2~=~\{\Sigma S(U\oplus V)^+\wedge X;V^{\infty} \wedge Y\}~,\\[1ex]
A_3~=~\{S(U\oplus V)^+ \wedge X;V^{\infty} \wedge Y\}~,~
A_4~=~\{\Sigma S(U)^+ \wedge X;Y\}~.
\end{array}$$
\end{proposition}
\begin{proof} These are the exact sequences
(\ref{BP2}) determined by the homotopy commutative braid of homotopy cofibrations

$$\xymatrix@C-30pt{
&U^{\infty} \wedge S(V)^+\ar[dr]^-{\di{1 \wedge s_V}}
\ar@/^2pc/[rr] && \Sigma S(U)^+ \ar[dr]&\\
S(U\oplus V)^+\ar[dr]^-{\di{s_{U\oplus V}}} \ar[ur] &&
U^{\infty}\wedge S^0 =U^{\infty}
\ar[dr]^-{\di{1\wedge 0_V}}\ar[ur]^-{\di{\alpha_V}}&&\Sigma S(U \oplus V)^+&\\
&S^0\ar@/_2pc/[rr]^-{\di{0_{U \oplus V}}}
\ar[ur]^-{\di{0_U}} &&U^{\infty} \wedge V^{\infty}=
(U \oplus V)^{\infty}\ar[ur]^-{\di{\alpha_{U \oplus V}}} &}$$

\noindent given by Proposition \ref{cofib2} (iv), involving two homotopy pushout squares.\\
\hfill\qed\end{proof}

We shall need the following version of exactness at $\{X;A\}$
in the sequence (\ref{BP2})
$$\begin{array}{l}
\xymatrix{\dots \ar[r]& \{\Sigma X;\Cc(F)\} \ar[r]^-{\di{H}} &
\{X;A\}\ar[r]^-{\di{F}} & \{X;B\} }\\
\hskip100pt
\xymatrix{\ar[r]^-{\di{G}} &\{X;\Cc(F)\} \ar[r]^-{\di{H}} & \{X;\Sigma A\} \ar[r] & \dots~.}
\end{array}$$

\begin{definition}~  \label{stablereldif1}
{\rm Let $X$ be a pointed space with a difference map
$\overline{\nabla}:X \to X \vee X$ (\ref{sumdifference} (iii)). Given maps $a_1,a_2:X \to A$, $F:A \to B$ and
a homotopy $b:Fa_1\simeq Fa_2:X \to B$ let
$$\overline{b}~:~Fa_2-Fa_1 ~\simeq~ *~:~X \to B$$
be the null-homotopy of $Fa_2-Fa_1=F(a_2-a_1):X \to B$
defined by the concatenation of the  homotopy
$$b-1~:~Fa_2-Fa_1~\simeq~ Fa_1-Fa_1~:~X \to B$$
and the standard null-homotopy $Fa_1-Fa_1 \simeq *:X \to B$.
The {\it stable relative difference}
$$\delta(a_1,a_2,b)~=~c \in \{\Sigma X;\Cc(F)\}$$
is the stable homotopy class of the map $c:\Sigma X \to \Cc(F)$
in the map of homotopy cofibrations
$$\xymatrix{X \ar[r] \ar[d]_-{\di{a_1-a_2}} & CX \ar[r]
\ar[d]^-{\di{\overline{b}}}
& \Sigma X \ar[d] ^-{\di{c}} \ar@{=}[r]& \Sigma X \ar[d]^-{\di{\Sigma (a_1-a_2)}}\\
A \ar[r]^-{\di{F}} & B \ar[r]^-{\di{G}} & \Cc(F) \ar[r]^-{\di{H}} & \Sigma A}$$
\hfill\qed}
\end{definition}

\begin{proposition}~\label{stablereldif2}
The stable relative difference $\delta(a_1,a_2,b) \in
\{\Sigma X;\Cc(F)\}$ of \ref{stablereldif1} is such that
$$\begin{array}{l}
a_1-a_2~=~H(\delta(a_1,a_2,b)) \in{\rm ker}(F:\{X;A\} \to \{X;B\})\\[1ex]
\hskip50pt =~{\rm im}(H:\{\Sigma X;\Cc(F)\}
\to \{X;A\}) \subseteq \{X;A\}~=~\{\Sigma X;\Sigma A\}~.
\end{array}$$
\end{proposition}
\begin{proof} By construction.\\
\hfill\qed\end{proof}

The homotopy cofibration sequence
$$\xymatrix{S(V)^+ \wedge Y \ar[r]^-{\di{s_V}}&
Y \ar[r]^-{\di{0_V}} & V^{\infty} \wedge Y \ar[r] &
\Cc(0_V) \simeq \Sigma S(V)^+ \wedge Y \ar[r]& \dots}$$
induces a long exact sequence of stable homotopy groups
$$\begin{array}{l}
\xymatrix{\dots \ar[r] & \{V^{\infty} \wedge X;S(V)^+ \wedge Y\}
\ar[r]^-{\di{s_V}} & \{V^{\infty} \wedge X;Y\}}\\[1ex]
\xymatrix{
\ar[r]^-{\di{0_V}}&\{V^{\infty} \wedge X;V^{\infty} \wedge Y\}=\{X;Y\} \ar[r] &
\{V^{\infty} \wedge X;\Sigma S(V)^+ \wedge Y\}\ar[r] &\dots}
\end{array}$$

\begin{definition}~\label{stablereldif3}
{\rm Given a map $p:V^{\infty} \wedge X \to Y$ define the map
$$f~=~0_Vp~:~V^{\infty} \wedge X \to V^{\infty} \wedge Y~;~(v,x) \mapsto
(0,p(v,x))~.$$
Given also a map $q:V^{\infty} \wedge X \to Y$ which agrees with $p$
on $0^+ \wedge X \subset V^{\infty} \wedge X$ let
$g=0_Vq:V^{\infty} \wedge X \to V^{\infty} \wedge Y$, and
define a homotopy
$$h~:~f \simeq g~:~V^{\infty} \wedge X \to V^{\infty} \wedge Y$$
by
$$\begin{array}{l}
h~:~I \times V^{\infty} \wedge X \to V^{\infty} \wedge Y~;\\[1ex]
\hskip50pt
(t,v,x) \mapsto \begin{cases}
(2tv,q((1-2t)v,x))&{\it if}~0 \leqslant t \leqslant 1/2\\
((2-2t)v,p((2t-1)v,x))&{\it if}~1/2 \leqslant t \leqslant 1~.
\end{cases}
\end{array}$$
The {\it stable relative difference} of $p,q$ is the stable
relative difference (\ref{stablereldif1})
$$\delta'(p,q)~=~\delta(f,g,h) \in \{\Sigma(V^{\infty} \wedge X);\Cc(0_V)\}~=~
\{V^{\infty} \wedge X;S(V)^+ \wedge Y\}~.$$}
\hfill\qed
\end{definition}

\begin{proposition}~ \label{stablereldif4}
The stable relative difference of \ref{stablereldif3} is such that
$$\begin{array}{l}
p-q~=~s_V\delta'(p,q)
\in {\rm im}(s_V:\{V^{\infty} \wedge X;S(V)^+ \wedge Y\} \to \{V^{\infty} \wedge X;Y\})\\[1ex]
\hskip120pt =~{\rm ker}(0_V:\{V^{\infty} \wedge X;Y\} \to \{X;Y\})~.
\end{array}$$
\end{proposition}
\begin{proof} Immediate from Proposition \ref{stablereldif2}.
\hfill\qed\end{proof}

The Umkehr maps of open embeddings $V \times X \emb V \times Y$
are stable maps. Specifically, given a map $f:X \to Y$, an inner product space $V$
and a map $g:V \times X \to V$ such that
$$e~=~(g,f)~:~V \times X \to V \times Y~;~(v,x) \mapsto (g(v,x),f(x))$$
is an open embedding then Definition
\ref{umkehrcompact} gives a compactification Umkehr stable map
$$F~:~(V\times Y)^{\infty}~=~V^{\infty} \wedge Y^{\infty} \to
(V\times X)^{\infty}~=~V^{\infty} \wedge X^{\infty}~.$$
For compact $X,Y$ this is an adjunction Umkehr stable map
$F:V^{\infty} \wedge Y^+ \to V^{\infty} \wedge X^+$.

\begin{example} \label{fincover}
{\rm For any finite cover $f:\widetilde{K} \to K$ of a $CW$ complex $K$
there is defined a transfer map
$$f^!~:~C(K) \to C(\widetilde{K})~;~x \mapsto \sum\limits_{y \in f^{-1}(x)}y~.$$
For finite $K$ there exist a finite-dimensional inner product space $V$ and a
map $g:V \times \widetilde{K} \to V$ such that
$$e~=~(g,f)~:~V \times \widetilde{K} \to V \times K~;~(v,x) \mapsto (g(v,x),f(x))$$
is an open embedding, with compactification (= adjunction) Umkehr map
$$F~:~(V\times K)^{\infty}~=~V^{\infty} \wedge K^+ \to
(V\times \widetilde{K})^{\infty}~=~V^{\infty} \wedge \widetilde{K}^+$$
inducing the transfer map $F=f^!:C(K) \to C(\widetilde{K})$
(cf. Adams \cite[\S 4.2]{adams2}, Proposition \ref{Z2cover}
and Example \ref{ad} below).\\
\hfill\qed}
\end{example}

We shall need a version of the Umkehr construction for pointed maps.

\begin{definition}~ \label{pointedumkehr}
{\rm Let $f:X \to Y$ be a pointed map, with $f(x_0)=y_0 \in Y$ for
base points $x_0 \in X$, $y_0 \in Y$. Let $V$ be an inner product space,
and suppose given a map $g:V \times X \to V$ such that
$$e~=~(g,f)~:~V \times X \to V \times Y~;~(v,x) \mapsto (g(v,x),f(x))$$
restricts to an open embedding
$$e\vert~:~V\times (X\backslash \{x_0\}) \to V \times Y$$
and also
$$g(v,x_0)~=~0\in V~(v \in V)~.$$
The {\it Umkehr stable map} of $e$ is\index{Umkehr!stable}
$$\begin{array}{l}
F~:~V^{\infty} \wedge Y \to V^{\infty} \wedge X~;\\[1ex]
\hskip50pt
(w,y) \mapsto \begin{cases} (v,x)&{\rm if}~(w,y)=(g(v,x),f(x)) \in V \times Y \\
\infty&{\rm otherwise}~.
\end{cases}
\end{array}$$
\hfill\qed}
\end{definition}

\begin{example} \label{pointedfincover}
{\rm
(i) Let $f:X \to Y$ be a pointed map of finite $CW$ complexes such that
the restriction $f\vert:X \backslash \{x_0\} \to Y \backslash \{y_0\}$
is a finite cover, so that there is defined a transfer chain map
$$f^!~:~\dot C(Y) \to \dot C(X)~;~ y \mapsto \sum\limits_{x \in f^{-1}(y)}x$$
of the reduced singular chain complexes. As in Adams \cite[pp.\ 511--512]{adams3} there
exist a finite-dimensional inner product space $V$ and a
map $g:V \times X \to V$ such that
$$e~=~(g,f)~:~V \times X \to V \times Y~;~(v,x) \mapsto (g(v,x),f(x))$$
satisfies the hypothesis of Definition \ref{pointedumkehr}, giving an
Umkehr stable map
$$F~:~V^{\infty} \wedge Y \to V^{\infty} \wedge X$$
inducing $f^!$ on the chain level.\\
(ii) If $f:\widetilde{K} \to K$ is a finite cover of a finite $CW$ complex
then for any finite-dimensional inner product space $W$
$$1\wedge f^+~:~X~=~W^{\infty} \wedge \widetilde{K}^+ \to Y~=~W^{\infty} \wedge K^+$$
is a pointed map as in (i), and the Umkehr stable map
$$F~:~V^{\infty} \wedge Y~=~(V\oplus W)^{\infty} \wedge K^+ \to
V^{\infty} \wedge X~=~(V\oplus W)^{\infty} \wedge \widetilde{K}^+$$
is just the compactification Umkehr map of Example \ref{fincover} above,
inducing the transfer chain map
$$F~=~f^!~:~\dot C(Y)~=~C(K)_{*-{\rm dim}(W)} \to \dot C(X)~=~C(\widetilde{K})_{*-{\rm dim}(W)}~.$$
\hfill\qed}
\end{example}

\begin{definition}~{\rm
(i) The {\it stable homotopy} and {\it cohomotopy} groups
of a space $X$ are
$$\omega_n(X)~=~\{S^n;X^+\}~,~\omega^n(X)~=~\{X^+;S^n\}~~(n \in \ZZ)~.$$
(ii) The {\it reduced stable homotopy} and {\it cohomotopy} groups of  a pointed  space $X$ are
$$\widetilde{\omega}_n(X)~=~\{S^n;X\}~,~\widetilde{\omega}^n(X)~=~\{X;S^n\}~.$$
In particular, for an unpointed space $X$
$$\omega_n(X)~=~\widetilde{\omega}_n(X^+)~,~\omega^n(X)~=~
\widetilde{\omega}^n(X^+)~.$$
\hfill\qed}
\end{definition}

\begin{example} {\rm The stable homotopy and cohomotopy groups of $\{*\}$
are the stable homotopy groups of spheres
$$\omega_n(\{*\})~=~\omega^{-n}(\{*\})~=~\{S^n;S^0\}~(n \geqslant 0)$$
which we write as $\omega_n$.\\
\hfill\qed}
\end{example}

\begin{definition}~{\rm
(i) A {\it spectrum} $\underline{X}=\{X(V)\}$ is a
sequence of pointed spaces $X(V)$ indexed by finite-dimensional
inner product spaces $V$, with structure maps
$$(V^{\perp})^{\infty} \wedge X(V) \to X(W)$$
defined whenever $V \subseteq W$, where $V^{\perp} \subseteq W$ is
the orthogonal complement of $V$ in $W$. For $n \in \ZZ$ let
$$\pi_n(\underline{X})~=~
\begin{cases}
\varinjlim\limits_V [\Sigma^nV^{\infty},X(V)]&{\rm if}~n \geqslant 0\\[1ex]
\varinjlim\limits_V [V^{\infty},\Sigma^{-n}X(V)]&{\rm if}~n \leqslant -1~.
\end{cases}$$
(ii)  A spectrum $\underline{X}$ is {\it connective} if
$$\pi_n(\underline{X})~=~0~{\rm for}~n \leqslant -1~.$$
In particular, this is the case if each $X(V)$ is $({\rm dim}(V)-1)$-connected,
i.e.
$$\pi_n(X(V))~=~0~{\rm for}~n \leqslant ({\rm dim}(V)-1)~.$$
(iii) The {\it $\underline{X}$-coefficient homology} groups of a space $Y$ are
$$\underline{X}_n(Y)~=~\pi_n(\underline{X} \wedge Y^+)~~(n \in\ZZ)~.$$
The {\it reduced $\underline{X}$-coefficient homology} groups of a
pointed space $Y$ are
$$\widetilde{\underline{X}}_n(Y)~=~\pi_n(\underline{X} \wedge Y)~~(n \in\ZZ)~.$$
If $\underline{X}$ is connective then $\underline{X}_n(Y)=0$ for $n \leqslant -1$,
and similarly in the reduced case.\\
\hfill\qed}
\end{definition}

\begin{example} {\rm
(i) The suspension spectrum of a pointed space $X$ is
the connective spectrum
$$\underline{X}~=~\{V^{\infty} \wedge X\,\vert\, {\rm dim}(V)<\infty\}~,$$
such that for any pointed space $Y$
$$\underline{X}_n(Y)~=~\begin{cases}
\widetilde{\omega}_n(X\wedge Y)&{\rm if}~n \geqslant 0\\[1ex]
0&{\rm if}~n \leqslant -1~.
\end{cases}$$
In particular, for $Y=S^0$
$$\pi_n(\underline{X})~=~\underline{X}_n(S^0)=~\begin{cases}
\widetilde{\omega}_n(X)&{\rm if}~n \geqslant 0\\[1ex]
0&{\rm if}~n \leqslant -1~.
\end{cases}$$
(ii) The sphere spectrum $\underline{S}=\{V^{\infty}\}$ (= the suspension
spectrum of $S^0$) is  such that
$$\pi_n(\underline{S})~=~\underline{S}_n(S^0)~=~\omega_n~,~
\underline{S}_n(X)~=~\widetilde{\omega}_n(X)~.$$
\hfill\qed}
\end{example}

\section{Vector bundles}\label{vectorbundles}

We shall only be considering vector bundles and
spherical fibrations over $CW$ complexes.

\begin{definition}~{\rm
Let $V$ be an inner product space.\\
(i) A {\it $V$-bundle} over a space $X$ is a vector bundle
$$\xymatrix{\xi~:~V \ar[r] & E(\xi) \ar[r]^-{\di{p_{\xi}}}&X}$$
with each fibre $p^{-1}(x)=\xi_x$ ($x \in X$) an inner product
space isomorphic to $V$, and the transition functions
linear isometries. \\
(ii) The {\it Thom space} of $\xi$ is the pointed space
\index{Thom space, $T(\xi)$}
$$T(\xi)~=~D(\xi)/S(\xi)~,$$
with $(D(\xi),S(\xi))$ the total space of the $(D(V),S(V))$-bundle of $\xi$
$$(D(V),S(V)) \xymatrix{\ar[r] &} (D(\xi),S(\xi))
\xymatrix@C+30pt{\ar[r]^-{\di{p_{\xi}\vert}} & X~.}$$
(iii) The {\it trivial $V$-bundle} $\epsilon_V$ over $X$ has
$$\begin{array}{l}
p~=~{\rm projection}~:~E(\epsilon_V)~=~V \times X \to X~,\\[1ex]
D(\epsilon_V)~=~D(V)\times X ~,~S(\epsilon_V)~=~S(V)\times X~,\\[1ex]
T(\epsilon_V)~=~(D(V)\times X)/(S(V) \times X)~=~V^\infty\wedge X^+~.
\end{array}$$
In particular, for $V=0$
$$D(\epsilon_V)~=~X~,~S(\epsilon_V)~=~\emptyset~,~T(\epsilon_V)~=~X^+~.$$
(iv) The {\it Thom spectrum} of a $V$-bundle $\xi$ is the suspension
spectrum of $T(\xi)$
$$\underline{T}(\xi)~=~\{T(\xi \oplus \epsilon_V)\,\vert\,{\rm dim}(V)<\infty\}$$
with $T(\xi \oplus \epsilon_V)=V^{\infty} \wedge T(\xi)$.\\
\hfill\qed}
\end{definition}

\begin{remark}{\rm The canonical map
$$E(\xi)^{\infty} \to T(\xi)~;~(x,v) \mapsto (x,\dfrac{v}{1+\Vert v\Vert})~
(x \in X,v \in V)$$ is a bijection which is a homeomorphism if and
only if $X$ is compact.\\}
\hfill\qed
\end{remark}

\begin{definition}~{\rm
(i) The {\it pullback} of a $V$-bundle $\xi$ over $X$ along a map $f:Y \to X$
is the $V$-bundle over $Y$
$$\xymatrix{f^*\xi~:~V \ar[r] & E(f^*\xi) \ar[r]^-{\di{p_{f^*\xi}}}&Y}$$
with $E(f^*\xi)$ fitting into a pullback square
$$\xymatrix{
E(f^*\xi) \ar[r] \ar[d]_-{\di{p_{f^*\xi}}} &
E(\xi) \ar[d]^-{\di{p_\xi}}\\
Y \ar[r]^-{\di{f}} & X~.}$$
(ii) The {\it product} of a $V$-bundle over $X$ and a $W$-bundle over
$Y$
$$\xymatrix{\xi~:~V \ar[r] & E(\xi) \ar[r]^-{\di{p_\xi}}&X}~,~
\xymatrix{\eta~:~W \ar[r] & E(\eta) \ar[r]^-{\di{p_\eta}}&Y}$$
is the $V\oplus W$-bundle over $X \times Y$
$$\xymatrix{\xi\times \eta~:~V\oplus W \ar[r] & E(\xi\times \eta)
\ar[r]^-{\di{p_{\xi\times \eta}}}&X\times Y}$$
with
$$p_{\xi\times \eta}~=~p_\xi \times p_\eta~:~
E(\xi \times \eta)~=~E(\xi) \times  E(\eta) \to X \times Y~,~
T(\xi \times \eta)~=~T(\xi) \wedge  T(\eta)~.$$
(iii) The {\it Whitney sum} of a $V$-bundle over $X$ and a
$W$-bundle over $X$
$$\xymatrix{\xi~:~V \ar[r] & E(\xi) \ar[r]^-{\di{p_\xi}}&X}~,~
\xymatrix{\eta~:~W \ar[r] & E(\eta) \ar[r]^-{\di{p_\eta}}&X}$$
is the pullback $V \oplus W$-bundle over $X$
$$\xymatrix{\xi \oplus \eta~=~\Delta^*(\xi \times \eta)~ :~
V \oplus W \ar[r] & E(\xi \oplus \eta)
\ar[r]^-{\di{p_{\xi \oplus \eta}}} & X}$$
with $\Delta:X \to X \times X;x \mapsto (x,x)$, such that
$$T(\xi\oplus \eta)~=~T(\xi)\wedge T(\eta)~.$$
The total space $E(\xi \oplus \eta)$ fits into a pullback square
$$\xymatrix{
E(\xi \oplus \eta) \ar[r] \ar[d]\ar[dr]^-{\di{p_{\xi\oplus \eta}}}
 & E(\xi) \ar[d]^-{\di{p_\xi}}\\
E(\eta) \ar[r]^-{\di{p_\eta}} & X}$$
and the Thom space of $\xi \oplus \eta$ is
$$T(\xi \oplus \eta)~=~T(p_{D(\xi)}^*\eta)/T(p_{S(\xi)}^*\eta)~.$$}
\hfill\qed
\end{definition}

\begin{proposition}~ \label{adjunct} Let
$$\xi~:~V \to E(\xi) \to M~,~\eta~:~W \to E(\eta) \to N$$
be vector bundles.
Regard $E(\xi)$ as a subspace of the disk space $D(\xi)$ via the open embedding
$$E(\xi) \emb D(\xi)~;~(m,v) \mapsto (m,\dfrac{v}{1+\Vert v \Vert})~(m \in M,v \in V)$$
with $D(\xi)\backslash E(\xi)=S(\xi)$,
and similarly for $S(\eta) \subset D(\eta)$. If
$$f~:~X~=~E(\xi) \emb Y~=~E(\eta)$$
is an open embedding which extends to an embedding
$\overline{f}:D(\xi) \emb D(\eta)$ there is defined a homeomorphism
$$T(\xi)~=~D(\xi)/S(\xi) \to Y/(Y\backslash f(X))~;~x \mapsto \overline{f}(x)$$
and the adjunction Umkehr of $f$ with respect to $E(\eta)\subset D(\eta)$
is a map
$$F~:~D(\eta)/S(\eta)~=~T(\eta) \to Y/(Y \backslash f(X))~=~T(\xi)~.$$
\end{proposition}
\begin{proof} By construction.
\hfill\qed\end{proof}

\begin{example}~\label{adex}
{\rm Let $V,W$ be finite dimensional inner product spaces.
For an open embedding of the type
$$e~:~E(\epsilon_W)~=~W \times M \emb E(\epsilon_V)~=~V\times N$$
which extends to an embedding
$$f~:~D(\epsilon_W)~=~D(W) \times M \emb D(\epsilon_V)~=~D(V) \times N$$
the adjunction Umkehr of \ref{adjunct} is a  map
$$F~:~T(\epsilon_V)~=~V^{\infty} \wedge N^\infty \to
T(\epsilon_W)~=~W^{\infty} \wedge M^+~.$$}
\hfill\qed
\end{example}

Given a (topological) group $G$ with a right action on a space $X$
and a left action on a space $Y$ let $X \times_GY$ be the quotient
of $X \times Y$ be the equivalence relation
$$(xg,y) \sim (x,gy)~{\rm for~all}~g \in G,x\in X,y \in Y~.$$
Similarly for pointed spaces, with $X \wedge_G Y$ the corresponding quotient
of $X \wedge Y$.

\begin{definition}~ {\rm Let $V,W$ be inner product spaces.\\
(i) The {\it orthogonal group} $O(V)$ is the group of linear isometries
$g:V \to V$, regarded as a pointed space with base point $1:V \to V$.
Define a left $O(V)$-action on $V$ by\index{orthogonal group, $O(V)$}
$$O(V) \times V \to V~;~(g,x) \mapsto g(x)~.$$
(ii) Let $O(V,W)$ be the {\it Stiefel manifold} of linear isometries
$i:V \emb W$. Define a right $O(V)$-action on $O(V,W)$ by
\index{Stiefel manifold, $O(V,W)$}
$$O(V,W) \times O(V) \to O(V,W)~;~(i,g) \mapsto ig~.$$
(iii) Let $G(V,W)$ be the {\it Grassmann manifold} of the subspaces
$U \subseteq W$ which are isomorphic to $V$.\index{Grassmann manifold, $G(V,W)$} \\
\hfill\qed}
\end{definition}

Here are some standard properties of the Stiefel and Grassmann manifolds:

\begin{proposition}~ {\rm (i)} For any finite-dimensional
inner product space $W$ and
subspace $V \subseteq W$ there is defined a homeomorphism
$$G(V,W) \to G(V^{\perp},W)~;~U \mapsto U^{\perp}~.$$
{\rm (ii)} For any two inner product spaces $U,V$
$O(U)\times O(V)$ has a right action on $O(U \oplus V)$
$$O(U \oplus V)\times (O(U)\times O(V)) \to O(U\oplus V)~;~
(f,g,h) \mapsto f(g \oplus h)~.$$
Regard $O(V,U \oplus V)$ as a pointed space
with base point the inclusion
$$j~:~ V \emb U \oplus V~;~x \mapsto (0,x)~.$$
The homeomorphisms
$$\begin{array}{l}
O(U \oplus V)/O(V) \to O(U,U\oplus V)~;~h \mapsto h\vert_U~,\\[1ex]
O(U \oplus V)/O(U) \to O(V,U\oplus V)~;~h \mapsto h\vert_V~,\\[1ex]
O(U \oplus V)/(O(U)\times O(V)) \to G(V,U\oplus V)~;~h \mapsto h(V)
\end{array}$$
give identifications
$$\begin{array}{l}
O(U,U\oplus V)~=~O(U\oplus V)/O(V)~,\\[1ex]
O(V,U\oplus V)~=~O(U\oplus V)/O(U)~,\\[1ex]
G(V,U\oplus V)~=~O(U \oplus V)/(O(U) \times O(V))~=~G(U,U\oplus V)~.
\end{array}$$
{\rm (iii)} For any finite-dimensional inner product spaces $U,V,W$
there is defined a commutative braid of fibre bundles

$$\xymatrix@C-40pt{
O(V)\ar[dr] \ar@/^2pc/[rr] && O(V,U\oplus V \oplus W)
\ar[dr]\ar@/^2pc/[rr] &&O(U,V,W)\\
&O(V,U\oplus V) \ar[ur] \ar[dr]&&G(V,U \oplus V \oplus W)\ar[ur]&\\
&&G(V,U \oplus V)\ar[ur]&&}$$
where $O(U,V,W)=O(U\oplus V \oplus W)/\langle O(U \oplus V),O(U\oplus W)\rangle$.\\
\hfill\qed
\end{proposition}

\begin{definition}~ {\rm
The {\it canonical $U$- and $V$-bundles} over the
Grassmann manifold $G=G(U,U\oplus V)=G(V,U\oplus V)$ are
$$\begin{array}{l}
\xi(U)~:~U \to E(\xi(U))~=~O(U,U\oplus V)\times_{O(U)}U \to
O(U,U \oplus V)/O(U)~=~G~,\\[1ex]
\xi(V)~:~V \to E(\xi(V))~=~O(V,U\oplus V) \times_{O(V)}V\to
O(V,U \oplus V)/O(V)~=~G
\end{array}$$
with Thom spaces
$$\begin{array}{l}
T(\xi(U))~=~O(U,U \oplus V)^+\wedge_{O(U)}U^{\infty}~,\\[1ex]
T(\xi(V))~=~O(V,U \oplus V)^+\wedge_{O(V)}V^{\infty}~.
\end{array}$$
An element $(W,x) \in E(\xi(U))$ (resp. $(W,y) \in E(\xi(V))$)
can be regarded as a subspace $W \subseteq U \oplus V$ isomorphic to $U$,
together with an element $x \in W$ (resp. $y \in W^{\perp}$).\\
\hfill\qed}
\end{definition}

Let $p_U:U \oplus V \to U$, $p_V:U \oplus V \to V$ be the projections.
For any $W \in G(U,U\oplus V)$
there exists $h \in O(U \oplus V)$ such that
$$h(U)~=~W~,~h(V)~=~W^{\perp} \subseteq U \oplus V~,$$
The maps
$$\begin{array}{l}
G \times (U \oplus V) \to E(\xi(U))= O(U\oplus V)/O(V)\times_{O(U)}U~;\\[1ex]
\hskip150pt (W,u,v) \mapsto (h,p_Uh^{-1}(u,v))\\[1ex]
G \times (U \oplus V) \to E(\xi(V))=O(U\oplus V)/O(U)\times_{O(V)}V~;\\[1ex]
\hskip150pt (W,u,v) \mapsto (h,p_Vh^{-1}(u,v))
\end{array}$$
fit into a pullback square
$$\xymatrix{G \times (U \oplus V) \ar[r] \ar[d] & E(\xi(U)) \ar[d] \\
E(\xi(V)) \ar[r] & G}$$
defining the canonical isomorphism
$$\xi(U) \oplus \xi(V)~\cong~\epsilon_{U \oplus V}~.$$

\begin{definition}~ {\rm  Let $V$ be a finite-dimensional inner
product space.\\
{\rm (i)} The {\it classifying space} for $V$-bundles is
$$BO(V)~=~\varinjlim\limits_U G(V,U\oplus V)~=~
\varinjlim\limits_U O(U\oplus V)/(O(U)\times O(V))$$
with $U$ running over finite-dimensional inner product spaces.\\
{\rm (ii)} The {\it universal $V$-bundle} is
$$\xi(V)~=~\varinjlim\limits_U \xi(V)~
:~V \to EO(V)~=~\varinjlim\limits_U O(V,U\oplus V)\times_{O(V)}V \to BO(V)~.$$
{\rm (iii)} The {\it universal Thom spectrum} is
$$\underline{MO}~=~\{MO(V)\,\vert\,{\rm dim}(V)<\infty\}$$
with $EO(V)$ a contractible space with free $O(V)$-action
$$MO(V)~=~T(\xi(V))~=~\varinjlim\limits_U O(V,U\oplus V)^+\wedge_{O(V)}V^{\infty}~.$$
\hfill\qed}
\end{definition}

\begin{terminology}~{\rm The Stiefel manifold
$$V_{j+k,k}~=~O(\R^k,\R^{j+k})~=~O(\R^{j+k})/O(\R^j)$$
of linear isometries $\R^k \emb \R^{j+k}$ is (homeomorphic to)
the space of orthogonal $k$-frames in $\R^{j+k}$.
The Grassmann manifold
$$G(\R^k,\R^{j+k})~=~V_{j+k,k}/O(\R^k)~=~O(\R^{j+k})/(O(\R^j) \times O(\R^k))$$
is (homeomorphic to)  the space of $k$-dimensional subspaces
$$W~=~{\rm im}(\R^k \emb \R^{j+k}) \subseteq \R^{j+k}~.$$
The canonical $\R^j$-bundle $\xi(\R^j)$ and the canonical
$\R^k$-bundle $\xi(\R^k)$ over the Grassmann manifold
$$\begin{array}{l}
\R^j \to E(\xi(\R^j)) \to G(\R^k,\R^{j+k})~,\\[1ex]
\R^k \to E(\xi(\R^k))\to G(\R^k,\R^{j+k})
\end{array}$$
are such that
$$\begin{array}{l}
E(\xi(\R^k))~=~V_{j+k,k}\times_{O(\R^k)}\R^k ~=~
\{(W \subseteq \R^{j+k},x \in W)\,\vert\,{\rm dim}(W)=k\}~,\\[1ex]
E(\xi(\R^j))~=~\{(W \subseteq \R^{j+k},y\in W^{\perp})
\,\vert\,{\rm dim}(W)=k\}~,\\[1ex]
\xi(\R^j) \oplus \xi(\R^k)~=~\epsilon_{\R^{j+k}}~.
\end{array}$$
Passing to the limit as $j \to \infty$ gives the universal $\R^k$-bundle
$$\begin{array}{l}
\R^k \to EO(\R^k)~=~\varinjlim\limits_j V_{j+k,k}\times_{O(\R^k)}\R^k\\[2ex]
\hskip100pt \to BO(\R^k)~=~\varinjlim\limits_j G(\R^k,\R^{j+k})~.
\end{array}$$
with Thom space $MO(\R^k)$ and Thom spectrum $\underline{MO}=\{MO(\R^k)\}$.\\
}
\hfill\qed
\end{terminology}

\begin{remark}~ \label{stiefel0}
{\rm  We recall some more standard facts about the Stiefel spaces.
(i) $V_{j+k,k}=O(\R^{j+k})/O(\R^j)$  is $(j-1)$-connected, with
$$\pi_i(V_{j+k,k})~=~\begin{cases} 0&{\rm if}~i<j\\
\ZZ&{\rm if}~i=j~\hbox{\rm is even, or if}~i=j~,~k=1\\
\ZZ_2&{\rm if}~i=j~\hbox{\rm is odd and}~k \geqslant 2
\end{cases}$$
(ii) $V_{j+k,k}$ fits into a fibre bundle
$$V_{j+k,k} \to BO(\R^j) \to BO(\R^{j+k})~.$$
The homotopy classes of maps $c:X \to V_{j+k,k}$ are in one-one correspondence
with the equivalence classes of pairs $(\xi,\delta\xi)$ with
$\xi:X \to BO(\R^j)$ an $\R^j$-bundle over $X$ and
$\delta\xi:\xi\oplus \epsilon_{\R^k} \cong \epsilon_{\R^{j+k}}$
a $k$-stable isomorphism, with
$$E(\xi)~=~\{(x \in X,y \in c(x)^{\perp} \subset \R^{j+k})\} \subset
E(\epsilon_{\R^{j+k}})~=~X \times \R^{j+k}~.$$
(iii) $V_{j+1,1}=S^j$, with $1:S^j \to V_{j+1,1}$ classifying the
tangent $\R^j$-bundle $\tau_{S^j}:S^j \to BO(\R^j)$ with the stable
isomorphism $\tau_{S^j} \oplus \epsilon_{\R} \cong \epsilon_{\R^{j+1}}$
determined by the standard framed embedding $S^j \subset \R^{j+1}$.\\
(iv) Let $j \geqslant 0$, $k \geqslant 1$.
For any $\R^{j+k}$-bundle $\xi:X \to BO(\R^{j+k})$ over a $CW$ complex $X$
it is possible to split
$\xi\vert_{X^{(j)}} \cong \xi' \oplus \epsilon_{\R^k}$
for some $\R^j$-bundle $\xi'$ over $X^{(j)}$.
The {\it $(j+1)$-th Stiefel-Whitney class} of $\xi$
$$w_{j+1}(\xi) \in H^{j+1}(X;\{\pi_j(V_{j+k,k})\})$$
is the primary obstruction to splitting
$\xi \cong \xi'' \oplus \epsilon_{\R^k}$
for some $\R^j$-bundle $\xi''$ over $X$. \\
(v) The {\it Euler class} of an $\R^{j+1}$-bundle $\xi:X \to BO(\R^{j+1})$
is the $(j+1)$th Stiefel-Whitney class
$$\gamma(\xi)~=~ w_{j+1}(\xi)\in H^{j+1}(X;\{\ZZ\})~,$$
the primary obstruction to splitting $\xi=\xi'\oplus \epsilon_\R$.\\
(vi) The $\bmod\,2$ Stiefel-Whitney class $w_{j+1}(\xi) \in H^{j+1}(X;\ZZ_2)$
of $\xi:X \to BO(\R^{j+k})$ is given by the evaluation of the Steenrod square
$$Sq^{j+1}~:~\dot H^{j+k}(T(\xi);\ZZ_2) \to \dot H^{2j+k+1}(T(\xi);\ZZ_2)$$
on the $\bmod\,2$ Thom class $U_\xi \in\dot H^{j+k}(T(\xi);\ZZ_2)$, with
$$Sq^{j+1}(U_\xi)~=~U_\xi \cup w_{j+1}(\xi) \in\dot H^{2j+k+1}(T(\xi);\ZZ_2)~.$$
\hfill\qed}
\end{remark}

The projection
$$p~:~O(V,U\oplus V) \to G(U,U\oplus V)~=~O(V,U\oplus V)/O(V)$$
is such that there is defined a pullback square
$$\xymatrix{
O(V,U \oplus V)  \times V\ar[r] \ar[d]^-{\di{p^*\xi(V)}} &
E(\xi(V))=O(V,U\oplus V)\times_{O(V)}V \ar[d]^-{\di{\xi(V)}}\\
O(V,U \oplus V) \ar[r]^-{\di{p}} & O(V,U\oplus V)/O(V)}$$
defining an isomorphism of $V$-bundles over $O(V,U \oplus V)$
$$p^*\xi(V)~\cong~\epsilon_V~.$$
The pullback
$$p^*\xi(U)~:~
U \to O(U \oplus V)\times_{O(U)}U \to O(U \oplus V)/O(U)~=~O(V,U\oplus V)$$
is such that
$$\begin{array}{l}
O(U \oplus V)\times_{O(U)}(U \oplus V) \to
O(V,U\oplus V)\times (U\oplus V)~;\\[1ex]
\hskip125pt (h,u,v) \mapsto (h\vert_V,p_Uh^{-1}(u),p_Vh^{-1}(u)+v)
\end{array}$$
defines an isomorphism of $U\oplus V$-bundles over $O(V,U \oplus V)$
$$p^*\xi(U)\oplus \epsilon_V~\cong~ \epsilon_{U \oplus V}~,$$
corresponding to the fibration sequence
$$\xymatrix{
O(U) \ar[r] &O(U \oplus V) \ar[r] &
O(V,U \oplus V) \ar[r]^-{\di{p}} &BO(U) \ar[r] &BO(U\oplus V)~.}$$
\begin{proposition}~Let $X$ be a reasonable space, such as a finite
$CW$ complex.\\
{\rm (i)} The isomorphism classes of $V$-bundles over $X$
$$\xymatrix{\xi~:~V \ar[r] & E(\xi) \ar[r]&X}$$
are in bijective correspondence with the set $[X,BO(V)]$ of
homotopy classes of maps
$$f~:~X \to G(V,U\oplus V)\subset BO(V)~~({\rm dim}(U)~{\it large})$$
with
$$\begin{array}{l}
\xi~=~f^*\xi(V)~,\\[1ex]
E(\xi)~=~f^*E(\xi(V))~=~\{(u,v,x) \in U\oplus V \times X\,\vert\,
(u,v) \in f(x) \subseteq U\oplus V\}~.
\end{array}$$
The function
$$\begin{array}{l}
[X,BO(V)] \to \{\hbox{isomorphism classes of $V$-bundles over $X$}\}~;\\[1ex]
\hskip150pt (f:X \to BO(V)) \mapsto f^*\xi
\end{array}$$
is a bijection.\\
{\rm (ii)} The isomorphism classes of pairs $(\xi,\delta\xi)$
with $\xi$ a $V$-bundle over $X$ and
$$\delta\xi~:~\xi \oplus \epsilon_U~\cong~\epsilon_{U \oplus V}$$
a $U \oplus V$-bundle isomorphism are in bijective
correspondence with the homotopy classes of maps
$g:X \to O(U,U \oplus V)$, with
$$\begin{array}{l}
\xymatrix{\xi=(pg)^*\xi(V)~:~X \ar[r]^-{\di{g}} &
O(U,U \oplus V) \ar[r]^-{\di{p}} & G(V,U\oplus V) \subset BO(V)}~,\\[1ex]
E(\xi)~=~(pg)^*E(\xi(V))~=~\{(u,v,x) \in U\oplus V \times X \,\vert\,
(u,v) \in pg(x) \subseteq U \oplus V\}~,\\[1ex]
\delta\xi ~:~E(\xi \oplus \epsilon_U)~=~U \times E(\xi)
 \to E(\epsilon_{U \oplus V})~=~U \oplus V \times X ~;\\[1ex]
\hskip75pt (t,u,v,x) \mapsto (t+u,v,x)~((u,v) \in g(x)(U)^{\perp} \subseteq U \oplus V)~.
\end{array}$$
\hfill\qed
\end{proposition}
\begin{example}
{\rm Let $M^m$ an $m$-dimensional manifold with an embedding
$M \subset \R^m\oplus V$, for some inner product space $V$.
The tangent $\R^m$-bundle $\tau_M$ and the normal $V$-bundle
$\nu_M=\nu_{M \subset \R^m \oplus V}$ are classified by
$$\begin{array}{l}
\tau_M~:~M \to G(\R^m,\R^m\oplus V) \subset BO(\R^m)~,\\[1ex]
\nu_M~=~\tau^{\perp}_M~:~M \to G(V,\R^m\oplus V) \subset BO(V)
\end{array}$$
with
$$\tau_M \oplus \nu_M~=~\tau_{\R^m\oplus V}\vert_M~=~\epsilon_{\R^m\oplus V}~:~
M \to G(\R^m\oplus V,\R^m\oplus V) \subset BO(\R^m\oplus V)~.$$
\hfill\qed}
\end{example}

\begin{example}~ \label{tangentsphere}
{\rm Let $V=U \oplus \R$, for some finite-dimensional inner product space $U$.
The homeomorphism
$$S(V) \to O(\R,V)~;~ x \mapsto (t \mapsto tx)$$
is the clutching function of the tangent $U$-bundle $\tau_{S(V)}:S(V) \to BO(U)$
$$U \to E(\tau_{S(V)}) \to S(V)$$
with sphere bundle
$$S(\tau_{S(V)})~=~\{(x,y) \in V \oplus V\,\vert\, \Vert x \Vert =1,\,
\langle x,y \rangle=0\}~.$$
The open embedding
$$S(V)\times \R \emb V~;~(x,t) \mapsto e^tx$$
corresponds to the $V$-bundle isomorphism
$$\delta\tau_{S(V)}~:~\tau_{S(V)}\oplus \nu_{S(V) \subset V}~=~
\tau_{S(V)}\oplus\epsilon_{\R}~\cong~\tau_V\vert~=~\epsilon_V$$
classified by $S(V) \to O(\R,V)=O(V)/O(U)$.}\\
\hfill\qed
\end{example}

\begin{definition}~\label{virtual}
{\rm Let $\xi:X \to BO(V)$, $\xi':X \to BO(V')$
be vector bundles over a space $X$, for some finite-dimensional
inner product spaces $V,V'$.\\
(i) The {\it dimension} of the virtual bundle $\xi-\xi'$ is
$${\rm dim}(\xi-\xi')~=~{\rm dim}(V) - {\rm dim}(V') \in \ZZ~.$$
(ii) The {\it virtual Thom space} of $\xi-\xi'$ is the spectrum
$$\underline{T}(\xi-\xi')~=~\{T(\xi-\xi')_U\,\vert\,U\}$$
defined by
$$T(\xi-\xi')_U~=~T(\xi \oplus \eta)$$
for finite-dimensional inner product spaces $U$ with an isometry $V'\to U$ and a vector bundle
$\eta$ over $X$ such that $\xi'\oplus \eta=\epsilon_U$. Note that
$${\rm dim}(\xi-\xi')~=~{\rm dim}(\xi \oplus \eta) - {\rm dim}(U) \in \ZZ~.$$
(iii) Write
$$\omega_*(X;\xi-\xi')~=~\widetilde{\omega}_*(\underline{T}(\xi-\xi'))~,~
\omega^*(X;\xi-\xi')~=~\widetilde{\omega}^*(\underline{T}(\xi-\xi'))$$
for the reduced stable homotopy and cohomotopy of the virtual Thom space
of $\xi-\xi'$, so that
$$\begin{array}{l}
\omega_*(X;\xi-\xi')~=~\widetilde{\omega}_{*+{\rm dim}(U)}(T(\xi\oplus \eta))~,\\[1ex]
\omega^*(X;\xi-\xi')~=~\widetilde{\omega}^{*+{\rm dim}(U)}(T(\xi\oplus \eta))~.
\end{array}$$
}
\hfill\qed
\end{definition}

\begin{example} {\rm
(i) In the special case $\xi'=\epsilon_{\R^{\ell}}$ ($\ell \geqslant 0$)
the virtual Thom space $\underline{T}(\xi-\xi')$
is just the $\ell$-fold desuspension of the
suspension spectrum $\underline{T}(\xi)$ of the actual Thom space $T(\xi)$
$$\begin{array}{ll}
\underline{T}(\xi-\xi')&=~
\Sigma^{-\ell}\{T(\xi\oplus \epsilon_U)\,\vert\, U\}~,\\[1ex]
&=~\{T(\xi\oplus \epsilon_{U/\R^{\ell}})\,\vert\, U, \R^{\ell} \subseteq U\}~,
\end{array}$$
and
$$\begin{array}{l}
\omega_*(X;\xi-\epsilon_{\R^{\ell}})~=~\widetilde{\omega}_{*+\ell}(T(\xi))~,\\[1ex]
\omega^*(X;\xi-\epsilon_{\R^{\ell}})~=~\widetilde{\omega}^{*+\ell}(T(\xi))~.
\end{array}$$
(iii) For a virtual trivial bundle $\xi-\xi'=\epsilon_{\R^k}-\epsilon_{\R^{\ell}}$
over a compact space $X$
$$\begin{array}{l}
\omega_*(X;\epsilon_{\R^k}-\epsilon_{\R^{\ell}})~=~\omega_{*-k+\ell}(X)~,\\[1ex]
\omega^*(X;\epsilon_{\R^k}-\epsilon_{\R^{\ell}})~=~\omega^{*-k+\ell}(X)~.
\end{array}$$
}\hfill\qed
\end{example}

\begin{proposition}~ {\rm (The tubular neighbourhood theorem).}\\
{\rm (i)} Every embedding of manifolds $M^m \emb N^n$ has a normal  $\R^{n-m}$-bundle
$$\nu_{M \subset N}~:~M \to BO(n-m)$$
such that
$$\tau_M \oplus \nu_{M \subset N}~=~\tau_N\vert_M~:~M \to BO(\R^n)~,$$
with a codimension 0 embedding $E(\nu_{M \subset N}) \subset N$
and hence an Umkehr map
$$F~:~N^{\infty} \to N/(N \backslash E(\nu_{M \subset N}))~=~T(\nu_{M \subset N})~.$$
{\rm (ii)} Every immersion of manifolds $f:M^m \imm N^n$ has a normal $\R^{n-m}$-normal $W$-bundle
$$\nu_{M \imm N}~:~\R^n \to E(\nu_{M\imm N}) \to M~,$$
with a codimension 0 immersion $E(\nu_{M \imm N}) \imm N$.
For any inner product space $V$ with ${\rm dim}(V)=j$ the immersion
$$(0,f)~:~M \imm V \times N~;~x \mapsto (0,f(x))$$
has codimension $j+n-m$, and normal bundle
$\nu_{(0,f)}=\nu_{M \imm N}\oplus \epsilon_V$.
For $j \geqslant 2m-n+1$ $(0,f)$ is regular homotopic to an
embedding of the form
$$(e,f)~:~M \emb V \times N~;~x \mapsto (e(x),f(x))$$
for some $e:M \to V$, with $\nu_{(e,f)}=\nu_{M \imm N}\oplus \epsilon_V$.
There is an extension of $(e,f)$ to an open embedding
$$E(\nu_{M \imm N}\oplus \epsilon_V) \emb E(\epsilon_V)~=~V \times N$$
with adjunction Umkehr map
$$F~:~T(\epsilon_V)~=~V^{\infty} \wedge N^\infty \to
T(\nu_{M \imm N}\oplus \epsilon_V)~=~V^{\infty} \wedge T(\nu_{M \imm N})~.$$
\hfill\qed
\end{proposition}

\begin{example} {\rm For any $m$-dimensional
manifold $M$ there exists an embedding $M \subset \R^m\oplus V$,
for some finite-dimensional inner product space $V$.
By the tubular neighbourhood theorem the embedding has a
normal $V$-bundle  $\nu_M$ with an embedding $E(\nu_M) \subset \R^m\oplus V$.
The composite of the Umkehr map
$$\begin{array}{l}
\alpha_{\R^m \oplus V}~:~(\R^m\oplus V)^{\infty}\\[1ex]
\hskip50pt
\to (\R^m\oplus V)/(\R^m\oplus V \backslash E(\nu_M))~=~D(\nu_M)/S(\nu_M)~=~T(\nu_M)
\end{array}$$
and the diagonal map $\Delta:T(\nu_M) \to M^+ \wedge T(\nu_M)$
$$\Delta \alpha_{\R^m \oplus V}~:~(\R^m\oplus V)^{\infty} \to M^+ \wedge T(\nu_M)$$
represents an element
$$(M,1)~=~\Delta \alpha_{\R^m \oplus V} \in \omega_m(M,\nu_M-\epsilon_V)~=~
\widetilde{\omega}_{m+{\rm dim}(V)}(T(\nu_M))~.$$
}
\hfill\qed
\end{example}

\section{Bordism}\label{bordism}

Every $m$-dimensional manifold $M$ admits an embedding
$M \subset \R^{m+n}$ ($n$ large). The tangent
$\R^m$-bundle $\tau_M:M \to BO(\R^m)$ and the normal $\R^n$-bundle
$\nu_{M \subset \R^{m+n}}:M \to BO(\R^n)$ are such that
$$\tau_M \oplus \nu_M~=~\tau_{\R^{m+n}}\vert_M~=~\epsilon_{\R^{m+n}}~:~
M \to BO(\R^{m+n})~.$$
The Pontrjagin-Thom Umkehr map $\alpha_M:S^{m+n} \to T(\nu_M)$
is transverse regular at the zero section $M \subset T(\nu_M)$ with
$$(\alpha_M)^{-1}(M)~=~M \subset \R^{m+n} \subset (\R^{m+n})^{\infty}~=~S^{m+n}~.$$

\begin{definition}~ {\rm
Let $\xi-\xi'$ be a virtual bundle over a space $X$, represented by
bundles $\xi,\xi'$ with ${\rm dim}(\xi)=i$, ${\rm dim}(\xi')=i'$.\\
(i) A {\it $m$-dimensional normal map}
$$(f,b)~:~(M,\nu_M-\epsilon_{\R^n}) \to (X,\xi-\xi')$$
consists of an $m$-dimensional manifold $M$ with
an embedding $M \subset \R^{m+n+i-i'}$ ($n$ large), a map $f:M \to X$
and a bundle map over $f$
$$b~:~\nu_M\oplus f^*\xi' \to \xi\oplus\epsilon_{\R^n}$$
with
$$\nu_M~=~\nu_{M \subset \R^{m+n+i-i'}}~:~M \to BO(n+i-i')~.$$
Equivalently, $b$ can be regarded as a bundle map
$$b~:~\tau_M \oplus f^*\xi\oplus \epsilon_{\R^n} \to
\xi'\oplus \epsilon_{\R^{m+n+i-i'}}$$
and as a virtual bundle map $b:\nu_M-\epsilon_{\R^n}\to \xi-\xi'$.
Let $\eta:M \to BO(j)$ be a bundle such that $\xi' \oplus \eta=\epsilon_{\R^{i'+j}}$,
so that $M \subset \R^{m+n+i+j}$ has normal bundle
$$\begin{array}{ll}
\nu_{M \subset \R^{m+n+i+j}}&=~\nu_M \oplus \epsilon_{\R^{i'+j}}\\[1ex]
&=~\nu_M \oplus f^*\xi' \oplus f^*\eta\\[1ex]
&=~f^*(\xi\oplus \eta \oplus \epsilon_{\R^n})~:~M \to BO(n+i+j)
\end{array}$$
allowing the {\it Pontrjagin-Thom map} to be defined by
$$\begin{array}{l}
\xymatrix{\alpha(f,b)~:~S^{m+n+i+j} \ar[r]^-{\di{\alpha_M}}&
\Sigma^{i'+j}T(\nu_M)~=~T(\nu_{M \subset \R^{m+n+i+j}})}\\[1ex]
\hskip175pt \xymatrix{\ar[r]^-{\di{T(b)}} &T(\xi\oplus \eta \oplus \epsilon_{\R^n})~.}
\end{array}$$
(ii) Let $\Omega_m(X,\xi-\xi')$ be the bordism group of $m$-dimensional normal maps
$(f,b):(M,\nu_M-\epsilon_{\R^n}) \to (X,\xi-\xi')$, with {\it Pontrjagin-Thom} map
$$\begin{array}{l}
PT~:~\Omega_m(X,\xi-\xi') \to\omega_{m+i-i'}(X;\xi-\xi')=
\{S^{m+n+i+j};T(\xi\oplus \eta \oplus \epsilon_{\R^n})\}~;\\[1ex]
\hskip 175pt (f,b) \mapsto \alpha(f,b)~.
\end{array}$$
(iii) Let $\fN_m(X)$ be the bordism group of $m$-dimensional manifolds $M$ with a map $f:M \to X$.\\
\hfill\qed}
\end{definition}

\begin{proposition}~ \label{Poincare}
{\rm (i)} The Pontrjagin-Thom maps
$$PT~:~\Omega_m(X,\xi-\xi') \to \omega_{m+i-i'}(X;\xi-\xi')~(m \geqslant 0)$$
are isomorphisms, with inverses given by the transversality construction.
Every element
$$F \in \omega_{m+i-i'}(X;\xi-\xi')~=~\{S^{m+n+i+j};T(\xi\oplus \eta \oplus \epsilon_{\R^n})\}$$
is represented by a map $F:S^{m+n+i+j}\to T(\xi\oplus \eta\oplus \epsilon_{\R^n})$
which is transverse regular at the zero section
$X \subset T(\xi\oplus \eta\oplus \epsilon_{\R^n})$ with the restriction
$$f~=~F\vert~:~M^m~=~F^{-1}(X) \to X$$
such that
$$\begin{array}{ll}
\nu_{M \subset \R^{m+n+i+j}} \oplus f^*\xi'&=~
f^*(\xi\oplus \xi' \oplus \eta\oplus \epsilon_{\R^n})\\[1ex]
&=~f^*\xi\oplus \epsilon_{\R^{n+i'+j}}~:~M \to BO(\R^{n+i+j})~,
\end{array}$$
and
$$T(b)\alpha_M~=~F \in \omega_{m+i-i'}(X;\xi-\xi')~.$$
{\rm (ii)} The maps
$$\begin{array}{l}
\fN_m(X) \to \varinjlim\limits_V\,\Omega_m(X\times BO(V),0\times \xi(V)-\epsilon_V)~;\\[1ex]
(f:M \to X) \mapsto (f \times \nu_M:(M,\nu_M-\epsilon_V) \to (X\times BO(V),\xi(V)-\epsilon_V))
\end{array}$$
are isomorphisms. The Pontrjagin-Thom maps
$$\begin{array}{l}
PT~:~\fN_m(X)~=~\varinjlim\limits_V \Omega_m(X\times BO(V),0 \times \xi(V)) \to\\[1ex]
\hskip100pt MO_m(X)~=~\varinjlim\limits_V [\Sigma^mV^{\infty};X^+ \wedge MO(V)]
\end{array}$$
are isomorphisms.\\
\hfill\qed
\end{proposition}

\begin{definition} {\rm
(i) A {\it framed $m$-dimensional manifold} is an $m$-dimensional
manifold $M$ with an embedding $M \subset \R^{m+n}$ and a trivialization
$\nu_M\cong \epsilon_{\R^n}$ of the normal bundle $\nu_M:M \to BO(n)$.
The  {\it Pontrjagin-Thom Umkehr map}
$$\alpha_M~:~(\R^{m+n})^{\infty}~=~S^{m+n} \to T(\nu_M)~=~\Sigma^nM^+$$
represents an element
$$\alpha_M \in \{S^{m+n};\Sigma^nM^+\}~=~\omega_m(M)~.$$
{\rm (ii)} The {\it $m$-dimensional framed  bordism
group} of a space $X$
$$\Omega^{fr}_m(X)~=~\Omega_m(X,0)~~(\xi=\epsilon_V)$$
is the bordism group of $m$-dimensional framed manifolds $M$
with a map $f:M \to X$, and the {\it framed Pontrjagin-Thom map}  is
$$PT^{fr}~:~\Omega^{fr}_m(X) \to\omega_m(X)~;~
((f,b):M \to X) \mapsto f_*\alpha_M~.$$
\hfill\qed}
\end{definition}

\begin{proposition}~
The framed Pontrjagin-Thom maps are isomorphisms
$$PT^{fr}~:~\Omega^{fr}_m(X) \xymatrix{\ar[r]^-{\di{\cong}}&}\omega_m(X)~.$$
\end{proposition}
\begin{proof} This is just the special case $\xi-\xi'=0$ of Proposition
\ref{Poincare}.\\
\hfill\qed\end{proof}

\begin{example} {\rm In particular, the Pontrjagin-Thom isomorphism
identifies $\Omega_m^{fr}$ with $\omega_m=\{S^m;S^0\}$.}\\
\hfill\qed
\end{example}

As usual, there are also relative bordism groups $\fN_m(X,Y)$ for a
pair of spaces $(X,Y \subseteq X)$, with elements the bordism classes
of maps $f:(M,\partial M) \to (X,Y)$ from $m$-dimensional manifolds
with boundary. The relative bordism groups fit into an exact sequence
$$\dots \to \fN_m(Y) \to \fN_m(X) \to \fN_m(X,Y) \to \fN_{m-1}(Y) \to \dots$$
and the natural maps
$$\begin{array}{l}
\fN_m(X,Y)~=~MO_m(X,Y) \to\\[1ex]
\hskip50pt
\widetilde{\fN}_m(X/Y)~=~
\widetilde{MO}_m(X/Y)~=~\varinjlim\limits_V [\Sigma^mV^{\infty};X/Y \wedge MO(V)]
\end{array}$$
are isomorphisms. For any $\R^i$-bundle $\xi:X \to BO(i)$ there are
defined Thom isomorphisms
$$\begin{array}{l}
\fN_m(X) \to \fN_{m+i}(D(\xi),S(\xi))~=~\widetilde{\fN}_{m+i}(T(\xi))~;\\[1ex]
\hskip100pt (M,f:M \to X) \mapsto (D(f^*\xi),S(f^*\xi))~.
\end{array}$$

\begin{example} {\rm
(i) For pointed spaces $X,Y$ a stable map
$F:V^{\infty}\wedge X \to V^{\infty}\wedge Y$ induces morphisms
$$F_*~:~\widetilde{\fN}_m(X)~=~\widetilde{MO}_m(X) \to \widetilde{\fN}_m(Y)~=~\widetilde{MO}_m(Y)$$
(ii) For a map $F:X \to Y$ of unpointed spaces the corresponding pointed map
$F^+:X^+ \to Y^+$ is such that
$$F^+_*~:~\fN_m(X) \to \fN_m(Y)~;~(M,f:M \to X) \mapsto (M,Ff:M \to Y)~.$$
(iii) The Umkehr map $F:V^{\infty} \wedge N^\infty \to V^{\infty} \wedge T(\nu_f)$
of an immersion $f:M^m \imm N^n$ is such that
$$\begin{array}{l}
F_*~:~\fN_n(N) \to \widetilde{\fN}_n(T(\nu_f))~=~\fN_n(D(\nu_f),S(\nu_f))~=~
\fN_m(M)~;\\[1ex]
\hskip100pt (N,1) \mapsto (N,\emptyset,z)~=~(M,1)
\end{array}
$$
with $z:N \to D(f)$ the zero section.\\
\hfill\qed}
\end{example}

\section{$S$-duality}\label{Sduality}

\begin{definition}~{\rm
Let $X,Y$ be pointed spaces.\\
(i) Define the {\it slant products} in stable (co)homotopy
$$\begin{array}{l}
\widetilde{\omega}_N(X \wedge Y) \otimes \widetilde{\omega}^i(X)
\to \widetilde{\omega}_{N-i}(Y)~;\\[1ex]
\hskip20pt (\sigma:S^N \to X \wedge Y)\otimes
(f:X \to S^i) \mapsto ((1\wedge f)\sigma:S^N \to \Sigma^iY)~,\\[1ex]
\widetilde{\omega}^N(X \wedge Y) \otimes \widetilde{\omega}_i(X)
\to \widetilde{\omega}^{N-i}(Y)~;\\[1ex]
\hskip20pt (\sigma^*:X \wedge Y \to S^N)\otimes
(f:S^i \to X) \mapsto (\sigma^*(f\wedge 1):S^i \wedge Y \to S^N)
\end{array}$$
for $i,N \in \ZZ$.\\
(ii) An element $\sigma \in \widetilde{\omega}_N(X \wedge Y)$ is an
{\it $S$-duality} if the products
$$\sigma\otimes-~:~\widetilde{\omega}^i(X)\to \widetilde{\omega}_{N-i}(Y)~~(i \in\ZZ)$$
are isomorphisms.\\
(iii) An element $\sigma^* \in \widetilde{\omega}^N(X \wedge Y)$ is a
{\it reverse $S$-duality} if the products
$$\sigma^*\otimes-~:~\widetilde{\omega}_i(X)\to \widetilde{\omega}^{N-i}(Y)~~(i \in\ZZ)$$
are isomorphisms.\\
\hfill\qed}
\end{definition}

\begin{proposition}\label{Sdual} ~{\rm (i)}
If $\sigma \in \widetilde{\omega}_N(X\wedge Y)$ is an
$S$-duality there are induced isomorphisms
$$\begin{array}{l}
\sigma~:~\{X \wedge A;B\} \xymatrix{\ar[r]^-{\di{\cong}}&}
\{\Sigma^NA;B \wedge Y \}~;~F \mapsto (F \wedge 1_Y)(\sigma\wedge 1_A)~,\\[1ex]
\sigma~:~\{A\wedge Y;B\} \xymatrix{\ar[r]^-{\di{\cong}}&}
\{\Sigma^NA;X \wedge B \}~;~
G \mapsto (1_X \wedge G)(1_A\wedge \sigma)
\end{array}$$
for any pointed $CW$ complexes $A,B$.\\
{\rm (ii)} If $\sigma^* \in \widetilde{\omega}^N(X\wedge Y)$ is a
reverse $S$-duality there are induced isomorphisms
$$\begin{array}{l}
\sigma^*~:~\{Z;X \wedge A\} \xymatrix{\ar[r]^-{\di{\cong}}&}
\{Z \wedge Y;\Sigma^NA\}~;~
F \mapsto (\sigma^*\wedge 1_A)(F\wedge 1_Y)~,\\[1ex]
\sigma^*~:~\{Z;A\wedge Y\} \xymatrix{\ar[r]^-{\di{\cong}}&}
\{X \wedge Z;\Sigma^NA\}~;~G \mapsto (\sigma^* \wedge 1_A)(1_X \wedge G)
\end{array}$$
for any pointed $CW$ complexes $A,Z$.\\
{\rm (iii)} If $X,Y$ are finite pointed $CW$ complexes
an element $\sigma \in \widetilde{\omega}_N(X \wedge Y)$ is
an $S$-duality if and only if the products
$$[\sigma] \otimes - ~:~\widetilde{H}^*(X) \to \widetilde{H}_{N-*}(Y)$$
are isomorphisms, with $[\sigma] \in \widetilde{H}_N(X \wedge Y)$ the
Hurewicz image. Similarly for reverse $S$-duality.\\
{\rm (iv)} If $X,Y$ are finite pointed $CW$ complexes there exists
an $S$-duality $\sigma \in \widetilde{\omega}_N(X \wedge Y)$ if
and only if there exists a reverse
$S$-duality $\sigma^* \in \widetilde{\omega}^N(X \wedge Y)$,
with $\sigma^*\sigma = 1 \in \{X;X\}$.\\
{\rm (v)} For any finite pointed $CW$ complex $X$ there exists
a finite pointed $CW$ complex $Y$ with an $S$-duality
$\sigma \in \widetilde{\omega}_N(X \wedge Y)$ and a reverse
$S$-duality $\sigma^* \in \widetilde{\omega}^N(X \wedge Y)$,
with $\sigma^*\sigma = 1 \in \{X;X\}$.\\
{\rm (vi)} Let $M$ be an $m$-dimensional manifold with
an embedding $M \subset \R^m\oplus V$ for $V=\R^i$, and let $\nu_M$ be the normal $V$-bundle.
The composite of the Umkehr map of $E(\nu_M) \subset \R^m \oplus V$
$$\begin{array}{l}
\alpha_{\R^m \oplus V}~:~(\R^m\oplus V)^{\infty}\\[1ex]
\hskip50pt
\to (\R^m\oplus V)/(\R^m\oplus V \backslash E(\nu_M))~=~D(\nu_M)/S(\nu_M)~=~T(\nu_M)
\end{array}$$
and the diagonal map $\Delta:T(\nu_M) \to M^+ \wedge T(\nu_M)$
represents an $S$-duality map
$$\Delta\alpha_{\R^m \oplus V} \in
\omega_m(M \times M;0 \times \nu_M - \epsilon_{\R^i})~=~
\widetilde{\omega}_{m+i}(M^+\wedge T(\nu_M))~.$$
\end{proposition}
\begin{proof} Standard.
\hfill\qed\end{proof}

\begin{proposition}~ \label{manifoldSduality}
Let $M$ be an $m$-dimensional manifold with an embedding $M
\subset \R^{m+n}$ with normal $\R^n$-bundle $\nu_M$.
For any bundles $\xi:M \to BO(\R^i)$, $\eta:M \to BO(\R^{n-i})$ such that
$$\xi \oplus \eta~=~\nu_M~:~M \to BO(\R^n)$$
the composite
$$\xymatrix{
\sigma~=~\Delta \alpha_M~:~S^{m+n} \ar[r]^-{\di{\alpha_M}} &
T(\nu_M) \ar[r]^-{\di{\Delta}} &T(\xi) \wedge T(\eta)}$$
is an $S$-duality map such that products with the $S$-duality
$$\sigma \in \omega_m(M \times M;\xi \times\eta - \epsilon_{\R^n})~=~
\widetilde{\omega}_{m+n}(T(\xi)\wedge T(\eta))$$
give Poincar\'e duality isomorphisms in stable homotopy
$$\sigma~:~\omega^*(M;\eta-\epsilon_{\R^{n-i}})~
\cong~\omega_{m-*}(M;\xi-\epsilon_{\R^i})~=~\Omega_{m-*}(M,\xi)~.$$
\hfill\qed
\end{proposition}

\begin{example} {\rm
Take $\xi=\nu_M$, $\eta=0$ in Proposition \ref{manifoldSduality},
so that $i=n$ and the $S$-duality map is
$$\xymatrix{
\sigma~=~\Delta \alpha_M~:~S^{m+n} \ar[r]^-{\di{\alpha_M}} &
T(\nu_M) \ar[r]^-{\di{\Delta}} &M^+ \wedge T(\nu_M)}$$
Products with the $S$-duality $\sigma \in \widetilde{\omega}_{m+n}
(M^+ \wedge T(\nu_M))$ define Poincar\'e duality isomorphisms
$$\sigma~:~\omega^*(M)~\cong~\omega_{m-*}(M;\nu_M-\epsilon_{\R^n})~=~
\Omega_{m-*}(M,\nu_M)~.$$
In particular, $1 \in \omega^0(M)$ is the Poincar\'e dual of
$$\sigma(1)~=~\alpha_M~=~(M,1) \in \omega_m(M;\nu_M-\epsilon_{\R^n})~=~\Omega_m(M,\nu_M)~.$$
\hfill\qed}
\end{example}

\begin{example}{\rm
(i) For any finite $CW$ complex $L$ there exists
an embedding $L \subset S^N$ ($N \geqslant 2\,{\rm dim}(L)+1$).
Let $W \subset S^N$ be a closed regular neighbourhood of $L\subset S^N$,
so that the projection $p:W \to L$ is a homotopy equivalence with
contractible point inverses.  For any subcomplex $K \subseteq L$ the
inverse image
$$V~=~p^{-1}(K) \subset W$$
is a closed regular neighbourhood of $K\subset S^N$, with $p\vert:V \to
K$ a homotopy equivalence with contractible point inverses.
The codimension 0 submanifolds
$$U~=~\overline{W \backslash V}~,~V \subset W$$
have boundaries
$$\partial U~=~\partial_VU \cup \partial_WU~,~\partial V~=~\partial_VU \cup \partial_WV$$
with
$$\begin{array}{l}
\partial_WU~=~U \cap \partial W \subseteq U~,~
\partial_WV~=~V \cap \partial W \subseteq V~,\\[1ex]
\partial_VU~=~U \cap V~,~U/\partial_VU~=~W/V~\simeq~L/K~.
\end{array} $$
The composite
$$\xymatrix@C-5pt{\sigma~:~S^N \ar[r] & S^N/(S^N\backslash
U)=U/\partial U \ar[r]^-{\di{\Delta}} & U/\partial_VU \wedge
U/\partial_WU~\simeq~L/K \wedge U/\partial_WU}$$
represents an $S$-duality
$\sigma \in \widetilde{\omega}_N((L/K) \wedge (U/\partial_WU))$.\\
(ii) The special case
$K=\{*\} \subset L$ in (i) gives an $S$-duality
$\sigma\in \widetilde{\omega}_N(L \wedge (U/\partial_WU)$
for any finite pointed $CW$ complex $L$.\\
(iii) For a finite unpointed $CW$ complex $L$ the special case
$K=\emptyset \subset L$ (with $V=\emptyset$,
$(U,\partial_WU)=(W,\partial W)$) in (i) gives an $S$-duality map
$$\sigma~:~S^N \to L^+ \wedge W/\partial W~.$$
If $L$ is an $n$-dimensional geometric Poincar\'e complex
(e.g. an $n$-dimensional manifold) then
$$S^{N-n-1} \to \partial W \to W~\simeq~L$$
is the Spivak normal fibration $\nu_L$, and
$$\sigma\in \widetilde{\omega}_N(L^+ \wedge W/\partial W)~=~
\widetilde{\omega}_N(L^+\wedge T(\nu_L))$$
gives the Atiyah-Wall $S$-duality between $L^+$ and the Thom space
$T(\nu_L)=W/\partial W$.\\
\hfill\qed}
\end{example}

\begin{proposition}~\label{S-dual} \label{difcon5}
Let $V$ be a finite-dimensional inner product space.\\
{\rm (i)} The composite of the projection
$$\begin{array}{ll}
\alpha_V~:&V^{\infty} \to V^{\infty}/0^+~=~\Sigma S(V)^+~;\\[1ex]
&v=[t,u]=\di{tu \over 1-t} \mapsto
\begin{cases}
(t,u)=(\di{\Vert v \Vert \over 1+\Vert v\Vert},
\di{v \over \Vert v\Vert})&\hbox{\it if}~v \neq 0,\infty\\
\infty&\hbox{\it if}~v=0~{\it or}~\infty
\end{cases}
\end{array}$$
and the diagonal map
$$\Delta~:~\Sigma S(V)^+ \to S(V)^+ \wedge \Sigma S(V)^+~;~
(t,u) \mapsto (u,(t,u))$$
is an $S$-duality map
$$\begin{array}{l}
\sigma_V~=~\Delta\alpha_V~:~V^{\infty}
\xymatrix{\ar[r]^-{\di{\alpha_V}}&}
\Sigma S(V)^+\xymatrix{\ar[r]^-{\di{\Delta}}&}
S(V)^+ \wedge \Sigma S(V)^+~;\\[1ex]
\hskip100pt v=[t,u]\mapsto
\begin{cases}
(u,(t,u))&\hbox{\it if}~v \neq 0,\infty\\
\infty&\hbox{\it if}~v=0~{\it or}~\infty~.
\end{cases}
\end{array}$$
with $S$-duality isomorphisms
$$\begin{array}{l}
\sigma_V~:~\{\Sigma S(V)^+\wedge X;Y\} \xymatrix{\ar[r]^-{\di{\cong}}&}
\{V^{\infty} \wedge X;S(V)^+ \wedge Y\}~;~F \mapsto F \sigma_V~,\\[1ex]
\sigma_V~:~\{S(V)^+\wedge X;Y\} \xymatrix{\ar[r]^-{\di{\cong}}&}
\{V^{\infty} \wedge X;\Sigma S(V)^+ \wedge Y\}~;~G \mapsto G\sigma_V
\end{array}$$
for any pointed spaces $X,Y$.
The corresponding reverse $S$-duality map
$\sigma^*_V:S(V)^+ \wedge \Sigma S(V)^+ \to V^{\infty}$
is given by the composite
$$\xymatrix{\sigma^*_V~:~S(V)^+ \wedge \Sigma S(V)^+
\ar[r]^-{\di{\Sigma F}} & \Sigma T(\tau_{S(V)})=V^{\infty} \wedge S(V)^+ \ar[r]
&V^{\infty}}$$
with $F:S(V)^+ \wedge S(V)^+ \to T(\tau_{S(V)})$ the Pontrjagin-Thom map
of the diagonal embedding $\Delta:S(V) \emb S(V) \times S(V)$ with
normal bundle $\nu_{\Delta}=\tau_{S(V)}$ such that
$\nu_{\Delta}\oplus \epsilon_{\R}=\epsilon_V$.\\
{\rm (ii)} The $S$-duality isomorphism
$$\{X;S(V)^+ \wedge Y\} \xymatrix{\ar[r]^-{\di{\cong}}&}
\{\Sigma S(V)^+ \wedge X;V^{\infty} \wedge Y\} ~;~
f \mapsto \Delta \alpha_V \wedge f$$
sends the stable relative difference $\delta'(p,q) \in \{X;S(V)^+ \wedge Y\}$
(\ref{stablereldif3}) of maps $p,q:V^{\infty} \wedge X \to V^{\infty} \wedge Y$
which agree on $0^+ \wedge X\subset V^{\infty} \wedge X$
to the stable homotopy class of the relative  difference
$\delta(p,q) \in \{\Sigma S(V)^+ \wedge X;V^{\infty} \wedge Y\}$
(\ref{difcon1}), with
$$\begin{array}{l}
s_V\delta'(p,q)~=~s_V\delta(p,q)~=~q-p\\[1ex]
\in {\rm ker}(0_V:\{X;Y\} \to \{X;V^{\infty} \wedge Y\})~=~
{\rm im}(s_V:\{X;S(V)^+ \wedge Y\} \to\{X;Y\})~.
\end{array}$$
{\rm (iii)} The inclusion
$S(V) \subset V$ is an embedding with trivial normal $\R$-bundle
$\nu_{S(V) \subset V}=\epsilon_{S(V)}$. By {\rm (i)}
the composite of
$$\alpha_{V}~:~V^{\infty} \to T(\epsilon_{S(V)})~=~\Sigma S(V)^+~;~[t,u] \mapsto (t,u)$$
and the diagonal map $\Delta:\Sigma S(V)^+ \to S(V)^+ \wedge \Sigma S(V)^+$
is an $S$-duality map
$$\sigma_{V}~=~\Delta \alpha_{V}~:~V^{\infty} \to S(V)^+\wedge \Sigma S(V)^+~.$$
Thus for any pointed $CW$ complexes $A,B$ there is
defined an $S$-duality isomorphism
$$\{\Sigma S(V)^+ \wedge A; V^{\infty} \wedge B\}
\to \{A;S(V)^+ \wedge B\}~;~
F \mapsto (1 \wedge F)(\Delta \alpha_{V} \wedge 1)~.$$
\hfill\qed
\end{proposition}

\begin{example} {\rm
Let $M$ be an $m$-dimensional manifold, with tangent bundle
$\tau_M:M \to BO(\R^m)$. Let $\nu_M:M \to BO(k)$
be the normal bundle of an embedding $M \subset \R^{m+k}$, so that
$$\tau_M \oplus \nu_M~=~\tau_{\R^{m+k}}\vert_M~=~\epsilon^{m+k}~,$$
and by the tubular neighbourhood theorem there is defined
a codimension 0 embedding $E(\nu_M) \subset S^{m+k}$.
By Proposition \ref{Sdual} (vi) the composite
$$\begin{array}{l}
\sigma_M~:~S^{m+k} \to S^{m+k}/(S^{m+k}\backslash E(\nu_M))
~=~E(\nu_M)/S(\nu_M)~=~T(\nu_M) \\[1ex]
\hskip25pt\xymatrix{\ar[r]^-{\di{\Delta}}&}
E(\nu_M) \times E(\nu_M)/(E(\nu_M) \times S(\nu_M))~\simeq~
M^+ \wedge T(\nu_M)
\end{array}$$
is an $S$-duality map. Let $z:M \emb E(\nu_M)$ be the zero section. The embedding
$$(1 \times z)\Delta~:~M \emb M \times E(\nu_M)~;~x \mapsto (x,z(x))$$
has normal bundle
$$\nu_{(1 \times z)\Delta}~=~\tau_M \oplus \nu_M~=~\epsilon^{m+k}~:~M
\to BO(\R^{m+k})$$
with an adjunction Umkehr map (\ref{umkehradjunct})
$$M^+\wedge T(\nu_M) \to T(\nu_{(1 \times z)\Delta})~=~\Sigma^{m+k}M^+~.$$
The composite
$$\sigma^*_M~:~M^+\wedge T(\nu_M) \xymatrix{\ar[r] &}
\Sigma^{m+k}M^+ \xymatrix{\ar[r] &} S^{m+k}$$
is a reverse $S$-duality map
such that $\sigma^*_M\sigma_M=1 \in \{M^+;M^+\}$, with $S$-duality isomorphisms
$$\begin{array}{l}
\sigma^*_M~:~
\{X;T(\nu_M)\} \xymatrix{\ar[r]^-{\di{\cong}}&} \{M^+ \wedge X ;S^{m+k}\}~;~
F \mapsto \sigma^*_M  F~,\\[1ex]
\sigma^*_M~:~\{X;M^+\}
\xymatrix{\ar[r]^-{\di{\cong}}&} \{X \wedge T(\nu_M);S^{m+k}\}~;~
G \mapsto \sigma^*_M G~.
\end{array}$$
}
\hfill\qed
\end{example}

\section{The stable cohomotopy Thom and Euler classes}

\begin{proposition}~ \label{thomsum} Let
$$\xymatrix{\xi~:~U \ar[r] & E(\xi) \ar[r]^-{\di{p_\xi}}&X}~,~
\xymatrix{\eta~:~V \ar[r] & E(\eta) \ar[r]^-{\di{p_\eta}}&X}$$
be a $U$- and  a $V$-bundle over a space $X$.\\
{\rm (i)} There is defined a homotopy cofibration sequence
$$\xymatrix{
T(\eta\vert_{S(\xi)})  \ar[r] & T(\eta) \ar[r]^-{\di{z}} &
T(\xi \oplus \eta) \ar[r] & \Sigma T(\eta\vert_{S(\xi)}) \ar[r]& \dots}$$
with
$$z~:~T(\eta) \to T(\xi \oplus \eta)~;~x \mapsto (0,x)~.$$
{\rm (ii)} There is defined an isomorphism of exact sequences
$$\xymatrix@C-10pt{
\dots \ar[r] & \{T(\xi \oplus \eta);(U \oplus V)^{\infty}\}
\ar[r]^-{\di{z^*}} \ar[d]^-{\di{\cong}} &
\{T(\eta);(U \oplus V)^{\infty}\} \ar[r] \ar[d]^-{\di{\cong}} &
\{T(\eta\vert_{S(\xi)});(U \oplus V)^{\infty}\}\ar[r] \ar[d]^-{\di{\cong}}
&\dots \\
\dots \ar[r] & \omega^0(D(\xi),S(\xi);\eta) \ar[r] &
\omega^0(X;\eta) \ar[r]^-{\di{p^*_{S(\xi)}}}
&\omega^0(S(\xi);\eta) \ar[r] &\dots}$$
using the terminology of Definition \ref{virtual} (iii).\\
{\rm (iii)} If $\xi \oplus \eta=\epsilon_{U \oplus V}$ and
${\rm dim}(U \oplus V)=m$
$$\begin{array}{l}
E(\xi \oplus \eta)~=~(U \oplus V)\times X~=~\R^m \times X ~,\\[1ex]
T(\xi \oplus \eta)~=~T(\epsilon_{U \oplus V})~=~
(U \oplus V)^{\infty}\wedge X^+~=~\Sigma^mX^+~,\\[1ex]
\omega^n(D(\xi);-\xi)~=~\widetilde{\omega}^{m+n}(T(\eta))~=~
\{T(\eta);\Sigma^n(U \oplus V)^{\infty}\}~=~\widetilde{\omega}^{m+n}(T(\eta))~,\\[1ex]
\omega^n(D(\xi),S(\xi);-\xi)~=~\omega^{m+n}(T(\eta),T(\eta\vert_{S(\xi)}))\\[1ex]
\hphantom{\omega^n(D(\xi),S(\xi),-\xi)}~=~
\{T(\xi \oplus \eta);\Sigma^n(U \oplus V)^{\infty}\}~=~
\{X^+;S^n\}~=~\omega^n(X)~.
\end{array}$$
\hfill\qed
\end{proposition}

\begin{definition}~ \label{eulerthom}
 {\rm (Crabb \cite[\S\S2,5]{crabb})\\
Let $\xi$ be a $U$-bundle over a space $X$, and let $-\xi=\eta$ be a $V$-bundle over $X$ such that $\xi\oplus \eta=\epsilon_{U \oplus V}$.\\
(i) The {\it stable cohomotopy Thom class}
$$u(\xi) \in \omega^0(D(\xi),S(\xi);-\xi)~\cong~\{X^+;S^0\}$$
is  represented by
$$1~:~X^+ \to S^0~;~x \mapsto -1~,~\infty \mapsto 1~.$$
(ii) The {\it stable cohomotopy Euler class} of $\xi$
$$\gamma(\xi)~=~z^*u(\xi) \in \omega^0(X;-\xi)~\cong~\omega^0(D(\xi);-\xi)~\cong~
\{T(\eta);(U \oplus V)^{\infty}\}$$
is represented by
$$T(\eta) \xymatrix{\ar[r]^-{\di{z}}&} T(\xi \oplus \eta)~=~
(U \oplus V)^{\infty} \wedge X^+ \to (U \oplus V)^{\infty}~.$$
\hfill\qed}
\end{definition}

\begin{example} {\rm For
the trivial $U$-bundle $\xi=\epsilon_U$ over $X$ we can take
$V=\{0\}$, $\eta=0$, so that $\gamma(\epsilon_U)$ is represented by the constant map
$$T(\eta)~=~X^+ \to U^{\infty}~;~x \mapsto 0$$
and
$$\gamma(\epsilon_U)~=~0 \in \omega^0(X;-\epsilon_U)~=~\{X^+;U^{\infty}\}~.$$
}
\hfill\qed
\end{example}

\begin{example} {\rm Let $M$ be an $m$-dimensional manifold, with
tangent bundle $\tau_M:M \to BO(m)$, and let $M \subset \R^{m+k}$ be an embedding
with normal bundle $\nu_M:M \to BO(k)$, so that
$$\tau_M\oplus \nu_M~=~\epsilon_{\R^{m+k}}~:~M \to BO(m+k)~.$$
The stable cohomotopy Euler class of $\tau_M$
$$\gamma(\tau_M) \in \omega^0(M;-\tau_M)~=~\{T(\nu_M);S^{m+k}\}$$
is represented by the composite
$$T(\nu_M) \xymatrix{\ar[r]^-{\di{z}}&} T(\tau_M \oplus \nu_M)~=~\Sigma^{m+k}M^+
\to S^{m+k}~.$$
}
\hfill\qed
\end{example}

\begin{proposition}~ \label{difcon3}
 {\rm (Crabb \cite[\S\S2,5]{crabb})}\\
{\rm (i)} The stable cohomotopy Euler class of the Whitney sum
$V\oplus W$-bundle $\xi\oplus \eta$ of a $V$-bundle $\xi$ over $M$ and a
$W$-bundle $\eta$ over $X$ is the product of the stable cohomotopy Euler classes of $\xi,\eta$
$$\gamma(\xi \oplus \eta)~=~\gamma(\xi) \gamma(\eta)
\in \omega^0(X;-(\xi \oplus \eta))~.$$
{\rm (ii)} For any $V$-bundle $\xi$ over $X$ there
is a bijective correspondence between
the splittings $\xi \cong \xi_1 \oplus \epsilon_{\R}$ and sections
$s:X \to S(\xi)$ of $p_{S(\xi)}:S(\xi) \to X$. If there exists such
a splitting (or section) then
$$\gamma(\xi)~=~\gamma(\xi_1)\gamma(\epsilon_{\R})~=~0 \in \omega^0(X;-\xi)$$
(since $\gamma(\epsilon_{\R})=0$).\\
{\rm (iii)} If $\xi$ is a $V$-bundle over $X$ and
$s:Y \to S(\xi\vert_Y)$ is a section of
$p_{S(\xi\vert_Y)}:S(\xi\vert_Y) \to Y$ (for some subcomplex $Y \subseteq X$ of the $CW$ complex $X$)
then using any extension $\widetilde{s}:(X,Y) \to (D(\xi),S(\xi))$ of $s$
there is defined an Euler class rel $Y$
$$\gamma(\xi,s)~=~\widetilde{s}^*(u) \in \omega^0(X,Y;-\xi)~,$$
with image
$$[\gamma(\xi,s)]~=~\gamma(\xi) \in \omega^0(X;-\xi)~.$$
{\rm (iv)} If $\xi$ is a $V$-bundle over $X$ and
$s_0,s_1:Y \to S(\xi\vert_Y)$ are sections of
$p_{S(\xi\vert_Y)}:S(\xi\vert_Y) \to Y$ which
agree on $Z \subseteq Y$ there is defined a difference class
$$\delta(s_0,s_1) \in \omega^{-1}(Y,Z;-\xi)$$
with image $\gamma(\xi,s_0)-\gamma(\xi,s_1) \in \omega^0(X,Y;-\xi)$.\\
{\rm (v)} Let $X$ be an $n$-dimensional manifold, and let $X \subset
\R^n\oplus U$ be an embedding, with normal $U$-bundle $\nu_X$.
Given a $V$-bundle $\xi$ over $X$ let $g:X \to E(\xi)$ be a generic section
transverse at the zero section $X \subset E(\xi)$, so that the restriction
$$(f,b)~=~g\vert~:~M^m~=~g^{-1}(X)  \to X$$
is the normal map defined by the inclusion of
an $m$-dimensional submanifold $M^m \subset X$ with
$$\begin{array}{l}
m~=~n-{\rm dim}(V)~,~\nu_{M \subset X}~=~\xi\vert_M~,\\[1ex]
b~:~\nu_M~=~\nu_{M \subset X \subset\R^n\oplus U}~=~(\xi \oplus \nu_X)\vert_M
 \to \xi\oplus \nu_X~.
\end{array} $$
For any $W$-bundle $\eta$ over $X$ such that
$\xi\oplus \eta=\epsilon_{V \oplus W}$ the restriction
$\eta\vert_M$ is trivial, and the Euler class
$$\gamma(\xi) \in \omega^0(X;-\xi)~=~\{T(\eta);(V \oplus W)^{\infty}\}$$
is the stable homotopy class of the composite
$$\xymatrix{
T(\eta) \ar[r]^-{\di{F}} &
T(\eta\vert_M)~=~M^+\wedge (V \oplus W)^{\infty}
\ar[r] & (V \oplus W)^{\infty}}$$
with $F$ the adjunction Umkehr map (\ref{umkehradjunct})
of the inclusion
$$E(\eta\vert_M)~=~M \times (V \oplus W) \subset E(\eta)~.$$
In terms of the $S$-dual formulation
$$\begin{array}{l}
\gamma(\xi)~=~T(b)_*(M,1)\\[1ex]
 \in \omega^0(X;-\xi)~=~\omega_{m+{\rm dim}(U)+{\rm dim}(V)}(X;\xi \oplus \nu_X)~=~\Omega_m(X;\xi \oplus \nu_X)
\end{array}$$
with $(M,1) \in \omega_{m+{\rm dim}(U)+{\rm dim}(V)}(M;\nu_M)$.\\
{\rm (vi)} Let $n,U,V,W,X,\nu_X,\xi,\eta$ be as in {\rm (v)}. The fibre product bundles
$$\begin{array}{l}
S(U) \times S(U) \to S(\xi)\times_X S(\xi) \to X~,\\[1ex]
S(U \oplus U) \to S(\xi \oplus \xi)~=~S(\xi)\times_X D(\xi)\cup D(\xi)\times_X S(\xi)
\to X
\end{array}$$
are such that there are defined bundle isomorphisms
$$\begin{array}{l}
\nu_{S(\xi)\times_X S(\xi)\emb S(\xi \oplus \xi)}~\cong~\epsilon_\R~,\\[1ex]
\nu_{\Delta:S(\xi) \emb S(\xi)\times_X S(\xi)} \oplus \epsilon_\R~\cong~
\nu_{\Delta:S(\xi) \emb S(\xi\oplus \xi)}~\cong~\epsilon_U~.
\end{array}$$
Given two sections $s_0,s_1:X \to S(\xi)$ of $\xi$  which agree on a submanifold
$Y \subseteq X$ the section $(s_0,-s_1):X \to S(\xi)\times_X S(\xi)$
is homotopic to a section $(t_0,-t_1):X \to S(\xi)\times_X S(\xi)$
transverse at the diagonal
$\Delta_{S(\xi)} \subset S(\xi)\times_X S(\xi)$. The inverse image
$$C~=~(t_0,-t_1)^{-1}(\Delta_{S(\xi)}) \subseteq X\backslash Y$$
is an $(m+1)$-dimensional submanifold with
$$m~=~n-{\rm dim}(V)~,~(\nu_{C \subset X})\oplus \epsilon_\R \cong \xi\vert_C~.$$
Inclusion defines a normal map $(f,b):C \to X\backslash Y$ with
$$b~:~\nu_C~=~\nu_{C \subset X \subset \R^n \oplus U}\oplus \epsilon_\R~=~
\nu_{C \subset X}\oplus \epsilon_\R\oplus
\nu_X\vert_C \to (\xi \oplus \nu_X)\vert_{X \backslash Y}~.$$
The rel $Y$ difference class
$$\delta(s_0,s_1)\in \omega^{-1}(X,Y;-\xi)~=~
\{\Sigma (T(\eta)/T(\eta\vert_Y));(U \oplus V)^{\infty}\}$$
is the stable homotopy class of the composite
$$\begin{array}{l}
\xymatrix{\Sigma (T(\eta)/T(\eta\vert_Y))
~=~T(\eta\oplus \epsilon_\R)/T(\eta\oplus \epsilon_\R\vert_Y)
\ar[r]^-{\di{F}} &}\\
\hskip75pt
\xymatrix{T(\xi \oplus \eta\vert_C)~=~C^+\wedge (U \oplus V)^{\infty}
\ar[r]&(U \oplus V)^{\infty}}
\end{array}$$
with $F$ the adjunction Umkehr map of the inclusion
$$E((\eta\oplus \epsilon_\R)\vert_C)~=~
C \times (U \oplus V) \subset E(\eta\oplus \epsilon_\R)~.$$
In terms of the $S$-dual formulation
$$\begin{array}{l}
\delta(s_0,s_1)~=~T(b)_*(C,1)\in \omega^{-1}(X,Y;-\xi)\\[1ex]
\hskip5pt =~
\omega_{m+{\rm dim}(U)+{\rm dim}(V)+1}(X\backslash Y;(\xi \oplus \nu_X)\vert_{X \backslash Y})~=~
\Omega_{m+1}(X\backslash Y;(\xi \oplus \nu_X)\vert_{X \backslash Y})
\end{array}$$
with $(C,1) \in \omega_{m+{\rm dim}(U)+{\rm dim}(V)+1}(C;\nu_C)$.\\
\hfill\qed
\end{proposition}

\begin{example} {\rm
(i) A non-zero section of the trivial $V$-bundle $\epsilon_V$
over a space $X$ is essentially the same as a map of the type
$$X \to S(\epsilon_V)~=~X \times S(V)~;~x \mapsto (x,s(x))~,$$
as determined by a map $s:X \to S(V)$. Any such map $s$ determines a map
$$p_s~:~CX^+ \to V^{\infty}~;~(t,x) \mapsto
\frac{t s(x)}{1-t}$$
sending $X=\{0\}\times X \subset CX^+$ to $0 \in V^{\infty}$.
The rel $X\cup CY$ difference class of the non-zero sections (\ref{difcon3} (iv))
given by two such maps $s_0,s_1:X \to S(V)$ which agree on $Y \subseteq X$
is just the rel $X\cup CY$ difference
(\ref{main} (iii)) of $p_{s_0},p_{s_1}:CX^+ \to V^{\infty}$
$$\delta(s_0,s_1)~=~\delta(p_{s_0},p_{s_1}) \in
\omega^{-1}(X,Y;-\epsilon_V)~=~\{\Sigma (X/Y);V^{\infty}\}~.$$
(ii) If $X$ is an $n$-dimensional manifold and
$s_0,s_1:X \to S(V)$ are maps such that
$$d~:~X \to S(V) \times S(V)~;~x \mapsto (s_0(x),-s_1(x))$$
is transverse regular at $\Delta_{S(V)} \subset
S(V) \times S(V)$ then
$$C~=~d^{-1}(\Delta_{S(V)})~=~\{x \in X\,\vert\, s_0(x)=-s_1(x)
\in S(V)\} \subset X$$
is an $(m+1)$-dimensional submanifold with
$$m~=~n-{\rm dim}(V)~,~
\nu_{C \subset X}~=~(d\vert_C)^*\nu_{\Delta:S(V) \emb S(V) \times S(V)}~=~
(d\vert_C)^*\tau_{S(V)}$$
such that $\nu_{C \subset X}\oplus \epsilon_{\R}
\cong \epsilon_V$, and
$$\xymatrix{
\delta(s_0,s_1)~:~\Sigma X^+ \ar[r]^-{\di{\Sigma F}} &
\Sigma T(\nu_{C \subset X})~=~
C^+ \wedge V^{\infty} \ar[r] & V^{\infty}}$$
with $F:X^+ \to T(\nu_{C \subset X})$ the adjunction Umkehr map of
$E(\nu_{C \subset X}) \subset X$.\\
(iii) For any maps $s_0,s_1:X=S(V) \to S(V)$ with
$$C~=~\{x \in X\,\vert\, s_0(x)=-s_1(x) \in S(V)\}$$
a 0-dimensional submanifold (= finite subset) of $X$,
the difference class of the corresponding sections
$$s_i~:~X \to S(\epsilon_V)~=~X \times S(V)~;~x \mapsto (x,s_i(x))~(i=0,1)$$
of $p_{S(\epsilon_V)}$ is
$$\delta(s_0,s_1)~=~\vert C \vert\in \{\Sigma X^+;V^{\infty}\}~=~\ZZ~,$$
counting the points of $C$ algebraically. In particular, if
$$s_0(x)~=~x~,~s_1(x)~=~x_0 \in S(V)~(x \in S(V))$$
for some point $x_0 \in S(V)$ then $C=\{-x_0\} \subset X$ and
$$\delta(s_0,s_1)~=~\vert C \vert~=~1\in \{\Sigma X^+;V^{\infty}\}~=~\ZZ~.$$
}\hfill\qed
\end{example}

\begin{definition} \label{adjoint} {\rm
An element $c\in O(V,U \oplus V)$ is a linear isometry $c:V \to U \oplus V$.
The {\it adjoint} of a map $c:X \to O(V,U\oplus V)$ is the pointed map
$$F_c~:~V^{\infty}\wedge X^+  \to (U\oplus V)^{\infty}~;~(v,x) \mapsto c(x)(v)~.$$
\hfill\qed}
\end{definition}

\begin{proposition}~\label{euler1}
{\rm (i)} The stable cohomotopy Euler class of the canonical
$V$-bundle $\xi(V):G=G(V,U \oplus V) \to BO(V)$
$$\gamma(\xi(V)) \in \omega^0(G;-\xi(V))~=~\{T(\xi(V));(U \oplus V)^{\infty}\}$$
is the stable homotopy class of the adjoint map
$$\gamma(\xi(V))~:~T(\xi(U))~=~O(U,U \oplus V) \wedge_{O(U)}U^{\infty} \to
(U \oplus V)^{\infty}~;~(i,x) \mapsto i(x)~.$$
{\rm (ii)} Any $V$-bundle $\xi:X \to BO(V)$ over a finite $CW$ complex $X$
is isomorphic to the pullback $f^*\xi(V)$ of $\xi(V)$ along a map
$f:X \to G=G(V,U \oplus V)$, for sufficiently large ${\rm dim}(U)$.
The pullback $\eta=f^*\xi(U):X \to BO(U)$ is a $U$-bundle over $X$
such that $\xi \oplus \eta \cong \epsilon_{U \oplus V}$, and
the stable cohomotopy Euler class
$$\gamma(\xi)~=~f^*\gamma(\xi(V)) \in \omega^0(X;-\xi)~=~\{T(\eta);(U \oplus V)^{\infty}\}$$
is the stable homotopy class of the composite
$$\xymatrix{\gamma(\xi)~:~T(\eta) \ar[r]^-{\di{f}} &T(\xi(U))~=~
O(U,U \oplus V)^+\wedge_{O(U)}U^{\infty} \ar[r]^-{\di{\gamma(\xi(U))}}
& (U \oplus V)^{\infty}~.}$$
A $U \oplus V$-bundle
isomorphism $\delta\xi:\xi \oplus \epsilon_U \cong \epsilon_{U \oplus V}$
corresponds to a lift of $\xi:X \to BO(V)$ to a map $\delta\xi:X \to O(U,U \oplus V)$,
in which case $\eta\cong \epsilon_U$ and
$\gamma(\xi) \in \omega^0(X;-\xi)~=~\{X^+;V^{\infty}\}$
is the stable homotopy class of the adjoint of $\delta \xi$
$$\gamma(\xi)~:~U^{\infty}\wedge X^+  \to (U \oplus V)^{\infty}~;~
(u,x) \mapsto \delta \xi(x)(u)~.$$
{\rm (iii)} A section $s:X \to S(\xi)$ of $\xi:X \to BO(V)$ determines
a null-homotopy
$$\gamma(s)~:~\gamma(\xi)~\simeq~*~:~T(\eta) \to (U \oplus V)^{\infty}~,$$
as given by an extension of $\gamma(\xi)$ to a map
$\gamma(s):CT(\eta) \to (U \oplus V)^{\infty}$. The rel $Y$
difference class (\ref{difcon3} (iv)) of sections $s_0,s_1:X \to S(\xi)$
which agree on $Y \subseteq X$ is just the
rel $T(\eta)\cup CT(\eta\vert_Y)$ difference (\ref{main} (iii))
$$\delta(s_0,s_1)~=~\delta(\gamma(s_0),\gamma(s_1))\in
\omega^{-1}(X,Y;-\xi)~=~\{\Sigma (T(\eta)/T(\eta\vert_Y));
(U \oplus V)^{\infty}\}~.$$
\end{proposition}
\begin{proof} (i) The composite
$$\begin{array}{l}
\xymatrix@C+20pt{T(\xi(U))~=~O(U,U \oplus V)^+\wedge_{O(U)}U^{\infty}
\ar[r]^-{\di{z=1 \wedge j}} &}\\[1ex]
\hskip50pt  \xymatrix@C-10pt{T(\xi(U) \oplus \xi(V))~=~
O(U,U \oplus V)^+\wedge_{O(U)}(U \oplus V)^{\infty}&}\\[1ex]
\hskip100pt
\xymatrix{\cong G^+ \wedge (U \oplus V)^{\infty} \ar[r] & (U \oplus V)^{\infty}}
\end{array}$$
is given by $(h,x) \mapsto h(x)$.\\
(ii) Since $X$ is a finite $CW$ complex
there exists a sufficiently high-dimensional $U$
such that the classifying map $\xi:X \to BO(V)$ factors up to homotopy as
$$\xymatrix{\xi~:~X \ar[r]^-{\di{f}} & G=G(V,U \oplus V)\ar[r] & BO(V)}$$
with $\xi \cong f^*\xi(V)$, $\xi \oplus f^*\xi(U) \cong \epsilon_{U \oplus V}$.\\
(iii) The stable map representing $\gamma(\xi)$ is the composite
$$\xymatrix{\gamma(\xi)~:~T(\eta) \ar[r]^-{\di{z}}
& T(\xi \oplus \eta)~\cong~(U \oplus V)^{\infty}\wedge X^+  \ar[r] &
(U \oplus V)^{\infty}~.}$$
The Thom space of $\xi \oplus \eta$ fits into a homotopy cofibration sequence
$$\xymatrix{T(p^*\eta) \ar[r]^-{\di{T(p)}} & T(\eta) \ar[r]^-{\di{z}}
& T(\xi \oplus \eta) \ar[r] & \dots}$$
with $p=p_{S(\xi)}:S(\xi) \to X$, and $T(p):T(p^*\eta) \to T(\eta)$ the
map of Thom spaces induced by the map of $V$-bundles $p:p^*\eta \to \eta$
over $p:S(\xi) \to X$. For any section $s:X \to S(\xi)$ of $p:S(\xi) \to X$
it follows from $p\circ s=1:X \to X$ that $s^*(p^*\eta)=\eta$, so that
$s$ lifts to a map of $U$-bundles $s:\eta \to p^*\eta$ inducing
a section $T(s):T(\eta) \to T(p^*\eta)$ of $T(p)$, giving the null-homotopies
$z \simeq *$, $\gamma(\xi) \simeq *$.\\
\hfill\qed\end{proof}

\chapter{$\ZZ_2$-equivariant homotopy and bordism theory}\label{z2homotopy}

Intersections and self-intersections of maps of manifolds are
expressed in terms of algebraic topology by means of the $\ZZ_2$-equivariant
homotopy properties of the diagonal maps
$$\begin{array}{l}
\Delta~:~X \to X \times X~;~x \mapsto (x,x)~(\hbox{$X$ unpointed})~,\\[1ex]
\Delta~:~X \to X \wedge X~;~x \mapsto (x,x)~(\hbox{$X$ pointed})~.
\end{array}$$

\section{$\pi$-equivariant homotopy theory}

Let $\pi$ be a group. We shall be concerned with pointed spaces $X$ with a
$\pi$-action
$$\pi \times X \to X~;~(g,x) \mapsto gx$$
fixing the base point, particularly for $\pi=\ZZ_2$.

\begin{definition}~
{\rm (i) Given a (pointed) $\pi$-space $X$ let $\vert X \vert$ be the (pointed) space
defined by $X$ with the $\pi$-action forgotten.\\
(ii) Given pointed $\pi$-spaces $X,Y$ let $[X,Y]_{\pi}$ be the set of
(pointed) $\pi$-equivariant homotopy classes of $\pi$-equivariant maps.
For any pointed $\pi$-spaces $X,Y$ there is defined a forgetful
function
$$[X,Y]_{\pi} \to [\vert X\vert,\vert Y \vert]~.$$
(iii) The {\it fixed point set} of a $\pi$-space $X$ is\index{fixed point set}
$$X^{\pi}~=~\{x \in X\,\vert\, gx=x \in X~\hbox{for all}~g \in \pi\}~.$$
A $\pi$-equivariant map $f:X \to Y$ restricts to a map of
the fixed point sets
$$\rho(f)~:~X^{\pi} \to Y^{\pi}~.$$
Similarly for $\pi$-equivariant homotopies, with a fixed point function
$$\rho~:~[X,Y]_{\pi} \to [X^{\pi},Y^{\pi}]$$
for any pointed $\pi$-spaces $X,Y$.\\
(iv) A $\pi$-space $X$ has the  {\it trivial} $\pi$-action if\index{trivial $\pi$-action}
$$gx~=~x~~(x \in X,~g\in \pi)~,$$
or equivalently $X^{\pi}=X$. Given a (pointed) space $X$ let
$X$ also denote the (pointed) $\pi$-space defined by the trivial
$\pi$-action on $X$.\\
(v) A $\pi$-space $X$ is {\it free} if\index{free $\pi$-space} for each $x \in X$ the only $g \in \pi$ with $gx=x$ is $g=1 \in \pi$. For $\pi=\ZZ_2$ this is equivalent to
$$X^{\ZZ_2}~=~\emptyset~,$$
or equivalently if the quotient map $X \to X/\pi$ is a regular covering projection.\\
(vi) A pointed $\pi$-space $X$ is {\it semifree} if\index{semifree pointed $\pi$-space} if $X \backslash \{*\}$ is a free $\pi$-space.}\\
\hfill\qed
\end{definition}

\begin{example}~ {\rm
(i) If $X$ is a free $\pi$-space then $X^+$ is a semifree $\pi$-space.\\
(ii) For any $\pi$-spaces $X,Y$
$$(X \times Y)^{\pi}~=~X^{\pi} \times Y^{\pi}$$
so that if $X$ is free then so is $X \times Y$.\\
(iii) For any pointed $\ZZ_2$-spaces $X,Y$
$$(X \wedge Y)^{\pi}~=~X^{\pi} \wedge Y^{\pi}$$
so that if $X$ is semifree then so is $X \wedge Y$.\\
\hfill\qed}
\end{example}

\begin{definition}~
{\rm (i) For any pointed spaces $A,B$ let $\map_*(A,B)$ be the
space of pointed maps $A \to B$, with path component set
$$\pi_0(\map_*(A,B))~=~[A,B]~.$$
(ii) For any pointed $\pi$-spaces $A,B$ let ${\rm map}^{\pi}_*(A,B)$ be the
space of $\pi$-equivariant pointed maps $A \to B$, with path component set
$$\pi_0({\rm map}_*^{\pi}(A,B))~=~[A,B]_{\pi}~.$$}
\hfill\qed
\end{definition}

\begin{proposition}~{\rm (i)} For any space $X$ and $\pi$-space $Y$
a $\pi$-equivariant map $f:X \to Y$ is the composite
$$\xymatrix{f~:~X \ar[r] &} Y^{\pi} \xymatrix{\ar[r] &Y}$$
of a map $X \to Y^{\pi}$ and the inclusion $Y^{\pi} \to Y$.\\
{\rm (ii)} For a pointed space $X$ and pointed $\pi$-space $Y$
$${\rm map}_*^{\pi}(X,Y)~=~{\rm map}_*(X,Y^{\pi})~,~[X,Y]_{\pi}~=~[X,Y^{\pi}]~.$$
\end{proposition}
\begin{proof} By construction.\\
\hfill\qed\end{proof}

We now specialize to the case $\pi=\ZZ_2$. In Chapter \ref{pi-equivariant} we shall deal with arbitrary $\pi$.

\begin{example}~ {\rm
{\rm (i)} A (pointed) $\ZZ_2$-space $X$ is a (pointed) space with an involution
$$T~:~X \to X~.$$
{\rm (ii)} For any pointed space $Y$ regard $Y \wedge Y$ as a pointed $\ZZ_2$-space
by the transposition
$$T~:~ Y \wedge Y \to Y \wedge Y~;~(y_1,y_2) \mapsto (y_2,y_1)~.$$
The inclusion of the fixed point set is just the diagonal map
$$\Delta~:~(Y \wedge Y)^{\ZZ_2}~=~Y \to Y \wedge Y~;~y \mapsto (y,y)~,$$
so that for any pointed space $X$
$$[X,Y \wedge Y]_{\ZZ_2}~=~[X,Y]~.$$
\hfill\qed}
\end{example}

\begin{definition}~{\rm
(i) An {\it inner product $\ZZ_2$-space} $V$ is an inner product space with
a $\ZZ_2$-action $T:V \to V$ which is an isometry. We shall write
$\vert V \vert$ for the underlying inner product space. Let
$S(V)$ and $P(V)$ denote the $\ZZ_2$-spaces defined by the
unit sphere $S(\vert V \vert)$ and the projective space $P(\vert V \vert)$
with the $\ZZ_2$-action induced by $T$.\index{inner product $\ZZ_2$-space}\\
(ii) An inner product space $V$ can be regarded as an inner product $\ZZ_2$-space
with the trivial $\ZZ_2$-action $T=1:V \to V$, so that $\vert V \vert = V$.
The $\ZZ_2$-actions on $S(V)$ and $P(V)$ are trivial.\\
(iii) Let $L$ be the $\R[\ZZ_2]$-module $\R$ with the involution $-1$.
For any inner product space $V$ there is then defined an inner product
$\ZZ_2$-space
$$LV~=~L\otimes_{\R} V~,~T_{LV}~:~LV \to LV~;~v \mapsto -v~,$$
with $\vert LV \vert = \vert V \vert$. The $\ZZ_2$-action on $S(LV)$
is non-trivial, and the $\ZZ_2$-action on $P(LV)$ is trivial. For $V=\R^n$
write
$$LV^{\infty}~=~LS^n~.$$
\hfill\qed}
\end{definition}

\begin{proposition}~Let $V$ be an inner product $\ZZ_2$-space.\\
{\rm (i)} The inner product spaces
$$\begin{array}{l}
V_+~=~V^{\ZZ_2}~=~\{x \in V\,\vert\, Tx=x\}~,\\[1ex]
V_-~=~(LV)^{\ZZ_2}~=\{x \in V\,\vert\, Tx=-x\} \subseteq V
\end{array}$$
are such that as an inner product $\ZZ_2$-space $V$ has a decomposition
as a sum of inner product $\ZZ_2$-spaces
$$V~=~V_+ \oplus LV_-~,~T(x_+,x_-)~=~(x_+,-x_-)$$
with $LV_-$ short for $L(V_-)$.\\
{\rm (ii)} The one-point compactification $V^{\infty}$
of $V=V_+\oplus LV_-$ is a pointed $\ZZ_2$-space
$$V^{\infty}~=~V_+^{\infty} \wedge LV_-^{\infty}$$
with
$$\begin{array}{l}
T~:~V^{\infty} \to V^{\infty}~;~(x_+,x_-) \mapsto (x_+,-x_-)~,\\[1ex]
(V^{\infty})^{\ZZ_2}~=~(V_+^{\infty})^{\ZZ_2} \wedge (LV_-^{\infty})^{\ZZ_2}~=~V_+^{\infty} \wedge \{0\}^{\infty}~=~V_+^{\infty}~,\\[1ex]
V^{\infty}/\ZZ_2~=~V_+^{\infty}/\ZZ_2\wedge LV_-^{\infty}/\ZZ_2~=~
V_+^{\infty}\wedge s P(V_-)~.
\end{array}$$
{\rm (iii)} The unit sphere $S(V)$ of $V=V_+\oplus LV_-$
is a $\ZZ_2$-space such that the homeomorphism of Proposition \ref{pushout2}
$$\begin{array}{l}
\lambda_{V_+,LV_-}~:~S(V_+)*S(LV_-) \to S(V_+ \oplus LV_-)~;\\[1ex]
\hskip50pt (t,v_+,v_-) \mapsto \big(v_+\cos(\pi t/2),v_-\sin(\pi t/2)\big)
\end{array}$$
is $\ZZ_2$-equivariant, with
$$\begin{array}{l}
T~:~S(V_+)*S(LV_-) \to S(V_+)*S(LV_-)~;~(t,v_+,v_-) \mapsto (t,v_+,-v_-)~,\\[1ex]
S(V)^{\ZZ_2}~=~S(V_+)*\emptyset ~=~S(V_+)~,~S(V)/\ZZ_2~=~S(V_+) * P(V_-)~.
\end{array}$$
The induced $\ZZ_2$-action on the projective space $P(V)$
$$T~:~P(V) \to P(V)~;~[v_+,v_-] \mapsto [v_+,-v_-]$$
is such that
$$P(V)^{\ZZ_2}~=~P(V_+) \sqcup P(V_-)~.$$
{\rm (iv)} A subspace $W \subseteq V$ is $\ZZ_2$-invariant
(i.e. $TW=W$) if and only if
$W=W_+ \oplus W_-$ for subspaces $W_+ \subseteq V_+$, $W_- \subseteq V_-$.\\
\end{proposition}
\begin{proof}
(i) It is clear that $V_+ \cap V_-=\{0\}$, and every $x \in V$ can be written as
$$x~=~\frac{x+Tx}{2}+ \frac{x-Tx}{2} \in V_+ + V_-~.$$
(ii)+(iii)+(iv) By construction.\\
\hfill\qed\end{proof}

\begin{proposition}\label{forget}~
Let $X,Y$ be pointed $\ZZ_2$-spaces.\\
{\rm (i)} The pointed $\ZZ_2$-space
$$S(L\R)^+ \wedge X~=~X \vee X~,~T(\pm x)~=~\mp Tx~(x \in X)$$
is the one-point union of two copies of $X$ which are transposed by the
$\ZZ_2$-action. The map
$$S(L\R)^+ \wedge X \to S(L\R)^+ \wedge \vert X \vert~;~
\begin{cases}
(+,x) \mapsto (+,x)\\[1ex]
(-,x) \mapsto (-,Tx)
\end{cases}$$
is a $\ZZ_2$-equivariant homeomorphism.\\
{\rm (ii)} A $\ZZ_2$-equivariant map $S(L\R)^+ \wedge X \to Y$ is
essentially the same as a map $\vert X \vert \to \vert Y\vert$, so that
$$[S(L\R)^+\wedge X,Y]_{\ZZ_2}~=~[\vert X\vert,\vert Y\vert]~.$$
{\rm (iii)} The $\ZZ_2$-equivariant homotopy cofibration sequence
$$S(L\R)^+ \to S^0 \to L\R^{\infty} \to \Sigma S(L\R)^+ \to \dots $$
induces the exact sequence
$$\begin{array}{l}\dots \to
[\Sigma S(L\R)^+ \wedge X,Y]_{\ZZ_2}=[\Sigma X,Y]\to
[L\R^{\infty} \wedge X,Y]_{\ZZ_2}\\[1ex]
\hskip75pt \to [X,Y]_{\ZZ_2} \to [S(L\R)^+ \wedge X,Y]_{\ZZ_2}=[\vert X \vert,\vert Y \vert]
\end{array}$$
with $[X,Y]_{\ZZ_2} \to [\vert X \vert,\vert Y \vert]$ the function which forgets $\ZZ_2$-equivariance,
and
$$[L\R^{\infty} \wedge X,Y]_{\ZZ_2}\to [X,Y]_{\ZZ_2}~;~F \mapsto (x \mapsto F(0,x))~.$$
\hfill\qed
\end{proposition}

\begin{example} {\rm The zero map $0_{L\R}:S^0 \to L\R^{\infty}$
induces a bijection
$$0_{L\R}^*~:~[L\R^{\infty},L\R^{\infty}]_{\ZZ_2} \to [S^0,L\R^{\infty}]_{\ZZ_2}$$
which sends $1:L\R^{\infty} \to L\R^{\infty}$ to $0_{L\R}$.\\
\hfill\qed}
\end{example}

\begin{proposition}~
{\rm (i)} For any inner product $\ZZ_2$-spaces $U,V$
$$P(U \oplus V)/P(V)~=~S(U)^+ \wedge_{\ZZ_2} V^{\infty}~.$$
{\rm (ii)} The continuous bijection
$V^{\infty} \to sS(V)$ of Proposition \ref{sphere}
(which is a homeomorphism for finite-dimensional $V$) is
$\ZZ_2$-equivariant, allowing the identifications
$$V^{\infty}~=~sS(V)~,~(V^{\infty})^{\ZZ_2}~=~sS(V_+)~,~
V^{\infty}/\ZZ_2~=~S(V_+)*sP(V_-)~.$$
{\rm (iii)} For an inner product space $V$ with the trivial $\ZZ_2$-action
$$S(LV)^{\ZZ_2}~=~\emptyset~,~S(LV)/\ZZ_2~=~P(V)~.$$
For any $\ZZ_2$-space $X$ the $\ZZ_2$-space $S(LV) \times X$ is free; for any
pointed $\ZZ_2$-space $X$ the pointed $\ZZ_2$-space
$S(LV)^+ \wedge X$ is semifree.\\
\end{proposition}
\begin{proof} (i) Immediate from the $\ZZ_2$-equivariant homotopy cofibration sequence
$$S(V) \to S(U \oplus V) \to S(U)^+ \wedge V^{\infty}~.$$
(ii)+(iii) By construction.\\
\hfill\qed\end{proof}

\begin{proposition}~\label{sequence}
For any pointed $\ZZ_2$-spaces $X,Y$ and inner product space $V$
there is defined a long exact sequence of abelian groups/pointed sets
$$\begin{array}{l}
\dots \xymatrix{\ar[r]&} [\Sigma X,Y]_{\ZZ_2}
\xymatrix{\ar[r]^-{\di{s_{LV}^*}} &}
[\Sigma S(LV)^+  \wedge X,Y]_{\ZZ_2}
\xymatrix{\ar[r]^-{\di{\alpha^*_{LV}}}&}\\
\hskip175pt
[LV^{\infty} \wedge X,Y]_{\ZZ_2} \xymatrix{\ar[r]^-{\di{0_{LV}^*}}&}
[X,Y]_{\ZZ_2}~.
\end{array}$$
\end{proposition}
\begin{proof} Immediate from the $\ZZ_2$-equivariant homotopy cofibration sequence
$$\xymatrix{
S(LV)^+ \ar[r]^-{\di{s_{LV}}}& S^0 \ar[r]^-{\di{0_{LV}}}&
LV^{\infty} \ar[r]^-{\di{\alpha_{LV}}}& \Sigma S(LV)^+  \ar[r]& \dots~.}$$
of Proposition \ref{cofib2} (iii).\\
\hfill\qed\end{proof}

\begin{terminology}~{\rm The infinite-dimensional inner product space
$$\R(\infty)~=~\varinjlim\limits_k \R^k$$
is denoted by $\R(\infty)$ to avoid confusion with the one-point
compactification $\R^{\infty}$ of $\R$. The unit sphere
$$S(\infty)~=~S(L\R(\infty))~=~\mathop{\varinjlim}\limits_k S(L\R^k)$$
is a contractible space with a free $\ZZ_2$-action, with quotient
the infinite-dimensional real projective space
$$P(\infty)~=~S(\infty)/\ZZ_2~.$$
The unreduced suspension
$$sS(\infty)~=~\varinjlim\limits_k sS(L\R^k)~=~
\varinjlim\limits_k (L\R^k)^{\infty}$$
is a contractible space with a non-free $\ZZ_2$-action, such that there
are two fixed points
$$sS(\infty)^{\ZZ_2}~=~S^0~.$$
In dealing with inner product spaces $V$ we interpret $LV^{\infty}$
in the infinite-dimensional case $V=\R(\infty)$ to be
$$sS(\infty)~=~\varinjlim\limits_{U \subset \R(\infty)~{\rm finite-dimensional}}LU^{\infty}$$
(which is not compact, and in particular not
the actual one-point compactification $L\R(\infty)^{\infty}$ of
$L\R(\infty)=\varinjlim\limits_U LU$). \\
\hfill\qed}
\end{terminology}

\begin{proposition}~ \label{fixed}
For any pointed $\ZZ_2$-spaces $X,Y$ and inner product space $V$ let
$$\rho_V~:~\map_*^{\ZZ_2}(X,LV^{\infty} \wedge Y) \to
\map_*(X^{\ZZ_2},Y^{\ZZ_2})$$
be the fixed point map, and let
$$\begin{array}{l}
\rho_\infty~=~\varinjlim\limits_V \rho_V~:~
\varinjlim\limits_V \map_*^{\ZZ_2}(X,LV^{\infty} \wedge Y)~=~\map_*^{\ZZ_2}(X,L\R(\infty)^{\infty} \wedge Y)\\[1ex]
\hskip150pt  \to  \map_*(X^{\ZZ_2},Y^{\ZZ_2})
\end{array}$$
be the map obtained by passing to the limit over finite-dimensional $V$.\\
{\rm (i)} If $X$ is a $CW$ $\ZZ_2$-complex then $\rho_{\infty}$
is a fibration with contractible point inverses, and so induces
isomorphisms in the homotopy groups. In particular, there is induced a
bijection
$$\begin{array}{l}
\rho_{\infty}~:~\pi_0(\map_*^{\ZZ_2}(X,L\R(\infty)^{\infty} \wedge Y))~=~
[X,L\R(\infty)^{\infty} \wedge Y]_{\ZZ_2}\\[1ex]
\hskip100pt
\xymatrix{\ar[r]^-{\di{\cong}}&} \pi_0(\map_*(X^{\ZZ_2},Y^{\ZZ_2}))~=~[X^{\ZZ_2},Y^{\ZZ_2}]~.$$
\end{array}$$
{\rm (ii)} If the $\ZZ_2$-action on $X$ is trivial then
$\rho_V$ is a homeomorphism, with the inclusion
$$\sigma_V~:~Y^{\ZZ_2}~=~(LV^{\infty} \wedge Y)^{\ZZ_2} \to LV^{\infty} \wedge Y~;~
y \mapsto (0,y)$$
such that
$$\rho_V^{-1}~=~\sigma_V~:~\map_*(X,Y^{\ZZ_2}) \to \map_*^{\ZZ_2}(X,LV^{\infty} \wedge Y)~,$$
and
$$[X,LV^{\infty}\wedge Y]_{\ZZ_2}~=~[X,Y^{\ZZ_2}]~.$$
The $\ZZ_2$-equivariant homotopy group
$[LV^{\infty} \wedge  X,LV^{\infty} \wedge Y]_{\ZZ_2}$ fits into a
direct sum system
$$[\Sigma S(LV)^+ \wedge  X,LV^{\infty}\wedge Y]_{\ZZ_2}
\xymatrix{\ar@<1ex>[r]^-{\di{\gamma}}\ar@<-1ex>@{<-}[r]_-{\di{\delta}} & }
[LV^{\infty} \wedge  X,LV^{\infty} \wedge Y]_{\ZZ_2}
\xymatrix{\ar@<1ex>[r]^-{\di{\rho}}\ar@<-1ex>@{<-}[r]_-{\di{\sigma}}& }
[ X, Y^{\ZZ_2}]$$
with
$$\delta \gamma~=~1~,~\rho\sigma~=~1~,~\gamma\delta+\sigma\rho~=~1~.$$
Here $\gamma=\alpha^*_{LV}$, and $\sigma$ is defined by
$$\sigma~:~[ X, Y^{\ZZ_2}] \to [LV^{\infty} \wedge  X,LV^{\infty}\wedge Y]_{\ZZ_2}~;~
G \mapsto 1_{LV^{\infty}}\wedge \sigma_VG~.$$
A $\ZZ_2$-equivariant map
$F:LV^{\infty} \wedge  X \to LV^{\infty} \wedge Y$
with fixed point map $G=\rho(F): X \to  Y^{\ZZ_2}$
is such that $F$ and $\sigma(G)$ agree on $0^+ \wedge  X$,
and $\delta(F)$ is defined to be
the relative difference $\ZZ_2$-equivariant map
$$\delta(F)~=~\delta(F,\sigma(G))~:~\Sigma S(LV)^+ \wedge  X \to
LV^{\infty} \wedge  Y~,$$
such that
$$F-\sigma(G)~=~\alpha^*_{LV}\delta(F,\sigma(G))
\in {\rm im}(\alpha^*_{LV})~=~{\rm ker}(\rho)~.$$
\end{proposition}
\begin{proof} (i) We repeat the argument of Sinha \cite[p.277]{sinha}
(a space level version of Crabb \cite[Lemma (A.1), p.60]{crabb}).
The fibre over a component of $\map_*(X^{\ZZ_2},Y^{\ZZ_2})$ is the
space of $\ZZ_2$-maps $f:X\to L\R(\infty)^{\infty}\wedge Y$ which are
specified on $X^{\ZZ_2}$ by $\rho(f):X^{\ZZ_2} \to Y^{\ZZ_2}$.
We consider the effect on this mapping space of attaching
$\ZZ_2$-cells to $X$. There are two types: single
cells with trivial $\ZZ_2$-action $D^n$ and pairs of cells with
free $\ZZ_2$-action $\ZZ_2 \times D^n$. For the single cells
the extension of $f$ is itself specified by the extension of
$\rho(f)$. For the pairs of cells note that for any space $W$
the $\ZZ_2$-maps $g:\ZZ_2 \times W \to L\R(\infty)^{\infty}\wedge Y$
are just the maps $g:W \to L\R(\infty)^{\infty}\wedge Y$, and
that $L\R(\infty)^{\infty}\wedge Y$ is contractible, so that
the space of $\ZZ_2$-maps $g:\ZZ_2 \times D^n \to L\R(\infty)^{\infty}\wedge Y$
extending a given $\ZZ_2$-map $\partial g:\ZZ_2 \times S^{n-1} \to L\R(\infty)^{\infty}\wedge Y$
is contractible.\\
(ii) By Proposition \ref{sequence} there is defined an exact sequence
$$\begin{array}{l}
\dots \to
[\Sigma S(LV)^+ \wedge  X,LV^{\infty} \wedge Y ]_{\ZZ_2}
\xymatrix{\ar[r]^-{\di{\alpha^*_{LV}}}&}
[LV^{\infty} \wedge  X,LV^{\infty}\wedge Y ]_{\ZZ_2}\\[1ex]
\hskip120pt \xymatrix{\ar[r]^-{\di{0^*_{LV}}}&}
[ X,LV^{\infty} \wedge Y]_{\ZZ_2}~=~[ X,Y^{\ZZ_2}]
\end{array}$$
\hfill\qed\end{proof}

\section{The bi-degree}

The degree is a $\ZZ$-valued homotopy invariant
of a pointed map $V^{\infty}\to V^{\infty}$ for a
 finite-dimensional inner product space $V$, which defines an
isomorphism
$$[V^{\infty},V^{\infty}] ~\cong~ \ZZ~.$$
Similarly, the bi-degree is a $\ZZ\oplus \ZZ$-valued $\ZZ_2$-homotopy invariant
of a pointed $\ZZ_2$-equivariant map $V^{\infty}\to V^{\infty}$
for a finite-dimensional inner product $\ZZ_2$-space $V$, which defines an
isomorphism
$$[V^{\infty},V^{\infty}]_{\ZZ_2} ~\cong~ \ZZ \oplus \ZZ~.$$

We start by recollecting the main properties of the degree.

\begin{definition}~{\rm Let $V$ be a finite-dimensional inner product
space. The {\it degree} of a pointed map $F:V^{\infty} \to V^{\infty}$ is the homotopy
invariant defined by\index{degree}
$${\rm degree}(F)~=~F_*(1) \in \ZZ~,$$
with
$$F_*~:~\dot H_{{\rm dim}(V)}(V^{\infty}) \to
\dot H_{{\rm dim}(V)}(V^{\infty})$$
and $1 \in \dot H_{{\rm dim}(V)}(V^{\infty}) \cong \ZZ$ either generator.}
\hfill\qed
\end{definition}

\begin{example}~{\rm
(i) The identity $1:V^\infty \to V^\infty;v \mapsto v$ has ${\rm degree}(1)=1$.\\
(ii) The constant map $\infty:V^\infty \to V^\infty;v \mapsto \infty$
has ${\rm degree}(\infty)=0$.\\
\hfill\qed}
\end{example}

\begin{proposition}~Let $V$ be a finite-dimensional inner product space.\\
{\rm (i)} The degree function
$$[V^{\infty},V^{\infty}] \to
\begin{cases}
\ZZ&{\rm if}~{\rm dim}(V)>0\\
\{0,1\}&{\rm if}~{\rm dim}(V)=0
\end{cases}~;~F \mapsto {\rm degree}(F)$$
is a bijection of pointed sets,
which is an isomorphism of abelian groups for ${\rm dim}(V)>0$.\\
{\rm (ii)} For any pointed map $F:V^{\infty} \to V^{\infty}$ and finite-dimensional
inner product space $W$ the pointed map
$$F'~=~F\wedge 1 ~:~(V \oplus W)^{\infty} ~=~V^{\infty}\wedge W^{\infty} \to
(V \oplus W)^{\infty} ~=~V^{\infty} \wedge W^{\infty}$$
has
$${\rm degree}(F')~=~{\rm degree}(F) \in \ZZ~.$$
\hfill\qed
\end{proposition}

We now move on to the bi-degree itself.

\begin{proposition}~ \label{half}
Let $V$ be a non-zero inner product space.\\
{\rm (i)} The cellular $\ZZ[\ZZ_2]$-module chain complexes of $S(LV)$
and $LV^{\infty}$ are given by
$$\begin{array}{l}
C^{cell}(S(LV))~:~C(S(LV))_{{\rm dim}(V)-1}=\ZZ[\ZZ_2]
\xymatrix@C+20pt{\ar[r]^-{1+(-)^{{\rm dim}(V)-1}T}&}\\[1ex]
\hskip25pt C(S(LV))_{{\rm dim}(V)-2}=\ZZ[\ZZ_2]
\xymatrix@C+20pt{\ar[r]^-{1+(-)^{{\rm dim}(V)}T}&}
\dots \xymatrix{\ar[r]^-{1-T}&} C(S(LV))_0=\ZZ[\ZZ_2]~,\\[2ex]
\dot C^{cell}(LV^{\infty})~=~S^{{\rm dim}(V)}(\ZZ,(-1)^{{\rm dim}(V)})
\end{array}$$
where $(\ZZ,\pm 1)$ denotes the $\ZZ[\ZZ_2]$-module $\ZZ$ with $T \in \ZZ_2$
acting by $\pm 1$.\\
{\rm (ii)} A $\ZZ[\ZZ_2]$-module chain map
$$e~:~\dot C^{cell}(\Sigma S(LV)^+)~=~SC^{cell}(S(LV))\to\dot C^{cell}(LV^{\infty})$$
induces a $\ZZ[\ZZ_2]$-module morphism
$$\begin{array}{l}
e_*~:~H_{{\rm dim}(V)-1}(S(LV))~=~
\begin{cases}
(\ZZ,(-)^{{\rm dim}(V)})&{\rm if}~{\rm dim}(V)>1\\
\ZZ[\ZZ_2]&{\rm if}~{\rm dim}(V)=1
\end{cases} \\[3ex]
\hskip125pt
\to \dot H_{{\rm dim}(V)}(LV^{\infty})~=~(\ZZ,(-)^{{\rm dim}(V)})
\end{array}$$
with $\begin{cases} e_*(1) \in 2 \ZZ\\ e_*(1)=-e_*(T)\in \ZZ.\end{cases}$
\end{proposition}
\begin{proof} (i) By construction.\\
(ii) The generator $1 \in H_{{\rm dim}(V)-1}(S(LV))=\ZZ$
is represented by
$$1+(-)^{{\rm dim}(V)}T \in C(S(LV))_{{\rm dim}(V)-1}~=~\ZZ[\ZZ_2]~.$$
A $\ZZ[\ZZ_2]$-module morphism
$$e~:~C(S(LV))_{{\rm dim}(V)-1}~=~\ZZ[\ZZ_2]  \to (\ZZ,(-)^{{\rm dim}(V)})$$
is of the form
$$e(a+bT)~=~na + (-)^{{\rm dim}(V)} nb~~(a,b \in \ZZ)$$
for some $n \in \ZZ$, and
$$\begin{cases}
e(1+(-)^{{\rm dim}(V)}T)=2n&{\rm if}~{\rm dim}(V)>1\\
e(1)=-e(T)=n&{\rm if}~{\rm dim}(V)=1~.
\end{cases}$$
\hfill\qed\end{proof}

\begin{definition} (i) The {\it semidegree} of a $\ZZ[\ZZ_2]$-module
chain map\index{semidegree!chain map}
$$e~:~\dot C^{cell}(\Sigma S(LV)^+)\to\dot C^{cell}(LV^{\infty})$$
is
$$\hbox{\rm semidegree}(e)~=~\begin{cases}
e_*(1)/2 \in \ZZ& {\rm if}~{\rm dim}(V)>1\\
e_*(1)=-e_*(T)\in \ZZ&{\rm if}~{\rm dim}(V)=1.
\end{cases}$$
(ii) Let $V,W$ be finite-dimensional inner product spaces.
The {\it semidegree} of a $\ZZ_2$-equivariant pointed map
\index{semidegree!pointed map}
$$E~:~\Sigma S(LV)^+\wedge W^{\infty} \to LV^{\infty}\wedge W^{\infty}$$
is the semidegree of the induced $\ZZ[\ZZ_2]$-module chain map, that is
$$\hbox{\rm semidegree}(E)~=~\begin{cases}
E_*(1)/2 \in \ZZ& {\rm if}~{\rm dim}(V)>1\\
E_*(1)=-E_*(T)\in \ZZ&{\rm if}~{\rm dim}(V)=1
\end{cases}$$
with
$$\begin{array}{l}
E_*~:~\dot H_{{\rm dim}(V)+{\rm dim}(W)}(\Sigma S(LV)^+\wedge W^{\infty})~=~
\begin{cases} \ZZ&{\rm if}~{\rm dim}(V)>1  \\
\ZZ[\ZZ_2]&{\rm if}~{\rm dim}(V)=1 \end{cases}\\[1ex]
\hskip100pt
\to \dot H_{{\rm dim}(V)+{\rm dim}(W)}(LV^{\infty}\wedge W^{\infty})~=~\ZZ~.
\end{array}$$
\hfill\qed
\end{definition}

\begin{proposition}~ {\rm (i)} The function
$$H_0({\rm Hom}_{\ZZ[\ZZ_2]}(\dot C^{cell}(\Sigma S(LV)^+),\dot C^{cell}(LV^{\infty})))
\to \ZZ~;~e \mapsto \hbox{\rm semidegree}(e)$$
is an isomorphism of abelian groups. \\
{\rm (ii)} The semidegree of a $\ZZ_2$-equivariant pointed map
$$E~:~\Sigma S(LV)^+\wedge W^{\infty} \to LV^{\infty}\wedge W^{\infty}$$
is a $\ZZ_2$-equivariant homotopy invariant, such that for any $j,k \geqslant 0$
the semidegrees of $E$ and the $\ZZ_2$-equivariant pointed map defined by
$$\begin{array}{l}
\Sigma^{j,k}E~:~\Sigma S(LV\oplus L\R^k)^+\wedge (W\oplus\R^j)^{\infty} \\[1ex]
\xymatrix@C+20pt{\ar[r]^-{\di{j_{LV,L\R} \wedge 1}}&}
\Sigma S(LV)^+ \wedge (L\R^k)^{\infty}\wedge W^{\infty} \wedge (\R^j)^{\infty}\\[1ex]
\xymatrix@C+10pt{\ar[r]^-{\di{E \wedge 1}}&}
LV^{\infty} \wedge (L\R^k)^{\infty}\wedge W^{\infty}
\wedge (\R^j)^{\infty} =~(LV\oplus L\R^k)^{\infty}\wedge (W\oplus \R^j)^{\infty}
\end{array}$$
are related by
$$\hbox{\rm semidegree}(\Sigma^{j,k}E)~=~\hbox{\rm semidegree}(E) \in \ZZ~.$$
\end{proposition}
\begin{proof}
(i) For each $n \in \ZZ$ construct a chain map $e$ with semidegree $n$
as in the proof of Proposition \ref{half}.\\
(ii) The only non-trivial case is for $V=\R$, $W=\{0\}$, $j=0$, $k=1$.
Consider the commutative diagram
$$\xymatrix@C+10pt{\dot H_2(\Sigma S(L\R^2)^+)=\ZZ \ar[r]^-{\di{\Sigma^{0,1}E_*}}
\ar[d]^-{\di{1-T}} & \dot H_2((L\R^2)^{\infty})=\ZZ \ar[d]^-{\di{1}}\\
\dot H_1(\Sigma S(L\R)^+)=\ZZ[\ZZ_2] \ar[r]^-{\di{E_*}}
& \dot H_1(L\R^{\infty})=\ZZ}$$
The evaluations of the two composites on the generator $1 \in \dot H_2(\Sigma S(L\R^2)^+)=\ZZ$
are
$$\begin{array}{l}
E_*(1-T)(1)~=~2\,\hbox{\rm semidegree}(E)~,\\[1ex]
\Sigma^{1,0}E_*(1)~=~2\,\hbox{\rm semidegree}(\Sigma^{1,0}E) \in \ZZ~,
\end{array}$$
so
$$\hbox{\rm semidegree}(\Sigma^{0,1}E)~=~\hbox{\rm semidegree}(E) \in \ZZ~.$$
\hfill\qed\end{proof}

\begin{definition}~
{\rm Given a $\ZZ_2$-equivariant pointed map
$$F~:~LV^{\infty}\wedge W^{\infty} \to LV^{\infty}\wedge W^{\infty}$$
define the pointed map of the $\ZZ_2$-fixed point sets
$$\begin{array}{l}
G~=~\rho(F)~:~(LV^{\infty}\wedge W^{\infty})^{\ZZ_2}~=~W^{\infty} \to W^{\infty}~;\\[2ex]
\hskip125pt w \mapsto \begin{cases} w'&{\rm if}~F(0,w)=(0,w')\\[1ex]
\infty&{\rm if}~F(0,w)=\infty~.
\end{cases}
\end{array}$$
The $\ZZ_2$-equivariant pointed map
$$\sigma(G)~=~1 \wedge G~:~LV^{\infty}\wedge W^{\infty} \to
LV^{\infty}\wedge W^{\infty}~;~(v,w) \mapsto (v,G(w))$$
is such that
$$F(0,w)~=~(0,G(w))~=~\sigma(G)(0,w) \in LV^{\infty} \wedge W^{\infty}~~(w \in W)~,$$
so that the relative difference is defined, a $\ZZ_2$-equivariant pointed map
$$\delta(F,\sigma(G))~:~\Sigma S(LV)^+ \wedge W^{\infty} \to LV^{\infty} \wedge W^{\infty}~.$$
The {\it bi-degree} of $F$ is defined by\index{bi-degree}
$$\begin{array}{ll}
\hbox{\rm bi-degree}(F)&=~(\hbox{\rm semidegree}(\delta(F,\sigma(G))),{\rm degree}(G))\\[1ex]
&\in \begin{cases}\ZZ \oplus \ZZ& {\rm if}~{\rm dim}(W)>0\\
\ZZ \times \{0,1\} & {\rm if}~{\rm dim}(W)=0~.
\end{cases}
\end{array}$$
}
\hfill\qed
\end{definition}

\begin{proposition}~ \label{00}
{\rm (i)} The bi-degree of a $\ZZ_2$-equivariant pointed map
$F:LV^{\infty}\wedge W^{\infty} \to LV^{\infty}\wedge W^{\infty}$
is
$$\begin{array}{ll}
\hbox{\rm bi-degree}(F)&=~\big(
\dfrac{{\rm degree}(F)-{\rm degree}(G)}{2},{\rm degree}(G)) \\[1ex]
&\in \begin{cases}\ZZ \oplus \ZZ&{\it if}~{\rm dim}(W)>0\\
\ZZ \times \{0,1\} &{\it if}~{\rm dim}(W)=0
\end{cases}
\end{array}$$
with $G=\rho(F):W^{\infty} \to W^{\infty}$. In particular
$${\rm deg}(F) \equiv {\rm deg}(G)~ (\bmod\, 2)~.$$
{\rm (ii)} The degree of a $\ZZ_2$-equivariant pointed map
$F:LV^{\infty} \to LV^{\infty}$ is odd
$${\rm deg}(F) \equiv 1~ (\bmod\, 2)~,$$
and for any pointed map $G:W^{\infty} \to W^{\infty}$ the bi-degree of the
$\ZZ_2$-equivariant pointed map
$$F \wedge G~:~LV^{\infty}\wedge W^{\infty} \to LV^{\infty}\wedge W^{\infty}$$
is
$$\begin{array}{ll}
\hbox{\rm bi-degree}(F \wedge G)&=~\big(
\dfrac{({\rm degree}(F)-1)}{2}{\rm degree}(G),{\rm degree}(G)) \\[1ex]
&\in \begin{cases}\ZZ \oplus \ZZ& {\rm if}~{\rm dim}(W)>0\\
\ZZ \times \{0,1\} & {\rm if}~{\rm dim}(W)=0~.
\end{cases}
\end{array}$$
{\rm (iii)} If $F_1,F_2:LV^{\infty}\wedge W^{\infty} \to LV^{\infty}\wedge W^{\infty}$
are $\ZZ_2$-equivariant pointed maps with
$$\rho(F_1)~=~\rho(F_2)~=~G~:~W^{\infty} \to W^{\infty}$$
then
$$\hbox{\rm bi-degree}(F_1)-\hbox{\rm bi-degree}(F_2)~=~\big(
\dfrac{{\rm degree}(F_1)-{\rm degree}(F_2)}{2},0) \in \ZZ \oplus \ZZ~.$$
In particular
$${\rm degree}(F_1)~\equiv~{\rm degree}(F_2)~(\bmod\, 2)~,$$
with
$$\hbox{\rm degree}(F_1)-\hbox{\rm degree}(F_2)~=~2\,\hbox{\rm semidegree}(\delta(F_1,F_2)) \in \ZZ~.$$
{\rm (iv)} The degree, semidegree and bi-degree functions define an
isomorphism of direct sum systems of abelian groups
$$\xymatrix@C+35pt{
[\Sigma S(LV)^+ \wedge W^{\infty},LV^{\infty} \wedge W^{\infty}]_{\ZZ_2}
\ar[r]^-{\di{\hbox{\rm semidegree}}}_-{\di{\cong}}
\ar[d]^-{\di{\gamma=\alpha_{LV}^*}} & \ZZ\ar[d]^-{\begin{pmatrix} 1 \\ 0 \end{pmatrix}} \\
[LV^{\infty} \wedge W^{\infty},LV^{\infty} \wedge W^{\infty}]_{\ZZ_2}
\ar[r]^-{\di{\hbox{\rm bi-degree}}}_-{\di{\cong}}
\ar@<5pt>[u]^-{\di{\delta}} \ar[d]^-{\di{\rho=s_{LV}^*}}&
\ZZ \oplus \ZZ\ar[d]^-{\begin{pmatrix} 0 & 1 \end{pmatrix}}
\ar@<5pt>[u]^-{\begin{pmatrix} 1 & 0 \end{pmatrix}} \\
[W^{\infty},LV^{\infty} \wedge W^{\infty}]_{\ZZ_2}=
[W^{\infty},W^{\infty}] \ar@<5pt>[u]^-{\di{\sigma}}
\ar[r]^-{\di{\rm degree}}_-{\di{\cong}}&
\ZZ \ar@<5pt>[u]^-{\begin{pmatrix} 0 \\ 1 \end{pmatrix}}}$$
for non-zero $W$. The forgetful map
$$[LV^{\infty} \wedge W^{\infty},LV^{\infty} \wedge W^{\infty}]_{\ZZ_2}~=~
\ZZ \oplus \ZZ \to [LV^{\infty} \wedge W^{\infty},LV^{\infty} \wedge W^{\infty}]~=~\ZZ$$
sends a $\ZZ_2$-equivariant map
$F:LV^{\infty}\wedge W^{\infty} \to LV^{\infty}\wedge W^{\infty}$
with bi-degree $(a,b)\in \ZZ\oplus \ZZ$ to a map $F$ with
$${\rm degree}(F)~=~2a+b \in \ZZ~.$$
For $W=\{0\}$ the functions define a bijection of
products of pointed sets
$$\xymatrix@C+35pt{
[\Sigma S(LV)^+ ,LV^{\infty} ]_{\ZZ_2}
\ar[r]^-{\di{\hbox{\rm semidegree}}}_-{\di{\cong}}
\ar[d]^-{\di{\alpha_{LV}^*}} & \ZZ\ar[d]\\
[LV^{\infty} ,LV^{\infty} ]_{\ZZ_2}
\ar[r]^-{\di{\hbox{\rm bi-degree}}}_-{\di{\cong}}
\ar@<5pt>[u]^-{\di{\delta}} \ar[d]^-{\di{s_{LV}^*}}&
\ZZ \times \{0,1\}\ar[d]\ar@<5pt>[u] \\
[S^0,LV^{\infty}]_{\ZZ_2}=[S^0,S^0] \ar@<5pt>[u]^-{\di{\sigma}}
\ar[r]^-{\di{\rm degree}}_-{\di{\cong}}& \{0,1\} \ar@<5pt>[u]}$$
\end{proposition}
\begin{proof} (i)
The map $\alpha_{LV}:LV^{\infty} \to \Sigma S(LV)^+$ induces
$$\begin{array}{l}
(\alpha_{LV})_*~=~\begin{cases} 1 \\ 1-T \end{cases}~:~
\dot H_{{\rm dim}(V)}(LV^{\infty})~=~ \ZZ\\[1ex]
\to \dot H_{{\rm dim}(V)}(\Sigma S(LV)^+)~=~
H_{{\rm dim}(V)-1}(S(LV)) ~=~\begin{cases}
\ZZ \\ \ZZ[\ZZ_2] \end{cases}~{\rm if}~
\begin{cases}
{\rm dim}(V)>1\\
{\rm dim}(V)=1~.
\end{cases}
\end{array}$$
The composite
$$F-\sigma(G)~:~\xymatrix@C+10pt{LV^{\infty} \wedge W^{\infty} \ar[r]^-{\di{\alpha_{LV}\wedge 1}}&
\Sigma S(LV)^+ \wedge W^{\infty} \ar[r]^-{\di{\delta(F,\sigma(G))}} &
LV^{\infty} \wedge W^{\infty}}$$
has degree
$${\rm degree}(\delta(F,\sigma(G))(\alpha_{LV} \wedge 1))~
=~{\rm degree}(F)-{\rm degree}(\sigma(G))$$
so
$$\begin{array}{ll}
\hbox{\rm semidegree}(\delta(F,\sigma(G)))&=~
\dfrac{{\rm degree}(F)-{\rm degree}(\sigma(G))}{2}\\[1ex]
&=~\dfrac{{\rm degree}(F)-{\rm degree}(G)}{2} \in \ZZ~.
\end{array}$$
(ii) Immediate from (i).\\
(iii) Taking $X=Y=S^0$ in Proposition \ref{sequence} (ii), we have a
direct sum system
$$\begin{array}{l}
[\Sigma S(LV)^+ \wedge W^{\infty},LV^{\infty} \wedge W^{\infty}]_{\ZZ_2}
\xymatrix{\ar@<1ex>[r]^-{\di{\gamma}}
\ar@<-1ex>@{<-}[r]_-{\di{\delta}}&}\\[1ex]
\hskip50pt
[LV^{\infty} \wedge W^{\infty},LV^{\infty} \wedge W^{\infty}]_{\ZZ_2}
\xymatrix{\ar@<1ex>[r]^-{\di{\rho}} \ar@<-1ex>@{<-}[r]_-{\di{\sigma}}&}
[W^{\infty},W^{\infty}]~.
\end{array}$$
The function
$$\begin{array}{l}
(\alpha_{LV}^*,\sigma)~:~[\Sigma S(LV)^+ \wedge W^{\infty},
LV^{\infty} \wedge W^{\infty}]_{\ZZ_2} \times [W^{\infty},W^{\infty}]\\[1ex]
\hskip125pt
\xymatrix{\ar[r]&}
[LV^{\infty} \wedge W^{\infty},LV^{\infty} \wedge W^{\infty}]_{\ZZ_2}
\end{array}$$
is an isomorphism such that
$$\begin{array}{l}
F~=~(\alpha^*_{LV},\sigma)(\delta(F,\sigma(G)),G)\\[1ex]
\in {\rm im}((\alpha_{LV}^*,\sigma):[\Sigma S(LV)^+ \wedge W^{\infty},
LV^{\infty} \wedge W^{\infty}]_{\ZZ_2} \times [W^{\infty},W^{\infty}]\\[1ex]
\hskip125pt
\xymatrix{\ar[r]^-{\di{\cong}}&}
[LV^{\infty} \wedge W^{\infty},LV^{\infty} \wedge W^{\infty}]_{\ZZ_2})~.
\end{array}$$
\hfill\qed\end{proof}

\begin{example} {\rm
Given a $\ZZ_2$-equivariant pointed map $f:LV^{\infty} \to LV^{\infty}$
and a pointed map $g:W^{\infty} \to W^{\infty}$ define the $\ZZ_2$-equivariant
pointed map
$$F~=~f\wedge g~:~LV^{\infty}\wedge W^{\infty} \to LV^{\infty}\wedge W^{\infty}~;~
(v,w) \mapsto (f(v),g(w))~.$$
The $\ZZ_2$-fixed point map is
$$G~=~\rho(F)~:~W^{\infty} \to W^{\infty}~;~w \mapsto \begin{cases}
g(w)&{\rm if}~f(0)=0 \\
\infty&{\rm if}~f(0)=\infty
\end{cases}$$
and
$$\sigma(G)~:~LV^{\infty}\wedge W^{\infty} \to LV^{\infty} \wedge W^{\infty}~;~
(v,w) \mapsto \begin{cases}
(v,g(w))&{\rm if}~f(0)=0 \\
\infty&{\rm if}~f(0)=\infty~,
\end{cases}$$
so that
$$\begin{array}{ll}
\hbox{\rm bi-degree}(F)&=~(\hbox{\rm semidegree}(\delta(F,\sigma(G))),
\hbox{\rm degree}(G))\\[2ex]
&=~(\dfrac{{\rm degree}(F)-{\rm degree}(G)}{2},{\rm degree}(G))\\[2ex]
&=~\begin{cases}
\big(\dfrac{({\rm degree}(f)-1){\rm degree}(g)}{2},{\rm degree}(g)\big)&{\rm if}~f(0)=0 \\[1ex]
\big(\dfrac{{\rm degree}(f){\rm degree}(g)}{2},0)&{\rm if}~f(0)=\infty~.
\end{cases}
\end{array}$$
}
\hfill\qed
\end{example}

\begin{example} {\rm
The $\ZZ_2$-equivariant pointed map defined for any $\lambda\in \R$ by
$$F~=~\lambda~:~L\R^{\infty} \to L\R^{\infty}~;~x \mapsto x/\lambda$$
has
$$\begin{array}{l}
G~=~\rho(F)~:~(L\R^{\infty})^{\ZZ_2}~=~\{0,\infty\} \to \{0,\infty\}~;~
\infty \mapsto \infty~,~
0 \mapsto \begin{cases}
0&\hbox{if $\lambda \neq 0$}\\
\infty&\hbox{if $\lambda =0$,}
\end{cases}\\[3ex]
\hbox{\rm bi-degree}(F)~=~(\hbox{\rm semidegree}\,\delta(F,\sigma(G)),{\rm degree}\,G)~
=~\begin{cases}
(0,1)&\hbox{if $\lambda >0$}\\
(0,0)&\hbox{if $\lambda =0$}\\
(-1,1)&\hbox{if $\lambda <0$}
\end{cases}\\[5ex]
\hskip25pt
\in [L\R^{\infty},L\R^{\infty}]_{\ZZ_2}~=~
[\Sigma S(L\R)^+,L\R^{\infty}]_{\ZZ_2} \times [S^0,S^0]~=~\ZZ \times \{0,1\}~.
\end{array}$$
}
\hfill\qed
\end{example}

\section{Stable $\ZZ_2$-equivariant homotopy theory}

\begin{definition}\label{stabdef}~{\rm
The {\it stable $\ZZ_2$-equivariant homotopy group} of pointed
$\ZZ_2$-spaces $X,Y$ is defined by
\index{stable $\ZZ_2$-equivariant homotopy!group, $\{X;Y\}_{\ZZ_2}$}
$$\{X;Y\}_{\ZZ_2}~=~\mathop{\varinjlim}\limits_V\,
[V^{\infty} \wedge X,V^{\infty}  \wedge Y]_{\ZZ_2}$$
with the direct limit running over all the finite-dimensional
inner product $\ZZ_2$-spaces $V$.\\
\hfill\qed}
\end{definition}

\begin{remark} {\rm In dealing with
$$\{X;Y\}_{\ZZ_2}~=~\mathop{\varinjlim}\limits_j\,
\mathop{\varinjlim}\limits_k\,[(L\R^j\oplus \R^k)^{\infty} \wedge X,
(L\R^j \oplus \R^k)^{\infty} \wedge Y]_{\ZZ_2}$$
there is no loss of generality in only considering $j=k$. Noting that
$(L\R^k \oplus \R^k)^{\infty} \wedge X$ is
$\ZZ_2$-equivariantly homeomorphic to
$$\Sigma^{k,k}X~=~S^k \wedge S^k \wedge X~,~
T~:~\Sigma^{k,k}X \to \Sigma^{k,k}X~;~(s,t,x) \mapsto (t,s,Tx)$$
 (see Definition \ref{transposition} below for details) we can identify
$$\{X;Y\}_{\ZZ_2}~=~\mathop{\varinjlim}\limits_k\,
[\Sigma^{k,k}X,\Sigma^{k,k}Y]_{\ZZ_2}~.$$
\hfill\qed}
\end{remark}

\begin{example} \label{00stable} {\rm
For any non-zero inner product space $V$ Proposition \ref{00} (ii)
gives an isomorphism of direct sum systems of abelian groups
$$\xymatrix@C+35pt{
\{\Sigma S(LV)^+;LV^{\infty}\}_{\ZZ_2}
\ar[r]^-{\di{\hbox{\rm semidegree}}}_-{\di{\cong}}
\ar[d]^-{\di{\gamma=\alpha_{LV}^*}} & \ZZ\ar[d]^-{\begin{pmatrix} 1 \\ 0 \end{pmatrix}}\\
\{LV^{\infty};LV^{\infty}\}_{\ZZ_2}=\{S^0;S^0\}_{\ZZ_2}
\ar[r]^-{\di{\hbox{\rm bi-degree}}}_-{\di{\cong}}
\ar@<5pt>[u]^-{\di{\delta}} \ar[d]^-{\di{\rho=s_{LV}^*}}&
\ZZ \oplus \ZZ\ar[d]^-{\begin{pmatrix} 0 & 1 \end{pmatrix}}
\ar@<5pt>[u]^-{\begin{pmatrix} 1 & 0 \end{pmatrix}} \\
\{S^0;LV^{\infty}\}_{\ZZ_2}=\{S^0;S^0\} \ar@<5pt>[u]^-{\di{\sigma}}
\ar[r]^-{\di{\rm degree}}_-{\di{\cong}}& \ZZ \ar@<5pt>[u]^-{\begin{pmatrix} 0 \\ 1 \end{pmatrix}}\\
}$$
\hfill\qed
}
\end{example}

\begin{proposition}\label{slr}~Let $X,Y$ be pointed $\ZZ_2$-spaces.\\
{\rm (i)} Given a pointed map $F:\vert X \vert \to \vert Y \vert$ define
the $\ZZ_2$-equivariant map
$$F\vee TFT~:~S(L\R)^+\wedge X \to Y~;~
\begin{cases}
(+,x) \mapsto F(x)\\
(-,x) \mapsto TFT(x)~.
\end{cases}$$
The functions
$$\begin{array}{l}
\{\vert X\vert ;\vert Y\vert\} \to \{S(L\R)^+\wedge X;Y\}_{\ZZ_2}~;~F \mapsto F\vee TFT~,\\[1ex]
\{S(L\R)^+\wedge X;Y\}_{\ZZ_2}\to \{\vert X\vert ;\vert Y\vert\} ~;~G \mapsto G\vert_{(+,X)}
\end{array}$$
are inverse isomorphisms of abelian groups. The $\ZZ_2$-equivariant homotopy cofibration sequence
$$S(L\R)^+ \to S^0 \to L\R^{\infty} \to \Sigma S(L\R)^+ \to \dots $$
induces the exact sequence
$$\begin{array}{l}
\dots \to \{\Sigma S(L\R)^+ \wedge X,Y\}_{\ZZ_2}=\{\Sigma X,Y\}\to
\{L\R^{\infty} \wedge X,Y\}_{\ZZ_2}\\[1ex]
\hskip50pt \to \{X,Y\}_{\ZZ_2} \to \{S(L\R)^+ \wedge X,Y\}_{\ZZ_2}=\{\vert X \vert,\vert Y \vert\}
\to \dots~.
\end{array}$$
{\rm (ii)} Given a pointed map $F:\vert X \vert \to \vert Y \vert$ define
the $\ZZ_2$-equivariant map
$$\begin{array}{l}
F+TFT~:~L\R^{\infty}\wedge X \to L\R^{\infty} \wedge (S(L\R)^+ \wedge Y)~;\\[1ex]
\hskip50pt
(t,x) \mapsto \begin{cases}
(2t,(+,F(x)))&{\it if}~0 \leqslant t \leqslant 1/2\\
(2t-1,(-,TFT(x)))&{\it if}~1/2 \leqslant t \leqslant 1~.
\end{cases}
\end{array}$$
The function
$$q~:~\{\vert X\vert ;\vert Y\vert\} \to \{X;S(L\R)^+\wedge Y\}_{\ZZ_2}~;~F \mapsto F+TFT$$
is an isomorphism of abelian groups.
\end{proposition}
\begin{proof} (i) By Proposition \ref{forget} (i)
$$[S(L\R)^+\wedge V^{\infty} \wedge X,V^{\infty} \wedge Y]_{\ZZ_2}
~=~[\vert V^{\infty}\vert  \wedge \vert X\vert ,\vert V^{\infty}\vert  \wedge \vert Y\vert ]~.$$
The exact sequence is the stable version of \ref{forget} (iii).\\
(ii) Use the projection $p:S(L\R)^+\wedge Y \to S(L\R)^+\wedge_{\ZZ_2} Y$
and the homeomorphism
$$\vert Y \vert \to S(L\R)^+\wedge_{\ZZ_2} Y~;~y \mapsto (+,y)$$
to define a morphism
$$p_*~:~\{X;S(L\R)^+\wedge Y\}_{\ZZ_2} \to
\{\vert X\vert;S(L\R)^+\wedge _{\ZZ_2}Y\}~=~\{\vert X\vert ;\vert Y\vert\}$$
such that $p_*q=1$, $qp_*=1$.\\
\hfill\qed\end{proof}

Given pointed $\ZZ_2$-spaces $A,B$ and a $\ZZ_2$-equivariant pointed
map $F:A \to B$ let
$G:B \to \Cc(F)$ be the inclusion in the mapping cone, and let
$H:\Cc(F) \to \Sigma A$ be the projection, so that
$$\xymatrix@C+10pt{A \ar[r]^-{\di{F}} &
B \ar[r]^-{\di{G}} &
\Cc(F)\ar[r]^-{\di{H}}& \Sigma A \ar[r]^-{\di{\Sigma F}} & \Sigma B \ar[r] &\dots}$$
is a $\ZZ_2$-equivariant homotopy cofibration sequence.
By analogy with the nonequivariant case (\ref{BP2}):

\begin{proposition}~\label{BP3}
For any pointed $\ZZ_2$-space $X$ there is induced a
Barratt-Puppe exact sequence of stable $\ZZ_2$-equivariant homotopy groups
$$\begin{array}{l}
\xymatrix{\dots \ar[r]& \{\Sigma X;\Cc(F)\}_{\ZZ_2} \ar[r]^-{\di{H}} &
\{X;A\}_{\ZZ_2}\ar[r]^-{\di{F}} & \{X;B\}_{\ZZ_2} }\\
\hskip100pt
\xymatrix{\ar[r]^-{\di{G}} &\{X;\Cc(F)\}_{\ZZ_2} \ar[r]^-{\di{H}} & \{X;\Sigma A\}_{\ZZ_2} \ar[r] & \dots~.}
\end{array}$$
\hfill\qed
\end{proposition}

By analogy with the braid of stable homotopy groups of Proposition
\ref{stablesequence0} :

\begin{proposition}~ \label{stablesequenceZ2}
For any inner product $\ZZ_2$-spaces $U,V$ and pointed $\ZZ_2$-spaces $X,Y$
there is defined a commutative braid of exact sequences of stable
homotopy groups
$$\xymatrix@C-30pt{
A_1 \ar[dr] \ar@/^2pc/[rr]^-{\di{\alpha_V}} && \{X;S(V)^+ \wedge Y\}_{\ZZ_2}
\ar[dr]^-{\di{s_V}}\ar@/^2pc/[rr] &&
\{S(U)^+ \wedge X;Y\}_{\ZZ_2}\ar[dr]&\\
&~~~A_2~~~ \ar[ur] \ar[dr] &&
\{X;Y\}_{\ZZ_2} \ar[ur]^-{\di{s^*_U}} \ar[dr]^-{\di{0_{V}}}&&A_3\\
A_4\ar[ur] \ar@/_2pc/[rr]^-{\di{\alpha^*_U}}
&&\{U^{\infty} \wedge X;Y\}_{\ZZ_2}
\ar@/_2pc/[rr]\ar[ur]^-{\di{0^*_U}} &&
\{X;V^{\infty}\wedge Y\}_{\ZZ_2}\ar[ur]&}$$

\bigskip

\noindent with
$$\begin{array}{l}
A_1~=~\{\Sigma X;LV^{\infty} \wedge Y\}_{\ZZ_2}~,~
A_2~=~\{\Sigma S(U\oplus V)^+\wedge X;V^{\infty} \wedge Y\}_{\ZZ_2}~,\\[1ex]
A_3~=~\{S(U\oplus V)^+ \wedge X;V^{\infty} \wedge Y\}_{\ZZ_2}~,~
A_4~=~\{\Sigma S(U)^+ \wedge X;Y\}_{\ZZ_2}~.
\end{array}$$
\end{proposition}
\begin{proof} These are the Barratt-Puppe exact sequences
(\ref{BP3}) determined by the homotopy commutative braid of $\ZZ_2$-equivariant homotopy cofibrations

$$\xymatrix@C-50pt{
&U^{\infty} \wedge S(V)^+\ar[dr]^-{\di{1 \wedge s_V}}
\ar@/^2pc/[rr] && \Sigma S(U)^+ \ar[dr]&\\
S(U\oplus V)^+\ar[dr]_-{\di{s_{U\oplus V}}} \ar[ur] &&
U^{\infty}\wedge S^0 =U^{\infty}
\ar[dr]^-{\di{1\wedge 0_V}}\ar[ur]^-{\di{\alpha_V}}&&\Sigma S(U \oplus V)^+&\\
&S^0\ar@/_2pc/[rr]^-{\di{0_{U \oplus V}}}
\ar[ur]^-{\di{0_U}} &&U^{\infty} \wedge V^{\infty}=
(U \oplus V)^{\infty}\ar[ur]^-{\di{\alpha_{U \oplus V}}} &}$$
\noindent with two homotopy pushouts.\\
\hfill\qed\end{proof}

\begin{definition} A {\it $\ZZ_2$-spectrum}\index{$\ZZ_2$-spectrum} $\underline{X}=\{X(V)\,\vert\, V\}$ is a
sequence of pointed $\ZZ_2$-spaces $X(V)$ indexed by finite-dimensional
inner product $\ZZ_2$-spaces $V$, with structure maps
$$(V^{\perp})^{\infty} \wedge X(V) \to X(W)$$
defined whenever $V \subseteq W$, where $V^{\perp} \subseteq W$ is
the orthogonal complement of $V$ in $W$. For $n \in \ZZ$ let
$$\pi^{\ZZ_2}_n(\underline{X})~=~
\begin{cases}
\varinjlim\limits_V [\Sigma^nV^{\infty},X(V)]_{\ZZ_2}&{\rm if}~n \geqslant 0\\[1ex]
\varinjlim\limits_V [V^{\infty},\Sigma^{-n}X(V)]_{\ZZ_2}&{\rm if}~n \leqslant -1~.
\end{cases}$$
\end{definition}

\begin{proposition}~\label{stablefixed}
Let $Y$ be a pointed $\ZZ_2$-space.\\
{\rm (i)}  For any pointed $\ZZ_2$-space $X$ the forgetful map
$\{X;Y\}_{\ZZ_2} \to \{\vert X\vert ;\vert Y\vert \}$ fits into a long exact sequence
$$\dots \to \{L\R^{\infty}\wedge X;Y\}_{\ZZ_2} \to
\{X;Y\}_{\ZZ_2} \to \{\vert X\vert ;\vert Y\vert \} \to \{L\R^{\infty}\wedge X;\Sigma Y\}_{\ZZ_2} \to \dots$$
{\rm (ii)} For any $CW$ $\ZZ_2$-spectrum $X$ the fixed point map
$$\rho~:~\{X;L\R(\infty)^{\infty}\wedge Y\}_{\ZZ_2} \to \{X^{\ZZ_2};Y^{\ZZ_2}\}$$
is an isomorphism, and the fixed point map
$$\begin{array}{l}
\rho~:~\{X;Y\}_{\ZZ_2} \to \{X^{\ZZ_2};Y^{\ZZ_2}\}~;\\[1ex]
(F:V^{\infty} \wedge X \to V^{\infty} \wedge Y) \mapsto
(\rho(F):V_+^{\infty} \wedge X^{\ZZ_2} \to V_+^{\infty} \wedge Y^{\ZZ_2})
\end{array}$$
fits into a long exact sequence
$$\begin{array}{l}
\dots \to \{X;S(\infty)^+\wedge Y\}_{\ZZ_2} \to
\{X;Y\}_{\ZZ_2} \xymatrix{\ar[r]^-{\di{\rho}}&} \{X^{\ZZ_2};Y^{\ZZ_2}\} \\[1ex]
\hskip100pt
\to \{X;\Sigma S(\infty)^+\wedge Y\}_{\ZZ_2} \to \dots~.
\end{array}$$
{\rm (iii)} For any semifree $CW$ $\ZZ_2$-spectrum $X$
$$\{X;S(\infty)^+\wedge Y\}_{\ZZ_2}~=~\{X;Y\}_{\ZZ_2}~,~
\{X^{\ZZ_2};Y^{\ZZ_2}\}~=~0~.$$
{\rm (iv)} For $CW$ $\ZZ_2$-spectrum $X$ with the trivial $\ZZ_2$-action
the fixed point map $\rho:\{X;Y\}_{\ZZ_2} \to \{X;Y^{\ZZ_2}\}$ is a surjection
which is split by the map $\sigma:\{X;Y^{\ZZ_2}\}\to \{X;Y\}_{\ZZ_2}$ induced by
$$\sigma~:~Y^{\ZZ_2}~=~(L\R(\infty)^{\infty} \wedge Y)^{\ZZ_2} \to L\R(\infty)^{\infty} \wedge Y~;~y \mapsto (0,y)~.$$
Thus  $\{X;Y\}_{\ZZ_2}$ fits into a direct sum system
$$\{X;S(\infty)^+\wedge Y\}_{\ZZ_2}
\xymatrix{\ar@<1ex>[r]^-{\di{\gamma}}\ar@<-1ex>@{<-}[r]_-{\di{\delta}} & }
\{X;Y\}_{\ZZ_2}
\xymatrix{\ar@<1ex>[r]^-{\di{\rho}}\ar@<-1ex>@{<-}[r]_-{\di{\sigma}}& }
\{X;Y^{\ZZ_2}\}$$
{\rm (v)} The direct sum system in {\rm (iv)}
is natural with respect to $\ZZ_2$-equivariant maps. A stable
$\ZZ_2$-equivariant map $F:V^{\infty} \wedge Y \to V^{\infty} \wedge Y'$
induces
$$\begin{array}{l}
F~=~\begin{pmatrix} 1 \wedge F & \delta(Fi_Y,i_{Y'}\rho(F))\\
0 & \rho(F) \end{pmatrix}~:~\\[3ex]
\{X;Y\}_{\ZZ_2}~=~\{X;S(\infty)^+\wedge Y\}_{\ZZ_2} \oplus \{X;Y^{\ZZ_2}\}\\[1ex]
\hskip50pt
\to
\{X;Y'\}_{\ZZ_2}~=~\{X;S(\infty)^+\wedge Y'\}_{\ZZ_2} \oplus \{X;{Y'}^{\ZZ_2}\}
\end{array}$$
with $i_Y:Y^{\ZZ_2} \to Y$, $i_{Y'}:{Y'}^{\ZZ_2} \to Y'$ the inclusions.
\end{proposition}
\begin{proof} (i) Immediate from the $\ZZ_2$-equivariant homotopy cofibration sequence
$$\xymatrix{
S(L\R)^+ \ar[r]^-{\di{s_{L\R}}}& S^0 \ar[r]^-{\di{0_{L\R}}}&
L\R^{\infty} \ar[r]^-{\di{\alpha_{L\R}}}& \Sigma S(L\R)^+  \ar[r]& \dots}$$
of Proposition \ref{cofib2} (iii) and the isomorphism
$$\{X;Y\} \to \{S(L\R)^+ \wedge X;Y\}_{\ZZ_2}~;~F \mapsto F \vee TF$$
given by Proposition \ref{forget}.\\
(ii) Combine Proposition \ref{fixed} (i) with the long exact sequence
$$\begin{array}{l}
\dots \to \{X;S(\infty)^+\wedge Y\}_{\ZZ_2} \to
\{X;Y\}_{\ZZ_2} \to  \{X;L\R(\infty)^{\infty} \wedge Y\}_{\ZZ_2} \\[1ex]
\hskip100pt
\to \{X;\Sigma S(\infty)^+\wedge Y\}_{\ZZ_2} \to \dots~.
\end{array}$$
induced by the $\ZZ_2$-equivariant homotopy cofibration sequence
$$\xymatrix@C+10pt{
S(\infty)^+ \ar[r]^-{\di{s_{S(\infty)}}}& S^0 \ar[r]^-{\di{0_{L\R(\infty)}}}&
L\R(\infty)^{\infty} \ar[r]^-{\di{\alpha_{L\R(\infty)}}}& \Sigma S(\infty)^+  \ar[r]& \dots}$$
given by Proposition \ref{cofib2} (iii).\\
(iii) This is the stable version of Proposition \ref{fixed}, with
each of the fixed point maps
$$\rho~:~[LV^{\infty} \wedge W^{\infty} \wedge X,
L\R(\infty)^{\infty} \wedge LV^{\infty} \wedge W^{\infty}\wedge  Y]_{\ZZ_2}
\to [W^{\infty} \wedge X^{\ZZ_2},W^{\infty} \wedge Y^{\ZZ_2}]$$
is an isomorphism.\\
(iv)+(v) By construction.\\
\hfill\qed\end{proof}

Proposition \ref{sequence} has an analogue for the stable $\ZZ_2$-equivariant
homotopy groups:

\begin{definition}~{\rm
(i) The {\it reduced stable $\ZZ_2$-equivariant homotopy} and {\it cohomotopy} groups
of a pointed $\ZZ_2$-space $X$ are
\index{reduced!stable $\ZZ_2$-equivariant!homotopy, $\widetilde{\omega}^{\ZZ_2}_n(X)$}
\index{reduced!stable $\ZZ_2$-equivariant!cohomotopy, $\widetilde{\omega}^n_{\ZZ_2}(X)$}
$$\widetilde{\omega}^{\ZZ_2}_n(X)~=~\{S^n;X\}_{\ZZ_2}~,~
\widetilde{\omega}^n_{\ZZ_2}(X)~=~\{X;S^n\}_{\ZZ_2}~(n \in \ZZ)~.$$
(ii) The {\it absolute stable $\ZZ_2$-equivariant homotopy} and {\it cohomotopy} groups of a $\ZZ_2$-space $X$ are
\index{absolute stable $\ZZ_2$-equivariant!homotopy, $\omega^{\ZZ_2}_n(X)$}
\index{absolute stable $\ZZ_2$-equivariant!cohomotopy, $\omega^n_{\ZZ_2}(X)$}
$$\omega^{\ZZ_2}_n(X)~=~\widetilde{\omega}^{\ZZ_2}_n(X^+)~,~\omega^n_{\ZZ_2}(X)~=~
\widetilde{\omega}^n_{\ZZ_2}(X^+)~(n \in \ZZ)~.$$
{\rm (iii)} The {\it stable $\ZZ_2$-equivariant homotopy groups of spheres}
 are\index{stable $\ZZ_2$-equivariant homotopy group!of sphere, $\omega^{\ZZ_2}_n$}
$$\omega^{\ZZ_2}_n~=~\omega^{\ZZ_2}_n(\{*\})~=~\omega^{-n}_{\ZZ_2}(\{*\})~=~
\{S^n;S^0\}_{\ZZ_2}~(n \in \ZZ)~.$$
\hfill\qed}
\end{definition}

Stable $\ZZ_2$-equivariant maps arise from double covers:

\begin{proposition}~ \label{Z2cover}\label{Adams}
{\rm (i)}  {\rm (Adams \cite[p.511]{adams3})} For any semifree
pointed $\ZZ_2$-space $X$ the projection
$f:X \to X/\ZZ_2$ induces a transfer $\ZZ[\ZZ_2]$-module chain map
$$f^!~:~\dot C(X/\ZZ_2) \to \dot C(X)~;~ y \mapsto \sum\limits_{x \in f^{-1}(y)}x$$
of the reduced singular chain complexes.
If $X$ is a finite pointed $CW$ $\ZZ_2$-complex the Umkehr map $f^!$ is induced by a
stable $\ZZ_2$-equivariant Umkehr  map
$$F~:~(LV \oplus W)^{\infty} \wedge X/\ZZ_2 \to(LV \oplus W)^{\infty} \wedge X~.$$
{\rm (ii)} Let $A$ be a pointed $\ZZ_2$-space with a $\ZZ_2$-equivariant map
$\epsilon:A \to [0,\infty)$ such that $\epsilon^{-1}(0)=\{*\} \subset A$,
with $* \in A$  the base point. For any non-zero inner product space $V$ let $X=S(LV)^+ \wedge A$.
The maps
$$\begin{array}{l}
f~=~projection~:~X \to X/\ZZ_2~=~S(LV)^+\wedge_{\ZZ_2}A~,\\[1ex]
g~:~LV \times X \to LV~;~(u,(v,a)) \mapsto \epsilon(a)(v+\dfrac{u}{1+\Vert u \Vert})
\end{array}$$
are such that there are defined an open $\ZZ_2$-equivariant embedding
$$e~=~(g,f)~:~LV \times X \emb LV \times X/\ZZ_2$$
and a $\ZZ_2$-equivariant Umkehr map
$$F~:~LV^{\infty} \wedge X/\ZZ_2 \to LV^{\infty} \wedge X~.$$
{\rm (iii)} {\rm (Adams \cite[Thms. 5.3,5.4,5.5]{adams3})}
Let $X$ be a semifree finite pointed $CW$ $\ZZ_2$-complex, let
$f:X \to X/\ZZ_2$ be the projection, and let
$F \in \{X/\ZZ_2;X\}_{\ZZ_2}$ be the stable $\ZZ_2$-equivariant Umkehr map of
{\rm (i)}. For any pointed finite $CW$ complex $Y$ the functions
$$\begin{array}{l}
\{X/\ZZ_2;Y\}\to \{X;Y\}_{\ZZ_2}~;~G \mapsto Gf~,\\[1ex]
\{Y;X/\ZZ_2\} \to \{Y;X\}_{\ZZ_2}~;~H \mapsto FH~,\\[1ex]
\mathop{\varinjlim}\limits_V [V^{\infty}\wedge X,V^{\infty}\wedge Y]_{\ZZ_2}
\to \{X;Y\}_{\ZZ_2}
\end{array}$$
are isomorphisms, with $V$ running over all finite-dimensional inner product
spaces (with the trivial $\ZZ_2$-action).
\end{proposition}
\begin{proof}
(i) There exist finite-dimensional inner product spaces $V,W$ and a
$\ZZ_2$-equivariant map
$$g~:~LV \oplus W \times X \to LV \oplus W$$
with the property that
$$e~=~(g,f)~:~LV \oplus W \times X \to LV \oplus W \times X/\ZZ_2~;~
(v,w,x) \mapsto (g(v,w),f(x))$$
restricts to an open $\ZZ_2$-equivariant embedding
$$e\vert~:~LV \oplus W \times (X \backslash \{x_0\})
 \to LV \oplus W \times X/\ZZ_2$$
and as $x$ approaches $x_0$ then
$g(LV \oplus W \times \{x\})$ approaches $\{0\} \subset LV \oplus W$.
The construction of Definition \ref{pointedumkehr}
gives a stable $\ZZ_2$-equivariant Umkehr  map
$$F~:~(LV \oplus W)^{\infty} \wedge X/\ZZ_2 \to(LV \oplus W)^{\infty} \wedge X$$
representing a stable $\ZZ_2$-equivariant homotopy class
$F \in \{X/\ZZ_2;X\}_{\ZZ_2}$ inducing the transfer chain map
$F=f^!:\dot C(X/\ZZ_2) \to \dot C(X)$.\\
(ii) The map $g$ has the property that as $a\in A$ approaches the base point
$*\in A$ then $g(LV \times (S(LV)^+ \wedge \{a\}))$ approaches
$0\in LV$. Then
$$e~=~(g,f)~:~LV \times X \emb LV \times X/\ZZ_2~;~
(u,(v,a)) \mapsto (g(u,(v,a)),f(v,a))$$
is an open $\ZZ_2$-equivariant embedding with $\ZZ_2$-equivariant Umkehr map
$$\begin{array}{l}
F~:~LV^{\infty} \wedge X/\ZZ_2\to LV^{\infty} \wedge X~;\\[1ex]
(u,(v,a)) \mapsto
\begin{cases}
(\dfrac{u- v}{1-\Vert u - v \Vert},(v,a))
&{\rm if}~(v,a) \in X,~\Vert u - v\Vert <1\\[2ex]
(\dfrac{u+v}{1-\Vert u + v \Vert},(-v,a))
&{\rm if}~(v,a) \in X,~\Vert u + v\Vert <1\\[2ex]
*&{\rm otherwise}~.
\end{cases}
\end{array}$$
(iii) See \cite{adams3}.\\
\hfill\qed\end{proof}

\begin{example} {\rm
The special case $Y=S^i$ of \ref{Adams} (iii) gives identifications
$$\widetilde{\omega}_i^{\ZZ_2}(X)~=~\widetilde{\omega}_i(X/\ZZ_2)~,~
\widetilde{\omega}^i_{\ZZ_2}(X)~=~\widetilde{\omega}^i(X/\ZZ_2)$$
for any finite pointed semifree $CW$ $\ZZ_2$-complex.\\
\hfill\qed}
\end{example}

\begin{definition} \label{transposition} {\rm
For any inner product space $V$ regard $V \oplus V$ as an inner product
$\ZZ_2$-space via the {\it transposition involution}
$\ZZ_2$-action\index{transposition!$T$}
$$T~:~ V \oplus V \to V \oplus V~;~(u,v) \mapsto (v,u)~.$$
The $\ZZ_2$-equivariant linear isomorphism\index{transposition!$\kappa_V$}
$$\kappa_V~:~LV\oplus V \iso  V\oplus V ~; ~(u,v)\mapsto (u+v,-u+v)$$
has inverse
$$\kappa^{-1}_V~:~V\oplus V \iso LV\oplus V ~; ~
(x,y)\mapsto \bigg(\dfrac{x-y}{2},\dfrac{x+y}{2}\bigg)$$
and induces a $\ZZ_2$-equivariant homeomorphism
$$\kappa_V~:~LV^{\infty} \wedge V^{\infty}\iso V^{\infty} \wedge V^{\infty}~.$$
The restriction of $\kappa_V$ is a $\ZZ_2$-equivariant homeomorphism
$$\kappa_V\vert ~:~(LV\backslash \{0\}) \times V \iso
(V\oplus V)\backslash \Delta(V)$$
which induces a $\ZZ_2$-equivariant homeomorphism
$$\kappa_V\vert~:~(LV\backslash \{0\})^{\infty}
\wedge V^{\infty} \iso ((V\oplus V)\backslash \Delta(V))^{\infty}$$
where $(LV\backslash \{0\})^{\infty}\cong \Sigma S(LV)^+$.
\hfill\qed}
\end{definition}

\begin{proposition}~ \label{stablesequence1}
{\rm (i)}  For any pointed $\ZZ_2$-spaces $X,Y$ and
inner product spaces $U,V$ there is defined
a commutative braid of exact sequences of stable $\ZZ_2$-equivariant
homotopy groups

$$\xymatrix@C-35pt{
A_1 \ar[dr] \ar@/^2pc/[rr]^-{\di{\alpha_{LV}}} && \{X;S(LV)^+ \wedge Y\}_{\ZZ_2}
\ar[dr]^-{\di{s_{LV}}}\ar@/^2pc/[rr] &&
\{S(LU)^+ \wedge X;Y\}_{\ZZ_2}\ar[dr]&\\
&~~~~~A_2~~~~~ \ar[ur] \ar[dr] &&
\{X;Y\}_{\ZZ_2} \ar[ur]^-{\di{s^*_{LU}}} \ar[dr]^-{\di{0_{LV}}}&&A_3\\
A_4\ar[ur] \ar@/_2pc/[rr]^-{\di{\alpha^*_{LU}}}
&&\{LU^{\infty} \wedge X;Y\}_{\ZZ_2}
\ar@/_2pc/[rr]\ar[ur]^-{\di{0^*_{LU}}} &&
\{X;LV^{\infty}\wedge Y\}_{\ZZ_2}\ar[ur]&}$$

\bigskip

\noindent with
$$\begin{array}{l}
A_1~=~\{\Sigma X;LV^{\infty} \wedge Y\}_{\ZZ_2}~,~
A_2~=~\{\Sigma S(LU\oplus LV)^+\wedge X;LV^{\infty} \wedge Y\}_{\ZZ_2}~,\\[1ex]
A_3~=~\{S(LU\oplus LV)^+ \wedge X;LV^{\infty} \wedge Y\}_{\ZZ_2}~,~
A_4~=~\{\Sigma S(LU)^+ \wedge X;Y\}_{\ZZ_2}~.
\end{array}$$
The morphism
$$0_{LV}~:~\{X;X \wedge X\}_{\ZZ_2}
 \to \{X;LV^{\infty} \wedge X \wedge X\}_{\ZZ_2}~=~
 \{V^{\infty} \wedge X;V^{\infty} \wedge LV^{\infty} \wedge X \wedge X\}_{\ZZ_2}$$
 sends
 $\Delta_X$ to $0_{LV} \wedge \Delta_X=(\kappa^{-1}_V \wedge 1)\Delta_{V^{\infty} \wedge X}$.\\
{\rm (ii)}
For any pointed space $X$ and any pointed $\ZZ_2$-space $Y$
$$\{X;Y\}_{\ZZ_2}~=~\{X;S(\infty)^+\wedge Y\}_{\ZZ_2}\oplus\{X;Y^{\ZZ_2}\}~,$$
with a direct sum system
$$\{X;S(\infty)^+ \wedge Y\}_{\ZZ_2}
\xymatrix{\ar@<1ex>[r]^-{\di{\gamma}}\ar@<-1ex>@{<-}[r]_-{\di{\delta}} & }
\{X;Y\}_{\ZZ_2}
\xymatrix{\ar@<1ex>[r]^-{\di{\rho}}\ar@<-1ex>@{<-}[r]_-{\di{\sigma}}& }
\{X;Y^{\ZZ_2}\}~.$$
The map $\gamma$ is induced by the projection
$s_{L\R(\infty)}:S(\infty)^+\to S^0$.
The map $\delta$ is defined by sending a $\ZZ_2$-equivariant map
$$F~:~U^{\infty} \wedge LV^{\infty} \wedge X \to U^{\infty} \wedge LV^{\infty} \wedge Y$$
to the stable $\ZZ_2$-equivariant homotopy class of the relative difference
of $F$ and
$$\begin{array}{l}
\sigma\rho(F)~:~
U^{\infty} \wedge LV^{\infty} \wedge X \to U^{\infty} \wedge LV^{\infty} \wedge Y~;\\[1ex]
\hskip100pt (u,v,x) \mapsto (w,v,y)~~(F(u,0,x)=(w,0,y))~,
\end{array}$$
that is
$$\begin{array}{ll}
\delta(F)~=~\delta(F,\sigma\rho(F))&\in
\varinjlim\limits_{U,V}\{U^{\infty} \wedge \Sigma S(LV)^+ \wedge X;U^{\infty} \wedge LV^{\infty} \wedge Y\}_{\ZZ_2}\\[1ex]
&=~\varinjlim\limits_V\{X;S(LV)^+\wedge Y\}_{\ZZ_2}~=~\{X;S(\infty)^+\wedge Y\}_{\ZZ_2}~.
\end{array}$$
{\rm (iii)}
For any inner product space $U$, any pointed space $X$
and any pointed $\ZZ_2$-space $Y$ there is defined a long exact sequence
$$\begin{array}{l}
\xymatrix{
\dots \ar[r]&\{\Sigma S(LU \oplus L\R(\infty))^+\wedge
X;L\R(\infty)^{\infty} \wedge Y\}_{\ZZ_2}&}\\
\hskip25pt \xymatrix{\ar[r] &\{LU^{\infty} \wedge X;Y\}_{\ZZ_2}
\ar[r]^-{\di{\rho}} & \{X;Y^{\ZZ_2}\}&}\\
\hskip25pt \xymatrix{\ar[r] &\{S(LU \oplus L\R(\infty))^+\wedge
X;L\R(\infty)^{\infty} \wedge Y\}_{\ZZ_2} \ar[r] &\dots}
\end{array}
$$
with $\rho$ defined by the fixed points of the $\ZZ_2$-action, and
$$\begin{array}{l}
\{S(LU \oplus L\R(\infty))^+\wedge
X;L\R(\infty)^{\infty} \wedge Y\}_{\ZZ_2}\\[1ex]
\hskip100pt =~ \{LU^\infty \wedge X;\Sigma S(LU \oplus L\R(\infty))^+
\wedge Y\}_{\ZZ_2}~.
\end{array}$$
\end{proposition}
\begin{proof}
(i) Immediate from the $\ZZ_2$-equivariant homotopy cofibration sequences
$$\begin{array}{l}
S(LU)^+ \to S^0 \to LU^{\infty} \to \Sigma S(LU)^+ \to \dots~,\\[1ex]
S(LV)^+ \to S^0 \to LV^{\infty} \to \Sigma S(LV)^+ \to \dots~.
\end{array}$$
(ii) This is the $\ZZ_2$-equivariant version of the braid of
Proposition \ref{stablesequence0}, noting that the fixed point
function
$$\rho~:~\{X;L\R(\infty)^{\infty} \wedge Y\}_{\ZZ_2} \to \{X;Y^{\ZZ_2}\}$$
is an isomorphism by Proposition \ref{stablefixed}. The identification
$$0_{LV} \wedge \Delta_X~=~(\kappa^{-1}_V \wedge 1)\Delta_{V^{\infty} \wedge X}$$
is immediate from the commutative triangle of $\ZZ_2$-equivariant maps
$$\xymatrix{V^{\infty} \ar[dr]_-{\di{\Delta_V}}
\ar[rr]^-{\di{0_{LV}}} && LV^{\infty} \wedge V^{\infty}\\
&V^{\infty} \wedge V^{\infty} \ar[ur]_-{\di{\kappa_V}} &}$$
(iii) The $\ZZ_2$-equivariant homotopy cofibration sequence
$$\xymatrix{S(\infty)^+ \ar[r]^-{\di{s_{L\R(\infty)}}} &S^0
\ar[r]^-{\di{0_{L\R(\infty)}}} & L\R(\infty)^{\infty}}$$
induces an exact sequence
$$\begin{array}{l}
\dots \to \{X;S(\infty)^+ \wedge Y\}_{\ZZ_2}
\xymatrix{\ar[r]^-{\di{\gamma}}&} \{X;Y\}_{\ZZ_2} \\[1ex]
\hskip100pt \xymatrix{\ar[r]^-{\di{\rho}}&}
\{X;L\R(\infty)^{\infty} \wedge Y\}_{\ZZ_2}=\{X;Y^{\ZZ_2}\}\to \dots
\end{array}$$
\hfill\qed\end{proof}

\begin{example} {\rm Let $Y$ be a pointed $\ZZ_2$-space.\\
(i) The splitting of Proposition \ref{stablesequence1} (ii)
$$\{X;Y\}_{\ZZ_2}~=~\{X;S(\infty)^+\wedge Y\}_{\ZZ_2}\oplus  \{X;Y^{\ZZ_2}\}$$
decomposes the stable $\ZZ_2$-equivariant homotopy theory of $Y$
according to the actual fixed point space $Y^{\ZZ_2}$ and the
$\ZZ_2$-homotopy orbit space $S(\infty)^+\wedge_{\ZZ_2} Y$.
In particular,
$$\widetilde{\omega}^{\ZZ_2}_*(Y)~=~\widetilde{\omega}_*(Y^{\ZZ_2})
\oplus \widetilde{\omega}_*(S(\infty)^+\wedge_{\ZZ_2}Y)~.$$
(ii) If $Y$ is semifree
$$\begin{array}{l}
Y^{\ZZ_2}~=~\{*\}~,~\widetilde{\omega}_*(Y^{\ZZ_2})~=~0~,~
S(\infty)^+\wedge_{\ZZ_2} Y~\simeq~Y/\ZZ_2~,\\[1ex]
\widetilde{\omega}^{\ZZ_2}_*(Y)~=~
\widetilde{\omega}_*(S(\infty)^+\wedge_{\ZZ_2}Y)~=~\widetilde{\omega}_*(Y/\ZZ_2)~.
\end{array}$$
while if the $\ZZ_2$-action on $Y$ is trivial then
$$\begin{array}{l}
Y^{\ZZ_2}~=~Y~,~S(\infty)^+\wedge_{\ZZ_2} Y~=~P(\infty)^+ \wedge Y~,\\[1ex]
\widetilde{\omega}^{\ZZ_2}_*(Y)~=~\widetilde{\omega}_*(Y) \oplus
\widetilde{\omega}_*(P(\infty)^+\wedge Y)~.
\end{array}
$$}
\hfill\qed
\end{example}

\begin{remark} {\rm
The decomposition of the stable $\ZZ_2$-equivariant homotopy groups
of a $\ZZ_2$-space $X$ as a sum of ordinary stable homotopy groups
of the fixed point space $X^{\ZZ_2}$ and the homotopy orbit space
$S(\infty)\times_{\ZZ_2}X$ of the $\ZZ_2$-action
$$\omega^{\ZZ_2}_*(X)~=~\omega_*(X^{\ZZ_2})\oplus \omega_*(S(\infty)\times_{\ZZ_2}X)$$
is the special case $G=\ZZ_2$ of a general decomposition of the
stable $G$-equivariant homotopy groups of a $G$-space
 for an arbitrary compact Lie (e.g.  finite) group $G$
-- see Segal \cite{segal1}, tom Dieck \cite{tomdieck}, Hauschild
\cite{hauschild2}, May et al.  \cite{may}.\\
\hfill\qed}
\end{remark}

\begin{example} \label{stablesequence2} {\rm
(i) The braid of Proposition \ref{stablesequence1} (i) for $U=\R^{i-j}$,
$X=S^j$, $V=\R(\infty)$ can be written

$$\xymatrix@C-35pt{
A_1 \ar[dr] \ar@/^2pc/[rr]^-{\di{0}} && \widetilde{\omega}_j(S(\infty)^+ \wedge_{\ZZ_2} Y)
\ar[dr]\ar@/^2pc/[rr] &&
\{S^j \wedge S(L\R^{i-j})^+;Y\}_{\ZZ_2}\ar[dr]&\\
&~~~~~A_2~~~~~ \ar[ur] \ar[dr] &&
\widetilde{\omega}^{\ZZ_2}_j(Y) \ar[ur] \ar[dr]&&A_3\\
A_4\ar[ur] \ar@/_2pc/[rr]
&&\widetilde{\omega}_{i,j}(Y)
\ar@/_2pc/[rr]^-{\di{\rho}}\ar[ur] && \widetilde{\omega}_j(Y^{\ZZ_2})\ar[ur]&}$$

\bigskip

\noindent with
$$\begin{array}{l}
\widetilde{\omega}_{i,j}(Y)~=~\{S^j \wedge LS^{i-j};Y\}_{\ZZ_2}~,\\[1ex]
A_1~=~\widetilde{\omega}_{j+1}(Y^{\ZZ_2})~,~
A_2~=~\widetilde{\omega}_{i,j}(S(\infty)^+\wedge Y)~,\\[1ex]
A_3~=~\widetilde{\omega}_{i-1,j-1}(S(\infty)^+\wedge Y)~,~
A_4~=~\{\Sigma^{j+1}S(L\R^{i-j})^+;Y\}_{\ZZ_2}~.
\end{array}$$
For $i=j$ $\widetilde{\omega}_{j,j}(Y)=\widetilde{\omega}^{\ZZ_2}_j(Y)$.\\
(ii) The braid of Proposition \ref{stablesequence1} (i) for $U=\R$
can be written

$$\xymatrix@C-35pt{
A_1 \ar[dr] \ar@/^2pc/[rr]^-{\di{\alpha_{LV}}} && \{X;S(LV)^+ \wedge Y\}_{\ZZ_2}
\ar[dr]^-{\di{s_{LV}}}\ar@/^2pc/[rr] &&
\{X;Y\}\ar[dr]&\\
&~~~~~A_2~~~~~ \ar[ur] \ar[dr] &&
\{X;Y\}_{\ZZ_2} \ar[ur]^-{\di{s^*_{L\R}}} \ar[dr]^-{\di{0_{LV}}}&&A_3\\
A_4\ar[ur] \ar@/_2pc/[rr]^-{\di{\alpha^*_{L\R}}}
&&\{LS^1 \wedge X;Y\}_{\ZZ_2}
\ar@/_2pc/[rr]\ar[ur]^-{\di{0^*_{L\R}}} &&
\{X;LV^{\infty}\wedge Y\}_{\ZZ_2}\ar[ur]&}$$

\bigskip

\noindent with
$$\begin{array}{l}
A_1~=~\{\Sigma X;LV^{\infty} \wedge Y\}_{\ZZ_2}~,~
A_2~=~\{\Sigma S(LV\oplus L\R)^+\wedge X;LV^{\infty} \wedge Y\}_{\ZZ_2}~,\\[1ex]
A_3~=~\{S(LV\oplus L\R)^+ \wedge X;LV^{\infty} \wedge Y\}_{\ZZ_2}~,~
A_4~=~\{\Sigma X;Y\}~.
\end{array}$$
The braid includes the exact sequence
$$\xymatrix{\dots \to\{\Sigma X;Y\}_{\ZZ_2}
\ar[r]^-{\di{\alpha_{L\R}^*}}&}
\xymatrix{\{L\R^{\infty}\wedge X;Y\}_{\ZZ_2} \ar[r]^-{\di{0_{L\R}^*}}&
\{X;Y\}_{\ZZ_2} \ar[r]^-{\di{s_{L\R}^*}}& \{X;Y\} \to \dots}$$
with $s_{L\R}^*:\{X;Y\}_{\ZZ_2} \to
 \{S(L\R)^+\wedge X;Y\}_{\ZZ_2}=\{X;Y\}$ the forgetful map.\\
(iii) The special case $X=S^j \wedge LS^{i-j}$, $Y=S^0$, $V=\R(\infty)$
of the braid in (ii) is

$$\xymatrix@C-30pt{
\omega_{j+1} \ar[dr] \ar@/^2pc/[rr]^-{\di{\alpha_{L\R(\infty)}}} &&
\omega_{i,j}(S(L(\R(\infty)))
\ar[dr]^-{\di{s_{L\R(\infty)}}}\ar@/^2pc/[rr] &&
~~~\omega_i~~~\ar[dr]&\\
&\widetilde{\omega}_j(P(\infty,j-i-1)) \ar[ur] \ar[dr] &&
~~\omega_{i,j}~~ \ar[ur]_-{\di{s^*_{L\R}}} \ar[dr]^-{\di{0_{L\R(\infty)}}}&&
\widetilde{\omega}_{j-1}(P(\infty,j-i-1))\\
\omega_{i+1}\ar[ur] \ar@/_2pc/[rr]^-{\di{\alpha^*_{L\R}}}
&&\omega_{i+1,j}\ar@/_2pc/[rr]^-{\di{\rho}}\ar[ur]^-{\di{0^*_{L\R}}} &&
~~~\omega_j~~~\ar[ur]&}$$

\bigskip

\noindent (iv) Proposition \ref{00} gives the bi-degree isomorphism
$$\begin{array}{l}
\omega_{0,0} \cong \ZZ \oplus \ZZ~;\\[1ex]
(F:S^j \wedge LS^{i-j} \to S^j \wedge LS^{i-j}) \mapsto
\big(\dfrac{{\rm degree}(F)-{\rm degree}(G)}{2},{\rm degree}(G))
\end{array}$$
with $G=\rho(F):S^j \to S^j$. If ${\rm degree}(F)=0 \in \ZZ$ then
$${\rm degree}(G) \equiv 0 (\bmod\,2)~.$$
The function
$$\begin{array}{l}
\omega_{1,0} \to \ZZ~;~(H:S^j \wedge LS^{i-j+1} \to S^j \wedge LS^{i-j}) \\[1ex]
\hphantom{\omega_{1,0} \to \ZZ~;~}
\mapsto {\rm degree}(H(1 \wedge 0_{L\R}):S^j \wedge LS^{i-j} \to S^j \wedge LS^{i-j})/2
\end{array}$$
is an isomorphism, and the special case $i=j=0$ of the braid in (iii)
is given by

$$\xymatrix@C-35pt{
\omega_1=\ZZ_2 \ar[dr] \ar@/^2pc/[rr]^-{\di{\alpha_{L\R(\infty)}}} &&
\omega_{0,0}(S(L(\R(\infty)))=\ZZ
\ar[dr]^-{\di{s_{L\R(\infty)}}}\ar@/^2pc/[rr]^-{\di{2}} &&
\omega_0=\ZZ \ar[dr]&\\
&\widetilde{\omega}_0(P(\infty,-1))=0 \ar[ur] \ar[dr] &&
\omega_{0,0}=\ZZ \oplus \ZZ \ar[ur]^-{\di{s^*_{L\R}}} \ar[dr]^-{\di{0_{L\R(\infty)}}}&&
\widetilde{\omega}_{-1}(P(\infty,-1))=\ZZ_2\\
\omega_1=\ZZ_2\ar[ur] \ar@/_2pc/[rr]^-{\di{\alpha^*_{L\R}}}
&&\omega_{1,0}=\ZZ\ar@/_2pc/[rr]^-{\di{\rho=2}}\ar[ur]^-{\di{0^*_{L\R}}} &&
\omega_0=\ZZ \ar[ur]&}$$
\hfill\qed}
\end{example}

\begin{proposition}~
For any pointed $\ZZ_2$-space $X$ the $\ZZ_2$-equivariant homotopy cofibration sequence
$$X^{\ZZ_2} \to X \to X/X^{\ZZ_2}$$
induces long exact sequences in stable $\ZZ_2$-homotopy and cohomotopy
$$\begin{array}{l}
\dots\to \widetilde{\omega}^{\ZZ_2}_i(X^{\ZZ_2})\to \widetilde{\omega}^{\ZZ_2}_i(X)
\to \widetilde{\omega}^{\ZZ_2}_i(X/X^{\ZZ_2}) \to \widetilde{\omega}^{\ZZ_2}_{i-1}(X^{\ZZ_2})\to \dots~,\\[1ex]
\dots\to \widetilde{\omega}_{\ZZ_2}^i(X/X^{\ZZ_2})\to
\widetilde{\omega}_{\ZZ_2}^i(X)
\to \widetilde{\omega}_{\ZZ_2}^i(X^{\ZZ_2}) \to \widetilde{\omega}_{\ZZ_2}^{i+1}(X/X^{\ZZ_2})\to \dots
\end{array}
$$
with
$$\begin{array}{l}
\widetilde{\omega}^{\ZZ_2}_i(X^{\ZZ_2})~=~\widetilde{\omega}_i(X^{\ZZ_2})
\oplus \widetilde{\omega}_i(P(\infty)^+ \wedge X^{\ZZ_2})~,\\[1ex]
\widetilde{\omega}^{\ZZ_2}_i(X)~=~\widetilde{\omega}_i(X^{\ZZ_2})
\oplus \widetilde{\omega}_i(S(\infty)^+ \wedge_{\ZZ_2} X)~,\\[1ex]
\widetilde{\omega}_i^{\ZZ_2}(X/X^{\ZZ_2})~=~\widetilde{\omega}_i((X/X^{\ZZ_2})/\ZZ_2)~,\\[1ex]
\widetilde{\omega}_{\ZZ_2}^i(X^{\ZZ_2})~=~\widetilde{\omega}^i(X^{\ZZ_2})
\oplus \widetilde{\omega}^i(P(\infty)^+ \wedge X_{\ZZ_2})~,\\[1ex]
\widetilde{\omega}_{\ZZ_2}^i(X)~=~\widetilde{\omega}^i(X^{\ZZ_2})
\oplus \widetilde{\omega}^i(S(\infty)^+ \wedge_{\ZZ_2} X)~,\\[1ex]
\widetilde{\omega}^i_{\ZZ_2}(X/X_{\ZZ_2})~=~\widetilde{\omega}^i((X/X^{\ZZ_2})/\ZZ_2)~.
\end{array}$$
\end{proposition}
\begin{proof}
These are $\ZZ_2$-equivariant analogues of the Barratt-Puppe
exact sequence.
Proposition \ref{stablesequence1} (ii)  (= Proposition 3.13 of
Crabb \cite{crabb}) gives the direct sum system
$$\widetilde{\omega}^{\ZZ_2}_i(S(\infty)^+\wedge X)
\xymatrix{\ar@<1ex>[r]^-{\di{\gamma}}\ar@<-1ex>@{<-}[r]_-{\di{\delta}} & }
\widetilde{\omega}^{\ZZ_2}_i(X)
\xymatrix{\ar@<1ex>[r]^-{\di{\rho}}\ar@<-1ex>@{<-}[r]_-{\di{\sigma}}& }
\widetilde{\omega}_i(X^{\ZZ_2})$$
and Proposition \ref{Adams} gives the identification
$$\widetilde{\omega}^{\ZZ_2}_i(S(\infty)^+\wedge X^{\ZZ_2} )~=~
\widetilde{\omega}_i(P(\infty)^+ \wedge X^{\ZZ_2})~.$$
Similarly for $\widetilde{\omega}^i$.\\
\hfill\qed\end{proof}

\begin{proposition}~ \label{stablesequence3}
For any space $X$ and pointed $\ZZ_2$-space $Y$ there are defined
direct sum systems
$$\begin{array}{l}
\{X;S(\infty)^+ \wedge_{\ZZ_2} Y\}
\xymatrix{\ar@<1ex>[r]^-{\di{\gamma}}\ar@<-1ex>@{<-}[r]_-{\di{\delta}} & }
\{X;Y \}_{\ZZ_2}
\xymatrix{\ar@<1ex>[r]^-{\di{\rho}}\ar@<-1ex>@{<-}[r]_-{\di{\sigma}}& }
\{X;Y^{\ZZ_2}\}~,\\[1ex]
\{X;P(\infty)^+ \wedge Y^{\ZZ_2}\}
\xymatrix{\ar@<1ex>[r]^-{\di{\gamma}}\ar@<-1ex>@{<-}[r]_-{\di{\delta}} & }
\{X;Y^{\ZZ_2}\}_{\ZZ_2}
\xymatrix{\ar@<1ex>[r]^-{\di{\rho}}\ar@<-1ex>@{<-}[r]_-{\di{\sigma}}& }
\{X;Y^{\ZZ_2}\}
\end{array}$$
which fit together in a commutative braid of exact sequences

$$\xymatrix@C-40pt{
\{\Sigma X;(Y/Y^{\ZZ_2})/\ZZ_2\}
 \ar[dr] \ar@/^2pc/[rr] && \{X;Y^{\ZZ_2}\}_{\ZZ_2}
\ar[dr]\ar@/^2pc/[rr]^-{\di{\rho}} &&
\{X;Y^{\ZZ_2}\}\\
&\{X;P(\infty)^+\wedge Y^{\ZZ_2}\} \ar[ur] \ar[dr] &&
\{X;Y\}_{\ZZ_2} \ar[ur]^-{\di{\rho}}\ar[dr]&\\
&&\{X;S(\infty)^+\wedge_{\ZZ_2}Y\}
\ar@/_2pc/[rr]\ar[ur] && \{X;(Y/Y^{\ZZ_2})/\ZZ_2\}}$$

\bigskip

\end{proposition}
\begin{proof} The direct sum systems are given by
 Proposition \ref{stablesequence1} (ii).\\
\hfill\qed\end{proof}

\begin{example} {\rm
In the special case $Y=A \wedge A$ for a pointed space $A$
the commutative braid in Proposition \ref{stablesequence3} is given by
$$\xymatrix@C-40pt{
\{\Sigma X;B\}
\ar[dr] \ar@/^2pc/[rr] && \{X;A\}_{\ZZ_2}\ar[dr]^-{\di{\Delta_A}}
\ar@/^2pc/[rr]^-{\di{\rho}} &&
\{X;A\}\\
&\{X;P(\infty)^+\wedge A\} \ar[ur] \ar[dr]^-{\di{1\wedge \Delta_A}} &&
\{X;A\wedge A\}_{\ZZ_2} \ar[ur]^-{\di{\rho}} \ar[dr]&\\
&&\{X;S(\infty)^+\wedge_{\ZZ_2}(A \wedge A)\}
\ar@/_2pc/[rr]\ar[ur] && \{X;B\}}$$

\medskip

\noindent with $B=((A \wedge A)/\Delta_A)/\ZZ_2$ and
$$\begin{array}{l}
\{X;A\}_{\ZZ_2}~=~\{X;A\}\oplus \{X;P(\infty)^+\wedge A\}~,\\[1ex]
\{X;A\wedge A\}_{\ZZ_2}~=~\{X;A\}\oplus \{X;S(\infty)^+\wedge_{\ZZ_2}(A\wedge A)\}~.
\end{array}$$
\hfill\qed}
\end{example}

\section{$\ZZ_2$-equivariant bundles}\label{Z2bundles}

We shall be mainly concerned with $\ZZ_2$-equivariant $U$-bundles for
finite-dimensional inner product spaces $U$.

\begin{definition}~{\rm
{\rm (i)} Let $U$ be an inner product space, and $X$ a $\ZZ_2$-space.
A {\it $\ZZ_2$-equivariant $U$-bundle $\xi$ over $X$} is a $U$-bundle
$$\xymatrix{\xi~:~U \ar[r] & E(\xi) \ar[r]^-{\di{p}}&X}$$
such that $E(\xi)$ is a $\ZZ_2$-space, $p$ is $\ZZ_2$-equivariant,
and the restrictions of $T:E(\xi)\to E(\xi)$ to the fibres
$$T\vert~:~\xi_x~=~p^{-1}(x) \to \xi_{Tx}~=~p^{-1}(Tx)~=~T\xi_x~(x \in X)$$
are isometries, with
$${\rm dim}(\xi_x)~=~{\rm dim}(\xi_{T(x)})~=~{\rm dim}(U)$$
and the Thom space $T(\xi)$ is a pointed $\ZZ_2$-space.
For $x \in X^{\ZZ_2}$ the fibre $\xi_x$ is an inner product $\ZZ_2$-space, and
for finite-dimensional $U$ the dimension of the fixed point subspace is
a continuous function
$${\rm dim}\,\xi_+~:~X^{\ZZ_2} \to \N ~;~x \mapsto {\rm dim}(\xi_x)_+$$
which is constant on each component of $X^{\ZZ_2}$. \\
{\rm (ii)} Let $V$ be an inner product $\ZZ_2$-space.
The {\it trivial $\ZZ_2$-equivariant $\vert V\vert$-bundle} $\epsilon_V$ over
a $\ZZ_2$-space $X$ is defined by
$$\epsilon_V~:~\vert V\vert \to E(\epsilon_V)~=~X \times V \to X~,$$
with
$$\begin{array}{l}
X \times V \to X~;~(x,v) \mapsto x~,~(\epsilon_V)_x~=~V~(x \in X)~,\\[1ex]
T~:~E(\epsilon_V)~=~X \times V \to X \times V~;~(x,v) \mapsto (Tx,Tv)~.
\end{array}$$
{\rm (iii)} The {\it Thom $\ZZ_2$-spectrum} of a $\ZZ_2$-equivariant
$U$-bundle $\xi$ is the suspension $\ZZ_2$-spectrum of $T(\xi)$
$$\underline{T}(\xi)~=~\{T(\xi \oplus \epsilon_V)\,\vert\,V\}$$
with $V$ running over finite-dimensional inner product $\ZZ_2$-spaces,
and
$$T(\xi \oplus \epsilon_V)~=~V^{\infty} \wedge T(\xi)~.$$
For the trivial $\ZZ_2$-equivariant $V$-bundle $\epsilon_V$
over $X$
$$\begin{array}{l}
E(\epsilon_V)^{\ZZ_2}~=~V_+ \times X^{\ZZ_2}~,~
T(\epsilon_V)~=~V^{\infty} \wedge X^{\infty}~,\\[1ex]
T(\epsilon_V)^{\ZZ_2}~=~(V_+)^{\infty} \wedge (X^{\ZZ_2})^{\infty}~.
\end{array}$$}
\hfill\qed
\end{definition}

\begin{example}~{\rm Let $U$ be an inner product space.\\
(i) A $\ZZ_2$-equivariant $U$-bundle $\xi$ over $\{*\}$ is a
$\ZZ_2$-action on $E(\xi)=U$, so that
$$\xi~=~\epsilon_U~=~\epsilon_{U_+} \oplus \epsilon_{U_-}~.$$
(ii) For a $\ZZ_2$-equivariant $U$-bundle $\xi$ over a $\ZZ_2$-space $X$,
forgetting the $\ZZ_2$-structure determines a  $U$-bundle $\xi$ over $X$.\\
(iii) A $U$-bundle $\xi$ over a space $X$ can be
regarded as  $\ZZ_2$-equivariant $U$-bundle $\xi$ over $X$ (with the
trivial $\ZZ_2$-action).\\
(iv) A  $U$-bundle $\xi$ over a space $X$
determines a $\ZZ_2$-equivariant $U$-bundle $L\xi$ over $X$ with
$$E(L\xi)~=~E(\xi)~,~T~:~E(L\xi) \to E(L\xi)~;~(x,u) \mapsto (x,-u)~,$$
such that
$$(L\xi)_+~=~0~,~(L\xi)_-~=~\xi~.$$
In particular, for the trivial $U$-bundle $\xi=\epsilon_U$ the
construction gives the $\ZZ_2$-equivariant $U$-bundle $L\xi=\epsilon_{LU}$.\\
\hfill\qed}
\end{example}

\begin{proposition}~ Let $\xi$ be a $\ZZ_2$-equivariant $\R^i$-bundle
over a $\ZZ_2$-space $X$. \\
{\rm (i)} Each fibre $\xi_x$  $(x \in X^{\ZZ_2})$
is an inner product $\ZZ_2$-space, and
$$E(\xi)^{\ZZ_2}~=~\bigcup\limits_{x \in X^{\ZZ_2}}(\xi_x)^{\ZZ_2}~.$$
The subspaces
$$X(\xi,k)~=~\{x \in X^{\ZZ_2}\,\vert\,
{\rm dim}((\xi_x)_+)=k\}\subseteq X^{\ZZ_2}~~(0 \leqslant k \leqslant i)$$
are unions of components of $X^{\ZZ_2}$ such that
$$X^{\ZZ_2}~=~\coprod\limits^i_{k=0}X(\xi,k)~.$$
The restriction $\ZZ_2$-equivariant bundle $(\xi,k)=\xi\vert_{X(\xi,k)}$
over $X(\xi,k)$ splits as a Whitney sum
$$(\xi,k)~=~(\xi,k)_+ \oplus L(\xi,k)_-$$
with
$$(\xi,k)_+~:~X(\xi,k) \to BO(\R^k)~,~(\xi,k)_-~:~X(\xi,k) \to BO(\R^{i-k})$$
nonequivariant bundles such that
$$\begin{array}{l}
E((\xi,k)_+)~=~\bigcup\limits_{x \in X(\xi,k)}(\xi_x)_+~,~
E((\xi,k)_-)~=~\bigcup\limits_{x \in X(\xi,k)}(\xi_x)_-~,\\[2ex]
E(\xi)^{\ZZ_2}~=~\coprod\limits^i_{k=0}E((\xi,k)_+)~,~
T(\xi)^{\ZZ_2}~=~\bigvee\limits^i_{k=0}T((\xi,k)_+)~.
\end{array}$$
{\rm (ii)} For a trivial $\ZZ_2$-equivariant bundle $\epsilon_V$
over a $\ZZ_2$-space $X$
$$X(\epsilon_V,k)~=~\begin{cases} X^{\ZZ_2}&{\rm if}~k={\rm dim}(V_+)\\
\emptyset&{\rm if}~k \neq {\rm dim}(V_+)~.
\end{cases}$$
{\rm (iii)} For the Whitney sum  $\xi \oplus \eta$ of
$\ZZ_2$-equivariant bundles $\xi,\eta$ over a $\ZZ_2$-space $X$
$$X(\xi \oplus \eta,k)~=~\coprod\limits_{i+j=k}X(\xi,i) \cap X(\eta,j)~.$$
In particular, if $\eta=\epsilon_V$
$$X(\xi \oplus \epsilon_V,k)~=~\begin{cases} X(\xi,k-{\rm dim}(V_+))&
{\rm if}~k\geqslant {\rm dim}(V_+)\\
\emptyset&{\rm if}~k < {\rm dim}(V_+)~.
\end{cases}$$
{\rm (iv)} If $\xi \oplus \eta=\epsilon_V$ for some inner product $\ZZ_2$-space $V$,
then
$$(\xi_x)_+ \oplus (\eta_x)_+~=~V_+~~(x \in X^{\ZZ_2})~,$$
so that
$$\begin{array}{ll}
X(\xi,k)&=~\{x \in X^{\ZZ_2}\,\vert\, {\rm dim}((\xi_x)_+)=k\}\\[1ex]
&=~\{x \in X^{\ZZ_2}\,\vert\, {\rm dim}((\eta_x)_+)={\rm dim}(V_+)-k)\}\\[1ex]
&=~X(\eta,{\rm dim}(V_+)-k)~.
\end{array}$$
\end{proposition}
\begin{proof} By construction.\\
\hfill\qed\end{proof}

\begin{proposition}
Let $U$ be a finite-dimensional inner product space,
and let $X$ be a $\ZZ_2$-space, with $p:X \to X/\ZZ_2$ the projection.\\
{\rm (i)} For any $U$-bundle $\xi$ over $X/\ZZ_2$ the pullback $p^*\xi$
is a $\ZZ_2$-equivariant $U$-bundle over $X$, with the trivial $\ZZ_2$-action
on each fibre $p^*\xi_x=\xi_{p(x)}$ ($x \in X$).\\
{\rm (ii)} If $X$ is a free $\ZZ_2$-space the
projection $p:X \to X/\ZZ_2$ is a double covering, and the function
$$\{U\hbox{-bundles over}~X/\ZZ_2\} \to
\{\ZZ_2\hbox{-equivariant}~U\hbox{-bundles over}~X\}~;~
\xi \mapsto p^*\xi$$
is a bijection of the sets of isomorphism classes, with
$$\begin{array}{l}
E(p^*\xi)^{\ZZ_2}~=~\emptyset~,~T(p^*\xi)^{\ZZ_2}~=~\{\infty\}~,\\[1ex]
E(p^*\xi)/\ZZ_2~=~E(\xi)~,~T(p^*\xi)/\ZZ_2~=~T(\xi)~.
\end{array}
$$
{\rm (iii)} For any $\ZZ_2$-equivariant $U$-bundle $\xi$
over a $\ZZ_2$-space $X$ the function
$$\begin{array}{l}
\{(\hbox{$\ZZ_2$-equivariant map $f:M \to X$ with $M$ $\ZZ_2$-free,}\\[1ex].
\hskip100pt
\hbox{$\ZZ_2$-equivariant bundle map $b:\nu_M \to \xi$})\}\\[1ex]
\to \{(\hbox{$g:M/\ZZ_2 \to S(\infty)\times_{\ZZ_2} X$,
bundle map $c:\nu/\ZZ_2 \to 0\times_{\ZZ_2}\xi$})\}
\end{array}$$
is a bijection of homotopy classes.
\end{proposition}
\begin{proof} (i) By construction.\\
(ii)  For any $\ZZ_2$-equivariant $U$-bundle $\Xi$ over $X$ the $\ZZ_2$-action
on the total $\ZZ_2$-space $E(\Xi)$ is free, and the quotient is the
total space
$$E(\Xi)/\ZZ_2~=~E(\xi)$$
of a $U$-bundle $\xi$ over $X/\ZZ_2$ such that $\Xi=p^*\xi$.\\
(iii) For a $\ZZ_2$-free space $M$ there is a $\ZZ_2$-equivariant
map $a:M \to S(\infty)$, so $(f:M \to X,b:\nu_M \to \xi)$
induces $(g:M/\ZZ_2 \to S(\infty)\times_{\ZZ_2}X,c:\nu_M/\ZZ_2 \to 0\times_{\ZZ_2}\xi)$ with $g([x])=[a(x),f(x)]$ ($x \in M$), $c=[b]$. Conversely, given $(g,c)$ let $(f,b)=(\bar{g},\bar{c})$
be the double cover.\\
\hfill\qed\end{proof}

\begin{example}~ {\rm Let $V$ be an inner product $\ZZ_2$-space.
The trivial $\ZZ_2$-equivariant $\vert V
\vert$-bundle $\epsilon_V$ over a free $\ZZ_2$-space $X$ is the pullback
$\epsilon_V=p^*(\epsilon_V/\ZZ_2)$ along the projection $p:X \to X/\ZZ_2$
of the $\vert V\vert$-bundle $\epsilon_V/\ZZ_2$ over $X/\ZZ_2$ with
$$V \to E(\epsilon_V/\ZZ_2)~=~E(\epsilon_V)/\ZZ_2~=~X\times_{\ZZ_2}V
\to X/\ZZ_2~.$$
\hfill\qed}
\end{example}

\begin{proposition}~
{\rm (i)} An $m$-dimensional $\ZZ_2$-manifold $M$ has a tangent
$\ZZ_2$-equivariant $\R^m$-bundle $\tau_M$,
namely the nonequivariant tangent $\R^m$-bundle
$\tau_M$ with the $\ZZ_2$-action $T:E(\tau_M) \to E(\tau_M)$ defined by
the differentials of the $\ZZ_2$-action $T:M \to M$.\\
{\rm (ii)} The tangent $\R^m$-bundle of an $m$-dimensional
free $\ZZ_2$-manifold $M$ is the pullback
$\tau_M=p^*\tau_{M/\ZZ_2}$ along the projection $p:M \to M/\ZZ_2$ of
$\tau_{M/\ZZ_2}$, with
$$E(\tau_{M/\ZZ_2})~=~E(\tau_M)/\ZZ_2~,~\tau_{M/\ZZ_2}~=~\tau_M/\ZZ_2~.$$
\hfill\qed
\end{proposition}

\begin{example} {\rm Let $U$ be a finite-dimensional inner product space.\\
(i) The $\ZZ_2$-manifold $LU$ has tangent $\ZZ_2$-equivariant $U$-bundle
$$\tau_{LU}~=~\epsilon_{LU}~.$$
(ii) The free $\ZZ_2$-manifold $S(LU) \subset LU$ is such that
there is defined a $\ZZ_2$-equivariant homeomorphism
$$S(LU) \times \R \to LU \backslash \{0\}~;~(u,x) \mapsto e^xu~,$$
so that the tangent $\ZZ_2$-equivariant bundle $\tau_{S(LU)}$ is  such that
$$\tau_{S(LU)} \oplus \epsilon_\R~=~\epsilon_{LU}~.$$
(iii) The tangent bundle of the projective space $P(U)=S(LU)/\ZZ_2$
$$\tau_{P(U)}~=~\tau_{S(LU)}/\ZZ_2$$
As usual, it is possible to identify
$$\begin{array}{l}
P(U)~=~G(\R,U)~=~\{K \subseteq U\,\vert\,{\rm dim}(K)=1\}~,\\[1ex]
E(\tau_{S(LU)})~=~\{(x\in LU,y \in S(LU))\,\vert\, \langle x,y \rangle =0\}~,\\[1ex]
E(\tau_{P(U)})~=~E(\tau_{S(LU)})/\ZZ_2~=~
\{(K,x\in K^{\perp})\,\vert\, K \subseteq U,\,{\rm dim}(K)=1\}~.
\end{array}$$
\hfill\qed}
\end{example}

\begin{definition}~ \label{Hopfbundle}
{\rm Let $U,V$ be finite-dimensional inner product spaces.
The $\ZZ_2$-equivariant embedding $S(LU) \subset S(LU \oplus LV)$
has normal $\ZZ_2$-equivariant $V$-bundle
$$\nu_{S(LU) \subset S(LU \oplus LV)}~=~\epsilon_{LV}~.$$
(i) The {\it Hopf $V$-bundle} over the projective space $P(U)$
$$H_V~=~\epsilon_{LV}/\ZZ_2~:~V \to E(H_V)~=~S(LU)\times_{\ZZ_2}LV \to P(U)$$
is the normal $V$-bundle of the embedding $P(U) \subset P(U \oplus V)$
$$H_V~=~\nu_{P(U) \subset P(U \oplus V)}~:~P(U) \to BO(V)~.$$
Also, for $V=U$
$$H_U~=~\tau_{P(U)} \oplus \epsilon_\R~=~\tau_{P(U \oplus \R)}\vert~
:~P(U) \to BO(U)~.$$
(ii) The Thom space of $H_V$ is the {\it stunted projective space}
$$T(H_V)~=~S(LU)^+ \wedge_{\ZZ_2} LV^{\infty}~.$$
For $U=\R^m$, $V=\R^n$ write
$$P(\R^{m+n})/P(\R^n)~=~P(m+n,n)~=~S(L\R^m)^+ \wedge_{\ZZ_2} LS^n~,$$
with $P(m,0)=P(\R^m)^+$.\\
(iii) The {\it infinite stunted projective space} $P(\infty,n)$ is defined
for any $n \in \ZZ$ to be the virtual Thom space (\ref{virtual})
of the virtual bundle ${\rm sgn}(n)\,H_{\R^{\vert n\vert}}$
over $P(\R^m)$ for sufficiently large $m \geqslant 0$.
For $n \geqslant 0$ this is the actual Thom space
$$\begin{array}{ll}
P(\infty,n)~=~T(H_{\R^n})~=~P(\R(\infty)\oplus \R^n)/P(\R^n)
&\,=~S(\infty)^+ \wedge_{\ZZ_2} LS^n\\[1ex]
&(=~P(\infty)^+~\hbox{for}~n=0)~,
\end{array}$$
so that
$$\widetilde{\omega}_j(P(\infty,n))~=~\{S^j;S(\infty)^+\wedge LS^n\}_{\ZZ_2}~.$$
For $n<0$ this is the `space' $P(\infty,n)$ with
$$\widetilde{\omega}_j(P(\infty,n))~=~\{S^j \wedge LS^{-n};S(\infty)^+\}_{\ZZ_2}~.$$
\hfill\qed}
\end{definition}

\begin{example}~{\rm
The Hopf $\R$-bundle $H_{\R}$ over $P(U)$ is such that
$$\begin{array}{l}
E(H_{\R})~=~\{(K,x \in K)\,\vert\, K \subseteq U,\,{\rm dim}(K)=1\}~=~
S(LU) \times_{\ZZ_2}L\R~,\\[1ex]
S(H_{\R})~=~\{(K,x \in K)\,\vert\, K \subseteq U,\,{\rm dim}(K)=1, \Vert x \Vert =1\}~=~S(LU)~,\\[1ex]
H_{\R}~=~\nu_{P(U) \subset P(U \oplus \R)}~:~
P(U) \to \varinjlim\limits_pP(U\oplus \R^p)~=~BO(\R)~~,\\[1ex]
\tau_{P(U)}\oplus \epsilon_{\R}~=~\tau_{P(U \oplus \R)}\vert~=~H_U~:~
P(U) \to BO(U)~,\\[1ex]
T(H_{\R})~=~P(U \oplus \R)/P(\R)~=~S(LU)^+\wedge_{\ZZ_2}L\R^{\infty}~,$$
\end{array}$$
and for any $n \geqslant 1$
$$\begin{array}{l}
\nu_{P(U) \subset P(U \oplus \R^n)}~=~H_{\R^n}~=~\bigoplus\limits_n H_{\R}~:~
P(U) \to BO(\R^n)~,\\[1ex]
\tau_{P(U)} \oplus \epsilon_{\R} \oplus H_{\R^{n-1}}~=~
\tau_{P(U \oplus \R^n)}\vert~=~H_{U\oplus \R^{n-1}}~:~P(U) \to BO(U \oplus \R^{n-1})~.
\end{array}$$
\hfill\qed}
\end{example}

\begin{definition}\label{Z2omega}~
{\rm Let $\xi-\xi'$ be a virtual $\ZZ_2$-equivariant bundle over
a $\ZZ_2$-space $X$, and let $\eta$ be a $\ZZ_2$-equivariant bundle over $X$
such that $\xi'\oplus \eta=\epsilon_V$ for some finite-dimensional
inner product $\ZZ_2$-space $V$, so that $\xi-\xi'=\xi\oplus \eta - \epsilon_V$.\\
(i) The {\it virtual Thom $\ZZ_2$-space} of $\xi-\xi'$ is the $\ZZ_2$-spectrum
$$\underline{T}(\xi-\xi')~=~\{T(\xi-\xi')_U\,\vert\,U\}$$
defined by
$$T(\xi-\xi')_U~=~T(\xi \oplus \eta\oplus \epsilon_{U/V})$$
for finite-dimensional inner product $\ZZ_2$-spaces $U$ such that
$V \subseteq U$. \\
(ii) Write
$$\begin{array}{l}
\omega^*_{\ZZ_2}(X;\xi-\xi')~=~\widetilde{\omega}^*_{\ZZ_2}(\underline{T}(\xi-\xi'))~,\\[1ex]
\omega_*^{\ZZ_2}(X;\xi-\xi')~=~\widetilde{\omega}_*^{\ZZ_2}(\underline{T}(\xi-\xi'))
\end{array}$$
so that
$$\begin{array}{l}
\omega_{\ZZ_2}^m(X;\xi-\xi')~=~\{T(\xi\oplus \eta);\Sigma^mV^{\infty}\}_{\ZZ_2}~,\\[1ex]
\omega^{\ZZ_2}_m(X;\xi-\xi')~=~\{\Sigma^mV^{\infty};T(\xi\oplus \eta)\}_{\ZZ_2}~.
\end{array}$$
\hfill\qed}
\end{definition}

\begin{proposition}\label{Z2free}~
{\rm (i)} Let $X$ be a pointed $\ZZ_2$-space
and let $p:X \to X/\ZZ_2$ be the projection, so that
there are induced morphisms
$$\begin{array}{l}
p^*~:~\widetilde{\omega}^m(X/\ZZ_2) \to \widetilde{\omega}^m_{\ZZ_2}(X)~;\\[1ex]
\hskip25pt
(F:V^{\infty} \wedge X/\ZZ_2 \to \Sigma^mV^{\infty})
\mapsto (F(1 \wedge p):V^{\infty} \wedge X \to \Sigma^mV^{\infty})~,\\[1ex]
p_*~:~\widetilde{\omega}_m^{\ZZ_2}(X)\to \widetilde{\omega}_m(X/\ZZ_2) ~;\\[1ex]
\hskip25pt
(F:\Sigma^mV^{\infty} \to V^{\infty}\wedge X)
\mapsto ((1 \wedge p)F:\Sigma^mV^{\infty} \to V^{\infty}\wedge X/\ZZ_2)~.
\end{array}$$
If $X$ is semifree then $p^*$ is an isomorphism,
and if in addition $X$ is a finite $CW$ $\ZZ_2$-complex then $p_*$
is also an isomorphism.\\
{\rm (ii)} For a virtual $\ZZ_2$-equivariant bundle $\xi-\xi'$
over a free $\ZZ_2$-space $X$
$$\omega^*_{\ZZ_2}(X;\xi-\xi')~=~\omega^*(X/\ZZ_2;\xi/\ZZ_2-\xi'/\ZZ_2)$$
and if in addition $X$ is a finite $CW$ $\ZZ_2$-complex then
$$\omega_*^{\ZZ_2}(X;\xi-\xi')~=~\omega_*(X/\ZZ_2;\xi/\ZZ_2-\xi'/\ZZ_2)~.$$
\end{proposition}
\begin{proof} (i) For any finite-dimensional inner product $\ZZ_2$-space $V$
there is defined a $V$-bundle $\epsilon_V/\ZZ_2$ over $X/\ZZ_2$ with
$$E(\epsilon_V/\ZZ_2)~=~V \times_{\ZZ_2}X~.$$
For sufficiently large $n \geqslant 0$ there exists a bundle $\lambda_V$
over $X/\ZZ_2$ such that
$$(\epsilon_V/\ZZ_2)\oplus \lambda_V~=~\epsilon_{\R^n}~.$$
An element $G \in \omega^m_{\ZZ_2}(X)$ is represented by a
$\ZZ_2$-equivariant map $G:V^{\infty} \wedge X \to \Sigma^mV^{\infty}$,
and the composite
$$\begin{array}{l}
H~:~(V^{\infty} \wedge_{\ZZ_2} X)\wedge_{X/\ZZ_2}T(\lambda_V)~=~
\Sigma^nX^+ \xymatrix{\ar[r]^-{\di{G \wedge 1}}&}\\[1ex]
\hskip100pt  \Sigma^m(V^{\infty} \wedge_{\ZZ_2} X)\wedge_{X/\ZZ_2}T(\lambda_V)~=~
\Sigma^{m+n}X \to S^{m+n}
\end{array}$$
is a stable map such that
$$\omega^m_{\ZZ_2}(X) \to \omega^m(X/\ZZ_2)~;~G \mapsto H$$
is an inverse to $p^*$.\\
For $X$ a finite $CW$ $\ZZ_2$-complex the stable $\ZZ_2$-equivariant
Umkehr map $F \in \{X/\ZZ_2;X\}_{\ZZ_2}$ of Proposition \ref{Z2cover} (ii)
induces a morphism
$$\omega_m(X/\ZZ_2) \to \omega_m^{\ZZ_2}(X)~;~G \mapsto FG$$
inverse to $p_*$, by \ref{Z2cover} (iii).\\
(ii) Let $\eta$ be a bundle over $X$
such that $\xi'/\ZZ_2\oplus \eta=\epsilon_V$ for some finite-dimensional
inner product space $V$, so that
$\xi' \oplus p^*\eta=\epsilon_V$,
$$\begin{array}{l}
\omega_{\ZZ_2}^m(X;\xi-\xi')~=~\{T(\xi\oplus p^*\eta);\Sigma^mV^{\infty}\}_{\ZZ_2}~,\\[1ex]
\omega^m(X/\ZZ_2;\xi/\ZZ_2-\xi'/\ZZ_2)~=~\{T(\xi/\ZZ_2\oplus \eta);\Sigma^mV^{\infty}\}~,\\[1ex]
\omega^{\ZZ_2}_m(X;\xi-\xi')~=~\{\Sigma^mV^{\infty};T(\xi\oplus \eta)\}_{\ZZ_2}~,\\[1ex]
\omega_m(X/\ZZ_2;\xi/\ZZ_2-\xi'/\ZZ_2)~=~\{\Sigma^mV^{\infty};T(\xi/\ZZ_2\oplus \eta)\}~.
\end{array}$$
The $\ZZ_2$-action on $T(\xi\oplus p^*\eta)$ is semifree,
with quotient
$$T(\xi\oplus p^*\eta)/\ZZ_2~=~T(\xi/\ZZ_2\oplus \eta)~,$$
so that by (i)
$$\begin{array}{ll}
\omega_{\ZZ_2}^m(X;\xi-\xi')
&=~\{T(\xi\oplus p^*\eta);\Sigma^mV^{\infty}\}_{\ZZ_2}\\[1ex]
&=~\{T(\xi/\ZZ_2\oplus \eta);\Sigma^mV^{\infty}\}~=~\omega^m(X/\ZZ_2;\xi/\ZZ_2-\xi'/\ZZ_2)~,\\[1ex]
\omega^{\ZZ_2}_m(X;\xi-\xi')
&=~\{\Sigma^mV^{\infty};T(\xi\oplus p^*\eta)\}_{\ZZ_2}\\[1ex]
&=~\{\Sigma^mV^{\infty};T(\xi/\ZZ_2\oplus \eta)\}~=~
\omega_m(X/\ZZ_2;\xi/\ZZ_2-\xi'/\ZZ_2)~.
\end{array}$$
\hfill\qed\end{proof}

\begin{terminology}~ {\rm For $i \geqslant j \geqslant 0$
the $i$-dimensional inner product $\ZZ_2$-space $V=\R^j\oplus L\R^{i-j}$
has a $j$-dimensional fixed point space $V_+=\R^j$. For any
$\ZZ_2$-space $X$ write
$$\begin{array}{l}
\omega^{i,j}(X)~=~\omega_{\ZZ_2}^0(X;-\epsilon_V)~=~\{X^+;V^{\infty}\}_{\ZZ_2}~,\\[1ex]
\omega_{i,j}(X)~=~\omega^{\ZZ_2}_0(X;-\epsilon_V)~=~\{V^{\infty};X^+\}_{\ZZ_2}
\end{array}$$
exactly as in Crabb \cite[pp. 28,29]{crabb}.}\\
\hfill\qed
\end{terminology}

\begin{example} {\rm If $X$ is a free $\ZZ_2$-space and $V=\R^j\oplus L\R^{i-j}$
then
$$\begin{array}{l}
\omega^{i,j}(X)~=~\omega_{\ZZ_2}^0(X;-\epsilon_V)~=~
\omega^0(X/\ZZ_2;-\epsilon_V/2)~,\\[1ex]
\omega_{i,j}(X)~=~\omega^{\ZZ_2}_0(X;-\epsilon_V)~=~
\omega_0(X/\ZZ_2;-\epsilon_V/2)
\end{array}$$
by Proposition \ref{Z2free}, assuming $X$ is a finite $CW$ $\ZZ_2$-complex
for $\omega_{i,j}(X)$.\\
\hfill\qed}
\end{example}

\begin{example} {\rm The $\ZZ_2$-equivariant homotopy cofibration sequence
$$S(\infty)^+ \to S^0 \to L\R(\infty)^{\infty} \to \Sigma S(\infty)^+ \to \dots $$
induces the exact sequence
$$\xymatrix@C-5pt{\dots \ar[r]&
\widetilde{\omega}_j(P(\infty,j-i))\ar[r] & \omega_{i,j} \ar[r]^-{\di{\rho}} &\omega_j \ar[r] &
\widetilde{\omega}_{j-1}(P(\infty,j-i)) \ar[r] &\dots}$$
of Crabb \cite[Proposition 4.6]{crabb}, with $\omega_{i,j}=\omega_{i,j}(\{*\})$.\\
\hfill\qed}
\end{example}

\begin{example}~{\rm (i) For $i=j$
$$\begin{array}{l}
\omega^{i,i}(X)~=~\omega_{\ZZ_2}^i(X)~=~\{X^+;S^i\}_{\ZZ_2}~,\\[1ex]
\omega_{i,i}(X)~=~\omega^{\ZZ_2}_i(X)~=~\{S^i;X^+\}_{\ZZ_2}~.
\end{array}$$
(ii) For $V=\R^j\oplus L\R^{i-j}$, $V'=\R^{j'}\oplus L\R^{i'-j'}$,
$$\begin{array}{l}
\omega^n_{\ZZ_2}(X;\epsilon_V-\epsilon_{V'})~=~
\{V^{\infty} \wedge X^+;(V' \oplus \R^n)^{\infty}\}_{\ZZ_2}~=~\omega^{n+i'-i,n+j'-j}(X)~,\\[1ex]
\omega_n^{\ZZ_2}(X;\epsilon_V-\epsilon_{V'})~=~
\{(V' \oplus \R^n)^{\infty};V^{\infty} \wedge X^+\}_{\ZZ_2}~=~\omega_{n+i'-i,n+j'-j}(X)~.
\end{array}$$
\hfill\qed}
\end{example}

\begin{proposition}\label{bundlesequence}~
Let $\xi$ be a $\ZZ_2$-equivariant $\R^i$-bundle over a $\ZZ_2$-space $X$,
so that $0 \times \xi$ is a $\ZZ_2$-equivariant $\R^i$-bundle over the
free $\ZZ_2$-space $S(\infty) \times X$ with
$$E(0 \times \xi)~=~S(\infty) \times E(\xi)~,~T(0 \times \xi)~=~S(\infty)^+ \wedge T(\xi)~,$$
and $0 \times_{\ZZ_2} \xi$ is an $\R^i$-bundle over $S(\infty) \times_{\ZZ_2} X$ with
$$E(0 \times_{\ZZ_2} \xi)~=~S(\infty) \times_{\ZZ_2} E(\xi)~,~T(0 \times_{\ZZ_2} \xi)~=~S(\infty)^+ \wedge_{\ZZ_2}T(\xi)~.$$
{\rm (i)} For any inner product $\ZZ_2$-space $V$
the stable $\ZZ_2$-equivariant homotopy groups
$\omega^{\ZZ_2}_*(X;\xi-\epsilon_V)$
fit into a long exact sequence
$$\begin{array}{l}
\dots \to \omega^{\ZZ_2}_m(S(\infty)\times X;0\times(\xi-\epsilon_V))
\to \omega^{\ZZ_2}_m(X;\xi-\epsilon_V)\\[1ex]
\xymatrix{\ar[r]^-{\di{\rho}}&} \omega_m(\underline{T}(\xi-\epsilon_V)^{\ZZ_2})
\to \omega^{\ZZ_2}_{m-1}(S(\infty)\times X;0\times(\xi-\epsilon_V))
\to \dots
\end{array}$$
with $\rho$ the fixed point map, and
$$\begin{array}{l}
\omega^{\ZZ_2}_m(X;\xi-\epsilon_V)~=~\{(V \oplus \R^m)^{\infty};T(\xi)\}_{\ZZ_2}~,\\[1ex]
\omega_m(\underline{T}(\xi-\epsilon_V)^{\ZZ_2})~=~
\{(V_+\oplus \R^m)^{\infty};T(\xi)^{\ZZ_2}\}\\[1ex]
\hphantom{\omega_m(\underline{T}(\xi-\epsilon_V)^{\ZZ_2})~}
=~\bigoplus\limits^i_{k=0}\omega_{m+{\rm dim}(V_+)}(X(\xi,k);(\xi,k)_+)~.
\end{array}$$
Furthermore, if $X$ is a finite $CW$ $\ZZ_2$-complex
$$\omega^{\ZZ_2}_m(S(\infty)\times X;0\times(\xi-\epsilon_V))~=~
\omega_m(S(\infty)\times_{\ZZ_2}X;0\times_{\ZZ_2}(\xi-\epsilon_V))~.$$
{\rm (ii)} For $V=\{0\}$ the sequence in {\rm (i)} breaks up into split short exact
sequences: the stable $\ZZ_2$-equivariant homotopy groups
$\omega^{\ZZ_2}_*(X;\xi)$ split as
$$\begin{array}{ll}
\omega^{\ZZ_2}_m(X;\xi)&=~\widetilde{\omega}^{\ZZ_2}_m(T(\xi))\\[1ex]
&=~\widetilde{\omega}_m(S(\infty)^+ \wedge_{\ZZ_2} T(\xi))\oplus \widetilde{\omega}_m(T(\xi)^{\ZZ_2})\\[1ex]
&=~\omega_m(S(\infty) \times_{\ZZ_2}X;0\times_{\ZZ_2}\xi)
\oplus \bigoplus\limits^i_{k=0}\omega_m(X(\xi,k);(\xi,k)_+)~.
\end{array}$$
\end{proposition}
\begin{proof} These are special cases of Proposition \ref{stablefixed} (ii)+(iii),
combined with Proposition \ref{Z2free}.\\
\hfill\qed\end{proof}

The stable cohomotopy Thom and Euler classes (\ref{difcon3}) have
$\ZZ_2$-equivariant versions:

\begin{definition}~
{\rm
(i) The {\it stable $\ZZ_2$-equivariant cohomotopy Thom class} of
a $\ZZ_2$-equivariant $V$-bundle $\xi$ over $X$
$$u^{\ZZ_2}(\xi)~=~1 \in \omega^0_{\ZZ_2}(D(\xi),S(\xi);-\xi)~=~\omega^0_{\ZZ_2}(X)$$
is represented by the pointed $\ZZ_2$-map $1:X^+ \to S^0$
sending $X$ to the non-base point. \\
(ii) The {\it stable $\ZZ_2$-equivariant cohomotopy Euler class} of
a $\ZZ_2$-equivariant $V$-bundle $\xi$ over $X$
$$\gamma^{\ZZ_2}(\xi)~=~z^*u^{\ZZ_2}(\xi)
 \in \omega^0_{\ZZ_2}(X;-\xi)~=~\{T(\eta);(V\oplus W)^{\infty}\}_{\ZZ_2}$$
is represented by the $\ZZ_2$-equivariant map
$$\xymatrix{T(\eta) \ar[r]^-{\di{z}} & T(\xi \oplus \eta)~=~
(V \oplus W)^{\infty} \wedge X^+ \ar[r] & (V\oplus W)^{\infty}}$$
with $\xi \oplus \eta=\epsilon_{V \oplus W}$ and
$$z~:~E(\eta) \to E(\xi \oplus \eta)~;~(w,x) \mapsto ((0,w),x)~.$$}
\hfill\qed
\end{definition}

\begin{example}~ {\rm Let $V$ be a finite-dimensional inner product $\ZZ_2$-space.
The stable $\ZZ_2$-equivariant cohomotopy Euler class of
the trivial $\ZZ_2$-equivariant $V$-bundle $\epsilon_V:X\to BO^{\ZZ_2}(V)$
over a $\ZZ_2$-space $X$
$$\gamma^{\ZZ_2}(\epsilon_V)~=~0_V
 \in \omega^0_{\ZZ_2}(X;-\epsilon_V)~=~\{X^+;V^{\infty}\}_{\ZZ_2}$$
is represented by the pointed $\ZZ_2$-map $0_V:X^+ \to V^{\infty}$
sending $X$ to $0_V \in V^{\infty}$.\\
\hfill\qed}
\end{example}

\begin{example} \label{cohomotopy} {\rm
(i) The special case
$$X~=~(\R^j \oplus L\R^{i-j})^{\infty}~=~\Sigma^j(LS^{i-j})~,~Y~=~S^0$$
of Proposition \ref{stablesequence1} (i) is the exact sequence
$$\xymatrix{\dots \ar[r]&
\omega_{i+1,j} \ar[r]^{\di{b}} & \omega_{i,j} \ar[r] &\omega_i \ar[r] &\omega_{i,j-1}
\ar[r] &\dots}$$
of Crabb \cite[Lemma (4.3)]{crabb}, with
$$\begin{array}{l}
\omega_i~=~\omega_i(\hbox{\rm pt.})~=~\{S^i;S^0\}~=~\pi_i^S~,\\[1ex]
\omega_{i,j}~=~\omega_{i,j}(\hbox{\rm pt.})~
=~\{\Sigma^j(LS^{i-j});S^0\}_{\ZZ_2}\\[1ex]
\hphantom{\omega_{i,j}~=~\omega_{i,j}(\hbox{\rm pt.})}~
=~\begin{cases}
0&{\rm if}~i<0~{\rm and}~j<0\\
\widetilde{\omega}_j(P(\infty,j-i))&{\rm if}~j<-1\\
\widetilde{\omega}_j(P(\infty,j-i))\oplus \omega_j&{\rm if}~i\leqslant j
\end{cases}
\end{array}$$
with $P(\infty,j-i)=P(\infty)/P(\R^{j-i})$ the
infinite stunted projective space (\ref{Hopfbundle}), and
$$b~=~\gamma^{\ZZ_2}(\epsilon_{L\R})~=~
1 \in \omega_{-1,0}~=~\{S^0;LS^1\}_{\ZZ_2}~=~\omega_0~=~\ZZ$$
the $\ZZ_2$-equivariant Euler class of $\epsilon_{L\R}$,
represented by $0:S^0 \to LS^1$ (\cite[Remarks 4.7]{crabb}).\\
(ii) The degree defines an isomorphism
$$\omega_0~=~\{S^0;S^0\} \xymatrix{\ar[r]^-{\di{\cong}}&} \ZZ~;~
(G:W^{\infty} \to W^{\infty}) \mapsto {\rm degree}(G)~.$$
By Example \ref{00stable} there is defined an isomorphism
$$\hbox{\rm bi-degree}~:~
\omega_{0,0}~=~\{S^0;S^0\}_{\ZZ_2} \xymatrix{\ar[r]^-{\di{\cong}}&}
\omega_0(P(\infty))\oplus \omega_0~=~\ZZ \oplus \ZZ~.$$
\hfill\qed}
\end{example}

\begin{definition} \label{Z2adjoint1} {\rm
An element $c\in O(V,U \oplus V)$ is a linear isometry $c:V \to U \oplus V$.
The {\it $\ZZ_2$-equivariant adjoint} of a map $c:X \to O(V,U\oplus V)$ is the $\ZZ_2$-equivariant pointed map
$$F_c~:~LV^{\infty}\wedge X^+  \to (LU\oplus LV)^{\infty}~;~(v,x)\mapsto c(x)(v)~.$$
\hfill\qed}
\end{definition}

\begin{proposition}~
{\rm (i)} Let $V,W$ be finite-dimensional inner product $\ZZ_2$-spaces.
For a stably trivial $\ZZ_2$-equivariant $V$-bundle $\xi$ and a
$\ZZ_2$-equivariant $V\oplus W$-bundle isomorphism
$\delta\xi:\xi \oplus \epsilon_W \cong \epsilon_{V \oplus W}$
$$\gamma^{\ZZ_2}(\xi) \in \omega^0_{\ZZ_2}(X;-\xi)~=~\{X^+;V^{\infty}\}_{\ZZ_2}$$
is the stable $\ZZ_2$-equivariant homotopy class of the adjoint of $\delta \xi$
$$\gamma^{\ZZ_2}(\xi)~:~V^{\infty}\wedge X^+  \to (V \oplus W)^{\infty}~;~
(v,x) \mapsto \delta \xi(x)(v)~,$$
exactly as in the nonequivariant case (\ref{euler1}).\\
{\rm (ii)} If $s:Y \to S(\xi\vert_Y)$ is a $\ZZ_2$-equivariant section of
$p_{S(\xi\vert_Y)}:S(\xi\vert_Y) \to Y$ (for some $Y \subseteq X$)
then
$$\gamma^{\ZZ_2}(\xi) \in
{\rm ker}(\omega^0_{\ZZ_2}(X;-\xi) \to \omega^0_{\ZZ_2}(Y;-\xi\vert_Y))
~=~{\rm im}(\omega^0_{\ZZ_2}(X,Y;-\xi) \to \omega^0_{\ZZ_2}(X;-\xi))$$
and there is defined a rel $Y$ Euler class
$\gamma^{\ZZ_2}(\xi,s) \in \omega^0_{\ZZ_2}(X,Y;-\xi)$
with image $\gamma^{\ZZ_2}(\xi) \in \omega_{\ZZ_2}^0(X;-\xi)$.\\
{\rm (iii)} The rel $Y$ Euler classes of $\ZZ_2$-equivariant sections $s_0,s_1:Y \to S(\xi\vert_Y)$
which agree on a subspace $Z \subseteq Y$ are such that
$$\begin{array}{l}
\gamma^{\ZZ_2}(\xi,s_0)-\gamma^{\ZZ_2}(\xi,s_1) \in
{\rm ker}(\omega^0_{\ZZ_2}(X,Y;-\xi) \to \omega^0_{\ZZ_2}(X,Z;-\xi))\\[1ex]
\hphantom{\gamma^{\ZZ_2}(\xi,s_0)-\gamma^{\ZZ_2}(\xi,s_1) \in}
=~{\rm im}(\omega^{-1}_{\ZZ_2}(Y,Z;-\xi) \to \omega^0_{\ZZ_2}(X,Y;-\xi))
\end{array}$$
and there is a $\ZZ_2$-equivariant difference class
$$\delta(s_0,s_1) \in \omega^{-1}_{\ZZ_2}(Y,Z;-\xi)~=~\{\Sigma (Y/Z);V^{\infty}\}_{\ZZ_2}$$
with image
$$\gamma^{\ZZ_2}(\xi,s_0)-\gamma^{\ZZ_2}(\xi,s_1) \in \omega^0_{\ZZ_2}(X,Y;-\xi)~=~
\{X/Y;V^{\infty}\}_{\ZZ_2}~.$$
\hfill\qed
\end{proposition}

\begin{example}~ \label{z2euler}
{\rm (i) Let $V$ be a finite-dimensional inner product space.
Given a $V$-bundle $\xi:X \to BO(V)$ let $\eta:X \to BO(W)$
be a $W$-bundle such that
$$\xi \oplus \eta~=~\epsilon_{V \oplus W}~:~X \to BO(V \oplus W)~,$$
so that
$$L\xi \oplus L\eta~=~\epsilon_{LV \oplus LW}~:~X \to BO^{\ZZ_2}(V \oplus W)~.$$
The $\ZZ_2$-equivariant Euler class of $L\xi:X \to BO^{\ZZ_2}(V)$
$$\gamma^{\ZZ_2}(L\xi) \in \omega^0_{\ZZ_2}(X;-L\xi)~
=~\{T(L\eta);(LV \oplus LW)^{\infty}\}_{\ZZ_2}$$
is represented by the composite $\ZZ_2$-equivariant map
$$T(L\eta) \xymatrix{\ar[r]^-{\di{z}}&} T(L\xi \oplus L\eta)~=~
(LV \oplus LW)^{\infty}\wedge X^+ \to (LV \oplus LW)^{\infty}$$
with
$$z~:~E(L\eta) \to E(L\xi \oplus L\eta)~=~(LV \oplus LW) \times X~;~(v,x) \mapsto ((0,v),x)~.$$
(ii) For $\xi=\epsilon_V$ can take $W=\{0\}$, $\eta=0$ in (i), so that
$$\gamma^{\ZZ_2}(\epsilon_{LV}) \in \omega^0_{\ZZ_2}(X;-\epsilon_{LV})~
=~\{X^+;LV^{\infty}\}_{\ZZ_2}$$
is represented by the $\ZZ_2$-equivariant map
$$X^+ \to LV^{\infty} ~;~x \mapsto 0~,~\infty \mapsto \infty~.$$
If $V$ is non-zero then $\epsilon_{LV}$ does not admit a $\ZZ_2$-equivariant
section, and $\gamma^{\ZZ_2}(\epsilon_{LV})$ is non-zero -- see
Example \ref{cohomotopy} above.\\
\hfill\qed}
\end{example}

\section{$\ZZ_2$-equivariant $S$-duality}

We recall the $\ZZ_2$-equivariant $S$-duality theory of Wirthm\"uller \cite{wirth}. (The theory is for $G$-spaces with $G$ a compact Lie group, although we shall only be concerned with the case $G=\ZZ_2$).  The theory deals with the stable $\ZZ_2$-equivariant homotopy groups
$$\{X;Y\}_{\ZZ_2}~=~\mathop{\varinjlim}\limits_U\,
[U^{\infty} \wedge X,U^{\infty}  \wedge Y]_{\ZZ_2}$$
as in Definition \ref{stabdef}, with the direct limit running over all the finite-dimensional
inner product $\ZZ_2$-spaces $U$.

\begin{definition}~{\rm
Let $X,Y$ be pointed $\ZZ_2$-spaces, and let $U-V$
be a formal difference of finite-dimensional inner product $\ZZ_2$-spaces.\\
(i) Define {\it cap} products
$$\begin{array}{l}
\widetilde{\omega}^{\ZZ_2}_0(X \wedge Y;\epsilon_U-\epsilon_V) \otimes
\widetilde{\omega}^i_{\ZZ_2}(X)
\to \widetilde{\omega}^{\ZZ_2}_{-i}(Y;\epsilon_U-\epsilon_V)~;\\[1ex]
(\sigma:V^{\infty} \to U^{\infty} \wedge X \wedge Y)\otimes
(f:X \to S^i) \mapsto ((1\wedge f)\sigma:V^{\infty} \to \Sigma^iU^{\infty} \wedge Y)~.
\end{array}$$
(ii) An element
$\sigma \in \widetilde{\omega}^{\ZZ_2}_0(X \wedge Y;\epsilon_U-\epsilon_V)$ is a
{\it $\ZZ_2$-equivariant $S$-duality} if the products
$$\sigma\otimes-~:~\widetilde{\omega}^i_{\ZZ_2}(X)\to \widetilde{\omega}^{\ZZ_2}_{-i}(Y;\epsilon_U-\epsilon_V)~~(i \in\ZZ)$$
are isomorphisms.\\
}\hfill\qed
\end{definition}

\begin{proposition}~ {\rm (Wirthm\"uller \cite{wirth})}\label{Z2S-dual} \\
{\rm (i)}
If $\sigma \in \widetilde{\omega}^{\ZZ_2}_0(X\wedge Y;\epsilon_U-\epsilon_V)$ is a
$\ZZ_2$-equivariant $S$-duality there are induced isomorphisms
$$\begin{array}{l}
\sigma~:~\{X \wedge A;B\}_{\ZZ_2} \xymatrix{\ar[r]^-{\di{\cong}}&}
\{V^{\infty}\wedge A;U^{\infty}\wedge B \wedge Y \}_{\ZZ_2}~;~F \mapsto (F \wedge 1_Y)(\sigma\wedge 1_A)~,\\[1ex]
\sigma~:~\{A\wedge Y;B\}_{\ZZ_2} \xymatrix{\ar[r]^-{\di{\cong}}&}
\{V^{\infty}\wedge A;U^{\infty}\wedge X \wedge B \}_{\ZZ_2}~;~
G \mapsto (1_X \wedge G)(1_A\wedge \sigma)
\end{array}$$
for any pointed $CW$ $\ZZ_2$-complexes $A,B$.\\
{\rm (ii)} For any finite pointed $CW$ $\ZZ_2$-complex $X$ there exist a
finite-dimensional inner product $\ZZ_2$-space $V$, a
finite pointed $CW$ $\ZZ_2$-complex $Y$ and a $\ZZ_2$-equivariant
map $\sigma:V^{\infty} \to X \wedge Y$ such that
$\sigma \in \widetilde{\omega}^{\ZZ_2}_0(X \wedge Y;-\epsilon_V)$
is a $\ZZ_2$-equivariant $S$-duality, and for any pointed $CW$ $\ZZ_2$-complex $B$
$$\{X;B\}_{\ZZ_2}~\cong~\{V^\infty;B \wedge Y\}_{\ZZ_2}~,~\{Y;B\}_{\ZZ_2}~\cong~\{V^\infty;B \wedge X\}_{\ZZ_2}~.$$
\hfill\qed
\end{proposition}

\begin{example}
{\rm  Let $X$ be a pointed $CW$ complex, let $Y$ a finite pointed $CW$ $\ZZ_2$-complex,
and let $i:Y^{\ZZ_2} \to Y$ be the inclusion.
For any stable $\ZZ_2$-equivariant map $F:X \to Y$ the fixed point stable map
$G=\rho(F):X \to Y^{\ZZ_2}$ is such that $F$ and $\sigma(G)=iG:X\to Y$
agree on the fixed points, with the relative difference $\ZZ_2$-equivariant map
$$\delta(F,\sigma(G))~:~\Sigma S(\infty)^+ \wedge X \to  L\R(\infty)^{\infty}\wedge Y$$
$\ZZ_2$-equivariantly $S$-dual (by \ref{Z2S-dual} (ii), with $V=\R(\infty)$)
to a  stable $\ZZ_2$-equivariant map
$$\delta'(F,\sigma(G))~:~X \to S(\infty)^+ \wedge Y$$
such that
$$F-\sigma(G)~=~\delta'(F,\sigma(G))  \in {\rm im}(\delta)~=~{\rm ker}(\rho)$$
in the direct sum system of Proposition \ref{stablesequence2}
$$\{X;S(\infty)^+ \wedge Y\}_{\ZZ_2}
\xymatrix{\ar@<1ex>[r]^-{\di{\gamma}}\ar@<-1ex>@{<-}[r]_-{\di{\delta}} & }
\{X;Y\}_{\ZZ_2}
\xymatrix{\ar@<1ex>[r]^-{\di{\rho}}\ar@<-1ex>@{<-}[r]_-{\di{\sigma}} & }
\{X;Y^{\ZZ_2}\}~.$$
\hfill\qed}
\end{example}

\begin{proposition}\label{Z2manifold}~Let $V$
be a finite-dimensional inner product $\ZZ_2$-space, and let $M$ be an
$m$-dimensional $V$-restricted $\ZZ_2$-manifold, so that there
exists a $\ZZ_2$-equivariant embedding $M \subset V \oplus \R^m$
with a $\ZZ_2$-equivariant normal $\vert V\vert$-bundle $\nu_M$ such that
$$\tau_M \oplus \nu_M~=~\epsilon_{V \oplus \R^m}~.$$
The composite of the $\ZZ_2$-equivariant Pontrjagin-Thom map
$\alpha:(V \oplus \R^m)^{\infty} \to T(\nu_M)$
and the diagonal map $\Delta:T(\nu_M) \to M^+ \wedge T(\nu_M)$ is a $\ZZ_2$-equivariant
map
$$\sigma~=~\Delta \alpha~:~(V \oplus \R^m)^{\infty} \to M^+ \wedge T(\nu_M)$$
such that
$$\begin{array}{ll}
\sigma \in \omega^{\ZZ_2}_0(M \times M,0\times -\tau_M)~
&=~\omega^{\ZZ_2}_m(M \times M,0\times (\nu_M-\epsilon_V))\\[1ex]
&=~\omega^{\ZZ_2}_0(M^+ \wedge T(\nu_M);-\epsilon_{V\oplus\R^m})\\[1ex]
&=~\{(V \oplus \R^m)^{\infty};M^+ \wedge T(\nu_M)\}_{\ZZ_2}
\end{array}$$
is a $\ZZ_2$-equivariant $S$-duality between
$M^+$ and $T(\nu_M)$ inducing Poincar\'e duality isomorphisms
$$\sigma \cap -~:~\omega^*_{\ZZ_2}(M)~\cong~
\omega^{\ZZ_2}_{m-*}(M,\nu_M-\epsilon_V)~.$$
The isomorphism
$\omega^0_{\ZZ_2}(M) \cong \omega^{\ZZ_2}_m(M;\nu_M-\epsilon_V)$
sends $1 \in \omega_{\ZZ_2}^0(M)=\{M^+;S^0\}_{\ZZ_2}$ to
$$\alpha \in \omega^{\ZZ_2}_0(M,-\tau_M)~=~\omega^{\ZZ_2}_m(M,\nu_M-\epsilon_V)~=~
\{(V \oplus \R^m)^{\infty};T(\nu_M)\}_{\ZZ_2}~.$$
\hfill\qed
\end{proposition}

\begin{example} {\rm Let $M$ be an $m$-dimensional free $\ZZ_2$-manifold,
with a $\ZZ_2$-equivariant embedding $M \subset V \oplus \R^m$ as in
\ref{Z2manifold}, and let $p:M \to M/\ZZ_2$ be the double covering projection,
so that
$$\begin{array}{l}
\tau_M~=~p^*\tau_{M/\ZZ}~=~\epsilon_{V\oplus \R^m}-\nu_M~=~
p^*(\epsilon_{V \oplus \R^m}/\ZZ_2)
-p^*(\nu_M/\ZZ_2)~,\\[1ex]
\tau_{M/\ZZ}~=~\epsilon_{V \oplus \R^m}/\ZZ_2-\nu_M/\ZZ_2~.
\end{array}$$
The normal bundle of an embedding $M/\ZZ_2 \subset U \oplus \R^m$ (for some
finite-dimensional inner product space $U$) is such that
$$\nu_{M/\ZZ_2}~=~\epsilon_{U \oplus \R^m}-\tau_{M/\ZZ_2}~
=~\epsilon_U -\epsilon_V/\ZZ_2 + \nu_M/\ZZ_2~,$$
and
$$\begin{array}{ll}
\alpha \in \omega^{\ZZ_2}_0(M,-\tau_M)&=~
\{(V \oplus \R^m)^{\infty};T(\nu_M)\}_{\ZZ_2}\\[1ex]
&=~\omega_0(M/\ZZ_2,-\tau_M/\ZZ_2)~=~\{(U \oplus \R^m)^{\infty};T(\nu_{M/\ZZ_2})\}
\end{array}$$
is the Pontrjagin-Thom map of $M/\ZZ_2 \subset U \oplus \R^m$.
The $\ZZ_2$-equivariant $S$-duality between $M^+$ and $T(\nu_M)$ can
also be regarded as a nonequivariant $S$-duality between $(M/\ZZ_2)^+$ and
$T(\nu_{M/\ZZ_2})$
$$\begin{array}{ll}
\sigma \in \omega^{\ZZ_2}_0(M \times M,0\times -\tau_M)~
&=~\omega^{\ZZ_2}_m(M \times T(\nu_M),-\epsilon_V)\\[1ex]
&=~\omega_0(M/\ZZ_2 \times M/\ZZ_2,0\times_{\ZZ_2} -\tau_M)\\[1ex]
&=~\omega_m((M/\ZZ_2)^+ \wedge T(\nu_{M/\ZZ_2});-\epsilon_U)~,
\end{array}$$
with the $\ZZ_2$-equivariant Poincar\'e duality isomorphisms
$$\sigma \cap -~:~\omega^*_{\ZZ_2}(M)~\cong~
\omega^{\ZZ_2}_{m-*}(M,\nu_M-\epsilon_V)$$
regarded as nonequivariant Poincar\'e duality isomorphisms
$$\sigma \cap -~:~\omega^*(M/\ZZ_2)~\cong~
\omega_{m-*}(M/\ZZ_2,\nu_{M/\ZZ_2}-\epsilon_U)~.$$
\hfill\qed}
\end{example}

\begin{example} \label{Z2adjoint2} {\rm
(i) For any finite-dimensional inner product space $V$
the inclusion $S(LV) \subset LV$ is a $\ZZ_2$-equivariant embedding with
trivial normal $\ZZ_2$-equivariant $\R$-bundle
$\nu_{S(LV) \subset LV}=\epsilon_{\R}$.
The composite of
$$\alpha_{LV}~:~LV^{\infty} \to T(\epsilon_{\R})~=~\Sigma S(LV)^+~;~[t,u] \mapsto (t,u)$$
and the diagonal map $\Delta:\Sigma S(LV)^+ \to S(LV)^+ \wedge \Sigma S(LV)^+$
$$\sigma_{LV}~=~\Delta \alpha_{LV}~:~LV^{\infty} \to S(LV)^+\wedge \Sigma S(LV)^+$$
represents a $\ZZ_2$-equivariant $S$-duality
$$\sigma_{LV} \in \omega_0(S(LV)^+\wedge \Sigma S(LV)^+;-\epsilon_{LV})~.$$
Thus for any pointed $CW$ $\ZZ_2$-complexes $A,B$ there is
defined an $S$-duality isomorphism
$$\{\Sigma S(LV)^+ \wedge A; LV^{\infty} \wedge B\}_{\ZZ_2}
\to \{A;S(LV)^+ \wedge B\}_{\ZZ_2}~;~
F \mapsto (1 \wedge F)(\Delta \alpha_{LV} \wedge 1)~.$$
The $\ZZ_2$-equivariant $S$-duality map $\sigma_{LV}$ induces a
$\ZZ[\ZZ_2]$-module chain map
$$\sigma_{LV}~:~
\dot C^{cell}(LV^{\infty}) \to C^{cell}(S(LV))\otimes_{\ZZ}SC^{cell}(S(LV))$$
with adjoint an isomorphism of $\ZZ[\ZZ_2]$-module chain complexes
$$\begin{array}{l}
\sigma_{LV}~:~{\rm Hom}_{\ZZ[\ZZ_2]}(SC^{cell}(S(LV)),\dot C^{cell}(LV^{\infty}))\\[1ex]
\xymatrix{\ar[r]^-{\di{\cong}}&}
{\rm Hom}_{\ZZ[\ZZ_2]}(\dot C^{cell}(LV^{\infty}),
C^{cell}(S(LV))\otimes_{\ZZ}\dot C^{cell}(LV^{\infty}))~=~C^{cell}(S(LV))~.
\end{array}$$
(ii) For any $\ZZ[\ZZ_2]$-module chain complexes $D,E$ the isomorphism
of (i) determines a chain level $\ZZ_2$-equivariant $S$-duality isomorphism
$$\begin{array}{l}
\sigma_{LV}~:~{\rm Hom}_{\ZZ[\ZZ_2]}(SC^{cell}(S(LV))\otimes_{\ZZ}D,\dot C^{cell}(LV^{\infty})\otimes_{\ZZ}E)\\[1ex]
\hskip100pt \xymatrix{\ar[r]^-{\di{\cong}}&}
{\rm Hom}_{\ZZ[\ZZ_2]}(D,C^{cell}(S(LV))\otimes_{\ZZ}E)~.
\end{array}$$
(iii) The $\ZZ_2$-equivariant $S$-duality isomorphism
$$\sigma_{LV}~:~\{\Sigma S(LV)^+ \wedge A;LV^{\infty} \wedge B\}_{\ZZ_2}
\xymatrix{\ar[r]^-{\di{\cong}}&} \{A;S(LV)^+\wedge B\}_{\ZZ_2} $$
induces the chain level $\ZZ_2$-equivariant $S$-duality isomorphism
of (ii) with $D=C(A)$, $E=C(B)$. A stable $\ZZ_2$-equivariant map
$F:\Sigma S(LV)^+ \wedge A
\to LV^{\infty} \wedge B$ as in \ref{Z2S-dual} (ii) induces a
$\ZZ[\ZZ_2]$-module chain map
$$f~:~S C^{cell}(S(LV))\otimes_{\ZZ}\dot C(A) \to \dot C^{cell}(LV^{\infty})
\otimes_{\ZZ}\dot C(B)~.$$
The $S$-dual stable $\ZZ_2$-equivariant map
$$G~=~(1 \wedge F)(\Delta \alpha_{LV} \wedge 1)~:~A \to S(LV)^+ \wedge B$$
induces a $\ZZ[\ZZ_2]$-module chain map
$$g~:~\dot C(A) \to C^{cell}(S(LV))\otimes_{\ZZ}\dot C(B)$$
which is just the adjoint of $f$, i.e. $f$ and $g$ correspond under the
isomorphism of (ii). Thus there is defined a commutative diagram
$$\xymatrix@C+10pt{
\{\Sigma S(LV)^+\wedge A;LV^{\infty} \wedge B\}_{\ZZ_2}
\ar[d]_-{\di{\sigma_{LV}}}^-{\di{\cong}} \ar[r]^-{Hurewicz} & H'
 \ar[d]_-{\di{\sigma_{LV}}}^-{\di{\cong}}\\
\{A,S(LV)^+ \wedge B\}_{\ZZ_2}\ar[r]^-{Hurewicz} & H''}$$
where
$$\begin{array}{l}
H'~=~H_0({\rm Hom}_{\ZZ[\ZZ_2]}(S C^{cell}(S(LV)) \otimes_{\ZZ} \dot C(A),
 C^{cell}(LV^{\infty})\otimes_{\ZZ} \dot C(B)))~,\\[1ex]
H''~=~H_0({\rm Hom}_{\ZZ[\ZZ_2]}(\dot C(A), C^{cell}(S(LV))\otimes_{\ZZ} \dot C(B)))~.
\end{array}$$
(iv) For $A=B=S^0$ in (iii) there is defined a commutative diagram of
isomorphisms
$$\xymatrix@C+10pt{
\{\Sigma S(LV)^+;LV^{\infty}\}_{\ZZ_2}
\ar[d]_-{\di{\sigma_{LV}}}^-{\di{\cong}} \ar[r]^-{Hurewicz}_-{\di{\cong}} & H'
 \ar[d]_-{\di{\sigma_{LV}}}^-{\di{\cong}}\\
\{S^0;S(LV)^+\}_{\ZZ_2}\ar[r]^-{Hurewicz}_-{\di{\cong}} & H''}$$
with
$$\begin{array}{l}
\{\Sigma S(LV)^+;LV^{\infty}\}_{\ZZ_2} \to H'=\ZZ~;\\[1ex]
\hskip50pt
(F':\Sigma S(LV)^+\to LV^{\infty}) \mapsto
\hbox{\rm semidegree}(F')=(F')_*(1)/2~,\\[1ex]
\{S^0;S(LV)^+\}_{\ZZ_2} \to H''=H_0(P(V))=\ZZ~;\\[1ex]
\hskip50pt
(F'':S^0 \to S(LV)^+) \mapsto (F'')_*(1)/2~.
\end{array}$$
\hfill\qed
}
\end{example}

\chapter{The geometric Hopf invariant}\label{geohopf}

We shall now construct the geometric Hopf invariant of a stable map
$F:V^{\infty} \wedge X \to V^{\infty} \wedge Y$.

\section{The $Q$-groups}\label{qgroups}

In the first instance we recall the definition of the various
$\ZZ_2$-hypercohomology $Q$-groups required for the quadratic
construction.

Let $A$ be a ring with an involution
$$A \to A~;~a \mapsto \overline{a}~,$$
a function such that
$$\overline{a+b}~=~\overline{a}+\overline{b}~,~\overline{ab}~=~\overline{b}.\overline{a}~,~
\overline{\overline{a}}~=~a~,~\overline{1}~=~1 \in A~.$$
The two main examples are:
\begin{itemize}
\item[1.] A commutative ring $A$ with the identity involution
$$A \to A~;~a \mapsto \overline{a}=a~.$$
\item[2.] A group ring $A=\ZZ[\pi]$ with the involution
$$\ZZ[\pi] \to \ZZ[\pi]~;~\sum\limits_{g \in \pi}n_gg \mapsto
\sum\limits_{g \in \pi}n_gw(g)g^{-1}$$
for a group morphism $w:\pi \to \ZZ_2=\{\pm 1\}$.
\end{itemize}

The involution on $A$ can be used to define the dual of a (left)
$A$-module $M$ to be the $A$-module
$$M^*~=~{\rm Hom}_A(M,A)~,~A \times M^* \to M^*~;~(a,f) \mapsto (x \mapsto f(x)\overline{a})~,$$
with a duality morphism
$${\rm Hom}_A(M,N) \to {\rm Hom}_A(N^*,M^*)~;~f \mapsto (f^*:g \mapsto gf)~.$$
The involution can also be used to define a right $A$-module structure
on the additive group of an $A$-module $M$
$$M \times A \to M~;~(x,a) \mapsto \overline{a}x~.$$
For $A$-modules $M,N$ there is defined a $\ZZ$-module
$$M \otimes_A N~=~M\otimes_{\ZZ}N/\{ax \otimes y - x\otimes \overline{a}y\,\vert\,
a \in A,\, x\in M,\,y\in N\}$$
and the natural $\ZZ$-module morphism
$$M\otimes_A N \to {\rm Hom}_A(M^*,N)~;~x\otimes y \mapsto (f \mapsto f(x)y)$$
is an isomorphism for f.g. projective $M$. In particular, for $N=A$
and a f.g. projective $A$-module $M$ there is defined a natural isomorphism
$$M \to M^{**}~;~x \mapsto (f \mapsto f(x))~.$$
\indent Given $A$-module chain complexes $C,D$ let $C\otimes_AD$,
${\rm Hom}_A(C,D)$ be the $\ZZ$-module chain complexes with
$$\begin{array}{l}
(C\otimes_AD)_n~=~\sum\limits_{p+q=n}C_p\otimes_AD_q~,~
d(x\otimes y)~=~x \otimes d_D(y)+(-)^qd_C(x)\otimes y~,\\[1ex]
{\rm Hom}_A(C,D)_n~=~\sum\limits_{q-p=n}{\rm Hom}_A(C_p,D_q)~,~d(f)~=~d_Df+(-)^qfd_C.
\end{array}$$
A cycle $f \in {\rm Hom}_A(C,D)_0$ is a chain map $f:C \to D$,
and $H_0({\rm Hom}_A(C,D))$ is the abelian group of chain homotopy
classes of chain maps $f:C \to D$.

The {\it suspension} of an $A$-module chain complex $C$ is the $A$-module
chain complex $SC$ with\index{suspension!chain complex, $SC$}
$$d_{SC}~=~d_C~:~(SC)_r~=~C_{r-1} \to (SC)_{r-1}~=~C_{r-2}~.$$
Let $C^{-*}$ be the $A$-module chain complex with
$$(C^{-*})_r~=~C^{-r}~=~{\rm Hom}_A(C_{-r},A)~,~d_{C^{-*}}~=~(d_C)^*~.$$
The natural $\ZZ$-module chain map
$$C\otimes_AD \to {\rm Hom}_A(C^{-*},D)~;~x \otimes y \mapsto (f \mapsto
\overline{f(x)}y)$$
is an isomorphism if $C$ is a bounded f.g. projective
$A$-module chain complex, in which case a homology
class $f \in H_n(C\otimes_AD)=H_0({\rm Hom}_A(C^{n-*},D))$ is a chain homotopy
class of $A$-module chain maps $f:C^{n-*} \to D$, with
$$d_{C^{n-*}}~=~(-)^rd^*_C~:~(C^{n-*})_r~=~C^{n-r} \to
(C^{n-*})_{r-1}~=~C^{n-r+1}~.$$
\indent
For an $A$-module chain complex $C$ let the generator $T \in\ZZ_2$
act on the $\ZZ$-module chain complex $C\otimes_AC$ by the signed
transposition
$$T~:~C_p\otimes_AC_q \to C_q\otimes_AC_p~;~x\otimes y \mapsto (-)^{pq}y \otimes x~,$$
so that $C\otimes_AC$ is a $\ZZ[\ZZ_2]$-module chain complex.
If $C$ is a bounded f.g. projective
$A$-module chain complex transposition corresponds to duality under the
natural isomorphism $C\otimes_AC \cong {\rm Hom}_A(C^{-*},C)$, with
$$T~:~{\rm Hom}_A(C_p^*,C_q) \to {\rm Hom}_A(C_q^*,C_p)~;~
f \mapsto (-)^{pq}f^*~.$$

\begin{definition}~ {\rm
 For $-\infty \leqslant i \leqslant j \leqslant \infty$ let
$W[i,j]$ be the $\ZZ[\ZZ_2]$-module chain complex with
$$W[i,j]_r~=~\begin{cases} \ZZ[\ZZ_2]&\hbox{\rm if $i \leqslant r \leqslant j$}\\
0&\hbox{\rm otherwise}\end{cases}~,~
d~=~1+(-)^rT~:~W[i,j]_r \to W[i,j]_{r-1}~.$$
(i) The {\it $[i,j]$-symmetric $Q$-groups} of $C$ are
\index{$Q$-groups!$[i,j]$-symmetric, $Q_{[i,j]}^n(C)$}
$$Q_{[i,j]}^n(C)~=~H_n({\rm Hom}_{\ZZ[\ZZ_2]}(W[i,j],C\otimes_AC))~.$$
An element $\phi \in Q_{[i,j]}^n(C)$ is an equivalence class of collections
$$\phi~=~\{\phi_s \in (C\otimes_AC)_{n+s}\,\vert\, i \leqslant s \leqslant j\}$$
satisfying
$$\begin{array}{l}
(d \otimes 1)\phi_s+(-)^r\phi_s(1\otimes d)+(-)^{n+s-1}
(\phi_{s-1}+(-)^sT\phi_{s-1})~=~0\\[1ex]
\hskip90pt  \in (C\otimes_AC)_{n+s-1}~=~\sum\limits_r C_{n-r+s-1}\otimes_AC_r~~(\phi_{i-1}=0)~.
\end{array} $$
(ii) The {\it $[i,j]$-quadratic $Q$-groups} of $C$ are
\index{$Q$-groups!$[i,j]$-quadratic, $Q^{[i,j]}_n(C)$}
$$Q^{[i,j]}_n(C)~=~H_n(W[i,j]\otimes_{\ZZ[\ZZ_2]}(C\otimes_AC))~.$$
An element $\psi \in Q^{[i,j]}_n(C)$ is an equivalence class of collections
$$\psi~=~\{\psi_s \in (C\otimes_AC)_{n-s}\,\vert\, i \leqslant s \leqslant j\}$$
satisfying
$$\begin{array}{l}
(d \otimes 1)\psi_s+(-)^r\psi_s(1\otimes d)+(-)^{n-s-1}(\psi_{s+1}+(-)^sT\psi_{s+1})~=~0 \\[1ex]
\hskip90pt \in (C\otimes_AC)_{n-s-1}~=~\sum\limits_r C_{n-r-s-1}\otimes_AC_r~~(\psi_{j+1}=0)~.
 \end{array} $$
\hfill\qed}
\end{definition}

\begin{proposition}~ \label{les1} {\rm (\cite[1.1]{ranicki1})}
Let $-\infty \leqslant i \leqslant j \leqslant k \leqslant \infty$.\\
{\rm (i)} The functions
$$\begin{array}{l}
Q^n_{[i,j]}(C) \to Q_{n-1}^{[-1-j,-1-i]}(C)~;~\phi \mapsto \psi~,~
\psi_s~=~\phi_{-1-s}~,\\[1ex]
Q^n_{[i,j]}(C) \to Q^{n+1}_{[i+1,j+1]}(SC)~;~\phi \mapsto S\phi~,~
S\phi_s~=~\phi_{s-1}
\end{array}$$
are isomorphisms. \\
{\rm (ii)} There are defined long exact sequences of $Q$-groups
$$\begin{array}{l}
\dots \to Q^n_{[j+1,k+1]}(C) \to Q^n_{[i,k+1]}(C) \to
Q^n_{[i,j]}(C) \to Q^{n-1}_{[j+1,k+1]}(C) \to\dots\\[1ex]
\dots \to Q_n^{[i,j]}(C) \to Q_n^{[i,k+1]}(C) \to
Q_n^{[j+1,k+1]}(C) \to Q_{n-1}^{[i,j]}(C) \to\dots
\end{array}$$
\end{proposition}
\begin{proof} (i) The functions are already isomorphisms on the chain level.\\
(ii) Immediate from the short exact sequence of $\ZZ[\ZZ_2]$-module
chain complexes
$$0 \to W[i,j] \to  W[i,k+1] \to W[j+1,k+1] \to 0~.$$
\hfill\qed\end{proof}

For an inner product space $V$ with ${\rm dim}(V)=k$
($1 \leqslant k \leqslant \infty$) give $S(LV)$ the standard
$\ZZ_2$-equivariant $CW$ structure with cells
$$e_0~,~Te_0~,~e_1~,~Te_1~,~\dots~,~e_{k-1}~,~Te_{k-1}$$
and $(k-1)$-dimensional free $\ZZ[\ZZ_2]$-module cellular chain complex
$$C^{cell}(S(LV))~=~W[0,k-1]~.$$

\begin{definition}~\label{iso} {\rm (\cite[pp. 100-101]{ranicki1})}
(i) Given $\ZZ[\ZZ_2]$-modules
$M,N$ define a $\ZZ[\ZZ_2]$-module structure on ${\rm Hom}_{\ZZ}(M,N)$ by
$$T~:~{\rm Hom}_{\ZZ}(M,N) \to {\rm Hom}_{\ZZ}(M,N)~;~f \mapsto T_NfT_M~.$$
(ii) Given an inner product space $V$ and
 $\ZZ[\ZZ_2]$-module chain complexes $C,D$ define a
{\it $V$-coefficient $\ZZ_2$-isovariant chain map} $f:C \to D$ to be a cycle
\index{$\ZZ_2$-isovariant!$V$-coefficient chain map}
$$f~=~\{f_s\,\vert\, 0 \leqslant s \leqslant {\rm dim}(V)-1\}
\in {\rm Hom}_{\ZZ[\ZZ_2]}(C^{cell}(S(LV)),{\rm Hom}_{\ZZ}(C,D))_0~,$$
with $f_s \in {\rm Hom}_{\ZZ}(C_r,D_{r+s})$ such that
$$\begin{array}{l}
d_Df_s+(-)^{s-1}f_sd_C+(-)^{s-1}(f_{s-1}+(-)^sT_Df_{s-1}T_C)~=~0~:\\[1ex]
\hskip100pt C_r \to D_{r+s-1}~~(f_{-1}=0)~,
\end{array}$$
which can be regarded as a $\ZZ[\ZZ_2]$-module chain map $f:C^{cell}(S(LV)) \otimes_{\ZZ} C \to D$.
Thus $f_0:C \to D$ is a $\ZZ$-module chain map, $f_1:f_0 \simeq T_Df_0T_C$
is a chain homotopy, $f_2$ is a higher chain homotopy, etc.\\
(iii) There is a corresponding notion of {\it $V$-coefficient
$\ZZ_2$-isovariant chain homotopy}\index{$\ZZ_2$-isovariant!$V$-coefficient chain homotopy}
$h:f\simeq g:C \to D$, and
$$\begin{array}{l}
H_0({\rm Hom}_{\ZZ[\ZZ_2]}(C^{cell}(S(LV)),{\rm Hom}_{\ZZ}(C,D)))\\[1ex]
\hskip100pt =~
H_0({\rm Hom}_{\ZZ[\ZZ_2]}(C^{cell}(S(LV))\otimes_{\ZZ}C,D))
\end{array}$$
is the abelian group of $V$-coefficient
$\ZZ_2$-isovariant chain homotopy classes of $V$-coefficient
$\ZZ_2$-isovariant chain maps.\\
(iv) A $V$-coefficient $\ZZ_2$-isovariant chain map $f:C \to D$
{\it induces} a $\ZZ$-module chain map
$$\begin{array}{l}
f^{\%}~:~{\rm Hom}_{\ZZ[\ZZ_2]}(C^{cell}(S(LV)),C) \to
{\rm Hom}_{\ZZ[\ZZ_2]}(C^{cell}(S(LV)),D)~;\\[1ex]
\hskip100pt g=\{g_s\} \mapsto (f\otimes g)\Delta
\end{array}$$
with
$$\Delta~:~C^{cell}(S(LV)) \to C^{cell}(S(LV))\otimes_{\ZZ}C^{cell}(S(LV))~;~
1_r \mapsto \sum\limits^r_{s=0}1_s\otimes (T_{r-s})^s$$
the {\it cellular diagonal chain approximation} $\ZZ[\ZZ_2]$-module chain map.
\index{cellular diagonal chain approximation, $\Delta$}\\
(v) The {\it composite} of $V$-coefficient $\ZZ_2$-isovariant chain maps
$f:C \to D$, $g:D \to E$ is the $V$-coefficient $\ZZ_2$-isovariant chain map
$gf:C\to E$ with\index{$\ZZ_2$-isovariant!composite chain map}
$$\begin{array}{l}
gf~=~f^{\%}(g)~:~C^{cell}(S(LV))\otimes_{\ZZ}C
\xymatrix{\ar[r]^-{\di{\Delta\otimes 1}} & } \\[1ex]
C^{cell}(S(LV))\otimes_{\ZZ}C^{cell}(S(LV))\otimes_{\ZZ}C
\xymatrix{\ar[r]^-{\di{1\otimes f}}&}
C^{cell}(S(LV))\otimes_{\ZZ}D \xymatrix{\ar[r]^-{\di{g}}&}  E~.
\end{array}$$
(vi) The {\it cup product} of $V$-coefficient $\ZZ_2$-isovariant chain maps
$f:C \to D$, $g:E \to F$ is the $V$-coefficient $\ZZ_2$-isovariant chain map
with\index{$\ZZ_2$-isovariant!cup product}
$$\begin{array}{l}
f \cup g~=~(f\otimes g)\Delta:~C^{cell}(S(LV))\otimes_{\ZZ}C \otimes_{\ZZ}E
\xymatrix{\ar[r]^-{\di{\Delta\otimes 1}} & } \\[1ex]
C^{cell}(S(LV))\otimes_{\ZZ}C^{cell}(S(LV))\otimes_{\ZZ}C\otimes_{\ZZ}E
\xymatrix{\ar[r]^-{\di{f\otimes g}}&}  D\otimes_{\ZZ}F~.
\end{array}$$
\hfill\qed
\end{definition}

\begin{example} {\rm
(i) For $V=\R$ a $V$-coefficient $\ZZ_2$-isovariant chain map $f:C \to D$
is just a $\ZZ$-module chain map $f_0:C \to D$.\\
(ii) For $V=\R^2$ a $V$-coefficient $\ZZ_2$-isovariant chain map $f:C \to D$
is a $\ZZ$-module chain map $f_0:C \to D$ together with a $\ZZ$-module
chain homotopy $f_1:f_0 \simeq T_Df_0T_C:C \to D$.}\\
\hfill\qed
\end{example}

\begin{proposition}~
Let $f:C \to D$ be a $V$-coefficient $\ZZ_2$-isovariant chain map
of $\ZZ[\ZZ_2]$-module chain complexes which are bounded below.
If $f_0:C \to D$ induces isomorphisms $(f_0)_*:H_*(C) \cong H_*(D)$
the $\ZZ$-module chain map $f^{\%}$ induces isomorphisms
$$(f^{\%})_*~:~H_*({\rm Hom}_{\ZZ[\ZZ_2]}(C^{cell}(S(LV)),C))~\cong~
H_*({\rm Hom}_{\ZZ[\ZZ_2]}(C^{cell}(S(LV)),D))~.$$
\end{proposition}
\begin{proof} Standard homological algebra.\\
\hfill\qed\end{proof}

Let $\widehat{S}(LV)$ be the suspension $\ZZ_2$-spectrum with
$$\widehat{S}(LV)_j~=~\Sigma^{j-k}S(LV \oplus L\R^j)^+~~(j\geqslant k)$$
and $\ZZ[\ZZ_2]$-module cellular chain complex
$$\dot C^{cell}(\widehat{S}(LV))~=~C^{cell}(S(LV \oplus LV))_{*+k}~=~W[-k,k-1]~.$$

\begin{definition}~ {\rm Let $C$ be an $A$-module chain complex, and
let $V$ be an inner product space with ${\rm dim}(V)=k$ ($1 \leqslant k \leqslant\infty$).\\
(i) The {\it $V$-coefficient symmetric $Q$-groups} of $C$ are
\index{$Q$-groups!$V$-coefficient symmetric, $Q_V^n(C)$}
$$Q_V^n(C)~=~H_n({\rm Hom}_{\ZZ[\ZZ_2]}(C^{cell}(S(LV)),C\otimes_AC))~=~Q_{[0,k-1]}^n(C)~.$$
(ii) The {\it $V$-coefficient quadratic $Q$-groups} of $C$ are
\index{$Q$-groups!$V$-coefficient quadratic, $Q^V_n(C)$}
$$Q^V_n(C)~=~H_n(C^{cell}(S(LV))\otimes_{\ZZ[\ZZ_2]}(C\otimes_AC))~=~Q^{[0,k-1]}_n(C)~.$$
(iii) The {\it $V$-coefficient hyperquadratic $Q$-groups} of $C$ are
\index{$Q$-groups!$V$-coefficient hyperquadratic, $\widehat{Q}_V^n(C)$}
$$\widehat{Q}^n_V(C)~=~H_n({\rm Hom}_{\ZZ[\ZZ_2]}(C^{cell}(\widehat{S(LV)}),C\otimes_AC))~=~Q_{[-k,k-1]}^n(C)~.$$}
\hfill\qed
\end{definition}

\begin{proposition}~
{\rm (i)} An $A$-module chain map $f:C \to D$ determines a
$\ZZ[\ZZ_2]$-module chain map
$$f \otimes f~:~C \otimes_AC \to D \otimes_A D$$
which can be regarded as a $V$-coefficient $\ZZ_2$-isovariant chain map
with
$$(f\otimes f)_s~=~\begin{cases}
f\otimes f&{\it if}~s=0\\[1ex]
0&{\it if}~s\geqslant 1
\end{cases}~:~(C\otimes_AC)_r \to (D\otimes_AD)_{r+s}~.$$
{\rm (ii)} An $A$-module chain homotopy $g:f \simeq f':C \to D$ determines a
$V$-coefficient $\ZZ_2$-isovariant chain homotopy
$$g\otimes g~:~f\otimes f~\simeq~f'\otimes f'~:~C\otimes_AC \to  D \otimes_AD$$
with
$$\begin{array}{l}
(g\otimes g)_s~=~\begin{cases}
f\otimes g + (-)^qg \otimes f'&{\it if}~s=0\\
(-)^qg\otimes g&{\it if}~s=1\\
0&{\it if}~s\geqslant 2
\end{cases}~:\\[5ex]
\hskip75pt (C\otimes_AC)_r \to (D\otimes_AD)_{r+s+1}~=~
\sum\limits_q D_{-q+r+s+1}\otimes_AD_q~.
\end{array}$$
{\rm (iii)} An $A$-module chain map $f:C \to D$ induces
morphisms in the $Q$-groups
$$\begin{array}{l}
f^{\%}~:~Q^n_V(C) \to Q^n_V(D)~,\\[1ex]
f_{\%}~:~Q_n^V(C) \to Q_n^V(D)~,\\[1ex]
\widehat{f}^{\%}~:~\widehat{Q}^n_V(C) \to \widehat{Q}^n_V(D)
\end{array}$$
which depend only on the chain homotopy class of $f$, and are isomorphisms
if $f$ is a chain equivalence.\\
{\rm (iv)} The chain level cup products
$$\begin{array}{l}
{\rm Hom}_{\ZZ}(C^{cell}S(LV)),C\otimes_AC)
\otimes_{\ZZ} {\rm Hom}_{\ZZ}(C^{cell}S(LV)),D\otimes_BD)\\[1ex]
\to {\rm Hom}_{\ZZ}(C^{cell}S(LV)),(C\otimes_{\ZZ}D)\otimes_{A\otimes_{\ZZ}B}(C\otimes_{\ZZ}D))~;\\[1ex]
\hskip100pt \phi \otimes \theta \mapsto \phi \cup \theta=(\phi\otimes \theta)\Delta
\end{array}$$
induce cup products in the $V$-coefficient symmetric $Q$-groups
$$\cup~:~Q_V^m(C)\otimes_{\ZZ}Q_V^n(D) \to Q_V^{m+n}(C\otimes_{\ZZ}D)~;~
\phi \otimes \theta \mapsto \phi \cup \theta~.$$
Similarly for products involving
the $V$-coefficient quadratic and hyperquadratic $Q$-groups
$$\begin{array}{l}
Q_V^m(C)\otimes_{\ZZ}Q^V_n(D) \to Q^V_{m+n}(C\otimes_{\ZZ}D)~,\\[1ex]
\widehat{Q}^V_m(C)\otimes_{\ZZ}\widehat{Q}^n_V(D) \to \widehat{Q}_V^{m+n}(C\otimes_{\ZZ}D)~.
\end{array}$$
\end{proposition}
\begin{proof} See \cite[pp. 100-101,\,\S 8]{ranicki1}.\\
\hfill\qed\end{proof}

\begin{definition}~ {\rm
Let $\sigma \in {\rm Hom}_{\ZZ[\ZZ_2]}(W[0,1],S\ZZ\otimes_AS\ZZ)_1$ be
the cycle defined by
$$\sigma_1~=~1~:~(S\ZZ)^1~=~\ZZ \to (S\ZZ)_1~=~\ZZ~.$$
The {\it suspension} chain map is defined for any $A$-module chain complex $C$
to be the evaluation of the cup product on $\sigma$
$$\begin{array}{l}
S~=~\sigma \cup -~:\\[1ex]
{\rm Hom}_{\ZZ[\ZZ_2]}(W[0,k],C\otimes_AC)\to
{\rm Hom}_{\ZZ[\ZZ_2]}(W[0,k+1],SC\otimes_ASC)_{*+1}~; \phi \mapsto S\phi~,\\[1ex]
S\phi_s=\phi_{s-1} \in (SC\otimes_A SC)_{n+s+1}=(C\otimes_A C)_{n+s-1}~
(0 \leqslant s \leqslant k+1,~\phi_{-1}=0)
\end{array}$$
inducing suspension maps in the symmetric $Q$-groups.
$$S~:~Q^n_V(C) \to Q^{n+1}_{V \oplus \R}(SC)~.$$
Similarly, product with $\sigma$ induces suspension maps in the
quadratic and hyperquadratic $Q$-groups
$$\begin{array}{l}
S~:~Q_n^V(C) \to Q_{n+1}^{V \oplus \R}(SC)~,\\[1ex]
S~:~\widehat{Q}^n_V(C) \to \widehat{Q}^{n+1}_{V \oplus \R}(SC)~.
\end{array}$$}
\hfill\qed
\end{definition}

\begin{proposition}~ \label{les2} {\rm (\cite[1.3]{ranicki1})}\\
{\rm (i)} For any inner product spaces $U,V$ there is defined a
long exact sequence of $Q$-groups
$$\dots \to Q_n^V(C) \xymatrix{\ar[r]^-{1+T}&} Q^n_U(C)
\xymatrix{\ar[r]^-{S^{{\rm dim}(V)}}&}
Q^{n+{\rm dim}(V)}_{U \oplus V}(S^{{\rm dim}(V)}C) \to Q^V_{n-1}(C) \to \dots$$
with $S^{{\rm dim}(V)}C=C_{*-{\rm dim}(V)}$ and
$$\begin{array}{l}
1+T~:~Q^U_n(C) \to Q_V^n(C)~;~\psi \mapsto (1+T)\psi~,\\[1ex]
\hskip100pt (1+T)\psi_s~=~
\begin{cases}
\psi_0+T\psi_0&\hbox{\rm if $s=0$}\\
0&\hbox{\rm if $s\geqslant 1$~,}
\end{cases}\\[1ex]
Q^{n+{\rm dim}(V)}_{U \oplus V}(S^{{\rm dim}(V)}C) \to Q^V_{n-1}(C)~;\\[1ex]
\hskip50pt \phi \mapsto \psi~,~
\psi_s~=~\phi_{{\rm dim}(V)-s-1}~(0 \leqslant s \leqslant {\rm dim}(V)-1)
\end{array}$$
{\rm (ii)} If $U=\{0\}$ then $Q_U^*(C)=0$ and the sequence in {\rm (i)}
gives the isomorphisms
$$Q^{n+{\rm dim}(V)}_V(S^{{\rm dim}(V)}C) \to Q^V_{n-1}(C)~;~
\phi \mapsto \psi~,~\psi_s~=~\phi_{{\rm dim}(V)-s-1}$$
of {\rm \ref{les1} (i)}, which are induced by the
chain level $\ZZ_2$-equivariant $S$-duality isomorphism of
Example \ref{Z2adjoint2} (ii)
$$\begin{array}{l}
{\rm Hom}_{\ZZ[\ZZ_2]}(SC^{cell}(S(LV)),\dot C^{cell}(LV^{\infty})\otimes_{\ZZ}(C \otimes_AC))\\[1ex]
\hskip100pt \xymatrix{\ar[r]^-{\di{\cong}}&}
C^{cell}(S(LV))\otimes_{\ZZ}(C \otimes_AC)~.
\end{array}$$
{\rm (iii)} If $U=V$ there are defined isomorphisms
$$\begin{array}{l}
Q^{n+{\rm dim}(V)}_{V \oplus V}(S^{{\rm dim}(V)}C) \to \widehat{Q}^n_V(C)~;\\[1ex]
\hskip50pt \phi \mapsto \widehat{\phi}~,~
\widehat{\phi}_s~=~\phi_{s+{\rm dim}(V)}~(-{\rm dim}(V) \leqslant s \leqslant {\rm dim}(V)-1)
\end{array}$$
and the sequence in {\rm (i)} is
$$\xymatrix@C+5pt{
\dots \ar[r] & Q_n^V(C) \ar[r]^-{\di{1+T}} & Q^n_V(C) \ar[r]^-{\di{J}} &
\widehat{Q}^n_V(C) \ar[r]^-{\di{H}} & Q^V_{n-1}(C) \ar[r] & \dots}$$
with
$$\begin{array}{l}
J~:~Q_V^n(C) \to \widehat{Q}_V^n(C)~;~\phi \mapsto J\phi~,~J\phi_s~=~
\begin{cases}
\phi_s&\hbox{\rm if $0 \leqslant s\leqslant {\rm dim}(V)-1$}\\
0&\hbox{\rm if $-{\rm dim}(V) \leqslant s\leqslant -1$~,}
\end{cases}\\[1ex]
H~:~\widehat{Q}^n_V(C) \to Q_{n-1}^V(C)~;~\widehat{\phi} \mapsto H\widehat{\phi}~,~
H\widehat{\phi}_s~=~\widehat{\phi}_{-s-1}~.
\end{array}$$
{\rm (iv)} The suspension maps in the hyperquadratic $Q$-groups
$$S~:~\widehat{Q}^n_V(C) \to \widehat{Q}^{n+1}_{V \oplus \R}(SC)$$
are isomorphisms, and there is defined a commutative braid of exact
sequences of $Q$-groups

$$\xymatrix@C-20pt{
\widehat{Q}^{n+1}_V(C)\ar[dr]^-{\di{H}} \ar@/^2pc/[rr]^-{\di{H}}  &&
Q_n^{V \oplus \R}(SC)\ar[dr]^-{\di{1+T}}\ar@/^2pc/[rr]&& H_{n-1}(C \otimes_A C)\\
&Q^V_n(C) \ar[ur]^-{\di{S}} \ar[dr]^-{\di{1+T}}&&
Q^n_{V \oplus \R}(SC)\ar[ur] \ar[dr]^-{\di{J}} &\\
H_n(C\otimes_AC)\ar@/_2pc/[rr]\ar[ur] &&
Q^n_V(C)\ar@/_2pc/[rr]^-{\di{J}}\ar[ur]^-{\di{S}} &&\widehat{Q}_V^n(C)}$$

\end{proposition}
\begin{proof} (i) This is a special case of Proposition \ref{les1}
$$\begin{array}{l}
\dots \to Q_n^V(C)=Q^{n+1}_{[-{\rm dim}(V),-1]}(C)
\to Q^n_U(C)=Q^n_{[0,{\rm dim}(U)-1]}(C) \to\\[1ex]
Q^{n+{\rm dim}(V)}_{U \oplus V}(C)=Q^n_{[-{\rm dim}(V),{\rm dim}(U)-1]}(C)
\to Q_{n-1}^V(C)=Q^n_{[-{\rm dim}(V),-1]}(C) \to \dots~.
\end{array}$$
Specifically, let ${\rm dim}(U)=j$, ${\rm dim}(V)=k$, and note that
the homotopy cofibration sequence of pointed $CW$-$\ZZ_2$-complexes
$$S(LV)^+ \to S(LU\oplus LV)^+ \to S(LU)^+ \wedge LV^{\infty}$$
induces the short exact sequence of cellular free $\ZZ[\ZZ_2]$-module
chain complexes
$$\begin{array}{l}
0 \to C^{cell}(S(LV))=W[0,k-1]\to C^{cell}(S(LU\oplus LV))=W[0,j+k-1]\\[1ex]
\hskip100pt
\to \dot C^{cell}(S(LU)^+\wedge LV^{\infty})=W[k,j+k-1] \to 0
\end{array}$$
inducing a long exact sequence
$$\begin{array}{l}
\dots \to Q_n^V(C)=Q_n^{[0,k-1]}(C)=Q^{n+k+1}_{[0,k-1]}(S^kC)\\[1ex]
\hskip50pt \to Q^n_U(C)=Q^n_{[0,j-1]}(C)=Q^{n+k}_{[k,j+k-1]}(S^kC)\\[1ex]
\hskip50pt \to Q^{n+{\rm dim}(V)}_{U \oplus V}(S^{{\rm dim}(V)}C)=Q^{n+k}_{[0,j+k-1]}(S^kC)\\[1ex]
\hskip100pt \to Q_{n-1}^V(C)=Q^{n+k}_{[0,k-1]}(S^kC) \to \dots~.
\end{array}$$
(ii) Immediate from (i).\\
(iii) Immediate from (i), noting that the long exact sequence
is induced from the short exact sequence of $\ZZ[\ZZ_2]$-module
cellular chain complexes
$$0 \to C^{cell}(S(LV))^{-*-1} \to  \dot C^{cell}(\widehat{S}(LV)) \to C^{cell}(S(LV)) \to 0~.$$
(iv) Immediate from (ii) and (iii).\\
\hfill\qed\end{proof}

\begin{terminology}~ \label{term} {\rm (\cite[\S1]{ranicki1})\\
The infinite-dimensional inner product space
$$V~=~\R(\infty)~=~\varinjlim\limits_k\R^k$$
has unit sphere
$$S(LV)~=~S(\infty)~=~E\ZZ_2$$
a contractible space with a free $\ZZ_2$-action, and
$$\begin{array}{l}
C^{cell}(S(\infty))~=~W[0,\infty]~=~W~:\\[1ex]
\hskip25pt
\dots \to W_2=\ZZ[\ZZ_2]\xymatrix{\ar[r]^-{\di{1+T}}&} W_1=\ZZ[\ZZ_2]
\xymatrix{\ar[r]^-{\di{1-T}}&} W_0=\ZZ[\ZZ_2]
\end{array}$$
is the standard free $\ZZ[\ZZ_2]$-module resolution of $\ZZ$.\\
(i) The {\it symmetric $Q$-groups} of an $A$-module chain complex $C$ are
$$Q^n(C)~=~Q_{\R(\infty)}^n(C)~=~Q_{[0,\infty]}^n(C)~.$$
(ii) The {\it quadratic $Q$-groups} of an $A$-module chain complex $C$ are
$$Q_n(C)~=~Q^{\R(\infty)}_n(C)~=~Q^{[0,\infty]}_n(C)~.$$
(ii) The {\it hyperquadratic $Q$-groups} of an $A$-module chain complex $C$ are
$$\widehat{Q}^n(C)~=~\widehat{Q}_{\R(\infty)}^n(C)~=~Q_{[-\infty,\infty]}^n(C)~.$$
}
\hfill\qed
\end{terminology}

\begin{proposition}~ \label{les3} {\rm (\cite[1.2]{ranicki1})}
Let $C$ be an $A$-module chain complex.\\
{\rm (i)} The symmetric and quadratic $Q$-groups are related by
a long exact sequence
$$\xymatrix{
\dots \ar[r] & Q^{[0,k-1]}_n(C) \ar[r]^-{\di{1+T}} & Q^n(C) \ar[r]^-{\di{S^k}} &
Q^{n+k}(S^kC) \ar[r]^-{\di{H}} & Q^{[0,k-1]}_{n-1}(C) \ar[r] & \dots}$$
for any $k \geqslant 0$.\\
{\rm (ii)} The symmetric, quadratic and hyperquadratic $Q$-groups are related by
a long exact sequence
$$\xymatrix{
\dots \ar[r] & Q_n(C) \ar[r]^-{\di{1+T}} & Q^n(C) \ar[r]^-{\di{J}} &
\widehat{Q}^n(C) \ar[r]^-{\di{H}} & Q_{n-1}(C) \ar[r] & \dots}$$
and
$$\widehat{Q}^n(C)~=~\varinjlim_k Q^{n+k}(S^kC)$$
is the direct limit of the suspension maps
$$\xymatrix{Q^n(C) \ar[r]^-{\di{S}} &Q^{n+1}(SC) \ar[r]^-{\di{S}} &Q^{n+2}(S^2C) \ar[r] &\dots~.}$$
\end{proposition}
\begin{proof}
(i) This is the special case $U=\R(\infty)$, $V=\R^k$ of \ref{les2}.\\
(ii) This is the special case $U=V=\R(\infty)$ of \ref{les2}.\\
\hfill\qed\end{proof}

\begin{definition}~{\rm Let $C$ be an $A$-module chain complex,
and  $1 \leqslant k \leqslant \infty$.\\
{\rm (i)} A {\it symmetric $k$-class} is an element $\phi \in Q_{[0,k-1]}^n(C)$.
A symmetric $\infty$-class is an element $\phi \in Q^n_{[0,\infty]}(C)=Q^n(C)$,
and this is called a {\it symmetric class}. A symmetric class $\phi$ determines a
symmetric $k$-class $\phi[0,k-1]$ for every $k$, via the morphism
$$Q^n(C) \to Q^n_{[0,k-1]}(C)~;~\phi=\{\phi_s\,\vert\,0 \leqslant s <\infty\}
\mapsto \phi[0,k-1]=\{\phi_s\,\vert\,0 \leqslant s <k\}~.$$
{\rm (ii)} A {\it quadratic $k$-class} is an element $\psi \in Q^{[0,k-1]}_n(C)$.
For $k=\infty$ this is an element $\psi \in Q_n(C)$, and this is called a
{\it quadratic class}.\\
{\rm (iii)} A {\it quadratic refinement} of a symmetric $k$-class $\phi$
is a quadratic $k$-class $\psi$ such that $(1+T)\psi=\phi$.}\\
\hfill\qed
\end{definition}

\begin{proposition}~\label{sym55}
{\rm (i)} For any $A$-module chain complex $C$ and $0 \leqslant k \leqslant \infty$
there is defined a commutative braid of exact sequences

$$\xymatrix@C-20pt{
Q^{[0,k-1]}_n(C)\ar[dr] \ar@/^2pc/[rr]  &&
Q^n_{[0,k-1]}(C)\ar[dr]^-{\di{S^k}}\ar@/^2pc/[rr]&& Q^{n-1}_{[k,\infty]}(C)\\
&Q^n(C) \ar[ur] \ar[dr]^-{\di{S^k}}&&
Q^{n+k}_{[-k,k-1]}(S^kC)\ar[ur] \ar[dr] &\\
Q^n_{[k,\infty]}(C) \ar@/_2pc/[rr]\ar[ur] &&
Q^{n+k}(S^kC)\ar@/_2pc/[rr]\ar[ur] &&Q_{n-1}^{[0,k-1]}(C)}$$

\noindent with
$$\begin{array}{c}
Q^n_{[-k,k-1]}(C)~=~Q^{n+k}_{[-k,k-1]}(S^kC)~,\\[1ex]
Q^n_{[-k,\infty]}(C)~=~Q^{n+k}(S^kC)~,~
Q^n_{[k,\infty]}(C)~=~Q^{n+k}_{[0,k-1]}(S^kC)~.
\end{array}$$
{\rm (ii)} The following conditions on symmetric class $\phi \in Q^n(C)$
are equivalent:
\begin{itemize}
\item[{\rm (a)}]~there exists a quadratic refinement
$\psi \in Q^{[0,k-1]}_n(C)$
\item[{\rm (b)}]~$S^k\phi=0 \in Q^{n+k}(S^kC)$,
\item[{\rm (c)}]~$S^k\phi[0,k-1]=0 \in Q_{[0,k-1]}^{n+k}(S^kC)$,
\end{itemize}
In particular, a chain level solution of $S^k\phi=0$ gives a particular $\psi$,
in which case
$$\begin{array}{l}
(1+T)\psi~=~\phi\\[1ex]
\in {\rm im}(1+T:Q^{[0,k-1]}_n(C) \to Q^n(C))~=~
{\rm ker}(S^k:Q^n(C) \to Q^{n+k}(S^kC))~.
\end{array}$$
This applies also for $k=\infty$, with $Q^{[0,k-1]}_n(C)=Q_n(C)$,
$S^k=J$.\\
{\rm (iii)} Let $A=\ZZ[\pi]$ be a group ring with an involution, so that for any $A$-modules
$M,N$  the diagonal $\pi$-action can be used to regard
$M \otimes_{\ZZ}N$ as an $A$-module with
$$\ZZ\otimes_A(M \otimes_{\ZZ}N)~=~M \otimes_AN~.$$
For any $A$-module chain complexes $C,D$ and $k \geqslant 1$ there is
a natural isomorphism of short exact sequences of $\ZZ$-module chain complexes
$$\xymatrix@C-15pt{
0 \ar[d] & 0 \ar[d] \\
{\rm Hom}_{A[\ZZ_2]}(W[k,\infty]\otimes_{\ZZ}S^kC,S^kD\otimes_{\ZZ}S^kD)
\ar[d] \ar[r]^-{\di{\cong}}&
{\rm Hom}_{A[\ZZ_2]}(W\otimes_{\ZZ}C,D \otimes_{\ZZ}D) \ar[d]^-{\di{S^k}} \\
{\rm Hom}_{A[\ZZ_2]}(W\otimes_{\ZZ}S^kC,S^kD\otimes_{\ZZ}S^kD)
\ar[d] \ar[r]^-{\di{\cong}}&
{\rm Hom}_{A[\ZZ_2]}(W\otimes_{\ZZ}S^kC,S^kD \otimes_{\ZZ}S^kD)\ar[d]\\
{\rm Hom}_{A[\ZZ_2]}(W[0,k-1]\otimes_{\ZZ}S^kC,S^kD\otimes_{\ZZ}S^kD)
\ar[d] \ar[r]^-{\di{\cong}}&
{\rm Hom}_{A[\ZZ_2]}(C,SW[0,k-1]\otimes_{\ZZ}(D \otimes_{\ZZ}D))\ar[d]\\
0 & 0}$$
as well as an isomorphism
$${\rm Hom}_{A[\ZZ_2]}(SW[0,k-1]\otimes_{\ZZ}S^kC,S^kD\otimes_{\ZZ}S^kD)~
\cong~{\rm Hom}_{A[\ZZ_2]}(C,W[0,k-1]\otimes_{\ZZ}(D \otimes_{\ZZ}D))$$
Suppose given $A$-module chain complexes $C,C',D,D'$ with $A$-module chain maps
$$E_C~:~S^kC \to C'~,~E_D~:~S^kD \to D'~,~F~:~C \to D~,~F'~:~C' \to D'~,$$
such that $E_C,E_D$ are chain equivalences with
$$F' E_C~=~E_D F~:~C \to D'~.$$
Suppose given also $A[\ZZ_2]$-module chain maps
$$\begin{array}{l}
\phi_C~:~W\otimes_{\ZZ}C \to C \otimes_{\ZZ} C~,~
\phi_{C'}~:~W\otimes_{\ZZ}C' \to C' \otimes_{\ZZ} C'~,\\[1ex]
\phi_D~:~W\otimes_{\ZZ}D \to D \otimes_{\ZZ} D~,~
\phi_{D'}~:~W\otimes_{\ZZ}D' \to D' \otimes_{\ZZ} D'
\end{array}$$
such that
$$(F' \otimes F')\phi_{C'}~=~\phi_{D'}(1\otimes F')~:~W \otimes_{\ZZ}C' \to D'\otimes_{\ZZ}D'$$
and a chain homotopy commutative diagram of $A[\ZZ_2]$-module chain complexes and chain maps
$$\xymatrix@C+20pt{
W\otimes_{\ZZ}S^kC\ar@{~>}[drr]^-{\delta\phi_C}
\ar[ddd]^-{1\otimes F}
\ar[dr]_-{1\otimes  E_C}^-{\simeq}\ar[rrr]^-{S^k \phi_C} &&&
S^kC\otimes_{\ZZ}S^k C\ar[ddd]^-{F\otimes F}
\ar[dl]_-{E_C \otimes E_C}^-{\simeq} \\
&W\otimes_{\ZZ}C'\ar[r]^-{\phi_{C'}}\ar[d]^-{1\otimes F'}&
C'\otimes_{\ZZ}C' \ar[d]^-{F'\otimes F'}&\\
&W\otimes_{\ZZ}D'\ar[r]^-{\phi_{D'}}& D'\otimes_{\ZZ}D' &\\
W\otimes_{\ZZ}S^kD\ar@{~>}[urr]^-{\delta\phi_D}
\ar[ur]^-{1\otimes E_D}_-{\simeq}
\ar[rrr]^-{S^k\phi_D} & &&
\ar[ul]^-{E_D \otimes E_D}_-{\simeq} S^kD\otimes_{\ZZ}S^kD}$$
with chain homotopies
$$\begin{array}{l}
\delta\phi_C~:~(E_C\otimes E_C)(S^k\phi_C)~\simeq~\phi_{C'}(1\otimes E_C)~:~
W\otimes_{\ZZ}S^kC \to C' \otimes_{\ZZ} C'~,\\[1ex]
\delta\phi_D~:~(E_D\otimes E_D)(S^k\phi_D)~\simeq~\phi_{D'}(1\otimes E_D)~:~
W\otimes_{\ZZ}S^kD \to D' \otimes_{\ZZ} D'~.
\end{array}$$
The data determines an $A$-module chain map
$$\psi~=~\delta(F,\phi_C,\phi_D)~:~C \to W[0,k-1]\otimes_{\ZZ[\ZZ_2]}(D \otimes_{\ZZ}D)$$
such that
$$(1+T)\psi\simeq~(F\otimes F)\phi_C - \phi_DF~:~
C \to {\rm Hom}_{\ZZ[\ZZ_2]}(W,D \otimes_{\ZZ} D)~.$$
The induced $\ZZ$-module chain maps
$$\begin{array}{l}
\phi_C~:~\ZZ\otimes_AC \to {\rm Hom}_{\ZZ[\ZZ_2]}(W,C \otimes_A C)~,\\[1ex]
\phi_D~:~\ZZ\otimes_AD \to  {\rm Hom}_{\ZZ[\ZZ_2]}(W,D \otimes_A D)~,\\[1ex]
\psi~=~\delta(f,\phi_C,\phi_D)~:~\ZZ\otimes_AC \to  W[0,k-1]\otimes_{\ZZ[\ZZ_2]}(D \otimes_A D)
\end{array}$$
are such that
$$(1+T)\psi~\simeq ~(F\otimes F)\phi_C - \phi_DF~:~
\ZZ\otimes_AC \to {\rm Hom}_{\ZZ[\ZZ_2]}(W,D \otimes_A D)~.$$
\end{proposition}
\begin{proof}
(i) The short exact sequence of $\ZZ$-module chain complexes
underlying the  proof of Proposition \ref{les3}
$$\begin{array}{l}
0 \to {\rm Hom}_{\ZZ[\ZZ_2]}(W[0,j-1],C\otimes_A C)\\
\hskip50pt \xymatrix{\ar[r]^-{\di{S^k}}&}
{\rm Hom}_{\ZZ[\ZZ_2]}(W[0,j+k-1],S^kC\otimes_AS^kC)_{*+k}\\[1ex]
\hskip75pt \to S(W[0,k-1]\otimes_{\ZZ[\ZZ_2]}(C\otimes_A C))\to 0
\end{array}$$
gives an identification
$$W[0,k-1]\otimes_{\ZZ[\ZZ_2]}(C\otimes_A C)~=~ \Cc(S^k)_{*+1}$$
and a chain equivalence
$$\begin{array}{l}
\Cc(1+T:W[0,k-1]\otimes_{\ZZ[\ZZ_2]}(C\otimes_A C) \to
{\rm Hom}_{\ZZ[\ZZ_2]}(W[0,j-1],C\otimes_A C))\\[1ex]
\hskip100pt \simeq~{\rm Hom}_{\ZZ[\ZZ_2]}(W[0,j+k-1],S^kC\otimes_AS^kC)_{*+k}~.
\end{array}$$
An $n$-cycle $\psi \in (W[0,k-1]\otimes_{\ZZ[\ZZ_2]}(C\otimes_A C))_n$ is thus
essentially the same (up to homology) as an $n$-cycle
$\phi \in {\rm Hom}_{\ZZ[\ZZ_2]}(W[0,j-1],C\otimes_A C)_n$ together with
an $(n+k+1)$-chain
$$\delta\phi \in {\rm Hom}_{\ZZ[\ZZ_2]}(W[0,j+k-1],S^kC\otimes_AS^kC)_{n+k+1}$$
such that $S^k\phi = d(\delta \phi)$, so that up to signs
$$\begin{array}{l}
\phi_r~=~(d\otimes 1 +1\otimes d)\delta \phi_{r+k}+\delta\phi_{r+k-1}+T\delta\phi_{r+k-1}\\[1ex]
\hskip100pt\in (C\otimes_AC)_{n+r}~=~(S^kC\otimes_AS^kC)_{n+2k+r}\\[1ex]
\hskip125pt
(-k \leqslant r \leqslant j-1,~\phi_q=0~{\rm for}~q<0)~,\\[1ex]
\psi_s~=~\delta\phi_{k-s-1} \in
(C\otimes_AC)_{n-s}~=~(S^kC\otimes_AS^kC)_{n+2k-s}~(0 \leqslant s \leqslant k-1)~.
\end{array}$$
and
$$\begin{array}{l}
\phi~=~(1+T)\psi \in
{\rm im}(1+T:Q^{[0,k-1]}_n(C) \to Q_{[0,j-1]}^n(C))\\[1ex]
\hskip100pt =~
{\rm ker}(S^k:Q_{[0,j-1]}^n(C) \to Q_{[0,j+k-1]}^{n+k}(S^kC))~.
\end{array}$$
For $j=k=\infty$ this is just
the short exact sequence of $\ZZ$-module chain complexes
underlying the  proof of Proposition \ref{les3}
$$\begin{array}{l}
0 \to {\rm Hom}_{\ZZ[\ZZ_2]}(W,C\otimes_A C)\\[1ex]
\hskip30pt \xymatrix{\ar[r]^-{\di{J}}&}
\varinjlim\limits_k {\rm Hom}_{\ZZ[\ZZ_2]}(W,S^kC\otimes_AS^kC)_{*+k}=
{\rm Hom}_{\ZZ[\ZZ_2]}(\widehat{W},C\otimes_AC)\\[1ex]
\hskip120pt \xymatrix{\ar[r]^-{\di{H}}&}
S(W\otimes_{\ZZ[\ZZ_2]}(C\otimes_A C))\to 0
\end{array}$$
which gives an identification
$$W\otimes_{\ZZ[\ZZ_2]}(C\otimes_A C)~=~ \Cc(J)_{*+1}$$
and a chain equivalence
$$\begin{array}{l}
\Cc(1+T:W\otimes_{\ZZ[\ZZ_2]}(C\otimes_A C) \to
{\rm Hom}_{\ZZ[\ZZ_2]}(W,C\otimes_A C))\\[1ex]
\hskip100pt \simeq~{\rm Hom}_{\ZZ[\ZZ_2]}(\widehat{W},C\otimes_AC)~.
\end{array}$$
An $n$-cycle $\psi \in (W\otimes_{\ZZ[\ZZ_2]}(C\otimes_A C))_n$ is thus
essentially the same (up to homology) as an $n$-cycle
$\phi \in {\rm Hom}_{\ZZ[\ZZ_2]}(W,C\otimes_A C)_n$ together with
an $(n+k+1)$-chain
$$\delta\phi \in {\rm Hom}_{\ZZ[\ZZ_2]}(W,S^kC\otimes_AS^kC)_{n+k+1}$$
such that $S^k\phi = d(\delta \phi)$,  and
$$\phi~=~(1+T)\psi \in
{\rm im}(1+T:Q^{[0,k-1]}_n(C) \to Q^n(C))~=~{\rm ker}(S^k:Q^n(C) \to \widehat{Q}^n(C))~.$$
(ii) The $A$-module chain map
$$\begin{array}{l}
\theta~=~(F'E_C\otimes F'E_C)(\phi_C)- (E_D\otimes E_D)(\phi_D)F~:\\[1ex]
\hskip100pt C \to {\rm Hom}_{\ZZ[\ZZ_2]}(W,S^{-k}D' \otimes_{\ZZ} S^{-k}D')
\end{array}$$
and the $A$-module chain null-homotopy
$$(F'\otimes F')\delta\phi_C - \delta\phi_D(1\otimes F)~:~
S^k\theta~\simeq~0~:~S^kC \to {\rm Hom}_{\ZZ[\ZZ_2]}(W,D' \otimes_{\ZZ} D')$$
correspond by (i) to an $A$-module chain map
$$\psi'~:~C \to W[0,k-1]\otimes_{\ZZ[\ZZ_2]}(S^{-k}D' \otimes_{\ZZ}S^{-k}D')$$
such that $(1+T)\psi'\simeq \theta$. Use any chain homotopy inverse
$(E_D)^{-1}:S^{-k}D' \to D$ of $E_D$ to define an $A$-module chain map
$$\psi~=~((E_D)^{-1}\otimes (E_D)^{-1})\psi'~:~C \to W[0,k-1]\otimes_{\ZZ[\ZZ_2]}(D \otimes_{\ZZ}D)$$
such that $(1+T)\psi\simeq (F\otimes F)\phi_C - \phi_DF$.\\

\hfill\qed\end{proof}

\begin{remark}{\rm For a bounded f.g. projective $A$-module chain complex $C$ the forgetful map
$$Q^n(C) \to H_n(C\otimes_AC)~;~\phi \mapsto \phi_0$$
sends a symmetric class $\phi \in Q^n(C)$ to a chain homotopy class
$$\phi_0 \in H_n(C\otimes_AC)~=~H_0({\rm Hom}_A(C^{-*},C))$$
of $A$-module chain maps $\phi_0:C^{n-*} \to C$ with
$$d_{C^{n-*}}~=~(-)^rd_C^*~:~
(C^{n-*})_r~=~C^{n-r}\to (C^{n-*})_{r-1}~=~C^{n-r+1}~.$$
\hfill\qed}
\end{remark}

The symmetric $L$-groups $L^n(A)$ of Mishchenko \cite{mishchenko}
are the cobordism groups of
$n$-dimensional symmetric Poincar\'e complexes $(C,\phi)$ over $A$,
as defined by an $n$-dimensional f.g. free $A$-module chain complex $C$
together with a symmetric class $\phi \in Q^n(C)$ such that $\phi_0:C^{n-*} \to C$
is a chain homotopy class of chain equivalences.
The quadratic $L$-groups $L_n(A)$ of Wall \cite{wall2}
were expressed in Ranicki \cite{ranicki1} as the cobordism groups of
$n$-dimensional quadratic Poincar\'e complexes $(C,\psi)$ over $A$,
as defined by an $n$-dimensional f.g. free $A$-module chain complex $C$
together with a quadratic class $\psi \in Q_n(C)$ such that $(1+T)\psi_0:C^{n-*} \to C$
is a chain homotopy class of chain equivalences.
The surgery obstruction of an $n$-dimensional normal map
$(f,b):M \to X$ was expressed in Ranicki \cite{ranicki2} as the cobordism class
$$\sigma_*(f,b)~=~(C,\psi) \in L_n(\ZZ[\pi_1(X)])$$
with $C$ a $\ZZ[\pi_1(X)]$-module chain complex such that
$$H_*(C)~=~K_*(M)~=~{\rm ker}(\widetilde{f}_*:H_*(\widetilde{M}) \to
H_*(\widetilde{X}))~,$$
with $\widetilde{X}$ the universal cover of $X$ and
$\widetilde{M}=f^*\widetilde{X}$ the pullback cover of $M$.  The above
chain level mechanism was used in \cite{ranicki2} to refine the
symmetric structure $\phi \in Q^n(C)$ on the chain complex kernel $C$
of a normal map $(f,b)$ to a quadratic structure $\psi\in Q_n(C)$ using
a stable $\pi_1(X)$-equivariant map $F:\Sigma^k\widetilde{X}^+ \to
\Sigma^k\widetilde{M}^+$ ($k$ large) $S$-dual to the induced
$\pi_1(X)$-equivariant map of Thom spaces
$T(\widetilde{b}):T(\widetilde{\nu}_M) \to T(\widetilde{\nu}_X)$
to provide a chain level solution of $S^k\phi=0 \in Q^{n+k}(S^kC)$.
In the following sections we shall use the geometric Hopf invariant
to construct a homotopy theoretic mechanism inducing this chain level
mechanism. In the first instance we ignore the fundamental group $\pi_1(X)$,
returning to it in Chapter \ref{pi-equivariant}.

\section{The symmetric construction $\phi_V(X)$}\label{sym}

As before, let $C(X)$ denote the singular chain complex of a space $X$, with
$$H_*(C(X))~=~H_*(X)~,~H^*(C(X))~=~H^*(X)$$
the singular homology and cohomology groups of $X$.

\begin{proposition}~ {\rm (Eilenberg-Zilber)}\\
{\rm (i)} For any spaces $X,Y$ there exist
inverse $\ZZ$-module chain equivalences
$$\begin{array}{l}
E(X,Y)~:~C(X) \otimes_{\ZZ} C(Y)\to C(X \times Y)~,\\[1ex]
F(X,Y)~:~C(X \times Y) \to C(X) \otimes_\ZZ C(Y)
\end{array}$$
which are natural in $X,Y$, and such that
$$E(X,Y)(a \otimes b)~=~a\times b~,~F(X,Y)(a \times b)~=~a \otimes b$$
for 0-simplexes $a:\Delta^0 \to X$, $b:\Delta^0 \to Y$.\\
{\rm (ii)} The diagrams
$$\xymatrix@C+10pt{C(X) \otimes_{\ZZ} C(Y)\ar[d]^-{\di{T}}
\ar[r]^-{\di{E(X,Y)}} & C(X \times Y)\ar[d]^-{\di{T}} \\
C(Y) \otimes_{\ZZ} C(X)\ar[r]^-{\di{E(Y,X)}} &C(Y \times X)}
\hskip15pt
\xymatrix@C+10pt{C(X\times Y)\ar[d]^-{\di{T}}
\ar[r]^-{\di{F(X,Y)}} & C(X)\otimes_{\ZZ}C(Y)\ar[d]^-{\di{T}} \\
C(Y\times X)\ar[r]^-{\di{F(Y,X)}} &C(Y) \otimes_{\ZZ} C(X)}
$$
are naturally chain homotopy commutative, with
$$\begin{array}{l}
T~:~X \times Y \to Y \times X~;~(x,y) \mapsto (y,x)~,\\[1ex]
T~:~C(X)_p \otimes_{\ZZ} C(Y)_q \to C(Y)_q \otimes_{\ZZ} C(X)_p~;~
x \otimes y \mapsto (-)^{pq}y \otimes x~.
\end{array}$$
{\rm (iii)} There exist chain homotopy inverse $\R(\infty)$-coefficient
$\ZZ_2$-isovariant chain equivalences
$$\begin{array}{l}
E(X)~=~\{E(X)_s\,\vert\,s \geqslant 0\}~:~C(X)\otimes_{\ZZ}C(X) \to C(X\times X)~,\\[1ex]
F(X)~=~\{F(X)_s\,\vert\,s \geqslant 0\}~:~C(X\times X)\to C(X)\otimes_{\ZZ}C(X)
\end{array}$$
which are natural in $X$, with $E(X)_0=E(X,X)$, $F(X)_0=F(X,X)$.
\end{proposition}
\begin{proof} Standard acyclic model theory.\\
\hfill\qed\end{proof}

\begin{example} {\rm
(i) The diagonal chain approximation for $X$
$$\Delta_{C(X)}~=~F(X)_0\Delta_X~:~
C(X) \xymatrix{\ar[r]^-{\di{\Delta_X}}&}
C(X \times X) \xymatrix{\ar[r]^-{\di{F(X)_0}}&} C(X)\otimes_{\ZZ}C(X)$$
is used to define the cup products
$$\cup ~:~ H^p(X) \times H^q(X) \to H^{p+q}(X)~;~
(x,y) \mapsto (z \mapsto \langle x \otimes y,\Delta_{C(X)}(z) \rangle)~.$$
(ii) The diagonal chain approximation of (i) is the degree 0 component
of an $\R(\infty)$-coefficient $\ZZ_2$-isovariant chain map
$$F(X)\Delta_X~:~C(X)  \xymatrix{\ar[r]^-{\di{\Delta_X}}&}
C(X \times X) \xymatrix{\ar[r]^-{\di{F(X)}}&} C(X)\otimes_{\ZZ}C(X)$$
which is used to define the Steenrod squares --
see Example \ref{steenrod} below.}\\
(iii) The restriction of the
diagonal chain approximation $\Delta_{C(S(LV))}$ to the chain
homotopy deformation retract
$C^{cell}(S(LV)) \subset C(S(LV))$ is chain homotopic
to the cellular diagonal chain approximation $\Delta$ of \ref{iso} (iv).\\
\hfill\qed
\end{example}

\begin{definition}~\label{gamma}
{\rm Let $V$ be an inner product space, and $X$ a space.\\
{\rm (i)} The {\it space level $V$-coefficient symmetric construction} is
the $\ZZ_2$-equivariant map
$$\phi_V(X)~:~S(LV) \times X \to X \times X~;~(v,x) \mapsto (x,x)~.$$
{\rm (ii)} The {\it chain level $V$-coefficient symmetric construction}
is the $V$-coefficient $\ZZ[\ZZ_2]$-module chain map
$$\phi_V(X)~:~C^{cell}(S(LV))\otimes_{\ZZ}C(X)\to C(X)\otimes_{\ZZ} C(X)~,$$
given by the composite of the $\ZZ_2$-isovariant chain maps
$$\begin{array}{l}
C^{cell}(S(LV))\otimes_{\ZZ} C(X) \subset C(S(LV))\otimes_{\ZZ} C(X) \\[1ex]
\hskip90pt
\xymatrix@C+35pt{\ar[r]^-{\di{E(S(LV),X)}}_-{\di{\simeq}}&} C(S(LV)\times X)
\xymatrix{\ar[r]^-{\di{\phi_V(X)}}&} C(X \times X)~,\\[1ex]
C^{cell}(S(LV)) \otimes_{\ZZ}C(X \times X) \subset
C(S(\infty)) \otimes_{\ZZ}C(X \times X) \\[1ex]
\hskip90pt
\xymatrix@C+70pt{\ar[r]^-{\di{F(X)}}_-{\di{\simeq}}&} C(X)\otimes_{\ZZ} C(X)
\end{array}$$
(using an embedding $V \subseteq \R(\infty)$), inducing morphisms
$$\phi_V(X)~:~H_n(X) \to Q^n_V(C(X))~.$$
\hfill\qed}
\end{definition}

\begin{proposition}~\label{delta}
For any inner product space $V$ and spaces $X,Y$ there exists a natural
$\ZZ[\ZZ_2]$-module chain homotopy
$$\begin{array}{l}
\delta\phi_V(X,Y)~:\\[1ex]
\phi_V(X \times Y)E(X,Y)\simeq (E(X,Y) \otimes E(X,Y))(\phi_V(X) \cup \phi_V(Y))~:\\[1ex]
\hskip70pt
C^{cell}(S(LV))\otimes_{\ZZ}C(X)\otimes_{\ZZ}C(Y) \to C(X \times Y)\otimes_{\ZZ}C(X \times Y)
\end{array}$$
in the diagram
$$\xymatrix@C+20pt{
C^{cell}(S(LV))\otimes_{\ZZ}C(X)\otimes_{\ZZ} C(Y)
\ar@{~>}[dr]^-{\delta\phi_V(X,Y)}
\ar[d]^-{1\otimes E(X,Y)}_{\simeq}
\ar[r]^-{\phi_V(X) \cup \phi_V(Y)} &
C(X)\otimes_{\ZZ}C(Y)\otimes_{\ZZ}C(X)\otimes_{\ZZ}C(Y)
\ar[d]^-{E(X,Y) \otimes E(X,Y)}_{\simeq} \\
C^{cell}(S(LV))\otimes_{\ZZ}C(X \times Y)
\ar[r]^-{\phi_V(X \times Y)}& C(X \times  Y)\otimes_{\ZZ}C(X \times Y)}$$
\end{proposition}
\begin{proof} Standard acyclic model theory.\\
\hfill\qed\end{proof}

\begin{terminology}~{\rm For $V=\R(\infty)$ write
$$\begin{array}{l}
\phi(X)~=~\phi_{\R(\infty)}(X)~:~S(\infty) \times X \to X \times X~,\\[1ex]
\phi(X)~=~\phi_{\R(\infty)}(X)~:~W\otimes_{\ZZ} C(X) \to C(X)\otimes_{\ZZ} C(X)~,\\[1ex]
\delta\phi(X,Y)~=~\delta\phi_{\R(\infty)}(X,Y)~:\\[1ex]
\phi(X \times Y)E(X,Y)\simeq (E(X,Y) \otimes E(X,Y))(\phi(X) \cup \phi(Y))~:\\[1ex]
\hskip120pt
C(X)\otimes_{\ZZ}C(Y) \to C(X \times Y)\otimes_{\ZZ}C(X \times Y)~.
\end{array}$$}
\hfill\qed
\end{terminology}

For $V=\R^k$ ($1 \leqslant k \leqslant \infty$) the adjoint of the
natural $\ZZ[\ZZ_2]$-module chain map
$$\phi_V(X)~:~W[0,k-1] \otimes_{\ZZ} C(X)  \to C(X)\otimes_{\ZZ}C(X)$$
is a natural $\ZZ$-module chain map
$$\phi_V(X)~:~C(X) \to {\rm Hom}_{\ZZ[\ZZ_2]}(W[0,k-1],C(X)\otimes C(X))$$
inducing a natural transformation of homology groups
$$\phi_V(X)~:~H_*(X) \to Q^*_V(C(X))~=~Q^*_{[0,k-1]}(C(X))~.$$
By Proposition \ref{delta} there is defined a commutative diagram
$$\xymatrix@C+50pt{
H_m(X) \otimes_{\ZZ} H_n(Y) \ar[r]^-{\di{E(X,Y)}}
\ar[d]_-{\di{\phi_V(X) \cup \phi_V(Y)}} & H_{m+n}(X \times Y)
\ar[d]_-{\di{\phi_V(X \times Y)}}\\
Q^{m+n}_V(C(X)\otimes_{\ZZ}C(Y))
\ar[r]^-{\di{E(X,Y)^{\%}}}_-{\di{\cong}} & Q^{m+n}_V(C(X\times Y))~.
}$$
The $V$-coefficient symmetric constructions $\phi_V(X)$ are the
restrictions of $\phi(X)$
$$\begin{array}{l}
\phi_V(X)~:~S(LV) \times X  \subseteq S(\infty) \times X
\xymatrix{\ar[r]^-{\di{\phi(X)}}&} X \times X~,\\[1ex]
\phi_V(X)~:~W[0,k-1] \otimes_{\ZZ}C(X)  \subseteq W \otimes_{\ZZ}C(X)
\xymatrix{\ar[r]^-{\di{\phi(X)}}&} C(X) \otimes_{\ZZ}C(X)
\end{array}$$
and similarly for $\delta\phi_V(X,Y)$.

\begin{remark}~\label{steenrod} {\rm
(i) The chain level symmetric construction
$$\phi(X)~:~C(X) \to {\rm Hom}_{\ZZ[\ZZ_2]}(W,C(X)\otimes C(X))$$
is the symmetric construction in Ranicki \cite{ranicki2},
with a natural transformation
$$\phi(X)~:~H_*(X) \to Q^*(C(X))~.$$
The components of $\phi(X)$ are natural transformations
$$\phi(X)_s~=~\phi(X)(1_s\otimes -)~:~
C(X)_s \to (C(X)\otimes_{\ZZ}C(X))_{n+s}~~(0 \leqslant s \leqslant k-1)$$
such that
$$\begin{array}{l}
d \phi(X)_s+(-)^r\phi(X)_s d
+(-)^{n+s-1}(\phi(X)_{s-1}+(-)^sT\phi(X)_{s-1})~=~0~:\\[1ex]
C(X)_n \to (C(X) \otimes C(X))_{n+s-1}=
\sum\limits_r C(X)_{n-r+s-1} \otimes C(X)_r~(\phi(X)_{-1}=0)~.
\end{array}$$
(ii) The Steenrod squares of $X$ are given by
$$\begin{array}{l}
Sq^i~:~H^r(X;\ZZ_2) \to H^{r+i}(X;\ZZ_2)~;\\[1ex]
\hskip50pt x \mapsto (y \mapsto\langle x\otimes x,\phi(X)_{r-i}(y) \rangle)~(=0~{\rm for}~r<i)~.
\end{array}
$$
The symmetric Poincar\'e complex associated by Mishchenko
\cite{mishchenko} to an $n$-dimensional geometric Poincar\'e complex $X$  is
$$\sigma^*(X)~=~(C(X),\phi(X)[X]\in Q^n(C(X)))~,$$
with $[X] \in H_n(X)$ the fundamental class.\\
(iii) The $\R(\infty)$-coefficient $\ZZ_2$-isovariant chain homotopy $\delta\phi(X,Y)$
is used in the chain level proof of the Cartan product formula for the Steenrod squares
$$\begin{array}{l}
Sq^k~=~\sum\limits_{i+j=k}Sq^i \otimes Sq^j~:\\[2ex]
H^r(X\times Y;\ZZ_2)~=~\sum\limits_{p+q=r}H^p(X;\ZZ_2)\otimes H^q(Y;\ZZ_2)\\[2ex]
\hskip25pt \to
H^{r+k}(X\times Y;\ZZ_2)~=~\sum\limits_{p+q+i+j=r+k}H^{p+i}(X;\ZZ_2)\otimes H^{q+j}(Y;\ZZ_2)~.
\end{array}$$
(iv) As in Steenrod \cite{steenrod} the {\it Hopf invariant} ${\rm Hopf}(F) \in \ZZ$
of a map $F:S^{2n-1} \to S^n$ is defined by
$${\rm Hopf}(F)~:~H^n(S^n \cup_F D^{2n})~=~\ZZ \to
H^{2n}(S^n \cup_F D^{2n})~=~\ZZ~;~1 \mapsto 1 \cup 1$$
and the {\it mod 2 Hopf invariant} ${\rm Hopf}_2(F) \in \ZZ_2$
of a map $F:S^{n+k} \to S^n$ can be defined by
$$\begin{array}{l}
Sq^{k+1}~=~{\rm Hopf}_2(F)~:\\[1ex]
H^n(S^n \cup_F D^{n+k+1};\ZZ_2)~=~\ZZ_2 \to
H^{n+k+1}(S^n \cup_F D^{n+k+1};\ZZ_2)~=~\ZZ_2~.
\end{array}$$
For $n=k+1$ ${\rm Hopf}_2(F) \in \ZZ_2$ is the mod 2 reduction of
${\rm Hopf}(F) \in \ZZ$.
By Adams \cite{adams1} ${\rm Hopf}(F) \in 2\ZZ\subset \ZZ$ for $n-1 \neq 1,3,7$,
and ${\rm Hopf}_2(F)=0 \in \ZZ_2$ for $k \neq 1,3,7$.
\\
}\hfill\qed
\end{remark}

\begin{definition}~ {\rm
{\rm (i)} The relative singular chain complex of a pair of spaces $(X,A \subseteq X)$
is defined by
$$C(X,A)~=~C(X)/C(A)~,$$
with $H_*(C(X,A))=H_*(X,A)$ the relative singular homology of $(X,A)$.\\
{\rm (ii)} For a pointed space $X$ let
$$\dot C(X)~=~C(X,\{*\})~=~C(X)/C(\{*\})$$
be the reduced singular chain complex, with $H_*(\dot C(X))=\dot H_*(X)$
the reduced singular homology of $X$.\\
{\rm (iii)} For a non-empty subspace $A \subseteq X$ let
$$p_{X/A}~:~C(X,A)~=~C(X)/C(A) \to \dot C(X/A)$$
be the projection.}
\hfill\qed
\end{definition}

In our applications $p_{X/A}$ will be a chain equivalence, e.g.
if $X$ is a $CW$ complex and $A \subset X$ is a subcomplex.

Given pointed spaces $X,Y$ let
$$Z~=~X \times \{*\} \cup \{*\} \times Y \subseteq X \times Y$$
so that
$$(X\times Y)/Z~=~X \wedge Y~,$$
and there are defined $\ZZ$-module chain maps
$$\begin{array}{l}
E(X,Y)~:~\dot C(X)\otimes_{\ZZ}\dot C(Y)~=~
C(X,\{*\}) \otimes_{\ZZ} C(Y,\{*\}) \to C(X\times Y,Z)~,\\[1ex]
\dot E(X,Y)~=~p_{X \wedge Y}E(X,Y)~:~\dot C(X)\otimes_{\ZZ}\dot C(Y) \to
\dot C(X \wedge Y)
\end{array}$$
with $E(X,Y)$ a chain equivalence, and
$p_{X\wedge Y}:C(X \times Y,Z) \to \dot C(X \wedge Y)$ the
projection. As noted above, in our applications $p_{X\wedge Y}$
will be a chain equivalence, e.g. for $CW$ complexes $X,Y$, so that
$\dot E(X,Y)$ will also be a chain equivalence.

Definition \ref{gamma} and Proposition \ref{delta} have pointed versions:

\begin{definition}~ {\rm Let $V$ be an inner product space, and $X$ a pointed
space.\\
{\rm (i)} The {\it reduced space level $V$-coefficient symmetric construction} is
the $\ZZ_2$-equivariant map
$$\dot\phi_V(X)~:~S(LV)^+\wedge X \to X \wedge X~;~(v,x) \mapsto (x,x)~.$$
{\rm (ii)} The {\it chain level $V$-coefficient symmetric construction}
is the $V$-coefficient $\ZZ[\ZZ_2]$-module chain map
$$\dot\phi_V(X)~:~C^{cell}(S(LV))\otimes_{\ZZ}\dot C(X)\to \dot C(X)\otimes_{\ZZ} \dot C(X)~,$$
given by the composite of the $\ZZ_2$-isovariant chain maps
$$\begin{array}{l}
C^{cell}(S(LV))\otimes_{\ZZ} \dot C(X) \subset C(S(LV))\otimes_{\ZZ} \dot C(X) \\[1ex]
\hskip90pt
\xymatrix@C+35pt{\ar[r]^-{\di{E(S(LV),X)}}_-{\di{\simeq}}&} \dot C(S(LV)^+\wedge X)
\xymatrix{\ar[r]^-{\di{\dot\phi_V(X)}}&} \dot C(X \wedge X)~,\\[1ex]
C^{cell}(S(LV)) \otimes_{\ZZ}\dot C(X \wedge X) \subset
C(S(\infty)) \otimes_{\ZZ}\dot C(X \wedge X) \\[1ex]
\hskip90pt
\xymatrix@C+70pt{\ar[r]^-{\di{F(X)}}_-{\di{\simeq}}&} \dot C(X)\otimes_{\ZZ} \dot C(X)
\end{array}$$
(using an embedding $V \subseteq \R(\infty)$), inducing morphisms
$$\dot\phi_V(X)~:~\dot H_n(X) \to Q^n_V(\dot C(X))~.$$
\hfill\qed}
\end{definition}

\begin{proposition}~\label{reduceddelta}
For any inner product space $V$ and pointed spaces $X,Y$ there exists a natural
$\ZZ[\ZZ_2]$-module chain homotopy
$$\begin{array}{l}
\delta\dot\phi_V(X,Y)~:\\[1ex]
\dot\phi_V(X \times Y)E(X,Y)\simeq (E(X,Y) \otimes E(X,Y))(\dot\phi_V(X) \cup \dot\phi_V(Y))~:\\[1ex]
\hskip70pt
C^{cell}(S(LV))\otimes_{\ZZ}\dot C(X)\otimes_{\ZZ}\dot C(Y) \to \dot C(X \wedge Y)\otimes_{\ZZ}\dot C(X \wedge Y)
\end{array}$$
in the diagram
$$\xymatrix@C+20pt{
C^{cell}(S(LV))\otimes_{\ZZ}\dot C(X)\otimes_{\ZZ} \dot C(Y)
\ar@{~>}[dr]^-{\delta\dot\phi_V(X,Y)}
\ar[d]^-{1\otimes E(X,Y)}_{\simeq}
\ar[r]^-{\dot\phi_V(X) \cup \dot\phi_V(Y)} &
\dot C(X)\otimes_{\ZZ}\dot C(Y)\otimes_{\ZZ}\dot C(X)\otimes_{\ZZ}\dot C(Y)
\ar[d]^-{E(X,Y) \otimes E(X,Y)}_{\simeq} \\
C^{cell}(S(LV))\otimes_{\ZZ}\dot C(X \wedge Y)
\ar[r]^-{\dot\phi_V(X \wedge Y)}& \dot C(X \wedge  Y)\otimes_{\ZZ}\dot C(X \wedge Y)}$$
\end{proposition}
\hfill\qed

\begin{example} {\rm
The $\ZZ[\ZZ_2]$-module chain homotopy
$$\begin{array}{l}
\delta\dot\phi_V(S^1,X)~:\\[1ex]
\dot\phi_V(S^1 \wedge X)\dot E(S^1,X)~\simeq~
(\dot E(S^1,X)\otimes \dot E(S^1,X))(\dot \phi_V(S^1) \cup \dot \phi_V(X))~:\\[1ex]
\hskip100pt
\dot C(S^1)\otimes_{\ZZ}\dot C(X) \to \dot C(S^1 \wedge X)\otimes_{\ZZ}\dot C(S^1 \wedge X)
\end{array}$$
for $V=\R(\infty)$
gives the stability of the Steenrod squares, i.e. that the diagrams
$$\xymatrix{\dot H^r(X;\ZZ_2) \ar[d]^-{\di{Sq^i}} \ar[r]_-{\di{\cong}}^-{\di{S}} &
\dot H^{r+1}(\Sigma X;\ZZ_2)\ar[d]^-{\di{Sq^i}} \\
\dot H^{r+i}(X;\ZZ_2) \ar[r]_-{\di{\cong}}^-{\di{S}} & \dot H^{r+i+1}(\Sigma X;\ZZ_2)}$$
commute, with $S$ the suspension isomorphisms induced by the natural chain equivalence
$$S~=~E(S^1,X)\vert~:~
S\dot C(X)~=~\dot C^{cell}(S^1)\otimes_{\ZZ}\dot C(X) \to \dot C(\Sigma X)~=~
\dot C(S^1 \wedge X)~.$$
For $i=r+1$ this is just the vanishing of cup products in a suspension.}\hfill\qed
\end{example}

\begin{example} {\rm
It follows from the stability of the Steenrod squares and from
$$Sq^i~=~0~:~H^r(X;\ZZ_2) \to H^{r+i}(X;\ZZ_2)~\text{for}~i>r$$
that if $X=\Sigma^kY$ is a $k$-fold suspension then
$$Sq^i~=~0~:~\dot H^r(X;\ZZ_2) \to \dot H^{r+i}(X;\ZZ_2)~\text{for}~i>r-k~.$$
\hfill\qed}
\end{example}

Proposition \ref{stablesequence1} for $U=V\subseteq \R(\infty)$ and
$X=Y$ compact gives a commutative braid of exact sequences of stable
$\ZZ_2$-equivariant homotopy groups:

\newpage

\bigskip

$$\xymatrix@C-60pt{
\{X;S(LV)^+ \wedge X \wedge X\}_{\ZZ_2}
\ar[dr]^-{\di{s_{LV}}}\ar@/^2pc/[rr] &&
\{S(LV)^+ \wedge X;X \wedge X\}_{\ZZ_2}\ar[dr]&&&\\
&\{X;X \wedge X\}_{\ZZ_2} \ar[ur]^-{\di{s^*_{LV}}} \ar[dr]^-{\di{0_{LV}}}&&
\{S(LV \oplus LV)^+\wedge X;LV^{\infty} \wedge X\wedge X\}_{\ZZ_2}
\ar[dr] &\\
\{LV^{\infty} \wedge X;X \wedge X\}_{\ZZ_2}
\ar@/_2pc/[rr]\ar[ur]^-{\di{0^*_{LV}}} &&
\{X;LV^{\infty}\wedge X \wedge X\}_{\ZZ_2}\ar[ur]\ar@/_2pc/[rr]&&A
}$$

\bigskip

\noindent with
$$\begin{array}{l}
\Delta_X~=~(0,1)\in
\{X;X \wedge X\}_{\ZZ_2}~=~\{X;S(\infty)^+\wedge(X \wedge X)\}_{\ZZ_2} \oplus \{X;X\}~,\\[1ex]
s_{LV}^*(\Delta_X)~=~\dot\phi_V(X) \in \{S(LV)^+ \wedge X;X \wedge X\}_{\ZZ_2}~,\\[1ex]
[\dot\phi_V(X)]~=~\dot\phi_{V \oplus V}(V^{\infty} \wedge X)
\in \{S(LV \oplus LV)^+\wedge X;LV^{\infty} \wedge X\wedge X\}_{\ZZ_2} \\[1ex]
\hskip50pt =~\{S(LV \oplus LV)^+\wedge (V^{\infty} \wedge X);
(V^{\infty} \wedge X)\wedge (V^{\infty} \wedge X)\}_{\ZZ_2}~,\\[1ex]
0_{LV}(\Delta_X)~=~\Delta_{V^{\infty} \wedge X}~=~(0,1)\\[1ex]
\hskip25pt \in
\{V^{\infty} \wedge X;(V^{\infty} \wedge X) \wedge (V^{\infty} \wedge X)\}_{\ZZ_2}~=~
\{X;LV^{\infty} \wedge (X \wedge X)\}_{\ZZ_2}\\[1ex]
\hskip50pt =~
\{X;S(\infty)/S(LV) \wedge (X \wedge X)\}_{\ZZ_2} \oplus \{X;X\}~,\\[1ex]
A~=~\{S(LV)^+\wedge X;LV^{\infty} \wedge X\wedge X\}_{\ZZ_2}~=~
\{X;\Sigma S(LV)^+ \wedge (X \wedge X)\}_{\ZZ_2}~.
\end{array}$$
If $X=V^{\infty} \wedge X_0$ for some pointed space $X_0$ then
$$\begin{array}{l}
\Delta_{X_0}\in
\{X_0;X_0 \wedge X_0\}_{\ZZ_2}~=~\{LV^{\infty} \wedge X;X \wedge X\}_{\ZZ_2}~,\\[1ex]
\Delta_X~=~0^*_{LV}(\Delta_{X_0}) \in
{\rm im}(0^*_{LV}:\{LV^{\infty} \wedge X;X \wedge X\}_{\ZZ_2} \to
\{X;X \wedge X\}_{\ZZ_2})\\[1ex]
\hskip50pt
=~{\rm ker}(s^*_{LV}:\{X;X \wedge X\}_{\ZZ_2} \to \{S(LV)^+\wedge X;X\wedge X\}_{\ZZ_2})~,\\[1ex]
\dot\phi_V(X)~=~0 \in \{S(LV)^+\wedge X;X\wedge X\}_{\ZZ_2}~,
\end{array}$$
generalizing the vanishing of the cup product in a suspension
(the special case $V=\R$, on the chain level). More precisely~:

\begin{proposition}~ \label{sym4}
{\rm (i)} For any pointed spaces $X,Y$ and inner product space $V$
$$\begin{array}{l}
\dot\phi_V(X \wedge Y)~=~\dot\phi_V(X) \wedge \Delta_Y~=~
\Delta_X \wedge \dot\phi_V(Y)~:\\[1ex]
S(LV)^+ \wedge X \wedge Y \to (X \wedge Y) \wedge (X \wedge Y)~=~
(X \wedge X) \wedge (Y \wedge Y)~.
\end{array}$$
{\rm (ii)} The space level symmetric construction $\dot\phi_V(V^{\infty})$
has a natural $\ZZ_2$-equivariant null-homotopy
$$\delta\dot\phi_V(V^{\infty})~:~\dot\phi_V(V^{\infty})~\simeq~*~:~
S(LV)^+\wedge V^{\infty} \to V^{\infty} \wedge V^{\infty}~.$$
{\rm (iii)} The space level $V$-coefficient symmetric construction $\dot\phi_V(V^{\infty}\wedge X)$
has a natural $\ZZ_2$-equivariant null-homotopy
$$\begin{array}{l}
\delta\dot\phi_V(V^{\infty}\wedge X)~=~\delta\dot\phi_V(V^{\infty})\wedge \Delta_X~:\\[1ex]
\hskip25pt
\dot\phi_V(V^{\infty}\wedge X)~\simeq~*~:~S(LV)^+\wedge V^{\infty}\wedge X \to V^{\infty}\wedge X \wedge V^{\infty}\wedge X
\end{array}$$
inducing a natural $\ZZ[\ZZ_2]$-module null-chain homotopy
$$\begin{array}{l}
\delta\dot\phi_V(V^{\infty}\wedge X)~:~\dot\phi_V(V^{\infty}\wedge X)~\simeq~0~:\\[1ex]
\hskip50pt
C^{cell}(S(LV))\otimes_{\ZZ} \dot C(V^{\infty}\wedge X) \to \dot C(V^{\infty}\wedge X)\otimes_{\ZZ}\dot C(V^{\infty}\wedge X)~.
\end{array}$$
\end{proposition}
\begin{proof} (i) By construction.\\
(ii) Apply the construction of \ref{difcon6} to the $\ZZ_2$-equivariant map
$$\begin{array}{l}
p~=~\kappa_V~:~LV^{\infty} \wedge V^{\infty}
\to V^{\infty}\wedge V^{\infty}~;~(r,s) \mapsto (r+s,-r+s)
\end{array}$$
to obtain a $\ZZ_2$-equivariant null-homotopy
$$\delta\dot\phi_V(V^{\infty})~=~\delta p~:~CS(LV)^+\wedge V^{\infty} \to
V^{\infty} \wedge V^{\infty}~;~(t,u,v) \mapsto p([t,u],v)$$
of
$$\delta\dot\phi_V(V^{\infty})\vert~=~\dot\phi_V(V^{\infty})~:~S(LV)^+ \wedge V^{\infty} \to V^{\infty} \wedge V^{\infty}~;~
(u,v) \mapsto (v,v)~,$$
using the projection
$$CS(LV)^+ \to LV^{\infty}~;~(t,u) \mapsto [t,u]~=~{tu \over 1-t}~.$$
(iii) Combine (i) and (ii).
\hfill\qed\end{proof}

More generally~:

\begin{proposition}~ \label{sym44}
Let $U,V$ be inner product spaces.\\
{\rm (i)} There is defined a $\ZZ_2$-equivariant homotopy
$$\begin{array}{l}
\delta\dot\phi_{U,V}~:~
\kappa_V\dot\phi_{U \oplus V}(V^{\infty})~\simeq~
(\dot\phi_U(S^0) \wedge 1_{LV^{\infty} \wedge V^{\infty}})(j_{LU,LV} \wedge 1_{V^{\infty}})\\[1ex]
\hskip75pt :~S(LU\oplus LV)^+\wedge V^{\infty} \to LV^{\infty} \wedge V^{\infty}~.
\end{array}$$
in the diagram
$$\xymatrix@C+40pt{
S(LU \oplus LV)^+ \wedge V^{\infty} \ar[r]^-{j_{LU,LV}\wedge 1_{V^{\infty}}}
\ar[d]_-{\dot\phi_{U \oplus V}(V^{\infty})} &
S(LU)^+\wedge LV^{\infty} \wedge V^{\infty}
\ar[d]^-{\dot\phi_U(S^0)\wedge  1_{ LV^{\infty}\wedge V^{\infty} } } \\
V^{\infty} \wedge V^{\infty} \ar[r]^-{\kappa_V}
&LV^{\infty} \wedge V^{\infty} }$$
with
$$\begin{array}{l}
j_{LU,LV}~=~{\rm projection}~:\\[1ex]
S(LU \oplus LV)^+~=~(D(LU) \times S(LV) \cup S(LU)\times D(LV))^+\\[1ex]
\to S(LU) \times D(LV)/S(LU) \times S(LV)~=~S(LU)^+ \wedge D(LV)/S(LV)\\[1ex]
\hphantom{\to S(LU) \times D(LV)/S(LU) \times S(LV)~}=~S(LU)^+\wedge LV^{\infty}~.
\end{array}$$
{\rm (ii)} For any pointed space $X$ there is defined a natural
$\ZZ_2$-equivariant homotopy
$$\begin{array}{l}
\delta\dot\phi_{U,V}(X)~:\\[1ex]
(\kappa_V\wedge 1_{X \wedge X})\dot\phi_{U \oplus V}(V^{\infty}\wedge X)~\simeq~
(\dot\phi_U(X) \wedge 1_{LV^{\infty} \wedge V^{\infty}})(j_{LU,LV} \wedge 1_{V^{\infty}})\\[1ex]
\hskip100pt :~S(LU\oplus LV)^+\wedge V^{\infty}\wedge X \to
LV^{\infty} \wedge V^{\infty}\wedge X\wedge X~.
\end{array}$$
in the diagram
$$\xymatrix@C+40pt{
S(LU \oplus LV)^+ \wedge V^{\infty}\wedge X \ar[r]^-{j_{LU,LV} \wedge 1_{V^{\infty}\wedge X}}
\ar[d]_-{\dot\phi_{U \oplus V}(V^{\infty}\wedge X)} &
S(LU)^+\wedge LV^{\infty} \wedge V^{\infty}\wedge X
\ar[d]^-{\dot\phi_U(X)\wedge  1_{ LV^{\infty}\wedge V^{\infty} } } \\
V^{\infty}\wedge V^{\infty}\wedge X \wedge  X \ar[r]^-{\kappa_V\wedge 1_{X \wedge X}}
&LV^{\infty} \wedge V^{\infty}\wedge X\wedge X }$$
{\rm (iii)} The natural $\ZZ_2$-equivariant homotopy in {\rm (ii)}
induces a natural $U \oplus V$-coefficient $\ZZ_2$-isovariant
chain homotopy
$$\begin{array}{l}
\delta\dot\phi_{U,V}~:~\dot\phi_{U \oplus V}(V^{\infty}\wedge X)~\simeq~
\dot \phi_U(X)j_{LU,LV}\\[1ex]
\hskip100pt :~\dot C(V^{\infty}\wedge X) \to \dot C(V^{\infty}\wedge X) \otimes_{\ZZ}
\dot C(V^{\infty}\wedge X)~.
\end{array}$$
\end{proposition}
\begin{proof} (i) Use the pushout square of $\ZZ_2$-spaces given
by Proposition \ref{pushout2} (ii)
$$\xymatrix{S(LU \oplus LV)^+ \ar[r]^-{\di{j_{LU,LV}}} \ar[d] & S(LU)^+ \wedge LV^{\infty}
\ar[d]^-{\di{k_{LU,LV}}} \\
CS(LU \oplus LV)^+ \ar[r]^-{\di{\delta j_{LU,LV}}} & \Sigma S(LV)^+}$$
to define
$$\delta\dot\phi_{U,V}~:~I \times
S(LU\oplus LV)^+\wedge V^{\infty} \to LV^{\infty} \wedge V^{\infty}~;~
(t,(u,v,w)) \mapsto (tv,w)~.$$
(For $U=\{0\}$  $\delta\dot\phi_{U,V}=\delta\dot\phi_V$ was already
defined in the proof of Proposition \ref{sym4} (ii)).\\
(ii)+(iii) Immediate from (i), with $\delta\dot\phi_{U,V}(X)=\delta\dot\phi_{U,V}\wedge \Delta_X$.
\hfill\qed\end{proof}

We shall also need a relative version of the symmetric construction:

\begin{definition}~ \label{relsym1}
{\rm Given a pair of pointed spaces $(A,B \subseteq A)$ there is
defined a {\it relative space level reduced $V$-coefficient
symmetric construction} $\ZZ_2$-equivariant map
$$\dot\phi_V(A,B)~:~S(LV)^+\wedge A/B \to (A \wedge A)/(B \wedge B)$$
which fits into a natural transformation of homotopy cofibration sequences of
$\ZZ_2$-spaces
$$\xymatrix@C-5pt@R+10pt{
S(LV)^+\wedge B \ar[r]
\ar[d]_-{\dot\phi_V(B)}&
S(LV)^+\wedge A \ar[r]
\ar[d]_-{\dot\phi_V(A)}&
S(LV)^+ \wedge A/B  \ar[r]
\ar[d]_-{\dot\phi_V(A,B)}& S(LV)^+\wedge \Sigma B
\ar[d]_-{\Sigma\dot\phi_V(B)} \\
B \wedge B \ar[r] &
A \wedge A \ar[r] &
(A \wedge A)/(B \wedge B) \ar[r]&\Sigma(B\wedge B)}$$}
\hfill\qed
\end{definition}

\begin{proposition}~ \label{relsym2}
Given a pointed map $F:X \to Y$ let
$$(A,B \subseteq A)~=~(\Mm(F),\{1\} \times X)~,~Z~=~A/B~=~\Cc(F)$$
so that there is defined a homotopy cofibration sequence of pointed maps
$$\xymatrix@C+10pt{X \ar[r]^-{\di{F}} &  Y \ar[r]^-{\di{G}}&
Z\ar[r]^-{\di{H}}& \Sigma X\ar[r]^-{\di{\Sigma F}} &\Sigma Y \ar[r]&\dots}$$
with $G$ the inclusion, $H$ the projection. The
space level reduced $V$-coefficient symmetric construction of $Z$
is given by
$$\begin{array}{l}
\xymatrix@C+20pt{\dot \phi_V(Z)~:~S(LV)^+ \wedge Z
\ar[r]^-{\di{\dot\phi_V(A,B)}}&}\\[1ex]
\hskip50pt
\xymatrix{A \wedge A/B \wedge B \ar[r] &
A\times A/(A \times B \cup B \times A)~=~Z \wedge Z}
\end{array}$$
\end{proposition}
\begin{proof}
By construction.\\
\hfill\qed\end{proof}

\section{The geometric Hopf invariant $h_V(F)$}\label{hopf}

Let $G:V^{\infty} \wedge X \to V^{\infty} \wedge Y$ be a stable
$\ZZ_2$-map, with $X,Y$ pointed $\ZZ_2$-spaces and $V$ an inner product $\ZZ_2$-space, with
$$V~=~V_- \oplus V_+~,~V_{\pm}~=~\{v \in V\st Tv=\pm v\}~.$$
(We shall be mainly concerned with the case when $G$ is the square
of a stable nonequivariant map $F$).
The geometric Hopf invariant uses the difference construction to
associate to $G$ a $\ZZ_2$-map $\delta(p,q)$ whose $\ZZ_2$-homotopy class
is the primary obstruction to the desuspension of $G$. The diagram
$$\xymatrix@R+20pt@C+40pt{
V^{\infty} \wedge X^{\ZZ_2} \ar[d]_-{\di{1 \wedge \rho(G)}}
\ar[r]^-{\di{1 \wedge i_X}} &
V^{\infty} \wedge X
\ar[d]^-{\di{G}}\\
V^{\infty} \wedge Y^{\ZZ_2}\ar[r]^-{\di{1\wedge i_Y}}
&V^{\infty} \wedge Y}$$
does not commute in general, with $i_X:X^{\ZZ_2} \to X$, $i_Y:Y^{\ZZ_2} \to Y$ the inclusions
and $\rho(G):V_+^{\infty} \wedge X^{\ZZ_2} \to V_+^{\infty} \wedge Y^{\ZZ_2}$ the restriction of $G$ to the fixed point sets.
However, the $\ZZ_2$-maps
$$p~=~G(1 \wedge i_X)~,~q~=~(1 \wedge i_Y) (1 \wedge \rho(G))~:~
V^{\infty} \wedge X^{\ZZ_2} \to V^{\infty} \wedge Y$$
are such that
$$p(0,v,x)~=~G(0,v,x)~=~(0,v,\rho(G)(x))~=~q(0,v,x)~~(v \in V_+,x \in X^{\ZZ_2})$$
so the relative difference (\ref{difcon1} (ii)) $\ZZ_2$-map
$$\delta(p,q)~:~\Sigma S(V_-)^+ \wedge V_+^{\infty} \wedge X^{\ZZ_2}
\to V^{\infty} \wedge Y$$
is defined, with
$$p-q~:~V^{\infty} \wedge X^{\ZZ_2}
\xymatrix{\ar[r]^-{\di{\alpha_{V_-}\wedge i_X}}&}
\Sigma S(V_-)^+ \wedge V_+^{\infty} \wedge X
\xymatrix{\ar[r]^-{\di{\delta(p,q)}}&} V^{\infty} \wedge Y~.$$
If $G$ is $\ZZ_2$-homotopic to
$1 \wedge G_0: V^{\infty} \wedge X \to V^{\infty} \wedge Y$
for some $\ZZ_2$-map $G_0:X \to Y$ then
$$\delta(p,q) ~\simeq~\delta(1\wedge G_0,1 \wedge G_0)~\simeq~*~.$$
\indent Now let $X,Y$ be spaces, and $V$ an inner product space.
The square of a nonequivariant stable map $F:V^{\infty} \wedge X \to V^{\infty} \wedge Y$
is a stable $\ZZ_2$-map
$$F \wedge F~:~ V^{\infty} \wedge V^{\infty} \wedge X \wedge X
\to V^{\infty} \wedge V^{\infty} \wedge Y \wedge Y~.$$
The geometric Hopf invariant of $F$ is the above difference construction
$$h_V(F)~=~\delta(p,q)~:~ \Sigma S(LV)^+ \wedge V^{\infty} \wedge X \to
LV^{\infty}\wedge V^{\infty} \wedge Y\wedge Y$$
applied to the stable $\ZZ_2$-map
$$\begin{array}{l}
G~=~(\kappa^{-1}_V \wedge 1)(F \wedge F)(\kappa_V \wedge 1)~:\\[1ex]
\hskip50pt
LV^{\infty} \wedge V^{\infty} \wedge X \wedge X
\to LV^{\infty} \wedge V^{\infty} \wedge Y \wedge Y~;\\[1ex]
\hskip75pt (u,v,x_1,x_2) \mapsto ((w_1-w_2)/2,(w_1+w_2)/2,y_1,y_2)\\[1ex]
\hskip75pt (F(u+v,x_1)=(w_1,y_1),F(-u+v,x_2)=(w_2,y_2))
\end{array}$$
with
$$\rho(G)~=~F~:~V^{\infty} \wedge X \to V^{\infty} \wedge Y~.$$
In essence, the geometric Hopf invariant measures the difference
$(F\wedge F)\Delta_X - \Delta_YF$, given that
$(F\wedge F)\Delta_{V^{\infty} \wedge X}=\Delta_{V^{\infty} \wedge Y}F$.
The diagram of $\ZZ_2$-equivariant maps
$$\xymatrix@R+20pt@C+40pt{
LV^{\infty} \wedge V^{\infty} \wedge X \ar[r]^-{\di{1 \wedge \Delta_X}}
\ar[d]_-{\di{1 \wedge F}} &
LV^{\infty} \wedge V^{\infty} \wedge X \wedge X
\ar[d]^-{\di{G}}\\
LV^{\infty}\wedge V^{\infty} \wedge Y
\ar[r]^-{\di{1 \wedge \Delta_Y}}
&LV^{\infty} \wedge V^{\infty} \wedge Y \wedge Y}$$
does not commute in general.
However, the $\ZZ_2$-equivariant maps defined by
$$\begin{array}{l}
p~=~G(1\wedge \Delta_X)~:~LV^{\infty} \wedge V^{\infty} \wedge X \to LV^{\infty} \wedge V^{\infty}\wedge Y \wedge Y~;\\[1ex]
\hskip75pt (u,v,x) \mapsto ((w_1-w_2)/2,(w_1+w_2)/2,y_1,y_2)\\[1ex]
\hskip75pt (F(u+v,x)=(w_1,y_1),F(-u+v,x)=(w_2,y_2))~,\\[1ex]
q~=~(1 \wedge \Delta_Y)(1 \wedge F)~:~
LV^{\infty} \wedge V^{\infty} \wedge X \to
LV^{\infty} \wedge V^{\infty} \wedge Y \wedge Y~;\\[1ex]
\hskip75pt (u,v,x) \mapsto (u,w,y,y)~~(F(v,x)=(w,y))
\end{array}$$
agree on $0^+ \wedge V^{\infty} \wedge X =V^{\infty} \wedge X \subset
LV^{\infty} \wedge V^{\infty} \wedge X$, with
$$\begin{array}{l}
p\vert~=~q\vert~=~(\kappa^{-1}_V \wedge 1)\Delta_{V^{\infty} \wedge Y}
F~=~(\kappa^{-1}_V \wedge 1)(F \wedge F)(\kappa_V \wedge 1)
\Delta_{V^{\infty} \wedge X}~:\\[1ex]
\hskip100pt V^{\infty} \wedge X
    \to LV^{\infty} \wedge V^{\infty} \wedge Y\wedge Y
    \end{array}$$
on account of the commutative diagram
$$\xymatrix@C+30pt@R+20pt{
V^{\infty} \wedge X \ar[r]^-{\di{\Delta_{V^{\infty} \wedge X}}}
\ar[d]_{\di{F}}
&(V^{\infty} \wedge X )\wedge (V^{\infty} \wedge X)
\ar[d]_{\di{F\wedge F}} \ar[r]^-{\di{\kappa^{-1}_V\wedge 1}}&
LV^{\infty} \wedge V^{\infty} \wedge X \wedge X
\ar[d]_-{\di{(\kappa^{-1}_V \wedge 1)(F \wedge F)(\kappa_V \wedge 1)}}\\
V^{\infty} \wedge Y\ar[r]^-{\di{\Delta_{V^{\infty} \wedge Y}}}&
(V^{\infty} \wedge Y )\wedge (V^{\infty} \wedge Y)
\ar[r]^-{\di{\kappa^{-1}_V\wedge 1}}
&LV^{\infty} \wedge V^{\infty} \wedge Y \wedge Y~.}$$
\begin{definition}~ \label{hopf1}
{\rm  The {\it geometric Hopf invariant} of a map $F:V^{\infty} \wedge X \to
V^{\infty} \wedge Y$ is the $\ZZ_2$-equivariant map given by the relative
difference of the $\ZZ_2$-equivariant maps
$p=(\kappa^{-1}_V\wedge 1) (F \wedge F)(\kappa_V \wedge \Delta_X)$,
$q=(1 \wedge \Delta_Y)(1 \wedge F)$
\index{geometric Hopf invariant!$h_V(F)$}
$$\begin{array}{ll}
h_V(F)~=~\delta(p,q)~:&\Sigma S(LV)^+ \wedge V^{\infty} \wedge X \to
LV^{\infty} \wedge V^{\infty} \wedge Y\wedge Y~;\\[1ex]
&(t,u,v,x) \mapsto \begin{cases}
q([ 1-2t,u],v,x)&\hbox{if $0 \leqslant  t \leqslant  1/2$}\\[1ex]
p([ 2t-1,u],v,x)&\hbox{if $1/2 \leqslant  t
\leqslant  1$}
\end{cases}\\[3ex]
&(t \in I,u \in S(LV),v \in V,x \in X)~.
\end{array}$$
\hfill\qed
}\end{definition}

\begin{remark} \label{hopf2}
{\rm $h_V(F)$ is the geometric Hopf invariant of
Crabb \cite[p.\ 61]{crabb}, Crabb and James \cite[p.\ 306]{crabbjames}.
\hfill\qed}
\end{remark}

\begin{example} {\rm Suppose that $V=\{0\}$, so that
$$V^{\infty}~=~S^0~,~S(LV)~=~\emptyset~,~S(LV)^+~=~\{*\}~.$$
The geometric Hopf invariant of a map
$$F~:~V^{\infty} \wedge X ~=~X \to V^{\infty} \wedge Y~=~Y$$
is
$$h_V(F)~=~*~:~\Sigma S(LV)^+ \wedge V^{\infty} \wedge X ~=~\{*\}\to
LV^{\infty} \wedge V^{\infty} \wedge Y \wedge Y~=~Y \wedge Y~.$$
\hfill\qed}
\end{example}

The geometric Hopf invariant $h_V(F)$ of a map
$F:V^{\infty} \wedge X \to V^{\infty} \wedge Y$ has the following properties.

\begin{proposition}~ \label{hopf7}{\it
{\rm (i)} {\rm (}Naturality{\rm )} If $f:X \to X'$, $g:Y \to Y'$,
$F:V^{\infty} \wedge X \to V^{\infty} \wedge Y$,
$F':V^{\infty} \wedge X' \to V^{\infty} \wedge Y'$ are maps
such that there is defined a commutative diagram
$$\xymatrix{
V^{\infty} \wedge X \ar[r]^-{\di{F}} \ar[d]_-{\di{1\wedge f}} &
V^{\infty} \wedge Y  \ar[d]^-{\di{1\wedge g}} \\
V^{\infty} \wedge X' \ar[r]^-{\di{F'}} & V^{\infty} \wedge Y'}$$
then there is defined a commutative diagram
$$\xymatrix@C+15pt@R+5pt{
\Sigma S(LV)^+\wedge V^{\infty} \wedge X \ar[r]^-{\di{h_V(F)}}
\ar[d]_-{\di{1\wedge f}}
& LV^{\infty} \wedge V^{\infty} \wedge Y \wedge Y  \ar[d]^-{\di{1\wedge g \wedge g}} \\
\Sigma S(LV)^+\wedge V^{\infty} \wedge X' \ar[r]^-{\di{h_V(F')}} &
LV^{\infty} \wedge V^{\infty} \wedge Y' \wedge Y'}$$
{\rm (ii)} {\rm (}Homotopy invariance{\rm )}
The $\ZZ_2$-equivariant homotopy class of
$h_V(F)$ depends only on the homotopy class of $F$.\\
{\rm (iii)} {\rm (}Suspension formula{\rm )}
Let $U$ be another inner product space.
The geometric Hopf invariant $h_{U\oplus V}(1_{U^{\infty}}\wedge F)$ of the
stabilization of $F:V^{\infty} \wedge X \to V^{\infty} \wedge Y$
$$1_{U^{\infty}} \wedge F~:~(U \oplus V)^{\infty} \wedge X ~=~
U^{\infty} \wedge V^{\infty} \wedge X
\to (U \oplus V)^{\infty} \wedge Y ~=~U^{\infty} \wedge V^{\infty} \wedge Y$$
is determined by $h_V(F)$, with a commutative diagram
$$\xymatrix@R+10pt@C-45pt
{U^{\infty} \wedge V^{\infty} \wedge \Sigma S(LU \oplus LV)^+ \wedge X
 \ar[dr]^-{\di{1 \wedge \Sigma j_{LV,LU}}}
\ar[dd]_-{\di{~~~~~~~h_{U\oplus V}(1_{U^{\infty}}\wedge F)}} &\\
&\ar[dl]^-{\di{~~~~~~~1_{(U\oplus LU)^{\infty}}\wedge h_V(F)}}
(U \oplus LU)^{\infty} \wedge V^{\infty} \wedge \Sigma S(LV)^+\wedge X\\
(U \oplus LU)^{\infty} \wedge (V \oplus LV)^{\infty} \wedge Y \wedge Y&}$$
with $j_{LV,LU}$ the $\ZZ_2$-equivariant adjunction
Umkehr map of the  $\ZZ_2$-equivariant embedding
$$S(LU) \times D(LV) \subset
D(LU)\times S(LV)\cup S(LU)\times D(LV)~=~S(LU \oplus LV)~,$$
that is
$$\begin{array}{l}
j_{LV,LU}~=~projection~:~S(LU \oplus LV)^+ \to\\[1ex]
\hskip50pt S(LU \oplus LV)/(S(LU) \times D(LV))~=~
LU^{\infty}\wedge S(LV)^+~.
\end{array}$$
{\rm (iv)} {\rm (}Symmetrization{\rm )}
Write
$$A~=~V^{\infty} \wedge X~,~B~=~LV^{\infty} \wedge V^{\infty} \wedge Y \wedge Y~.$$
The images of the $\ZZ_2$-equivariant homotopy class of the geometric
Hopf invariant $h_V(F) \in [\Sigma S(LV)^+ \wedge A;B]_{\ZZ_2}$
under the maps in the commutative braid of exact sequences of
$\ZZ_2$-equivariant homotopy groups
$$\xymatrix@C-60pt{
[\Sigma S(LV)^+ \wedge A;B]_{\ZZ_2}
\ar[dr]^-{\di{\alpha^*_{LV}}}\ar@/^2pc/[rr]^-{\di{k^*_{LV,LV}}} &&
[S(LV)^+ \wedge LV^{\infty} \wedge A;B]_{\ZZ_2}
\ar[dr]^-{\di{j_{LV,LV}^*}}&&\\
&[LV^{\infty} \wedge A;B]_{\ZZ_2}
\ar[ur]^-{\di{s^*_{LV}}}
\ar[dr]^-{\di{0^*_{LV}}}&&
[S(LV \oplus LV)^+\wedge A;B]_{\ZZ_2}\\
[LV^{\infty} \wedge LV^{\infty} \wedge A;B]_{\ZZ_2}
\ar@/_2pc/[rr]^-{\di{0^*_{LV\oplus LV}}}\ar[ur]^-{\di{0^*_{LV}}} &&
[A;B]_{\ZZ_2}\ar[ur]^-{\di{s^*_{LV\oplus LV}}}&
}$$
are given by
$$\begin{array}{l}
\alpha^*_{LV}h_V(F)~=~(F \wedge F)\Delta_X -\Delta_Y (1 \wedge F)
\in [LV^{\infty} \wedge A;B]_{\ZZ_2}~,\\[1ex]
k^*_{LV,LV}h_V(F)~=~(F \wedge F)\dot\phi_V(X) -
\dot\phi_V(Y) F\in [S(LV)^+\wedge LV^{\infty} \wedge A;B]_{\ZZ_2}~.
\end{array}$$
In particular, the diagram of $\ZZ_2$-equivariant maps
$$\xymatrix@C+30pt@R+20pt
{LV^{\infty} \wedge V^{\infty} \wedge X \ar[r]^-{\di{\alpha_{LV} \wedge 1}}
\ar[dr]_-{\di{(F \wedge F)\Delta_X -\Delta_Y (1\wedge F)~~~~~~~~~~}} &
\ar[d]^-{\di{h_V(F)}}\Sigma S(LV)^+ \wedge V^{\infty} \wedge X\\
& LV^{\infty}\wedge V^{\infty} \wedge Y \wedge Y}$$
commutes up to $\ZZ_2$-equivariant homotopy.\\
{\rm (v)} {\rm (}Composition formula{\rm )}
The geometric Hopf invariant of the composite
$$GF~:~V^{\infty} \wedge X \to V^{\infty}\wedge Z$$
of maps $F:V^{\infty} \wedge X \to V^{\infty} \wedge Y$,
$G:V^{\infty} \wedge Y \to V^{\infty}\wedge Z$
is given up to $\ZZ_2$-equivariant homotopy by
$$\begin{array}{ll}
&h_V(GF)~=~h_V(G)(1\wedge F)+(\kappa^{-1}_V\wedge 1)(G \wedge G)(\kappa_V\wedge 1)h_V(F)~:\\[1ex]
&\hskip50pt
\Sigma S(LV)^+ \wedge V^{\infty} \wedge X \to
LV^{\infty} \wedge V^{\infty} \wedge Z \wedge Z~.
\end{array}$$
{\rm (vi)} {\rm (}Product formula{\rm )}
Let $F_i :V_i^{\infty}\wedge X_i\to V_i^{\infty}\wedge Y_i$, $i=1,\, 2$,
be two maps.  The geometric Hopf invariant $h_{V_1\oplus V_2} (F_1\wedge F_2)$
of the smash product
$$F_1\wedge F_2~:~(V_1\oplus V_2)^{\infty} \wedge (X_1\wedge X_2)
\to (V_1\oplus V_2)^{\infty} \wedge (Y_1\wedge Y_2)$$
is the composite
$$\begin{array}{ll}
&\big((\Delta_{Y_1}F_1\wedge h_{V_2}(F_2))\vee (h_{V_1}(F_1)\wedge h_{V_2}(F_2))
\vee (h_{V_1}(F_1)\wedge \Delta_{Y_2}F_2)\big) \circ \Sigma c~:\\[1ex]
&\hskip25pt \Sigma S(LV_1 \oplus LV_2)^+ \wedge
(V_1 \oplus V_2)^{\infty} \wedge X_1 \wedge X_2\\[1ex]
&\hskip50pt \to ((LV_1\oplus LV_2)^{\infty} \wedge
(V_1\oplus V_2)^{\infty}\wedge (Y_1\wedge Y_2) \wedge(Y_1\wedge Y_2))
\end{array}$$
where
$$c~:~A~=~S(LV_1\oplus LV_2)^+ \to B_1 \vee B_2 \vee B_3$$
is the Umkehr map for the inclusion of the submanifold
$$S(LV_1) \sqcup S(LV_2)\sqcup (S(LV_1)\times S(LV_2))
\emb S(LV_1\oplus LV_2)~.$$
The $\ZZ_2$-equivariant stable homotopy class
$$c~=~(c_1,c_2,c_3) \in \{A;B_1\vee B_2\vee B_3\}_{\ZZ_2}~=~
\{A;B_1\}_{\ZZ_2} \oplus\{A;B_2\}_{\ZZ_2} \oplus\{A;B_3\}_{\ZZ_2}$$
has components
$$\begin{array}{l}
c_1~=~j_{LV_2,LV_1}~:~
A~=~S(LV_1 \oplus LV_2)\\[1ex]
\hskip20pt  \to B_1~=~S(LV_1 \oplus LV_2)/(S(LV_1)\times D(LV_2))~=~
LV_1^{\infty} \wedge S(LV_2)^+~,\\[1ex]
c_2~=~j_{LV_1,LV_2}~:~
A~=~S(LV_1 \oplus LV_2)\\[1ex]
\hskip20pt \to B_2~=~S(LV_1 \oplus LV_2)/(D(LV_1)\times S(LV_2))~=~
S(LV_1)^+ \wedge LV_2^{\infty}~,\\[1ex]
c_3~=~projection~:~A~=~S(LV_1 \oplus LV_2)\\[1ex]
\hskip20pt \to B_3~=~S(LV_1 \oplus LV_2)/(S(LV_1)\times D(LV_2) \sqcup
D(LV_1)\times S(LV_2))\\[1ex]
\hskip150pt =~\Sigma (S(LV_1)^+\wedge S(LV_2)^+)~.\end{array}$$
{\rm (vii)} {\rm (}One-point union formula{\rm )}
The geometric Hopf invariant of a map from a one-point union
$$F_1\vee F_2~:~V^{\infty}\wedge (X_1\vee X_2) ~=~
(V^{\infty}\wedge X_1)\vee (V^{\infty}\wedge X_2)\to V^{\infty}\wedge Y$$
is given by
$$\begin{array}{l}
h_V(F_1\vee F_2)~=~h_V(F_1)\vee h_V(F_2)~:\\[1ex]
\Sigma S(LV)^+ \wedge V^{\infty} \wedge
(X_1 \vee X_2)\to LV^{\infty} \wedge V^{\infty} \wedge Y \wedge  Y~.
\end{array}$$
{\rm (viii)} {\rm (}Sum map{\rm )}
For $V \neq \{0\}$ the geometric Hopf invariant of the sum map
$$\nabla_V~:~V^{\infty} \to V^{\infty} \vee V^{\infty}$$
is given up to $\ZZ_2$-equivariant homotopy by the composite
$$h_V(\nabla_V)~=~i_Vj_V~:~\Sigma S(LV)^+ \wedge V^{\infty}
\to (V^{\infty} \vee V^{\infty}) \wedge (V^{\infty} \vee V^{\infty})$$
with
$$\begin{array}{l}
i_V~:~(V^{\infty} \wedge V^{\infty}) \vee (V^{\infty} \wedge V^{\infty}) \to
(V^{\infty} \vee V^{\infty}) \wedge (V^{\infty} \vee V^{\infty})~;\\[1ex]
\hskip50pt
(u,v)_1 \mapsto (u_1,v_2)~~,~~(u,v)_2 \mapsto (u_2,v_1)\end{array}$$
and
$$\begin{array}{l}j_V~:~
\Sigma S(LV)^+ \wedge V^{\infty} \to (V^{\infty} \wedge
V^{\infty}) \vee (V^{\infty} \wedge V^{\infty})~;\\[1ex]
\hskip40pt(t,u,v) \mapsto \begin{cases}(\alpha_V(2t,u),v)_1&
\hbox{if $0
\leqslant  t \leqslant  1/2$} \\[1ex] (v,\alpha_V(2t-1,u))_2& \hbox{if $1/2
\leqslant  t \leqslant  1$~.}
\end{cases}
\end{array}$$
{\rm (ix)} {\rm (}Sum formula{\rm )}
For $V \neq \{0\}$ the geometric Hopf invariant of the sum of maps
$F,G:V^{\infty} \wedge X \to V^{\infty} \wedge Y$
$$F+G~:~V^{\infty} \wedge X \to V^{\infty}\wedge Y$$
is given up to $\ZZ_2$-equivariant homotopy by
$$\begin{array}{l}
h_V(F+G)~=~(h_V(F) \vee h_V(G) \vee
((F \wedge G)\vee (G \wedge F))(1\wedge \Delta_X))\circ (d \wedge 1)~:\\[1ex]
\hskip100pt
\Sigma S(LV)^+ \wedge V^{\infty} \wedge X \to
LV^{\infty} \wedge V^{\infty} \wedge Y \wedge Y\end{array}$$
with
$$\begin{array}{l}
d~:~A~=~\Sigma S(LV)^+ \wedge V^{\infty} \to \\[1ex]
B_1 \vee B_2 \vee B_3~=~
(\Sigma S(LV)^+ \wedge V^{\infty})\vee
(\Sigma S(LV)^+ \wedge V^{\infty})\\[1ex]
\hphantom{B_1 \vee B_2 \vee B_3~=~(\Sigma S(LV)^+ \wedge V^{\infty})}
\vee ((V^{\infty} \wedge V^{\infty})\vee
(V^{\infty} \wedge V^{\infty}))\end{array}$$
a $\ZZ_2$-equivariant map with stable components
$$d~=~(d_1,d_2,d_3) \in
\{A;B_1\vee B_2\vee B_3\}_{\ZZ_2}~=~ \{A;B_1\}_{\ZZ_2}
\oplus\{A;B_2\}_{\ZZ_2} \oplus\{A;B_3\}_{\ZZ_2}$$ given by
$$\begin{array}{l}
d_1~=~d_2~=~1 ~:~
A~=~\Sigma S(LV)^+ \wedge V^{\infty} \to B_1~=~B_2~=~
\Sigma S(LV)^+ \wedge V^{\infty}~,\\[1ex]
d_3~=~j_V~:~A~=~\Sigma S(LV)^+ \wedge V^{\infty} \to B_3~=~
(V^{\infty} \wedge V^{\infty})\vee (V^{\infty} \wedge V^{\infty})\end{array}$$
Thus up to stable $\ZZ_2$-equivariant homotopy}
$$\begin{array}{l}
h_V(F+G)~=~h_V(F)+h_V(G)+[(F\wedge G)\vee
(G \wedge F)](j_V \wedge \Delta_X)~:\\[1ex]
\hskip50pt
\Sigma S(LV)^+ \wedge V^{\infty} \wedge X \to
LV^{\infty} \wedge V^{\infty} \wedge Y \wedge Y~.\end{array}$$
\end{proposition}
\begin{proof}
(i) By construction.\\
(ii) A homotopy $G:F \simeq F': V^{\infty} \wedge X \to
V^{\infty} \wedge Y$ determines $\ZZ_2$-equivariant homotopies
$$\begin{array}{l}
1 \wedge G~:~1 \wedge F~\simeq~1 \wedge F'~:~
LV^{\infty} \wedge V^{\infty} \wedge X \to
LV^{\infty} \wedge V^{\infty} \wedge Y~,\\[1ex]
G \wedge G~:~F \wedge F~\simeq~F' \wedge F'~:~
(V^{\infty} \wedge X) \wedge (V^{\infty} \wedge X) \to
(V^{\infty} \wedge Y) \wedge (V^{\infty} \wedge Y) \end{array}$$
and hence by \ref{difcon2} (i) a $\ZZ_2$-equivariant homotopy
$$\begin{array}{l}
\delta((\kappa_V^{-1} \wedge 1)(G \wedge G)(\kappa_V \wedge \Delta_X),
(1\wedge \Delta_Y)(1\wedge G))~:\\[1ex]
h_V(F)~=~\delta((\kappa_V^{-1} \wedge 1)(F \wedge F)(\kappa_V \wedge \Delta_X),
(1\wedge \Delta_Y)(1\wedge F))\\[1ex]
\hphantom{h_V(F)~} \simeq~
h_V(F')~=~\delta((\kappa_V^{-1} \wedge 1)(F' \wedge F')(\kappa_V \wedge \Delta_X),
(1\wedge \Delta_Y)(1\wedge F'))\\[1ex]
\hskip50pt:~
\Sigma S(LV)^+ \wedge V^{\infty} \wedge X \to
LV^{\infty}  \wedge V^{\infty} \wedge Y\wedge Y~.\end{array}$$
(iii) The identity $h_{U \oplus V}(1_{U^{\infty}}\wedge F)=
(1_{(U\oplus LU)^{\infty}}\wedge h_V(F))(1 \wedge \Sigma j_{LU,LV})$
follows from the commutative diagram
$$\xymatrix@C+30pt@R+20pt
{(LU \oplus LV)^{\infty}=LU^{\infty} \wedge LV^{\infty}
\ar[r]^-{\di{\alpha_{L(U\oplus V)}}}
\ar[dr]_-{\di{1_{LU^{\infty}}\wedge\alpha_{LV}}} &
\ar[d]^-{\di{\Sigma j_{LU,LV}}} \Sigma S(L(U\oplus V))^{\infty}\\
& LU^{\infty} \wedge \Sigma S(LV)^+~.}$$
(iv) This is a special case of \ref{difcon2} (i). \\
(v) By construction $h_V(F)=\delta(p,q)$, $h_V(G)=\delta(r,s)$ with
$$\begin{array}{l}
p~=~(\kappa^{-1}_V\wedge 1) (F \wedge F)(\kappa_V \wedge \Delta_X)~,~
q~=~(1 \wedge \Delta_Y)(1 \wedge F)~:\\[1ex]
\hskip100pt
LV^{\infty} \wedge V^{\infty} \wedge X \to LV^{\infty} \wedge V^{\infty}\wedge Y \wedge Y~,\\[1ex]
r~=~(\kappa^{-1}_V\wedge 1) (G \wedge G)(\kappa_V \wedge \Delta_Y)~,~s~=~
(1 \wedge \Delta_Z)(1 \wedge G)~:\\[1ex]
\hskip100pt
LV^{\infty} \wedge V^{\infty} \wedge Y \to LV^{\infty} \wedge V^{\infty}\wedge Z \wedge Z~.
\end{array}$$
Now
$$\begin{array}{ll}
(\kappa^{-1}_V \wedge 1)(G \wedge G)(\kappa_V \wedge 1)q&
=~(\kappa^{-1}_V \wedge 1)(G \wedge G)(\kappa_V \wedge 1)\Delta_Y(1\wedge F)\\[1ex]
&=~r(1 \wedge F)~:\\[1ex]
&LV^{\infty} \wedge
V^{\infty} \wedge X \to LV^{\infty} \wedge V^{\infty} \wedge Z \wedge Z~,
\end{array}$$
so Proposition \ref{difcon2} (v) can be applied to give a
($\ZZ_2$-equivariant) homotopy
$$\begin{array}{ll}
h_V(GF)&=~\delta((\kappa^{-1}_V\wedge 1)
(G \wedge G)(\kappa_V\wedge 1)p,s(1 \wedge F))\\[1ex]
&\simeq~\delta(r(1\wedge F),s(1 \wedge F),)\\[1ex]
&\hskip25pt +
\delta(
(\kappa_V^{-1}\wedge 1)(G\wedge G)(\kappa_V\wedge 1)p,
(\kappa_V^{-1}\wedge 1)(G \wedge G)(\kappa_V\wedge 1)q)\\[1ex]
&=~h_V(G)(1\wedge F)+(\kappa_V^{-1}\wedge 1)(G \wedge G)(\kappa_V\wedge 1)h_V(F)~:\\[1ex]
&\hskip40pt \Sigma S(LV)^+ \wedge V^{\infty} \wedge X
\to LV^{\infty} \wedge V^{\infty} \wedge Z \wedge Z ~.
\end{array}$$
(vi) Consider first the special case of a product map of the type
$$F \wedge 1_Z~:~
V^{\infty} \wedge X \wedge Z \to V^{\infty} \wedge Y \wedge Z$$
when the product formula to be obtained is
$$\begin{array}{l}
h_V(F\wedge 1_Z)~=~h_V(F)\wedge \Delta_Z~:\\[1ex]
\hskip25pt
\Sigma S(LV)^+ \wedge V^{\infty} \wedge X \wedge Z
\to LV^{\infty} \wedge V^{\infty}\wedge ((Y \wedge Z) \wedge (Y \wedge Z))~.
\end{array}$$
By construction $h_V(F)=\delta(p_F,q_F)$ with
$$\begin{array}{l}
p_F~=~(\kappa^{-1}_V\wedge 1) (F \wedge F)(\kappa_V \wedge \Delta_X)~,~
q_F~=~(1 \wedge \Delta_Y)(1 \wedge F)~:\\[1ex]
\hskip50pt \Sigma S(LV)^+ \wedge V^{\infty} \wedge X
\to LV^{\infty} \wedge V^{\infty} \wedge Y\wedge Y~.\end{array}$$
It now follows from
$$\Delta_{X \wedge Z}~=~\Delta_X \wedge \Delta_Z~~,~~
\Delta_{Y \wedge Z}~=~\Delta_Y \wedge \Delta_Z$$
that
$$\begin{array}{l}
h_V(F \wedge 1)~=~\delta(p_{F \wedge 1},q_{F \wedge 1})~
=~\delta(p_F \wedge \Delta_Z,q_F \wedge \Delta_Z)\\[1ex]
\hphantom{h_V(F \wedge 1)~}=~\delta(p_F,q_F) \wedge \Delta_Z~=~h_V(F)\wedge \Delta_Z~:\\[1ex]
\Sigma S(LV)^+ \wedge V^{\infty} \wedge X\wedge Z
\to LV^{\infty} \wedge V^{\infty} \wedge ((Y \wedge Z) \wedge (Y\wedge Z))~.\end{array}$$
The general case follows from the suspension formula (iii), the composition
formula (v)  and the special case, because
$$F_1\wedge F_2~=~(1\wedge F_2)\circ (F_1\wedge 1)~:~
(V_1^{\infty}\wedge X_1)\wedge (V_2^{\infty}\wedge X_2)\to
(V_1^{\infty}\wedge Y_1)\wedge (V_2^{\infty}\wedge Y_2)~.$$
(vii) By construction.\\
(viii) The Hurewicz map
$$\begin{array}{l}
[\Sigma S(LV)^+\wedge V^{\infty},(V^{\infty} \vee V^{\infty})\wedge
(V^{\infty} \vee V^{\infty})]_{\ZZ_2} \to \\[1ex]
\to H_0({\rm Hom}_{\ZZ[\ZZ_2]}(\dot C(\Sigma S(LV)^+\wedge V^{\infty}),
\dot C((V^{\infty} \vee V^{\infty})\wedge (V^{\infty} \vee V^{\infty})))\\[1ex]
=~Q_0^V(\ZZ \oplus \ZZ)~=~Q_0^V(\ZZ) \oplus Q_0^V(\ZZ)
\oplus H_0((\ZZ \oplus \ZZ)\otimes (\ZZ \oplus \ZZ))~=~\ZZ \oplus \ZZ \oplus \ZZ^4
\end{array}$$
is an isomorphism and the symmetrization map
$$\begin{array}{l}
1+T~=~2 \oplus 2 \oplus 1~:~Q_0^V(\ZZ \oplus \ZZ)~=~\ZZ \oplus \ZZ \oplus \ZZ^4\\[1ex]
\to
Q^0_V(\ZZ \oplus \ZZ)~=~Q^0_V(\ZZ) \oplus Q^0_V(\ZZ) \oplus H_0((\ZZ \oplus \ZZ)\otimes
(\ZZ \oplus \ZZ))~=~\ZZ \oplus \ZZ \oplus \ZZ^4
\end{array}$$
is an injection with
$$\begin{array}{ll}
(1+T)h_V(\nabla)&=~\nabla^{\%}\dot\phi_V(S^0)-\dot\phi_V(S^0\vee S^0)\nabla\\[1ex]
&=~(0,0,(1,0)\otimes (0,1)+(0,1)\otimes (1,0))\\[1ex]
&=~(1+T)i_Vj_V \in Q^0_V(\ZZ \oplus \ZZ)~.
\end{array}$$
(In Example \ref{pinch} we shall give a more geometric proof of
$h_V(\nabla)=i_Vj_V$ using the fact that $\nabla$ is an
Umkehr map of an immersion $\{*_1,*_2\} \to \{*\}$ of 0-dimensional manifolds
with a single double point.)\\
(ix) By definition, $F+G$ is the composite
$$\begin{array}{l}F+G~=~(F \vee G)(\nabla_V \wedge 1)~:\\[1ex]
V^{\infty} \wedge X \xymatrix@C+10pt{\ar[r]^-{\di{\nabla_V\wedge 1}}&}
(V^{\infty} \vee V^{\infty}) \wedge X ~=~V^{\infty} \wedge(X \vee X)
\xymatrix@C+10pt{\ar[r]^-{\di{F \vee G}}&} V^{\infty} \wedge Y\end{array}$$
with $\nabla_V:V^{\infty} \to V^{\infty} \vee V^{\infty}$ a sum map
as in Definition \ref{sum-map} (iii).
By the composition formula (v), the product formula (vi),
the one-point union formula (vii) and the sum map formula (viii)
$$\begin{array}{l}h_V(F+G)\\[1ex]
=~h_V(F \vee G)(\nabla_V \wedge 1)+
((F \vee G) \wedge (F \vee G))h_V(\nabla_V \wedge 1)\\[1ex]
=~(h_V(F) \vee h_V(G))(\nabla_V\wedge 1)+
((F \vee G) \wedge (F \vee G))(h_V(\nabla_V) \wedge \Delta_X)\\[1ex]
=~h_V(F)+h_V(G)+((F \vee G) \wedge (F \vee G))(h_V(\nabla_V) \wedge \Delta_X)\\[1ex]
=~h_V(F)+h_V(G)+((F\wedge G)\vee (G\vee F))(j_V \wedge \Delta_X)~:\\[1ex]
\hskip100pt
\Sigma S(LV)^+ \wedge V^{\infty}\wedge X \to
LV^{\infty}\wedge V^{\infty}\wedge (Y  \wedge Y)~.\end{array}$$
\hfill\qed\end{proof}

\begin{remark} \label{hopf8}
{\rm As a special case of the composition formula \ref{hopf7} (v),
we have the naturality of the geometric Hopf invariant for $X$ and
$Y$ in the homotopy category, that is
$$h_V((1\wedge g)\circ F\circ (1\wedge f))~=~
(1\wedge (g\wedge g))\circ h_V(F)\circ (1\wedge f)$$ for maps
$F:V^{\infty} \wedge X \to V^{\infty} \wedge Y$, $f: X'\to X$ and
$g:Y\to Y'$. \hfill\qed}
\end{remark}

In the case $X=Y=S^0$ we have:

\begin{proposition}~  \label{X=Y=S0}
The stable $\ZZ_2$-equivariant homotopy class of the
geometric Hopf invariant of a stable map $F:V^{\infty} \to V^{\infty}$
of degree $d \in \ZZ$ is
$$h_V(F)~=~\frac{d(d-1)}{2} \in \{\Sigma S(LV)^+\wedge V^{\infty};
LV^{\infty} \wedge V^{\infty}\}_{\ZZ_2}~=~\{S^0;P(V)^+\}~=~\ZZ~.$$
\end{proposition}
\begin{proof} We offer three proofs:\\
(I) The maps
$p,q:LV^{\infty}\wedge V^{\infty} \to LV^{\infty}\wedge V^{\infty}$
used to define $h_V(F)=\delta(p,q)$ have
$$\hbox{\rm degree}(p)~=~d^2~,~\hbox{\rm degree}(q)~=~d \in \ZZ$$
so that by Proposition \ref{00}
$$\begin{array}{ll}
h_V(F)&=~\hbox{\rm semidegree}(h_V(F))\\[1ex]
&=~\dfrac{\hbox{\rm degree}(p)-\hbox{\rm degree}(q)}{2}\\[1ex]
&=~\dfrac{d^2-d}{2} \in
[\Sigma S(LV)^+\wedge V^{\infty},LV^{\infty} \wedge V^{\infty}]_{\ZZ_2}~=~\ZZ~.
\end{array}$$
(II) As in the proof of Proposition \ref{hopf7} (ix)
note that the Hurewicz map
$$\begin{array}{l}
[\Sigma S(LV)^+\wedge V^{\infty},V^{\infty}\wedge V^{\infty}]_{\ZZ_2} \to \\[1ex]
\hskip25pt
H_0({\rm Hom}_{\ZZ[\ZZ_2]}(
\dot C(\Sigma S(LV)^+\wedge V^{\infty}),
\dot C(V^{\infty}\wedge V^{\infty}) ))~=~Q_0^V(\ZZ)~=~\ZZ
\end{array}$$
is an isomorphism and the symmetrization map
$$1+T~=~2~:~Q^V_0(Z)~=~\ZZ \to Q^0_V(\ZZ)~=~\ZZ$$
is an injection with
$$(1+T)h_V(d)~=~d^{\%}\dot\phi_V(S^0)-\dot\phi_V(S^0)d~
=~d^2-d\in Q^0_V(\ZZ)~=~\ZZ~.$$
(III) Proposition \ref{stablesequence1} (i)
gives a split short exact sequence
$$\begin{array}{l}
\{S^1;LV^{\infty}\}_{\ZZ_2}=0 \to
\{\Sigma S(LV)^+;LV^{\infty}\}_{\ZZ_2}=\{S^0;P(V)^+\}_{\ZZ_2}=\ZZ\\[1ex]
\hskip50pt
\to \{LV^{\infty};LV^{\infty}\}_{\ZZ_2}=\{S^0;S^0\}_{\ZZ_2}=A(\ZZ_2)\\[1ex]
\hskip100pt \xymatrix{\ar[r]^-{\di{\rho}}&} \{S^0;LV^{\infty}\}_{\ZZ_2}=\{S^0;S^0\}=\ZZ \to 0
\end{array}$$
with $A(\ZZ_2)$ the Burnside ring of finite $\ZZ_2$-sets and $\rho$ the
$\ZZ_2$-fixed point map.
Assume that $F$ is transverse regular at $0 \in V^{\infty}$, so that
$\{x \in V\,\vert\, F(x)=0\}$ is a finite set with $d$ points (counted
algebraically), with
$$F~=~d \in [V^{\infty},V^{\infty}]~=~\ZZ~.$$
The image of $h_V(F)$ in $A(\ZZ_2)$
is the class of the finite free $\ZZ_2$-set
$$S~=~\{(x,y) \in V \times V\,\vert\,F(x)=F(y)=0,x \neq y\}$$
with $d^2-d=d(d-1)$ points, so that
$$h_V(F)~=~\frac{\vert S \vert}{2}~=~\frac{d(d-1)}{2} \in \ZZ~.$$
\hfill\qed\end{proof}

\begin{example} \label{twist2} {\rm
(i) Given $c \in O(V)$ define
$$d~=~{\rm degree}(c^{\infty}:V^{\infty} \to V^{\infty})~=~
\begin{cases}
+1&{\rm if}~{\rm det}(c)>0\\
-1&{\rm if}~{\rm det}(c)<0~.
\end{cases}$$
By Proposition \ref{X=Y=S0} the stable $\ZZ_2$-equivariant homotopy class
of the geometric Hopf invariant of $c^{\infty}:V^{\infty} \to V^{\infty}$
is given by
$$\begin{array}{l}
h_V(c^{\infty}) ~=~\dfrac{d(d-1)}{2}~=~\begin{cases}
0&{\rm if}~{\rm det}(c)>0\\
1&{\rm if}~{\rm det}(c)<0
\end{cases}\\[3ex]
\hskip75pt
\in \omega_0(P(\R(\infty)))~=~{\rm ker}(\rho:\omega_{0,0}\to \omega_0)~=~ \ZZ
\end{array}
$$
(using the terminology of Example \ref{cohomotopy}). \\
(ii) In the special case of (i)
$$\begin{array}{l}
V~=~\R~,~c~=~-1 \in GL(V)~=~\R\backslash \{0\}~,\\[1ex]
c^{\infty}~=~-1~:~V^{\infty}~=~S^1 \to V^{\infty}~=~S^1~;~t \mapsto 1-t
\end{array}$$
with degree $d=-1$ the formula in (i) gives that
$$h_{\R}(-1)~=~1 \in \ZZ~.$$
\hfill\qed}
\end{example}

\begin{example} \label{hopf5}
{\rm Let $V=\R$, $X=S^{2n}$, $Y=S^n$, for any $n \geqslant 0$.
The $\ZZ_2$-equivariant homotopy class of the geometric Hopf invariant of a map
$$F~:~V^{\infty} \wedge X~=~S^{2n+1} \to V^{\infty} \wedge Y~=~S^{n+1}$$
is an integer
$$\begin{array}{l}
h_{\R}(F) \in
[\Sigma S(LV)^+ \wedge V^{\infty}\wedge X,LV^{\infty} \wedge V^{\infty}\wedge Y \wedge Y]_{\ZZ_2}\\[1ex]
\hskip35pt =~[V^{\infty}\wedge V^{\infty}\wedge X,LV^{\infty} \wedge V^{\infty}\wedge Y \wedge Y]~=~
\pi_{2n+2}(S^{2n+2})~=~\ZZ
\end{array}$$
(using $\Sigma S(LV)^+=V^{\infty} \vee LV^{\infty}$), namely the degree of the map
$$h_{\R}(F)\vert~:~V^{\infty} \wedge V^{\infty} \wedge X~=~S^{2n+2} \to
LV^{\infty} \wedge V^{\infty} \wedge Y \wedge Y~=~S^{2n+2}~.$$
(i) For $n=0$ and any $d \in \ZZ$ the stable $\ZZ_2$-equivariant homotopy class
of the geometric Hopf invariant of the map
$$d~:~\Sigma X~=~S^1 \to \Sigma Y~=~S^1~;~t \mapsto dt$$
is given by Proposition \ref{X=Y=S0} to be
$$h_{\R}(d)~=~\frac{d(d-1)}{2} \in \ZZ~,$$
See \ref{double7} below for a geometric interpretation in the case $d
\geqslant 1$.\\
(ii) See \ref{link6} below for the identification of
the $\ZZ_2$-equivariant homotopy class of the
geometric Hopf invariant of a map $F:\Sigma X=S^{2n+1} \to \Sigma
Y=S^{n+1}$ in the case $n \geqslant 1$ with the classical Hopf
invariant of $F$.  \hfill\qed}
\end{example}

The geometric Hopf invariant is closely related to the relative version
of the symmetric construction, for maps:

\begin{proposition}~ \label{sym3}
{\rm (i)} Given a map $F:X \to Y$ let
$$Z~=~\Cc(F)~=~Y \cup_F CX$$
be the mapping cone, so that there is defined a homotopy cofibration sequence
$$\xymatrix@C+10pt{X \ar[r]^-{\di{F}} &  Y \ar[r]^-{\di{G}}&
Z\ar[r]^-{\di{H}}& \Sigma X}$$
with $G$ the inclusion, $H$ the projection.
Let $E:CX \to Z$ be the inclusion, a canonical null-homotopy of
$GF:X \to Z$ . For any inner product space $V$ the $V$-coefficient
symmetric construction $\dot\phi_V(Z)$ is the composite
$\ZZ_2$-equivariant map
$$\dot\phi_V(Z)~:~S(LV)^+ \wedge Z
\xymatrix{\ar[r]^-{1\wedge H} &} \Sigma S(LV)^+ \wedge X
\xymatrix{\ar[r]^-{\delta(r,Tr)}&}  Z\wedge Z$$
with
$$r~=~(E \wedge GF)\phi_V(X)~,~Tr~=~(GF \wedge E)\phi_V(X)~:~
CS(LV)^+ \wedge X \to Z\wedge Z$$
two null-homotopies of
$$(GF \wedge GF)\phi_V(X)~:~S(LV)^+ \wedge X \to Z \wedge Z~.$$
{\rm (ii)}  The geometric Hopf invariant of a map $F:V^{\infty} \wedge X \to V^{\infty}
\wedge Y$ is the relative difference (\ref{reldif1}) of the null-homotopies of
$$\begin{array}{l}
(F\wedge F)\dot\phi_V(V^{\infty} \wedge X)~=~
\dot\phi_V(V^{\infty} \wedge Y)(1 \wedge F)~:\\[1ex]
\hskip50pt
S(LV)^+ \wedge V^{\infty} \wedge X \to V^{\infty} \wedge V^{\infty} \wedge Y
\wedge Y
\end{array}
$$
given by \ref{sym4}
$$\begin{array}{l}
\delta \dot\phi_V(V^{\infty}\wedge Y)(1 \wedge F)~,~
(F \wedge F)\delta \dot\phi_V(V^{\infty}\wedge X)~:\\[1ex]
\hskip50pt
C S(LV)^+ \wedge V^{\infty} \wedge X \to
V^{\infty} \wedge V^{\infty} \wedge Y \wedge Y\end{array}$$
that is
$$\begin{array}{l}
h_V(F)~=~(\kappa_V^{-1}\wedge 1)
\delta\big(\delta \dot\phi_V(V^{\infty}\wedge Y)(1 \wedge F),(F \wedge F)\delta \dot\phi_V(V^{\infty}\wedge X)\big)~:\\[1ex]
\hskip75pt
\Sigma S(LV)^+ \wedge V^{\infty} \wedge X \to
LV^{\infty} \wedge V^{\infty} \wedge Y \wedge Y~.
\end{array}$$
{\rm (iii)} Given a map $F:V^{\infty} \wedge X \to V^{\infty} \wedge Y$  let
$$\xymatrix@C+10pt{V^{\infty}\wedge X \ar[r]^-{\di{F}} &  V^{\infty} \wedge Y \ar[r]^-{\di{G}}&
Z\ar[r]^-{\di{H}}& \Sigma V^{\infty} \wedge X}$$
be as in {\rm (i)} with $Z=\Cc(F)$ the mapping cone, $G$ the inclusion and $H$
the projection. The $V$-coefficient
symmetric construction $\dot\phi_V(Z)$ is given up to $\ZZ_2$-equivariant
homotopy by
$$\dot\phi_V(Z)~=~(G \wedge G)h_V(F)H~:~S(LV)^+\wedge Z \to Z \wedge Z~.$$
\end{proposition}
\begin{proof}
(i) By construction.\\
(ii) By Proposition \ref{difcon7} (i).\\
(iii) Combine (ii) and Proposition \ref{reldif2} (ii).\\
\hfill\qed\end{proof}

The geometric Hopf invariant $h_V(F)$ has a particularly simple
description in the special case $V=\R$.

\begin{definition}~ \label{plus} {\rm For $V = \R$
$$\begin{array}{l}
V^{\infty}~=~S^1~,~S(LV)~=~S^0~,~
LV\backslash \{0\}~=~\R_+ \sqcup \R_-~,\\[1ex]
\Sigma S(LV)^+~=~
(LV\backslash \{0\})^{\infty}~=~S^1 \vee S^1
\end{array}$$
with $\ZZ_2$ acting on $S^0$ and $S^1\vee S^1$ by transposition, and
$$\R_+~=~\{x \in \R\,\vert\,x >0\}~,~\R_-~=~\{x \in \R\,\vert\,x <0\}~.$$
The geometric Hopf invariant of a map
$$F~:~\R^{\infty} \wedge X~=~\Sigma X \to \R^{\infty} \wedge Y~=~\Sigma Y$$
is the $\ZZ_2$-equivariant map
$$h_\R(F)~=~h_\R(F)_+\vee h_\R(F)_-~:~
(S^1 \vee S^1)\wedge \Sigma X \to
L\R^{\infty}\wedge \R^{\infty} \wedge Y \wedge Y~=~\Sigma Y \wedge \Sigma Y
$$
which is entirely determined by the non-equivariant
{\it $\R$-coefficient geometric Hopf invariant}
$$h_\R(F)_+~:~(S^1 \vee \{*\})\wedge \Sigma X~=~\Sigma^2X
 \to \Sigma Y \wedge \Sigma Y~.$$
 }\hfill\qed
\end{definition}

\begin{proposition}~ \label{hopfr}
{\rm (i)} The geometric Hopf invariant of the sum
of maps $F,G:\Sigma X \to \Sigma Y$
$$F+G~:~\Sigma X \to \Sigma Y~;~(s,x) \mapsto
\begin{cases}
F(2s,x)&{\it if}~
0 \leqslant s \leqslant 1/2\\
G(2s-1,x)&{\it if}~
1/2 \leqslant s \leqslant 1
\end{cases}$$
is given up to homotopy by
$$h_{\R}(F+G)_+~=~h_{\R}(F)_++h_{\R}(G)_++F \wedge G~:~\Sigma^2X
    \to \Sigma Y \wedge \Sigma Y~.$$
{\rm (ii)} The geometric Hopf invariant of the  map defined for
any pointed spaces $A_1,A_2$ by
$$F~~:~\Sigma (A_1 \times A_2) \to \Sigma (A_1 \vee A_2)~;~
(s,a_1,a_2) \mapsto \begin{cases}
(2s,a_1)&{\it if}~
0 \leqslant s \leqslant 1/2\\
(2s-1,a_2)&{\it if}~
1/2 \leqslant s \leqslant 1
\end{cases}$$
is given up to homotopy by
$$h_{\R}(F)_+~:~\Sigma^2(A_1\times A_2)
\xymatrix{\ar[r]^-{\Sigma^2 p}&} \Sigma^2(A_1 \wedge A_2) \subset
\Sigma (A_1\vee A_2) \wedge \Sigma (A_1 \vee A_2)$$
with $p:A_1 \times A_2 \to A_1 \wedge A_2$ the projection.
\end{proposition}
\begin{proof} (i) This is the special case $V=\R$ of the sum formula of
Proposition \ref{hopf7} (ix).\\
(ii) Apply (i) to $F=\Sigma F_1+\Sigma F_2$, with
$$F_i~:~A_1 \times A_2 \to A_1 \vee A_2~;~(a_1,a_2) \mapsto a_i~.$$
\hfill\qed\end{proof}

\begin{example} {\rm
For any pointed space $X$ define the map
$$F~=~\Sigma p_1+\Sigma p_2~:~\Sigma(X \times X) \to \Sigma X$$
with
$$p_i~:~X\times X \to X~;~(x_1,x_2) \mapsto x_i~~(i=1,2)~.$$
By Proposition \ref{hopfr} the geometric Hopf invariant of $F$ is
given up to homotopy by
$$h_{\R}(F)_+~=~\Sigma^2p~:~\Sigma^2(X \times X) \to \Sigma^2(X \wedge X)$$
with $p:X\times X \to X \wedge X$ the projection.
}\hfill\qed
\end{example}

\section{The stable geometric Hopf invariant $h'_V(F)$}\label{h'}

The stable relative difference $\delta'(p,q)$ defined in \ref{stablereldif3}
above will now be used to construct a second form of the geometric Hopf
invariant of a map $F:V^{\infty} \wedge X \to V^{\infty} \wedge Y$, to
be a $\ZZ_2$-equivariant map
$$h'_V(F)~:~LV^{\infty} \wedge V^{\infty} \wedge X \to
LV^{\infty} \wedge V^{\infty} \wedge S(LV)^+ \wedge Y \wedge Y$$
with $LV$ as defined in section \ref{hopf}. We shall regard $h'_V(F)$
as a stable $\ZZ_2$-equivariant map
$$h'_V(F)~:~ X \sto S(LV)^+ \wedge Y \wedge Y~.$$

The symmetric construction on the mapping
cone of a map $F:V^{\infty} \wedge X \to V^{\infty} \wedge Y$ will be
shown to factorize (up to $\ZZ_2$-equivariant homotopy) through the
stable geometric Hopf invariant $h'_V(F)$ of \S\ref{h'}.  Thus
$h'_V(F)$ is a geometric version of the functional Steenrod squares
(cf. Example \ref{quad2} below).

As before, define the $\ZZ_2$-equivariant maps
$$\begin{array}{l}
p~=~(\kappa_V^{-1}\wedge 1)(F \wedge F)(\kappa_V \wedge \Delta_X)~,~
q~=~(1 \wedge \Delta_Y)(1 \wedge F)~:\\[1ex]
\hskip110pt
LV^{\infty} \wedge V^{\infty} \wedge X \to LV^{\infty} \wedge V^{\infty}
\wedge Y \wedge Y\end{array}$$
which agree on $0^+ \wedge V^{\infty} \wedge X =V^{\infty} \wedge X \subset
LV^{\infty} \wedge V^{\infty} \wedge X$. The geometric Hopf invariant
(\ref{hopf1}) is the $\ZZ_2$-equivariant map defined by
$$h_V(F)~=~\delta(p,q)~:~\Sigma S(LV)^+ \wedge V^{\infty} \wedge X \to
LV^{\infty} \wedge V^{\infty} \wedge Y \wedge Y~.$$
\begin{definition}~ \label{hopf9}
{\rm The {\it stable geometric Hopf invariant} of a
map $F:V^{\infty} \wedge X \to V^{\infty} \wedge Y$ is
the stable relative difference of $p$ and $q$ (\ref{stablereldif3}),
the $\ZZ_2$-equivariant map
\index{geometric Hopf invariant!stable, $h'_V(F)$}
$$h'_V(F)~=~\delta'(p,q)~:~LV^{\infty} \wedge V^{\infty} \wedge X \to
LV^{\infty} \wedge V^{\infty} \wedge S(LV)^+ \wedge Y \wedge Y$$
which we shall regard as a stable $\ZZ_2$-equivariant map
$$h'_V(F)~:~X \sto S(LV)^+ \wedge Y \wedge Y~.$$}
\hfill\qed
\end{definition}

\begin{proposition}~ \label{hopf10}
{\rm (i)} The $\ZZ_2$-equivariant $S$-duality isomorphism
$$\begin{array}{ll}
\{\Sigma S(LV)^+\wedge V^{\infty} \wedge X;
LV^{\infty}\wedge V^{\infty} \wedge Y\wedge Y\}_{\ZZ_2}
&\to \{X; S(LV)^+\wedge Y\wedge Y\}_{\ZZ_2}~;\\[1ex]
&f \mapsto (1\wedge f)(\Delta\alpha_{LV}\wedge 1)
\end{array}$$
sends the geometric Hopf invariant map
$$h_V(F)~:~\Sigma S(LV)^+ \wedge V^{\infty} \wedge X \to
LV^{\infty}\wedge V^{\infty} \wedge Y\wedge Y$$
to the stable geometric Hopf invariant map $h'_V(F)$.\\
{\rm (ii)} The stable geometric Hopf invariant is such that there is a
commutative diagram of $\ZZ_2$-equivariant maps
$$\xymatrix@C+60pt@R+20pt
{LV^{\infty} \wedge V^{\infty} \wedge X \ar[r]^-{\di{h'_V(F)}}
\ar[dr]_-{\di{1\wedge ((\kappa^{-1}_V \wedge 1)(F \wedge F)
(\kappa_V\wedge\Delta_X) -\Delta_Y (1\wedge F))~~~~~~~~~~~~~}} &
\ar[d]^-{\di{s_{LV} \wedge 1}}LV^{\infty}
\wedge V^{\infty} \wedge S(LV)^+ \wedge Y\wedge Y\\
& LV^{\infty} \wedge V^{\infty} \wedge Y \wedge Y}$$
with
$$s_{LV}~:~S(LV)^+ \to S(LV)^+/S(LV)~=~0^+~;~u \mapsto 0$$
the $\ZZ_2$-equivariant pinch map.\\
{\rm (iii)} The stable geometric Hopf invariant defines a function
$$h~:~[V^{\infty}  \wedge X;V^{\infty} \wedge Y] \to
\{X;S(LV)^+ \wedge Y \wedge Y\}_{\ZZ_2}~;~F \mapsto h'_V(F)$$
such that for any $F:V^{\infty}  \wedge X \to V^{\infty} \wedge Y$,
$G:V^{\infty}  \wedge Y \to V^{\infty} \wedge Z$
$$h'_V(GF)~=~h'_V(G)F+(G \wedge G)h'_V(F)~:~X \sto
S(LV)^+ \wedge Z \wedge Z~.$$
\end{proposition}
\begin{proof}
(i) This is a special case of \ref{difcon5} (i). (See also Example \ref{Z2adjoint2}).\\
(ii) This is a special case of \ref{difcon5} (ii).\\
(iii) By construction and the composition formula of \ref{hopf7}
(v). \\
\hfill\qed\end{proof}

Proposition \ref{stablesequence1} gives a commutative braid of
exact sequences of stable $\ZZ_2$-equivariant homotopy groups

$$\xymatrix@C-60pt{
\{X;S(LV)^+ \wedge Y \wedge Y\}_{\ZZ_2}
\ar[dr]^-{\di{s_{LV}}}\ar@/^2pc/[rr] &&
\{S(LV)^+ \wedge X;Y \wedge Y\}_{\ZZ_2}\ar[dr]&&&\\
&\{X;Y \wedge Y\}_{\ZZ_2} \ar[ur]^-{\di{s^*_{LV}}} \ar[dr]^-{\di{0_{LV}}}&&
\{S(LV \oplus LV)^+\wedge X;LV^{\infty} \wedge Y\wedge Y\}_{\ZZ_2}
\ar[dr] &\\
\{LV^{\infty} \wedge X;Y \wedge Y\}_{\ZZ_2}
\ar@/_2pc/[rr]\ar[ur]^-{\di{0^*_{LV}}} &&
\{X;LV^{\infty}\wedge Y \wedge Y\}_{\ZZ_2}\ar[ur]\ar@/_2pc/[rr]&&A
}$$

\bigskip

\noindent and for a pointed map $F:V^{\infty} \wedge X \to V^{\infty} \wedge Y$
with $V \subseteq \R(\infty)$
$$\begin{array}{l}
p~=~(F\wedge F)\Delta_X~=~(h'_{\R(\infty)}(F),F)~,~q~=~\Delta_YF~=~(0,F)\\[1ex]
\hskip25pt \in
\{X;Y \wedge Y\}_{\ZZ_2}~=~\{X;S(\infty)^+\wedge(Y \wedge Y)\}_{\ZZ_2} \oplus \{X;Y\}~,\\[1ex]
s_{LV}^*(p)~=~(F\wedge F)\dot\phi_V(X)~,~s_{LV}^*(q)~=~\dot\phi_V(Y)F
\in \{S(LV)^+ \wedge X;Y \wedge Y\}_{\ZZ_2}~,\\[1ex]
0_{LV}(p)~=~0_{LV}(q)~=~(0,F)\\[1ex]
\hskip25pt \in
\{V^{\infty} \wedge X;(V^{\infty} \wedge Y) \wedge (V^{\infty} \wedge Y)\}_{\ZZ_2}~=~
\{X;LV^{\infty} \wedge Y \wedge Y\}_{\ZZ_2}\\[1ex]
\hskip50pt =~
\{X;S(\infty)/S(LV) \wedge (Y \wedge Y)\}_{\ZZ_2} \oplus \{X;Y\}~,\\[1ex]
p-q~=~s_{LV}h'_V(F)
\in {\rm ker}(0_{LV}:\{X;Y\wedge Y\}_{\ZZ_2} \to \{X;LV^{\infty} \wedge Y\wedge Y\}_{\ZZ_2})\\[1ex]
\hskip50pt
=~{\rm im}(s_{LV}:\{X;S(LV)^+ \wedge Y\wedge Y\}_{\ZZ_2} \to \{X;Y\wedge Y\}_{\ZZ_2})~,\\[1ex]
\{LV^{\infty} \wedge X;Y \wedge Y\}_{\ZZ_2}~=~
\{S(\infty)/S(LV) \wedge X;Y \wedge Y\}_{\ZZ_2} \oplus \{X;Y\}~,\\[1ex]
A~=~\{S(LV)^+\wedge X;LV^{\infty} \wedge Y\wedge Y\}_{\ZZ_2}~=~
\{X;\Sigma S(LV)^+ \wedge Y \wedge Y\}_{\ZZ_2}~.
\end{array}$$

For any pointed space $X$ there is defined a split short exact sequence
$$\xymatrix{0 \ar[r] & \{X;P(\R(\infty))^+\}\ar[r] &
\widetilde{\omega}^0_{\ZZ_2}(X)\ar[r]^-{\di{\rho}} &
\widetilde{\omega}^0(X) \ar[r] & 0}$$
with $\rho$ the $\ZZ_2$-fixed point map and
$$\begin{array}{l}
\sigma~:~\widetilde{\omega}^0(X) \to \widetilde{\omega}^0_{\ZZ_2}(X)~;\\[1ex]
(F:V^{\infty} \wedge X \to V^{\infty}) \mapsto
(1\wedge F:LV^{\infty} \wedge V^{\infty} \wedge X \to LV^{\infty} \wedge
V^{\infty})
\end{array}$$
such that $\rho \sigma=1:\widetilde{\omega}^0(X) \to \widetilde{\omega}^0(X)$
(Example \ref{stablesequence2} (iv)).

\begin{definition} \label{square} (Crabb \cite[pp.32-33]{crabb})\\{\rm
(i) The {\it squaring operation} from stable cohomotopy to stable
$\ZZ_2$-cohomotopy
$$\begin{array}{l}
P^2~:~\omega^0(X) \to \omega^0_{\ZZ_2}(X)~;~
(F:V^{\infty} \wedge X^+ \to V^{\infty}) \mapsto \\[1ex]
((\kappa^{-1}_V\wedge 1)(F \wedge F)(\kappa_V\wedge 1)(1\wedge \Delta_X)
:LV^{\infty} \wedge V^{\infty} \wedge X^+ \to LV^{\infty} \wedge V^{\infty})
\end{array}$$
is a non-additive function such that
$\rho P^2=1:\omega^0(X) \to \omega^0(X)$.\\
(ii) The {\it reduced squaring operation} is defined by
$$\begin{array}{l}
\overline{P}^2~:~\omega^0(X) \to
{\rm ker}(\rho:\omega^0_{\ZZ_2}(X\times X)\to
\omega^0(X))~=~\{X;P(\R(\infty))^+\}~;\\[1ex]
\hspace*{100pt}
(F:V^{\infty} \wedge X^+ \to V^{\infty}) \mapsto P^2(F)-\sigma(F)=h'_V(F)~.
\end{array}$$
\hfill\qed}
\end{definition}

\begin{example} {\rm (i)
$\overline{P}^2$ agrees with the Hopf invariant function $\theta^2$ of
Segal \cite{segal2} (cf. Crabb \cite[p.35]{crabb}).\\
(ii) For $j \equiv 1(\bmod 2)$
the composite of $\overline{P}^2$ for $X=S^j$ and the Hurewicz map
$$\xymatrix{\omega_j=\widetilde{\omega}^0(S^j)
\ar[r]^-{\di{\overline{P}^2}} & \{S^j;P(\R(\infty))^+\}=\omega_j(P(\R(\infty)))
\ar[r]& H_j(P(\R(\infty)))=\ZZ_2}$$
is the classical mod\,2 Hopf invariant (\cite[4.16]{crabb} and Example \ref{steenrod} (iv)).\\
(iii) For $j \equiv 3(\bmod 4)$
the composite of $\overline{P}^2$ for $X=S^j$ and the real $KO_*$-theory
Hurewicz map
$$\begin{array}{l}
\xymatrix{\omega_j=\widetilde{\omega}^0(S^j)
\ar[r]^-{\di{\overline{P}^2}} & \{S^j;P(\R(\infty))^+\}=\omega_j(P(\R(\infty)))}\\[1ex]
\hskip150pt \xymatrix{\ar[r]& KO_j(P(\R(\infty)))=\ZZ[1/2]/\ZZ}
\end{array}$$
is the 2-primary $e$-invariant (\cite[4.17]{crabb}).\\
}\hfill\qed
\end{example}

\section{The quadratic construction $\psi_V(F)$}\label{quad}

The $V$-coefficient `quadratic construction' on a pointed space $X$ is the
pointed space $S(LV)^+\wedge_{\ZZ_2} (X\wedge X)$.

\begin{definition}~ \label{quad1}
{\rm Let $V$ be an inner product space, and let
$F:V^{\infty} \wedge X \to V^{\infty} \wedge Y$ be a pointed map.\\
(i) The {\it space level $V$-coefficient quadratic construction} on $F$
is the stable map
\index{quadratic construction!space level, $\psi_V(F)$}
$$\psi_V(F)~:~X\sto S(LV)^+\wedge_{\ZZ_2} (Y\wedge Y)$$
given by the image of the $\ZZ_2$-equivariant stable geometric Hopf invariant (\ref{hopf9})
$$h'_V(F)~=~\delta'(p,q)~:~ LV^{\infty} \wedge V^{\infty} \wedge X \to
LV^{\infty} \wedge V^{\infty} \wedge S(LV)^+ \wedge Y \wedge Y$$
under the isomorphism given by Proposition \ref{Adams}
$$\{X;S(LV)^+ \wedge_{\ZZ_2} (Y \wedge Y)\}
\iso \{X;S(LV)^+ \wedge Y \wedge Y\}_{\ZZ_2}~.$$
(ii) The {\it chain level $V$-coefficient quadratic construction} on $F$
is the $\ZZ$-module chain map induced by the space level $\psi_V(F)$
$$\psi_V(F)~:~\dot C(X)\to C^{cell}(S(LV))\otimes_{\ZZ[\ZZ_2]}(\dot C(Y)\otimes_{\ZZ} \dot C(Y))~,$$
inducing morphisms
$$\psi_V(F)~:~\dot H_n(X) \to Q^V_n(\dot C(Y))~.$$
More directly, note that the chain level quadratic construction $\psi_V(F)$
is the adjoint of the $\ZZ[\ZZ_2]$-module chain map
$$h_V(F)~:~SC^{cell}(S(LV))\otimes_{\ZZ}\dot C(X)
\to C^{cell}(LV^{\infty})\otimes_{\ZZ} (\dot C(Y)\otimes_{\ZZ} \dot C(Y))$$
induced by the geometric Hopf invariant $\ZZ_2$-equivariant map
$$h_V(F):\Sigma S(LV)^+ \wedge V^{\infty} \wedge X \to
LV^{\infty}\wedge V^{\infty} \wedge Y\wedge Y~,$$
with respect to the chain level $\ZZ_2$-equivariant $S$-duality isomorphism
of \ref{Z2adjoint1} (ii) and  \ref{les2} (ii)
$$\begin{array}{l}
C^{cell}(S(LV))\otimes_{\ZZ[\ZZ_2]}(\dot C(Y)\otimes_{\ZZ} \dot C(Y))\\[1ex]
\hskip50pt
\cong~{\rm Hom}_{\ZZ[\ZZ_2]}(SC^{cell}(S(LV)),\dot C(LV^{\infty})\otimes_{\ZZ}
(\dot C(Y)\otimes_{\ZZ} \dot C(Y)))~.
\end{array}$$
\hfill\qed}
\end{definition}

\begin{remark} {\rm
We refer to Crabb and James \cite{crabbjames} for an extended treatment of fibrewise homotopy theory - there is a brief account in Appendix
\ref{appendix1} below.
Here is an explicit formula for the space level quadratic construction
$\psi_V(F)$ in terms of fibrewise $\ZZ_2$-equivariant homotopy theory.
Think of the Hopf invariant map
$$h_V(F)~:~\Sigma S(LV)^+ \wedge V^{\infty} \wedge X
\to LV^{\infty} \wedge V^{\infty}\wedge Y \wedge Y$$
as a fibrewise $\ZZ_2$-equivariant map
$$S(LV)\times (\Sigma V^{\infty}\wedge X) \to
S(LV)\times (LV^{\infty} \wedge V^{\infty}\wedge Y\wedge Y)$$
over $S(LV)$. By dividing out the free $\ZZ_2$-action, we get
a fibrewise map
$$P(V)\times (\Sigma V^{\infty}\wedge X) \to
S(LV)\times_{\ZZ_2} (LV^{\infty} \wedge V^{\infty}\wedge Y\wedge Y)\eqno{(*)}$$
over the real projective space $P(V)$.\\
\indent
We now use fibrewise duality theory as in Crabb and James \cite[Part
II, (14.43)]{crabbjames}.  This requires an embedding of $P(V)$ in some Euclidean
space.  To be definite, we use the embedding
$$P(V) \emb M~=~{\rm End}(V)~;~[x] \mapsto (v \mapsto
\langle x, v \rangle x )~~(x \in S(LV))$$
mapping a line to the orthogonal projection onto that subspace.
Use the inner product to construct a tubular neighbourhood $\nu
\emb M$, where $\nu$ is the (total space of) the normal
bundle. Notice that the direct sum $\R\oplus\tau$ of a trivial line
and the tangent bundle $\tau$ of $P(V)$ is just $S(LV)\times_{\ZZ_2}
LV$.  And $\tau\oplus \nu$ is the trivial bundle $P(V)\times M$.\\
\indent
So the fibrewise smash product of (*) with the identity on the
fibrewise one-point compactification $\nu^{\infty}_{P(V)}$ of $\nu$
gives a fibrewise map
$$\begin{array}{l}
\nu^{\infty}_{P(V)} \wedge_{P(V)} (P(V) \times (\Sigma V^{\infty}\wedge X))
\to \\[1ex]
\hskip100pt (P(V)\times M^+)\wedge_{P(V)} (S(LV)\times_{\ZZ_2}
(\Sigma V^{\infty}\wedge Y\wedge Y)
\end{array}
$$
over $P(V)$. Collapsing the basepoint section $P(V)$, we get a map
$$
P(V)^\nu \wedge (\Sigma V^{\infty}\wedge X) \to (M^+\wedge \Sigma V^{\infty})
\wedge (S(LV)^+ \wedge_{\ZZ_2}Y \wedge Y)$$
Compose this with the Pontrjagin-Thom map $M^{\infty} \to P(V)^\nu$, and
we have produced an explicit map
$$\psi_V(F)~:~(M^+\wedge \Sigma V^{\infty})\wedge X \to
(M^+ \wedge \Sigma V^{\infty}) \wedge(S(LV)^+ \wedge_{\ZZ_2}Y \wedge Y)~.$$
\hfill\qed}
\end{remark}

\begin{proposition}\label{quadquad}
Let $F:V^{\infty} \wedge X \to V^{\infty} \wedge Y$
be a stable map, with $V=\R^k$. The quadratic construction $\psi_V(F):X \sto
S(LV)^+\wedge_{\ZZ_2} (Y\wedge Y)$ induces the chain level
quadratic construction of Ranicki \cite[\S1]{ranicki2}, \cite[p.29]{ranicki3}
$$\begin{array}{l}
\psi_V(F)~:~\dot C(X) \to \dot C^{cell}(S(LV)^+\wedge_{\ZZ_2} (Y\wedge Y))\\[1ex]
\hphantom{\psi_V(F)~:~\dot C(X) \to} =~
W[0,k-1]\otimes_{\ZZ[\ZZ_2]}(\dot C(Y)\otimes_{\ZZ}\dot C(Y))
\end{array}$$
and hence also the quadratic construction on the level of homology groups
$$\psi_V(F)~:~\dot H_*(X) \to \dot H_*(S(LV)^+\wedge_{\ZZ_2} (Y\wedge Y))~=~
Q^{[0,k-1]}_*(\dot C(Y))~,$$
identifying $C^{cell}(S(LV))=W[0,k-1]$.
\end{proposition}
\begin{proof} From Proposition \ref{sym44} we have a diagram of
$\ZZ[\ZZ_2]$-module chain complexes, chain maps and chain homotopies
$$\xymatrix@C+25pt{
C^{cell}(S(LV))\otimes_{\ZZ}\dot C(V^{\infty})\otimes_{\ZZ} \dot C(X)
\ar@{~>}[dr]^-{\delta\dot\phi_V(V^{\infty},X)}
\ar[d]_-{1\otimes \dot E(V^{\infty},X)}^-{\simeq}
\ar[r]^-{\dot\phi_V(V^{\infty}) \cup \dot\phi_V(X)} &
\dot C(V^{\infty})\otimes_{\ZZ}\dot C(X)\otimes_{\ZZ}\dot C(V^{\infty})\otimes_{\ZZ}\dot C(X)
\ar[d]_-{\dot E(V^{\infty},X) \otimes \dot E(V^{\infty},X)}^-{\simeq}
\\
C^{cell}(S(LV))\otimes_{\ZZ}\dot C(V^{\infty} \wedge X)
\ar[r]^-{\dot \phi_V(V^{\infty} \wedge X)}
\ar[d]_-{1\otimes F}&
\dot C(V^{\infty} \wedge X)\otimes_{\ZZ}\dot C(V^{\infty} \wedge X)
\ar[d]^-{F\otimes F}&\\
C^{cell}(S(LV))\otimes_{\ZZ}\dot C(V^{\infty} \wedge Y)
\ar[r]^-{\dot \phi_V(V^{\infty} \wedge Y)}&
\dot C(V^{\infty} \wedge Y)\otimes_{\ZZ}\dot C(V^{\infty} \wedge Y)\\
C^{cell}(S(LV))\otimes_{\ZZ}\dot C(V^{\infty})\otimes_{\ZZ} \dot C(Y)
\ar@{~>}[ur]^-{\delta\dot\phi_V(V^{\infty},Y)}
\ar[u]^-{1\otimes \dot E(V^{\infty},Y)}_-{\simeq}
\ar[r]^-{\dot\phi_V(V^{\infty}) \cup \dot\phi_V(Y)} &
\dot C(V^{\infty})\otimes_{\ZZ}\dot C(Y)\otimes_{\ZZ}\dot C(V^{\infty})\otimes_{\ZZ}\dot C(Y)
\ar[u]^-{\dot E(V^{\infty},Y) \otimes \dot E(V^{\infty},Y)}_-{\simeq} }$$
Evaluation on a generator of $\dot C(V^{\infty})=S^k\ZZ$ now gives
a diagram of $\ZZ[\ZZ_2]$-module chain complexes, chain maps and chain homotopies
$$\xymatrix@C+100pt{
W[0,k-1]\otimes_{\ZZ}S^k\dot C(X)
\ar@{~>}[dr]^-{\delta\dot\phi_{\R^k}(S^k,X)}
\ar[d]_-{1\otimes \dot E(S^k,X)}^-{\simeq}
\ar[r]^-{S^k \dot\phi_{\R^k}(X)} &
S^k\dot C(X)\otimes_{\ZZ}S^k\dot C(X)
\ar[d]_-{\dot E(S^k,X) \otimes \dot E(S^k,X)}^-{\simeq} \\
W[0,k-1]\otimes_{\ZZ}\dot C(\Sigma^k X)
\ar[r]^-{\dot \phi_{\R^k}(\Sigma^k X)}
\ar[d]_-{1\otimes F}&
\dot C(\Sigma^k X)\otimes_{\ZZ}\dot C(\Sigma^k X)
\ar[d]^-{F\otimes F}&\\
W[0,k-1]\otimes_{\ZZ}\dot C(\Sigma^k Y)
\ar[r]^-{\dot \phi_{\R^k}(\Sigma^k Y)}&
\dot C(\Sigma^k Y)\otimes_{\ZZ}\dot C(\Sigma^k Y)\\
W[0,k-1]\otimes_{\ZZ}S^k\dot C(Y)
\ar@{~>}[ur]^-{\delta\dot\phi_{\R^k}(S^k,Y)}
\ar[u]^-{1\otimes \dot E(S^k,Y)}_-{\simeq}
\ar[r]^-{S^k\dot\phi_{\R^k}(Y)} &
S^k\dot C(Y)\otimes_{\ZZ}S^k\dot C(Y)
\ar[u]^-{\dot E(S^k,Y) \otimes \dot E(S^k,Y)}_-{\simeq} }$$
The stable geometric Hopf invariant $h'_{\R^k}(F)$ of
Proposition \ref{hopf9} (ii) induces the same chain map
as the quadratic construction $\psi_{\R^k}(F)$, which is thus the
chain map given by Proposition \ref{sym55} (ii)
$$\psi_{\R^k}(F)~=~\delta(F,\dot\phi_{\R^k}(X),\dot\phi_{\R^k}(Y))~
:~\dot C(X) \to W[0,k-1]\otimes_{\ZZ}(\dot C(Y)\otimes_{\ZZ} \dot C(X))~.$$
But this is exactly the chain level quadratic construction of \cite{ranicki2,ranicki3}.
\\
\hfill\qed\end{proof}

The quadratic construction $\psi_V(F)$ has the following properties~:

\begin{proposition}~ \label{quad3} {\it
{\rm (i)} The stable $\ZZ_2$-equivariant homotopy class of $\psi_V(F):X \sto
S(LV)^+\wedge (Y \wedge Y)$ depends only on the homotopy
class of $F:V^{\infty} \wedge X \to V^{\infty}\wedge Y$.  The function
$$\psi_V~ :~[V^{\infty} \wedge X,V^{\infty} \wedge Y]\to
\{ X;\, S(LV)^+\wedge (Y\wedge Y)\}_{\ZZ_2}~;~F \mapsto \psi_V(F)$$
is such that
$$\psi_V(F+G)~ =~ \psi_V(F)+\psi_V(G) + i((F\wedge G)\circ \Delta_X)~,$$
where
$$i~:~\{ X;\, Y\wedge Y\}~=~\{ X;\, (S^0)^+ \wedge (Y\wedge Y)\}_{\ZZ_2}
\to    \{ X;\, S(LV)^+\wedge(Y\wedge Y)\}_{\ZZ_2}$$
is induced by the map
$$Y\wedge Y~=~(S^0)^+\wedge (Y\wedge Y) \to
    S(LV)^+\wedge_{\ZZ_2}(Y\wedge Y)~;~ (\pm 1,y_1,y_2) \mapsto (\pm v,y_1,y_2)$$
for any $v \in S(LV)$.\\
{\rm (ii)} The symmetrization of the quadratic construction is the difference
of the symmetric constructions, in the sense that
$\psi_V(F) \in \{X;S(LV)^+\wedge Y \wedge Y\}_{\ZZ_2}$
has images
$$\begin{array}{l}
s_{LV}(\psi_V(F))~=~(F \wedge F)\Delta_X -\Delta_YF \in
\{X;Y\wedge Y\}_{\ZZ_2}~,\\[1ex]
s^*_{LV}s_{LV}(\psi_V(F))~=~(F \wedge F)\dot\phi_V(X) -\dot\phi_V(Y)F \in
\{S(LV)^+ \wedge X;Y\wedge Y\}_{\ZZ_2}~.
\end{array}$$
{\rm (iii)} The quadratic construction on the suspension
$$\Sigma F~:~(V \oplus \R)^{\infty} \wedge X~=~\Sigma(V^{\infty} \wedge X) \to
(V \oplus \R)^{\infty} \wedge Y~=~\Sigma(V^{\infty} \wedge Y)$$
is the composite
$$\psi_{V \oplus \R}(\Sigma F)~:~
X\xymatrix{\ar[r]^{\di{\psi_V(F)}}&} S(LV)^+\wedge_{\ZZ_2}(Y\wedge Y)
\xymatrix{\ar[r]^{\di{e \wedge 1}}&} S(L(V \oplus \R))^{\infty}\wedge_{\ZZ_2}Y \wedge Y$$
with $e:S(LV) \to S(L(V \oplus \R))$ the inclusion induced by
$$e~:~V \to V \oplus \R~;~x \mapsto (x,0)~.$$
{\rm (iv)} If $V=\{0\}$ then $\psi_V(F)=0$, i.e. the quadratic construction is 0
for an unstable map $F:X \to Y$.\\
{\rm (v)} {\rm (}Sum formula{\rm )}
The sum of maps
$F_1,F_2:V^{\infty} \wedge X \to V^{\infty} \wedge Y$
(with ${\rm dim}(V)>0$) is a map
$$F_1+F_2~:~V^{\infty} \wedge X \to V^{\infty} \wedge Y$$
with quadratic construction
$$\begin{array}{ll}
\psi_V(F_1+F_2)&=~\psi_V(F_1) + \psi_V(F_2) +
(F_1 \wedge F_2+F_2 \wedge F_1) \Delta_X\\[1ex]
&:~X \sto S(LV)^+ \wedge_{\ZZ_2}(Y\wedge Y)\end{array}$$
with
$$\begin{array}{l}
F_1 \wedge F_2+F_2 \wedge F_1~:~X \wedge X
\xymatrix@C+15pt{\ar[r]^-{\begin{pmatrix} F_1 \wedge F_2 \\ F_2 \wedge F_1 \end{pmatrix}}&}\\[2ex]
(Y \wedge Y)\vee (Y\wedge Y)~=~(S^0)^+\wedge_{\ZZ_2}(Y\wedge Y)
\xymatrix{\ar@{^{(}->}[r]&} S(LV)^+\wedge_{\ZZ_2}(Y\wedge Y)
\end{array}$$
where $S^0 \subset S(LV)\,;\, \pm \mapsto \pm v$ for any $v \in S(LV)$.\\
{\rm (vi)} {\rm (}Product formula{\rm )}
Let $F_i :V_i^{\infty}\wedge X_i\to V_i^{\infty}\wedge Y_i$, $i=1,\, 2$,
be two maps.  The quadratic construction $\psi_{V_1\oplus V_2} (F_1\wedge F_2)$
of the smash product
$$F_1\wedge F_2~:~(V_1\oplus V_2)^{\infty} \wedge (X_1\wedge X_2)
\to (V_1\oplus V_2)^{\infty} \wedge (Y_1\wedge Y_2)$$
is the sum
$$\begin{array}{l}
\psi_{V_1\oplus V_2} (F_1\wedge F_2)\\[1ex]
=~
(a_1\wedge 1)(\psi_{V_1}(F_1)\wedge \Delta_{Y_2}F_2)+
(a_2\wedge 1)(\Delta_{Y_1}F_1\wedge \psi_{V_2}(F_2))\\[1ex]
\hskip150pt + (a_3 \wedge 1)(\psi_{V_1}(F_1)\wedge \psi_{V_2}(F_2))~:\\[1ex]
\hskip50pt
X_1 \wedge X_2 \sto S(LV_1\oplus LV_2)^+ \wedge (Y_1 \wedge Y_1)
\wedge (Y_2 \wedge Y_2)
\end{array}$$
with
$$\begin{array}{l}
a_1~:~S(LV_1) \to S(LV_1\oplus LV_2)~;~x_1 \mapsto (x_1,0)~,\\[1ex]
a_2~:~S(LV_2) \to S(LV_1\oplus LV_2)~;~x_2 \mapsto (0,x_2)~,\\[1ex]
a_3~:~S(LV_1) \times S(LV_2) \to S(LV_1\oplus LV_2)~;~
(x_1,x_2) \mapsto (x_1/\sqrt{2},x_2/\sqrt{2})~.
\end{array}$$
{\rm (vii)} {\rm (}Composition formula{\rm )} The composite of maps
$$F~:~V^{\infty} \wedge X \to V^{\infty} \wedge Y~~,~~
G~:~V^{\infty} \wedge Y \to V^{\infty} \wedge Z$$
is a map $GF:V^{\infty} \wedge X \to V^{\infty} \wedge Z$
with quadratic construction}
$$\psi_V(GF)~=~(G \wedge G)\psi_V(F) + \psi_V(G)F~:~
X \sto S(LV)^+ \wedge_{\ZZ_2} (Z \wedge Z)~.$$
\end{proposition}
\begin{proof} Immediate from the corresponding properties of the
geometric Hopf invariant $h_V(F)$ (\ref{hopf7}).\\
\hfill\qed\end{proof}

For any pointed space $X$ write
$$\Sigma^{\infty}X~=~\R(\infty)^{\infty}\wedge X$$
with $\R(\infty)=\varinjlim\limits_k\R^k$.
The quadratic construction $\psi_V(F)$ also has a stable version,
for $V=\R(\infty)$ :

\begin{definition} {\rm
The {\it quadratic construction} on a stable map
$F:\Sigma^{\infty}X \to \Sigma^{\infty}Y$ is the stable homotopy class
$$\psi(F)~=~\psi_{\R(\infty)}(F)
\in \{ X;\, S(\infty)^+\wedge (Y\wedge Y)\}_{\ZZ_2}~=~
\{ X;\, S(\infty)^+\wedge_{\ZZ_2} (Y\wedge Y)\}$$
with $S(\infty)$ a contractible space with a free $\ZZ_2$-action.  \\
\hfill\qed}
\end{definition}

\begin{remark} \label{quad2} {\rm
(i) The stable homotopy group $\pi^S_{2i}(K(\ZZ_2,i))=\ZZ_2$ is generated by
the Hopf construction $S^{2i+1} \to \Sigma K(\ZZ_2,i)$ on the map
$$S^i \times S^i \to K(\ZZ_2,i) \times K(\ZZ_2,i) \to K(\ZZ_2,i)~,$$
with the $\bmod\,2$ coefficient quadratic construction  defining an isomorphism
$$\begin{array}{l}
\pi^S_{2i}(K(\ZZ_2,i)) \xymatrix{\ar[r]^-{\di{\cong}}&}
\pi_{2i}^S(S(\infty)^+ \wedge_{\ZZ_2}(K(\ZZ_2,i)\wedge K(\ZZ_2,i)))~=~\ZZ_2~;\\[1ex]
\hskip50pt (F:V^{\infty} \wedge S^{2i} \to V^{\infty} \wedge K(\ZZ_2,i))\mapsto \psi(F)~.
\end{array}$$
(ii) The $\bmod\,2$ coefficient quadratic construction of a stable map
$F:V^{\infty} \wedge X \to V^{\infty} \wedge Y$
can be expressed in terms of the functional Steenrod squares:
for any $y \in H^i(Y;\ZZ_2)=[Y,K(\ZZ_2,i)]$
$$\begin{array}{l}
y_{\%}\psi_V(F)~:~H_j(X;\ZZ_2) \xymatrix{\ar[r]^-{\di{\psi_V(F)}}&}Q_j(C(Y;\ZZ_2))
\xymatrix{\ar[r]^-{\di{y_{\%}}}&} Q_j(S^i\ZZ_2)~=~\ZZ_2~;\\[1ex]
\hskip100pt
x \mapsto \langle Sq^{j-i+1}_{yF}(\iota),x\rangle
\end{array}$$
with $\iota \in H^i(K(\ZZ_2,i)=\ZZ_2$ the generator
(cf. \cite[Proposition 1.6]{ranicki2}).  If $j=2i$ and
$x \in {\rm im}(\pi^S_{2i}(X) \to H_{2i}(X;\ZZ_2))$  then
$$y_{\%}\psi_V(x)~=~yFx \in \pi_{2i}^S(K(\ZZ_2,i))~=~\ZZ_2~,$$
with (i) the special case $Y=K(\ZZ_2,i)$.\\
(iii) The $\pi$-equivariant version of the quadratic construction
was used in \cite{ranicki2} to express the non-simply-connected
surgery obstruction of Wall \cite{wall2} as the cobordism class
of a quadratic Poincar\'e complex: see \S\ref{surgeryobstruction} below for a more
detailed discussion.}\hfill\qed
\end{remark}

\section{The ultraquadratic construction $\widehat{\psi}(F)$}\label{ultra}

The ultraquadratic construction is the quadratic construction
$\psi_V(F)$ on a map $F:V^{\infty}\wedge X \to V^{\infty}\wedge Y$ in
the special case $V=\R$. A knot determines such a map $F$, and the
ultraquadratic construction on $F$ is a homotopy theoretic version of
the Seifert form on the homology of a Seifert surface.  See Ranicki
\cite[pp. 814--842]{ranicki3} for an earlier account of the ultraquadratic
theory.

\begin{definition}~ \label{ultraquadratic}
{\rm The {\it ultraquadratic construction} on a map $F:\Sigma X \to \Sigma Y$
is the quadratic construction for the special case $V=\R$
\index{ultraquadratic construction, $\widehat{\psi}(F)$}
$$\widehat{\psi}(F)~=~\psi_{\R}(F)~:~X \sto
S(LV)^+\wedge_{\ZZ_2} (Y \wedge Y)~=~Y \wedge Y~,$$
identifying $S(LV)=S^0=\{1,-1\}$ with $\ZZ_2$ acting by permutation.\hfill\qed}
\end{definition}

\begin{proposition}~ The ultraquadratic construction defines a function
$$\begin{array}{l}
\widehat{\psi}~:~[\Sigma X,\Sigma Y]~=~[\R^{\infty}\wedge X,\R^{\infty} \wedge Y] \to\\[1ex]
\hskip50pt
\{X;S(L\R)^+ \wedge Y \wedge Y\}_{\ZZ_2}~=~\{X;Y \wedge Y\}~;~F \mapsto \psi_\R(F)~,
\end{array}$$
with a commutative braid of exact sequences of abelian groups
$$\xymatrix@C-40pt{
\{X;Y \wedge Y\}
\ar[dr]^-{\di{s_{L\R}}}\ar@/^2pc/[rr]^-{\di{1+T}} &&
\{X;Y \wedge Y\}\ar[dr]&\\
&\{X;Y \wedge Y\}_{\ZZ_2} \ar[ur]^-{\di{s^*_{L\R}}} \ar[dr]^-{\di{0_{L\R}}}&&
A\\
\{L\R^{\infty} \wedge X;Y \wedge Y\}_{\ZZ_2}
\ar@/_2pc/[rr]\ar[ur]^-{\di{0^*_{L\R}}} &&
\{X;L\R^{\infty}\wedge Y \wedge Y\}_{\ZZ_2}\ar[ur]&}$$
\bigskip

\noindent where
$A=\{S(L\R\oplus L\R)^+ \wedge X;L\R^{\infty} \wedge Y \wedge Y\}_{\ZZ_2}$,
$s^*_{L\R}:\{X;Y \wedge Y\}_{\ZZ_2} \to \{X;Y \wedge Y\}$ is the
forgetful map, and
$$s_{L\R}\widehat{\psi}(F)~=~(F \wedge F)\Delta_X-\Delta_YF
\in \{X;Y \wedge Y\}_{\ZZ_2}~.$$
\end{proposition}
\begin{proof} This is the special case $V=\R$ of Proposition \ref{quad3}
(i)+(ii).
\hfill\qed\end{proof}

\begin{example} {\rm Let $N^n \subset S^{n+1}$ be a codimension 1 framed
submanifold with boundary $\partial N$, with Pontrjagin-Thom Umkehr map
$$F~:~S^{n+1}~=~\Sigma S^n \to
S^{n+1}/{\rm cl.}(S^{n+1}-N \times [0,1])~=~\Sigma (N/\partial N)$$
(interpreting $N/\partial N$ as $N^\infty$ if $\partial N=\emptyset$).
The Hurewicz image
$$[\widehat{\psi}(F)] \in
H_n(\dot C(N/\partial N)\otimes_\ZZ \dot C(N/\partial N))$$
of the ultraquadratic construction
$\widehat{\psi}(F) \in \{S^n;(N/\partial N) \wedge (N/\partial N)\}$
is a $\ZZ$-module chain map (or rather a chain homotopy class)
$$[\widehat{\psi}(F)]~:~\dot C(N/\partial N)^{n-*} \to \dot C(N/\partial N)~,$$
defining an
`$n$-dimensional ultraquadratic complex over $\ZZ$' in the sense of
Ranicki \cite[p.814]{ranicki3}. If $\partial N=S^{n-1}$
then $N$ is a Seifert surface for the $(n-1)$-knot
$\partial N=S^{n-1} \subset S^{n+1}$
and $(\dot C(N/\partial N),[\widehat{\psi}(F)])$ is a chain complex
generalization of the Seifert form on $\dot H_*(N/\partial N)=H_*(N,\partial N)$ -- see
Example \ref{embed12} below for the connection with linking numbers and
the Hopf invariant.\\
\hfill \qed}
\end{example}

In dealing with the special case $V=\R$ we shall use the following
terminology~:

\begin{definition}~ \label{difcon8} {\rm
(i) The {\it positive and negative lines}
\index{negative!line, $V^-$}
\index{positive!line, $V^{\infty}$}
$$V^+~=~\{v \in V\st v\geqslant 0\}^{\infty}~,~
V^-~=~\{v \in V\st v\leqslant  0\}^{\infty} \subset V^{\infty}$$
are homeomorphic to $[ 0,1]$, and such that
$$V^{\infty}~=~V^+\cup_{0,\infty}V^-~.$$
(ii) The {\it positive and negative circles}
\index{negative!circle, $S^-(V)$}
\index{positive!circle, $S^+(V)$}
$$S^+(V)~=~\alpha_V(V^+)~,~S^-(V)~=~\alpha_V(V^-) \subset \Sigma S(V)^+$$
(with $\alpha_V:V^+ \to V^+/0^+=\Sigma S(V)^+$ the canonical projection)
are homeomorphic to $S^1$, and such that
$$\Sigma S(V)^+~=~S^+(V)\vee S^-(V)~.$$
(iii) The {\it positive and negative differences} of maps $p,q:V^{\infty} \wedge X \to Y$ such that
$p\vert=q\vert:0^+ \wedge X \to Y$ are the restrictions of the
relative difference $\delta(p,q):\Sigma S(V)^+ \wedge X \to Y$
(\ref{difcon1} (ii)) to the positive and negative circles
\index{negative!difference, $\delta^-(p,q)$}
\index{positive!difference, $\delta^+(p,q)$}
$$\begin{array}{l}
\delta^+(p,q)~=~\delta(p,q)\vert~:~S^+(V) \wedge X \to Y~,\\[1ex]
\delta^-(p,q)~=~\delta(p,q)\vert~:~S^-(V) \wedge X \to Y
\end{array}$$
with
$$\delta(p,q)~=~\delta^+(p,q)\vee\delta^-(p,q)~:~
(S^+(V)\wedge X) \vee (S^-(V)\wedge X) \to Y~.$$
\hfill\qed
}\end{definition}

More directly, working with suspension coordinates
the positive and negative differences of maps
$p,q:\Sigma X \to Y$ such that
$$p(1/2,x)~=~q(1/2,x)\in Y~~(x \in X)$$
are given (up to homotopy and rescaling) by
$$\begin{array}{l}
\delta^+(p,q)~:~\Sigma X \to Y~;~(t,x) \mapsto \begin{cases}q(1-t,x)&\hbox{if
$0 \leqslant  t \leqslant  1/2$}\\[1ex]
p(t,x)&\hbox{if $1/2 \leqslant t
\leqslant  1$~,}
\end{cases}\\[4ex]
\delta^-(p,q)~:~\Sigma X \to Y~;~(t,x)
\mapsto \begin{cases}
q(t,x)&\hbox{if $0 \leqslant  t \leqslant  1/2$}\\[1ex]
p(1-t,x)&\hbox{if $1/2 \leqslant  t \leqslant  1$~.}
\end{cases}
\end{array}$$

\begin{lemma} \label{difcon9}
For any map $q:V^{\infty} \wedge X \to Y$
such that $q(0^+ \wedge X)=*$ the positive difference
$\delta^+(*,q)$ is homotopic to the map
$$[q]^+~:~S^+(V) \wedge X \to Y~;~(v,x) \mapsto q(v,x)~.$$
Similarly for the negative difference.
\end{lemma}
\begin{proof}  The homotopies are the
restrictions of the homotopy of \ref{difcon2} (iv)
$$\begin{array}{l}
\delta(*,q)~ =~\delta^+(*,q)\vee \delta^-(*,q)~
\simeq~ [q]~=~[q]^+ \vee [q]^-~:\\[1ex]
\hskip100pt \Sigma S(V)^+ \wedge X~=~
(S^+(V) \wedge X) \vee (S^-(V) \wedge X) \to Y~.
\end{array}$$
\hfill\qed\end{proof}

\begin{proposition}~ \label{hopf3} {\it
Suppose that $V=\R$, so that (as in \ref{difcon8})
$$\Sigma S(LV)^+~=~S^+(LV)\vee S^-(LV)~.$$
The $\ZZ_2$-action $T:\Sigma S(LV)^+ \to \Sigma S(LV)^+$
interchanges the positive and negative circles
$$S^+(LV)~=~T(S^-(LV))~~,~~S^-(LV)~=~T(S^+(LV))~.$$
The Hopf invariant (\ref{hopf1}) of a map
$F:V^{\infty} \wedge X \to V^{\infty} \wedge Y$ is a $\ZZ_2$-equivariant map
$$\begin{array}{l}
h_V(F)~=~\delta(p,q)~=~\delta^+(p,q) \vee T\delta^+(p,q)~:\\[1ex]
\Sigma S(LV)^+ \wedge V^{\infty} \wedge X~=~
(S^+(LV) \wedge V^{\infty} \wedge X) \vee
T(S^+(LV) \wedge V^{\infty} \wedge X)\\[1ex]
\hskip175pt \to (V^{\infty} \wedge Y) \wedge (V^{\infty} \wedge Y)
\end{array}$$
which is determined by the positive difference
$$h_V(F)\vert~=~\delta^+(p,q)~:~
S^+(LV) \wedge V^{\infty} \wedge X \to
(V^{\infty} \wedge Y) \wedge (V^{\infty} \wedge Y)~.$$
The positive difference is homotopic to the Hopf map of Boardman and
Steer \cite{bs}}
$$\delta^+(p,q)~\simeq~\lambda_2(F)~:~
V^{\infty}\wedge V^{\infty} \wedge X \to
(V^{\infty} \wedge Y)\wedge (V^{\infty} \wedge Y)~.$$
\end{proposition}
\begin{proof} In the first instance,
we recall the maps $\mu_2$, $\lambda_2$ of \cite{bs}.
Express $V^{\infty} \wedge V^{\infty}$ as a union
$$V^{\infty} \wedge V^{\infty} ~=~H^+(V) \cup_{\Delta(V)^{\infty}} H^-(V)$$
of two contractible subspaces
$$\begin{array}{l}
H^+(V)~=~\{(v,w)\in V^{\infty} \wedge V^{\infty}\st v \geqslant w\}~,\\[1ex]
H^-(V)~=~\{(v,w)\in V^{\infty} \wedge V^{\infty}\st v \leqslant w\}~,
\end{array}$$ with
$$H^+(V) \cap H^-(V)~=~\Delta(V)^{\infty}~=~
\{(v,w)\in V^{\infty} \wedge V^{\infty}\st v=w\}~.$$
(Up to homeomorphism this is just $S^2=D^2\cup_{S^1}D^2$).
Given pointed spaces $X,Y_1,Y_2$ and a map
$$f~:~V^{\infty}\wedge X \to Y_1 \vee Y_2$$
let $\pi_i:Y_1 \vee Y_2 \to Y_i$ ($i=1,2$) be the projections, and define
$$f_i~=~\pi_if~:~V^{\infty} \wedge X \to Y_i~.$$
Note that for each $(v,x) \in V^{\infty} \wedge X$ either $f_1(v,x)=*$ or $f_2(v,x)=*$.
The $\mu_2$-function of \cite[5.1]{bs}
$$\mu_2~:~[V^{\infty}\wedge X,Y_1 \vee Y_2] \to [V^{\infty} \wedge V^{\infty} \wedge X,Y_1 \wedge Y_2]$$
is defined by
$$\begin{array}{l}
\mu_2(f)~=~(f_1 \wedge f_2)(1\wedge \Delta_X) \cup *~:\\[1ex]
V^{\infty} \wedge V^{\infty}\wedge X ~=~(H^+(V)\wedge X)
\cup_{\Delta(V)^{\infty}\wedge X} (H^-(V)\wedge X) \to Y_1 \wedge Y_2~.
\end{array}$$
The projection
$$\pi^+~:~V^{\infty} \wedge V^{\infty} \to H^+(V)/{\Delta(V)^{\infty}}~;~
(v,w) \mapsto \begin{cases}
(v,w)&\hbox{if $v \geqslant w$}\\[1ex]
*&\hbox{if $v \leqslant w$}
\end{cases}$$ is a homotopy equivalence such that there is defined a
commutative diagram
$$\xymatrix@R+20pt@C+20pt
{V^{\infty} \wedge V^{\infty}\wedge X \ar[dr]_-{\di{\mu_2(f)}}
\ar[r]^-{\di{\pi^+ \wedge 1}}_-{\di{\simeq}} &
(H^+(V)/\Delta(V)^{\infty})\wedge X
\ar[d]^-{\di{[(f_1 \wedge f_2)(1\wedge \Delta_X)]}}\\
&Y_1 \wedge Y_2~.}$$
Let $\nabla:V^{\infty} \to V^{\infty} \vee V^{\infty}$ be a sum map,
with the components
$$\nabla_i~=~\pi_i\nabla~:~V^{\infty} \xymatrix{\ar[r]^-{\di{\nabla}}&}
V^{\infty} \vee V^{\infty} \xymatrix{\ar[r]^-{\di{\pi_i}}&} V^{\infty}$$
both homotopic to the identity, and chosen such that
$$(\nabla_1 \wedge \nabla_2)\vert~=~*~:~
H^+(V)/\Delta(V)^{\infty} \to V^{\infty} \wedge V^{\infty}~.$$
(If $\dot\phi:(0,1)\to \R$ is an order-preserving homeomorphism, then
the sum map
$$\nabla~:~V^{\infty} \to V^{\infty} \vee V^{\infty}~;~\dot\phi(t) \mapsto
\begin{cases}\dot\phi(2t)_1&
\hbox{if $t \leqslant  1/2$}\\[1ex]
\dot\phi(2t-1)_2&\hbox{if $t \geqslant 1/2$}
\end{cases}
$$ has these properties, with $\dot\phi(0)=\dot\phi(1)=\infty$). The
Hopf map of \cite{bs} is given by the composite
$$\begin{array}{ll}
\lambda_2~=~&\mu_2(\nabla\wedge 1_Y) ~:\\[1ex]
&[V^{\infty}\wedge X,V^{\infty}\wedge Y]
\xymatrix{\ar[r]^-{\di{\nabla\wedge 1_Y}}&}
[V^{\infty}\wedge X,(V^{\infty}\wedge Y) \vee (V^{\infty}\wedge Y)] \\[1ex]
&\hskip25pt \xymatrix{\ar[r]^-{\di{\mu_2}}&}
[V^{\infty}\wedge V^{\infty} \wedge X,
(V^{\infty}\wedge Y) \wedge (V^{\infty}\wedge Y)]~.
\end{array}$$
The Hopf invariant of a map $F:V^{\infty}\wedge X \to V^{\infty}\wedge Y$
is thus given by
$$\begin{array}{ll}
\lambda_2(F)~~&=~\mu_2((\nabla \wedge 1_Y)F)\\[1ex]
&=~[((\nabla_1 \wedge 1_Y)F\wedge (\nabla_2 \wedge 1_Y)F))
(1\wedge \Delta_X)]\cup *~:\\[1ex]
&\hskip20pt
V^{\infty}\wedge V^{\infty} \wedge X ~=~(H^+(V)\wedge X)
\cup_{\Delta(V)^{\infty}\wedge X} (H^-(V)\wedge X)\\[1ex]
&\hskip120pt \to (V^{\infty}\wedge Y) \wedge (V^{\infty}\wedge Y)
\end{array}$$
with a commutative diagram
$$\xymatrix@R+20pt@C+20pt
{V^{\infty} \wedge V^{\infty}\wedge X \ar[dr]_-{\di{\lambda_2(F)}}
\ar[r]^-{\di{\pi^+ \wedge 1}}_-{\di{\simeq}} &
(H^+(V)/\Delta(V)^{\infty})\wedge X
\ar[d]^-{\di{[(\nabla_1\wedge \nabla_2 \wedge 1)
((F \wedge F)(1\wedge \Delta_X))]}}\\
&(V^{\infty} \wedge Y) \wedge (V^{\infty} \wedge Y)~.}$$
\indent Now return to the difference construction (\ref{hopf1}) of the Hopf invariant
$h_V(F)=\delta(p,q)$, with
$$\begin{array}{l}
p~=~(\kappa_V \wedge \Delta_Y)(1 \wedge F)~~,~~
q~=~(F \wedge F)(\kappa_V \wedge \Delta_X)~:\\[1ex]
\hskip75pt
LV^{\infty} \wedge V^{\infty} \wedge X \to
(V^{\infty} \wedge Y) \wedge (V^{\infty} \wedge Y)~.\end{array}$$
The composites
$$\begin{array}{l}
LV^{\infty} \wedge V^{\infty} \wedge X  \xymatrix{\ar[r]&}
LV^{\infty} \wedge V^{\infty} \wedge X \xymatrix{\ar[r]^-{\di{p}}&}
(V^{\infty} \wedge Y) \wedge (V^{\infty} \wedge Y)\\[1ex]
\hskip125pt \xymatrix@C+40pt
{\ar[r]^-{\di{\nabla_1 \wedge \nabla_2 \wedge 1}}&}
(V^{\infty} \wedge Y) \wedge (V^{\infty} \wedge Y)~,\\[1ex]
0^+ \wedge V^{\infty} \wedge X \xymatrix{\ar[r]&}
LV^{\infty} \wedge V^{\infty} \wedge X \xymatrix{\ar[r]^-{\di{q}}&}
(V^{\infty} \wedge Y) \wedge (V^{\infty} \wedge Y)\\[1ex]
\hskip125pt \xymatrix@C+40pt
{\ar[r]^-{\di{\nabla_1\wedge \nabla_2 \wedge 1}}&}
(V^{\infty} \wedge Y) \wedge (V^{\infty} \wedge Y)\end{array}$$
are both $*$. There exist homotopies
$$\nabla_1~\simeq~\nabla_2~\simeq~1~:~V^{\infty}\to V^{\infty}$$
and hence also a homotopy
$$\nabla_1\wedge \nabla_2~\simeq~1~:~
V^{\infty} \wedge V^{\infty} \to V^{\infty} \wedge V^{\infty}~.$$
Thus there exists a homotopy
$$\delta^+ (p,q)~\simeq~
(\nabla_1 \wedge \nabla_2\wedge 1) \delta^+(p,q)~=~
\delta^+(*,(\nabla_1\wedge \nabla_2 \wedge 1) q)~.$$
Furthermore, \ref{difcon9} gives a homotopy
$$\delta^+(*,(\nabla_1 \wedge \nabla_2 \wedge 1)q)~\simeq~
[(\nabla_1 \wedge \nabla_2 \wedge 1)q]^+~,$$
so that there exists a homotopy
$$\delta^+(p,q)~\simeq~[(\nabla_1 \wedge \nabla_2 \wedge 1)q]^+~:~
S^+(LV) \wedge V^{\infty} \wedge X \to
(V^{\infty} \wedge Y) \wedge (V^{\infty} \wedge Y)~.$$
The homeomorphism
$$\kappa^+_V~:~S^+(LV)\wedge V^{\infty} \to  H^+(V)/\Delta(V)^{\infty}~;~
(v,w) \mapsto (v+w,-v+w)$$
is such that there is defined a commutative diagram
$$\xymatrix@R+10pt@C+10pt
{S^+(LV) \wedge V^{\infty} \wedge X
\ar[dd]_-{\di{[(\nabla_1 \wedge \nabla_2 \wedge 1)q]^+}}
\ar[dr]^-{\di{\kappa_V^{\infty} \wedge 1}}_-{\di{\cong}} &\\
&(H^+(V)/\Delta(V)^{\infty})\wedge X
\ar[dl]^-{\di{~~~~~~~~~~~~~~~~~~~~
[(\nabla_1\wedge \nabla_2 \wedge 1)((F \wedge F)(1\wedge \Delta_X))]}}\\
(V^{\infty} \wedge Y) \wedge (V^{\infty} \wedge Y)~.&}$$
Putting this together, there are obtained a homotopy equivalence
$$\iota~=~(\kappa_V^{\infty})^{-1}\pi^+~:~V^{\infty} \wedge V^{\infty} \to
S^+(LV) \wedge V^{\infty}$$
and a homotopy
$$\lambda_2(F) ~\simeq~\delta^+(p,q)(\iota \wedge 1)~:~
V^{\infty} \wedge V^{\infty} \wedge X \to
(V^{\infty} \wedge Y) \wedge (V^{\infty} \wedge Y)~.$$
\hfill\qed\end{proof}

\begin{terminology}~ \label{hopf4}
{\rm In view of \ref{hopf3} we shall regard the Hopf
invariant of a map $F:\Sigma X \to \Sigma Y$ as a non-equivariant map
$$h_{\R}(F)~:~\Sigma^2X \to \Sigma Y \wedge \Sigma Y~.\eqno{\hbox{\qed}}$$}
\end{terminology}

\begin{definition} \label{james} {\rm
(i) The {\it James map} (\cite{imj2}) is
$$\begin{array}{l}
J~:~X \times X \to Q_\R(X)~=~\Omega \Sigma X~;\\[1ex]
(x_1,x_2) \mapsto \bigg(
s \mapsto \begin{cases} (2s,x_1)&{\rm if}~0 \leqslant s \leqslant 1/2\\
(2s-1,x_2)&{\rm if}~1/2 \leqslant s \leqslant 1\end{cases}\bigg)~.
\end{array}$$
(ii) Let $E:Q_\R(X) \sto X$ be the stable map defined by evaluation
$$E~:~\Sigma Q_\R(X) \to \Sigma X~;~(s,\omega) \mapsto \omega(s)~,$$
so that
$$EJ~=~\pi_1+\pi_2~:~X \times X \sto X$$
with $\pi_i:X \times X \to X;(x_1,x_2) \mapsto x_i$.
Note that the composite of the map
$$i~:~X \to X \times X~;~x \mapsto (x,*)$$
and $J$ is just the inclusion
$$Ji~:~X \to Q_\R(X)~;~x \mapsto (v \mapsto (v,x))~.$$
}
\hfill\qed
\end{definition}

\begin{proposition}~
{\rm (i)} The ultraquadratic construction
$$\widehat{\psi}(E)~:~Q_\R(X)\sto X \wedge X$$
is such that up to stable homotopy
$$\begin{array}{l}
\widehat{\psi}(E)J~=~\widehat{\psi}(EJ)~=~{\rm proj.}~:~X \times X  \sto X \wedge X~,\\[1ex]
\widehat{\psi}(E)Ji~=~\{*\}~:~X  \sto X \wedge X~.
\end{array}$$
{\rm (ii)} The adjoint of a stable map $F:\Sigma W \to \Sigma X$
$${\rm adj}(F)~:~W \to Q_\R(X)~;~x \mapsto (s \mapsto F(s,x))$$
is such that $F=E({\rm adj}(F))$ with
$$\widehat{\psi}(F)~=~\widehat{\psi}(E){\rm adj}(F)~:~W \sto (X \wedge X)~.$$
$$\xymatrix@C+30pt@R+20pt
{W \ar[r]^-{ \di{{\rm adj}(F)} }
\ar[dr]_-{\di{\widehat{\psi}(F)~~}} &
\ar[d]^-{\di{\widehat{\psi}(E)}} Q_\R(X)\\
& X \wedge X~.}$$
\end{proposition}
\begin{proof} By construction.
(This is just the special case $V=\R$ of \ref{quad4} (below).)\\
\hfill\qed\end{proof}

\section{The spectral quadratic construction $s\psi_V(F)$}\label{spectral}

The chain level `spectral quadratic construction' (Ranicki \cite[\S7.3]{ranicki3})
of  a `semistable' map $F:X \to V^{\infty} \wedge Y$ inducing the chain map
$f:\dot C(X)_{*+{\rm dim}(V)} \to \dot C(Y)$ is a natural transformation
$$\psi_F~:~\dot H_{*+{\rm dim}(V)}(X) \to Q_*({\mathcal C}(f))~.$$
We shall now show that this is induced by the following space level geometric Hopf invariant map~:

\begin{definition}~ \label{sym5}
{\rm  The {\it spectral Hopf invariant map} of
a map $F:X \to V^{\infty} \wedge Y$ is the $\ZZ_2$-equivariant map
given by the relative difference (\ref{reldif1})
\index{spectral Hopf invariant, $sh_V(F)$}
$$\begin{array}{l}sh_V(F)~=~
\delta\big((G \wedge G)\delta \dot\phi_V(V^{\infty}\wedge Y)(1 \wedge F),
\dot\phi_V(\Cc(F))(1\wedge GF)\big):\\[1ex]
\hskip150pt
\Sigma S(LV)^+\wedge X~\to \Cc(F)  \wedge \Cc(F)\end{array}$$
with $G:V^{\infty} \wedge Y \to \Cc(F)$ the inclusion in the mapping cone and
$$1 \wedge GF~:~CS(LV)^+ \wedge X ~=~S(LV)^+ \wedge CX\to S(LV)^+ \wedge \Cc(F)$$
the null-homotopy of $1\wedge GF:S(LV)^+ \wedge X \to S(LV)^+ \wedge \Cc(F)$
determined by the inclusion
$GF:CX \to\Cc(F)=(V^{\infty} \wedge Y)\cup_FCX$.\hfill\qed
}\end{definition}

The spectral Hopf invariant has the same properties as the Hopf invariant
obtained in \ref{hopf7}. Here are two particularly important special cases:

\begin{proposition}~  \label{sym7} {\it
{\rm (i)} The $\ZZ_2$-equivariant homotopy class of
$sh_V(F)$ depends only on the homotopy class of $F$.\\
{\rm (ii)} The spectral Hopf invariant of the composite $EF:X \to
V^{\infty}\wedge Z$ of maps $F:X \to V^{\infty} \wedge Y$, $E:V^{\infty} \wedge Y \to
V^{\infty}\wedge Z$ is given up to $\ZZ_2$-equivariant homotopy by
$$\begin{array}{l}
sh_V(EF)~=~(k \wedge \ell)sh_V(F)+(k \wedge \ell)h_V(E)(1 \wedge F)~:\\[1ex]
\hskip100pt \Sigma S(LV)^+ \wedge X \to \Cc(EF) \wedge \Cc(EF)\end{array}$$
with
$$k~=~E \cup 1~:~ \Cc(F)~=~V^{\infty} \wedge Y \cup_FCX \to  \Cc(EF)~=~V^{\infty} \wedge Z \cup_{GF}CX$$
and $\ell:V^{\infty} \wedge Z \to \Cc(GF)$ the inclusion.}
\end{proposition}
\begin{proof} (i) As for the homotopy
invariance \ref{hopf7} (ii)
of the Hopf invariant $h_V$. \\
(ii) As for the composition formula \ref{hopf7} (v) for the Hopf
invariant $h_V$.\\
\hfill\qed\end{proof}

The spectral Hopf invariant is closely related to the symmetric construction
on the mapping cone:

\begin{proposition}~ \label{sym8} {\it
Given a map $F:X \to V^{\infty} \wedge Y$ let
$$Z~=~\Cc(F)~=~V^{\infty} \wedge Y \cup_F CX$$
be the mapping cone, so that there is defined a homotopy cofibration sequence
$$\xymatrix@C+10pt{X \ar[r]^-{\di{F}} & V^{\infty} \wedge Y \ar[r]^-{\di{G}} &
Z\ar[r]^-{\di{H}} & \Sigma X}$$
with $G$ the inclusion and $H$ the projection.\\
{\rm (i)} The symmetric construction $\dot\phi_V(Z)$ is determined by the
spectral Hopf invariant $sh_V(F)$, with a $\ZZ_2$-equivariant homotopy
commutative diagram
$$\xymatrix@C+20pt@R+20pt{
S(LV)^+\wedge Z \ar[r]^-{\di{1 \wedge H}}
\ar[dr]_-{\di{\dot\phi_V(Z)}} & S(LV)^+\wedge \Sigma X=
\Sigma S(LV)^+\wedge X
\ar[d]^-{\di{sh_V(F)}}\\
& Z \wedge Z~.}$$
{\rm (ii)} If $X=V^{\infty} \wedge X_0$ the spectral Hopf invariant $sh_V(F)$
is determined by the Hopf invariant $h_V(F)$, with a commutative diagram
$$\xymatrix@C+20pt@R+20pt{\Sigma S(LV)^+\wedge X \ar[r]^-{\di{h_V(F)}}
\ar[dr]_-{\di{sh_V(F)}} & (V^{\infty} \wedge Y) \wedge (V^{\infty} \wedge Y)
\ar[d]^-{\di{G\wedge G}}\\
& Z \wedge Z~.}$$
{\rm (iii)} Suppose given a space $W$ and a homotopy equivalence
$Z\simeq V^{\infty} \wedge W$, so that the homotopy cofibration sequence can be written as
$$\xymatrix@C+10pt{X \ar[r]^-{\di{F}} & V^{\infty} \wedge Y \ar[r]^-{\di{G}}
&V^{\infty} \wedge  W \ar[r]^-{\di{H}}& \Sigma X~.}$$
The spectral Hopf invariant $sh_V(F)$ is determined by the Hopf invariant
$h_V(G)$, with a stable $\ZZ_2$-equivariant homotopy commutative diagram}
$$\xymatrix@C+20pt@R+20pt{\Sigma S(LV)^+\wedge X \ar[r]^-{\di{1 \wedge F}}
\ar[dr]_-{\di{sh_V(F)}} & \Sigma S(LV)^+\wedge V^{\infty} \wedge Y
\ar[d]^-{\di{-h_V(G)}}\\
& Z\wedge Z~.}$$
\end{proposition}
\begin{proof} (i) As in \ref{sym3} there is defined a
natural transformation of homotopy cofibration sequences
$$\xymatrix@C-10pt@R+10pt{
S(LV)^+ \wedge X \ar[r]^-{1 \wedge F}
\ar[d]_-{\di{\dot\phi_V(X)}}&
S(LV)^+ \wedge V^{\infty} \wedge Y \ar[r]^-{1 \wedge G}
\ar[d]_-{\di{\dot\phi_V(V^{\infty} \wedge Y)}}&
S(LV)^+ \wedge Z \ar[r]^-{1 \wedge H}
\ar[d]_-{\di{\dot\phi_V(F)}}&
S(LV)^+ \wedge \Sigma X \ar[d]_-{\di{\Sigma\dot\phi_V(X)}} \\
X \wedge X \ar[r]^-{F \wedge F} &
(V^{\infty} \wedge Y) \wedge (V^{\infty} \wedge Y) \ar[r] & \Cc(F \wedge F) \ar[r]& \Sigma(X\wedge X)}$$
and \ref{sym4} gives a $\ZZ_2$-equivariant null-homotopy
$\delta\dot\phi_V(V^{\infty} \wedge Y):\dot\phi_V(V^{\infty} \wedge Y) \simeq *$.
Now apply \ref{reldif2} (iii) to get a $\ZZ_2$-equivariant homotopy
$\dot\phi_V(F)\simeq \delta(i,j)$ (for the appropriate $i,j$) and
compose with the projection $ \Cc(F \wedge F) \to Z \wedge Z$ to
get a $\ZZ_2$-equivariant homotopy $\dot\phi_V(Z) \simeq sh_V(F)(1 \wedge H)$.\\
(ii) By construction.\\
(iii) The application of \ref{sym7} to $GF \simeq *$ shows that up to
stable $\ZZ_2$-equivariant homotopy
$$\begin{array}{ll}
sh_V(GF)&=~sh_V(F)+h_V(G)(1\wedge F)\\[1ex]
&=~sh_V(*)~=~0~:~\Sigma S(LV)^+\wedge X \to Z \wedge Z~.
\end{array}$$
\hfill\qed\end{proof}

\begin{remark} \label{sym9}
{\rm The expression in \ref{sym8} (i) of the symmetric
construction $\dot\phi_V(\Cc(F))$ on the mapping cone $\Cc(F)$ of a
map $F:X\to V^{\infty} \wedge Y$ in terms of the spectral Hopf
invariant $sh_V(F)$ is a generalization of the relationship between the
functional Steenrod squares of a stable map $F$ and the Steenrod
squares of $\Cc(F)$.\hfill\qed}
\end{remark}

\begin{definition}~
{\rm  The {\it spectral quadratic construction} on a  map
$F:X \to V^{\infty} \wedge Y$ is the stable $\ZZ_2$-equivariant map
\index{spectral quadratic construction, $s\psi_V(F)$}
$$s\psi_V(F)~:~LV^{\infty} \wedge X\sto S(LV)^+\wedge_{\ZZ_2} (\Cc(F)\wedge \Cc(F))$$
given by the image of the $\ZZ_2$-equivariant spectral Hopf invariant (\ref{hopf9})
$$sh_V(F)~:~\Sigma S(LV)^+  \wedge X \to \Cc(F) \wedge \Cc(F)$$
under the composite of the
$\ZZ_2$-equivariant $S$-duality isomorphism of Proposition \ref{Z2S-dual}
$$\{\Sigma S(LV)^+  \wedge X;\Cc(F) \wedge \Cc(F)\}_{\ZZ_2}~\cong~
\{LV^{\infty} \wedge X;S(LV)^+\wedge \Cc(F) \wedge \Cc(F)\}_{\ZZ_2}~.$$
\hfill\qed}
\end{definition}

\begin{example} {\rm The spectral quadratic construction on a map
$F:V^{\infty} \wedge X \to V^{\infty} \wedge Y$ is the composite
$$\begin{array}{l}
\xymatrix@C+35pt{s\psi_V(F)~:~V^{\infty} \wedge LV^{\infty} \wedge X
\ar[r]^-{\di{1 \wedge \psi_V(F)}} & S(LV)^+ \wedge LV^{\infty}
\wedge V^{\infty} \wedge Y \wedge Y}\\[1ex]
\hskip175pt \cong S(LV)^+ \wedge((V^{\infty} \wedge Y) \wedge (V^{\infty} \wedge Y))\\
\hskip125pt
\xymatrix@C+15pt
{\ar[r]^-{\di{1 \wedge G \wedge G}} & S(LV)^+ \wedge (\Cc(F) \wedge \Cc(F))}
\end{array}$$
with $G:V^{\infty} \wedge Y \to \Cc(F)$ the inclusion.\\
\hfill\qed}
\end{example}

\section{Stably trivialized vector bundles}

We now apply the geometric Hopf invariant to the classification of
pairs $(\delta\xi,\xi)$ consisting of a $U$-vector bundle $\xi:X \to
BO(U)$ and a $V$-stable trivialization
$$\delta\xi~:~\xi \oplus \epsilon_V~\cong~\epsilon_{U\oplus V}$$
with $U,V$ finite-dimensional inner product spaces.
We relate the classifying map
$$c~=~(\delta\xi,\xi)~:~X \to O(V,U \oplus V)$$
to the geometric Hopf invariant $h_V(F)$ of the $V$-stable map
$$\begin{array}{l}
F~:~V^{\infty} \wedge T(\xi)~=~T(\xi \oplus \epsilon_V)
\xymatrix{\ar[r]_-{\di{\cong}}^-{\di{T(\delta\xi)}}&}
T(\epsilon_{U \oplus V})~=~(U \oplus V)^{\infty} \wedge X^+\\[1ex]
\hskip175pt \xymatrix{\ar[r]^-{proj.}&} V^{\infty} \wedge U^{\infty}~,
\end{array}$$
with stable homotopy class
$$\begin{array}{ll}
h_V(F)~=~h'_V(F)~=~\psi_V(F)
\in& \{\Sigma S(LV)^+ \wedge T(\xi);
LV^{\infty} \wedge (U^{\infty}\wedge U^{\infty})\}_{\ZZ_2}\\[1ex]
&=~\{T(\xi);S(LV)^+\wedge (U^{\infty}\wedge U^{\infty})\}_{\ZZ_2}\\[1ex]
&=~\{T(\xi);S(LV)^+\wedge_{\ZZ_2} (U^{\infty}\wedge U^{\infty})\}\\[1ex]
&=~\{X^+;S(LV)^+\wedge_{\ZZ_2} LU^{\infty}\}~.
\end{array}$$
The stable homotopy relationship between the Stiefel space
$O(V,U\oplus V)$ and the stunted projective space
$$S(LV)^+\wedge_{\ZZ_2}LU^{\infty}~=~P(U\oplus V)/P(U)$$
studied by James \cite{imj5} and Crabb \cite{crabb} is  used
in Proposition \ref{twist} (v) to identify the quadratic construction $\psi_V(F)$
with the `local obstruction' $\theta(c)$ of \cite{crabb} (cf. Definition \ref{local} below)
$$\psi_V(F)~=~\theta(c) \in\{X^+;S(LV)^{\infty}\wedge_{\ZZ_2}LU^{\infty}\}~.$$
In particular, $c^{univ}=1:O(V,U \oplus V) \to O(V,U\oplus V)$ classifies
the universal $U$-bundle with a $V$-stable trivialization
$$(\eta(U)~:~O(V,U\oplus V) \to BO(U)~,~
\delta\eta(U)~:~\eta(U)\oplus \epsilon_V~\cong~\epsilon_{U \oplus V})$$
with a corresponding $V$-stable map
$$F^{univ}~:~V^{\infty} \wedge T(\eta(U))~\cong~
(V \oplus U)^{\infty}\wedge O(V,U\oplus V)^+ \to V^{\infty} \wedge U^{\infty}~.$$
The quadratic construction/local obstruction
$$\begin{array}{l}
\psi_V(F^{univ})~=~\theta(c^{univ})\\[1ex]
\in \{T(\eta(U));S(LV)^+\wedge_{\ZZ_2}(U^{\infty}\wedge U^{\infty})\}~=~
\{O(V,U\oplus V)^+;S(LV)^{\infty}\wedge_{\ZZ_2}LU^{\infty}\}
\end{array}$$
defines a stable map $\theta:O(V,U\oplus V)^+\sto S(LV)^{\infty}\wedge_{\ZZ_2}LU^{\infty}$
which is $2\,{\rm dim}(U)$-connected. It follows that the function
$$[X,O(V,U \oplus V)] \to \{X^+;S(LV)^+\wedge_{\ZZ_2}LU^{\infty}\}~;~
c \mapsto \psi_V(F)=\theta(c)$$
is a bijection if $X$ is a $CW$ complex with  ${\rm dim}(X)<2\,{\rm dim}(U)$.

\begin{definition}~ {\rm
(i) Given an inner product space $V$ and a 1-dimensional
subspace $L=\R x \subset V$ ($x \in S(V)$) let
$$L^{\perp}~=~\{v \in V\,\vert\, \langle v,x \rangle  = 0\}~,$$
and define the linear isometry reflecting $V$ in $L^{\perp} \subset V$
$$\begin{array}{l}
R_L~=~-1_L \oplus 1_{L^{\perp}}~:~V~=~L \oplus L^{\perp}
\to V~=~L \oplus L^{\perp}~;\\[1ex]
\hspace*{150pt} v=(\lambda,\mu) \mapsto v - 2\langle v,x \rangle x=(-\lambda,\mu)~.
\end{array}$$
(ii) The {\it reflection} map
$$\begin{array}{l}
R~:~P(U \oplus V)/P(U) \to O(V,U \oplus V)~=~O(U\oplus V)/O(U)~;\\[1ex]
\hspace*{25pt}
L=\R x \mapsto (R_L\vert:v \mapsto (0,v)-2 \langle (0,v),x \rangle x)
~(x \in S(U \oplus V),v \in V)
\end{array}$$
sends $L \subset U \oplus V$ to $R_L\vert:V \to U \oplus V$.\\
\hfill\qed}
\end{definition}

\begin{example} {\rm
(i) For $U=\{0\}$, $V=\R$ the reflection map is a homeomorphism
$$R~:~P(\R)^+~=~S^0 \to O(\R)~;~\pm 1 \mapsto \pm 1~.$$
(ii) For $U=\{0\}$, $V=\R^2$ the reflection map is the embedding
$$R~:~P(\R^2)^+~=~(S^1)^+ \emb O(\R^2)~;~
\begin{cases}
[\cos t,\sin t] \mapsto \begin{pmatrix} -\cos  2t & -\sin  2t\\
 -\sin  2t & \cos  2t \end{pmatrix}\\
* \mapsto \begin{pmatrix} 1 & 0 \\ 0 & 1 \end{pmatrix}
\end{cases}$$
with image consisting of all the reflections $\R^2 \to \R^2$ in lines
through the origin and the identity map.\\
\hfill\qed}
\end{example}

\begin{example}~{\rm
(i) The reflection map
$$\begin{array}{l}
R~:~P(U\oplus \R)/P(U)~=~U^{\infty} \to O(\R,U\oplus \R)~=~S(U\oplus \R)~;\\[1ex]
\hskip100pt
u \mapsto \big( \dfrac{-2u}{\vert u \vert^2+1},\dfrac{\vert u \vert^2-1}{\vert u \vert^2+1}\big)~,~
\infty \mapsto (0,1)
\end{array}$$
is a homeomorphism, with inverse
$$S(U\oplus \R) \to U^{\infty}~;~(x,y) \mapsto \dfrac{x}{y-1}~.$$
(ii) The composite
$$\begin{array}{l}
\xymatrix{\widetilde{R}~:~S(V) \ar[r] & P(V)
\subset P(V)^+~=~P(0\oplus V)/P(0) \ar[r]^-{\di{R}} &O(V)}~;\\[1ex]
\hskip100pt v \mapsto (w \mapsto w - 2\langle v,w \rangle v)
\end{array}$$
is the clutching function for the tangent $V$-bundle of $S(V\oplus \R)$
$$\tau_{S(V \oplus \R)}~:~S(V \oplus \R)~=~D(V)\cup_{S(V)}D(V) \to BO(V)~,$$
with
$$E(\tau_{S(V \oplus \R)})~=~D(V) \times V
\cup_{(x,v) \sim (x,\widetilde{R}(x)(v))}D(V) \times V~(x \in S(V),v \in V)~.$$
The classifying map for $\tau_{S(V \oplus \R)}$ fits into a fibration
$$\xymatrix@C+10pt{
S(V \oplus \R) \ar[r]^-{\di{\tau_{S(V \oplus \R)}}}& BO(V) \ar[r] & BO(V\oplus \R)}$$
corresponding to the stable isomorphism
$$\tau_{S(V \oplus \R)}\oplus \epsilon_{\R}~\cong~\epsilon_{V \oplus \R}$$
determined by the framed embedding $S(V \oplus \R) \subset V \oplus \R$.\\
}\hfill\qed
\end{example}

The reflection map $R:P(U\oplus V)/P(U) \to O(V,U\oplus V)$ is injective,
which we use to regard $P(U\oplus V)/P(U)$ as a subspace of $O(V,U\oplus V)$.

\begin{proposition}~ \label{james1}
{\rm (James \cite[Prop. 1.3 and Thm. 3.4]{imj5})}\\
The pair $(O(V,U \oplus V),P(U\oplus V)/P(U))$
is $2\,{\rm dim}(U)$-connected, so that
$$R_*~:~\pi_i(P(U \oplus V)/P(U)) \to \pi_i(O(V,U \oplus V))$$
is an isomorphism for $i<2\,{\rm dim}(U)$ and a surjection for
$i=2\,{\rm dim}(U)$.
\end{proposition}
\begin{proof} (Sketch) By induction on $\text{dim}(V)$, using the
Blakers-Massey theorem to compare the morphism induced by $R$ from the
homotopy exact sequence (in a certain range of dimensions)
of the homotopy cofibration sequence
$$P(U \oplus V)/P(U) \to P(U \oplus V \oplus \R)/P(U) \to (U \oplus V)^{\infty}$$
to the homotopy exact sequence of the fibration sequence
$$O(V,U\oplus V) \to O(V \oplus \R,U \oplus V \oplus \R) \to (U \oplus V)^{\infty}~.$$
\hfill\qed\end{proof}

For any finite-dimensional inner product space $V$ the $\ZZ_2$-action on
$LV$ restricts to
the antipodal involution on the unit sphere $S(LV)$, with
$$S(LV)/\ZZ_2~=~P(V)~,~LV^{\infty}/\ZZ_2~=~sP(V)~.$$
The homotopy cofibration sequence of Proposition \ref{cofib2} (iii)
$$\begin{array}{l}
S(LU) \to S(LU \oplus LV) \to LU^{\infty} \wedge S(LV)^+ \\[1ex]
\hskip50pt  \to LU^{\infty}\to
(LU \oplus LV)^{\infty} \to LU^{\infty} \wedge \Sigma S(LV)^+ \to \dots
\end{array}$$
is $\ZZ_2$-equivariant,  so that passing to the $\ZZ_2$-quotients
there is defined a homotopy cofibration sequence
$$\begin{array}{l}
P(U) \to P(U \oplus V) \to P(U \oplus V)/P(U)=
LU^{\infty}\wedge_{\ZZ_2}S(LV)^+\\[1ex]
\hskip25pt \to sP(U) \to sP(U\oplus V) \to
LU^{\infty}\wedge_{\ZZ_2}\Sigma S(LV)^+ \to \dots~.
\end{array}$$
If $U=\{0\}$ interpret $P(U \oplus V)/P(U)$ as $P(V)^+$.

Let $j={\rm dim}(U)$, $k={\rm dim}(V)$, so that
$$O(V,U \oplus V)~=~V_{j+k,k}~.$$
The reduced homology groups of the stunted projective space
$$P(U\oplus V)/P(U)~=~LU^{\infty} \wedge_{\ZZ_2}S(LV)^+$$
are the $Q$-groups
$$\begin{array}{ll}
\dot H_i(P(U\oplus V)/P(U))&=~Q^{[0,k-1]}_{i+j}(S^j\ZZ)\\[1ex]
&=~\begin{cases} 0&{\rm if}~i<j\\
\ZZ&{\rm if}~i=j~\hbox{\rm is even, or if}~i=j,~k=1\\
\ZZ_2&{\rm if}~i=j~\hbox{\rm is odd and}~k \geqslant 2~.
\end{cases}
\end{array}$$
The induced morphisms in homology
$$\begin{array}{l}
\theta_*~:~H_i(O(V,U\oplus V))~=~H_i(V_{j+k,k}) \\[1ex]
\hskip100pt \to H_i(P(U\oplus V)/P(U))~=~Q^{[0,k-1]}_{i+j}(S^j\ZZ)~(i>0)
\end{array}$$
are split surjections, which are isomorphisms for $i<2j$ by \ref{james1}
(cf. \ref{stiefel0} for $i \leqslant j$).

The Stiefel space
$$O(V,U\oplus V)~=~O(U\oplus V)/O(U)~=~V_{j+k,k}$$
fits into a fibration
$$O(V,U \oplus V) \to BO(U) \to BO(U\oplus V)~.$$
The canonical $U$-bundle $\eta(U)$ over $O(V,U \oplus V)$ is such that
$$\begin{array}{l}
E(\eta(U))~=~\{(f,x)\,\vert\, f \in O(V,U \oplus V),~x \in f(V)^{\perp}\subset U\oplus V\}\\[1ex]
\hphantom{E(\eta(U))~}=~O(U \oplus V)\times_{O(U)}U
\end{array}$$
with the canonical $U\oplus V$-bundle isomorphism
$$\delta\eta(U)~:~\eta(U) \oplus \epsilon_V~\cong~\epsilon_{U \oplus V}$$
defined by
$$\begin{array}{l}
\delta\eta(U)~: ~E(\eta(U) \oplus \epsilon_V)~=~
E(\eta(U)) \times V\\[1ex]
\hskip50pt \to E(\epsilon_{U \oplus V})~=~O(V,U\oplus V) \times U\oplus V~;~
(f,x,v) \mapsto (f,x+f(v))~.
\end{array}$$
A map $c:X \to O(V,U\oplus V)$ classifies
a $U$-bundle $\xi:X \to BO(U)$ with a stable isomorphism
$\delta\xi:\xi\oplus \epsilon_V \cong \epsilon_{U\oplus V}$, where
$$E(\xi)~=~\{(x \in X,y \in c(x)^{\perp} \subset U\oplus V)\} \subset
E(\epsilon_{U \oplus V})~=~X \times U\oplus V~.$$
The adjoint $\ZZ_2$-equivariant map
$$F_c~:~LV^{\infty} \wedge X^+ \to (LU \oplus LV)^{\infty}~;~(v,x) \mapsto
c(x)(v)$$
represents the $\ZZ_2$-equivariant Euler class of $L\xi:X \to BO^{\ZZ_2}(U)$
$$\gamma^{\ZZ_2}(L\xi)~=~F_c \in \omega^0_{\ZZ_2}(X;-L\xi)~=~
\{X;LU^{\infty}\}_{\ZZ_2}~.$$
The standard pair $(1,\epsilon_U)$ is classified by the constant map
$$0~:~X \to O(V,U\oplus V)~;~x \mapsto (v \mapsto (0,v))$$
with adjoint $\ZZ_2$-equivariant map
$$F_0~:~LV^{\infty} \wedge X^+ \to (LU \oplus LV)^{\infty}~;~(v,x) \mapsto (0,v)~.$$
The adjoint $\ZZ_2$-equivariant maps $F_c,F_0$ are such that
$$F_c(0,x)~=~F_0(0,x)~=~(0,0) \in (LU \oplus LV)^{\infty}~.$$
\begin{definition} ~\label{local}
{\rm (Crabb \cite[2.6]{crabb})\\
(i) The {\it local obstruction} of $c:X\to O(V,U\oplus V)$
is the relative difference $\ZZ_2$-equivariant map
$$\begin{array}{l}
\theta(c)~=~\delta(F_c,F_0)~:~
\Sigma S(LV)^+ \wedge X^+ \to (LU \oplus LV)^{\infty}~;\\[1ex]
\hspace*{50pt}
([t,v],x) \mapsto \begin{cases}
(0,[1-2t,v])&{\rm if}~0 \leqslant t \leqslant 1/2\\
c(x)[2t-1,v]&{\rm if}~1/2 \leqslant t \leqslant 1~.
\end{cases}
\end{array}$$
The stable $\ZZ_2$-equivariant homotopy class
$$\theta(c)\in \{\Sigma S(LV)^+ \wedge X^+;(LU \oplus LV)^{\infty}\}_{\ZZ_2}$$
has image
$$[\theta(c)]~=~\gamma^{\ZZ_2}(L\xi)-\gamma^{\ZZ_2}(L\epsilon_U)
\in \{X^+;LU^{\infty}\}_{\ZZ_2}~.$$
(ii) Let
$$\theta'(c) \in \{X^+;S(LV)^+\wedge LU^{\infty}\}_{\ZZ_2}~=~
\{X^+;S(LV)^+\wedge_{\ZZ_2}LU^{\infty}\}$$
be the $\ZZ_2$-equivariant $S$-dual of $\theta(c) \in
\{\Sigma S(LV)^+ \wedge X^+;(LU \oplus LV)^{\infty}\}_{\ZZ_2}$.\\
\hfill\qed}
\end{definition}

\begin{example}{\rm For any inner product space $U$ the homeomorphism
$$\begin{array}{c}
c~:~S(U \oplus \R) \to O(\R,U\oplus \R)~;~(u,v) \mapsto (t \mapsto (tu,tv))\\[1ex]
(t \in V=\R,(u,v) \in S(U\oplus \R))
\end{array}$$
classifies the tangent $U$-bundle
$$\xi~=~\tau_{S(U\oplus \R)}~:~S(U\oplus \R) \to BO(U)$$
with the $U\oplus \R$-bundle isomorphism
$$\delta\xi~:~\epsilon_\R \oplus \xi ~\cong~ \epsilon_{U\oplus \R}$$
determined by the embedding $S(U\oplus \R) \subset U \oplus \R$.
The local obstruction of $c$ is
$$\theta(c)~=~1 \in \{S(U\oplus \R)^+;S(L\R)^+\wedge LU^{\infty}\}_{\ZZ_2}~=~
\{S(U\oplus \R)^+;U^{\infty}\}~=~\ZZ~.$$
The map
$$\begin{array}{l}
F_c~:~\R^{\infty} \wedge S(U \oplus \R)^+ \to
\R^{\infty} \wedge U^{\infty} \wedge S(U\oplus \R)^+~;\\[1ex]
(t,u,v) \mapsto (c(u,v)(t),(u,v))~=~(tv,tu,u,v)~(t \in \R,(u,v) \in S(U\oplus \R))~.
\end{array}$$
has stable homotopy class
$$F_c~=~1 \in \{S(U\oplus \R)^+;U^{\infty} \wedge S(U\oplus \R)^+\}~=~
\{S(U\oplus \R)^+;U^{\infty}\}~=~\ZZ~,$$
and $\R$-coefficient quadratic construction
$$\begin{array}{l}
\psi_\R(F_c)~=~\theta(c) F_c~=~1\\[1ex]
\hskip25pt
\in \{S(U\oplus \R)^+;S(L\R)^+ \wedge U^{\infty}\wedge LU^{\infty} \wedge S(U\oplus \R)^+\wedge S(U\oplus \R)^+\}_{\ZZ_2}\\[1ex]
=~\{S(U\oplus \R)^+;U^{\infty}\wedge U^{\infty} \wedge S(U\oplus \R)^+\wedge S(U\oplus \R)^+\}~=~\ZZ~.
\end{array}$$}
\hfill\qed
\end{example}

\begin{example}~ \label{orthogonal} {\rm Suppose $U=\{0\}$,
and $c \in O(V)$.  The
$\ZZ_2$-equivariant homotopy class of the local obstruction
$\theta(c):\Sigma S(LV)^+ \to LV^{\infty}$ is given by
Proposition \ref{00} to be
$$\begin{array}{ll}
\theta(c)&=~\hbox{\rm bi-degree}(\delta(F_c,F_0))\\[1ex]
&=~(\dfrac{{\rm degree}(F_c)-{\rm degree}(F_0)}{2},
{\rm degree}(G_c)-{\rm degree}(G_0))\\[1ex]
&=~(\dfrac{{\rm degree}(c)-1}{2},1-1)\\[1ex]
&=~\begin{cases}
(0,0) &{\rm if}~{\rm det}(c)=1\\
(-1,0) &{\rm if}~{\rm det}(c)=-1
\end{cases} \in [\Sigma S(LV)^+,LV^{\infty}]_{\ZZ_2}~=~\ZZ\oplus \ZZ
\end{array}$$
with
$$G_c~=~G_0~=~1~:~(LV^{\infty})^{\ZZ_2}~=~\{0\}^+ \to \{0\}^+$$
the fixed point maps.\\
\hfill\qed}
\end{example}

\begin{proposition}~\label{twist}
Let $\xi:X \to BO(U)$ be a $U$-vector bundle with a $V$-stable trivialization
$\delta\xi: \xi \oplus \epsilon_V \cong \epsilon_{U \oplus V}$.\\
{\rm (i)} The map
$$T(\delta\xi)~:~T(\xi \oplus \epsilon_V)~=~V^{\infty} \wedge T(\xi)
\to T(\epsilon_{U \oplus V})~=~V^{\infty} \wedge U^{\infty} \wedge X^+$$
is a homeomorphism inducing an isomorphism
$$T(\delta\xi)^*~:~\{X;S^0\} \xymatrix{\ar[r]^{\di{\cong}}&}
\{T(\xi);U^{\infty}\}~.$$
The stable cohomotopy Thom class (\ref{eulerthom}) of $\xi$
$$u(\xi)~=~T(\delta\xi)^*(1) \in \{T(\xi);U^{\infty}\}~=~\omega^0(X;\xi-\epsilon_U)$$
is represented by the composite
$$\begin{array}{l}
F~=~(1\wedge p)T(\delta\xi)~:~T(\epsilon_V \oplus \xi)~=~V^{\infty} \wedge T(\xi)
\xymatrix{\ar[r]^-{\di{T(\delta\xi)}}&}\\[1ex]
\hskip75pt
T(\epsilon_{U \oplus V})~=~V^{\infty} \wedge U^{\infty} \wedge X^+
\xymatrix{\ar[r]^-{\di{1 \wedge p}}&}V^{\infty} \wedge U^{\infty}
\end{array}
$$
with $p:X^+ \to S^0$ the projection sending $X$ to the non-base-point of $S^0$.
The adjoint map
$$F_c~:~V^{\infty} \wedge X^+ \to V^{\infty}\wedge U^{\infty}~;~
(v,x) \mapsto c(x)(v)$$
is the composite
$$\xymatrix@C+30pt{F_c~:~V^{\infty} \wedge X^+ \ar[r]^-{\di{1 \wedge z_{\xi}}} &
V^{\infty} \wedge T(\xi) \ar[r]^-{\di{F}} & V^{\infty}\wedge U^{\infty}}$$
with $z_{\xi}:X^+ \to T(\xi)$ the zero section. The images of
$u(\xi) \in \{T(\xi);U^{\infty}\}$ in the commutative square
$$\xymatrix{\{T(\xi);U^{\infty}\} \ar[r] \ar[d]^-{\di{z_{\xi}^*}} &
\dot H^j(T(\xi)) \ar[d]^-{\di{z_{\xi}^*}} \\
\{X^+;U^{\infty}\} \ar[r]  &H^j(X)}$$
are the Thom class $[u(\xi)] \in \dot H^j(T(\xi))$,
the stable cohomotopy Euler class (\ref{eulerthom})
$$z_{\xi}^*u(\xi)~=~\gamma(\xi) \in \{X^+;U^{\infty}\}$$
and the Euler class $[\gamma(\xi)] \in H^j(X)$.
The adjoint $\ZZ_2$-equivariant maps
$$F_c~,~F_0~:~LV^{\infty} \wedge X^+ \to (LV \oplus LU)^{\infty}$$
with $0$ the constant map
$$0~:~X \to O(V,U \oplus V)~;~x \mapsto (v \mapsto (0,v))$$
represent the stable $\ZZ_2$-equivariant cohomotopy Euler classes
$$\gamma^{\ZZ_2}(L\xi)~=~F_c~,~\gamma^{\ZZ_2}(L\epsilon_U)~=~F_0
\in \{X^+;LU^{\infty}\}_{\ZZ_2}~.$$
{\rm (ii)} The geometric Hopf invariants of $F$ and $F_c$
$$\begin{array}{l}
h_V(F) \in
\{\Sigma S(LV)^+ \wedge T(\xi);LV^{\infty}\wedge(U^{\infty} \wedge U^{\infty})\}_{\ZZ_2} ~,\\[1ex]
h_V(F_c) \in
\{\Sigma S(LV)^+ \wedge X^+;LV^{\infty}\wedge (U^{\infty} \wedge U^{\infty})\}_{\ZZ_2}
\end{array}$$
are related by
$$h_V(F_c)~=~h_V(F)(1 \wedge z_{\xi})~.$$
{\rm (iii)} The $\ZZ_2$-equivariant stable homotopy class of the
geometric Hopf invariant of $F_c$  is given by
$$\begin{array}{l}
h_V(F_c)~=~(\theta(c) \wedge F_c)(1\wedge \Delta_X)\\[1ex]
\hskip25pt \in
\{\Sigma S(LV)^+ \wedge X^+;LV^{\infty}\wedge
(U^{\infty}\wedge U^{\infty})\}_{\ZZ_2}\\[1ex]
\hskip25pt
=~\{\Sigma S(LV)^+ \wedge V^{\infty} \wedge X^+;(LU \oplus LV)^{\infty} \wedge  (U\oplus V)^{\infty}\}_{\ZZ_2}
\end{array}$$
and is represented by the $\ZZ_2$-equivariant stable map
$$\begin{array}{l}
(\theta(c) \wedge F_c)(1\wedge \Delta_X)~:\\[1ex]
\Sigma S(LV)^+ \wedge V^{\infty} \wedge X^+  \to
(LU \oplus LV)^{\infty} \wedge  (U\oplus V)^{\infty}~;\\[1ex]
\hspace*{100pt} ([t,u],v,x) \mapsto (\theta(c)([t,u],x),c(x)(v))~.
\end{array}$$
{\rm (iv)} The isomorphism of exact sequences
$$\xymatrix@C-18pt{\{X^+;S(LV)^+\wedge LU^{\infty}\}_{\ZZ_2}
\ar[r] \ar[d]^-{\di{T(\delta\xi)^*\wedge \kappa_U}}_-{\di{\cong}} &
\{X^+; LU^{\infty}\}_{\ZZ_2}
\ar[r] \ar[d]^-{\di{T(\delta\xi)^*\wedge \kappa_U}}_-{\di{\cong}} &
\{X^+; (LU\oplus LV)^{\infty}\}_{\ZZ_2}
\ar[d]^-{\di{T(\delta\xi)^*\wedge \kappa_U}}_-{\di{\cong}} \\
\{T(\xi);S(LV)^+\wedge (U^{\infty} \wedge U^{\infty})\}_{\ZZ_2}
\ar[r]  &\{T(\xi); U^{\infty} \wedge U^{\infty}\}_{\ZZ_2}\ar[r] &
\{T(\xi); LV^{\infty}\wedge (U^{\infty} \wedge U^{\infty})\}_{\ZZ_2}}$$
sends $\theta'(c) \in \{X;S(LV)^+\wedge LU^{\infty}\}_{\ZZ_2}$
to the stable geometric Hopf invariant
$h'_V(F) \in \{T(\xi);S(LV)^+\wedge (U^{\infty} \wedge U^{\infty})\}_{\ZZ_2}$.
The stable geometric Hopf invariant of $F_c$ is given by
$$h'_V(F_c)~=~h'_V(F)z_{\xi} \in \{X^+;S(LV)^+\wedge (U^{\infty} \wedge U^{\infty})\}_{\ZZ_2}~.$$
{\rm (v)} The isomorphism
$$T(\delta\xi)^*\wedge \kappa_U~:~
\{X^+;S(LV)^+\wedge_{\ZZ_2}LU^{\infty}\} \xymatrix{\ar[r]^-{\di{\cong}}&}
\{T(\xi);S(LV)^+\wedge_{\ZZ_2}(U^{\infty}\wedge U^{\infty})\}$$
sends the local obstruction
$$\theta(c) \in
\{X^+;P(U\oplus V)/P(U)\}~=~\{X^+;S(LV)^+\wedge_{\ZZ_2}LU^{\infty}\}$$
to the $V$-coefficient quadratic construction $\psi_V(F)$, so we
can identify
$$\theta(c)~=~\psi_V(F) \in \{X^+;S(LV)^+\wedge_{\ZZ_2}LU^{\infty}\}~=~
\{T(\xi);S(LV)^+\wedge_{\ZZ_2}(U^{\infty}\wedge U^{\infty})\}~.$$
{\rm (vi)} The $V$-coefficient quadratic construction on $F_c$ is the product
of $F_c=\gamma(\xi) \in \{X^+;U^{\infty}\}$ and $\theta(c)$, that is
$$\psi_V(F_c)~=~(\theta(c)\times \gamma(\xi))\Delta_X \in
\{X^+;S(LV)^+ \wedge_{\ZZ_2} (U^{\infty}\wedge U^{\infty})\}$$
inducing
$$\psi_V(F_c)~:~H_i(X) \to Q^{[0,k-1]}_i(S^j\ZZ)~.$$
The local obstruction stable map
$$\theta~:~O(V,U\oplus V) \sto S(LV)^+\wedge_{\ZZ_2} LU^{\infty}~=~P(U\oplus V)/P(U)$$
is $2j$-connected (\ref{james1}), so that for $i<2j$
$$\theta(c)_*~=~c_*~:~H_i(X) \to H_i(O(V,U\oplus V))~=~Q_{i+j}^{[0,k-1]}(S^j\ZZ)~.$$
Now $F_c$ induces the $\ZZ$-module chain map
$$F_c~=~[\gamma(\xi)] \cap-~:~C(X) \to S^j\ZZ$$
and $\theta(c)$ induces a $\ZZ$-module chain map
$$\theta(c)~:~C(X) \to W[0,k-1]\otimes_{\ZZ[\ZZ_2]} S^j(\ZZ,-1)~,$$
so that $\psi_V(F_c)$ induces the $\ZZ$-module chain map
$$\begin{array}{l}
\psi_V(F_c)~:~C(X) \xymatrix{\ar[r]^-{\di{\Delta_X}} &}
C(X)\otimes_{\ZZ} C(X)\\[1ex]
\hskip100pt
\xymatrix@C+20pt{\ar[r]^-{\di{\theta(c)\otimes \gamma(\xi)}} &}
W[0,k-1]\otimes_{\ZZ[\ZZ_2]} (S^j\ZZ\otimes_{\ZZ}S^j\ZZ)~.
\end{array}$$
{\rm (vii)} Let $U=U' \oplus \R$. If $\xi=\xi' \oplus \epsilon_{\R}$
for some $U'$-bundle $\xi':X \to BO(U')$ then
$$\begin{array}{l}
F_c~\simeq~*~:~ V^{\infty}\wedge X^+ \to (U \oplus V)^{\infty}~
(\text{non-equivariantly})~,\\[1ex]
\gamma(\xi)~=~0 \in \{X^+;U^{\infty}\}~,\\[1ex]
h_V(F_c)~=~0 \in
\{\Sigma S(LV)^+ \wedge X^+;U^{\infty}\wedge (LU \oplus LV)^{\infty}\}_{\ZZ_2}~,\\[1ex]
\psi_V(F_c)~=~0 \in
\{X^+;S(LV)^+ \wedge_{\ZZ_2} (U^{\infty}\wedge U^{\infty})\}~.
\end{array}$$
If $X$ is a $j$-dimensional $CW$ complex then
$\xi=\xi' \oplus \epsilon_{\R}$ if and only if
$$\gamma(\xi)~=~0 \in \{X^+;U^{\infty}\}~=~H^j(X)~.$$
\end{proposition}
\begin{proof} (i)+(ii) By construction.\\
(iii) The geometric Hopf invariant of $F_c$ is the relative difference
$h_V(F_c)=\delta(p,q)$ of the $\ZZ_2$-equivariant maps
$$p,q~:~LV^{\infty} \wedge V^{\infty} \wedge X \to (LU \oplus LV)^{\infty}
    \wedge (U \oplus V)^{\infty}$$
defined by
$$\begin{array}{l}
p(v_1,v_2,x)\\[1ex]
=~( \di{c(x)(v_1+v_2)-c(x)(-v_1+v_2) \over 2},\di{c(x)(v_1+v_2)+c(x)(-v_1+v_2) \over 2})\\[1ex]
\hphantom{q(v_1,v_2,x)~}=~(c(x)(v_1),c(x)(v_2))~,\\[1ex]
q(v_1,v_2,x)~=~(v_1,c(x)(v_2))
\end{array}$$
with
$$p(0,v,x)~=~q(0,v,x)~=~(0,c(x)(v)) ~~(v \in V,x \in X)~.$$
By Example \ref{inj}
$$\begin{array}{l}
h_V(F_c)~=~\delta(p,q)~=~(\theta(c)\wedge F_c)(1 \wedge \Delta_X)~:\\[1ex]
\Sigma S(LV)^+ \wedge V^{\infty} \wedge X \to
LV^{\infty} \wedge V^{\infty}~;\\[1ex]
\hskip100pt
([t,u],v,x) \mapsto (\theta(c)([t,u],x),c(x)(v))~.
\end{array}$$
(iv) The isomorphism $T(\delta\xi)^*\wedge \kappa_U$ sends
$F_c,F_0 \in \{X^+;LU^{\infty}\}_{\ZZ_2}$ to
$$\begin{array}{l}
(T(\delta\xi)^*\wedge \kappa_U)(F_c)~=~(F \wedge F)\Delta_{T(\xi)}~=~p~,\\[1ex]
(T(\delta\xi)^*\wedge \kappa_U)(F_0)~=~\Delta_{U^{\infty}}F~=~q
\in \{T(\xi);U^{\infty} \wedge U^{\infty}\}_{\ZZ_2}~.
\end{array}$$
so that the stable relative differences are such that
$$\begin{array}{ll}
(T(\delta\xi)^*\wedge \kappa_U)(\theta'(c))&=~
(T(\delta\xi)^*\wedge \kappa_U)(\delta'(F_c,F_0))\\[1ex]
&=~\delta'(p,q)~=~h'_V(F) \in
\{T(\xi);S(LV)^+ \wedge U^{\infty} \wedge U^{\infty}\}_{\ZZ_2}~.
\end{array}$$
(v) Immediate from (iv).\\
(vi)+(vii) By construction.\\
\hfill\qed\end{proof}

Given a map $c:(X,Y) \to (O(V,U \oplus V),\{0\})$ define the
$\ZZ_2$-equivariant section of the trivial $\ZZ_2$-equivariant
bundle $\epsilon_{LU \oplus LV}$ over $S(LV)\times X$
$$c'~:~S(LV) \times X \to S(\epsilon_{LU \oplus LV})~=~
S(LV)\times X \times LU \oplus LV~;~(v,x) \mapsto (v,x,c(x)(v))~.$$
The local obstruction of $c$ is the rel $S(LV)\times Y$ $\ZZ_2$-equivariant
difference class
$$\begin{array}{l}
\theta(c)~=~\delta(c',0')\\[1ex]
\in \omega^{-1}_{\ZZ_2}(S(LV)\times (X,Y);-\epsilon_{LU \oplus LV})~
=~\{\Sigma S(LV)^+ \wedge X/Y;(LU\oplus LV)^{\infty}\}_{\ZZ_2}~.
\end{array}$$

\begin{proposition}~ \label{localR}
{\rm (Crabb \cite[5.13, 2.7]{crabb})}\\
{\rm (i)} The local obstruction determines a stable map
$$\theta~:~O(V,U\oplus V) \sto P(U \oplus V)/P(U)$$
inducing a function
$$\begin{array}{l}
\theta~:~[X/Y,O(V,U \oplus V)] \to
\{\Sigma S(LV)^+ \wedge X/Y;(LU \oplus LV)^{\infty}\}_{\ZZ_2}\\[1ex]
\hspace*{125pt}=~\{X/Y;S(LU \oplus LV)/S(LU)\}_{\ZZ_2}\\[1ex]
\hspace*{125pt}=~\{X/Y;S(LV)^+\wedge LU^{\infty}\}_{\ZZ_2}\\[1ex]
\hspace*{125pt}=~\{X/Y;P(U \oplus V)/P(U)\}\\[1ex]
\hspace*{125pt}=~\omega^{-1}_{\ZZ_2}(S(LV)\times (X,Y);-\epsilon_{LU \oplus LV})~.
\end{array}$$
The local obstruction of a map $c:X/Y \to O(V,U\oplus V)$ is a stable
$\ZZ_2$-equivariant map $\theta(c):X/Y \to S(LV)^+\wedge LU^{\infty}$ such that
the composite with the $\ZZ_2$-equivariant map
$$s_{LV} \wedge 1~:~S(LV)^+ \wedge LU^{\infty} \to S^0 \wedge
LU^{\infty}~=~LU^\infty$$
is $F_c-F_0 \in \{X/Y;LU^{\infty}\}_{\ZZ_2}$, with
$F_c:LV^{\infty} \wedge X/Y \to (LU \oplus LV)^{\infty}$ the adjoint map
(\ref{adjoint}) and $0$ the constant map
$$0~:~X/Y \to O(V,U\oplus V)~;~x \mapsto (v \mapsto (0,v))~.$$
{\rm (ii)} The composite stable map
$$\xymatrix{\theta \circ R~:~P(U \oplus V)/P(U)
\ar[r]^-{\di{R}} & O(V,U \oplus V)
 \ar[r]^-{\di{\theta}}  & P(U \oplus V)/P(U)}$$
is the identity, so that
$$\theta\circ R~:~[X/Y;P(U \oplus V)/P(U)] \to \{X/Y;P(U \oplus V)/P(U)\}$$
is the stabilization map. In particular, since
$(\,O(V,U \oplus V)\,,\, P(U\oplus V)/P(U)\,)$
is $2\,{\rm dim}(U)$-connected and $O(V,U \oplus V)$ is $({\rm dim}(U)-1)$-connected,
$$\theta~:~\pi_i(O(V,U \oplus V)) \to \pi_i(P(U \oplus V)/P(U))$$
is an isomorphism for $i \leqslant  {\rm dim}(U)$, and
$$\theta~:~\omega_i(O(V,U \oplus V)) \to \omega_i(P(U \oplus V)/P(U))$$
is an isomorphism for $i \leqslant  2\,{\rm dim}(U)-1$.\\
{\rm (iii)} If $c:(X,Y) \to (O(V),\{0\})$ is a map such that $X\backslash Y$ is a
manifold and the $\ZZ_2$-equivariant section of the trivial $S(LV)
\times S(LV)$-bundle over $S(LV)\times X$
$$\begin{array}{l}
d~=~(-j,c)~:~S(LV)\times X  \to S(LV) \times X  \times  S(LV) \times S(LV)~ ;\\[1ex]
\hskip175pt (v,x) \mapsto (v,x,-v,c(x)(v))
\end{array}$$
is transverse regular at $S(LV) \times X \times \Delta_{S(LV)}$ then
$$\begin{array}{ll}
C&=~d^{-1}(S(LV)\times X \times \Delta_{S(LV)})\\[1ex]
&=~\{(v,x) \,\vert\, c(x)(v)=-v \in S(LV)\} \subset S(LV)\times (X\backslash Y)
\end{array}$$
is a submanifold of codimension ${\rm dim}(V)-1$ with normal bundle
$$\nu_{C \subset S(LV)\times (X \backslash Y)}~=~
(d\vert_C)^*\nu_{\Delta_{S(LV)} \subset S(LV) \times S(LV)}~=~
(d\vert_C)^*\tau_{S(LV)}$$
such that $\nu_{C \subset S(LV)\times (X\backslash Y)}\oplus \epsilon_{\R}
\cong \epsilon_{LV}$, and
$$\begin{array}{l}
\xymatrix{
\theta(c)~:~\Sigma S(LV)^+\wedge X/Y \ar[r]^-{\di{\Sigma F}}&}\\[1ex]
\hskip75pt
\xymatrix{\Sigma T(\nu_{C \subset S(LV) \times (X\backslash Y)})~=~
C^+ \wedge LV^{\infty}\ar[r] &LV^{\infty}}
\end{array}
$$
with $F:S(LV)^+\wedge X/Y\to T(\nu_{C \subset S(LV) \times (X\backslash
Y)})$ the adjunction
Umkehr map of the codimension 0 embedding $E(\nu_{C \subset
S(LV) \times (X \backslash Y)}) \subset S(LV) \times (X \backslash Y)$.
\end{proposition}
\begin{proof} (i) The local obstruction of $c \in O(V,U\oplus V)$ is
the $\ZZ_2$-equivariant pointed map
$$\theta(c)~=~\delta(F_c,F_0)~:~\Sigma S(LV)^+ \to (LU\oplus LV)^{\infty}$$
with $F_c:LV^{\infty} \to (LU\oplus LV)^{\infty}$ the $\ZZ_2$-equivariant
adjoint (\ref{Z2adjoint1}) and $0\in O(V,U \oplus V)$ the linear isometry
$$0~:~V \to U \oplus V~;~v \mapsto (0,v)~.$$
For any $c:X/Y \to O(V,U \oplus V)$ the local obstruction
defines a $\ZZ_2$-equivariant map
$$\theta(c)~:~\Sigma S(LV)^+ \wedge X/Y \to
(LU \oplus LV)^{\infty}~;~(v,x) \mapsto c(x)(v)~.$$
The $\ZZ_2$-equivariant $S$-duality isomorphism
$$\{\Sigma S(LV)^+ \wedge X/Y;(LU \oplus LV)^{\infty}\}_{\ZZ_2} \to
\{X/Y;S(LU \oplus LV)/S(LU)\}_{\ZZ_2}$$
sends $\theta(c)$ to the composite stable $\ZZ_2$-equivariant map
$$\begin{array}{l}
\xymatrix@C+30pt{
(LU \oplus LV)^{\infty} \wedge X/Y
\ar[r]^-{\di{\alpha_{LU \oplus LV}\wedge 1}} &
\Sigma S(LU \oplus LV)^+ \wedge X/Y}\\[1ex]
\hspace*{100pt}\xymatrix@C+20pt{\ar[r]^-{\di{\Delta\wedge 1}} &
\Sigma S(LV)^+ \wedge  S(LU \oplus V))/S(LU) \wedge X/Y}\\[1ex]
\hspace*{100pt}
\xymatrix@C+20pt{\ar[r]^-{\di{1 \wedge \theta(c)}} &
(LU \oplus LV)^{\infty} \wedge   S(L(U \oplus V))/S(LU)}
\end{array}$$
with
$$\begin{array}{l}
\Delta~:~\Sigma S(LU \oplus LV)^+ \to
\Sigma S(LV)^+ \wedge S(LU \oplus LV)/S(LU)~;\\[1ex]
\hspace*{100pt}
(t,u,v) \mapsto ((t,\dfrac{v}{\Vert v \Vert}),[u,v])~.
\end{array}$$
(ii) Consider first the special case $U=\{0\}$, $V=\R$. The
reflection map is
$$R~:~P(\R)^+=\{*,[1]\} \to O(\R)=\{\pm 1\}~;~
* \mapsto 1~,~[1] \mapsto -1~.$$
By Example \ref{orthogonal}
the local obstruction stable map $\theta:O(\R) \to P(\R)^+$ induces
$$\begin{array}{l}
\theta~:~[S^0,O(\R)]~=~\{\pm 1\}
 \to \{\Sigma S(L\R)^+;L\R^{\infty}\}_{\ZZ_2}~=~\pi_0^S(P(\R)^+)~=~\ZZ~;\\[1ex]
 \hspace*{100pt} c \mapsto \delta(c,1)~=~\begin{cases} 0&{\rm if}~c=1\\
1&{\rm if}~c=-1~. \end{cases}
\end{array}$$
Thus $\theta \circ R \simeq {\rm id}:P(\R)^+ \sto P(\R)^+$.
For arbitrary $V$ and any $v \in P(V)$ there is defined a commutative
diagram
$$\xymatrix{
P(\R)^+ \ar[r]^-{\di{R}} \ar[d]^-{\di{v}} &
O(\R) \ar[r]^-{\di{\theta}} \ar[d]^-{\di{v}} &
P(\R)^+ \ar[d]^-{\di{v}} \\
P(V)^+ \ar[r]^-{\di{R}}  &
O(V) \ar[r]^-{\di{\theta}} & P(V)^+}$$
so that $\theta  \circ R\simeq {\rm id}:P(V)^+ \to P(V)^+$.
Similarly for arbitrary $U,V$.\\
(iii) This is a $\ZZ_2$-equivariant special case of \ref{difcon3} (vi).\\
\hfill\qed\end{proof}

\begin{proposition}~{\rm (Crabb \cite[2.9,2.10]{crabb})}
Let $U,V$ be finite-dimensional inner product spaces, and
let $\xi:X \to BO(U\oplus V)$ be a $U\oplus V$-bundle.\\
{\rm (i)} The $\ZZ_2$-equivariant homotopy cofibration sequence
$$S(LV)^+ \to S^0 \to LV^{\infty} \to \Sigma S(LV)^+$$
induces a long exact sequence of stable $\ZZ_2$-cohomotopy groups
$$\begin{array}{l}
\dots \to \omega^0_{\ZZ_2}(X;\epsilon_{LV}-L\xi)
\xymatrix@C+25pt{\ar[r]^-{\di{-\otimes \gamma^{\ZZ_2}(\epsilon_{LV})}}&}
\omega^0_{\ZZ_2}(X;-L\xi) \\[1ex]
\hskip125pt \to \omega^0(X \times P(V);-\xi\times H_\R) \to \dots
\end{array}
$$
with $H_\R$ the Hopf $\R$-bundle over $P(V)$ (\ref{Hopfbundle}) and
$$\omega^0(X \times P(V);-\xi\times H_\R)~=~\{T(L\eta);LU^{\infty} \wedge \Sigma S(LV)^+
\wedge LW^{\infty}\}_{\ZZ_2}$$
for any $W$-bundle $\eta:X \to BO(W)$ such that $\xi\oplus \eta=\epsilon_{U \oplus V\oplus W}$.\\
{\rm (ii)} If $\xi \cong \xi' \oplus \epsilon_V$ for some $U$-bundle
$\xi':X \to BO(U)$ then
$$\gamma^{\ZZ_2}(L\xi)~=~\gamma^{\ZZ_2}(L\xi')\otimes \gamma^{\ZZ_2}(\epsilon_{LV})
\in {\rm im}(\omega^0_{\ZZ_2}(X;\epsilon_{LV}-L\xi) \to \omega^0_{\ZZ_2}(X;-L\xi))$$
{\rm (iii)} If $X$ is an $m$-dimensional $CW$ complex and
$m<2{\rm dim}(U)$ then $\xi \cong \xi' \oplus \epsilon_V$ for some $U$-bundle
$\xi':X \to BO(U)$ if and only if
$$\begin{array}{l}
\gamma^{\ZZ_2}(L\xi)\in {\rm im}(-\otimes \gamma^{\ZZ_2}(\epsilon_{LV}):
\omega^0_{\ZZ_2}(X;\epsilon_{LV}-L\xi) \to \omega^0_{\ZZ_2}(X;-L\xi))\\[1ex]
\hskip50pt
=~{\rm ker}(\omega^0_{\ZZ_2}(X;-L\xi) \to \omega^0(X \times P(V);-\xi\times H_\R))~.
\end{array}$$
\hfill\qed
\end{proposition}

\begin{definition}~ {\rm
Let $H(V)$ be the space of homotopy equivalences $h:S(V) \to S(V)$,
regarded as a pointed space with base point $1:S(V) \to S(V)$.}
\hfill\qed
\end{definition}

Following Crabb \cite[\S3]{crabb} we shall now factor the local
obstruction stable map $\theta: O(V) \to P(V)^+$ through the inclusion
$J:O(V) \subset H(V)$.

\begin{definition}~ {\rm
(i) Let $H^{\ZZ_2}(LV\oplus W) \subset H(LV\oplus W)$ be the subspace of
$\ZZ_2$-equivariant homotopy equivalences $h:S(LV\oplus W) \to S(LV\oplus W)$.\\
(ii) Let $H^{\ZZ_2}(LV;W) \subset H^{\ZZ_2}(LV \oplus W)$ be the subspace of
$\ZZ_2$-equivariant homotopy equivalences $h:S(LV\oplus W) \to
S(LV\oplus W)$ which restrict to the identity $h\vert=1:S(W) \to S(W)$
on the $\ZZ_2$-fixed point set.}\\
\hfill\qed
\end{definition}

For any inner product spaces $V,W$ the homeomorphism defined in
Proposition \ref{pushout2} (ii)
$$\lambda_{LV,W}~:~S(LV)*S(W) \to S(LV \oplus W)$$
is $\ZZ_2$-equivariant. The spaces of
$\ZZ_2$-equivariant homotopy equivalences fit into a fibration sequence
$$\xymatrix{H^{\ZZ_2}(LV;W) \ar[r] & H^{\ZZ_2}(LV \oplus W) \ar[r]^-{\di{\rho}} &
H(W)}$$
with $\rho$ the $\ZZ_2$-fixed point map, which is split by
$$\sigma~:~H(W) \to H^{\ZZ_2}(LV \oplus W)~;~g \mapsto
\lambda_{LV,W}(1_{LV}*g)(\lambda_{LV,W})^{-1}~.$$
For any (reasonable) pointed space $X$ there is induced a split short exact sequence
of homotopy groups
$$\xymatrix{
0 \ar[r] & [X,H^{\ZZ_2}(LV;W)] \ar[r] &
[X,H^{\ZZ_2}(LV\oplus W)] \ar[r]^-{\di{\rho}} &[X,H(W)] \ar[r] &0}$$
Any map $g:X \to H^{\ZZ_2}(LV\oplus W)$ such that
$\rho(g):X \to H(W)$ is null-homotopic can be compressed (up to homotopy)
to a map
$$X \to H^{\ZZ_2}(LV;W) \subset H^{\ZZ_2}(LV\oplus W)~.$$
An element $h \in H^{\ZZ_2}(LV;W)$ suspends to a pointed $\ZZ_2$-equivariant
homotopy equivalence
$$sh~:~sS(LV\oplus W)~=~(LV\oplus W)^{\infty} \to sS(LV\oplus W)~=~
(LV\oplus W)^{\infty}$$
which is the identity on $W^{\infty} \subset (LV\oplus W)^{\infty}$.
Given a map $h:X \to H^{\ZZ_2}(LV;W)$ let
$$F_0~,~F_h~:~(LV\oplus W)^{\infty} \wedge X \to (LV\oplus W)^{\infty}$$
be the $\ZZ_2$-equivariant maps defined by
$$F_0(v,w,x)~=~(v,w)~,~F_h(v,w,x)~=~sh(x)(v,w)$$
such that $F_0(0,w,x)=F_h(0,w,x)=(0,w)$. Use the $\ZZ_2$-equivariant
difference map
$$\delta(F_h,F_0)~:~\Sigma S(LV)^+\wedge W^{\infty} \wedge X \to
LV^{\infty}\wedge W^{\infty}$$
to define a function
$$\begin{array}{l}
\zeta^{\ZZ_2}~:~[X,H^{\ZZ_2}(LV;W)] \to
\{\Sigma S(LV)^+\wedge X;LV^{\infty}\}_{\ZZ_2}~=~\{X;P(V)^+\}~;\\[1ex]
\hskip100pt h \mapsto \delta(F_h,F_0)~.
\end{array}$$
\indent The homeomorphism
$\lambda_{V,V}:S(V)*S(V) \cong S(V \oplus V)$ is $\ZZ_2$-equivariant
with respect to the transposition involutions
$$\begin{array}{l}
T~:~S(V)*S(V) \to S(V)*S(V)~;~(v,t,w) \mapsto (w,1-t,v)~,\\[1ex]
T~:~S(V \oplus V) \to S(V \oplus V)~;~(x,y) \mapsto (y,x)
\end{array}$$
with $\ZZ_2$-fixed point sets
$$\begin{array}{l}
(S(V)*S(V))^{\ZZ_2}~=~\{(v,1/2,v)\,\vert\, v \in S(V)~,\\[1ex]
S(V\oplus V)^{\ZZ_2}~=~\{(\dfrac{x}{\sqrt{2}},\dfrac{x}{\sqrt{2}})\,\vert\, x \in S(V)\}~.
\end{array}$$
The $\ZZ_2$-equivariant linear isometry
$$\kappa_V/\sqrt{2}~:~LV\oplus V \iso  V\oplus V ~; ~(v,w)\mapsto
(\dfrac{v+w}{\sqrt{2}},\dfrac{-v+w}{\sqrt{2}})$$
restricts to a $\ZZ_2$-equivariant homeomorphism
$$\kappa_V/\sqrt{2}~:~S(LV \oplus V) \to S(V \oplus V)~.$$
The composite
$$\begin{array}{l}
\mu_V~=~\lambda_{V,V}(\kappa_V/\sqrt{2})\lambda_{LV,V}^{-1}~:\\[1ex]
\hskip50pt S(LV)*S(V) \to S(LV \oplus V) \to S(V \oplus V) \to S(V)*S(V)
\end{array}$$
is a $\ZZ_2$-equivariant homeomorphism which is the identity on the
$\ZZ_2$-fixed point sets.
For any homotopy equivalence $h:S(V) \to S(V)$ the
diagram of $\ZZ_2$-equivariant maps
$$\xymatrix@C+35pt{
S(LV)*S(V) \ar[r]^-{\di{\mu_V}}
\ar[d]_-{\di{1 *h}} & S(V) * S(V) \ar[d]^-{\di{h *h}}\\
S(LV)*S(V) \ar[r]^-{\di{\mu_V}} & S(V) * S(V) }
$$
does not commute, but does restrict to a commutative diagram on the
$\ZZ_2$-fixed point subset $S(V) \subset S(LV)*S(V)$.
Choose a homotopy inverse  $h^{-1}:S(V) \to S(V)$ of $h$, and
define a $\ZZ_2$-equivariant homotopy equivalence
$$(1*h^{-1})\mu_V^{-1}(h*h)\mu_V~:~S(LV \oplus V) \to S(LV \oplus V)$$
which on the $\ZZ_2$-fixed point sets is homotopic to $1:S(V) \to S(V)$.

\begin{definition}~
{\rm (Crabb \cite[p.20]{crabb})\\
(i) The {\it doubling function} is defined by
$$S^2~:~[X,H(V)] \to [X,H^{\ZZ_2}(LV;V)]~;~h \mapsto \mu_V^{-1}(h*h)\mu_V~.$$
(ii) The {\it reduced doubling function}
$$\overline{S}^2~:~[X,H(V)] \to [X,H^{\ZZ_2}(LV;V)]~;~h \mapsto
(1*h^{-1})\mu_V^{-1}(h*h)\mu_V$$
is defined using a continuous choice of homotopy inverse
$h^{-1}:X \to H(V)$, with
$$\begin{array}{l}
\overline{S}^2(h)~ =~ (1*h^{-1})\mu_V^{-1}(h*h)\mu_V\\[1ex]
\hspace*{25pt}\in [X,H^{\ZZ_2}(LV;V)]~=~
{\rm ker}(\rho:[X,H^{\ZZ_2}(LV \oplus V)] \to [X,H(V)])~.
\end{array}$$
For $g:X \to O(V)$ and $h=Jg:X \to H(V)$ the actual inverse can be chosen
in (ii), and
$$\overline{S}^2(Jg)~=~g \oplus 1 \in {\rm im}(J:[X,O(V)] \to
[X,H^{\ZZ_2}(LV;V)])~.$$
(iii) Define the function
$$\eta~=~\zeta^{\ZZ_2}\overline{S}^2~:~[X,H(V)] \to [X,H^{\ZZ_2}(LV;V)]\to
\{X;P(V)^+\}~.$$
\hfill\qed}
\end{definition}

\begin{proposition}~ {\rm (Crabb \cite[pp.23,27]{crabb})}
Let $X$ be a pointed space.\\
{\rm (i)} There is defined a commutative diagram
$$\xymatrix@C+5pt@R+5pt{[X,O(V)] \ar[d]^-{\di{J}}
\ar[dr]^-{\di{\theta}} \ar[r]^-{\di{J_{\ZZ_2}F}}  &
[X,H^{\ZZ_2}(LV)]\ar[d]^-{\di{\zeta^{\ZZ_2}}} \\
[X,H(V)] \ar[r]^-{\di{\eta}} & \{X;P(V)^+\}}$$
with
$$J_{\ZZ_2}F~:~[X,O(V)] \to [X,H^{\ZZ_2}(LV)]~;~g \mapsto
g\vert_{S(LV)}$$
the forgetful map.\\
{\rm (ii)} Passing to the direct limit over all finite-dimensional
inner product spaces $V$ there is defined a commutative square
$$\xymatrix{\widetilde{KO}^{-1}(X)=\varinjlim\limits_V [X,O(V)]
\ar[d]_-{\di{\theta}}  \ar@<1ex>[r]^-{\di{J_{\ZZ_2}F}} &
\varinjlim\limits_V [X,H^{\ZZ_2}(LV\oplus V)]\ar[d] \\
\{X;P(\infty)^+\}=\varinjlim\limits_V \{X;P(V)^+\} \ar@<1ex>[r] &
\widetilde{\omega}_{\ZZ_2}^0(X)}$$
with
$$\begin{array}{l}
J_{\ZZ_2}F~:~\widetilde{KO}^{-1}(X)=\varinjlim\limits_V [X,O(V)] \to
\varinjlim\limits_V [X,H^{\ZZ_2}(LV\oplus V)]~;\\[2ex]
\hspace*{175pt} g \mapsto  (g \oplus 1)\vert_{S(LV \oplus V)}~,\\[1ex]
\varinjlim\limits_V [X,H^{\ZZ_2}(LV\oplus V)] \to
\widetilde{\omega}_{\ZZ_2}^0(X)=\varinjlim\limits_V \{(LV \oplus V)^{\infty} \wedge X;
(LV \oplus V)^{\infty}\}_{\ZZ_2}~;\\[1ex]
\hspace*{175pt} h \mapsto F_h-F_0~(F_0=1)
\end{array}$$
and $\{X;P(\infty)^+\} \to \widetilde{\omega}^0_{\ZZ_2}(X)$
the injection in the direct sum system
$$\begin{array}{l}
\{X;P(\infty)^+\}=\{X;S(\infty)^+\}_{\ZZ_2}
\xymatrix{\ar@<1ex>[r]^-{\di{\gamma}} \ar@<-1ex>@{<-}[r]_-{\di{\delta}}
& \{X;S^0\}_{\ZZ_2}=\widetilde{\omega}_{\ZZ_2}^0(X) }\\
\hspace*{150pt}
\xymatrix{\ar@<1ex>[r]^-{\di{\rho}} \ar@<-1ex>@{<-}[r]_-{\di{\sigma}}
& \{X;S^0\}=\widetilde{\omega}^0(X)~.}
\end{array}$$
\end{proposition}
\hfill\qed

\chapter{The double point theorem}\label{doublepoint}

We apply the geometric Hopf invariant to a homotopy theoretic treatment
of the double points of maps.

The connection between the Hopf invariant and double points has been much studied already, cf.  Boardman and Steer \cite{bs}, Dax \cite{dax}, Eccles \cite{eccles},
Haefliger \cite{haefliger}, Haefliger and Steer \cite{hs}, Hatcher and Quinn \cite{hatcherquinn}, Koschorke and Sanderson \cite{ks1},\cite{ks2}, Vogel \cite{vogel}, Wood \cite{wood} $\dots$~.

\section{Framed manifolds}

Every $m$-dimensional manifold $M$ admits an embedding
$M \subset V\oplus \R^m$ for some inner product space $V$
(and certainly if ${\rm dim}(V)>m$) with a normal vector
$V$-bundle $\nu_M:M \to BO(V)$. The embedding $M \subset V\oplus \R^m$
extends to an open embedding $E(\nu_M) \subset V \oplus \R^m$ with
compactification Umkehr map
$$(V\oplus \R^m)^{\infty}~=~V^{\infty} \wedge S^m \to T(\nu_M)~.$$

\begin{definition} {\rm
(i) A {\it framed $m$-dimensional manifold} $(M,b)$ is an $m$-dimensional
manifold $M$ together with an embedding $M \subset V \oplus \R^m$
and with an isomorphism $b:\nu_M\cong \epsilon_V$, or equivalently
with an extension to an open embedding $b:V \times M \subset V \oplus \R^m$.
\index{framed manifold, $(M,b)$}\\
(ii) The {\it Pontrjagin-Thom map} of  $(M,b)$ is the composite
\index{framed!manifold}\index{Pontrjagin-Thom!map, $\alpha_M$}
$$\alpha_M~:~(V\oplus \R^m)^{\infty}~=~V^{\infty}\wedge S^m \to
T(\nu_M)\xymatrix{\ar[r]^-{\di{T(b)}}_{\di{\cong}}&}
T(\epsilon_V)~=~ V^{\infty}\wedge M^+~,$$
representing an element
 $$\alpha_M \in \omega_m(M)~=~\varinjlim_V\,
 [V^{\infty} \wedge S^m,V^{\infty} \wedge M^+]~.$$
(iii) For any space $X$ let $\Omega^{fr}_m(X)$ be the
{\it framed bordism} group of pairs $(M,b,f)$
with $(M,b)$ a $m$-dimensional framed manifold and $f:M \to X$ a map.
\index{framed!bordism, $\Omega^{fr}_m(X)$}\\
 \hfill\qed
}
\end{definition}

As is well-known:

\begin{proposition}~
For any space $X$ the Pontrjagin-Thom construction defines an isomorphism
$$\Omega^{fr}_m(X) \xymatrix{\ar[r]^-{\di{\cong}}&} \omega_m(X)~;~
(M,b,f) \mapsto (1\wedge f)\alpha_M~.$$
In particular,
$$(M,1)~=~\alpha_M \in \Omega^{fr}_m(M)~=~\omega_m(M)~.$$
\hfill\qed
\end{proposition}

Let $\eta(\R^m)$ be the canonical $\R^m$-bundle over the
Stiefel space $V_{m+k,k}=O(\R^k,\R^{m+k})$ of linear isometries
$u:\R^k \to \R^{m+k}$, with
$$\begin{array}{l}
E(\eta(\R^m))~=~\{(u,x)\,\vert\,u \in O(\R^k,\R^{m+k}),x \in u(\R^k)^{\perp}\}~,\\[1ex]
\delta\eta(\R^m)~:~\eta(\R^m)\oplus\epsilon_{\R^k}~\cong~ \epsilon_{\R^{m+k}}~.
\end{array}$$

\begin{definition}~{\rm (Kervaire \cite{kerv1})\\
Let $(M,b)$ be a framed $m$-dimensional manifold, with an embedding $e:M \subset \R^{m+k}$ and an isomorphism $b:\nu_M \cong \epsilon_{\R^k}$. \\
{\rm (i)}  The {\it generalized Gauss map} of $(M,b)$\index{generalized Gauss map}
$$c~:~M \to V_{m+k,k}~;~x \mapsto (e_xb_x^{-1}:\R^k \cong
(\nu_M)_x \emb (\tau_{\R^{m+k}})_{e(x)}=\R^{m+k})$$
classifies the tangent bundle of $M$
$$c^*\eta(\R^m)~=~\tau_M~:~M \to BO(\R^m)$$
and the stable trivialization
$$\delta\tau_M~:~\tau_M \oplus \epsilon_{\R^k} \xymatrix{\ar[r]^-{\di{1 \oplus b^{-1}}}_-{\di{\cong}}&}
\tau_M \oplus \nu_M~=~\tau_{\R^{m+k}}\vert_M~=~\epsilon_{\R^{m+k}}~.$$
(ii) The {\it curvatura integra} of $(M,b)$ is the image of
the fundamental class $[M] \in H_m(M)$ under the generalized Gauss map
\index{curvatura integra}
$$c_*[M] \in H_m(V_{m+k,k})~=~Q_{(-)^m}(\ZZ)~=~\begin{cases}
\ZZ&{\rm if}~m \equiv 0(\bmod\,2)\\
\ZZ_2&{\rm if}~m \equiv 1(\bmod\,2)
\end{cases}$$
(assuming $k>1$).\\
(iii) The {\it Hopf invariant} of $(M,b)$ is \index{Hopf invariant!framed manifold}
$${\rm Hopf}(M,b)~=~\begin{cases}
0 \in \ZZ&{\rm if}~m \equiv 0(\bmod\,2)\\
{\rm Hopf}(p\alpha_M)\in \ZZ_2&{\rm if}~m \equiv 1(\bmod\,2)
\end{cases}$$
with $\alpha_M:S^{m+k} \to \Sigma^kM^+$ the Pontrjagin-Thom map
of $b:\R^k \times M \subset \R^{m+k}$, and $p:M \to \{*\}$ the projection.
Note that ${\rm Hopf}(M,b)=0$ for $m \neq 1,3,7$.\\
(iv) The {\it semicharacteristic} of an $m$-dimensional manifold $M$ is defined for odd $m$ by
\index{semicharacteristic, $\chi_{1/2}(M)$}
$$\chi_{1/2}(M)~=~ \sum\limits^{(m-1)/2}_{i=0}{\rm dim}_{\ZZ_2}H_i(M;\ZZ_2) \in \ZZ_2~.$$
\hfill\qed}
\end{definition}

By Proposition \ref{localR} (ii) the local obstruction function
$\theta$ (\ref{local}) defines an isomorphism
$$\theta~:~H_m(V_{m+k,k}) \xymatrix{\ar[r]^-{\di{\cong}}&}
H_m(P(\R^{m+k})/P(\R^m))~=~Q^{[0,k-1]}_{2m}(S^m\ZZ)$$
so that the {\it curvatura integra} is the image of $[M] \in H_m(M)$ under
the composite
$$\xymatrix{\theta(c)~:~M\ar[r]^-{\di{c}} & V_{m+k,k}
\ar[r]^-{\di{\theta}}&P(\R^{m+k})/P(\R^m)~=~S(L\R^k)^+\wedge_{\ZZ_2}(L\R^m)^{\infty}}~,$$
that is
$$\begin{array}{ll}
c_*[M]~=~\theta(c)_*[M] \in H_m(V_{m+k,k})
&=~Q^{[0,k-1]}_{2m}(S^m\ZZ)\\[1ex]
&=~\begin{cases} \ZZ
&{\rm if}~m~ {\rm is~even~or~if}~k=1\\
\ZZ_2&{\rm if}~m~ {\rm is~odd~and}~k \geqslant 2~.\end{cases}
\end{array}$$
The zero section maps
$$z_{\tau_M}~:~M \to T(\tau_M)~,~z_{\eta(\R^m)}~:~V_{m+k,k} \to
T(\eta(\R^m))$$
are such that (for any $m$) there is defined a commutative diagram
$$\xymatrix{H_m(M)\ar[d]^-{\di{z_{\tau_M}}} \ar[r]^-{\di{c_*}} &
H_m(V_{m+k,k}) \ar[d]^-{\di{z_{\eta(\R^m)}}}\\
\dot H_m(T(\tau_M))=\ZZ\ar[r]^-{\di{c_*}}_-{\cong} &
\dot H_m(T(\eta(\R^m)))=\ZZ}$$
with
$$z_{\tau_M}([M])~=~\chi(M)~,~z_{\eta(\R^m)}([S^m])~=~\chi(S^m)~=~1+(-)^m \in \ZZ~.$$
Thus for even $m$
$$c_*[M]~=~\chi(M)/2 \in H_m(V_{m+k,k})~=~\ZZ~,$$
as originally proved by Hopf \cite{hopf1}.  Kervaire \cite{kerv1,kerv2}
expressed the {\it curvatura integra} for odd $m$ in terms of the
semicharacteristic and the Hopf invariant:
$$c_*[M]~=~\chi_{1/2}(M)-{\rm Hopf}(M) \in H_m(V_{m+k,k})~.$$
We shall reprove this in Proposition \ref{curvint} below,
following the outline of Crabb \cite[Thm. 8.4]{crabb} and using
the quadratic construction on
$$\alpha_M~=~(M,1) \in \omega_m(M)~=~\Omega^{fr}_m(M)$$
which is a stable homotopy theoretic version of the quadratic
refinement $\mu$ of the intersection form $\lambda$ on $H^*(M)$.

\begin{proposition}~\label{mu}
{\rm (Pontrjagin \cite{pont} for $n=1$,
Kervaire and Milnor \cite{kervmiln} for $(n-1)$-connected $M^{2n}$,
Browder \cite{browder1}, \cite[\S III.4]{browder2} in general.)}\\
Let $(M,b)$ be an $m$-dimensional framed manifold,
with an embedding $\R^k \times M \subset \R^{m+k}$ and Pontrjagin-Thom
map $\alpha_M:S^{m+k}\to \Sigma^kM^+$.\\
{\rm (i)} Suppose that $m=2n$, and define the ring
$$A_n~=~\begin{cases} \ZZ&{\it if}~n \equiv 0(\bmod\,2)\\
\ZZ_2&{\it if}~n \equiv 1(\bmod\,2)
\end{cases}$$
so that there is a $(-)^n$-symmetric intersection pairing over $A_n$
$$\lambda~:~H^n(M;A_n) \times H^n(M;A_n) \to A_n~;~(x,y) \mapsto
\langle x \cup y,[M] \rangle$$
with
$$\lambda(x,y)~=~(-)^n\lambda(y,x) \in A_n~.$$
The framing of $M$ determines a $(-)^n$-quadratic refinement of $\lambda$
$$\begin{array}{l}
\mu~:~H^n(M;A_n) \to Q_{(-)^n}(A_n)~=~\ZZ/\{1-(-)^n\}~;\\[1ex]
\hskip100pt
x \mapsto \begin{cases} \lambda(x,x)/2&{\it if}~n \equiv 0(\bmod\,2)\\
\langle Sq^{n+1}_{xF}(\iota),[S^{2n}]\rangle&{\it if}~n \equiv 1(\bmod\,2)~,
\end{cases}
\end{array}$$
with $ \iota \in H^n(K(\ZZ_2,n);\ZZ_2)=\ZZ_2$ the generator,
$xF \in \pi^S_{2n}(K(\ZZ_2,n))$,
such that for all $x,y \in H^n(M;A_n)$, $a \in A_n$
$$\begin{array}{c}
\lambda(x,x)~=~\mu(x)+(-)^n\mu(x) \in A_n~,\\[1ex]
\mu(ax)~=~a^2\mu(x)~,~\lambda(x,y)~=~\mu(x+y)-\mu(x)-\mu(y) \in Q_{(-)^n}(A_n)~.
\end{array}$$
{\rm (ii)} Let $m=2n$ as in {\rm (i)}, and suppose that $x\in H^n(M;A_n)$
is the Poincar\'e dual of the homology class $x[S^n] \in H_n(M;A_n)$
represented by an embedding $x:S^n \subset M$ {\rm (}as is the case for all
$x \in H^n(M;A_n)$ if $M$ is $(n-1)$-connected with $n \geqslant 3${\rm )},
with normal bundle $\nu_x:S^n \to BO(\R^n)$ and  geometric Umkehr map $G:M^+ \to T(\nu_x)$.
The normal bundle of the composite $S^n \subset M^{2n} \subset \R^{2n+k}$ has a
canonical trivialization {\rm (}for large $k${\rm )}
$$\delta\nu_x~:~\nu_{S^n \subset \R^{2n+k}}~=~\nu_x \oplus \epsilon_{\R^k}~
\cong~\epsilon_{\R^{n+k}}~,$$
such that the corresponding Pontrjagin-Thom map
$$S^{2n+k} \xymatrix{\ar[r]^-{\di{F}}&} \Sigma^kM^+
\xymatrix{\ar[r]^-{\di{\Sigma^kG}}&} \Sigma^kT(\nu_x)
\xymatrix{\ar[r]^-{\di{\delta\nu_x}}_-{\di{\cong}}&}
\Sigma^k T(\epsilon^n) \to S^{n+k}$$
is null-homotopic. Then:\\
{\rm (a)} The evaluation of the quadratic form $\mu$ on $x \in H^n(M;A_n)$ is
$$\mu(x)~=~(\delta\nu_x,\nu_x)\in \pi_n(V_{n+k,k})~=~Q_{(-)^n}(A_n)~.$$
{\rm (b)} The exact sequence
$$\dots \to \pi_n(O(\R^n))
\to \pi_n(O) \to Q_{(-)^n}(A_n) \to \pi_n(BO(\R^n)) \to \pi_n(BO) \to \dots$$
is such that
$$Q_{(-)^n}(A_n)~=~\ZZ/\{1+(-)^{n+1}\}\to \pi_n(BO(\R^n))~;~1 \mapsto \tau_{S^n}$$
is injective for $n \notin \{1,3,7\}$ and is $0$ for $n \in \{1,3,7\}$. Thus if $n \notin \{1,3,7\}$
$$\mu(x)~=~\begin{cases} \chi(\nu_x)/2 \in \ZZ={\rm ker}(\pi_n(BO(\R^n)) \to \pi_n(BO))
&{\it if}~n \equiv 0(\bmod\,2)\\
\nu_x \in \ZZ_2={\rm ker}(\pi_n(BO(\R^n)) \to \pi_n(BO))&{\it if}~n \equiv 1(\bmod\,2)
\end{cases}$$
with $\mu(x)=0$ if and only if $\nu_x \cong \epsilon_{\R^n}$.\\
{\rm (c)} If $\nu_x=0 \in \pi_n(BO(\R^n))$ {\rm (}e.g. if $n \in \{1,3,7\}$
and $\nu_x$ is orientable{\rm )} then for any choice of trivialization
$\delta\nu'_x:\nu_x \cong \epsilon_{\R^n}$ the normal bundle of
$S^n \subset M^{2n} \subset \R^{2n+k}$ has a trivialization
$$\delta\nu'_x\oplus \epsilon_{\R^k}~:~\nu_{S^n \subset \R^{2n+k}}~=~
\nu_x \oplus \epsilon_{\R^k} ~\cong~\epsilon_{\R^{n+k}}~,$$
with the Pontrjagin-Thom map
$$S^{2n+k} \xymatrix{\ar[r]^-{\di{F}}&} \Sigma^kM^+
\xymatrix{\ar[r]^-{\di{\Sigma^kG}}&} \Sigma^kT(\nu_x)
\xymatrix{\ar[r]^-{\di{\Sigma^k\delta\nu'_x}}_-{\di{\cong}}&}
\Sigma^k T(\epsilon^n) \to S^{n+k}$$
such that
$$\mu(x)~=~{\rm Hopf}(S^n,\delta\nu'_x \oplus \epsilon_{\R^k}) \in
{\rm im}(\pi_n(BO) \to Q_{(-)^n}(A_n))~.$$
{\rm (d)} If $n \neq 2$ there exists a choice of trivialization
$\delta\nu'_x:\nu_x \cong \epsilon_{\R^n}$ such that $x[S^n] \in H_n(M;A_n)$
can be killed by a framed surgery on $M$ if and only if $\mu(x)=0$.\\
{\rm (iii)} For any $m$, the framing determines
a bundle map $b:\nu_M \to \epsilon_{\R^k}$
over a degree 1 map $f:M \to S^m$,  i.e. a normal map $(f,b):M \to S^m$.
The surgery obstruction of $(f,b)$ is
$$\begin{array}{l}
\sigma_*(f,b)~=~\begin{cases}
{\rm signature}(H^n(M),\lambda)/8\\[1ex]
{\rm Arf~invariant}(H^n(M;\ZZ_2),\lambda,\mu)\\[1ex]
0
\end{cases}\\[1ex]
\hskip50pt \in L_m(\ZZ)~=~
\begin{cases}
\ZZ&{\it if}~m=2n,~n \equiv 0(\bmod\,2)\\[1ex]
\ZZ_2&{\it if}~m=2n,~n \equiv 1(\bmod\,2)\\[1ex]
0&{\it if}~m \equiv 1(\bmod\,2)~.
\end{cases}
\end{array}$$
\end{proposition}
\hfill\qed

\begin{remark} {\rm
In Ranicki \cite{ranicki2} the quadratic function $\mu:H^n(M) \to Q_{(-)^n}(\ZZ)$ of a
$2n$-dimensional framed manifold $M$ was shown to be a special case of a
quadratic structure on the chain complex $C(M)$ of an
$m$-dimensional  framed manifold $M$ (for any $m$) obtained by the quadratic
construction on the Pontrjagin-Thom map $\alpha_M:V^{\infty} \wedge S^m
\to V^{\infty} \wedge M^+$ of an embedding $V \times M \subset V\oplus \R^m$,
namely
$$\psi~=~\psi(\alpha_M) \in H_m(S(\infty) \times_{\ZZ_2}(M \times M))~=~Q_m(C(M))~.$$
For any oriented $m$-dimensional manifold $M$ there is defined a degree 1 map
$$f~:~M \to M/(M \backslash D^m)~=~D^m/S^{m-1}~=~S^m$$
such that $f:C(M) \to C(S^m)$ is a chain homotopy surjection split by the
Umkehr chain map
$$f^!~:~C(S^m)~\simeq~C(S^m)^{m-*} \xymatrix{\ar[r]^-{\di{f^*}}&}
C(M)^{m-*}~\simeq~C(M)~,$$
so that there is defined a chain homotopy direct sum system
$$C(S^m) \xymatrix{\ar@<1ex>[r]^-{\di{f^!}}
\ar@<-1ex>@{<-}[r]_-{\di{f}}&} C(M)
\xymatrix{\ar@<1ex>[r]^-{\di{e}}
\ar@<-1ex>@{<-}[r]&} \Cc(f^!)$$
with $e:C(M) \to \Cc(f^!)$ the inclusion in the algebraic mapping cone.
A framing of $M$ determines a bundle map $b:\nu_M \to \epsilon_V$ over $f$,
i.e. a normal map $(f,b):M \to S^m$.
The $S$-dual of $T(b):T(\nu_M) \to T(\epsilon_V)$
is a geometric Umkehr map $F:\Sigma^{\infty}(S^m)^+ \to \Sigma^{\infty}M^+$
with $(\Sigma^{\infty}f)F \simeq 1$, such that
$$\alpha_M~:~\Sigma^{\infty}(S^m) \to \Sigma^{\infty}(S^m)^+
\xymatrix{\ar[r]^-{\di{F}}&} \Sigma^{\infty}M^+~.$$
The geometric Umkehr $F$ induces the chain Umkehr $f^!$.
The quadratic structure
$$\psi~=~\psi(\alpha_M) \in H_m(S(\infty) \times_{\ZZ_2}(M \times M))~=~Q_m(C(M))~,$$
is such that $(\Cc(f^!),e_{\%}(\psi))$ is an $m$-dimensional quadratic
Poincar\'e complex over $\ZZ$ representing the simply-connected surgery
obstruction
$$\sigma_*(f,b)~=~(\Cc(f^!),e_{\%}(\psi)) \in L_m(\ZZ)~.$$
See Chapter \ref{surgeryobstruction} below for more details.
\hfill\qed}
\end{remark}

For any space $M$ Proposition \ref{stablesequence1} (ii) gives a direct sum system
$$\omega_m(S(\infty) \times_{\ZZ_2}(M \times M))
\xymatrix{\ar@<1ex>[r]^-{\di{\gamma}}
\ar@<-1ex>@{<-}[r]_-{\di{\delta}}&}\omega_m^{\ZZ_2}(M \times M)
\xymatrix{\ar@<1ex>[r]^-{\di{\rho}}
\ar@<-1ex>@{<-}[r]_-{\di{\sigma}}&}\omega_m(M)~.$$
The quadratic construction function (\ref{quad1})
$$\psi~:~\omega_m(M) \to \omega_m(S(\infty) \times_{\ZZ_2}(M \times M))~;~
(\alpha:\Sigma^{\infty} \wedge S^m \to \Sigma^{\infty}M^+) \mapsto \psi(\alpha)$$
is such that
$$\begin{array}{l}
\gamma \psi(\alpha)~=~(\alpha \wedge \alpha)\Delta_{S^m} - \Delta_M(\alpha) \in
\omega^{\ZZ_2}_m(M \times M)~,\\[1ex]
\psi(\alpha+\beta)~=~\psi(\alpha)+\psi(\beta)+(\alpha \wedge \beta)\Delta_{S^m}
\in \omega_m(S(\infty) \times_{\ZZ_2}(M \times M))
\end{array}$$
with $(\alpha \wedge \beta)\Delta_{S^m}$ defined using the map induced by $S(L\R) \emb S(\infty)$
$$\omega_m(M\times M)~=~\omega_m(S(L\R)\times_{\ZZ_2}(M\times M))
 \to \omega_m(S(\infty)\times_{\ZZ_2}(M \times M))~.$$
 Let $p:M \to \{*\}$ be the projection, and let
$$q~=~(1 \times p \times p)~:~S(\infty)\times_{\ZZ_2}(M \times M) \to S(\infty)\times_{\ZZ_2}(
\{*\} \times \{*\})~=~P(\infty)~,$$
so that there is defined a commutative diagram
$$\xymatrix{
\omega_m(M) \ar[r]^-{\di{p_*}} \ar[d]^-{\di{\psi}} &
\omega_m \ar[d]^-{\di{\overline{P}^2}}\\
\omega_m(S(\infty) \times_{\ZZ_2}(M \times M))
\ar[r]^-{\di{q_*}} \ar[d]^-{\di{H}} &
\omega_m(P(\infty)) \ar[d]^-{\di{H}}\\
Q_m(C(M)) \ar[r]^-{\di{q_*}} & Q_m(\ZZ)=H_m(P(\infty))}$$
with $H$ the Hurewicz maps, and $\overline{P}^2$
the reduced squaring operation (\ref{square} (ii)).

\begin{proposition}~{\rm (Crabb \cite[p.45]{crabb})}\\
{\rm (i)} The quadratic construction on the Pontrjagin-Thom
map $\alpha_M:\Sigma^{\infty}S^m \to \Sigma^{\infty}M^+$ of an
$m$-dimensional framed manifold $(M,b)$ is a stable homotopy class
$$\psi(\alpha_M) \in \omega_m(S(\infty) \times_{\ZZ_2}(M \times M))$$
with Hurewicz image the quadratic construction
$\psi \in Q_m(C(M))$ of \cite{ranicki2}, defining a quadratic function
$$\begin{array}{l}
\psi~:~\Omega^{fr}_m(X)~=~\omega_m(X) \to \omega_m(S(\infty)\times_{\ZZ_2}(X \times X))~;\\[1ex]
\hskip100pt
(M,b,f:M \to X) \mapsto (1\times f \times f)_*\psi(\alpha_M)
\end{array}$$
for any space $X$.\\
{\rm (ii)} The bordism class $\alpha_M=(M,1)\in \omega_m(M)=\Omega^{fr}_m(M)$
has images
$$\begin{array}{l}
p\alpha_M~=~(M,b) \in \omega_m~=~\Omega^{fr}_m~,\\[1ex]
\psi(\alpha_M) \in \omega_m(S(\infty)\times_{\ZZ_2}(M \times M))~,\\[1ex]
H(\psi(\alpha_M)) \in Q_m(C(M))~,\\[1ex]
q_*H\psi(\alpha_M)~=~H\overline{P}^2(p\alpha_M)~=~{\rm Hopf}(M)\\[1ex]
\hskip50pt \in Q_m(\ZZ)~=~\begin{cases} \ZZ&{\it if}~m=0\\
0&{\it if}~m>0~{\it and}~m \equiv 0(\bmod\,2)\\
\ZZ_2&{\it if}~m \equiv 1(\bmod\,2)~.
\end{cases}
\end{array}$$
\end{proposition}
\begin{proof} By construction.\\
\hfill\qed\end{proof}

\begin{proposition} ~\label{semichar}
{\rm (Milnor and Stasheff \cite[\S11]{milnstas},
Crabb \cite[p.92]{crabb})}\\
Let $N$ be an $n$-dimensional manifold, so that
$N \times N$ is a $2n$-dimensional manifold, and $\Delta_N:N \emb N \times N$
is the embedding of an $n$-dimensional submanifold with normal bundle
$$\nu_{\Delta_N}~=~\tau_N~:~N \to BO(\R^n)~.$$
{\rm (i)} If $N$ is oriented then so is $N \times N$, and $\Delta_N$
represents a homology class $[N] \in H_n(N \times N)$ with a
Poincar\'e dual $[N]^* \in H^n(N \times N)$. The
$(-)^n$-symmetric intersection form $\lambda_{N \times N}$ on $H^n(N \times N)$ is such that
$$\lambda_{N \times N}([N]^*,[N]^*)~=~\chi(N) \in H_n(N)~=~\ZZ~,$$
with $\chi(N)=0$ for $n \equiv 1(\bmod\,2)$.\\
{\rm (ii)} If $N$ is framed then so is $N \times N$, and the $(-)^n$-quadratic
function $\mu_{N \times N}$ on $H^n(N \times N)$ is such that
$$\mu_{N \times N}([N]^*)~=~\begin{cases}
\chi(N)/2&{\it if}~n \equiv 0(\bmod\,2)\\
\chi_{1/2}(N)&{\it if}~n \equiv 1(\bmod\,2)
\end{cases} \in Q_{(-)^n}(\ZZ)$$
\end{proposition}
\begin{proof} (i) Let $F$ be a field, and work with $F$-coefficients.
By the K\"unneth formula
$$H^n(N \times N;F)~=~\sum\limits^n_{i=0}H^i(N;F)\otimes_FH^{n-i}(N;F)$$
with
$$\lambda_{N \times N}(\sum\limits_{j=1}^k x_j \otimes y_j,
\sum\limits_{j'=1}^{k'} x'_{j'} \otimes y'_{j'})~=~\sum\limits_{j=1}^k
\sum\limits_{j'=1}^{k'}\lambda_N(x_j,x'_{j'}) \lambda_N(y_j,y'_{j'}) \in F~.$$
Choose a basis $b_1,b_2,\dots,b_r$ for $H^*(N;F)$, and let
$b^*_1,b^*_2,\dots,b^*_r$ be the dual basis for $H^*(N;F)$ with
$$\lambda(b_p,b^*_q)~=~\begin{cases} 1&{\rm if}~p=q \\
0&{\rm if}~p\neq q~.
\end{cases}$$
See \cite[pp. 127--130]{milnstas} for the proof that
$$\begin{array}{l}
[N]^*~=~\sum\limits^r_{q=1}(-)^{{\rm dim}(b_q)}b_q \otimes b^*_q\\[1ex]
\hskip50pt
\in H^n(N \times N;F)~=~\sum\limits^n_{i=0}H^i(N;F)\otimes_FH^{n-i}(N;F)
\end{array}$$
with
$$\lambda_{N \times N}([N]^*,[N]^*)~=~\sum\limits^r_{q=1}(-)^{{\rm dim}(b_q)}
~=~\chi(N) \in {\rm im}(\ZZ \to F)~.$$
(ii) For $n \equiv 0(\bmod\,2)$ it is immediate from (i) and
$$\lambda_{N \times N}([N]^*,[N]^*)~=~2\mu_{N \times N}([N]^*) \in \ZZ$$
that $\mu_{N \times N}([N]^*)=\chi(N)/2$.\\
\indent So suppose that $n \equiv 1(\bmod\,2)$, and work with $\ZZ_2$-coefficients.
By the K\"unneth formula
$$H^n(N \times N;\ZZ_2)~=~\sum\limits^n_{i=0}H^i(N;\ZZ_2)\otimes_{\ZZ_2}H^{n-i}(N;\ZZ_2)$$
with
$$\begin{array}{l}
\lambda_{N \times N}(\sum\limits_{j=1}^k x_j \otimes y_j,
\sum\limits_{j'=1}^{k'} x'_{j'} \otimes y'_{j'})~=~\sum\limits_{j=1}^k
\sum\limits_{j'=1}^{k'}\lambda_N(x_j,x'_{j'}) \lambda_N(y_j,y'_{j'}) \in \ZZ_2~,\\[1ex]
\mu_{N \times N}(\sum\limits_{j=1}^k x_j \otimes y_j)~=~\sum\limits_
{1 \leqslant j < j' \leqslant k}\lambda_N(x_j,x_{j'}) \lambda_N(y_j,y_{j'}) \in
Q_{(-)^n}(\ZZ_2)~=~\ZZ_2~.
\end{array}$$
Continuing with the terminology of (i) (with $F=\ZZ_2$) we have
$$\mu_{N \times N}([N]^*)~=~\sum\limits_{{\rm dim}(b_q) <{\rm dim}(b^*_q)} 1
~=~\chi_{1/2}(N) \in \ZZ_2~.$$
\hfill\qed\end{proof}

\begin{proposition}\label{curvint}~
{\rm (Kervaire \cite{kerv1,kerv2}, Crabb \cite[Lemma 8.5]{crabb})} \\
{\rm (i)} Let $(M^{2n},b)$ be a $2n$-dimensional framed manifold,
so that there is given an embedding $M \subset \R^{2n+k}$ with a
trivialized normal $\R^k$-bundle $b:\nu_{M \subset \R^{2n+k}}\cong \epsilon_{\R^k}$.
Let $N^n \subset M^{2n}$ be an $n$-dimensional submanifold such that
the normal bundle of the composite embedding $N \subset  M \subset \R^{2n+k}$ is equipped with a trivialization
$a:\nu_{N \subset \R^{2n+k}}\cong \epsilon_{\R^{n+k}}$.
Let $c:N \to V_{n+k,k}$ classify the normal $\R^n$-bundle
$\nu_{N \subset M}:N \to BO(\R^n)$
with the corresponding stable trivialization
$$\delta\nu_{N \subset M}~:~
\nu_{N \subset M}\oplus \epsilon_{\R^k}~\cong~\epsilon_{\R^{n+k}}~.$$
The Hopf invariant of $(N,a)$ is given by
$${\rm Hopf}(N,a)~=~\mu([N]^*)-c_*[N] \in H_n(V_{n+k,k})~=~Q_{(-)^n}(\ZZ)$$
with $\mu([N]^*)$ the evaluation of the quadratic function $\mu$
on the Poincar\'e dual $[N]^* \in H^n(M)$ of $[N] \in H_n(M)$.
{\rm (}Note that ${\rm Hopf}(N,a)=0$ for $n \neq 1,3,7$,
or if $N=S^n$ with the canonical framing of the normal bundle of
$S^n \subset M^{2n} \subset \R^{2n+k}$.{\rm )}\\
{\rm (ii)} Let $(N,a)$ be an $n$-dimensional framed manifold,
so that there is given an embedding $N \subset \R^{n+k}$ with a
trivialized normal bundle. Then
$${\rm Hopf}(N,a)~=~\chi_{1/2}(N)-c_*[N] \in H_n(V_{n+k,k})~=~Q_{(-)^n}(\ZZ)$$
with $c:N \to V_{n+k,k}$ classifying the tangent $\R^n$-bundle
$\tau_N:N \to BO(\R^n)$ with the corresponding stable isomorphism
$$a~:~\tau_N\oplus \epsilon_{\R^k} \cong \epsilon_{\R^{n+k}}~.$$
{\rm (iii)} As in {\rm (i)} let $N^n \subset M^{2n}\subset \R^{2n+k}$ with
trivializations $b:\nu_{M \subset \R^{2n+k}}\cong \epsilon_{\R^k}$,
$a:\nu_{N \subset M^{2n}}\cong \epsilon_{\R^n}$, and let
$c:N \to V_{2n+k,n+k}$ classify $\tau_N:N \to BO(\R^n)$ with the corresponding
stable trivialization $\tau_N \oplus \epsilon_{\R^{n+k}}\cong \epsilon_{\R^{2n+k}}$.
The evaluation of $\mu:H^n(M) \to Q_{(-)^n}(\ZZ)$ on $[N]^* \in H^n(M)$ is
$$\begin{array}{l}
\mu([N]^*)~=~{\rm Hopf}(N,a)~=~\chi_{1/2}(N)-c_*[N]~(=0~{\rm for}~n \neq 1,3,7)\\[1ex]
\hskip150pt  \in H_n(V_{2n+k,n+k})~=~Q_{(-)^n}(\ZZ)~.
\end{array}$$
In particular, if $N=S^n$ then $\chi_{1/2}(N)=1$ and
$$\mu([S^n]^*)~=~{\rm Hopf}(S^n,\text{std.})~=~1-c_*[S^n] \in
H_n(V_{2n+k,n+k})~=~Q_{(-)^n}(\ZZ)~.$$
\end{proposition}
\begin{proof} (i) Let $V=\R^k$ and let
$$\alpha_M~:~V^{\infty}\wedge (\R^{2n})^{\infty} \to V^{\infty} \wedge M^+$$
be the Pontrjagin-Thom map of $V \times M \subset V \oplus \R^{2n}$,
so that $\psi(\alpha_M) \in \omega_{2n}(S(\infty)\times_{\ZZ_2}(M \times M))$
determines the $(-)^n$-quadratic form $\mu$ on $H^n(M)$
refining the $(-)^n$-symmetric intersection form $\lambda$. Let
$$\alpha_{N \subset M}~:~
M^+ \to M/(M \backslash E(\nu_{N\subset M}))~=~T(\nu_{N \subset M})$$
be the Umkehr map of $N \subset M$. Apply the composition formula for the quadratic construction
to the factorization
$$\xymatrix{ & V^{\infty} \wedge M^+ \ar[d]^-{\di{1\wedge \alpha_{N \subset M}}} \\
V^{\infty} \wedge S^{2n} \ar[ur]^-{\di{\alpha_M}} \ar[dr]^-{\di{\alpha_N}} &
V^{\infty} \wedge T(\nu_{N \subset M})
\ar[d]^-{\di{T(\delta\nu_{N \subset M})}} \\
& (V\oplus \R^n)^{\infty} \wedge N^\infty}$$
to obtain
$$\begin{array}{l}
\psi(\alpha_N)~=~(T(\delta\nu_{N \subset M})(1 \wedge \alpha_{N \subset M}))_{\%}\psi(\alpha_M)+
\psi(T(\delta\nu_{N \subset M}))(1\wedge \alpha_{N \subset M})\alpha_M\\[1ex]
\hskip150pt :~H_{2n}(S^{2n}) \to Q_{2n}(S^nC(M))~.
\end{array}$$
Let $p:N \to \{*\}$ be the projection. Now
$$\begin{array}{l}
p_{\%}\psi(\alpha_N)[S^{2n}]~=~{\rm Hopf}(N)~,\\[1ex]
p_{\%}(T(\delta\nu_{N \subset M})(1 \wedge \alpha_{N \subset M}))_{\%}\psi(\alpha_M)[S^{2n}]~=~\mu([N]^*)~,\\[1ex]
p_{\%}\psi(T(\delta\nu_{N \subset M}))(1\wedge \alpha_{N \subset M})\alpha_M[S^{2n}]~=~
-c_*[N] \in Q_{2n}(S^n\ZZ) ~~(\ref{twist})~,
\end{array}$$
so that
$${\rm Hopf}(N)~=~\mu([N]^*)-c_*[N] \in Q_{2n}(S^n\ZZ)~=~Q_{(-)^n}(\ZZ)~.$$
(ii) Apply (i) with $\Delta_N:N \subset M=N \times N$, $\nu_{\Delta_N}=\tau_N$,
using \ref{semichar} to identify $\mu([N]^*)=\chi_{1/2}(N)$.\\
(iii) Immediate from (i), (ii) and Proposition \ref{mu} (ii) (c).\\
\hfill\qed\end{proof}

\begin{example} {\rm The embeddings
$$x~:~S^n \times \{*\} \subset S^n \times S^n~,~y~:~\{*\}\times S^n \subset S^n \times S^n$$
have trivial normal $\R^n$-bundles $\nu_x=\nu_y=\epsilon^n$. For $n=1,3,7$
there exist framings $b$ of $S^n \times S^n \subset \R^{2n+k}$
such that the corresponding quadratic function
$$\mu_b~:~H^n(S^n \times S^n;\ZZ_2)~=~\ZZ_2 \oplus \ZZ_2 \to \ZZ_2$$
takes values
$$\mu_b(x)~=~\mu_b(y)~=~1 \in \ZZ_2$$
with the corresponding framings of $x,y:S^n \subset \R^{2n+k}$ such that
$${\rm Hopf}(S^n,x)~=~{\rm Hopf}(S^n,y)~=~1 \in \ZZ~.$$
\hfill\qed}
\end{example}

\begin{example} {\rm Given a map $\omega:S^m \to O(k)$ the embedding
$$e_{\omega}~:~S^m \times D^k \emb S^{m+k}~=~S^m \times D^k \cup
D^{m+1} \times S^{k-1}~;~(x,y) \mapsto (x,\omega(x)(y))$$
defines a submanifold $S^m \subset S^{m+k}$ with a framing $b_\omega$ such that
the generalized Gauss map is given by
$$c~:~S^m \to V_{m+k,k}~;~x \mapsto (y \mapsto (x,\omega(x)(y)))$$
and the Pontrjagin-Thom map is given by the $J$-homomorphism
$$J(\omega)~:~S^{m+k} \to S^k~;~
(x,y) \mapsto \begin{cases}\omega(x)(y)&{\rm if}~(x,y) \in S^m \times D^k\\
*&{\rm if}~(x,y) \in D^{m+1} \times S^{k-1}~.
\end{cases}$$
The {\it curvatura integra} is
$$\begin{array}{l}
c_*[S^m]~=~\begin{cases} \chi(S^m)/2~=~1 \\
\chi_{1/2}(S^m)-{\rm Hopf}(J(\omega))~=~1-{\rm Hopf}(J(\omega))
\end{cases}\\[1ex]
\hskip100pt \in H_m(V_{m+k,k})~=~\begin{cases}
\ZZ&{\rm if}~m~{\rm is~even}\\
\ZZ_2&{\rm if}~m~{\rm is~odd}
\end{cases}
\end{array}$$
with ${\rm Hopf}(J(\omega))=0$ for $m \neq 1,3,7$.\\
}\hfill\qed
\end{example}

An $m$-dimensional framed manifold $M$ determines a normal map
$(f,b):M \to S^m$ with a surgery obstruction $\sigma_*(f,b) \in L_m(\ZZ)$.
Following Crabb \cite{crabb} we shall now obtain a geometric
realization of the surgery obstruction by expressing the quadratic construction
$$\psi(\alpha_M) \in \omega_m(S(\infty)\times_{\ZZ_2}(M \times M))$$
(in Example \ref{Y} below) as a bordism class of an $m$-dimensional framed
manifold $N$ with a map
$g:N \to S(\infty)\times_{\ZZ_2}(M \times M)$ such that
$$\begin{array}{l}
(N,g)~=~(M,\Delta_M)-(S^m,(\alpha_M \wedge \alpha_M)\Delta_{S^m})~=~
-\psi(\alpha_M)\\[1ex]
\in {\rm ker}(\rho:\omega^{\ZZ_2}_m(M \times M) \to \omega_m(M))~=~
\omega_m(S(\infty)\times_{\ZZ_2}(M \times M))~,
\end{array}$$
and
$$(\overline{N},\overline{g})~=~(M,\Delta_M) \in \omega_m(M\times M)$$
with $\overline{N}=g^*(S(\infty)\times (M \times M))$ the canonical
double cover of $N$ and
$$\overline{g}~:~\overline{N} \to S(\infty)\times (M \times M) \to M\times M~.$$
The construction uses the $\ZZ_2$-equivariant $S$-duality analogue of:

\begin{proposition}~
An $m$-dimensional  framed  manifold $(M,b)$ is self-$S$-dual,
in the sense that the composite
$$\Delta \alpha_M~:~(V\oplus \R^m)^{\infty} \to V^{\infty} \wedge M^+
\to M^+ \wedge T(\epsilon_V)~=~
M^+ \wedge (V^{\infty} \wedge M^+)$$
is an $S$-duality map between $M^+$ and $V^{\infty} \wedge M^+$.
Cap product with the stable homotopy fundamental class
$$[M]~=~(M,1)~=~\alpha_M \in \omega_m(M)$$
defines Poincar\'e duality isomorphisms in stable homotopy
$$[M]\cap -~:~\omega^*(M)~\cong~\omega_{m-*}(M)~.$$
In particular $[M]\cap -:\omega^0(M)~\cong~\omega_m(M)$ sends
$1 \in \omega^0(M)$ to $[M] \in \omega_m(M).$\\
\hfill\qed
\end{proposition}

And now for the $\ZZ_2$-analogues.
Every $m$-dimensional $\ZZ_2$-manifold $M$ admits a $\ZZ_2$-equivariant
embedding $i:M \subset LV\oplus W\oplus R^m$
for some inner product spaces $V,W$ with a normal $\ZZ_2$-equivariant
vector $LV\oplus W$-bundle $\nu_M$, and $i$ extends to an open $\ZZ_2$-equivariant embedding
$i:E(\nu_M) \subset LV\oplus W \oplus  \R^m$.

\begin{definition} {\rm
(i) An $m$-dimensional $\ZZ_2$-manifold $M$ is {\it $\ZZ_2$-framed}
if there is given a $\ZZ_2$-embedding $i:M \subset LV \oplus W\oplus \R^m$
with an isomorphism $\nu_M\cong \epsilon_{LV\oplus W}$,
in which case $i$ extends to an open $\ZZ_2$-equivariant embedding
\index{$\ZZ_2$-framed!$\ZZ_2$-manifold}
$$i~:~LV\oplus W \times M \subset LV \oplus W \oplus \R^m~.$$
(ii) For any $\ZZ_2$-space let
$\Omega^{\ZZ_2\hbox{-}fr}_m(X)$ be the
{\it $\ZZ_2$-framed bordism} group of pairs $(M,f)$
with $M$ an $m$-dimensional $\ZZ_2$-framed $\ZZ_2$-manifold and $f:M \to X$ a
$\ZZ_2$-equivariant map.\index{$\ZZ_2$-framed!$\ZZ_2$-bordism} \\
(iii) The {\it $\ZZ_2$-equivariant Pontrjagin-Thom} map of a $\ZZ_2$-framed
$\ZZ_2$-manifold $M$ is
the compactification Umkehr map of $i:LV\oplus W \times M \subset
LV\oplus W \oplus \R^m$, a stable $\ZZ_2$-equivariant map
\index{$\ZZ_2$-equivariant Pontrjagin-Thom}
$$\begin{array}{l}
\alpha^{\ZZ_2}_M~:~(LV\oplus W\oplus \R^m)^{\infty}~=~
LV^{\infty}\wedge W^{\infty} \wedge S^m\\[1ex]
\hskip100pt
 \to T(\epsilon_{LV \oplus W})~=~LV^{\infty}\wedge W^{\infty}\wedge M^+
 \end{array} $$
 representing an element
 $$\alpha^{\ZZ_2}_M \in \omega^{\ZZ_2}_m(M)~=~
 \varinjlim_{V,W} [(LV\oplus W \oplus \R^m)^{\infty},
LV^{\infty}\wedge W^{\infty}\wedge M^+]_{\ZZ_2}$$
 \hfill\qed
}
\end{definition}

For any $\ZZ_2$-space $X$ the $\ZZ_2$-equivariant Pontrjagin-Thom
construction defines a morphism
$$
\Omega^{\ZZ_2\hbox{-}fr}_m(X) \to \omega^{\ZZ_2}_m(X)~;~
(M,f) \mapsto (1\wedge f)\alpha^{\ZZ_2}_M~.
$$
This morphism is surjective, but
it is NOT an isomorphism in general. We refer to
Appendix \ref{appendix2} for an account of $\ZZ_2$-equivariant bordism,
including the obstructions to the $\ZZ_2$-equivariant Pontrjagin-Thom map being an isomorphism.
The group $\Omega^{\ZZ_2\hbox{-}fr}_m(X)$ is written there, more systematically,
as $\Omega^{\ZZ_2}_m(X;\, 0,0)$.

\begin{proposition}
{\rm (i)} The $\ZZ_2$-equivariant Pontrjagin-Thom map
$\Omega^{\ZZ_2\hbox{-}fr}_m(X) \to \omega^{\ZZ_2}_m(X)$
sends $(M,1) \in \Omega^{\ZZ_2\hbox{-}fr}_m(M)$ to
 $\alpha^{\ZZ_2}_M \in \omega^{\ZZ_2}_m(M)$.\\
{\rm (ii)} A $\ZZ_2$-framed  $m$-dimensional $\ZZ_2$-manifold $M$ is
$\ZZ_2$-equivariantly self-$S$-dual, in the sense that the composite
$$\begin{array}{l}
\Delta \alpha^{\ZZ_2}_M~:~(LV\oplus W\oplus \R^m)^{\infty} \to
LV \oplus W^{\infty} \wedge M^+\\[1ex]
\hskip100pt
\to M^+ \wedge T(\epsilon_{LV\oplus W})~=~
M^+ \wedge (LV^{\infty} \wedge W^{\infty} \wedge M^+)
\end{array}$$
is a $\ZZ_2$-equivariant $S$-duality map. Cap product with the
$\ZZ_2$-equivariant stable homotopy fundamental class
$$[M]~=~(M,1)~=~\alpha^{\ZZ_2}_M\in \omega^{\ZZ_2}_m(M)$$
defines Poincar\'e duality isomorphisms in stable homotopy
$$[M]\cap -~:~\omega_{\ZZ_2}^*(M)~\cong~\omega^{\ZZ_2}_{m-*}(M)~.$$
In particular $[M]\cap -:\omega^0_{\ZZ_2}(M)~\cong~\omega^{\ZZ_2}_m(M)$ sends
$1 \in \omega^0_{\ZZ_2}(M)$ to $\alpha^{\ZZ_2}_M \in \omega^{\ZZ_2}_m(M).$\\
{\rm (iii)} For any space $X$ there is defined a direct sum system
$$\omega_m(S(\infty) \times_{\ZZ_2}(X \times X))
\xymatrix{\ar@<1ex>[r]^-{\di{\gamma}}
\ar@<-1ex>@{<-}[r]_-{\di{\delta}}&}\omega_m^{\ZZ_2}(X \times X)
\xymatrix{\ar@<1ex>[r]^-{\di{\rho}}
\ar@<-1ex>@{<-}[r]_-{\di{\sigma}}&}\omega_m(X)$$
with the various maps described in terms of representatives by
framed manifolds as
$$\begin{array}{l}
\rho~:~\omega_m^{\ZZ_2}(X \times X) \to \omega_m(X)~;\\[1ex]
\hskip100pt
(M,f:M \to X \times X) \mapsto (M^{\ZZ_2},f^{\ZZ_2}:M^{\ZZ_2} \to X)~,\\[1ex]
\gamma~:~\omega_m(S(\infty) \times_{\ZZ_2}(X \times X)) \to
\omega_m^{\ZZ_2}(X \times X)~;\\[1ex]
(N,g:N \to S(\infty) \times_{\ZZ_2}(X \times X))
\mapsto (\overline{N},\overline{g})~,~
\overline{N}=g^*(S(\infty)\times X\times X)~,\\[1ex]
\sigma~:~\omega_m(X) \to \omega_m^{\ZZ_2}(X \times X)~;\\[1ex]
\hskip100pt (P,h:P \to X) \mapsto (P,\Delta_Xh:P \to X \times X)~,\\[1ex]
\delta~:~\omega_m^{\ZZ_2}(X \times X)~=~\{S^m;X^+ \wedge X^+\}_{\ZZ_2}
\to \omega_m(S(\infty)\times_{\ZZ_2}(X \times X))~;\\[1ex]
(F:V^{\infty}\wedge LW^{\infty} \wedge S^m  \to
V^{\infty} \wedge LW^{\infty} \wedge X^+ \wedge X^+)
\mapsto \delta(F,\sigma \rho(F)) \, .
\end{array}$$
\end{proposition}
\begin{proof} (i)+(ii) By construction.\\
(iii) Immediate from Proposition \ref{stablesequence1} (ii).
\hfill\qed\end{proof}

For any framed $m$-dimensional manifold $(M,b)$ with embedding
$V \times M \subset V \oplus \R^m$ and Pontrjagin-Thom map
$\alpha_M:(V \oplus \R^m)^{\infty} \to (V \times M)^{\infty}$
the product $M \times M$ is a $2m$-dimensional $\ZZ_2$-manifold with a
$\ZZ_2$-equivariant embedding
$$\begin{array}{l}
(\kappa_{V \oplus \R^m})^{-1}(i \times i)(\kappa_V\times 1)~:~
(LV \oplus V) \times (M \times M)~\cong~V \times V\times M \times M\\[1ex]
\hskip75pt
\to (V \oplus \R^m) \times (V \oplus \R^m)~ \cong~(LV \oplus L\R^m)
\times (V \oplus \R^m)
\end{array}$$
and $\ZZ_2$-equivariant Pontrjagin-Thom map
$$\begin{array}{l}
(\kappa_V\wedge 1)(\alpha_M\wedge \alpha_M)(\kappa_{V \oplus \R^m})^{-1}~:~
(LV\oplus L\R^m)^{\infty} \wedge (V\oplus \R^m)^{\infty}\\[1ex]
\hskip175pt
\to LV^{\infty} \wedge V^{\infty} \wedge M^+\wedge M^+~.
\end{array}$$
The group $\omega_{2m,m}(M \times M)$ is defined in Example \ref{stablesequence2} (i). The image of
$$\begin{array}{ll}
[M\times M,b \times b]&=~(\kappa_V\wedge 1)(\alpha_M\wedge \alpha_M)(\kappa_{V \oplus \R^m})^{-1}\\[1ex]
&\in \omega_{2m,m}(M \times M)~=~\{(L\R^m\oplus \R^m)^{\infty},M^+\wedge M^+\}_{\ZZ_2}
\end{array}$$
under
$$0^*_{L\R^m}~:~\omega_{2m,m}(M \times M)\to
\omega^{\ZZ_2}_m(M \times M)~=~\omega_m(M) \oplus
\omega_m(S(\infty) \times_{\ZZ_2}(M \times M))$$
is of the form
$$\begin{array}{l}
0^*_{L\R^m}[M \times M,b \times b]~=~((M,b),(N,a)) \\[1ex]
\hskip50pt \in
\omega^{\ZZ_2}_m(M \times M)~=~
\omega_m(M)\oplus \omega_m(S(\infty) \times_{\ZZ_2}(M \times M))
\end{array}$$
with $(N,a)$ the framed $m$-dimensional manifold with a map
$g:N \to S(\infty) \times_{\ZZ_2}(M \times M)$
constructed from $M$ by Crabb \cite[p.47]{crabb}.
Note that there is defined an exact sequence
$$\begin{array}{l}
\to A_1=\omega_{2m,m}(S(L\R^m) \times M \times M) \to
A_2=\omega_{2m,m}(S(L\R^m \oplus L\R(\infty))\times M \times M) \\[1ex]
\to A_3=\omega_m(S(\infty)\times_{\ZZ_2} (M \times M))
\to A_4=\omega_{2m-1,m-1}(S(L\R^m) \times M \times M)
\end{array}$$
which fits into a commutative braid of exact sequences

$$\xymatrix@C-15pt{
A_1\ar[dr] \ar@/^2pc/[rr] && \omega_{2m,m}(M \times M)
\ar[dr]^-{\di{0^*_{L\R^m}}}\ar@/^2pc/[rr]^-{\di{\rho}} && \omega_m(M)\\
&A_2
\ar[ur] \ar[dr]&& \omega^{\ZZ_2}_m(M \times M) \ar[dr]
\ar@<1ex>[ur]^-{\di{\rho}} \ar@<-1ex>@{<-}[ur]_-{\di{\sigma}}&\\
\omega_{m+1}(M)\ar[ur] \ar@/_2pc/[rr]^-{\di{0}}
&&A_3
\ar@/_2pc/[rr]
\ar@<1ex>[ur]^-{\di{\gamma}} \ar@<-1ex>@{<-}[ur]_-{\di{\delta}}&&
A_4}$$

\begin{example} \label{Y} {\rm (Crabb \cite[pp. 44-48]{crabb})
If $(M,b)$ is a framed $m$-dimensional manifold with an open embedding
$i:V \times M \subset V \oplus \R^m$ the Pontrjagin-Thom map
$$\alpha_M~:~(V \oplus \R^m)^{\infty}~=~V^{\infty} \wedge S^m \to
(V\times M)^{\infty}~=~V^{\infty} \wedge M^+$$
defines the stable homotopy fundamental class
$$[M]~=~\alpha_M \in \omega_m(M)~=~\{S^m;M^+\}~.$$
The product $\ZZ_2$-manifold $M \times M$ has stable $\ZZ_2$-equivariant
homotopy fundamental class
$$[M\times M]~=~\alpha_M\wedge \alpha_M \in \omega^{\ZZ_2}_{2m,m}(M\times M)~=~
\{(L\R^m)^{\infty} \wedge S^m;M^+\wedge M^+\}_{\ZZ_2}~.$$
The quadratic construction on $\alpha_M$
$$\psi(\alpha_M) \in \omega_m(S(LV)\times_{\ZZ_2}(M \times M))$$
has the following geometric interpretation, with $n={\rm dim}(V)$. Define a
$(2m+n)$-dimensional manifold with boundary
$$(W,\partial W)~= (D(LV),S(LV))\times M \times M ~,$$
with $\ZZ_2$-action
$$T~:~W \to W~;~(v,x,y) \mapsto (-v,y,x)$$
(which is free on $\partial W$), with a stable $\ZZ_2$-equivariant homotopy
fundamental class
$$[W]_{\alpha_M \wedge \alpha_M} ~=~\alpha_M \wedge \alpha_M\in \omega^{\ZZ_2}_m(W,\partial W;
-\epsilon_{LV\oplus L\R^m})~=~\omega_{2m,m}(M \times M)~.$$
Consider the stable $\ZZ_2$-equivariant cohomotopy exact sequence
associated to the trivial $\ZZ_2$-bundle $\epsilon_{LV\oplus L\R^m}$ over $W$
$$\begin{array}{l}
\dots \to \omega^{-1}_{\ZZ_2}(\partial W;-\epsilon_{LV\oplus L\R^m})
\to \omega^0_{\ZZ_2}(W,\partial W;-\epsilon_{LV\oplus L\R^m})
\to \omega^0_{\ZZ_2}(W;-\epsilon_{LV\oplus L\R^m})\\[1ex]
\hspace*{150pt}
\to \omega^0_{\ZZ_2}(\partial W;-\epsilon_{LV\oplus L\R^m}) \to \dots~.
\end{array}$$
There are two reasons for the
stable $\ZZ_2$-equivariant cohomotopy Euler class
$$\gamma^{\ZZ_2}(\epsilon_{LV\oplus L\R^m}) \in
\omega^0_{\ZZ_2}(W;-\epsilon_{LV\oplus L\R^m})~=~\{(M\times M)^+;(LV\oplus L\R^m)^{\infty}\}_{\ZZ_2}$$
to restrict to 0 on $\partial W$, namely the $\ZZ_2$-equivariant sections
of the sphere bundle of $\epsilon_{LV\oplus L\R^m}$
$$s_1~,~s_2~:~\partial W \to S(\epsilon_{LV\oplus L\R^m})~=~
S(LV\oplus L\R^m)\times \partial W $$
defined by
$$s_1(v,x,y)~=~
\big(\dfrac{i(v,x)-i(-v,y)}{\Vert i(v,x)-i(-v,y) \Vert},(v,x,y)\big)~,
~s_2(v,x,y)~=~(v,(v,x,y))~.$$
The $\ZZ_2$-equivariant map
$$\begin{array}{l}
(s_1,-s_2)~:~\partial W \to S(\epsilon_{LV\oplus L\R^m})
\times_{\partial W}S(\epsilon_{LV\oplus L\R^m})\\[1ex]
\hspace*{100pt}=~
S(LV\oplus L\R^m)\times S(LV\oplus L\R^m) \times \partial W~;\\[1ex]
\hspace*{100pt} (v,x,y) \mapsto \big(s_1(v,x,y),-s_2(v,x,y),(v,x,y)\big)
\end{array}$$
is $\ZZ_2$-equivariant homotopic to a $\ZZ_2$-equivariant map
$$(t_1,-t_2)~:~\partial W \to S(\epsilon_{LV\oplus L\R^m})
\times_{\partial W}S(\epsilon_{LV\oplus L\R^m})$$
which is
$\ZZ_2$-equivariantly transverse regular at
$\Delta_{S(LV \oplus L\R^m)}\times\partial W$. The inverse image
$$Y~=~(t_1,-t_2)^{-1}(\Delta_{S(LV \oplus L\R^m)}\times\partial W)
\subset \partial W$$
is a framed $m$-dimensional manifold with a free $\ZZ_2$-action, and
with a $\ZZ_2$-equivariant map $Y \to \partial W$, such that $Y$ is
nonequivariantly framed cobordant to $M$. The $m$-dimensional manifold
$$N~=~Y/\ZZ_2$$
is equipped with a map $g:N \to \partial W/\ZZ_2=S(LV)\times_{\ZZ_2}(M \times M)$.
The $\ZZ_2$-equivariant
Poincar\'e duality defines an isomorphism of exact sequences
$$\xymatrix{
\omega^{-1}_{\ZZ_2}(\partial W;-\epsilon_{LV\oplus L\R^m})
\ar[r] \ar[d]^-{\di{\cong}}_{\di{[\partial W]_{\alpha_M \wedge \alpha_M}\cap -}} &
\omega^0_{\ZZ_2}(W,\partial W;-\epsilon_{V \oplus L\R^m})
\ar[r]
\ar[d]^-{\di{\cong}}_{\di{[W]_{\alpha_M \wedge \alpha_M}\cap -}}&
\omega^0_{\ZZ_2}(W;-\epsilon_{LV \oplus L\R^m})
\ar[d]^-{\di{\cong}}_{\di{[W]_{ \alpha_M \wedge \alpha_M}\cap -}}\\
\omega_m^{\ZZ_2}(\partial W)\ar[r] &
\omega_m^{\ZZ_2}(W)\ar[r]  & \omega_m^{\ZZ_2}(W,\partial W)}$$
and
$$\begin{array}{l}
\delta(s_1,s_2)~=~\psi_V(\alpha_M)~=~(N,g)\\[1ex]
\in \omega^{-1}_{\ZZ_2}(\partial W;-\epsilon_{LV \oplus L\R^m})
~=~\{\Sigma \partial W^+;(LV\oplus L\R^m)^{\infty}\}_{\ZZ_2}\\[1ex]
\hspace*{125pt} ~=~
\omega_m^{\ZZ_2}(\partial W)~=~\omega_m(S(LV)\times_{\ZZ_2}(M \times M)).
\end{array}$$
Note that for any $v \in S(LV)$ the $\ZZ_2$-equivariant map
$$v \times 1 \times 1~:~S(L\R)\times M \times M \to S(LV) \times M \times M$$
sends $\delta(s_1,s_2)$ to the element
$$(v \times 1 \times 1)^*\delta(s_1,s_2)
\in \omega^{-1}_{\ZZ_2}(S(L\R) \times M \times M;-\epsilon_{LV\oplus L\R^m})~=~
\{M \times M;S(V \oplus \R^m)\}$$
represented by the map
$$M \times M \to S(V \oplus \R^m)~=~S^{m+n-1}~;~
(x,y) \mapsto \dfrac{i(v,x)-i(-v,y)} {\Vert i(v,x)-i(-v,y) \Vert}~.$$
The rel $\partial W$ $\ZZ_2$-equivariant Euler classes are
$$\begin{array}{l}
\gamma^{\ZZ_2}(\epsilon_{LV \oplus L\R^m},s_1)~=~(\alpha_M\wedge \alpha_M)\Delta_{S^m}~,~
\gamma^{\ZZ_2}(\epsilon_{LV \oplus L\R^m},s_2)~=~\Delta_M\alpha_M\\[1ex]
\hspace*{25pt}\in \omega^0_{\ZZ_2}(W,\partial W;-\epsilon_{LV \oplus L\R^m})~
=~\omega^{\ZZ_2}_m(W)~=~\omega_m^{\ZZ_2}(M \times M)
\end{array}$$
and $\delta(s_1,s_2)$ has image
$$\begin{array}{l}
\gamma^{\ZZ_2}(\epsilon_{LV\oplus L\R^m},s_1)-
\gamma^{\ZZ_2}(\epsilon_{LV\oplus L\R^m},s_2)~
=~(\alpha_M\wedge \alpha_M)\Delta_{S^m}-\Delta_M\alpha_M\\[1ex]
\hspace*{15pt}
\in {\rm im}(\omega^{-1}_{\ZZ_2}(\partial W;-\epsilon_{LV\oplus L\R^m})
\to \omega^0_{\ZZ_2}(W,\partial W;-\epsilon_{LV\oplus L\R^m}))\\[1ex]
\hspace*{20pt}
={\rm ker}(\omega^0_{\ZZ_2}(W,\partial W;-\epsilon_{LV\oplus L\R^m})
\to\omega^0_{\ZZ_2}(W;-\epsilon_{LV\oplus L\R^m}))
\subset \omega_{2m,m}(M \times M)~.
\end{array}$$
Also, passing to the limit as ${\rm dim}(V) \to \infty$ note that
$$\begin{array}{l}
\varinjlim\limits_{{\rm dim}(V)} \psi_V(\alpha_M)~=~(\alpha_M \wedge \alpha_M)\Delta_{S^m} - \Delta_M \alpha_M\\[1ex]
\hspace*{50pt}
\in {\rm im}(\omega_m(S(\infty)\times_{\ZZ_2} (M \times M) )
\to \omega^{\ZZ_2}_{2m,m}(M \times M))\\[1ex]
\hspace*{150pt}
=~{\rm ker}(\rho:\omega^{\ZZ_2}_{2m,m}(M \times M)\to\omega_m(M))
\end{array}$$
is determined by
$(\alpha_M \wedge \alpha_M)\Delta_{S^m},\Delta_M \alpha_M\in \omega^{\ZZ_2}_m(M \times M)$,
since
$$\omega_m(S(\infty)\times_{\ZZ_2} (M \times M) ) \to
\omega^{\ZZ_2}_{2m,m}(M \times M)$$
is a split injection.\hfill\qed}
\end{example}

\section{Double points}

\begin{definition}~\label{double1} {\rm
{\rm (i)} The {\it fibrewise product} of maps $f_1:X_1 \to Y$, $f_2:X_2 \to Y$ is the
map
$$f_1 \times_Yf_2~:~X_1 \times_Y X_2 \to Y~;~(x_1,x_2) \mapsto f(x_1)=f(x_2)$$
with
$$X_1 \times_Y X_2~=~\{(x_1,x_2) \in X_1 \times X_2\st f_1(x_1)=f_2(x_2) \in Y\}$$
the pullback
$$\xymatrix{X_1 \times_Y X_2 \ar[r] \ar[d] \ar[dr]^-{f_1\times_Y f_2} &
X_1 \ar[d]^-{f_1} \\
X_2 \ar[r]^-{f_2} & Y}$$
{\rm (ii)} The {\it square} of a map $f:X \to Y$ is the $ \ZZ_2$-equivariant map\index{square!map $f \times f$}
$$f \times f~:~ X \times X \to Y \times Y~;~(x_1,x_2) \mapsto (f(x_1),f(x_2))$$
with $T \in \ZZ_2$ acting on $X \times X$ by
$$T~:~X\times X  \to X\times X~;~(x_1,x_2) \mapsto (x_2,x_1)$$
and similarly for $Y \times Y$.\\
{\rm (iii)} The {\it double point space} of $f:X \to Y$  is the $\ZZ_2$-space\index{double point set!$X\times_YX$}
$$X\times_YX~=~(f\times f)^{-1}(Y)~=~\{(x_1,x_2) \in X \times X\st f(x_1)=f(x_2) \in Y\} \subset X \times X~.$$
The diagonal embedding
$$\Delta_X~:~X \to X\times_YX~;~x \mapsto (x,x)$$
is $\ZZ_2$-equivariant, with image the fixed point set
$$(X\times_YX)^{\ZZ_2}~=~X~=~\{(x,x) \st x \in X\}~.$$
The map
$$f \times_Yf~=~(f\times f)\vert~:~X \times_YX \to Y~;~(x_1,x_2) \mapsto f(x_1)=f(x_2)~.$$
is also $\ZZ_2$-equivariant, with $(f \times_Y f)\Delta_X=f:X \to Y$.\\
{\rm (iv)}  Let $\widetilde{Y}\to Y$ be a regular cover of $Y$ with group of
covering translations $\pi$, and given a map $f:X \to Y$ let
$$\widetilde{X}~=~X\times_Y\widetilde{Y}~=~f^*\widetilde{Y}~=~\{(x,\widetilde{y}) \in X \times \widetilde{Y}\st f(x)=[\widetilde{y}] \in Y\}$$
be the pullback cover of $X$.
The {\it assembly map} is the inclusion\index{assembly map, $A$}
$$A~:~X \times_Y X \to \widetilde{X}\times_\pi \widetilde{X}~;~(x_1,x_2) \mapsto [(x_1,\widetilde{y}),(x_2,\widetilde{y})]$$
with $\widetilde{y} \in \widetilde{Y}$ any lift of $y=f(x_1)=f(x_2)\in Y$.\\
{\rm (v)} The {\it ordered double point space} of $f:X \to Y$ is the
free $\ZZ_2$-space\index{double point set!ordered, $D_2(f)$}
$$D_2(f)~=~\{(x_1,x_2) \in X \times X\st x_1 \neq x_2 \in
X,f(x_1)=f(x_2) \in Y\}$$
with
$$X\times_Y X~=~X \sqcup D_2(f)~.$$
{\rm (vi)} The {\it unordered double point set} of  $f:X \to Y$ is the quotient space\index{double point set!unordered, $D_2[f]$}
$$D_2[f]~=~D_2(f)/\ZZ_2~.$$
The projection $D_2(f) \to D_2[f]$ is a double cover.
\hfill\qed}
\end{definition}

\begin{remark} {\rm The diagonal $Y \subseteq Y \times Y$ is a closed
subspace (since $Y$ is Hausdorff), so that $X \times_YX=(f\times f)^{-1}(Y) \subseteq X \times X$ is also closed, and $X \subseteq X \times_YX$ is also closed (since $X$ is Hausdorff). We shall always assume that $X \subseteq X \times_Y X$ is
also open, so that
$$X \times_YX~=~X \sqcup D_2(f)$$
as a topological space, with one-point compactification
$$(X \times_YX)^{\infty}~=~X^{\infty} \vee D_2(f)^{\infty}~.$$
\hfill\qed}
\end{remark}

We shall be concerned with the double points of a map $f:X \to Y$
which extends to an embedding $V \times X \emb V \times Y$ (in the manner
of a codimension 0 immersion of manifolds) :

\begin{definition}~\label{embed0} {\rm
(i) An {\it embedding} of a map $f:X \to Y$ is an open embedding of the type
$$e~:~V \times X \emb V \times Y~;~(v,x) \mapsto (g(v,x),f(x))$$
for an inner product space $V$, so that the diagram\index{embedding of map}
$$\xymatrix@R+10pt{V \times X \ar@{^{(}->}[r]^-{\di{e}}
\ar[d]_-{\di{\rm proj.}}&
V \times Y \ar[d]^-{\di{\rm proj.}}\\
X \ar[r]^-{\di{f}} & Y}$$
commutes. For each $y \in Y$ the restriction of $e$ defines an embedding
$$e\vert~:~V \times f^{-1}(y) \to V~;~(v,x) \mapsto g(v,x)~.$$
The pair $e=(g,f)$ is an {\it embedded map}.\index{embedded map}\\
(ii) The {\it compactification Umkehr} of an embedded map $e=(g,f)$ is the
compactification Umkehr given by \ref{adjunct}\index{Umkehr!compactification}
$$\begin{array}{l}
F~:~V^{\infty} \wedge Y^{\infty}\to V^{\infty} \wedge X^{\infty}~;\\[1ex]
\hskip50pt (w,y) \mapsto
\begin{cases} (v,x)&{\rm if}~(w,y)=e(v,x) \in e(V \times X) \\
*&{\rm otherwise.}
\end{cases}
\end{array}
$$
\hfill\qed}
\end{definition}

Let $e=(g,f):V \times X \emb V \times Y$ be an embedded map, with
compactification Umkehr map $F:V^{\infty}\wedge Y^{\infty} \to
V^{\infty}\wedge X^{\infty}$.  We shall now construct a {\it local
geometric Hopf map}\index{local geometric Hopf map, $h_V(F)_Y$}
$$\begin{array}{l}
h_V(F)_Y~=~F_1 \vee \delta(F_2,F_3)~:~\Sigma S(LV)^+ \wedge V^{\infty} \wedge Y^{\infty}\\[1ex]
\hskip50pt
\to LV^{\infty} \wedge V^{\infty} \wedge (X \times_Y X)^{\infty}~=~
LV^{\infty}\wedge V^{\infty} \wedge (D_2(f)^{\infty}\vee X^{\infty})
\end{array}$$
with
$$F_1~:~\Sigma S(LV)^+ \wedge V^{\infty} \wedge Y^{\infty}
\to LV^{\infty} \wedge V^{\infty} \wedge D_2(f)^{\infty}$$
the Umkehr map of a $\ZZ_2$-equivariant embedding
$$e_1~:~LV \times V \times D_2(f) \emb
(LV\backslash \{0\}) \times V \times Y~,$$
and
$$\delta(F_2,F_3)~:~\Sigma S(LV)^+ \wedge V^{\infty} \wedge Y^{\infty}
\to LV^{\infty} \wedge V^{\infty} \wedge X^{\infty}$$
the relative difference of the Umkehr maps
$$F_2~,~F_3~:~LV^{\infty}\wedge V^{\infty} \wedge Y^{\infty} \to
LV^{\infty} \wedge V^{\infty} \wedge X^{\infty}$$
of two $\ZZ_2$-equivariant embeddings
$e_2,e_3:LV \times V \times X \emb LV \times V \times Y$ of $f:X \to Y$ such that
$$e_2\vert~=~e_3\vert~=~e~:~\{0\} \times V \times X \emb \{0\} \times V \times Y~.$$
The geometric Hopf invariant map is the assembly of the local geometric
Hopf map
$$\begin{array}{l}
h_V(F)~=~(1 \wedge A)h_V(F)_Y~:\\[1ex]
\hskip50pt \Sigma S(LV)^+ \wedge V^{\infty} \wedge Y^{\infty}
\to LV^{\infty} \wedge V^{\infty} \wedge (X \times X)^{\infty}
\end{array}$$
with $\widetilde{Y}=Y$, $\pi=\{1\}$, $A:X \times_Y X \subset X \times X$.

\begin{terminology}~ \label{Umkehrterm}
{\rm Given a map $f:X \to Y$ and an embedding $e=(g,f):V \times X\emb V\times Y$ define a $\ZZ_2$-equivariant embedding of $f \times f:X \times X \to Y \times Y$
$$\begin{array}{l}
e \times e~=~(g\times g,f\times f)~:\\[1ex]
(V \times V) \times (X \times X) \emb (V \times V) \times (Y \times Y)~;\\[1ex]
\hphantom{e \times e~=~(g\times g,f\times f)~:~}
(u,v,x,y) \mapsto (g(u,x),g(v,y),f(x),f(y))
\end{array}$$
with $\ZZ_2$-equivariant compactification Umkehr map
$$F \wedge F~:~
V^{\infty} \wedge V^{\infty} \wedge Y^{\infty} \wedge Y^{\infty} \to
V^{\infty} \wedge V^{\infty} \wedge X^{\infty} \wedge X^{\infty}~.$$
The conjugate of $e \times e$ by $\kappa_V \times 1$ is a $\ZZ_2$-equivariant
embedding of $f:X \to Y$
$$\begin{array}{l}
e'~=~(\kappa_V \times 1)^{-1}(e \times e)(\kappa_V \times 1)~:~
LV \times V \times X \times X \emb LV \times V \times Y\times Y~;\\[1ex]
 (u,v,x,y) \mapsto (g(u+v,y)-g(-u+v,x))/2,(g(u+v,y)+g(-u+v,x))/2,f(x),f(y))~.
\end{array}$$
with $\ZZ_2$-equivariant compactification Umkehr map
$$\begin{array}{l}
F'~=~(\kappa_V \wedge 1)(F \wedge F)(\kappa_V \wedge 1)^{-1}~:\\[1ex]
\hskip25pt LV^{\infty} \wedge V^{\infty} \wedge Y^{\infty} \wedge Y^{\infty} \to
LV^{\infty} \wedge V^{\infty} \wedge X^{\infty} \wedge X^{\infty}~.
\end{array}$$
The restriction of $e'$ defines a $\ZZ_2$-equivariant  open embedding
$$\begin{array}{l}
e_1~=~e'\vert~:~
LV \times V \times D_2(f) \emb  (LV \backslash \{0\}) \times V \times Y~;~(u,v,x,y) \mapsto\\[1ex]
 (g(u+v,y)-g(-u+v,x))/2,(g(u+v,y)+g(-u+v,x))/2,f(x)=f(y))
\end{array}$$
with a $\ZZ_2$-equivariant compactification Umkehr map
$$F_1~:~(LV \backslash \{0\})^{\infty} \wedge V^{\infty} \wedge Y^{\infty}~=~
\Sigma S(LV)^+ \wedge V^{\infty} \wedge Y^{\infty} \to
LV^{\infty} \wedge V^{\infty} \wedge D_2(f)^\infty~.$$
The components of the  $\ZZ_2$-equivariant embeddings of $f:X \to Y$
$$e_2~,~e_3~:~LV \times V \times X \emb  LV \times V \times Y$$
defined by
$$\begin{array}{l}
e_2(u,v,x)~=~(1\times e)(u,v,x)~=~(u,g(v,x),f(x))~,\\[1ex]
e_3(u,v,x)~=~e'(u,v,x,x)\\[1ex]
\hskip10pt =~(g(u+v,x)-g(-u+v,x))/2,(g(u+v,x)+g(-u+v,x))/2,f(x))
\end{array}$$
have $\ZZ_2$-equivariant compactification Umkehr maps
$$F_2~=~1 \wedge F~,~F_3~:~LV^{\infty} \wedge V^{\infty} \wedge Y^{\infty}\to
LV^{\infty} \wedge V^{\infty} \wedge X^{\infty}$$
with
$$\begin{array}{l}
e_2\vert~=~e_3\vert~=~e~:~
\{0\} \times V \times X \emb  \{0\} \times V \times Y~,\\[1ex]
F_2\vert~=~F_3 \vert~:~\{0\}^{\infty} \wedge V^{\infty} \wedge
Y^{\infty} \to LV^{\infty} \wedge V^{\infty} \wedge X^{\infty}~.
\end{array}$$
\hfill\qed}
\end{terminology}

\begin{theorem}  (Double Point Theorem) \label{doubleHopf}\\
For an embedding $e=(g,f):V \times X \emb V \times Y$ of a map
$f:X \to Y$ the geometric Hopf invariant $h_V(F)$ of the Umkehr map
$F:V^{\infty} \wedge Y^{\infty} \to V^{\infty} \wedge X^{\infty}$
is given up to natural $\ZZ_2$-equivariant homotopy by the composite
$$\begin{array}{l}
h_V(F)~=~(1 \wedge A)h_V(F)_Y~:~\Sigma S(LV)^+ \wedge V^{\infty} \wedge Y^{\infty} \\[1ex]
\hskip50pt \to LV^{\infty} \wedge V^{\infty} \wedge (X\times_YX)^{\infty} \to
LV^{\infty} \wedge V^{\infty} \wedge (X\times X)^{\infty}
\end{array}
$$
with
$$\begin{array}{l}
h_V(F)_Y~=~(1 \wedge i)F_1 \vee (1 \wedge \Delta_X)\delta(F_2,F_3)~:\\[1ex]
\Sigma S(LV)^+ \wedge V^{\infty} \wedge Y^{\infty} \to\\[1ex]
\hskip25pt
LV^{\infty} \wedge V^{\infty} \wedge (X\times_YX)^{\infty}~=~
LV^{\infty}\wedge V^{\infty} \wedge (D_2(f)^{\infty}\vee X^{\infty})~,\\[1ex]
i~=~{\rm inclusion}~:~D_2(f) \to X\times_YX ~;~(x_1,x_2) \mapsto (x_1,x_2)~,\\[1ex]
\Delta_X~=~{\rm diagonal}~:~X \to X \times_YX~;~x \mapsto (x,x)~,\\[1ex]
A~=~{\rm assembly}~=~\Delta_X \sqcup i~:~X \times_YX~=~X \sqcup D_2(f) \to X \times X~.
\end{array}$$
\end{theorem}
\begin{proof} By definition
$$h_V(F)~=~\delta(p,q)~:~\Sigma S(LV)^+ \wedge V^{\infty}
\wedge Y^{\infty} \to LV^{\infty} \wedge V^{\infty} \wedge
X^{\infty}\wedge X^{\infty}$$
with $p=(1\wedge \Delta_{X^{\infty}})F_2$, $q=F'(1 \wedge \Delta_{Y^{\infty}})$ the composites in the (noncommutative) square of $\ZZ_2$-equivariant maps
$$\xymatrix@R+20pt@C+20pt{
LV^{\infty} \wedge V^{\infty} \wedge Y^{\infty} \ar[r]^-{\di{1
\wedge \Delta_{Y^{\infty}}}} \ar[d]_-{\di{F_2}} &LV^{\infty}
\wedge V^{\infty} \wedge Y^{\infty} \wedge Y^{\infty}
\ar[d]^-{\di{F'}} \\
LV^{\infty} \wedge V^{\infty} \wedge X^{\infty} \ar[r]^-{\di{1
\wedge \Delta_{X^{\infty}}}} &LV^{\infty} \wedge V^{\infty} \wedge
    X^{\infty}\wedge X^{\infty} }$$
with
$$p\vert~=~q\vert~~:~\{0\}^{\infty}\wedge V^{\infty} \wedge
Y^{\infty} \to LV^{\infty}\wedge V^{\infty} \wedge X^{\infty}
\wedge X^{\infty}~.$$
 Consider the pullback commutative square of $\ZZ_2$-equivariant open embeddings and
 inclusions
$$\xymatrix@R+20pt@C+60pt{LV \times V \times( X\times_YX)
    \ar[d]_-{\di{1 \times  A}} \ar[r]^-{\di{e_1 \sqcup e_3}} &
    LV \times V  \times Y
    \ar[d]^-{\di{1 \times \Delta_{Y}}} \\
LV \times V \times X \times X \ar[r]^-{\di{e'}} & LV \times
    V \times Y \times Y}$$
which determines a commutative square of $\ZZ_2$-equivariant Umkehr maps
and inclusions
$$\xymatrix@R+20pt@C+30pt{
LV^{\infty} \wedge V^{\infty} \wedge Y^{\infty} \ar[d]_-{\di{1
\wedge \Delta_{Y^{\infty}}}} \ar[r]^-{\di{F_{13}}} &
LV^{\infty} \wedge V^{\infty} \wedge (X \times_YX)^\infty
\ar[d]^-{\di{1 \wedge A}} \\
LV^{\infty} \wedge V^{\infty} \wedge Y^{\infty} \wedge Y^{\infty} \ar[r]^-{\di{F'}}
&LV^{\infty} \wedge V^{\infty} \wedge
    X^{\infty}\wedge X^{\infty} }$$
with $F_{13}$ the Umkehr map of $e_1 \sqcup e_3$. Thus
$$\begin{array}{ll}
h_V(F)&=~\delta(p,q)~=~\delta((1\wedge \Delta_{X^{\infty}})F_2,
F'(1 \wedge \Delta_{Y^{\infty}}))\\[1ex]
&=~\delta((1\wedge \Delta_{X^{\infty}})F_2,
(1 \wedge (i \vee \Delta_{X^{\infty}}))F_{13})~:\\[1ex]
&\hskip20pt
\Sigma S(LV)^+ \wedge V^{\infty} \wedge Y^{\infty}
\to LV^{\infty} \wedge V^{\infty} \wedge X^{\infty} \wedge
X^{\infty}~.
\end{array}$$
Proposition \ref{difcon2} (iv) gives a natural $\ZZ_2$-equivariant homotopy
$$\begin{array}{l}
h_V(F) \simeq (1\wedge A\Delta_{X^{\infty}})\delta(F_2,F_3)+(1 \wedge A)
\delta((1 \wedge \Delta_{X^{\infty}})F_3,(1 \wedge (i \vee \Delta_{X^{\infty}}))F_{13})\\[1ex]
\hphantom{h_V(F) \simeq}~=~
(1\wedge A\Delta_{X^{\infty}})\delta(F_2,F_3)+
\delta(p',q')~:\\[1ex]
\hskip50pt
\Sigma S(LV)^+ \wedge V^{\infty} \wedge Y^{\infty}
\to LV^{\infty} \wedge V^{\infty} \wedge X^{\infty} \wedge
X^{\infty}
\end{array}$$
with
$$\begin{array}{l}
p'~=~(1 \wedge (\Delta_{X}\vert)^{\infty})F_3~,~q'~=~
(1 \wedge (i \sqcup \Delta_X \vert)^{\infty})F_{13}~:\\[1ex]
\hskip100pt
LV^{\infty} \wedge V^{\infty} \wedge Y^{\infty} \to
LV^{\infty} \wedge V^{\infty} \wedge (X\times_YX)^{\infty}
\end{array}$$
the composites in the (noncommutative) diagram of $\ZZ_2$-equivariant maps
$$\xymatrix@R+20pt@C+30pt{
LV^{\infty} \wedge V^{\infty} \wedge Y^{\infty} \ar[d]_-{\di{F_3}}
\ar[r]^-{\di{F_{13}}} &
LV^{\infty} \wedge V^{\infty} \wedge (D_2(f)^\infty \vee X^{\infty})
\ar[d]^-{1 \wedge(i \sqcup\Delta_X\vert)^{\infty}} \\
LV^{\infty} \wedge V^{\infty} \wedge X^{\infty}
\ar[r]^-{1 \wedge (\Delta_{X}\vert)^{\infty}} &LV^{\infty} \wedge V^{\infty} \wedge
    (X\times_YX)^{\infty} }$$
    with  $\Delta_X\vert : X \to X\times_YX;x \mapsto (x,x)$ the diagonal map. The neighbourhood
$$U~=~e_3(LV \times V \times X)\cup \{*\} \subset
LV^{\infty} \wedge V^{\infty} \wedge Y^{\infty}$$
of $0^+ \wedge V^{\infty} \wedge Y^{\infty}$ is such that
for all $t \in LV^{\infty}$, $x \in V^{\infty}\wedge Y^{\infty}$
$$p'(t,x)~=~\begin{cases} q'(t,x)&{\rm if}~(t,x) \in U\\
*&{\rm if}~(t,x) \in
\overline{(LV^{\infty} \wedge V^{\infty} \wedge Y^{\infty})\backslash U}~.
\end{cases}$$
Proposition \ref{difcon2} (v) gives a natural $\ZZ_2$-equivariant homotopy
$$\delta(p',q')~\simeq~q''~:~
\Sigma S(LV)^+ \wedge V^{\infty} \wedge Y^{\infty}
\to LV^{\infty} \wedge (X\times_YX)^{\infty}$$
with $q''$ defined by
$$\begin{array}{l}
q''~:~\Sigma  S(LV)^+\wedge V^{\infty} \wedge Y^{\infty}\to LV^{\infty}
\wedge (X\times_YX)^{\infty}~;\\[1ex]
\hphantom{q''~:~}(t,x) \mapsto \begin{cases}
*&{\rm if}~ (t,x) \in U\\
q'(t,x)&{\rm if}~ (t,x) \in \overline{
(LV^{\infty} \wedge V^{\infty} \wedge Y^{\infty})\backslash U}~.
\end{cases}
    \end{array}$$
But $q''=(1\wedge i)F_1$, giving the required result.
\hfill\qed\end{proof}

Similarly for the stable geometric Hopf invariant $h'_V(F)$~:

\begin{proposition}~ \label{doubleHopf2}
Let $e=(g,f):V \times X \emb V \times Y$ be an embedded map
$$\xymatrix@R+10pt{V \times X \ar@{^{(}->}[r]^-{\di{e}}
\ar[d]_-{\di{\rm proj.}}&
V \times Y \ar[d]^-{\di{\rm proj.}}\\
X \ar[r]^-{\di{f}} & Y}$$
with Umkehr map $F:V^{\infty} \wedge Y^{\infty} \to
V^{\infty} \wedge X^{\infty}$.
The stable geometric Hopf invariant of $F$ is given up to
natural $\ZZ_2$-equivariant homotopy by the sum
$$h'_V(F)~=~(1\wedge A)h'_V(F)_Y~:~
Y^{\infty} \sto S(LV)^+ \wedge X^{\infty} \wedge X^{\infty}$$
with
$$h'_V(F)_Y~=~cF'_1+(1\wedge \Delta_{X^{\infty}})\delta'(F_2,F_3)~:~
Y^{\infty} \sto S(LV)^+ \wedge (X\times_YX)^{\infty}$$
$$c F'_1:~Y^{\infty}
\xymatrix{\ar[r]^-{\di{F'_1}}&}
D_2(f)^\infty  \xymatrix{\ar[r]^-{\di{c}}&} S(LV)^+ \wedge (X\times_YX)^{\infty}$$
the composite of the $\ZZ_2$-equivariant Umkehr map
$$F'_1~:~LV^{\infty} \wedge V^{\infty} \wedge
Y^{\infty} \to LV^{\infty} \wedge V^{\infty} \wedge D_2(f)^\infty$$
of the $\ZZ_2$-equivariant open embedding
$$e'_1~=~(\kappa_V \times 1)^{-1}(e \times e)\vert(\kappa_V \times 1)~:~
LV \times V \times D_2(f) \emb LV \times V \times Y~,$$
and the $\ZZ_2$-equivariant map $c$ defined by
$$c~:~D_2(f)^\infty \to  S(LV)^+ \wedge (X\times_YX)^{\infty}~;~
(x_1,x_2) \mapsto \bigg(
\di{g(0,x_1)-g(0,x_2) \over \Vert g(0,x_1)-g(0,x_2) \Vert},x_1,x_2\bigg)~,$$
and $\delta'(F_3,F_2):Y^{\infty} \sto S(LV)^+ \wedge X^{\infty}$ the
stable relative difference of the $\ZZ_2$-equivariant Umkehr maps
$$F_2~,~F_3~:~LV^{\infty} \wedge V^{\infty} \wedge Y^{\infty}
\to LV^{\infty} \wedge V^{\infty} \wedge X^{\infty}$$
of \ref{Umkehrterm}.
\end{proposition}
\begin{proof} The $\ZZ_2$-equivariant open embedding
$e'_1$ is the composite
$$e'_1~:~LV \times V \times D_2(f) \xymatrix{\ar[r]^-{\di{e_1}}&}
(LV\backslash \{0\}) \times V \times Y  \emb LV \times V \times
Y$$ with $e_1=(\kappa_V \times 1)^{-1}(e \times e)\vert(\kappa_V
\times 1)$ as in  \ref{Umkehrterm}, and the Umkehr of
$e'_1$ is the composite
$$\begin{array}{l}
F'_1~=~F_1(\alpha_{LV} \wedge 1)~:~LV^{\infty}\wedge V^{\infty} \wedge Y^{\infty}\\[1ex]
\xymatrix@C+10pt{\ar[r]^-{\di{\alpha_{LV} \wedge 1}}&}
(LV\backslash \{0\})^{\infty}\wedge V^{\infty} \wedge Y^{\infty}
=\Sigma S(LV)^+\wedge V^{\infty} \wedge Y^{\infty}\\[1ex]
\hskip150pt \xymatrix{\ar[r]^-{\di{F_1}}&} LV^{\infty}\wedge
V^{\infty}\wedge D_2(f)^{\infty}
\end{array}$$
with $F_1$ the Umkehr map of $e_1$ as in \ref{Umkehrterm}.
Theorem \ref{doubleHopf} gives a natural $\ZZ_2$-equivariant homotopy
$$\begin{array}{l}
h_V(F)~\simeq~(1 \wedge i)F_1+(1\wedge \Delta_{X^{\infty}})\delta(F_2,F_3)~:\\[1ex]
\hskip50pt
\Sigma S(LV)^+\wedge V^{\infty} \wedge Y^{\infty} \to
LV^\infty \wedge V^{\infty} \wedge (X\times_YX)^{\infty}
\end{array}$$
with $i={\rm inclusion}:D_2(f) \to X \times_Y X$, so
that there is defined a natural $\ZZ_2$-equivariant homotopy
$$\begin{array}{l}
h'_V(F)_Y~=~cF'_1+(1 \wedge \Delta_{X^{\infty}})\delta'(F_2,F_3)\\[1ex]
\hskip50pt\simeq~((1 \wedge i)(1 \wedge F_1)(1\wedge \Delta_{X^{\infty}})+ \delta(F_2,F_3))(\Delta \alpha_{LV} \wedge 1)~:\\[1ex]
\hskip50pt
LV^{\infty} \wedge V^{\infty} \wedge Y^{\infty} \to S(LV)^+
\wedge LV^{\infty} \wedge V^{\infty} \wedge (X\times_YX)^{\infty}~.
\end{array}$$
Now
$$\begin{array}{l}
(1 \wedge i)(1 \wedge F_1)(\Delta \alpha_{LV} \wedge 1)~:\\[1ex]
LV^{\infty} \wedge V^{\infty} \wedge Y^{\infty} \to S(LV)^+
\wedge LV^{\infty} \wedge V^{\infty} \wedge (X\times_YX)^{\infty}~;\\[1ex]
((t,u),v,y) \mapsto
\begin{cases}
(u,v_1,v_2,x_1,x_2)\\[1ex]
\hskip10pt{\rm if}~[t,u]=(g(v_1+v_2,x_2)-g(-v_1+v_2,x_1))/2 \in
LV^{\infty}\backslash\{0\},\\[1ex]
\hskip10pt\hphantom{{\rm if}~} v=(g(v_1+v_2,x_2)+g(-v_1+v_2,x_1))/2 \in V,\\[1ex]
\hskip10pt\hphantom{{\rm if}~} y=f(x_1)=f(x_2)\\[1ex]
*\hskip15pt{\rm otherwise}
\end{cases}
\end{array}$$
and
$$\begin{array}{l}
(1 \wedge c)(1 \wedge F_1)(\alpha_{LV} \wedge 1)~:\\[1ex]
LV^{\infty} \wedge V^{\infty} \wedge Y^{\infty} \to S(LV)^+
\wedge LV^{\infty} \wedge V^{\infty} \wedge (X\times_YX)^{\infty}~;\\[1ex]
((t,u),v,y) \mapsto
\begin{cases}
(\di{g(0,x_2)-g(0,x_1) \over \Vert g(0,x_2)-g(0,x_1)\Vert},v_1,v_2,x_1,x_2)\\[2ex]
\hskip10pt {\rm if}~[t,u]=(g(v_1+v_2,x_2)-g(-v_1+v_2,x_1))/2 \in
LV^{\infty}\backslash\{0\},\\[1ex]
\hskip10pt \hphantom{{\rm if}~} v=(g(v_1+v_2,x_2)+g(-v_1+v_2,x_1))/2 \in V,\\[1ex]
\hskip10pt \hphantom{{\rm if}~} y=f(x_1)=f(x_2)\\[1ex]
*\hskip15pt {\rm otherwise}~.
\end{cases}
\end{array}$$
The map
$$\begin{array}{l}
h~:~I \times LV^{\infty} \wedge V^{\infty} \wedge Y^{\infty}
\to S(LV)^+
\wedge LV^{\infty} \wedge V^{\infty} \wedge (X\times_YX)^{\infty}~;\\[1ex]
(s,(t,u),v,y) \mapsto
\begin{cases}
(\di{g(sv_1+sv_2,x_2)-g(-sv_1+sv_2,x_1) \over
\Vert g(sv_1+sv_2,x_2)-g(-sv_1+sv_2,x_1)\Vert},v_1,v_2,x_1,x_2)\\[2ex]
\hskip10pt {\rm if}~[t,u]=(g(v_1+v_2,x_2)-g(-v_1+v_2,x_1))/2 \in
LV^{\infty}\backslash\{0\}\\[1ex]
\hskip10pt \hphantom{{\rm if}~} v=(g(v_1+v_2,x_2)+g(-v_1+v_2,x_1))/2 \in V,\\[1ex]
\hskip10pt \hphantom{{\rm if}~} y=f(x_1)=f(x_2)\\[1ex]
* \hskip15pt {\rm otherwise}
\end{cases}
\end{array}$$
defines a natural $\ZZ_2$-equivariant homotopy
$$\begin{array}{l}
h~:~(1 \wedge i)(1 \wedge F_1)(\Delta \alpha_{LV} \wedge 1)~
\simeq~(1 \wedge c)(1 \wedge F_1)(\alpha_{LV} \wedge 1)~:\\[1ex]
\hphantom{h~:~}
LV^{\infty} \wedge V^{\infty} \wedge Y^{\infty} \to S(LV)^+
\wedge LV^{\infty} \wedge V^{\infty} \wedge (X \times_YX)^{\infty}~.
\end{array}$$
The concatenation of the homotopies is a natural $\ZZ_2$-equivariant
homotopy
$$\begin{array}{l}
h'_V(F)_Y~\simeq~(1 \wedge c)(1 \wedge F_1)(\alpha_{LV} \wedge 1)+
\delta(F_2,F_3)(\Delta\alpha_{LV}\wedge 1)~=~
(1 \wedge c)(1 \wedge F'_1):\\[1ex]
\hphantom{h~:~}
LV^{\infty} \wedge V^{\infty} \wedge Y^{\infty} \to S(LV)^+
\wedge LV^{\infty} \wedge V^{\infty} \wedge (X\times_YX)^{\infty} ~.
\end{array}$$
Finally, identify
$$\delta(F_2,F_3)(\Delta \alpha_{LV} \wedge 1)~=~
\delta'(F_2,F_3)~:~LV^{\infty} \wedge V^{\infty} \wedge Y^{\infty}
\to LV^{\infty} \wedge V^{\infty} \wedge X^{\infty}~.$$
\hfill\qed\end{proof}

\section{Positive embeddings}

We shall now formulate a positivity condition on an embedded map
$$e~=~(g,f)~:~V \times X \emb V \times Y$$
which ensures that the embedding term $(1\wedge
\Delta_{X^{\infty}})\delta(F_2,F_3)$ in the expression given by Theorem
\ref{doubleHopf} for the geometric Hopf invariant $h_V(F)$ of the
compactification Umkehr map $F:V^{\infty} \wedge Y^{\infty} \to V^{\infty}
\wedge X^{\infty}$ is 0, so that $h_V(F)$ only depends on the double point
$\ZZ_2$-space $D_2(f)$ with $h_V(F)=(i \wedge 1)F_1$.

\begin{definition}~ \label{positive}
{\rm (i) A linear automorphism $a \in GL(V)$ of an
inner product space $V$ is {\it positive} if\index{positive!linear automorphism, $GL^+(V)$} if $a$ has no negative real eigenvalues,
i.e. $a+uI$ should be invertible for all $u>0$ ($u \in \R$).\\
(ii) Let $GL^+(V) \subseteq GL(V)$ be the subset of the positive linear
automorphisms.\\
(iii) Let
$$O^+(V)~=~O(V) \cap GL^+(V) \subseteq O(V)$$
 be the subset of the positive orthogonal automorphisms\index{positive! orthogonal automorphism, $O^+(V)$}.}
\end{definition}

\begin{proposition} \label{contract}
{\rm (i)} The space $GL^+(V)$ is contractible, with $I \in GL^+(V)$.\\
{\rm (ii)}  The following conditions on an orthogonal
automorphism $a:V \to V$ are equivalent:
\begin{itemize}
\item[{\rm (a)}]~$a$ is positive,
\item[{\rm (b)}]~$1+a:V \to V$ is a linear automorphism,
\item[{\rm (c)}]~$-1$ is not an eigenvalue of $a$.
\end{itemize}
{\rm (iii)}  The space $O^+(V)$ is contractible.
\end{proposition}
\begin{proof} (i) For any $a \in GL^+(V)$ and $t \in [0,1]$ the linear map
$$H_t(a)~:~V \to V~;~v \mapsto (1-t)a(v)+tv$$
defines a path in $GL^+(V)$ from $H_0(a)=a$ to $H_1(a)=I$: for
any $u>0$
$$H_t(a)+uI~:~V \to V~;~v \mapsto (1-t)a(v)+(u+t)v$$
is invertible because $u+t>0$.\\
{\rm (ii)} Immediate from the definitions.\\
{\rm (iii)} $O^+(V)$ is homeomorphic to the
contractible space ${\rm End}^-(V)$ of linear endomorphisms $b:V \to V$
such that
$$\langle b(u),v \rangle~=~-\langle u,b(v) \rangle  ~~(u,v \in V)~,$$
with the Cayley parametrization (Weyl \cite[II.10]{weyl})
of positive orthogonal automorphisms defining inverse homeomorphisms
$$\begin{array}{l}
O^+(V) \to {\rm End}^-(V)~;~a \mapsto (1+a)^{-1}(1-a)~,\\[1ex]
{\rm End}^-(V) \to O^+(V)~;~b \mapsto (1+b)^{-1}(1-b)~.
\end{array}$$
\hfill\qed\end{proof}

\begin{definition}~ {\rm An embedding $e=(g,f):V \times X \emb V \times Y$
of $f:X \to Y$ is
{\it positive} if for each $x \in X$ the embedding
\index{positive!embedding}
$$g_x~:~V \to V~;~u \mapsto g(u,x)$$
is differentiable and the differentials $dg_x(v):V \to V$ ($v \in
V$) are positive linear automorphisms.
\hfill\qed}
\end{definition}

\begin{example}{\rm For any map $c:X \to O^+(V)$ the embedding
$$V \times X \emb V \times X~;~(v,x) \mapsto (c(x)(v),x)$$
is positive, with each
$$g_x~:~V \to V~;~v \mapsto c(x)(v)~(x \in X)$$
a linear map with $dg_x=c(x) \in O^+(V)$.\\
\hfill\qed}
\end{example}

The embeddings $g_x:V \to V$ ($x \in X$) in a positive embedding
$e=(g,f):V \times X \to V \times Y$ are orientation-preserving.

\begin{theorem} \label{doubleHopfpositive}
Let $e=(g,f):V \times X \emb V \times Y$ be a positively embedded
map, with adjunction Umkehr map $F:V^{\infty}\wedge Y^{\infty} \to
V^{\infty}\wedge X^{\infty}$. \\
{\rm (i)} The local geometric Hopf invariant
is given up to natural $\ZZ_2$-equivariant homotopy by
$$h_V(F)_Y~=~(1 \wedge i)F_1~:~
 \Sigma S(LV)^+ \wedge V^{\infty} \wedge Y^{\infty} \to
LV^{\infty} \wedge V^{\infty} \wedge (X\times_YX)^{\infty}$$
with $i={\it inclusion}:D_2(f) \to X \times_Y X$, and
$$F_1~:~(LV \backslash \{0\})^{\infty} \wedge V^{\infty} \wedge Y^{\infty}~=~
\Sigma S(LV)^+\wedge V^{\infty} \wedge Y^{\infty} \to
LV^{\infty} \wedge V^{\infty} \wedge D_2(f)^{\infty}$$
the compactification Umkehr of the $\ZZ_2$-equivariant embedding
$$e_1~:~LV \times V \times D_2(f) \emb (LV\backslash \{0\}) \times V \times Y$$
of \ref{doubleHopf}.\\
{\rm (ii)} The geometric Hopf invariant
$h_V(F)$ is given up to natural $\ZZ_2$-equivariant homotopy by
$$\begin{array}{l}
h_V(F)~=~(1\wedge A)h_V(F)_Y~=~(1 \wedge Ai)F_1~:\\[1ex]
\Sigma S(LV)^+ \wedge V^{\infty} \wedge Y^{\infty}
\xymatrix@C+10pt{\ar[r]^-{\di{h_V(F)_Y}}&} \\[1ex]
\hskip50pt LV^{\infty} \wedge V^{\infty} \wedge (X\times_YX)^{\infty}
\xymatrix{\ar[r]^-{\di{1 \wedge A}}&}
LV^{\infty} \wedge V^{\infty} \wedge X^{\infty} \wedge X^{\infty}
\end{array}
$$
with $A:X \times_YX \subset X \times X$ the assembly.
\end{theorem}
\begin{proof} (i) We shall use the positive
property of the embedding $e=(g,f)$ to construct a
$\ZZ_2$-equivariant isotopy
$$E~:~e_3~\simeq~e_2~:~LV \times V \times X \to LV \times V \times Y$$
which is constant on $\{0\} \times V \times X$. This
will give a $\ZZ_2$-equivariant homotopy of the Umkehr maps
$$F_3~\simeq~F_2~:~
LV^{\infty} \wedge V^{\infty} \wedge Y^{\infty} \to LV^{\infty}
\wedge V^{\infty} \wedge X^{\infty}$$
which is constant on $\{0\}^+ \wedge V^{\infty} \wedge Y^{\infty}$,
so that $\delta(F_3,F_2) \simeq \{*\}$ by Proposition \ref{main} (i).\\
\indent In the first instance, note that for any $\epsilon >0$ and
any $\ZZ_2$-equivariant embedding $h:LV \times V \times X \emb LV
\times V \times Y$ there is defined a $\ZZ_2$-equivariant  isotopy
$$I \times LV \times V \times X \to LV \times V \times Y ~;~
(t,u,v,x) \mapsto h( \di{(1-t+t\epsilon) u \over 1+t\Vert u
\Vert},v,x)$$ from $h$ to the $\ZZ_2$-equivariant embedding
$$\overline{h}~:~LV \times V \times X \emb LV \times V \times Y ~;~
(u,v,x) \mapsto h( \di{\epsilon u \over 1+\Vert u \Vert},v,x)$$
which only involves the restriction
$$h\vert~:~D_{\epsilon}(0) \times V \times X \emb LV \times V \times Y~.$$
It therefore suffices to construct a $\ZZ_2$-equivariant  isotopy
$$e_3\vert~ \simeq~e_2\vert~:~D_{\epsilon}(0) \times V \times X \emb LV \times V \times Y$$
which is constant on $\{0\} \times V \times X$.  For sufficiently small
$\epsilon>0$ and $u \in D_{\epsilon}(0)$ we have
$$\begin{array}{l}
e_2(u,v,x)~=~((g(u+v,x)-g(-u+v,x))/2,(g(u+v,x)+g(-u+v,x))/2,f(x))\\[1ex]
\hskip150pt \approx~(dg_x(v)(u),g(v,x),f(x))~,\\[1ex]
e_3(u,v,x)~=~(u,g(v,x),f(x))
\end{array}$$
and
$$\begin{array}{l}
I \times D_{\epsilon}(0) \times V \times X \emb LV \times V \times Y~;\\[1ex]
(t,u,v,x) \mapsto (1-t)e_3(u,v,x)+te_2(u,v,x)~\approx~
((1-t)u+tdg_x(v)(u),g(v,x),f(x))
\end{array}$$
defines an isotopy $e_3\vert \simeq e_2\vert$
which is constant on $\{0\} \times V \times X$. The positive
condition on $e$ is used to ensure that the linear map
$$V \to V~;~u \mapsto tu+(1-t)dg_x(v)(u)$$
is an automorphism for each $t \in [0,1]$, exactly as in the
proof of Lemma \ref{contract}.\\
(ii) Immediate from (i) and \ref{doubleHopf}.
\hfill\qed\end{proof}

\begin{definition} {\rm Given maps $f_1:X_1 \to Y$, $f_2:X_2 \to Y$ write the
compactification of the  fibrewise product $X_1\times_YX_2$ of \ref{double1} as
$$(X_1 \times_Y X_2)^{\infty}~=~X_1^{\infty} \wedge_Y X_2^{\infty}~.$$}
\hfill\qed
\end{definition}

\begin{example} {\rm For compact $X_1,X_2$
the cartesian product of vector bundles $\xi_1:X_1 \to BO(m_1)$, $\xi_2:X_2 \to BO(m_2)$ is a vector bundle $\xi_1 \times \xi_2:X_1 \times X_2 \to BO(m_1+m_2)$
with total and Thom spaces
$$E(\xi_1 \times \xi_2)~=~E(\xi_1) \times E(\xi_2)~,~
T(\xi_1 \times \xi_2)~=~T(\xi_1) \wedge T(\xi_2)~.$$
For maps $f_1:X_1 \to Y$, $f_2:X_2 \to Y$
the restriction of $\xi_1 \times \xi_2$ to  $X_1 \times_YX_2\subseteq X_1 \times X_2$
$$\xi_1 \times_Y \xi_2~=~(\xi_1 \times \xi_2)\vert_{X_1\times_YX_2}~:~X_1\times_YX_2 \to BO(m_1+m_2)$$
is such that
$$E(\xi_1 \times_Y \xi_2)~=~E(\xi_1) \times_Y E(\xi_2)~,~
T(\xi_1 \times_Y \xi_2)~=~T(\xi_1)\wedge_YT(\xi_2)~.$$}
\hfill\qed
\end{example}

\begin{proposition} \label{immersions}
Given an immersion $f:M^m \imm N^n$ let
$\nu_f:M \to BO(n-m)$ be the normal $(n-m)$-plane bundle, so that $f$ extends to a codimension 0 immersion
$$f'~:~M'~=~E(\nu_f) \imm N~.$$
Assume that $f$ has no triple points, and only transverse double points,
so that the ordered double point space $M''=D_2(f)$ is a closed
$(2m-n)$-dimensional manifold.\\
{\rm (i)} The immersion
$$f\times_Nf~=~f \sqcup f''~:~M \times_N M~=~M \sqcup M'' \imm N$$
has normal bundle $\nu_f \sqcup \nu_{f''}$, with
$$\begin{array}{l}
f''~:~M'' \imm N~;~(x,y) \mapsto f(x)=f(y)~,\\[1ex]
\nu_{f''}~=~(\nu_f \times \nu_f)\vert_{M''}~:~M'' \to BO(2(n-m))~.
\end{array}$$
The ordered double point space of $f'$ is the $n$-dimensional manifold
$$D_2(f')^n~=~E(\nu_{f''})~=~E((\nu_f\times_N\nu_f)\vert_{M''})$$
and
$$T(\nu_f)\wedge_NT(\nu_f)~=~T(\nu_f \sqcup \nu_{f''})~=~
T(\nu_f) \vee D_2(f')^{\infty}~.$$
{\rm (ii)} For some finite-dimensional inner product space $V$
the codimension 0 immersion
$$1 \times f'~:~V \times M' \imm V \times N~;~(v,x) \mapsto (v,f'(x))$$
is regular homotopic to a positive embedding of $f'$
$$e~=~(g,f')~:~ V \times M' \emb V \times N~;~(v,x) \mapsto (g(v,x),f'(x))$$
with Umkehr map
$$F~:~V^{\infty} \wedge N^{\infty} \to V^{\infty} \wedge (M')^{\infty}~=~
V^{\infty} \wedge T(\nu_f)$$
and differentials of the form
$$\begin{array}{l}
de(v,x)~=~\begin{pmatrix} dg_x(v) & \alpha(v,x) \\ 0 & df(x) \end{pmatrix}~:~
\tau_{(v,x)}(V \times E(\nu_f))~=~V \oplus \tau_x(E(\nu_f)) \\[2ex]
\hskip50pt \to\tau_{(g(v,x),f'(x))}(V \times N)~=~V \oplus \tau_{f'(x)}(N)~(g_x(v)=g(v,x))
\end{array}$$
with $dg_x(v):V \to V$ positive linear automorphisms.
The fibrewise geometric Hopf invariant of $F$ is given up to natural
$\ZZ_2$-equivariant homotopy by
$$\begin{array}{l}
h_V(F)_N~=~(1 \wedge i)F_1~:~
\Sigma S(LV)^+ \wedge V^{\infty}\wedge N^\infty~=~(LV \backslash \{0\})^{\infty} \wedge V^{\infty}\wedge N^\infty \\[1ex]
\hskip50pt
\to LV^{\infty} \wedge V^{\infty} \wedge T(\nu_{f''}) \to LV^{\infty} \wedge V^{\infty} \wedge T(\nu_f\times_N\nu_f)~.
\end{array}$$
with
$$i~=~{\it inclusion}~:~T(\nu_{f''})~=~T((\nu_f \times_N \nu_f)\vert_{M''})
\to T(\nu_f \times_N \nu_f)~=~T(\nu_f) \wedge_N T(\nu_f)~.$$
The geometric Hopf invariant of $F$ is
$$h_V(F)~=~(1 \wedge Ai)F_1~:~\Sigma S(LV)^+ \wedge V^{\infty} \wedge N^{\infty} \to
LV^{\infty} \wedge V^{\infty} \wedge T(\nu_f)\wedge T(\nu_f)$$
with
$$A~=~{\it assembly}~:~T(\nu_f \times_N \nu_f)=T(\nu_f)\wedge_NT(\nu_f)
\to T(\nu_f \times \nu_f)~=~T(\nu_f)\wedge T(\nu_f)~.$$
\end{proposition}
\begin{proof} (i) By construction.\\
(ii) The double point space $D_2(f) \subset M \times M$ is a compact manifold with free involution,
and we can choose a smooth map $g_0 : M \to V$ for some $V$ of large
dimension  such that
$$g_0(x)~=~-g_0(y)\in S(V)~~((x,y) \in D_2(f))~.$$
We can take $g_0$ to be zero outside a neighbourhood of $D_2(f) \subset M \times M$,
so that there is defined an embedding
$$(g_0,f)~ ~: M \emb V\times N~;~x \mapsto (g_0(x),f(x))~.$$
Define
$$g~:~V \times M \to V~;~(v,x) \mapsto g(v,x)~=~g_0(x)+v~.$$
and let $D(V)=\{v \in V \st \Vert v \Vert \leqslant 1\}$ be the unit ball in $V$.
The map
$$(g,f)~:~V \times M \to V \times N~;~(v,x) \mapsto (g(v,x),f(x)) $$
 restricts to an embedding of $D(V)\times M \emb V \times N$, by the
following argument: the derivative has maximal rank  everywhere, and if
$$(g(v,x),f(x))~=~(g(w,y),f(y)) \in V \times N$$
then  for $(x,y) \in D_2(f)$ with $g(v,x)=g(w,y)$,
$$v-w=g_0(y)-g_0(x)~,~g_0(x)~=~-g_0(y)~,~\Vert v-w \Vert~=~2$$
and $(v,x)=(w,y)$. Composing with  the open embedding
$$j~:~V \to D(V)~;~ v \mapsto \dfrac{v}{1+\Vert v \Vert}$$
gives an embedding
$$e~:~V \times M \emb V \times N~;~(v,x) \mapsto (g(j(v),x),f(x))$$
which extends to a positive open embedding
$$e'~:~V \times M' \emb V \times N~;~(v,x) \mapsto (g'(v,x),f'(x))~.$$
Now apply Theorem \ref{doubleHopfpositive}.\\
\hfill\qed\end{proof}

\begin{example}\label{ad}
{\rm (i) Let $M=N=\{*\}$, and $V=\R^k$ with the standard inner product
$$\langle (u_1,u_2,\dots,u_k),(v_1,v_2,\dots,v_k) \rangle~=~
\sum\limits^k_{i=1}u_iv_i \in \R~.$$
For any
$\lambda_1,\lambda_2,\dots,\lambda_k> 0$ and $c=(c_1,c_2,\dots,c_k) \in V$ the embedding
$$g~:~V \emb V~;~v~=~(v_1,v_2,\dots,v_k) \mapsto
c+\di{(\lambda_1v_1,\lambda_2v_2,\dots,\lambda_kv_k)\over 1+\Vert v \Vert}$$
is such that $e=(g,f):V \times M \emb V \times N$ is a
positive embedding of the unique map $f:M \to N$, with
$$\begin{array}{l}
dg_*(v)~:~V \to V~;~u~=~(u_1,u_2,\dots,u_k) \mapsto
\di{(\lambda_1u_1,\lambda_2u_2,\dots,\lambda_ku_k)\over (1+\Vert v \Vert)^2}~,\\[1ex]
\langle dg_*(v)(u),u \rangle ~=~
\di{\sum\limits^k_{i=1}\lambda_iu_i^2\over (1+\Vert v \Vert)^2}>0~~
(u \neq 0 \in V)~,\\[3ex]
e(V \times M) \subseteq D_{\lambda}(c) \times N
\end{array}$$
where
$$D_{\lambda}(c)~=~
\{u=(u_1,u_2,\dots,u_k) \in V\,\vert\, \vert u_i-c_i \vert < \lambda_i,
1 \leqslant i \leqslant k\} \subset V~.$$
(ii) A finite cover $f:\widetilde{K} \to K$ admits positive embeddings
$e=(g,f):V \times \widetilde{K} \emb V \times K$ (cf. Example \ref{fincover}).
For $k> 2\,{\rm dim}(K)$ there exists an embedding
$d:\widetilde{K} \emb V=\R^k$, and for each $x \in \widetilde{K}$ it is
possible to make a continuous choice of $\lambda_1(x),\lambda_2(x),
\dots,\lambda_k(x)>0$ with
$$D_{\lambda(x)}(d(x)) \cap D_{\lambda(y)}(d(y))~=~\emptyset~
(x\neq y \in \widetilde{K},f(x)=f(y) \in K)~.$$
The map
$$g~:~V \times \widetilde{K} \to V~;~(v,x) \mapsto d(x)+ \di{(\lambda_1(x)v_1,
\lambda_2(x)v_2,\dots,\lambda_k(x)v_k) \over 1+\Vert v \Vert}$$
determines a positive embedding $e=(g,f)$ with
$$e(V \times \{x\})  \subseteq D_{\lambda(x)}(d(x)) \times \{f(x)\}~~
    (x \in \widetilde{K})~.$$
The (stable) homotopy class of the compactification Umkehr map $F:V^{\infty} \wedge
K^+ \to V^{\infty} \wedge \widetilde{K}^+$ of $e$
depends only on the homotopy class of the finite covering $f$.
\hfill\qed}
\end{example}

\begin{example} {\rm
(i) Given a map $c:X \to O(V)$ define an embedding of $f=1:X \to X$
$$e~=~(g,1)~:~V \times X \emb V \times X~;~(v,x) \mapsto (c(x)(v),x)$$
with adjunction Umkehr map
$$F~:~V^{\infty} \wedge X^+ \to V^{\infty} \wedge
    X^+~;~(v,x) \mapsto (c(x)^{-1}(v),x)~.$$
The $\ZZ_2$-equivariant embeddings $e_2,e_3:LV \times V \times X
\emb LV \times V \times X$ in Theorem \ref{doubleHopf} are given by
$$\begin{array}{l}
e_2(u,v,x) ~=~(u,g(v,x),f(x))~=~(u,c(x)(v),x)~,\\[1ex]
e_3(u,v,x) ~=~((g(u+v,x)-g(-u+v,x))/2,(g(u+v,x)+g(-u+v,x))/2,f(x))\\[1ex]
\hphantom{e_3(u,v,x)~}=~(c(x)(u),c(x)(v),x)~,
\end{array}$$
with compactification Umkehr maps
$$F_2~=~1 \wedge F~,~F_3~=~F \wedge F~:~LV^{\infty} \wedge V^{\infty} \wedge
X^+\to LV^{\infty} \wedge V^{\infty} \wedge X^+ $$
such that
$$h_V(F)~\simeq~(1\wedge \Delta_{X^+})\delta(F_2,F_3)~:~
\Sigma S(LV)^+ \wedge V^{\infty} \wedge
X^+\to LV^{\infty} \wedge V^{\infty} \wedge X^+\wedge X^+$$
(cf. Proposition \ref{twist}). Let $p:X^+ \to S^0$ be the projection, so that
$$pF~=~F_{c^{-1}}~:~LV^{\infty} \wedge X^+ \to LV^{\infty}$$
is the adjoint map of $c^{-1}$, and
$$-\theta(c)~=~\delta(1,F_{c^{-1}})\in \{\Sigma S(LV)^+ \wedge X^+;LV^{\infty}\}_{\ZZ_2}~=~
\omega^{-1}_{\ZZ_2}(S(LV)\times X;-\epsilon_{LV})$$
with
$$h_V(F_{c^{-1}})~=~(F_{c^{-1}}\wedge -\theta(c))\wedge \Delta_{X^+}
\in \{\Sigma S(LV)^+ \wedge X^+;LV^{\infty}\}_{\ZZ_2}$$
as in Proposition \ref{twist}.\\
(ii) As in Definition \ref{positive} let $O^+(V) \subset O(V)$ be
the subgroup of the positive orthogonal automorphisms $a:V \to V$.
If the map $c:X \to O(V)$ in (i) is such that $c(X) \subseteq
O^+(V)$ then $h_V(F) \simeq *$, by Theorem
\ref{doubleHopfpositive}. (Alternatively, note that the isotopy $e
\simeq {\rm id}$ defined by
    $$V \times X \times I \to V \times X~;~(v,x,t) \mapsto
    (tv+(1-t)c(x)(v),x)~,$$
induces a homotopy $F \simeq {\rm id}$.)\hfill\qed}
\end{example}

\begin{example} {\rm
(i) Given a manifold $N$, a space $X$, an inner product space $V$
and a map $\rho:V^{\infty} \wedge N^\infty \to V^{\infty} \wedge
X^+$ it is possible to make $\rho$ transverse at $\{0\}
\times X \subset V^{\infty} \wedge X^+$, so that there
exist a codimension 0 immersion $f:M=\rho^{-1}(\{0\} \times X)
\imm N$ with an embedding $e=(g,f):V \times M \emb V \times N$ and
a map $h=\rho\vert:M \to X$ so that up to homotopy
$$\rho~=~(1 \wedge h)F~:~
V^{\infty} \wedge N^\infty \xymatrix{\ar[r]^-{\di{F}}&}
V^{\infty} \wedge M^+ \xymatrix{\ar[r]^-{\di{1\wedge h}}&}
V^{\infty} \wedge X^+$$
with $F$ the adjunction Umkehr of $e$.
    (In general, it is not possible to choose $e$ to be a
    positive embedding -- see (ii) below for an explicit example).
The framing of $0 \times f$ determined by $g$ differs from the
canonical framing by a map $c:M \to O(V)$ such that
$$g(v,x)~=~g_0(c(x)(v),x) \in V~~(v \in V,x \in M)~.$$
The embeddings of $f$
$$e_0~=~(g_0,f)~,~ e~=~(g,f)~:~V \times M \emb V \times N$$
are such that $e=e_0(c,1)$ with
$$(c,1)~:~V \times M \to V \times M~;~(v,x) \mapsto
    (c(x)(v),x)~,$$
so that the adjunction Umkehr maps $F_0,F:V^{\infty} \wedge N^\infty \to
V^{\infty} \wedge M^+$ differ by the adjunction Umkehr of $(c,1)^{-1}$
$$C~:~V^{\infty} \wedge M^+ \to V^{\infty} \wedge
    M^+~;~(v,x) \mapsto (c(x)(v),x)$$
(cf. Proposition \ref{twist}) with
    $$F~:~V^{\infty} \wedge N^\infty
    \xymatrix{\ar[r]^-{\di{F_0}}&} V^{\infty} \wedge M^+
    \xymatrix{\ar[r]^-{\di{C}}&} V^{\infty} \wedge M^+~.$$
By Theorem \ref{doubleHopfpositive} and the composition formula of
Proposition \ref{hopf7} (v) the geometric Hopf invariant of $F$ is
given by
$$\begin{array}{ll}
h_V(F)&=~h_V(CF_0)~=~h_V(C)(1\wedge F_0)+(\kappa^{-1}_V\wedge 1)(C \wedge C)(\kappa_V\wedge 1)h_V(F_0)~:\\[1ex]
&\hskip50pt
\Sigma S(LV)^+ \wedge V^{\infty} \wedge N^\infty \to
LV^{\infty} \wedge V^{\infty} \wedge M^+ \wedge M^+~.
\end{array}$$
(ii) If $N=X=\{{\rm pt.}\}$, $V=\R$ then
$$\rho~=~-1~:~V^{\infty} \wedge N^\infty~=~S^1 \to V^{\infty} \wedge X^+~=~S^1~;~z \mapsto z^{-1}$$
is a map of degree $-1$ with $M=\rho^{-1}(\{0\}\times X)=\{*\}$,
such that the corresponding  embedding
$$e~=~(g,f)~:~V \times M~=~\R \emb V \times N~=~\R~;~v \mapsto
   -v$$
of the codimension 0 immersion $f:M \imm N$ is not positive. The
canonical (positive) embedding of $f$
$$e_0~=~(g_0,f)~:~V \times M~=~\R \emb V \times N~=~\R~;~v \mapsto
   v$$
has adjunction Umkehr
$$F_0~=~1~:~V^{\infty}\wedge N^\infty~=~S^1 \to V^{\infty}\wedge M^+~=~S^1~,$$
and $e=(g,f)$ has adjunction Umkehr
$$F~=~F_0C~=~\rho~=~-1~:~V^{\infty}\wedge N^\infty~=~S^1 \to
    V^{\infty}\wedge M^+~=~S^1~,$$
with $C:S^1 \to S^1$ the adjoint of $c=-1:M=\{*\} \to O(\R)=\{\pm
1\}$.  The geometric Hopf invariant of $C$ is
$h_{\R}(C)=1 \in \ZZ$, by Example \ref{twist2}.\\
\hfill\qed}
\end{example}

\begin{proposition}~ \label{RHopf}
Let $e=(g,f):\R \times M \emb \R \times N$ be a
positively embedded map, with adjunction Umkehr map
$$F~:~\R^{\infty}\wedge N^\infty~=~\Sigma N^\infty \to
\R^{\infty}\wedge M^+~=~\Sigma M^+~.$$
{\rm (i)} The transposition $\ZZ_2$-action on the set
$$E~=~\R \times \R \times D_2(f)$$
is free, with a decomposition
$$E~=~E_+\sqcup E_-$$
defined by
$$\begin{array}{l}
E_+~=~\{(v_1,v_2,x_1,x_2) \in \R \times \R \times M \times M\,\vert\,\\[1ex]
\hskip50pt g(v_1,x_1)>g(v_2,x_2) \in \R,x_1 \neq x_2 \in M,f(x_1)=f(x_2) \in N\}~,\\[1ex]
E_-~=~\{(v_1,v_2,x_1,x_2) \in \R \times \R \times M \times M\,\vert\,\\[1ex]
\hskip50pt g(v_1,x_1)<g(v_2,x_2) \in \R,x_1 \neq x_2 \in M,f(x_1)=f(x_2) \in N\}~.
\end{array}$$
{\rm (ii)} The geometric Hopf invariant of $F$ is given up to a natural homotopy by
$$\begin{array}{l}
h_{\R}(F)_+~=~ i_+F_+~:~(L\R\backslash \{0\})_+^{\infty}
\wedge \R^{\infty} \wedge N^\infty~=~\Sigma^2N^\infty\\[1ex]
\hskip100pt \to
L\R^{\infty} \wedge \R^{\infty}\wedge M^+ \wedge M^+~=~\Sigma^2(M^+\wedge M^+)
\end{array}$$
(using Terminology \ref{plus}) with
$$F_+~:~(L\R\backslash \{0\})_+^{\infty}\wedge \R^{\infty} \wedge N^\infty \to
E_+^{\infty}$$
the adjunction Umkehr of the embedding
$$\begin{array}{l}
e_+~:~E_+ \emb (L\R\backslash \{0\})_+ \times \R \times N~;~
(v_1,v_2,x_1,x_2) \mapsto (w_1,w_2,y)\\[1ex]
\hskip20pt (w_1=(g(v_1,x_1)-g(v_2,x_2))/2 \in(L\R\backslash \{0\})_+~,\\[1ex]
\hskip20pt \hphantom{(}w_2=(g(v_1,x_1)+g(v_2,x_2))/2\in \R~,~
y=f(x_1)=f(x_2)\in N)
\end{array}$$
and
$$i_+~=~{\it inclusion}~:~E_+ \to \R \times \R \times M\times M~.$$
\end{proposition}
\begin{proof} This is the special case $V=\R$ of
Theorem \ref{doubleHopfpositive}.
\hfill\qed\end{proof}

\begin{example} \label{pinch}
{\rm Let $f:M=\{1,2\}\to N=\{*\}$ be the unique map, with ordered
double point set
$$D_2(f)~=~(M\times M)\backslash \Delta(M)~=~\{(1,2),(2,1)\}~.$$
Define a positive embedding $(e:\R \times M \emb \R \times N,f:M
\to N)$ of $f$ by
$$e~:~\R \times \{1,2\} \emb \R~;~\begin{cases} (v,1) \mapsto e^v \\
(v,2) \mapsto -e^{-v}~.
\end{cases}$$
The compactification Umkehr of $e$ is a sum map
$$F~=~\nabla~=~h_1+h_2~:~
(\R \times N)^{\infty}~=~S^1 \to (\R \times M)^{\infty}~=~S^1\vee S^1$$
with
$$h_i~:~S^1 \to S^1 \vee S^1~;~t \mapsto t_i~~(i=1,2)~,$$
and by Proposition \ref{RHopf} the geometric Hopf invariant of $F$ is
determined up to $\ZZ_2$-equivariant homotopy by the map
    $$\begin{array}{l}
h_{\R}(F)\vert~=~h_1\wedge h_2~:\\[2ex]
    \Sigma^2N^\infty~=~S^2~=~S^1 \wedge S^1 \to \Sigma^2(M^+\wedge M^+)~=~
    (S^1 \vee S^1) \wedge (S^1 \vee S^1) ~.
    \end{array}$$
\hfill\qed}
\end{example}

\section{Finite covers}\label{finiteum}

In this section $K$ is a finite $CW$ complex.  Every finite covering
$f:\widetilde{K} \to K$ admits a positive embedding $e=(g,f):V \times
\widetilde{K} \emb V \times K$ (\ref{ad} (ii)), with a compactification
Umkehr map $F:V^{\infty} \wedge K^+ \to V^{\infty} \wedge \widetilde{K}^+$
(\ref{umkehradjunct}).

For a finite covering $f:\widetilde{K} \to K$ of degree $d$
define a commutative square of finite coverings
$$\xymatrix{\widetilde{K}' \ar[r]^-{\di{\widetilde{f}'}}
\ar[d]_-{\di{\widetilde{f}}}  & \widetilde{K}
\ar[d]^-{\di{f}} \\
K'\ar[r]^-{\di{f'}} & K}$$
with
$$\begin{array}{l}
\widetilde{K}'~=~D_2(f)~=~\{(x_1,x_2) \in \widetilde{K} \times \widetilde{K}\,\vert\,
x_1 \neq x_2 \in \widetilde{K}\,,\,f(x_1)=f(x_2) \in K\}~,\\[1ex]
K'~=~D_2[f]~=~\widetilde{K}'/\{(x_1,x_2) \sim (x_2,x_1)\}
\end{array}$$
so that the projections
$$\begin{array}{l}
\widetilde{f}~:~\widetilde{K}' \to K'~;~(x_1,x_2) \mapsto [x_1,x_2]~,\\[1ex]
f'~:~K'\to K~;~[x_1,x_2] \mapsto f(x_1)=f(x_2)~,\\[1ex]
\widetilde{f}'~:~\widetilde{K}' \to \widetilde{K}~;~(x_1,x_2) \mapsto x_1
\end{array}$$
are finite coverings of degree $2$, $d(d-1)/2$, $d-1$ respectively. The
composite
$$f\widetilde{f}'~=~f'\widetilde{f}~:~\widetilde{K}' \to K~;~
(x_1,x_2) \mapsto f(x_1)=f(x_2)$$
is a finite covering of degree $d(d-1)$.

\begin{proposition}~ \label{hopfhopf}
The geometric Hopf invariant of the compactification Umkehr map
$F:V^{\infty} \wedge K^+ \to V^{\infty} \wedge
\widetilde{K}^+$ of a positively embedded finite covering
$e=(g,f):V \times \widetilde{K} \emb V \times K$
is given up to natural $\ZZ_2$-equivariant homotopy by
$$\begin{array}{l}
h_V(F)~=~(1 \wedge i)F_1~:~
\Sigma S(LV)^+ \wedge V^{\infty}\wedge K^+\\[1ex]
\hskip100pt \to  LV^{\infty} \wedge V^{\infty} \wedge (\widetilde{K}')^+
\to LV^{\infty} \wedge V^{\infty} \wedge \widetilde{K}^+\wedge \widetilde{K}^+
\end{array}$$
with
$$F_1~:~\Sigma S(LV)^+\wedge V^{\infty}\wedge
K^+ ~=~
(LV\backslash \{0\})^+ \wedge V^{\infty}\wedge K^+
\to LV^{\infty}\wedge V^{\infty} \wedge (\widetilde{K}')^+$$
the Umkehr map of the $\ZZ_2$-equivariant embedding
$$\begin{array}{l}
e_1~:~LV \times V \times \widetilde{K}' \emb
(LV\backslash \{0\}) \times V \times K~;~(v_1,v_2,x_1,x_2) \mapsto \\[1ex]
((g(v_1+v_2,x_1)-g(-v_1+v_2,x_2))/2,(g(v_1+v_2,x_1)+g(-v_1+v_2,x_2))/2),y)\\[1ex]
\hskip200pt (y=f(x_1)=f(x_2))
\end{array}$$
and
$$i~=~{\it inclusion}~:~\widetilde{K}'~=~D_2(f) \to \widetilde{K}\times
\widetilde{K}~.$$
\end{proposition}
\begin{proof}
A direct application of Theorem \ref{doubleHopfpositive}.
\hfill\qed\end{proof}

\begin{example} \label{double7}
{\rm For any $d \geqslant 1$ consider the finite covering of degree $d$
$$f~:~\widetilde{K}~=~\{1,2,\dots,d\} \to K~=~\{*\}$$
with ordered double point set
$$\widetilde{K}'~=~D_2(f)~=~
\{(i,j)\,\vert\,1 \leqslant  i,j \leqslant  d, i \neq j\}$$
consisting of $d(d-1)$ points. For any positive embedding
$e=(g,f):\R \times \widetilde{K} \emb \R \times K=\R$
the Umkehr is given up to homotopy by the $d$-fold sum map
$$\begin{array}{ll}
F~=~h_1+h_2+\dots+h_d~:&\R^+~=~S^1~=~I/(0=1)  \to \R^+ \wedge \widetilde{K}^+~
=~\bigvee\limits_d{S^1}~;\\[1ex]
&t \mapsto (dt-i+1)_i ~{\rm if}~
(i-1)/d \leqslant  t \leqslant  i/d
\end{array}$$
where
$$h_i~:~S^1 \to \bigvee\limits_d{S^1}~;~t \mapsto t_i~(1 \leqslant i \leqslant d)~.$$
The geometric Hopf invariant of $F$ (\ref{RHopf}) is given up to
 homotopy by
$$h_{\R}(F)~=~\sum\limits_{1 \leqslant  i < j \leqslant  d}h_i\wedge h_j ~:~
S^1\wedge S^1 \to
(\bigvee\limits_d{S^1}) \wedge (\bigvee\limits_d{S^1})$$
The geometric Hopf invariant of the composite
$$d~:~S^1 \xymatrix@C+10pt{\ar[r]^-{\di{F}}&} \bigvee\limits_d{S^1}
\xymatrix@C+40pt{\ar[r]^-{\di{1\vee 1 \vee \dots \vee 1}}&} S^1$$
is thus
$$\begin{array}{ll}
h_{\R}(d)&=~\big((1\vee 1 \vee \dots \vee 1) \wedge
(1\vee 1 \vee \dots \vee 1)\big)h_{\R}(F)\\[1ex]
&=~d(d-1)/2~:~S^2 \to S^2~.\end{array}$$
\hfill\qed}
\end{example}

A double cover $f:\widetilde{K} \to K$ of a finite $CW$ complex $K$ is
classified by a map $c:K \to P(V)$ to the projective space
$P(V)=S(LV)/\ZZ_2$ of a inner product space $V$, for
sufficiently large ${\rm dim}(V)$.

\begin{proposition}~
For a finite cover $f:\widetilde{K} \to K$ of degree $d$ the quadratic
construction on the compactification Umkehr map $F:V^{\infty}\wedge K^+ \to
V^{\infty} \wedge \widetilde{K}^+$ of a positive embedding
$e=(g,f):V \times \widetilde{K} \emb V \times K$ is given by
$$\psi_V(F)~:~K^+ \xymatrix{\ar[r]^-{\di{F'}}&}
K^{\prime +} \xymatrix{\ar[r]^-{\di{c}}&}
S(LV)^+ \wedge_{\ZZ_2} (\widetilde{K} \times \widetilde{K})^+$$
with $F'$ the Umkehr map for the finite cover of degree $d(d-1)/2$
$$\begin{array}{l}
f'~:~K'~=~\{(x,y) \in \widetilde{K} \times \widetilde{K}\,\vert\,
x \neq y \in \widetilde{K},f(x)=f(y) \in K\}/\ZZ_2 \to K~;\\[1ex]
\hskip100pt [x,y] \mapsto f(x)=f(y)
\end{array}$$
and
$$c~:~K' \to S(LV) \times _{\ZZ_2} (\widetilde{K} \times \widetilde{K})~;~
[x,y] \mapsto \bigg[\di{ g(0,x) - g(0,y) \over \Vert g(0,x) -
    g(0,y) \Vert },x,y \bigg]~.$$
For $d=2$ $F'$ is a homeomorphism
$$F'~:~K^+ \to (K')^+~;~f(x) \mapsto [x,Tx]$$
with $T:\widetilde{K} \to \widetilde{K}$ the covering translation, and
$$K \to P(V)~;~f(x) \mapsto \di{ g(0,x) - g(0,Tx) \over \Vert g(0,x) -
    g(0,Tx) \Vert }$$
is a classifying map for the double cover $f:\widetilde{K} \to K$.\hfill$\qed$
\end{proposition}

\begin{proposition}~ \label{doublecover}
Let $c:K \to P(V)$ be a map, classifying the double cover of $K$
$$\widetilde{K}~=~c^*S(LV)~=~
\{(v,x) \in S(LV)\times K\,\vert\,[v]=c(x)\in P(V)\}$$
with covering projection
$$f~:~\widetilde{K} \to K~;~(v,x) \mapsto x$$
and covering translation
$$T~:~\widetilde{K} \to \widetilde{K}~;~(v,x) \mapsto (-v,x)~.$$
{\rm (i)} The map
$$g~:~V \times \widetilde{K} \to V~;~(u,(v,x)) \mapsto
v+\dfrac{u}{1+\Vert u \Vert}$$
determines a positive embedding
$$e~=~(g,f)~:~V \times \widetilde{K} \emb V \times K~;~
(u,(v,x)) \mapsto (v+\dfrac{u}{1+\Vert u \Vert},x)$$
with Umkehr map
$$\begin{array}{l}
F~:~V^{\infty} \wedge K^+ \to V^{\infty} \wedge \widetilde{K}^+~;\\[1ex]
(u,x) \mapsto
\begin{cases}
(\dfrac{u- v}{1-\Vert u - v \Vert},(v,x))
&{\it if}~(v,x) \in \widetilde{K},~\Vert u - v\Vert <1\\[2ex]
(\dfrac{u+v}{1-\Vert u + v \Vert},(-v,x))
&{\it if}~(v,x) \in \widetilde{K},~\Vert u + v\Vert <1\\[2ex]
*&{\it otherwise}~.
\end{cases}
\end{array}$$
{\rm (ii)} The restriction of $e$ is a $\ZZ_2$-equivariant embedding
$$e'~=~e\vert~:~LV \times \widetilde{K} \emb (LV \backslash \{0\}) \times K$$
with $\ZZ_2$-equivariant Umkehr map
$$F'~=~F\vert~:~(LV \backslash \{0\})^+ \wedge K^+~=~
\Sigma S(LV)^+ \wedge K^+ \to LV^{\infty} \wedge \widetilde{K}^+~.$$
{\rm (iii)} For any subspace $K_0 \subseteq K$ let
$$f_0~=~f\vert~:~\widetilde{K}_0~=~f^{-1}(K_0) \to K_0$$
be the pullback double cover of $K_0$. The maps $F,F'$ in {\rm (i)},
{\rm (ii)} induce maps
$$\begin{array}{l}
F~:~V^{\infty} \wedge K/K_0 \to V^{\infty} \wedge \widetilde{K}/\widetilde{K}_0~,\\[1ex]
F'~:~\Sigma S(LV)^+ \wedge K/K_0 \to LV^{\infty} \wedge \widetilde{K}/\widetilde{K}_0
\end{array}$$
with $F'$ $\ZZ_2$-equivariant, such that the geometric Hopf invariant of
$F$ is given up to natural $\ZZ_2$-homotopy by
$$\begin{array}{l}
\xymatrix@C+10pt{h_V(F)~:~
\Sigma S(LV)^+ \wedge V^{\infty} \wedge K/K_0
\ar[r]^-{\displaystyle{1\wedge F'}} &
LV^{\infty} \wedge V^{\infty} \wedge \widetilde{K}/\widetilde{K}_0}\\
\hspace*{100pt}
\xymatrix@C+60pt{
\ar[r]^-{\displaystyle{1 \wedge (1\wedge T)\Delta_{\widetilde{K}/\widetilde{K}_0}}}
&  LV^{\infty} \wedge V^{\infty} \wedge
\widetilde{K}/\widetilde{K}_0 \wedge \widetilde{K}/\widetilde{K}_0}~.
\end{array}$$
\end{proposition}
\begin{proof} (i)+(ii) By construction.\\
(iii) By Proposition \ref{hopfhopf} $h_V(F)$ is given up to natural
$\ZZ_2$-equivariant homotopy to
$$\begin{array}{l}
\xymatrix@C+10pt{h_V(F)~:~
\Sigma S(LV)^+ \wedge V^{\infty} \wedge K/K_0
\ar[r]^-{\displaystyle{F_1}} &
LV^{\infty} \wedge V^{\infty} \wedge D_2(f)/D_2(f_0)}\\[1ex]
\hspace*{100pt}
\xymatrix@C+10pt{
\ar[r]^-{\displaystyle{1 \wedge i}} & LV^{\infty} \wedge V^{\infty} \wedge
\widetilde{K}/\widetilde{K}_0 \wedge \widetilde{K}/\widetilde{K}_0}
\end{array}$$
with $F_1$ the Umkehr map of the $\ZZ_2$-equivariant open embedding
$$\begin{array}{l}
e_1~:~LV \times V \times D_2(f) \emb (LV \backslash \{0\}) \times V
\times K~;~(u,v,x_1,x_2) \mapsto\\[1ex]
((g(u+v,x_1)-g(-u+v,x_2))/2,(g(u+v,x_1)+g(-u+v,x_2))/2,y)\\[1ex]
\hskip150pt (y=f(x_1)=f(x_2) \in K)
\end{array}$$
and $i:D_2(f) \to \widetilde{K} \times \widetilde{K}$ the inclusion.
Use the $\ZZ_2$-equivariant homeomorphism
$$\widetilde{K} \to D_2(f)~;~(v,x) \mapsto ((v,x),(-v,x))$$
to identify
$$i~:~D_2(f)~=~\widetilde{K} \to \widetilde{K} \times \widetilde{K}~;~
(v,x)\mapsto ((v,x),(-v,x))$$
and to express $e_1$ as
$$\begin{array}{l}
e_1~:~LV \times V \times \widetilde{K} \emb
(LV\backslash \{0\}) \times V \times K~;~(u_1,u_2,v,x) \mapsto \\[1ex]
\hskip100pt
((g(u_1+u_2,v,x)-g(-u_1+u_2,-v,x))/2,\\[1ex]
\hskip100pt (g(u_1+u_2,v,x)+g(-u_1+u_2,-v,x))/2),f(x))~.
\end{array}$$
The $\ZZ_2$-equivariant embeddings $e_1,1\times e':
LV \times V \times \widetilde{K} \emb (LV\backslash \{0\}) \times V \times K$
are such that
$$e_1\vert~=~(1 \times e')\vert~:~LV \times \{0\} \times \widetilde{K}
\emb (LV\backslash \{0\}) \times V \times K~;~
(u,0,v,x) \mapsto (g(u,v,x),0,f(x))~,$$
so that by the $\ZZ_2$-equivariant tubular neighbourhood theorem
(applied in the special case $K=\{*\}$, $\widetilde{K}=S^0$)
they are related by a natural $\ZZ_2$-equivariant open isotopy
$e_1 \simeq (1\times e')$ inducing a natural $\ZZ_2$-equivariant
homotopy
$$F_1~\simeq~1 \wedge F'~:~
\Sigma S(LV)^+ \wedge V^{\infty} \wedge K/K_0 \to
LV^{\infty} \wedge V^{\infty} \wedge \widetilde{K}/\widetilde{K}_0
\wedge \widetilde{K}/\widetilde{K}_0~.$$
\hfill\qed\end{proof}

\begin{remark} {\rm (Segal \cite{segal2,segal3})
For any unpointed space $X$ let
$C(K;X)$ denote the set of isomorphism classes of pairs
$(f:\widetilde{K} \to K,g:\widetilde{K} \to X)$ with $f$ a finite cover
and $g \in [\widetilde{K},X]$.  The set $C(K;X)$ is an abelian semigroup with
respect to disjoint union of the spaces $\widetilde{K}$, which is a
contravariant homotopy-functor of $K$. The Umkehr construction defines
a natural transformation of semi-group valued functors $C(~;X) \to
\{~;X\}$ by
$$C(K;X) \to \{K^+;X^+\}~;~(f,g) \mapsto gF$$
which is a `group completion', i.e.  is universal among transformations
$T:C(~;X) \to F$ where $F$ is a representable abelian-group-valued
homotopy-functor and $T$ is a transformation of semigroup-valued
functors (Barratt, Priddy and Quillen).  As above, given a finite cover
$f:\widetilde{K} \to K$ of degree $d$ and a map $g:\widetilde{K} \to X$
there are defined a degree $d(d-1)/2$ cover $f':K' \to K$ and a degree 2
cover $\widetilde{f}':\widetilde{K}' \to K'$, with
$$\widetilde{K}'~=~\{(x_1,x_2) \in \widetilde{K} \times \widetilde{K}\,\vert\,
x_1 \neq x_2 \in \widetilde{K}\,,\,f(x_1)=f(x_2) \in K\}~.$$
Let $c:K' \to P(\R(\infty))$ be a classifying map for $\widetilde{f}'$,
with $\ZZ_2$-equivariant lift $\widetilde{c}:\widetilde{K}' \to S(\infty)$,
and let $g':K' \to S(\infty)\times_{\ZZ_2}(X \times X)$ be
the quotient of the $\ZZ_2$-equivariant map
$$\widetilde{g}'~:~\widetilde{K}' \to S(\infty) \times (X \times X)~;~
(x_1,x_2) \mapsto (\widetilde{c}(x_1,x_2),g(x_1),g(x_2))~.$$
The transformation of semigroup-valued functors
$$\theta^2~:~C(K;X) \to C(K;S(\infty)\times_{\ZZ_2}(X \times X))~;~
(f,g) \mapsto (f',g')$$
determines a transformation of abelian-group-valued functors
$$\theta^2~:~\{K^+;X^+\} \to \{K^+;(S(\infty)\times_{\ZZ_2}(X \times X))^+\}~.$$
In Example \ref{double11} below the corresponding stable map
$$\theta^2~:~X^+ \sto (S(\infty)\times_{\ZZ_2}(X \times X))^+$$
will be expressed in terms of the geometric Hopf invariant of the evaluation
stable map $\Omega^{\infty}\Sigma^{\infty}X^+\sto X^+$.\\
\hfil\qed}
\end{remark}

\section{Function spaces}

We shall now use the quadratic construction $\psi_{\R^k}(F)$ and the
double point theorem for finite covers to split off the quadratic
part of the function spaces $\Omega^k\Sigma^kX$ for $k=1,2,\dots,\infty$,
with $X$ a pointed (but not necessarily connected) space.

The geometric Hopf invariant defines a map of function spaces
$$h_V~:~{\rm map}(V^{\infty} \wedge X,V^{\infty} \wedge Y)
\to {\rm map}^{\ZZ_2}(\Sigma S(LV)^+ \wedge
V^{\infty} \wedge X,LV^{\infty} \wedge V^{\infty} \wedge Y\wedge Y)$$
for any pointed spaces $X,Y$ and inner product space $V$.

\begin{definition}~ \label{hopf12}
{\rm  The pointed space
$$Q_V(X)~=~{\rm map}(V^{\infty},V^{\infty} \wedge X)$$
is the function space of maps $\omega:V^{\infty} \to
V^{\infty}\wedge X$.}
\hfill\qed
\end{definition}

\begin{proposition}\label{eval}~
{\rm (i)} A pointed map $f:Y \to Q_V(X)$ is essentially the same as a pointed
map $F:V^{\infty} \wedge Y \to V^{\infty} \wedge X$, via the adjoint map
construction
$$F~=~{\rm adj}(f)~:~V^{\infty} \wedge Y \to V^{\infty} \wedge X~;~
(v,y) \mapsto f(y)(v)~.$$
The function
$$[Y,Q_V(X)] \to [V^{\infty} \wedge Y,V^{\infty} \wedge X]~;~f \mapsto F$$
is a bijection.\\
{\rm (ii)} The inclusion
$$i~:~X \to Q_V(X)~;~x \mapsto (v \mapsto (v,x))$$
has adjoint the identity map
$${\rm adj}(i)~=~1~:~V^{\infty} \wedge X \to V^{\infty} \wedge X~.$$
{\rm (iii)} The identity map $1_V:Q_V(X) \to Q_V(X)$
has adjoint the evaluation map
$$e_V~=~{\rm adj}(1_V)~:~V^{\infty} \wedge Q_V(X) \to V^{\infty}\wedge X~;~
(v,\omega) \mapsto \omega(v)~,$$
such that for $f,F$ as in {\rm (i)}
$$\xymatrix@C+10pt{F~:~V^{\infty} \wedge Y \ar[r]^-{\di{1 \wedge f}} &
V^{\infty} \wedge Q_V(X) \ar[r]^-{\di{e_V}} & V^{\infty} \wedge X~,}$$
and
$$\xymatrix@C+20pt{\psi_V(F)~:~Y \ar[r]^-{\di{f}} & Q_V(X)
\ar[r]^-{\di{\psi_V(e_V)}} & S(LV)^+ \wedge_{\ZZ_2} (X \wedge X)~.}$$
{\rm (iv)} For $V \neq \{0\}$ $Q_V(X)$ is an $H$-space via
$$Q_V(X) \times Q_V(X) \to Q_V(X)~;~(f,g) \mapsto f+g=(f \vee g)\nabla$$
with $\nabla:V^{\infty} \to V^{\infty} \vee V^{\infty}$ any sum map.
For ${\rm dim}(V) \geqslant 2$ $Q_V(X)$ is a homotopy commutative $H$-space.\\
{\rm (v)} For any pointed space $Y$ the quadratic construction defines a function
$$\psi_V~:~[Y,Q_V(X)]~=~[V^{\infty} \wedge Y,V^{\infty} \wedge X] \to
\{Y;S(LV)^+ \wedge_{\ZZ_2}(X \wedge X)\}~;~F \mapsto \psi_V(F)$$
which is induced by a stable map
$$\psi_V~:~Q_V(X) \sto S(LV)^+ \wedge_{\ZZ_2}(X \wedge X)~.$$
\hfill\qed
\end{proposition}

\begin{example} {\rm Let $V=\R^k$, with $k \geqslant 1$.\\
(i) If $X$ is connected then so is $Q_V(X)$, with
$$\pi_n(Q_V(X))~=~\pi_{n+k}(\Sigma^kX)~.$$
(ii) For $X=S^0$
$$Q_V(S^0)~=~{\rm map}(V^{\infty},V^{\infty})~=~\Omega^kS^k~,$$
with a bijection
$$\pi_0(Q_V(S^0)) \to \ZZ~;~(\omega:V^{\infty} \to V^{\infty}) \mapsto
{\rm degree}(\omega)~.$$
The quadratic construction is given by
$$\psi_V~:~\pi_0(Q_V(S^0))=\ZZ \to \pi^S_0(P(V)^+)=\ZZ~;~
d \mapsto \dfrac{d(d-1)}{2}~.$$
Proposition \ref{twist} gives a homotopy commutative diagram of stable maps
$$\xymatrix{O(V) \ar[r]^-{\di{\theta \wedge J}} \ar[d]_-{\di{J}}
& P(V)^+ \wedge Q_V(S^0)\ar[d]^-{\di{1\wedge e_V}}\\
Q_V(S^0) \ar[r]^-{\di{\psi_V}} & P(V)^+ \wedge S^0=P(V)^+}$$
given on the path-components by
$$\xymatrix{
\pi_0(O(V))=\{1,-1\} \ar[r]^-{\di{1}} \ar@{(->}[d]
& \pi_0(P(V)^+ \wedge Q_V(S^0))=\{1,-1\}
\ar@{(->}[d]\\
\pi_0(Q_V(S^0))=\ZZ \ar[r]^-{\di{\psi_V}} & \pi^S_0(P(V)^+)=\ZZ~.}$$
By Proposition \ref{localR} (ii) $\theta R \simeq {\rm id.}$, so
there is defined a homotopy commutative diagram of stable maps
$$\xymatrix@C+10pt{P(V)^+ \ar[r]^-{\di{1\wedge JR}} \ar[d]_-{\di{JR}}
& P(V)^+ \wedge Q_V(S^0)\ar[d]^-{\di{1\wedge e_V}}\\
Q_V(S^0) \ar[r]^-{\di{\psi_V}} & P(V)^+\wedge S^0=P(V)^+}$$
\hfill\qed
}
\end{example}

\begin{example}{\rm
(i) Let $K$ be a finite $CW$ complex with a map $c:K \to P(V)$,
classifying a double cover $f:\widetilde{K}=c^*S(LV) \to K$,
with Umkehr map $F:V^{\infty} \wedge K^+ \to V^{\infty} \wedge
\widetilde{K}^+$. The composite stable map
$$\xymatrix@C+20pt{K^+ \ar[r]^-{\di{{\rm adj}(F)}} &
Q_V(\widetilde{K}^+) \ar[r]^-{\di{\psi_V}} &
S(LV)^+ \wedge_{\ZZ_2} (\widetilde{K}^+ \wedge
\widetilde{K}^+)}$$
is represented by
$$K \to S(LV) \times_{\ZZ_2} (\widetilde{K} \times \widetilde{K})~;~
f(x) \mapsto [\widetilde{c}(x),x,Tx]$$
with $\widetilde{c}:\widetilde{K} \to S(LV)$ a $\ZZ_2$-equivariant
lift of $c$ and $T:\widetilde{K} \to \widetilde{K}$ the covering
translation (by Proposition \ref{quad3} (vii)).\\
(ii) In particular, (i) applies to the canonical double cover
$f:\widetilde{K}=S(LV) \to K=P(V)$ classified by $c=1:K \to P(V)$,
with Umkehr map
$$\begin{array}{l}
F~:~V^{\infty} \wedge P(V)^+ \to V^{\infty} \wedge S(LV)^+~;\\[1ex]
(u,[v]) \mapsto \begin{cases}
(\dfrac{u- v}{1-\Vert u - v \Vert},v)
&{\rm if}~v \in S(LV),~\Vert u - v\Vert <1\\[2ex]
(\dfrac{u+v}{1-\Vert u + v \Vert},-v)
&{\rm if}~v \in S(LV),~\Vert u + v\Vert <1\\[2ex]
*&{\rm otherwise}~.
\end{cases}
\end{array}$$
Let $p:S(LV)^+ \to S^0$ be the projection.
It follows from (i) that the map
$$\begin{array}{l}
D_V~=~(1\wedge p)F~:~V^{\infty} \wedge P(V)^+ \to V^{\infty}~;\\[1ex]
(u,[v]) \mapsto \begin{cases}
\dfrac{u- v}{1-\Vert u - v \Vert}
&{\rm if}~v \in S(LV),~\Vert u - v\Vert <1\\[2ex]
\dfrac{u+v}{1-\Vert u + v \Vert}
&{\rm if}~v \in S(LV),~\Vert u + v\Vert <1\\[2ex]
*&{\rm otherwise}
\end{cases}
\end{array}$$
is such that the composite stable map
$$\xymatrix@C+20pt{P(V)^+ \ar[r]^-{\di{{\rm adj}(D_V)}} &
Q_V(S^0) \ar[r]^-{\di{\psi_V}} & P(V)^+}$$
is the identity.\hfill\qed}
\end{example}

Given maps $f_1:A \to B_1$, $f_2:A \to B_2$ define the
{\it homotopy pushout} to be the double mapping cone
\index{homotopy pushout, $\Cc(f_1,f_2)$}
$$\Cc(f_1,f_2)~=~(A \times I \sqcup B_1 \sqcup B_2)/
\{(a,0) \sim f_1(a),\, (a,1) \sim f_2(a)\,\vert\, a\in A\}$$
which fits into a square
$$\xymatrix{A \ar[r]^-{\di{f_1}} \ar[d]_-{\di{f_2}}
& B_1 \ar[d]^-{\di{i_1}} \\
B_2 \ar[r]^-{\di{i_2}} &\Cc(f_1,f_2)}$$
with a canonical homotopy
$$j~:~i_1f_1 \simeq i_2f_2~:~A \to \Cc(f_1,f_2)~.$$
The homotopy pushout has the universal property that for any maps
$g_1:B_1 \to C$, $g_2:B_2 \to C$ and homotopy $h:g_1f_1 \simeq g_2f_2:A \to C$
there is a unique map $(g_1,g_2,h):\Cc(f_1,f_2) \to C$ such that
$$g_1~=~(g_1,g_2,h)i_1~,~g_2~=~(g_1,g_2,h)i_2~,~h~=~(g_1,g_2,h)j~.$$

\begin{definition}~ \label{hopf13}
{\rm (i) Given a pointed space $X$ and an inner product space $V$ let
$$J_{2,V}(X)~=~\Cc(f_1,f_2)$$
be the pointed space defined by the homotopy pushout of the maps
$$\begin{array}{l}
f_1~:~S(LV)^+\wedge_{\ZZ_2}(X \times \{*\} \cup \{*\} \times X) \to X~;~
(v,x,*)=(-v,*,x) \mapsto x~,\\[1ex]
f_2~=~{\rm inclusion}~:~
S(LV)^+\wedge_{\ZZ_2}(X \times \{*\} \cup \{*\} \times X) \to
S(LV)^+\wedge_{\ZZ_2}(X \times X)
\end{array}$$
with base point $* \in X$, so that there are defined a homotopy
commutative  diagram
$$\xymatrix{S(LV)^+\wedge_{\ZZ_2}(X \times \{*\} \cup \{*\} \times X)
\ar[d]_-{\di{f_2}} \ar[r]^-{\di{f_1}} &X \ar[d]^-{\di{i_1}} \\
S(LV)^+\wedge_{\ZZ_2}(X \times X) \ar[r]^-{\di{i_2}} & J_{2,V}(X)}$$
and a homotopy cofibration sequence
$$X \to J_{2,V}(X) \to S(LV)^+ \wedge_{\ZZ_2} (X \wedge X)~.$$
(It is assumed that the base point $* \in X$ is nondegenerate).\\
(ii) The {\it Dyer-Lashof map} is\index{Dyer-Lashof map, $d_V$}
$$d_V~:~J_{2,V}(X) \to Q_V(X)~;~
\begin{cases}
x \mapsto ((v,x) \mapsto (v,x))\\
[v,x,y] \mapsto (i(x)\vee i(y))\nabla_v
\end{cases}$$
with $\nabla_v:V^{\infty} \to V^{\infty} \vee V^{\infty}$ the
sum map defined by the Umkehr map of the open embedding
$$V \times \{v,-v\} \emb V~;~(u,\pm v) \mapsto \pm v+ \dfrac{u}{1+\Vert
u \Vert}$$
(as in \ref{sumdif} (ii)), so that
$$d_V[v,x,y]~:~V^{\infty} \to V^{\infty} \wedge X~;~
u \mapsto
\begin{cases}
(\di{u-v \over 1- \Vert u-v \Vert},x)&{\rm if}~ \Vert u-v \Vert <1\\[2ex]
(\di{u+v \over 1- \Vert u+v \Vert},y)&{\rm if}~ \Vert u+v \Vert <1\\[2ex]
*&{\rm otherwise}~.
\end{cases}$$
Write the adjoint map as
$$D_V~=~{\rm adj}(d_V)~:~V^{\infty} \wedge J_{2,V}(X) \to
V^{\infty} \wedge X~.$$
(iii) For an unpointed space $Y$ let
$$\begin{array}{l}
d^+_V~:~(S(LV)\times_{\ZZ_2}(Y \times Y))^+ \emb J_{2,V}(Y^+)
\xymatrix{\ar[r]^-{\di{d_V}}&} Q_V(Y^+)~,\\[1ex]
D^+_V~=~{\rm adj}(d^+_V)~:~V^{\infty} \wedge
(S(LV)\times_{\ZZ_2}(Y \times Y))^+ \to V^{\infty} \wedge Y^+~.
\end{array}$$
}\hfill$\qed$
\end{definition}

\begin{remark} {\rm
(i) The map
$$d_{\R^k}~:~J_{2,\R^k}(X) \to
Q_{\R^k}(X)~=~\Omega^k\Sigma^kX~~(1 \leqslant k \leqslant \infty)$$
is essentially the same as the James map $J$ of \ref{james} (for $k=1$)
and Dyer-Lashof \cite{dl} (for $k \geqslant 1$). \\
(ii) For a connected pointed space $X$ $Q_{\R^k}(X)$ is connected and
the inclusion $i:X \emb Q_{\R^k}(X)$
induces the $k$-fold suspension map in the homotopy groups
$$E^k~:~\pi_n(X) \to \pi_n(Q_{\R^k}(X))~=~\pi_{n+k}(\Sigma^kX)~;~
f \mapsto \Sigma^k f~.$$
$J_{2,\R^k}(X)$ is the second stage of the combinatorial model for $\Omega^k\Sigma^kX$
constructed by James \cite{imj2} for $k=1$ and by some combination
of Boardman-Vogt, Milgram, May, Segal, Snaith for $2 \leqslant k \leqslant \infty$,
where
$$\begin{array}{l}
Q_{\R^{\infty}}(X)~=~\varinjlim_kQ_{\R^k}X~=~\Omega^{\infty}\Sigma^{\infty}X~,\\[1ex]
\pi_n(Q_{\R^{\infty}}X)~=~\{S^n;X\}~=~\pi^S_n(X)~~(n \geqslant 0)~,\\[1ex]
\dot H_*(Q_{\R^k}(X))~=~\bigoplus\limits_{n=1}^\infty
\dot H_*(C(k,n)^+ \wedge_{\Sigma_n}(\bigwedge\limits_n X))~,\\[1ex]
\dot H_*(QX)~=~\bigoplus\limits_{n=1}^\infty
\dot H_*((E\Sigma_n)^+ \wedge_{\Sigma_n}(\bigwedge\limits_n X))
\end{array}$$
with $C(k,n)$ the $\Sigma_n$-equivariant configuration space of embeddings
$\{1,2,\dots,n\} \emb \R^k$, and $E\Sigma_n=\mathop{\varinjlim}_kC(k,n)$
a contractible space with a free $\Sigma_n$-action. In particular, for $k=1$
$$C(1,n)~\simeq_{\Sigma_n}\Sigma_n~,~
C(1,n)^+ \wedge_{\Sigma_n}(\bigwedge\limits_n X)~\simeq~
\bigwedge\limits_n X~,$$
and for $n=2$
$$C(k,2)~\simeq_{\Sigma_2} S^{k-1}~.$$
The quadratic construction on the evaluation map
$$e_{\R^k}~:~\Sigma^k(\Omega^k \Sigma^k X) \to \Sigma^k X~;~(s, \omega) \mapsto  \omega(s)$$
induces the projection
$$\psi_{\R^k}(e_{\R^k})~:~\dot H_*(\Omega^k\Sigma^k X) \to
\dot H_*((S^{k-1})^+ \wedge_{\ZZ_2}(X \wedge X))~.$$
The quadratic construction on a map $F:\Sigma^kY \to \Sigma^kX$ with adjoint
$${\rm adj}(F)~:~Y \to \Omega^k\Sigma^k X~;~y \mapsto (s \mapsto F(s,y))$$
is given by the composite
$$\xymatrix@C+15pt{
\psi_{\R^k}(F)~:~\dot H_*(Y) \ar[r]^-{\di{{\rm adj}(F)_*}} &
\dot H_*(\Omega^k \Sigma^k X) \ar[r]^-{\di{\psi_{\R^k}(e_{\R^k})}} &
\dot H_*((S^{k-1})^+ \wedge_{\ZZ_2}(X \wedge X))~.}$$
(iii) If $X$ is an $(m-1)$-connected pointed space the map
$d_{\R}:J_{2,\R}(X) \to \Omega \Sigma X$ is $(3m-2)$-connected, and for $n
\leqslant 3m-2$ the homotopy exact sequence
$$\begin{array}{l}
\xymatrix{\dots \ar[r] &\pi_n(X) \ar[r] &
\pi_n(J_{2,\R}(X)) \ar[r] & \pi_n(J_{2,\R}(X),X)}\\
\hspace*{100pt}
\xymatrix
{ \ar[r]& \pi_{n-1}(X) \ar[r] & \pi_{n-1}(J_{2,\R}(X)) \ar[r] & \dots}
\end{array}$$
becomes the $EHP$ exact sequence
$$\begin{array}{l}
\xymatrix{\dots \ar[r] &\pi_n(X) \ar[r]^-{\di{E}} &
\pi_{n+1}(\Sigma X) \ar[r]^-{\di{H}} & \pi_n(X \wedge X)}\\
\hspace*{100pt}
\xymatrix
{\ar[r]^-{\di{P}} & \pi_{n-1}(X)\ar[r]^-{\di{E}} & \pi_n(\Sigma X) \ar[r] & \dots}
\end{array}$$
with $H$ the Hopf invariant map (Whitehead \cite{wh}, Milgram \cite{milgram}).
The homotopy class of a map
$f:S^{n-1} \to X$ together with a null-homotopy $g:\Sigma f\simeq *:S^n
\to \Sigma X$ is an element $(f,g) \in \pi_n(X \wedge X)$.  The
Hurewicz image $h(f,g) \in \dot H_n(X \wedge X)$ has the following
description in terms of the quadratic construction.  Let $Y=X
\cup_fD^n$, and use $g$ to define a stable map
$$F~:~\Sigma S^n~=~ S^{n+1} \to
    \Sigma X \vee S^{n+1}\simeq \Sigma X \cup_{\Sigma f} D^{n+1}\simeq
    \Sigma Y~.$$
The quadratic construction $\psi_\R(F):\dot H_n(S^n)
\to \dot H_n(Y\wedge Y)$ sends $[S^n] \in \dot H_n(S^n)$ to
$$\psi_\R(F)[S^n]~=~h(f,g) \in \dot H_n(Y \wedge Y)~=~ \dot H_n(X \wedge X)~.$$
In particular, if $f=*:S^{n-1} \to X$ then a null-homotopy
$g:\Sigma f \simeq *$ is just a map $g:
\Sigma(S^n)=S^{n+1} \to \Sigma X$, and $F=g:S^{n+1} \to
\Sigma X$, with the Hurewicz image of $H(g) \in \pi_n(X \wedge X)$
given by $\psi_\R(g)[S^n] \in \dot H_n(X \wedge X)$. For $X=S^m$,
$n=2m$
$$\begin{array}{l}
H~=~{\rm Hopf~invariant}~:~\pi_{2m+1}(S^{m+1}) \to \pi_{2m}(S^m
\wedge S^m)~=~\ZZ~;\\[1ex]
\hspace*{100pt} (g:S^{2m+1}=\Sigma(S^{2m}) \to S^{m+1}=\Sigma(S^m))
\mapsto h_{\R}(g)~,\\[1ex]
P~:~\ZZ \to \pi_{2m-1}(S^m)~;~1 \mapsto J(\tau_{S^m})~=~[\iota,\iota]
\end{array}$$
with $\tau_{S^m} \in \pi_m(BSO(m))=\pi_{m-1}(SO(m))$.\\
(iv) For a connected unpointed space $Y$ there is defined a bijection
$$\pi_0(Q_{\R^k}(Y^+))=\pi_k(\Sigma^kY^+)
 \to \ZZ~;~(\omega:S^k \to \Sigma^kY^+) \mapsto{\rm degree}(\omega)~.$$
The Dyer-Lashof map
$$d^+_{\R^k}~:~J^+_{2,\R^k}(Y)~=~Y^+ \sqcup S^{k-1}\times_{\ZZ_2}(Y \times Y)
\emb\, J_{2,\R^k}(Y^+) \xymatrix{\ar[r]^-{\di{d_{\R^k}}}&} Q_{\R^k}(Y^+)$$
sends $Y$ to the degree 1 component of $Q_{\R^k}(Y^+)$
and $S^{k-1}\times_{\ZZ_2}(Y \times Y)$ to the degree 2 component.
$J^+_{2,\R^k}(Y)$ is the second stage of the
combinatorial model for $Q_{\R^k}(Y^+)$ provided by
the work of Barratt-Eccles, Priddy, Quillen etc., which expresses
$H_*(Q_{\R^k}(Y^+))$ for $1 \leqslant k < \infty$ (resp. $k=\infty$)
as the `group completion' of
$$\bigoplus\limits_{n=1}^\infty
\dot H_*(C(k,n) \times_{\Sigma_n}(\prod\limits_n Y))
~({\rm resp.}~\bigoplus\limits_{n=1}^\infty
\dot H_*(E\Sigma_n \times_{\Sigma_n}(\prod\limits_n Y)))~.$$
\hfill\qed}
\end{remark}

\begin{proposition}~ \label{hopf11}
The geometric Hopf invariant of the Dyer-Lashof stable map
$D_V:V^{\infty} \wedge J_{2,V}(X) \to V^{\infty} \wedge X$
$$h_V(D_V)~:~\Sigma S(LV)^+ \wedge V^{\infty} \wedge J_{2,V}(X) \to
LV^{\infty} \wedge V^{\infty} \wedge X \wedge X$$
is such that up to natural $\ZZ_2$-equivariant homotopy
$$\begin{array}{l}
h_V(D_V)(1\wedge i_1)~=~*~:~\Sigma S(LV)^+  \wedge V^{\infty}\wedge X
\to LV^{\infty}\wedge V^{\infty} \wedge X\wedge X~,\\[1ex]
h_V(D_V)(1\wedge i_2)~=~1\wedge E_V~:\\[1ex]
\hspace*{50pt}
\Sigma S(LV)^+  \wedge V^{\infty} \wedge
S(LV)^+\wedge_{\ZZ_2} (X \times X) \to LV^{\infty}\wedge V^{\infty} \wedge X\wedge X
\end{array}$$
with $i_1:X\to J_{2,V}(X)$, $i_2:S(LV)^+\times_{\ZZ_2}(X \times X) \to
J_{2,V}(X)$ as in \ref{hopf13}, and
$$\begin{array}{ll}
E_V~:&\Sigma S(LV)^+ \wedge S(LV)^+\times_{\ZZ_2}
(X \times X) \to LV^{\infty} \wedge X\wedge X~;\\[1ex]
&(u,[v,x,y]) \mapsto  \begin{cases}
(\di{u-v \over 1- \Vert u-v \Vert},x,y)&{\it if}~\Vert u-v \Vert <1\\[2ex]
(\di{u+v \over 1- \Vert u+v \Vert},y,x)&{\it if}~\Vert u+v \Vert <1\\[2ex]
*&{\it otherwise}~,
\end{cases}
\end{array}$$
using the canonical homeomorphism
$\Sigma S(LV)^+\cong (LV\backslash \{0\})^{\infty}$
to regard the elements $u\neq * \in \Sigma S(LV)^+$ as elements
$u \neq 0 \in LV$.
\end{proposition}
\begin{proof} The double cover
$$f~:~\widetilde{K}~=~S(LV) \times X \times X
\to K~=~S(LV)\times_{\ZZ_2}(X\times X)~;~(v,x,y) \mapsto [v,x,y]$$
is classified by the projection $c:K \to S(LV)/\ZZ_2=P(V)$.
As in Proposition \ref{doublecover} define a positive embedding of $f$
$$e~=~(g,f)~:~V \times \widetilde{K} \emb V \times K~;~
(u,(v,x,y)) \mapsto (v+\dfrac{u}{1+\Vert u \Vert},[v,x,y])~.$$
The Umkehr map of $f$ determined by $e$
$$\begin{array}{ll}
F~:&V^{\infty} \wedge K/K_0 \to V^{\infty} \wedge \widetilde{K}/\widetilde{K}_0~;\\
&(u,[v,x,y]) \mapsto
\begin{cases}
(\di{u-v \over 1- \Vert u-v \Vert},v,x,y)&{\rm if}~ \Vert u-v \Vert <1\\[2ex]
(\di{u+v \over 1- \Vert u+v \Vert},-v,y,x)&{\rm if}~ \Vert u+v \Vert <1\\[2ex]
*&{\rm otherwise}
\end{cases}
\end{array}$$
is such that there is defined a homotopy commutative diagram of stable maps
$$\xymatrix{
S(LV)^+\wedge_{\ZZ_2}(X \times \{*\} \cup \{*\} \times X)
\ar[d] \ar[r] & X \ar[d]_-{\di{i_1}} \ar[ddr]^-{\di{1}} &\\
S(LV)^+\wedge_{\ZZ_2}(X \times X)\ar[d]_-{\di{F}} \ar[r]^-{\di{i_2}} &
J_{2,V}(X) \ar[dr]_-{\di{D_V}}&\\
S(LV)^+\wedge (X \times X) \ar[rr]^{\di{p}} & & X}$$
with
$$p~:~S(LV)^+\wedge (X \times X)\to X~;~(v,x,y) \mapsto x~.$$
Now apply Proposition \ref{doublecover} with
$$K_0~=~P(V)~=~S(LV)\times_{\ZZ_2}(*,*)\subset K~,~
\widetilde{K}_0~=~S(LV)\times(*,*)\subset \widetilde{K}$$
such that
$$K/K_0~=~S(LV)^+\wedge_{\ZZ_2}(X \times X)~,~
\widetilde{K}/\widetilde{K}_0~=~S(LV)^+\wedge (X \times X)~.$$
The geometric Hopf invariant $h_V(F)$ fits into a
$\ZZ_2$-equivariant homotopy commutative diagram of stable $\ZZ_2$-equivariant
homotopy commutative maps
$$\xymatrix@C+10pt{
\Sigma S(LV)^+ \wedge K/K_0 \ar[d]_-{\di{h_V(F)}} \ar[r]^-{\di{1\wedge i_2}} &
\Sigma S(LV)^+ \wedge J_{2,V}(X) \ar[d]^-{\di{h_V(D_V)}}\\
LV^{\infty} \wedge \widetilde{K}/\widetilde{K}_0
\wedge \widetilde{K}/\widetilde{K}_0 \ar[r]^-{\di{1\wedge p\wedge p}} &
LV^{\infty} \wedge X  \wedge X }$$
with
$$\begin{array}{l}
\xymatrix@C+10pt@R-5pt{h_V(F)~:~
\Sigma S(LV)^+ \wedge V^{\infty} \wedge K/K_0
\ar[r]^-{\displaystyle{1\wedge F'}} &
LV^{\infty} \wedge V^{\infty} \wedge \widetilde{K}/\widetilde{K}_0}\\[1ex]
\hspace*{100pt}
\xymatrix@C+60pt{\ar[r]^-{\displaystyle{1 \wedge (1\wedge T)\Delta_{\widetilde{K}/\widetilde{K}_0}}}
&  LV^{\infty} \wedge V^{\infty} \wedge
\widetilde{K}/\widetilde{K}_0 \wedge \widetilde{K}/\widetilde{K}_0
}
\end{array}$$
where
$$F'~=~F\vert~:~(LV \backslash \{0\})^{\infty} \wedge K/K_0~=~
\Sigma S(LV)^+ \wedge K/K_0 \to V^{\infty} \wedge \widetilde{K}/\widetilde{K}_0$$
is the $\ZZ_2$-equivariant Umkehr map of the $\ZZ_2$-equivariant embedding
$$e'~=~e\vert~:~LV \times \widetilde{K} \emb (LV \backslash \{0\}) \times K~.$$
Finally, note that it follows from
$$\begin{array}{l}
D_Vi_1~=~1~:~V^{\infty} \wedge X \to V^{\infty} \wedge X~,\\[1ex]
D_Vi_2~=~(1\wedge p)F~:~V^{\infty} \wedge K/K_0 \to V^{\infty} \wedge X
\end{array}$$
that
$$\begin{array}{l}
h_V(D_V)(1\wedge i_1)~=~h_V(1)~=~*~:\\[1ex]
\hspace*{50pt} \Sigma S(LV)^+  \wedge V^{\infty}\wedge X
\to LV^{\infty}\wedge V^{\infty} \wedge X\wedge X~,\\[1ex]
h_V(D_V)(1\wedge i_2)~=~(1 \wedge p \wedge p)h_V(F)~=~1\wedge E_V~:\\[1ex]
\hspace*{50pt}
\Sigma S(LV)^+  \wedge V^{\infty} \wedge
S(LV)^+\wedge_{\ZZ_2} (X \times X) \to LV^{\infty}\wedge V^{\infty} \wedge X\wedge X~.
\end{array}$$
\hfill\qed\end{proof}

\begin{proposition}~ \label{quad4}
{\rm (i)} For a pointed space $X$ the quadratic constructions on the maps
$$\begin{array}{l}
D_V~:~V^{\infty} \wedge J_{2,V}(X) \to V^{\infty} \wedge X~,\\[1ex]
e_V~:~V^{\infty} \wedge Q_V(X) \to V^{\infty}\wedge X
\end{array}$$
are stable maps
$$\begin{array}{l}
\psi_V(D_V)~:~J_{2,V}(X) \sto S(LV)^+\wedge_{\ZZ_2} (X \wedge X)~,
\\[1ex]
\psi_V(e_V)~:~Q_V(X)\sto S(LV)^+\wedge_{\ZZ_2} (X\wedge X)
\end{array}$$
such that up to stable homotopy
$$\psi_V(D_V)~=~\psi_V(e_V)d_V~=~{\rm projection}~:~
J_{2,V}(X) \sto S(LV)^+\wedge_{\ZZ_2} (X \wedge X)~.$$
{\rm (ii)} For an unpointed space $Y$ the quadratic constructions on the maps
$$\begin{array}{l}
D^+_V~:~V^{\infty} \wedge (S(LV)^+ \wedge_{\ZZ_2}(Y \times Y)^+) \to V^{\infty} \wedge Y^+~,\\[1ex]
e_V~:~V^{\infty} \wedge Q_V(Y^+) \to V^{\infty}\wedge Y^+
\end{array}$$
are such that
$$\psi_V(D^+_V)~=~\psi_V(e_V)d^+_V~=~1~:~(S(LV)\times_{\ZZ_2} (Y \times Y))^+
\sto (S(LV)\times_{\ZZ_2} (Y \times Y))^+~.$$
\end{proposition}
\begin{proof} (i) By the proof of Proposition \ref{hopf11} there is defined a
commutative diagram of stable maps
$$\xymatrix{
S(LV)^+\wedge_{\ZZ_2}(X \times \{*\} \cup \{*\} \times X)
\ar[d] \ar[rrr] & & & X \ar[dll]_-{\di{i_1}} \ar@{=}[dd]^-{\di{1}}
\ar[dl]^-{\di{i}}\\
S(LV)^+\wedge_{\ZZ_2}(X \times X)\ar[r]^-{\di{i_2}}\ar[d]_-{\di{F}}
& J_{2,V}(X) \ar[r]^-{\di{d_V}}\ar[drr]_-{\di{D_V}} &
Q_V(X) \ar[dr]^-{\di{e_V}} &\\
S(LV)^+\wedge_{\ZZ_2}(X \wedge X) \ar[rrr]^-{\di{p}} &&& X}$$
with
$$F~:~ S(LV)^+ \wedge_{\ZZ_2}(X \times X) \sto  S(LV)^+ \wedge(X \times X)$$
the Umkehr for the double cover
$$f~:~S(LV) \times (X \times X) \to S(LV)\times_{\ZZ_2}(X \times X)$$
and
$$p~:~S(LV)^+ \wedge (X \times X) \to X~;~(v,x,y) \mapsto x~.$$
Now
$$\psi_V(D_V)i_1~=~\psi_V(1)~=~*~:~X \sto S(LV)^+\wedge_{\ZZ_2} (X \wedge X)~,$$
and by the finite cover formula of Proposition \ref{quad3} (vii) applied to $F$
$$\begin{array}{l}
\psi_V(D_V)i_2~=~(p \wedge p)\psi_V(F)~=~ {\rm projection}~:\\[1ex]
\hskip100pt S(LV)^+\wedge_{\ZZ_2}(X \times X) \sto S(LV)^+\wedge_{\ZZ_2} (X \wedge X)~,
\end{array}$$
so that by composition formula of \ref{quad3} (vi) applied to
$D_V=e_Vd_V$
$$\psi_V(D_V)~=~\psi_V(e_V)d_V~=~{\rm projection}~:~
J_{2,V}(X) \sto S(LV)^+\wedge_{\ZZ_2} (X \wedge X)~.$$
(ii) Immediate from (i), with
$$\begin{array}{l}
X~=~Y^+~=~Y \sqcup \{*\}~,\\[1ex]
X \wedge X~=~(Y \times Y)^+  \subset
X \times X~=~(Y \times Y)^+ \sqcup (Y \times \{*\} \sqcup \{*\} \times Y)~,\\[1ex]
S(LV)^+\wedge_{\ZZ_2}(X \wedge X)~=~
(S(LV)\times_{\ZZ_2}(Y \times Y))^+\subset S(LV)^+\wedge_{\ZZ_2}(X \times X)~.
\end{array}
$$
\hfill\qed\end{proof}

\begin{example} {\rm As in Proposition \ref{doublecover}
let $K$ be a space with a map $c:K \to P(V)$, classifying
the double cover
$$\widetilde{K}~=~c^*S(LV)~=~
\{(v,x) \in S(LV)\times K\,\vert\,[v]=c(x)\in P(V)\}$$
of $K$. The Umkehr map $F:V^{\infty} \wedge K^+ \to
V^{\infty} \wedge \widetilde{K}^+$ is such that there is
defined a commutative square
$$\xymatrix{K^+ \ar[r]^-{\di{c}} \ar[d]_-{\di{{\rm adj}(F)}}&
P(V)^+ \ar[d]^-{\di{D^+_V}} \\
Q_V(\widetilde{K}^+) \ar[r]^-{\di{p}} &Q_V(S^0)}$$
and up to stable homotopy
$$\xymatrix@C+30pt{c~:~K^+ \ar[r]^-{\di{p\,{\rm adj}(F)}}&
Q_V(S^0) \ar[r]^-{\di{\psi_V}} & P(V)^+~.}$$
\hfill\qed}
\end{example}

\begin{remark} {\rm The Dyer-Lashof map $d^+_V:P(V)^+ \to O(V)$
is related to the reflection map $R:P(V)^+ \to O(V)$ by the
identity (up to stable homotopy)
$$d^+_V~=~1-JR~:~P(V)^+ \to Q_V(S^0)~.$$
It is enough to verify the special case $V=\R$ directly, the general
case following by naturality, working as in the proof of $\theta \circ R=1$
in Proposition \ref{localR} (ii) with the commutative diagrams
$$\xymatrix@C+30pt{P(\R)^+ \ar[r]^-{\di{d^+_{\R}=1-JR}} \ar[d]_-{\di{v}} &
Q_{\R}(S^0) \ar[d]^-{\di{v}} \\
P(V)^+ \ar[r]^-{\di{d^+_V=1-JR}}  & Q_V(S^0)}$$
defined for $v \in P(V)$.}\\
\hfill\qed
\end{remark}

\section{Embeddings and immersions}\label{double}

We shall now consider the application of the geometric Hopf
invariant of \S\ref{hopf} to an embedding-immersion pair
$$(e\times f:M \emb V \times N~,~f:M \imm N)$$
of manifolds $M,N$, with an Umkehr map $F:V^{\infty}\wedge N^\infty \to V^{\infty}\wedge T(\nu_f)$.
The basic result (\ref{double6}) applies the Double Point Theorem
\ref{doubleHopf} to identify the quadratic construction
$$\psi_V(F)~:~N^\infty \sto S(LV)^+\wedge_{\ZZ_2}(T(\nu_f)\wedge T(\nu_f))$$
with the Pontrjagin-Thom map of the immersion
$D_2(f)^{2m-n} \imm N$ of the double point manifold.

\begin{definition}~ \label{double2}{\rm
An {\it $(m,n,j)$-dimensional embedding-immersion pair}
\index{embedding-immersion pair}
$$(\  e\times f:M \emb V \times N\ ,\  f:M \imm N)$$
consists of an embedding $e \times f$ and an immersion $f$, with
$M$ an $m$-dimensional manifold, $N$ an $n$-dimensional manifold, and
$V$ a $j$-dimensional inner product space.\hfill\qed
}\end{definition}

It will always be assumed that $f$ only has transverse
self-intersections, so that the double point sets $D_2(f)$, $D_2[f]$
(\ref{double1}) are $(2m-n)$-dimensional manifolds.

The immersion $f:M \imm N$ has a normal bundle $\nu_f:M \to BO(n-m)$,
and $f$ extends to a codimension 0 immersion $E(f):E(\nu_f) \imm N$
of the total space of $\nu_f$.
The embedding $e \times f:M \emb V \times N$
is regular homotopic to the composite $M
\imm N \emb V \times N$ of $f$ and the embedding
$$N \emb V \times N~;~x \mapsto (0,x)$$
which has trivial normal bundle $\epsilon^j:N \to
BO(j)$.  It follows that the normal bundle of $e\times f$ is
$$\nu_{e \times f}~=~\nu_f\oplus \epsilon^j~:~M \to BO(n-m+j)~.$$
The product immersion $1 \times E(f):V \times E(\nu_f) \imm V \times N$ is
regular homotopic to an embedding
$E(e \times f):V \times E(\nu_f) \emb V \times N$
 of the immersion $E(\nu_f) \imm N$ by Proposition \ref{immersions} (ii).

\begin{definition}~ \label{double3}
{\rm The {\it adjunction Umkehr map} (or {\it Pontrjagin-Thom map})
of an $(m,n,j)$-dimensional
embedding-immersion pair $(e \times f,f)$ is the adjunction Umkehr map of
$E(e\times f)$ \index{Umkehr!adjunction}\index{Pontrjagin-Thom}
$$F~:~V^{\infty} \wedge N^\infty \to V^{\infty} \wedge T(\nu_f)~;~
(a,x) \mapsto \begin{cases}(b,y)&\hbox{if $(a,x)=E(e \times f)(b,y)$} \\[1ex]
*&\hbox{otherwise}~.
\end{cases}$$
\hfill\qed
}\end{definition}

\begin{example} {\rm
A finite cover $f:M \to N$ of an $n$-dimensional manifold is an immersion
with normal bundle $\nu_f:M \to BO(0)=\{*\}$. An embedding
$e:V \times M \emb V \times N$ of $f$ in the sense of \ref{embed0}
with ${\rm dim}(V)=j$ determines an $(n,n,j)$-dimensional
embedding-immersion pair $(e_0 \times f,f)$ with
$$e_0~:~M \to V~;~x \mapsto v~{\rm if}~e(0,x)~=~(v,f(x))~.$$
The Umkehr map given by \ref{double3} is the same as the Umkehr map
given by \ref{embed0}
$$F~:~V^{\infty}\wedge N^\infty \to V^{\infty} \wedge T(\nu_f)~=~
V^{\infty} \wedge M^+~.$$
\hfill\qed}
\end{example}

The construction of \ref{double3} will also be applied to an embedding
$f: M \emb V \times N$ when $M$ is a disjoint union of manifolds of
different dimensions.

\begin{remark}\label{double4}
{\rm By the embedding theorem of Whitney \cite{whitney1} for any immersion $f:M^m
\imm N^n$ and $j > 2m-n$ there exists a map $e:M\to V=\R^j$ such
that $e(x)\not= e(y)$ whenever $f(x)=f(y)$ and $x \neq y$, i.e.
such that
$$e\times f ~:~ M\emb V\times N ~;~ x\mapsto (e(x),f(x))$$
is an embedding, and $(e \times f,f)$ is an $(m,n,j)$-dimensional
embedding-immersion pair.\hfill\qed}
\end{remark}

\begin{definition}~ \label{double5}{\rm
The {\it square} of an $(m,n,j)$-dimensional embedding-immersion pair
$(e\times f:M \emb V \times N,f:M \imm N)$
is the $\ZZ_2$-equivariant $(2m,2n,2j)$-dimensional
embedding-immersion pair
\index{square embedding-immersion}
$$(\  g:M\times M \emb (LV\oplus V)
\times N\times N\ ,\
f \times f:M\times M \imm N\times N\ )$$
with
$$\begin{array}{l}
g~=~(e\times f) \times (e \times f)~:~
M\times M \emb (LV \oplus V)\times (N\times N)~;\\[1ex]
\hskip50pt
(x,y)\mapsto ({1\over 2}(e(x)-e(y)),{1\over 2}(e(x)+e(y)),f(x),f(y))~.\end{array}$$
\hfill\qed
}\end{definition}

The embedding $g$ has normal bundle
$$\nu_g~=~(\nu_f\oplus \epsilon^j) \times (\nu_f\oplus \epsilon^j)~:~
M \times M \to BO(2(n-m+j))$$
and $\ZZ_2$-equivariant Umkehr map
$$F\wedge F~:~  LV^{\infty}\wedge V^{\infty}\wedge N^\infty\wedge N^\infty \to
T(\nu_g)~=~LV^{\infty}\wedge V^{\infty} \wedge T(\nu_f)\wedge T(\nu_f)~.$$
The $(2m-n)$-dimensional ordered double point manifold
$$\begin{array}{ll}
D_2(f)&=~\{(x,y) \in M \times M\st x \neq y \in M,f(x)=f(y) \in N\}\\[1ex]
&=~\{(x,y) \in M \times M\st e(x) \neq e(y) \in V,f(x)=f(y) \in N\}\end{array}$$
fits into a pullback square of $\ZZ_2$-equivariant embeddings
$$\xymatrix@C+20pt@R+20pt{
M \sqcup D_2(f)~  \ar@{^{(}->}[r]^-{\di{g_1 \sqcup g_2}} \ar@{^{(}->}[d]_-{\di{\Delta_M \sqcup i_2}}
&(LV \oplus V) \times N
\ar@{^{(}->}[d]^{\di{1\times \Delta_N}}\\
M \times M~ \ar@{^{(}->}[r]^-{\di{g}}
&(LV \oplus V) \times N\times N}$$
with
$$\begin{array}{l}
g_1~:~M \emb (LV \oplus V) \times N~;~x \mapsto (0,e(x),f(x))~,\\[1ex]
g_2~:~D_2(f) \emb (LV \oplus V) \times N~;~(x,y) \mapsto g(x,y)~,\\[1ex]
i_2~:~D_2(f) \emb M \times M~;~(x,y) \mapsto (x,y)~,\end{array}$$
and $g_2$ restricts to a $\ZZ_2$-equivariant embedding
$$g_3~:~D_2(f) \emb (LV\backslash \{0\} \times V) \times N~.$$
The normal bundles are given by
$$\begin{array}{l}
\nu_{g_1}~=~\nu_f \oplus \epsilon^j \oplus L\epsilon^j~:~
M \to BO(n-m+2j)~,\\[1ex]
\nu_{g_2}~=~\nu_{g_3}~=~
((\nu_f\oplus \epsilon^j)\times (\nu_f \oplus \epsilon^j))\vert~:~
D_2(f) \to BO(2n-2m+2j)~.\end{array}$$

\begin{proposition}~ \label{double6}{\it
Let $(e \times f:M \emb V \times N,f:M \imm N)$ be an
$(m,n,j)$-dimensional embedding-immersion pair.
Write the $\ZZ_2$-equivariant Umkehr maps of $g_1,g_2,g_3$ as
$$\begin{array}{l}G_1~:~LV^{\infty} \wedge V^{\infty} \wedge N^\infty \to
T(\nu_{g_1})~=~LV^{\infty} \wedge V^{\infty} \wedge T(\nu_f)~,\\[1ex]
G_2~:~LV^{\infty} \wedge V^{\infty} \wedge N^\infty \to
T(\nu_{g_2})~=~LV^{\infty} \wedge V^{\infty} \wedge
T(\nu_f\times \nu_f\vert_{D_2(f)})~,\\[1ex]
G_3~:~(LV\backslash \{0\})^{\infty} \wedge V^{\infty} \wedge N^\infty~
=~\Sigma S(LV)^+ \wedge V^{\infty} \wedge N^\infty \\[1ex]
\hskip100pt \to
T(\nu_{g_3})~=~LV^{\infty} \wedge V^{\infty} \wedge
T(\nu_f\times \nu_f\vert_{D_2(f)})~.\end{array}$$
{\rm (i)} The Umkehr maps for $g_1,g_2,g_3$ are such that
$$G_1~=~1\wedge F~~,~~G_2~=~G_3(\alpha_{LV} \wedge 1)~.$$
{\rm (ii)} There is defined a homotopy
$$h_V(F)~ \simeq~i_2G_3~:~
\Sigma S(LV)^+ \wedge V^{\infty} \wedge N^\infty  \to
LV^{\infty} \wedge V^{\infty} \wedge T(\nu_f)\wedge T(\nu_f)$$
with $h_V(F)=\delta(p,q)$ the geometric Hopf invariant of $F$,
i.e. the relative difference of the maps
$$\begin{array}{l}
p~=~(1 \wedge \Delta_{T(\nu_f)})(1 \wedge F)~,~
q~=~(F \wedge F)(1 \wedge \Delta_N)~:\\[1ex]
\hskip50pt LV^{\infty} \wedge V^{\infty} \wedge N^\infty \to
LV^{\infty} \wedge V^{\infty} \wedge T(\nu_f)\wedge T(\nu_f)~.\end{array}$$
which agree on $0^+ \wedge V^{\infty} \wedge N^\infty$.}\\
\end{proposition}
\noindent{\it Proof\  .}
(i) Immediate from the expressions of $g_1,g_2$ as composites
$$\begin{array}{l}
g_1~:~M \xymatrix@C+10pt{\ar@{^{(}->}[r]^-{\di{e \times f}}&}
V \times N \emb (LV \oplus V) \times N~,\\[1ex]
g_2~:~D_2(f)\xymatrix@C+10pt{\ar@{^{(}->}[r]^-{\di{g_3}}&}
(LV \backslash \{0\}) \times V \times N \emb (LV \oplus V) \times N~.
\end{array}$$
(ii) This is a special case of Theorem \ref{doubleHopfpositive}.
\hfill\qed

\begin{definition}~ \label{double8} {\rm
The {\it (second) extended power} of a $k$-plane bundle
$\alpha:M \to BO(k)$ with respect to an inner product space $V$
\index{second extended power, $e_2(\alpha)$}
$$e_2(\alpha)~:~S(LV)\times_{\ZZ_2}(M \times M) \to BO(2k)$$
is the $2k$-plane bundle with total space
$$E(e_2(\alpha))~=~S(LV) \times_{\ZZ_2} (E(\alpha) \times E(\alpha))~.$$
\hfill\qed
}\end{definition}

The Thom space of the extended power bundle is given by
$$T(e_2(\alpha))~=~S(LV)^+\wedge_{\ZZ_2} (T(\alpha) \wedge T(\alpha))~.$$

\begin{example} \label{double9}
{\rm For the trivial $k$-plane bundle over a space $M$
$$\alpha~=~\epsilon^k~~,~~E(\alpha)~=~V \times M~~,~~T(\alpha)~=~\Sigma^kM^+$$
with $V=\R^k$. The isomorphism of $\R[\ZZ_2]$-modules
$$\kappa_V~:~LV\oplus V \to V\oplus V ~; ~(v,w)\mapsto (v+w,-v+w)$$
can be used to identify
$$e_2(\epsilon^k)~=~k\lambda \oplus \epsilon^k~:~S(LV)\times_{\ZZ_2}(M \times M) \to BO(2k)$$
where $\lambda:S(LV)\times_{\ZZ_2}(M \times M) \to BO(1)$ is the line bundle
with
$$\begin{array}{l}E(\lambda)~=~S(LV) \times_{\ZZ_2}(L \R \times M \times M)~,\\[1ex]
T(\lambda)~=~(S(LV \oplus \R)/ S(LV)) \wedge_{\ZZ_2}(M^+ \wedge M^+)~,
\end{array}$$
and
$$\begin{array}{l}
E(e_2(\epsilon^k))~=~S(LV) \times_{\ZZ_2} ((V \times M) \times (V \times M))~,\\[1ex]
T(e_2(\epsilon^k))~=~\Sigma^k( (S(LV \oplus V)/ S(LV)) \wedge_{\ZZ_2}
(M^+ \wedge M^+))~.\end{array}$$
In the special case $k=1$, $V=\R$
$$\begin{array}{l}S(LV)~=~S^0~,~\lambda~=~\epsilon~,\\[1ex]
e_2(\epsilon^k)~=~\epsilon^{2k}~:~
S(LV) \times_{\ZZ_2} (M \times M)~=~M \times M \to BO(2k)~,\\[1ex]
T(e_2(\epsilon^k))~=~\Sigma^kM^+ \wedge
\Sigma^kM^+~.\end{array}$$ \hfill\qed}
\end{example}

\begin{proposition}~ \label{double10} {\it
Let $(e \times f:M \emb V \times N,f:M \imm N)$ be an
$(m,n,j)$-dimensional embedding-immersion pair, and write
the $(2m-n)$-dimensional unordered double point manifold as
$$M'~=~D_2[f]~.$$
{\rm (i)} The immersion
$$f'~:~ M' \imm N~;~[x,y] \mapsto f(x)=f(y)$$
has normal bundle
$$\nu_{f'}~=~c^*e_2(\nu_f)~:~M' \to BO(2(n-m))$$
with
$$c~:~M' \to S(LV) \times_{\ZZ_2}(M \times M)~;~[x,y] \mapsto
[{e(x)-e(y) \over \Vert e(x) - e(y) \Vert},x,y]~.$$
{\rm (ii)} The quadratic construction
{\rm (}\ref{quad1}{\rm )} on the Umkehr map
$F:V^{\infty}\wedge N^\infty \to V^{\infty}\wedge T(\nu_f)$ is given
by the composite
$$\psi_V(F)~:~N^\infty \xymatrix{\ar[r]^-{\di{F'}}&} T(\nu_{f'})
 \xymatrix{\ar[r]^-{\di{T(c)}}&} T(e_2(\nu_f))=
 S(LV)^+ \wedge_{\ZZ_2}(T(\nu_f) \wedge T(\nu_f))$$
with $F':N^\infty \sto T(\nu_{f'})$ the Umkehr stable map of $f'$.\\
{\rm (iii)} The quadratic construction $\psi_V(F)$
sends the fundamental class $[N]\in H_n(N)$ to the image of the fundamental class
$[M'] \in H_{2m-n}(M';\ZZ^{w'})$
$$\begin{array}{l}
\psi_V(F)_*[N]~=~c_*[M'] \\[1ex]
 \in\widetilde{H}_n(S(LV)^+\wedge_{\ZZ_2}(T(\nu_f) \wedge T(\nu_f)))=
H_{2m-n}(S(LV)\times_{\ZZ_2}(M \times M);\ZZ^{w'})\end{array}$$
where $\ZZ^{w'}$ refers to $\ZZ$ twisted by the orientation character
$w'=w_1(e_2(\nu_f))$.}
\end{proposition}
\noindent{\it Proof\  .}
This is just a matter of looking at the duality construction for
the sphere $S(LV)$ and the homotopy $h_V(F) \simeq i_2G_3$ of \ref{double6}
(ii).\hfill\qed

\begin{example} \label{double11}{\rm
(i) A finite covering $f:M \to N$ of an $n$-dimensional
manifold $N$ is a framed immersion of the $n$-dimensional manifold $M$,
and there exists a map
$e:M \to \R^j$ such that $e \times f:M \to \R^j \times N$ is an
embedding, defining an $(n,n,j)$-dimensional embedding-immersion pair
$$(\  e\times f:M \emb \R^j \times N\ ,\  f:M \imm N \ )$$
with an Umkehr stable map
$F:\Sigma^j N^\infty \to \Sigma^j M^+$.
The quadratic construction on $F$ is the composite
$$\psi_{\R^j}(F)~:~N^\infty \xymatrix{\ar[r]^-{\di{F'}}&}
{M'}^+ \xymatrix{\ar[r]^-{\di{c}}&}
(S(L\R^j) \times_{\ZZ_2}(M \times M))^+$$
with $F': N^\infty \sto {M'}^+$ the Umkehr stable map of the
finite covering map
$$f'~:~M'~=~D_2[f] \to N~;~[x,y]  \mapsto f(x)=f(y)$$
and
$$c~:~M' \to S(L\R^j) \times_{\ZZ_2}(M \times M)~;~
[x,y] \mapsto [{e(x)-e(y) \over \Vert e(x) - e(y) \Vert},x,y]~.$$
(ii) For a double covering $f:M \to N$ with covering translation
$T:M \to M$ there is defined a canonical homeomorphism
$$N \to M'~;~f(x) \mapsto [x,Tx]~~(x \in M)~,$$
with
$$f'~=~1~:~M'~=~N  \to N~.$$
For any $e:M \to \R^j$ with $e \times f:M \to \R^j \times N$ an embedding
the quadratic construction on the Umkehr map
$F:\Sigma^jN^\infty \to \Sigma^jM^+$ is given by
$$\psi_V(F)~=~c~:~N^\infty \to
(S^{j-1})^+\wedge_{\ZZ_2}(M^+ \wedge M^+)$$
with
$$c~:~M'~=~N \to S^{j-1}\times_{\ZZ_2}(M \times M)~;~
f(x) \mapsto [{e(x)-e(Tx) \over \Vert e(x) - e(Tx) \Vert},x,Tx]~.$$
(iii) The stable homotopy operation of Segal \cite{segal2}
$$\theta^2~:~\pi^0_S(X)~=~\{X;S^0\}
\to \pi^0_S(X;B\Sigma_2)~=~\{X;(B\Sigma_2)^+\}$$
is such that for any $F:\Sigma^jX \to S^j$
$$\theta^2(F)~:~X \xymatrix@C+15pt{\ar[r]^-{\di{\psi_{\R^j}(F)}}&}
(S(L\R^j)/\ZZ_2)^+~=~P(\R^j)^+ \emb (B\Sigma_2)^+$$
with
$$B\Sigma_2~=~P(\R(\infty))~=~\bigcup^{\infty}_{j=1}P(\R^j)~.$$
\hfill\qed}
\end{example}

\begin{remark}~{\rm
(i) A stable map $F:V^{\infty} \wedge X \to V^{\infty} \wedge Y$
induces a morphism of direct sum systems
$$\xymatrix@R+5pt{
\{X;S(\infty)^+\wedge_{\ZZ_2}(X \wedge X)\}
\ar@<1ex>[r]^-{\di{\gamma}}
\ar@<-1ex>@{<-}[r]_-{\di{\delta}}
\ar[d]^-{\di{(1\wedge F \wedge F)_*}}&
\{X;X \wedge X\}_{\ZZ_2} \ar[d]^-{\di{(F \wedge F)_*}}
\ar@<1ex>[r]^-{\di{\rho}}
\ar@<-1ex>@{<-}[r]_-{\di{\sigma}}& \{X;X\}\ar[d]^-{\di{F_*}}\\
\{X;S(\infty)^+\wedge_{\ZZ_2}(Y \wedge Y)\}
\ar@<1ex>[r]^-{\di{\gamma}}
\ar@<-1ex>@{<-}[r]_-{\di{\delta}}&
\{X;Y \wedge Y\}_{\ZZ_2}
\ar@<1ex>[r]^-{\di{\rho}}
\ar@<-1ex>@{<-}[r]_-{\di{\sigma}}&
\{X;Y\}
}$$
with
$$\begin{array}{l}
(F \wedge F)_*~=~\begin{pmatrix} (1\wedge F \wedge F)_* & h_{\infty}(F) \\
0 & F_* \end{pmatrix}~:\\[2ex]
\{X;X \wedge X\}_{\ZZ_2}~=~\{X;S(\infty)^+\wedge_{\ZZ_2}(X \wedge X)\}
\oplus \{X;X\}\to \\[1ex]
\hskip50pt
\{X;Y \wedge Y\}_{\ZZ_2}~=~\{X;S(\infty)^+\wedge_{\ZZ_2}(Y \wedge Y)\}
\oplus \{X;Y\}~,
\end{array}$$
where $h_{\infty}(F)=h_{\R^{\infty}}(F)\in \{X;S(\infty)^+\wedge_{\ZZ_2}(Y \wedge Y)\}$.
The image of
$$\Delta_X~=~(0,1) \in
\{X;X \wedge X\}_{\ZZ_2}~=~\{X;S(\infty)^+\wedge_{\ZZ_2}(X \wedge X)\}
\oplus \{X;X\}$$
is
$$\begin{array}{l}
(F \wedge F)_*(\Delta_X)~=~(h_{\infty}(F),F)\in \\[1ex]
\hskip75pt
\{X;Y \wedge Y\}_{\ZZ_2}~=~\{X;S(\infty)^+\wedge_{\ZZ_2}(Y \wedge Y)\}
\oplus \{X;Y\}~.
\end{array}$$
(ii) Now suppose given an $(m,n,j)$-dimensional embedding-immersion
pair $(e \times f:M \emb V \times N,f:M \imm N)$ with Umkehr map
$F:V^{\infty} \wedge N^\infty \to V^{\infty} \wedge T(\nu_f)$.
The induced morphism of direct sum systems
$$\xymatrix@R+5pt{
\omega_n(S(\infty)\times_{\ZZ_2}(N \times N))
\ar@<1ex>[r]^-{\di{\gamma}}
\ar@<-1ex>@{<-}[r]_-{\di{\delta}}
\ar[d]^-{\di{(1\wedge F \wedge F)_*}}&
\omega_{n,n}^{\ZZ_2}(N \times N) \ar[d]^-{\di{(F \wedge F)_*}}
\ar@<1ex>[r]^-{\di{\rho}}
\ar@<-1ex>@{<-}[r]_-{\di{\sigma}}& \omega_n(N)\ar[d]^-{\di{F_*}}\\
\widetilde{\omega}_n(S(\infty)^+\wedge_{\ZZ_2}(T(\nu_f) \wedge T(\nu_f)))
\ar@{=}[d]
\ar@<1ex>[r]^-{\di{\gamma}}
\ar@<-1ex>@{<-}[r]_-{\di{\delta}}&
\widetilde{\omega}^{\ZZ_2}_{n,n}(T(\nu_f) \wedge T(\nu_f))\ar@{=}[d]
\ar@<1ex>[r]^-{\di{\rho}}
\ar@<-1ex>@{<-}[r]_-{\di{\sigma}}&
\widetilde{\omega}_n(T(\nu_f))\ar@{=}[d]\\
\omega_{2m-n}(S(\infty) \times_{\ZZ_2}(M \times M))
\ar@<1ex>[r]^-{\di{\gamma}}
\ar@<-1ex>@{<-}[r]_-{\di{\delta}}&
\omega_{2m-n,m}^{\ZZ_2}(M \times M)
\ar@<1ex>[r]^-{\di{\rho}}
\ar@<-1ex>@{<-}[r]_-{\di{\sigma}}&
\omega_m(M)
}$$
is such that
$$\begin{array}{l}
(F \wedge F)_*~=~\begin{pmatrix} (1\wedge F \wedge F)_* & h_{\infty}(F) \\
0 & F_* \end{pmatrix}~:\\[2ex]
\omega^{\ZZ_2}_{n,n}(N \times N)~=~\omega_n(S(\infty)\times_{\ZZ_2}(N \times N))
\oplus \omega_n(N) \to \\[1ex]
\widetilde{\omega}^{\ZZ_2}_{n,n}(T(\nu_f) \wedge T(\nu_f))~=~
\widetilde{\omega}_n(S(\infty)^+\wedge_{\ZZ_2}(T(\nu_f) \wedge T(\nu_f)))
\oplus
\widetilde{\omega}_n(T(\nu_f))\\[1ex]
\hskip50pt =~\omega_{2m-n,m}^{\ZZ_2}(M \times M)~=~
\omega_{2m-n}(S(\infty) \times_{\ZZ_2}(M \times M)) \oplus \omega_m(M)~.
\end{array}$$
The image under $(F\wedge F)_*$ of the class in
$$
\omega_n^{\ZZ_2}(N\times N)=
\omega_{n,n}(N \times N)~=~\omega_n(S(\infty)\times_{\ZZ_2}(N \times N))
\oplus \omega_n(N)$$
represented by the framed manifold $(N,\Delta_N)\in\Omega_n^{\ZZ_2\hbox{-}fr}(N\times N)$,
corresponding to $(0,(N,1)) \in
\Omega^{fr}_n(S(\infty )\times_{\ZZ_2}(N\times N))
\oplus\Omega^{fr}_n(N)$,
is represented by
$$
((D_2[f],i_2),(M,1))\in
\Omega^{fr}_{2m-n}(S(\infty) \times_{\ZZ_2}(M \times M)) \oplus \Omega^{fr}_m(M)~.
$$
(iii) For any pointed spaces $X,Y$ there is defined a nonadditive function of direct sum systems
$$\xymatrix@R+5pt{
\{X;S(\infty)^+\wedge_{\ZZ_2}Y\}
\ar@<1ex>[r]^-{\di{\gamma}}
\ar@<-1ex>@{<-}[r]_-{\di{\delta}}
\ar[d]^-{\di{1\wedge P^2}}&
\{X;Y\}_{\ZZ_2} \ar[d]^-{\di{P^2}}
\ar@<1ex>[r]^-{\di{\rho}}
\ar@<-1ex>@{<-}[r]_-{\di{\sigma}}& \{X;Y\}\ar@{=}[d]\\
\{X;S(\infty)^+\wedge_{\ZZ_2}(Y \wedge Y)\}
\ar@<1ex>[r]^-{\di{\gamma}}
\ar@<-1ex>@{<-}[r]_-{\di{\delta}}&
\{X;Y \wedge Y\}_{\ZZ_2}
\ar@<1ex>[r]^-{\di{\rho}}
\ar@<-1ex>@{<-}[r]_-{\di{\sigma}}&
\{X;Y\}
}$$
with
$$\begin{array}{l}
P^2~:~\{X;Y\}_{\ZZ_2} \to \{X;Y \wedge Y\}_{\ZZ_2}~;~F \mapsto (F \wedge F)\Delta_X~,\\[1ex]
P^2~=~\begin{pmatrix} 1\wedge P^2 & h_{\infty}(E) \\
0 & 1 \end{pmatrix}~:\\[2ex]
\{X;Y\}_{\ZZ_2}~=~\{X;S(\infty)^+\wedge_{\ZZ_2}Y\}
\oplus \{X;Y\}\to \\[1ex]
\hskip50pt
\{X;Y \wedge Y\}_{\ZZ_2}~=~\{X;S(\infty)^+\wedge_{\ZZ_2}(Y \wedge Y)\}
\oplus \{X;Y\}~,
\end{array}$$
using $\{X;Y\}=[X,\Omega^{\infty} \Sigma^{\infty} Y]$
and the geometric Hopf invariant
$$h_{\infty}(E) \in \{\Omega^{\infty} \Sigma^{\infty} Y;
S(\infty)^+\wedge_{\ZZ_2}(Y \wedge Y)\}$$
of the evaluation stable map
$$E~:~\Sigma^{\infty}(\Omega^{\infty} \Sigma^{\infty} Y)\to \Sigma^{\infty}Y$$
(Proposition \ref{eval}).
For any stable map $F:\Sigma^{\infty}X \to \Sigma^{\infty}Y$ the image of
$$F~=~(0,F) \in
\{X;Y\}_{\ZZ_2}~=~\{X;S(\infty)^+\wedge_{\ZZ_2}Y\}  \oplus \{X;Y\}$$
is
$$\begin{array}{l}
P^2(F)~=~(h_{\infty}(E)(F),F)~=~(h_{\infty}(F),F) \\[1ex]
\hskip75pt\in\{X;Y \wedge Y\}_{\ZZ_2}~=~\{X;S(\infty)^+\wedge_{\ZZ_2}(Y \wedge Y)\}
\oplus \{X;Y\}~.
\end{array}$$
(iv) For an immersion $f:M^m \imm N^n$ consider the
nonadditive function of direct sum systems
$$\xymatrix@R+5pt@C-5pt{
\{N^\infty;S(\infty)\wedge_{\ZZ_2}T(\nu_f)\}
\ar@<1ex>[r]^-{\di{\gamma}}
\ar@<-1ex>@{<-}[r]_-{\di{\delta}}
\ar[d]^-{\di{1\wedge P^2}}&
\{N^\infty;T(\nu_f)\}_{\ZZ_2}
\ar[d]^-{\di{P^2}}\ar@<1ex>[r]^-{\di{\rho}}
\ar@<-1ex>@{<-}[r]_-{\di{\sigma}}&
\{N^\infty;T(\nu_f)\}\ar@{=}[d]\\
\{N^\infty;S(\infty)^+\wedge_{\ZZ_2}(T(\nu_f) \wedge T(\nu_f))\}
\ar@<1ex>[r]^-{\di{\gamma}}
\ar@<-1ex>@{<-}[r]_-{\di{\delta}}&
\{N^\infty;T(\nu_f) \wedge T(\nu_f)\}_{\ZZ_2}
\ar@<1ex>[r]^-{\di{\rho}}
\ar@<-1ex>@{<-}[r]_-{\di{\sigma}}&
\{N^\infty;T(\nu_f)\}
}$$
with
$$\begin{array}{l}
P^2~=~\begin{pmatrix} 1\wedge P^2 & h_{\infty}(E) \\
0 & 1 \end{pmatrix}~:\\[2ex]
\{N^\infty;T(\nu_f)\}_{\ZZ_2}~=~\{N^\infty;S(\infty)\wedge_{\ZZ_2}T(\nu_f)\}\oplus
\{N^\infty;T(\nu_f)\}\to \\[1ex]
\{N^\infty;T(\nu_f)\wedge T(\nu_f)\}_{\ZZ_2}~=~\{N^\infty;S(\infty)\wedge_{\ZZ_2}(T(\nu_f)\wedge T(\nu_f))\}\oplus
\{N^\infty;T(\nu_f)\}~.
\end{array}$$
The image of the Umkehr map $F:\Sigma^{\infty}N^\infty \to \Sigma^{\infty}T(\nu_f)$
$$F~=~(0,F) \in \{N^\infty;T(\nu_f)\}_{\ZZ_2} ~=~
\{N^\infty;S(\infty)\wedge_{\ZZ_2}T(\nu_f)\}\oplus \{N^\infty;T(\nu_f)\}$$
is
$$\begin{array}{l}
P^2(F)~=~(h_{\infty}(F),F)~=~(D_2[f],F) \in\\[1ex]
\{N^\infty;T(\nu_f)\wedge T(\nu_f)\}_{\ZZ_2} ~=~
\{N^\infty;S(\infty)^+\wedge_{\ZZ_2}(T(\nu_f) \wedge T(\nu_f))\}\oplus
\{N^\infty;T(\nu_f)\}~,
\end{array}
$$
with
$$h_{\infty}(F)~=~D_2[f] \in \{N^\infty;S(\infty)^+\wedge_{\ZZ_2}(T(\nu_f) \wedge T(\nu_f))\}$$
a regular homotopy invariant of the immersion $f:N \imm M$, and
$$\begin{array}{l}
[F]~=~[h_{\infty}(F)]~=~[D_2[f]] \\[1ex]
\hskip25pt
\in {\rm coker}\big(P^2:\{N^\infty;T(\nu_f)\}_{\ZZ_2} \to
\{N^\infty;T(\nu_f) \wedge T(\nu_f))\}_{\ZZ_2}\big)\\[1ex]
\hskip50pt
=~{\rm coker}\big(1\wedge P^2:\{N^\infty;S(\infty)^+\wedge_{\ZZ_2}T(\nu_f)\}\to\\[1ex]
\hskip125pt
\{N^\infty;S(\infty)^+\wedge_{\ZZ_2}(T(\nu_f) \wedge T(\nu_f))\}\big)
\end{array}$$
a homotopy invariant of $f$.}\\
\hfill\qed
\end{remark}

For any manifold $A$, bundle $\gamma:B \to BO(k)$ and integer
$\ell \geqslant 0$ let $\Omega_{\ell}(A,B,\gamma)$ be the bordism
group of pairs
$$\hbox{( immersion $L^{\ell} \imm A$,
stable bundle map $(L,\nu_{L \imm A}) \to (B,\gamma)$ )~.}$$
The Pontrjagin-Thom construction identifies this bordism group with
a stable homotopy group
$$\Omega_{\ell}(A,B,\gamma)~=~\{A^{\infty},T(\gamma)\}~.$$

\begin{definition}~ \label{double12}
{\rm  The {\it double point Hopf invariant}
of an $(m,n,j)$-dimensional embedding-immersion pair
$(e \times f :M^m \emb V \times N^n,f:M \imm N)$
is the bordism class
$$(M',f',c) \in \Omega_{2m-n}(N,S(LV)\times_{\ZZ_2}(M \times M),e_2(\nu_f))$$
of the immersion of the
$(2m-n)$-dimensional unordered double point manifold
$$f'~:~M'~=~D_2[f] \imm N~;~[x,y] \mapsto f(x)=f(y)$$
with
$$c~:~M' \to S(LV) \times_{\ZZ_2}(M \times M)~;~[x,y] \mapsto
[{e(x)-e(y) \over \Vert e(x) - e(y) \Vert},x,y]~.\eqno{\hbox{\qed}}$$
}\end{definition}

\begin{proposition}~\label{double13}
The double point Hopf invariant is the stable homotopy
class of the quadratic construction {\rm (}\ref{quad1}{\rm )} on the
Umkehr map $F:V^{\infty} \wedge N^\infty \to V^{\infty} \wedge T(\nu_f)$
$$\psi_V(F)~:~N^\infty \xymatrix{\ar[r]^-{\di{F'}}&}
T(\nu_{f'}) \xymatrix@C+10pt{\ar[r]^-{\di{T(c)}}&}
S(LV)^+\wedge_{\ZZ_2}(T(\nu_f) \wedge T(\nu_f))$$
with $F':N^\infty \sto T(\nu_{f'})$ the Umkehr stable map of $f'$
$$\begin{array}{ll}
(M',f',c)~=~\psi_V(F) \in &
\Omega_{2m-n}(N,S(LV)\times_{\ZZ_2}(M \times M),e_2(\nu_f))\\[1ex]
&=~\{N^\infty,S(LV)^+ \wedge_{\ZZ_2} (T(\nu_f)\wedge T(\nu_f))\}~.
\end{array}$$
\hfill \qed
\end{proposition}

For any $j,k \geqslant 1$ the normal bundle of the
standard embedding of real projective spaces $P(\R^j) \subset
P(\R^{j+k})$ is the Whitney sum of $k$ copies of the canonical
line bundle $\lambda:P(\R^j) \to BO(1)$
$$\nu_{P(\R^j) \subset P(\R^{j+k})}~=~k\lambda~:~P(\R^j) \to BO(k)$$
with bundle projection
$$E(k\lambda)~=~S(L\R^j)\times_{\ZZ_2}L\R^k \to S(L\R^j)/\ZZ_2~=~P(\R^j)$$
and Thom space the stunted projective space
$$T(k\lambda)~=~(S(L\R^{j+k})/S(L\R^k))/\ZZ_2~=~P(\R^{j+k})/P(\R^k)~.$$
The bordism group of triples
$$\hbox{\rm ( $\ell$-dimensional manifold $L^{\ell}$, map $L \to P(\R^j)$,
stable bundle map $\nu_L \to k\lambda$ )}$$
is given by the Pontrjagin-Thom isomorphism to be
$$\Omega_{\ell}(P(\R^j),k\lambda)~=~\pi^S_{k+\ell}(P(\R^{j+k})/P(\R^k))~.$$
In particular, for $j=0$ this is
$$\Omega^{fr}_{\ell}~=~\pi^S_{\ell}~.$$
\indent The Stiefel manifold $V_{j,k}$ ($1\leqslant k \leqslant j$)
of orthonormal $k$-frames in $\R^j$ fits into a fibration
$$V_{j,k}~=~O(j)/O(j-k) \to BO(j-k) \to BO(j)~.$$
The canonical embedding
$$P(\R^j)/P(\R^{j-k}) \emb V_{j,k}~;~x=[x_1,x_2,\dots,x_j]
\mapsto [{\bf v}_{j-k+1},{\bf v}_{j-k+2},\dots,{\bf v}_j]$$
is defined using the columns ${\bf v}_{\ell}$ of the orthogonal
$j\times j$-matrix
$$(\delta_{pq}-2x_px_q)_{1 \leqslant  p,q \leqslant  j}~~
(\sum\limits^j_{i=1}(x_i)^2=1)$$
of the reflection in the hyperplane $x^{\perp} \subset \P(\R^j)$.
The pair $(V_{j,k},P(\R^j)/P(\R^{j-k}))$ is $2(j-k)$-connected
(James \cite[p.\  5]{imj5})).

\begin{remark} \label{double14}{\rm Let $N=S^n$. \\
(i) The double point Hopf invariant (\ref{double12}) of an
$(m,n,j)$-dimensional embedding-immersion pair
$$(e \times f :M^m \emb \R^j \times S^n,f:M \imm S^n)$$
is the bordism class of the $(2m-n)$-dimensional unordered double point
manifold
$$\begin{array}{l}
(f':M' \imm S^n)~=~\psi_{\R^j}(F\vert)\\[1ex]
\in \{S^n;(S^{j-1})^{\infty}\wedge_{\ZZ_2}(T(\nu_f) \wedge T(\nu_f))\}~
=~\Omega_{2m-n}(S^{j-1}\times_{\ZZ_2}(M \times M),e_2(\nu_f))
\end{array}$$
with
$$f'~:~M'~=~D_2[f]^{2m-n} \imm S^n~;~[x,y] \mapsto f(x)=f(y)$$
and $F\vert:S^{n+j} \to \Sigma^jT(\nu_f)$ the restriction of the
Umkehr map of $e \times f$
$$F~:~(\R^j \times S^n)^{\infty}~=~S^{n+j} \vee S^j\to \Sigma^jT(\nu_f)~.$$
(ii) The bordism class of an
$(m,n,j)$-dimensional embedding-immersion pair
$(e \times f :M^m \emb \R^j \times S^n,f:M \imm S^n)$ with
a framing $\delta\nu_f:\nu_f \cong \epsilon^{n-m}$ is
the homotopy class
$$(e \times f :M^m \emb \R^j \times S^n,f:M \imm S^n,\delta\nu_f)~=~
\phi \in \pi_{n+j}(S^{n-m+j})$$
of the Pontrjagin-Thom map
$$\begin{array}{l}
\phi~:~S^{n+j} \xymatrix{\ar[r]&} (\R^j \times S^n)^{\infty}~=~
S^{n+j} \vee S^j\\[1ex]
\hskip75pt \xymatrix{\ar[r]^-{\di{F}}&} \Sigma^jT(\nu_f)
\xymatrix@C+5pt{\ar[r]^-{\di{\Sigma^jT(\delta\nu_f)}}&}
\Sigma^{n-m+j}M^+ \to S^{n-m+j}~.
\end{array}$$
As in \ref{double10} (i) the normal bundle of the immersion
$f:M'=D_2[f]^{2m-n} \imm S^n$ is given by
$$\nu_{f'}~=~d^*(e_2(\epsilon^{n-m}))~:~M' \to BO(2(n-m))$$
with
$$d~:~M' \to P(\R^j)~;~[x,y] \mapsto
[{e(x) - e(y) \over \Vert e(x) - e(y) \Vert }]$$
and
$$\begin{array}{l}
e_2(\epsilon^{n-m})~=~(n-m)\lambda \oplus \epsilon^{n-m}~:\\[1ex]
S(L\R^j) \times_{\ZZ_2}(\{{\rm pt.}\}\times \{{\rm pt.}\})~=~
P(\R^j) \to BO(2(n-m))\end{array}$$
with
$$\begin{array}{l}
E(e_2(\epsilon^{n-m}))~=~S(L\R^j)\times_{\ZZ_2}(L\R^{n-m} \oplus \R^{n-m})~,\\[1ex]
T(e_2(\epsilon^{n-m}))~=~\Sigma^{n-m}(P(\R^{n-m+j})/P(\R^{n-m}))~.\end{array}$$
The double point Hopf invariant of Koschorke and
Sanderson \cite{ks1} is the bordism class
$$\begin{array}{l}
(f':M'\imm S^n)~=~\psi_{\R^j}(\phi)\\[1ex]
\hskip50pt \in \Omega_{2m-n}(P(\R^j),(n-m)\lambda)~=~
\pi^S_m(P(\R^{j+n-m})/P(\R^{n-m}))
\end{array}$$
of the quadratic construction
$$\psi_{\R^j}(\phi)~:~S^n \sto (S^{j-1})^{\infty}\wedge_{\ZZ_2}(S^{n-m} \wedge
S^{n-m})~=~\Sigma^{n-m}(P(\R^{j+n-m})/P(\R^{n-m}))~.$$
For $n \leqslant 2(n-m)$ the double point Hopf invariant map
$$\begin{array}{l}
H^j~:~\pi_{n+j}(S^{n-m+j}) \to \Omega_{2m-n}(P(\R^j),(n-m)\lambda)~=~
\pi_m(V_{n-m+j,j})~;\\[1ex]
(e \times f :M^m \emb \R^j \times S^n,f:M \imm S^n,\delta\nu_f) \mapsto
(f':M'\imm S^n)~=~\psi_{\R^j}(\phi)
\end{array}$$
fits into the $EHP$ exact sequence of James \cite{imj1}
$$\pi_n(S^{n-m}) \xymatrix{\ar[r]^-{\di{E^j}}&} \pi_{n+j}(S^{n-m+j})
\xymatrix{\ar[r]^-{\di{H^j}}&}
\pi_m(V_{n-m+j,j}) \xymatrix{\ar[r]^-{\di{P^j}}&}
\pi_{n-1}(S^{n-m}) \to \dots~,$$
with
$$\begin{array}{ll}
\pi_m(V_{n-m+j,j})&=~\pi_m(P(\R^{j+n-m})/P(\R^{n-m}))\\[1ex]
&=~\pi^S_m(P(\R^{j+n-m})/P(\R^{n-m}))~=~\Omega_{2m-n}(P(\R^j),(n-m)\lambda)~.
\end{array}$$
The $J$-homomorphism defines a natural transformation of exact sequences
$$\xymatrix@C-15pt{
\pi_m(SO(n-m)) \ar[r] \ar[d]^-{\di{J}} &
\pi_m(SO(j+n-m)) \ar[r]  \ar[d]^-{\di{J}}&
\pi_m(V_{j+n-m,j}) \ar[r] \ar@{=}[d] &
\pi_{m-1}(SO(n-m)) \ar[d]^-{\di{J}} \\
\pi_n(S^{n-m})\ar[r]^-{\di{E^j}} & \pi_{j+n}(S^{j+n-m})
\ar[r]^-{\di{H^j}}& \pi_m(V_{j+n-m,j}) \ar[r]^-{\di{P^j}} &
    \pi_{n-1}(S^{n-m})}$$
(iii) For $n \leqslant 2(n-m)$, $j=1$ the
exact sequence in (ii) is just the EHP sequence
$$\pi_n(S^{n-m}) \xymatrix{\ar[r]^-{\di{E}}&} \pi_{n+1}(S^{n-m+1})
\xymatrix{\ar[r]^-{\di{H}}&} \pi_m(S^{n-m})
    \xymatrix{\ar[r]^-{\di{P}}&}\pi_{n-1}(S^{n-m}) \to \dots~,$$
with $V_{1+n-m,1}=S^{n-m}$, $E$ the suspension map and $H$ the
classical Hopf invariant, interpreted as sending an
$(m,n,1)$-dimensional embedding-immersion pair $(e \times f :M^m
\emb \R \times S^n,f:M \imm S^n)$ with a framing
$\delta\nu_f:\nu_f \cong \epsilon^{n-m}$ to the double point Hopf
invariant
$$(f':M'\imm S^n)~=~\psi_{\R}(\phi)
 \in \Omega^{fr}_{2m-n}~=~\pi_m(V_{n-m+1,1})~=~\pi_m(S^{n-m})~=~\pi^S_{2m-n}~.$$
For $n=2m\geqslant 2$ the EHP sequence is
$$\pi_{2m}(S^m) \xymatrix{\ar[r]^-{\di{E}}&} \pi_{2m+1}(S^{m+1})
\xymatrix{\ar[r]^-{\di{H}}&} \pi_m(S^m)=\ZZ
    \xymatrix{\ar[r]^-{\di{P}}&}\pi_{2m-1}(S^m) \to \dots~.$$
In particular, for $m=1$, $n=2$  the generator $1 \in \pi_3(S^2)$
is the cobordism class of the $(1,2,1)$-dimensional embedding-immersion pair
$$(e \times f :S^1 \emb \R \times S^2,f:S^1 \imm S^2)$$
with $f$ the figure 8 immersion and $\delta\nu_f \cong \epsilon$ one of
the two framings. The Hopf invariant 1 element
$$(e \times f :S^1 \emb \R \times S^2,f:S^1 \imm S^2)~=~\eta~=~1 \in
\pi_3(S^2)~=~\ZZ$$
is detected by the cobordism class of the 0-dimensional
unordered double point manifold $M'=\{*\}$ of the immersion $f:M=S^1 \imm S^2$
$$(f':M'\imm S^2)~=~\psi_{\R}(\phi)~=~ 1  \in \Omega^{fr}_0~=~\pi^S_0~=~\ZZ~.$$
\hfill\qed}
\end{remark}

\section{Linking and self-linking}\label{link}

The original invariant of Hopf [10] was detected by the linking of
$k$-dimensional submanifolds of $S^{2k+1}$
$$\begin{array}{l}H~:~\pi_{2k+1}(S^{k+1}) \to \ZZ~;\\[1ex]
(\phi:S^{2k+1} \to S^{k+1}) \mapsto
\hbox{linking number}((M_1)^k \cup (M_2)^k \subset S^{2k+1})\end{array}$$
with $M_1=\phi^{-1}(x_1)$, $M_2=\phi^{-1}(x_2) \subset S^{2k+1}$
for distinct regular values $x_1,x_2 \in S^{k+1}$ of $\phi$.
The homotopy-theoretic construction of the generalized Hopf invariant
$H$ of G.W. Whitehead \cite{wh} was interpreted by Kervaire \cite{kerv3} and
Haefliger and Steer \cite[p.\ 262]{hs} in terms of higher linking numbers.
The Hopf invariant constructions $\lambda_2,\mu_2$ of Boardman and
Steer \cite{bs} (recalled in \S\S\ref{hopf},4 above) were motivated by the linking and
self-linking of embedded submanifolds.  We shall now interpret the
constructions in terms of immersed submanifolds.

Let $M_1,M_2,N$ be closed manifolds such that
$${\rm dim}\,M_1~=~m_1~~,~~{\rm dim}\,M_2~=~m_2~~,~~{\rm dim}\,N~=~n~,$$
and suppose given disjoint embeddings
$$e_i \times f_i~:~M_i \emb \R \times N~~(i=1,2)$$
which intersect transversely, with $f_i:M_i \imm  N$ immersions which intersect transversely.
The normal bundle of $e_i \times f_i:M_i \hookrightarrow \R \times N$ is
$$\nu_{e_i \times f_i}~=~\nu_{f_i}\oplus\epsilon~:~M_i \to BO(n-m_i+1)$$
with $\nu_{f_i}:M_i \to BO(n-m_i)$ the normal bundle of $f_i:M_i \imm N$.

\begin{definition}~ \label{link1}
{\rm  The {\it linking manifold} of the
submanifolds $M_1,M_2 \subset \R \times N$ is the
$(m_1+m_2-n)$-dimensional submanifold of $M_1 \times M_2$
$$\begin{array}{l}L(M_1,M_2,N)\\[1ex]
=~\{(x_1,x_2) \in M_1 \times M_2\st
e_1(x_1) < e_2(x_2) \in \R,f_1(x_1)=f_2(x_2) \in N\}~.\end{array}$$
\hfill\qed
}\end{definition}

\begin{proposition}~ \label{link2}
{\rm (i)} The linking manifold $L(M_1,M_2,N)$ is homeomorphic to the
intersection in $\R \times N \times I$ of the tracks of isotopies
unlinking $M_1,M_2$.\\
{\rm (ii)} The embedding
$$\begin{array}{l}g~:~L(M_1,M_2,N) \emb \R \times \R \times N~;\\[1ex]
(x_1,x_2) \mapsto (e_1(x_1),e_2(x_2),y)~~(y=f_1(x_1)=f_2(x_2))\end{array}$$
has normal bundle
$$\nu_g~=~
((\nu_{f_1}\oplus \epsilon) \times (\nu_{f_2}\oplus\epsilon))\vert~:~
L(M_1,M_2,N) \to BO(2n-m_1-m_2+2)~.$$
\end{proposition}
\noindent{\it Proof\  .}
(i) Assume that $e_1(M_1),e_2(M_2) \subset (0,1)$, and define isotopies
$$\begin{array}{l}
h_1~:~M_1\times I \to \R\times N \times I~;~(x,s) \mapsto (e_1(x)+s,f_1(x),s)~,\\[1ex]
h_2~:~M_2\times I \to \R\times N \times I~;~(y,t) \mapsto (e_2(y),f_2(y),t)\end{array}$$
such that the submanifolds
$$h_1(M_1 \times \{1\})~~,~~
h_2(M_2 \times \{1\}) \subset  \R \times N \times \{1\}$$
isotopic to the submanifolds
$$M_1~=~h_1(M_1 \times \{0\})~~,~~
M_2~=~h_2(M_2 \times \{0\}) \subset \R \times N \times \{0\}$$
are disjoint and unlinked, with
$$\begin{array}{l}
h_1(M_1 \times \{1\}) \subset
\{(a,z,1)\in \R \times N \times \{1\}\st a>1,z \in N\}~,\\[1ex]
h_2(M_2 \times \{1\}) \subset
\{(b,z,1)\in \R \times N \times \{1\}\st b<1,z \in N\}\end{array}$$
separated by  $\{1\} \times N \times \{1\}$. The tracks of the
isotopies are $(m_i+1)$-dimensional submanifolds
$$L_i~=~h_i(M_i \times I) \subset \R \times N \times I~~(i=1,2)$$
and there is defined a homeomorphism
$$\begin{array}{l}L(M_1,M_2,N) \to L_1 \cap L_2~;\\[1ex]
(x_1,x_2) \mapsto h_1(x_1,e_2(x_2)-e_1(x_1))~=~h_2(x_2,e_2(x_2)-e_1(x_1))~.\end{array}$$
(ii) The submanifolds $L_i \subset \R \times \R \times N$ $(i=1,2)$
intersect transversely, with
$$\nu_{L_i \subset \R \times \R \times N}~:~L_i \xymatrix{\ar[r]&}
M_i \xymatrix@C+20pt{\ar[r]^-{\di{\nu_{f_i} \oplus \epsilon}}&} BO(n-m_i+1)$$
so that
$$\begin{array}{l}
\nu_g~=~((\nu_{L_1\subset \R \times \R \times N})\times
(\nu_{L_2 \subset \R \times \R \times N}))\vert\\[1ex]
=~((\nu_{f_1}\oplus \epsilon) \times (\nu_{f_2}\oplus\epsilon))\vert~:~
L(M_1,M_2,N) \to BO(2n-m_1-m_2+2)~.\end{array}$$
\hfill\qed

Let
$$F_i~:~ \Sigma N^\infty~=~(\R \times N)^{\infty} \to
T(\nu_{e_i\times f_i})~=~\Sigma T(\nu_{f_i})$$
be the Umkehr stable maps given by the Pontrjagin-Thom construction.
The embedding of the disjoint union
$$e_1 \times f_1 \sqcup e_2 \times f_2~:~M_1\sqcup M_2 \emb \R \times N$$
has a compactification Umkehr stable map
$$F~:~\Sigma N^\infty \to \Sigma T(\nu_{f_1}) \vee  \Sigma T(\nu_{f_2})~;~
(t,x) \mapsto \begin{cases}F_1(2t,x)&\hbox{if $0 \leqslant  t \leqslant 1/2$}
\\[1ex] F_2(2t-1,x)&\hbox{if $1/2 \leqslant  t \leqslant  1$~.}
\end{cases}$$

\begin{proposition}~\label{link3}
{\rm (i)} The map of Boardman and Steer \cite{bs}
$$
\mu_2(F)~:~\Sigma^2 N^\infty \to \Sigma T(\nu_{f_1}) \wedge
\Sigma T(\nu_{f_2})~;~(s,t,x) \mapsto
\begin{cases}(F(s,x),F(t,x))&
\hbox{if $s \leqslant  t$}\\[1ex] *&\hbox{otherwise}\end{cases}$$ is
the composite
$$\mu_2(F)~=~T(i)G~:~\Sigma^2 N^\infty \xymatrix{\ar[r]^{\di{G}}&} T(\nu_g)
\xymatrix@C+10pt{\ar[r]^{\di{T(i)}}&}
\Sigma T(\nu_{f_1}) \wedge \Sigma T(\nu_{f_2})$$
of the Umkehr stable map for $g:L(M_1,M_2,N) \emb \R\times \R \times N$
$$G~:~(\R \times \R \times N)^{\infty}~=~\Sigma^2 N^\infty \to T(\nu_g)$$
and the inclusion of Thom spaces
$$T(i)~:~T(\nu_g) \to T((\nu_{f_1}\oplus\epsilon) \times
(\nu_{f_2}\oplus \epsilon))~=~
\Sigma T(\nu_{f_1}) \wedge \Sigma T(\nu_{f_2})$$
induced by the inclusion
$$i~:~L(M_1,M_2,N) \emb M_1 \times M_2~;~(x_1,x_2) \mapsto (x_1,x_2)~.$$
{\rm (ii)}
The image of the fundamental class $[N] \in H_n(N)$ under $\mu_2(F)$ is
$$\mu_2(F)_*[N]~=~i_*[L] \in \widetilde{H}_n(T(\nu_{f_1})\wedge T(\nu_{f_2}))~
=~H_{m_1+m_2-n}(M_1\times M_2)$$
with $[L] \in H_{m_1+m_2-n}(L)$ the fundamental class of $L=L(M_1,M_2,N)$.
Viewing $i_*[L] \in H_{m_1+m_2-n}(M_1\times M_2)$ as a chain map
$$i_*[L]~:~C(M_1)^{m_1-*} \to C(M_2)_{*-n+m_2}$$
(up to chain homotopy) there is defined a commutative diagram
$$\xymatrix@R+10pt{
H^{m_1-*}(M_1) \ar[r]^-{\di{i_*[L]}}\ar[d]^-{\di{\cong}}_-{\di{[M_1]\cap -}}&
H_{*-n+m_2}(M_2)~\cong~H_{*+1}(\R \times N,\R \times N\backslash
(e_2 \times f_2)(M_2))\ar[d] \\
H_*(M_1)\ar[r]^-{\di{(e_1 \times f_1)_*}} & H_*(\R \times N\backslash
(e_2 \times f_2)(M_2))~.}$$
\hfill\qed
\end{proposition}

\begin{example} \label{link4}
{\rm The linking manifold $L=L(M_1,M_2,S^n)$ of
$M_1,M_2 \subset \R \times S^n$ with $m_1+m_2=n$ is 0-dimensional, and
$$\mu_2(F)_*[N]~=~i_*[L] \in H_0(M_1\times M_2)~=~\ZZ$$
may be identified with
$$\begin{array}{l}\hbox{\rm linking number}(M_1 \cup M_2 \subset \R \times S^n)~=~
(e_1 \times f_1)_*[M_1]\\[1ex]
\hskip75pt
\in H_{m_1}((\R \times S^n) \backslash (e_2 \times f_2)(M_2))~=~
H^{m_2}(M_2)~=~\ZZ~.\end{array}$$
\hfill\qed}
\end{example}

\begin{proposition}~ \label{link5} {\rm (Boardman and Steer \cite{bs})}
Let $N$ be an $n$-dimensional manifold. The Hopf map
$$H~:~[\Sigma N^\infty,S^{k+1}] \to [\Sigma^2 N^\infty,S^{2k+2}]~;~
\phi \mapsto h_{\R}(\phi)$$
sends $\phi:\Sigma N^\infty \to S^{k+1}$ to the cobordism class of
the framed submanifold
$$L(M_1,M_2,N)^{n-2k} \subset \R\times \R \times N$$
with
$$M_i~=~\phi^{-1}(x_i) \subset \R \times N~~(i=1,2)$$
for distinct regular values $x_1,x_2 \in S^{k+1}$, regarding
$[\Sigma^2 N^\infty,S^{2k+2}]$ as the bordism
group of framed embeddings $L^{n-2k} \subset \R\times \R\times N$.
\end{proposition}
\noindent{\it Proof\  .} It may be assumed that the
inclusions are framed embeddings $e_i \times f_i:M_i \emb
\R\times N$ with $f_i:M_i \imm N$ immersions. By \ref{link3}
$$H(\phi)~:~\Sigma^2 N^\infty \xymatrix@C+10pt{\ar[r]^-{\di{\mu_2(F)}}&}
\Sigma T(\nu_{f_1}) \wedge \Sigma T(\nu_{f_2})~=~
\Sigma^{2k+2}(M_1 \times M_2)^{\infty}$$
is such that
$$H(\phi)^{-1}(M_1 \times M_2)~=~L(M_1,M_2,N)^{n-2k} \subset \R \times \R \times N~.$$
It follows that
$$h_{\R}(\phi)~:~\Sigma^2N^\infty \xymatrix@C+10pt{\ar[r]^-{\di{H(\phi)}}&}
\Sigma^{2k+2}(M_1 \times M_2)^+\to S^{2k+2}$$
is such that
$$h_{\R}(\phi)^{-1}(*)~=~L(M_1,M_2,N)^{n-2k} \subset \R\times \R \times N~.$$
\hfill\qed

\begin{example} \label{link6} {\rm
(i) For $N=S^{2k}$ \ref{link5} gives the original invariant of Hopf \cite{hopf2}
$$\begin{array}{ll}
H~:&\pi_{2k+1}(S^{k+1}) \to \pi_{2k+2}(S^{2k+2})~=~\ZZ~;\\[1ex]
&\phi \mapsto H(\phi)~=~\hbox{\rm linking number}
((M_1)^k \cup (M_2)^k \subset S^{2k+1})\end{array}$$
with $(M_i)^k=\phi^{-1}(x_i)$.\\
(ii) For $N=S^n$ \ref{link5} gives the generalized Hopf invariant
$$H~:~\pi_{n+1}(S^{k+1}) \to \pi_{n+2}(S^{2k+2})~;~\phi \mapsto H(\phi)~.$$
of G.W. Whitehead \cite{wh} and  Hilton \cite{hilton}, with
$$\begin{array}{ll}
H(\phi)~=&
\hbox{framed cobordism class}(L(M_1,M_2,S^n)^{n-2k} \subset S^{n+2})
\in  \pi_{n+2}(S^{2k+2})\\[1ex]
&((M_i)^{n-k}~=~ \phi^{-1}(x_i)\subset \R \times S^n \subset
S^{n+1}~~(i=1,2))\end{array}$$
as in Kervaire \cite{kerv3}.\hfill\qed}
\end{example}

\begin{definition}~ \label{link7}
{\rm  Let $e \times f:M \emb \R \times N$
be an embedding of an $m$-dimensional manifold $M$, with $N$ an
$n$-dimensional manifold and $f:M \imm N$ an immersion.
The {\it self-linking manifold} is the $(2m-n)$-dimensional manifold
$$L(M,N)~=~\{(x,y) \in M \times M\st e(x) < e(y)\in \R, f(x)=f(y)\in N \}~.$$
\hfill\qed
}\end{definition}

\begin{proposition}~\label{link8}
{\rm (i)} The ordered double point manifold
$$D_2(f)~=~\{(x,y) \in M \times M\st x \neq y\in M, f(x)=f(y)\in N \}$$
is the disjoint union of the self-linking manifold and its transpose
$$D_2(f)~=~L(M,N) \sqcup TL(M,N)~,$$
with
$$TL(M,N)~=~\{(x,y) \in M \times M\st e(x) > e(y)\in \R, f(x)=f(y)\in N \}~.$$
{\rm (ii)} The embedding
$$g~:~L(M,N) \emb \R \times \R \times N~;~(x,y) \mapsto (e(x),e(y),f(x))$$
has normal bundle
$$\begin{array}{ll}
\nu_g&=~((\nu_f\oplus \epsilon) \times (\nu_f\oplus \epsilon))\vert\\[1ex]
&=~((\nu_f\times \nu_f)\oplus \epsilon^2)\vert ~:~L(M,N) \to BO(2n-2m+2)\end{array}$$
with $\nu_f:M \to BO(n-m)$ the normal bundle of $f$. \hfil\break
{\rm (iii)} The immersion
$$g_2~:~L(M,N) \imm N~;~(x,y) \mapsto f(x)$$
has normal bundle
$$\nu_{g_2}~=~(\nu_f \times \nu_f)\vert ~:~ L(M,N) \to BO(2n-2m)~.$$
{\rm (iv)} The immersion $f:M \imm N$ is an embedding if and only
if $L(M,N)=\emptyset$.\hfill\qed
\end{proposition}

Let $F:\Sigma N^\infty \to \Sigma T(\nu_f)$ be the Umkehr stable map of $e\times f$.

\begin{proposition}~\label{link9}
{\rm (i)} The Hopf invariant map of Boardman and Steer \cite{bs}
$$\lambda_2(F)~=~\mu_2(\nabla F)~:~\Sigma^2N^\infty \to \Sigma T(\nu_f)\wedge \Sigma T(\nu_f)$$
(identified with $h_{\R}(F)$ in \ref{hopf3}) is the composite
$$h_{\R}(F)~=~\lambda_2(F)~=~T(i)G~:~\Sigma^2 N^\infty \xymatrix{\ar[r]^{\di{G}}&} T(\nu_g)
\xymatrix@C+10pt{\ar[r]^{\di{T(i)}}&} \Sigma T(\nu_f) \wedge \Sigma T(\nu_f)$$
of the Umkehr stable map of $g$
$$G~:~(\R \times \R \times N)^{\infty}~=~\Sigma^2 N^\infty \to T(\nu_g)$$
and the inclusion of Thom spaces
$$T(i)~:~T(\nu_g) \to T((\nu_f \oplus \epsilon)\times (\nu_f\oplus\epsilon))~=~
\Sigma T(\nu_f) \wedge \Sigma T(\nu_f)$$
induced by the inclusion
$$i~:~L(M,N) \emb M \times M~;~(x_1,x_2) \mapsto (x_1,x_2)~.$$
{\rm (ii)} The double point Hopf invariant {\rm (}\ref{double12}{\rm )} is the bordism
class of the immersion of the self-linking manifold $g_2:L(M,N) \imm N^n$,
which is the stable homotopy class of the ultraquadratic construction
(\ref{ultraquadratic})
$$\psi_{\R}(F)~:~N^\infty \xymatrix{\ar[r]^-{\di{G_2}}&}
T(\nu_{g_2}) \xymatrix@C+10pt{\ar[r]^-{\di{T(i)}}&} T(\nu_f) \wedge T(\nu_f)$$
with $G_2$ the Umkehr stable map of $g_2$.\\
{\rm (iii)} The image of the fundamental class $[N] \in H_n(N)$ under $h_{\R}(F)$
is
$$h_{\R}(F)_*[N]~=~i_*[L(M,N)] \in \widetilde{H}_n(T(\nu_f)\wedge T(\nu_f))~=~
H_{2m-n}(M\times M)$$ with $[L(M,N)] \in H_{2m-n}(L(M,N))$ the
fundamental class of $L(M,N)$.
\end{proposition}
\noindent{\it Proof\  .} (i) Apply \ref{link3} with
$$M_1~=~M_2~=~M~~,~~f_1~=~f_2~=~f~:~M \to N~~,~~e_1~=~e+c~~,~~e_2~=~e$$
for some $c>0$ so small that
$$c < e(y)-e(x)~~((x,y) \in L(M,N))~,$$
in which case
$$e_1 \times f_1~:~M_1 \emb \R \times N~~,~~e_2 \times f_2~:~M_2 \emb \R \times N$$
are disjoint embeddings, and the linking manifold is
$$L(M_1,M_2,N)~=~L(M,N)~.$$
The embedding of the disjoint union
$$e_1 \times f \sqcup e_2 \times f~:~M \sqcup M \emb \R \times N$$
has compactification Umkehr stable map
$$\nabla F~:~\Sigma N^\infty \to \Sigma T(\nu_f) \vee  \Sigma T(\nu_f)~;~
(t,x) \mapsto \begin{cases}
F(2t,x)_1&\hbox{if $0 \leqslant  t \leqslant 1/2$}
\\[1ex] F(2t-1,x)_2&\hbox{if $1/2 \leqslant  t \leqslant  1$}
\end{cases}$$ with
projections
$$\pi_i (\nabla F)~ \simeq~F~:~ \Sigma N^\infty \to \Sigma T(\nu_f)~.$$
(ii) This is just the special case $j=1$ of \ref{double13} (i).\\
(iii) Immediate from (ii).\hfill\qed

The fibration
$$V_{k+1,1}~=~S^k \xymatrix@C+20pt{\ar[r]^-{\di{\tau_{S^k}}}&} BO(k)
\xymatrix{\ar[r]&} BO(k+1)$$
classifies the tangent bundle of $S^k$ and the stable trivialization
$$\delta\tau_{S^k}~:~\tau_{S^k}\oplus\epsilon~\cong~\epsilon^{k+1}$$
determined by the embedding $S^k \emb S^{k+1}$.

\begin{lemma} \label{link10} For any space $M$ there is a natural
one-one correspondence between the equivalence classes of pairs
$$\hbox{\rm ( \it $k$-plane bundle $\xi$ over $M$, stable isomorphism
$\delta\xi:\xi \oplus \epsilon \cong \epsilon^{k+1}$\rm )}$$
and the homotopy classes of maps $\gamma:M \to S^k$.
\end{lemma}
\begin{proof}
A map $\gamma:M \to V_{k+1,1}=S^k$ classifies a $k$-plane bundle
$$\xi~=~\gamma^*\tau_{S^k}~:~M \to BO(k)$$
together with a stable isomorphism
$$\delta\xi~=~\gamma^*\delta\tau_{S^k}~:~
\xi \oplus \epsilon ~\cong~ \epsilon^{k+1}~.$$
Conversely, given $(\xi,\delta\xi)$ the composite of the non-zero section
$$M \to E(\xi \oplus \epsilon)\backslash M~;~x \mapsto (x,(0,1))$$
and the map
$$E(\xi \oplus \epsilon)\backslash M \xymatrix{\ar[r]_-{\di{\cong}}^-{\di{\delta\xi}}&}
E(\epsilon^{k+1})\backslash M~=~
M \times (\R^{k+1} \backslash \{0\}) \to \R^{k+1} \backslash \{0\} ~
\simeq~S^k$$
is a map $\gamma:M \to S^k$ classifying $\xi$ and $\delta\xi$.\\
\hfill\qed\end{proof}
The Hopf invariant of
$$T(\delta\xi)~:~T(\xi \oplus \epsilon)~=~\Sigma T(\xi) \to
T(\epsilon^{k+1})~=~\Sigma^{k+1}M^+$$
can be regarded (as in \ref{hopf4}) as an inequivariant map
$$h_{\R}(T(\delta\xi))~:~\Sigma^2T(\xi)
\to \Sigma T(\epsilon^k)\wedge\Sigma T(\epsilon^k)~,$$
which has the following geometric interpretation for a manifold $M$.

\begin{definition}~\label{link11}
{\rm  Let $M$ be an $m$-dimensional manifold,
together with a $k$-plane bundle $\xi:M \to BO(k)$ and a stable
trivialization $\delta\xi:\xi \oplus \epsilon \cong \epsilon^{k+1}$.
A {\it linking manifold} for $(\xi,\delta\xi)$ is a framed
codimension $k$ submanifold
$$L(M,\xi,\delta\xi)^{m-k}~=~\gamma^{-1}(*) \subset M$$
with $\gamma:M \to S^k$ a classifying map for $(\xi,\delta\xi)$,
and $* \in S^k$ a regular value of $\gamma$.\hfill\qed
}\end{definition}

The homotopy class $\gamma \in [M,S^k]$ may be identified with the framed
bordism class of $L(M,\xi,\delta\xi)^{m-k} \subset M$.

\begin{proposition}~\label{link12}
{\rm (i)} A linking manifold for $(\xi,\delta\xi)$ is the linking manifold
$$L(M,\xi,\delta\xi)~=~L(M_1,M_2,M)$$
for the disjoint embeddings
$$\begin{array}{l}
e_1\times f_1~:~M_1~=~M \emb \R \times M~;~
x \mapsto (0,(0,x))~,\\[1ex]
e_2\times f_2~:~M_2~=~M \emb \R \times E(\xi)~;~
x \mapsto \delta\xi^{-1}(*,x)~,\end{array}$$
with
$$\delta\xi^{-1}~:~E(\epsilon^{k+1})~=~\R^{k+1} \times M \to
E(\xi \oplus \epsilon)~=~\R \times E(\xi)~.$$
{\rm (ii)} The bordism class of the immersion of the linking
manifold  $L(M,\xi,\delta\xi) \imm E(\xi)$
is the stable homotopy class of the Hopf invariant map
$$h_{\R}(T(\delta\xi))~:~\Sigma^2T(\xi) \to \Sigma^{k+1}M^+ \wedge
\Sigma^{k+1}M^+$$
with
$$T(\delta\xi)~:~T(\xi \oplus \epsilon)~=~\Sigma T(\xi) \to
T(\epsilon^{k+1})~=~\Sigma^{k+1} M^+~.$$
{\rm (iii)} The Hopf invariant of $T(\delta\xi)$ is given by
$$\begin{array}{l}
h_{\R}(T(\delta\xi))~:~\Sigma^2 T(\xi)~\cong~\Sigma^{k+2}M^+
\xymatrix@C+30pt{\ar[r]^-{\di{\Sigma^{k+2}\ell}}&}
\Sigma^{2k+2}L(M,\xi,\delta\xi)^{\infty} \\[1ex]
\hskip150pt \xymatrix@C+30pt{\ar[r]^-{\di{\Sigma^{2k+2}i}}&}
\Sigma^{2k+2}(M^+ \wedge M^+)\end{array}$$
with $\ell:M^+ \to \Sigma^kL(M,\xi,\delta\xi)^{\infty}$ a
Pontrjagin-Thom map for a linking manifold $L(M,\xi,\delta\xi)
\subset M$, and
$$i~:~L(M,\xi,\delta\xi) \emb M \times M~;~x \mapsto (x,x)~.$$
\hfill\qed
\end{proposition}

\begin{proposition}~\label{link13}
Let $(e\times f:M^m \emb \R \times N^n,f:M \imm N)$ be
an $(m,n,1)$-dimensional embedding-immersion pair with $e \times
f$ framed, so that $\nu_f:M \to BO(n-m)$ has a stable
trivialization
$$\delta\nu_f~:~\nu_f\oplus\epsilon~\cong~\epsilon^{n-m+1}~.$$
The Hopf invariant of the composite map
$$\phi~:~\Sigma N^\infty \xymatrix{\ar[r]^-{\di{F}}&}
\Sigma T(\nu_f)\xymatrix@C+20pt{\ar[r]^-{\di{T(\delta\nu_f)}}&}
T(\epsilon^{n-m+1})$$
is the sum
$$\begin{array}{l}
h_{\R}(\phi)~=~
(T(\delta\nu_f)\wedge T(\delta\nu_f))h_{\R}(F)+h_{\R}(T(\delta\nu_f))(\Sigma F)~:\\[1ex]
\hskip50pt
\Sigma^2N^\infty \to T(\epsilon^{n-m+1})\wedge T(\epsilon^{n-m+1})\end{array}$$
with
$$h_{\R}(F)~:~\Sigma^2N^\infty \to \Sigma T(\nu_f)\wedge \Sigma T(\nu_f)$$
the Pontrjagin-Thom map for the self-linking manifold
$L(M,N)\emb \R \times \R \times N$
and
$$\begin{array}{l}
h_{\R}(T(\delta\nu_f))~:~
\Sigma^2 T(\nu_f)~\cong~\Sigma^{n-m+2}M^+\\[1ex]
\xymatrix@C+30pt{\ar[r]^-{\di{\Sigma^{n-m+2}\ell}}&}
\Sigma^{2n-2m+2}L(M,\nu_f,\delta\nu_f)^{\infty}
\xymatrix@C+30pt{\ar[r]^-{\di{\Sigma^{2n-2m+2}i}}&}
\Sigma^{2n-2m+2}(M^+ \wedge M^+)
\end{array}$$
with $\ell:M^+ \to \Sigma^{n-m}L(M,\nu_f,\delta\nu_f)^+$
the Pontrjagin-Thom map for the linking manifold
$L(M,\nu_f,\delta\nu_f) \subset M$ and
$$i~:~L(M,\nu_f,\delta\nu_f) \emb M \times M~;~x \mapsto (x,x)~.$$
\hfill\qed
\end{proposition}

\begin{example} \label{link14}
{\rm (i) Given an $n$-dimensional manifold $N$ and a map
$$\phi~:~(\R \times N)^{\infty}~=~\Sigma N^\infty \to S^{k+1}$$
apply the Pontrjagin-Thom construction to obtain
an $(n-k,n,1)$-dimensional embedding-immersion pair
$$(e \times f:M^{n-k}=\phi^{-1}(*) \emb \R\times N~,~f:M \imm N)$$
with $e \times f$ framed, so that the normal bundle $\nu_f:M \to BO(k)$
has a stable trivialization
$$\delta\nu_f~:~\nu_f \oplus \epsilon~ \cong~ \epsilon^{k+1}~.$$
If $F:\Sigma N^\infty \to \Sigma T(\nu_f)$ is the Umkehr stable map then
$$\phi~:~\Sigma N^\infty \xymatrix{\ar[r]^-{\di{F}}&}
\Sigma T(\nu_f)~\cong~\Sigma^{k+1}M^+ \to S^{k+1}~.$$
By \ref{double6} the Hopf invariant map
$$h_{\R}(\phi)~:~\Sigma^2 N^\infty\to S^{2k+2}$$
is such that
$$h_{\R}(\phi)^{-1}(*)^{n-k}~=~L(M,N) \cup L(M,\nu_f,\delta\nu_f)
\subset \R \times \R \times N~.$$
(ii) Setting $N=S^n$ in (i) shows that the Hopf map (\ref{link6})
$$H~:~\pi_{n+1}(S^{k+1}) \to \pi_{n+2}(S^{2k+2})~;~
\phi \mapsto h_{\R}(\phi)$$
sends the bordism class of the $(n-k,n,1)$-dimensional embedding-immersion pair
$$(e\times f:M^{n-k}=\phi^{-1}(*) \emb \R \times S^n~,~
f:M \imm S^n)$$
with $e \times f$ framed to the bordism class of the framed submanifold
$$L(M,S^n)^{n-2k}\cup L(M,\nu_f,\delta\nu_f)^{n-2k} \subset S^{n+2}~.$$
\hfill\qed}
\end{example}

\begin{example} \label{link15} {\rm The homotopy group
$\pi_{m+k}(T(\tau_{S^k}))$ is the cobordism group of embeddings
$f:M^m \emb S^{m+k}$ which are framed in $S^{m+k+1}$. A
map $\rho:S^{m+k} \to T(\tau_{S^k})$ which is transverse at the
zero section $S^k\subset T(\tau_{S^k})$ determines an embedding
$$f~:~M^m~=~\rho^{-1}(S^k) \emb S^{m+k}$$
which is framed in $S^{m+k+1}$, with a map $\gamma=\rho\vert :M \to S^k$
such that
$$\nu_f~=~\gamma^*\tau_{S^k}~:~M \to BO(k)~,$$
with a stable trivialization
$$\delta\nu_f~=~\gamma^*\delta\tau_{S^k}~:~
\nu_f \oplus \epsilon~\cong~\epsilon^{k+1}~.$$
The linking manifold construction  (\ref{link11}) gives rise to the commutative square
of Wood \cite{wood}
$$\xymatrix@C+10pt@R+10pt{\pi_{m+k}(T(\tau_{S^k})) \ar[r]^-{\di{h}} \ar[d] &
\pi_{m+k}(S^{2k}) \ar[d]^-{\di{E}} \cr
\pi_{m+k+1}(S^{k+1}) \ar[r]^-{\di{H}} & \pi_{m+k+1}(S^{2k+1})}$$
with
$$\begin{array}{l}
h~:~\pi_{m+k}(T(\tau_{S^k})) \to \pi_{m+k}(S^{2k})~;~\rho \mapsto
(L(M,\nu_f,\delta\nu_f)^{m-k} \subset S^{m+k})~,\\[1ex]
\pi_{m+k}(T(\tau_{S^k})) \to \pi_{m+k+1}(S^{k+1})~;~\rho \mapsto (M \subset S^{m+k+1})~.\end{array}$$
The Thom space of $\tau_{S^k}:S^k \to BO(k)$ is
$$T(\tau_{S^k})~=~S^k \cup_{J(\tau_{S^k})}e^{2k}$$
with $J(\tau_{S^k}):S^{2k-1} \to S^k$ given by the $J$-homomorphism
$$J~:~\pi_k(BO(k))~=~\pi_{k-1}(O(k)) \to \pi_{2k-1}(S^k)~,$$
and $h$ is induced by $T(\tau_{S^k}) \to S^{2k}$. The Hopf map $H$ is given by
$$\begin{array}{l}
H~:~\pi_{m+k+1}(S^{k+1}) \to \pi_{m+k+1}(S^{2k+1})~;\\[1ex]
(M,e,f) \mapsto (L(M,S^{m+k})^{m-k}\subset S^{m+k+1})+
(L(M,\nu_f,\delta\nu_f)^{m-k}\subset S^{m+k+1})~,\end{array}$$
regarding $\pi_{m+k+1}(S^{k+1})$ as the cobordism group of
$(m,m+k,1)$-dimensional embedding-immersion pairs $(e\times f:M^m
\emb \R \times S^{m+k},f:M \imm S^{m+k})$ such that $e
\times f$ is framed.\hfill\qed}
\end{example}

\begin{example}\label{link16} {\rm Let $m \leqslant  2k-2$, so that
$$\pi_{m+k+1}(S^{2k+1})~=~\pi_m(S^k)$$
by the Freudenthal suspension theorem; equivalently, every
(framed) submanifold $L^{m-k} \subset S^{m+k+1}$ can be compressed
to $L^{m-k} \subset S^m$. \\
(i) Regard $\pi_{m+k+1}(S^{k+1})$ as the cobordism group of
$(m,m+k,1)$-dimensional embedding-immersion pairs
$$(e\times f:M^m \emb \R \times S^{m+k}~,~f:M \imm S^{m+k})$$
with $e \times f$ framed and $f$ an embedding, so that
$$(L(M,S^{m+k})^{m-k} \subset S^{m+k+1})~=~0 \in \pi_{m+k+1}(S^{2k+1})$$
and
$$\begin{array}{l}H~:~\pi_{m+k+1}(S^{k+1}) \to \pi_{m+k+1}(S^{2k+1})~;\\[1ex]
\hskip50pt
(M,e,f) \mapsto (L(M,\nu_f,\delta\nu_f)^{m-k} \subset S^{m+k+1})\end{array}$$
is the `singularity Hopf invariant' of Koschorke and Sanderson
\cite[p.\  201]{ks1}.\\
(ii) Regard $\pi_{m+k+1}(S^{k+1})$ as the cobordism group of
$(m,m+k,1)$-dimensional embedding-immersion pairs
$$(e\times f:M^m \emb \R \times S^{m+k}~,~f:M \imm S^{m+k})$$
with a framing $\nu_f \cong \epsilon^k$, so that
$$(L(M,\nu_f,\delta\nu_f)^{m-k} \subset S^{m+k+1})~=~0 \in
\pi_{m+k+1}(S^{2k+1})$$
and
$$\begin{array}{l}H~:~\pi_{m+k+1}(S^{k+1}) \to \pi_{m+k+1}(S^{2k+1})~;\\[1ex]
\hskip50pt (M,e,f) \mapsto (L(M,S^{m+k})^{m-k} \subset
S^{m+k+1})\end{array}$$ is the `double point Hopf invariant' of
\cite[p. 202]{ks1} (and \ref{double14} (ii)).\hfill\qed}
\end{example}

\section{Intersections and self-intersections for $M^m \to N^{2m}$}\label{self}

We now consider the application of the geometric Hopf invariant to the
intersection and self-intersection properties of embeddings and
immersions $f:M \to N$ of an $m$-dimensional manifold $M^m$ in a
$2m$-dimensional manifold $N^{2m}$.  It will be assumed that $M$ and
$N$ are oriented, so that the normal bundle is an oriented $m$-plane
bundle $\nu_f:M \to BSO(m)$.

\begin{definition}~ \label{embed1} {\rm
(i) The {\it intersection pairing} of a $2m$-dimensional manifold $N$ is
the $(-1)^m$-symmetric cup product pairing
$$\lambda~:~H_m(N) \times H_m(N) \to \ZZ~;~
(a,b) \mapsto \langle a^*\cup b^*,[N] \rangle$$
with $a^*,b^* \in H^m(N)$ the Poincar\'e duals of $a,b \in H_m(N)$.\\
(ii) Let $(x_1,x_2) \in M_1 \times M_2$ be a transverse intersection
point of embeddings
$$f_1~:~(M_1)^m \emb N^{2m}~~,~~f_2:(M_2)^m \emb N^{2m}$$
so that
there is defined an isomorphism of $2m$-dimensional vector spaces
$$\tau_{M_1}(x_1)\oplus \tau_{M_2}(x_2)~\cong~\tau_N(y)$$
with
$$y~=~f_1(x_1)~=~f_2(x_2) \in N~.$$
Assuming $M_1,M_2,N$ are oriented the {\it intersection index}
of $(x_1,x_2)$ is
$$I(x_1,x_2)~=~\begin{cases}+1&
\hbox{if $\tau_{M_1}(x_1)\oplus \tau_{M_2}(x_2)\cong
\tau_N(y)$ is orientation-preserving}\\[1ex]
-1&\hbox{otherwise}
\end{cases}$$
where $y=f_1(x_1)=f_2(x_2) \in N$. \hfill\qed
}\end{definition}

\begin{proposition}~\label{embed2}
{\rm (i)} The intersection pairing
$\lambda:H_m(N) \times H_m(N) \to \ZZ$ of a $2m$-dimensional manifold $N$
is such that
$$\lambda(a,b)~=~(-1)^m\lambda(b,a) \in \ZZ~.$$
{\rm (ii)} If $f_1:(M_1)^m\emb N^{2m}$, $f_2:(M_2)^m
\emb N^{2m}$ are transverse embeddings of oriented
manifolds
    $$\lambda(f_1[M_1],f_2[M_2])~=~\sum\limits_{f_1(x_1)=f_2(x_2)}
    I(x_1,x_2) \in \ZZ~.\eqno{\hbox{\qed}}$$
\end{proposition}

The homotopy exact sequence of the fibration
$$V_{\infty,\infty-m} \to BO(m) \to BO$$
is mapped by the $J$-homomorphism to the stable $EHP$ sequence of $S^m$
$$\xymatrix@R+10pt@C+20pt{
\pi_{m+1}(BO) \ar[r] \ar[d]^-{\di{J}} &
\pi_m(V_{\infty,\infty-m}) \ar[r] \ar[d]^-{\di{\cong}} &
\pi_m(BO(m)) \ar[r] \ar[d]^-{\di{J}} &
\pi_m(BO) \ar[d]^-{\di{J}} \\
\pi_m^S \ar[r]^-{\di{H}} & Q_{(-1)^m}(\ZZ) \ar[r]&
\pi_{2m-1}(S^m) \ar[r] &
\pi_{m-1}^S}$$
The isomorphism
$$\pi_m(V_{\infty,\infty-m}) \xymatrix{\ar[r]^-{\di{\cong}}&}
\Omega_0(P(\R(\infty)),m\lambda)~=~Q_{(-1)^m}(\ZZ)~;~
(\xi,\delta\xi) \mapsto D_2(S^m,\xi,\delta\xi)/\ZZ_2$$
is defined by the double point Hopf invariant (\ref{double12}), and
$$\pi_{m+1}(BO,BO(m))~=~Q_{(-1)^m}(\ZZ) \to \pi_m(BO(m))~;~
1 \mapsto \tau_{S^m}~.$$
The Euler number defines a morphism
$$\chi~:~\pi_m(BSO(m)) \to \ZZ~;~\xi \mapsto \chi(\xi)$$
such that
$$\begin{array}{l}\pi_{m+1}(BO,BO(m))~
=~Q_{(-1)^m}(\ZZ) \xymatrix{\ar[r]&} \pi_m(BSO(m))
\xymatrix{\ar[r]^{\di{\chi}}&} \ZZ~;\\[1ex]
\hskip150pt 1 \mapsto  \chi(\tau_{S^m})~=~1+(-1)^m~.\end{array}$$

\begin{proposition}~\label{embed3}
For any embedding $f:M^m \emb N^{2m}$
$$\lambda(f_*[M^m],f_*[M^m])~=~\chi(\nu_f) \in \ZZ~.$$
\end{proposition}
\noindent{\it Proof\  .}
There exists an isotopic embedding $f':M^m \emb N^{2m}$ which intersects
transversely with $f$ in $\chi(\nu_f)$ points (counted algebraically), and
$$\lambda(f_*[M],f_*[M])~=~\sum\limits_{(x,x') \in D_2(f,f')}I(x,x')~
=~\chi(\nu_f) \in \ZZ~.$$
\hfill\qed

\begin{example}\label{embed4} {\rm If $N^{4k}$ is a $(2k-1)$-connected
$4k$-dimensional manifold and $k \geqslant 2$ then every element
$x \in H_{2k}(N)$ is represented by an embedding $f:S^{2k}
\emb N$, with
$$\lambda(x,x)~=~\chi(\nu_f) \in \ZZ$$
as in Wall \cite{wall1}.\hfill\qed}\end{example}

Recall the definition of the double point set (\ref{double1}) of a map $f:M \to N$
$$D_2(f)~=~\{(x,y) \in M \times M\st x \neq y \in M,f(x)=f(y) \in N\}~.$$
Transposition defines a free $\ZZ_2$-action
$$T~:~D_2(f) \to D_2(f)~;~(x,y) \mapsto (y,x)~.$$
For any ordered double point $(x,y) \in D_2(f)$ let
$[x,y] \in D_2(f)/\ZZ_2$ be the unordered double point.

\begin{definition}~ \label{embed5}{\rm
Let $f:M^m \imm N^{2m}$ be an immersion, so that the double point set
$D_2(f)$ is a finite $0$-dimensional manifold with a free $\ZZ_2$-action.
\\
(i) The {\it self-intersection index}
of an unordered double point $[x,y] \in D_2(f)/\ZZ_2$ is
$$I[x,y]~=~[I(x,y)] \in \ZZ/\{1+(-1)^{m+1}\}~=~
\begin{cases}\hbox{$\ZZ$}&\hbox{if $m\equiv 0({\rm mod}\,2)$} \\[1ex]
\hbox{$\ZZ_2$}&\hbox{if $m\equiv 1({\rm mod}\,2)$}\end{cases}$$
which does not depend on the choice of lift of $[x,y]$ to an ordered
double point $(x,y) \in D_2(f)$, since $I(x,y)=(-1)^mI(y,x)\in \ZZ$.\\
(ii) The {\it geometric self-intersection number} of $f$ is
$$\mu(f)~=~\sum\limits_{[x,y] \in D_2(f)/\ZZ_2}I[x,y] \in
Q_{(-)^m}(\ZZ)~=~\ZZ/\{1+(-1)^{m+1}\}~.
\eqno{\hbox{\qed}}$$
}\end{definition}

\begin{proposition}~ \label{embed6} {\rm (Whitney \cite{whitney3}, Wall \cite[5.3]{wall2})}
The geometric self-intersection number $\mu(f)$ is a regular
homotopy invariant of $f:M^m \imm N^{2m}$ such that
$$\lambda(f,f)~=~(1+(-1)^m)\mu(f) +\chi(\nu_f) \in \ZZ$$
with $\mu(f)=0\in Q_{(-)^m}(\ZZ)$ if (and for $\pi_1(N)=\{1\}$, $m \geqslant 2$ only
if) $f$ is regular homotopic to an embedding.\hfill\qed
\end{proposition}

In particular, for even $m$
$$\mu(f)~=~(\lambda(f,f) - \chi(\nu_f))/2 \in \ZZ~.$$

\begin{definition}~ \label{embed7} {\rm  Let
$(e\times f:M^m \emb V \times N^{2m},f:M \imm N)$
be an $(m,2m,j)$-dimensional embedding-immersion pair, with Umkehr stable map
$$F~:~(V \times N)^{\infty}~=~V^{\infty} \wedge N^\infty \to
V^{\infty} \wedge T(\nu_f)~.$$
(i) The {\it homological self-intersection number}
$$\mu_V(f)~=~\psi_V(F)_*[N] \in
\begin{cases}\hbox{$\ZZ$}&
\hbox{if $j=1$}\\[1ex] \hbox{$\ZZ$}/\{1+(-1)^{m+1}\}&\hbox{if
${\rm dim}\,V\geqslant 2$}\end{cases}$$ is the image of the fundamental class
$[N] \in H_{2m}(N)$ under the induced map
$$\begin{array}{ll}
\psi_V(F)_*~:~H_{2m}(N)\to
&\widetilde{H}_{2m}(S(LV) ^{\infty} \wedge_{\ZZ_2}(T(\nu_f)\wedge
T(\nu_f)))\\[1ex]
&=~H_0(S(LV)\times_{\ZZ_2}(M \times M);\ZZ^w)\\[1ex]
&=~\begin{cases}\hbox{$\ZZ$}&\hbox{if $j=1$}\\[1ex]
\hbox{$\ZZ$}/\{1+(-1)^{m+1}\}&\hbox{if $j\geqslant 2$~.}
\end{cases}\end{array}$$
(ii) If
$j=1$ set $V=\R$, and use the self-linking manifold (\ref{link7})
$$L(M,N)~=~\{(x,y) \in M \times M\st e(x) < e(y)\in \R, f(x)=f(y)\in N \}$$
and the decomposition
$$D_2(f)~=~L(M,N) \cup TL(M,N)$$
to lift each unordered double point $(x,y)\in D_2(f)/\ZZ_2$ to an ordered double
point $(x,y) \in L(M,N)$, ordered according to $e(x) < e(y)$.
The {\it integral geometric self-intersection} of $f$ is the
self-linking number
$$\mu^{\ZZ}(f)~=~[L(M,N)]~=~\sum\limits_{(x,y)\in L(M,N)}I(x,y)\in \hbox{$\ZZ$}~.$$
For even $m$ $\mu^{\ZZ}(f)=\mu(f)$. \hfill\qed
}\end{definition}

\begin{proposition}~ \label{embed8} {\rm (i)} If $m$ is even or if
${\rm dim}(V)\geqslant 2$ the geometric self-intersection number
is just the homological self-intersection number
$$\mu(f)~=~\mu_V(f)\in Q_{(-1)^m}(\ZZ)~=~\hbox{$\ZZ$}/\{1+(-1)^{m+1}\}~.$$
{\rm (ii)} If ${\rm dim}(V)=1$ the integral geometric
self-intersection number is just the homological self-intersection
number
$$\mu^{\ZZ}(f)~=~\mu_{\R}(f)\in \hbox{$\ZZ$}~.$$
\end{proposition}
\noindent{\it Proof\  .} Immediate from \ref{double6} and
\ref{double13}.\hfill\qed

\begin{example}
{\rm We refer to Crabb and Ranicki \cite{crabbranicki} for an interpretation
in terms of the geometric Hopf invariant of the Smale-Hirsch-Haefliger regular homotopy classification of immersions $f:M^m \imm N^n$  in the metastable dimension range
$3m < 2n-1$ (when a generic $f$ has no triple points). In particular, this applies
to the case $n=2m$ with $m \geqslant 2$.}\\
\hfill\qed
\end{example}

\chapter{The $\pi$-equivariant geometric Hopf invariant}\label{pi-equivariant}

The stable homotopy constructions of \S\ref{geohopf} (in particular the geometric Hopf invariant) and the double point theorem of \S\ref{doublepoint} are so natural that they have $\pi$-equivariant versions, for any discrete  group $\pi$, inducing the corresponding chain level constructions of Ranicki \cite{ranicki1,ranicki2}, with applications to
non-simply-connected surgery obstruction theory  (\S\ref{non}).

\section{$\pi$-equivariant $S$-duality}

\begin{definition}~{\rm
{\rm (i)} Given pointed $\pi$-spaces $X,Y$ define the
{\it integral $\pi$-equivariant homotopy groups}
$$\{X;Y\}_{0,\pi}~=~\varinjlim\limits_U [U^\infty \wedge X;U^\infty \wedge Y]_{\pi}$$
with $U$ running over finite-dimensional inner product spaces.\\
{\rm (ii)} Given pointed $\pi$-spaces $X,Y,Z$, a finite-dimensional inner
product space $V$ and a map
$$\sigma~:~V^{\infty} \to X\wedge_{\pi}Y$$
define the {\it slant} products
$$\begin{array}{l}
\sigma \backslash -~:~\{X,Z\}_{0,\pi} \to \{V^{\infty},Z \wedge_{\pi}Y\}~;\\[1ex]
\hskip25pt
(f:U^{\infty} \wedge X \to U^{\infty} \wedge Z) \mapsto
((f \wedge 1)(1\wedge \sigma):(U\oplus V)^{\infty}\to  U^{\infty} \wedge Z \wedge_{\pi}Y)~,\\[1ex]
\sigma \backslash -~:~\{Y,Z\}_{0,\pi} \to \{V^{\infty},X \wedge_{\pi}Z\}~;\\[1ex]
\hskip25pt
(f:U^{\infty} \wedge Y \to U^{\infty} \wedge Z) \mapsto
((f \wedge 1)(1\wedge \sigma):(U\oplus V)^{\infty} \to  U^{\infty} \wedge X \wedge_{\pi}Z)~.
\end{array}$$
{\rm (iii)} A map $\sigma:V^{\infty} \to X\wedge_{\pi}Y$ is an
{\it integral $\pi$-equivariant $S$-duality} map if the slant products
$\sigma \backslash -$ in {\rm (ii)} are isomorphisms for every
$\pi$-space $Z$.\\
}\hfill\qed
\end{definition}

\begin{proposition}~ {\rm (Ranicki \cite[\S3]{ranicki2})}\label{piS-dual}
For any semifree finite pointed $CW$ $\pi$-complex $X$ there exist a
finite-dimensional inner product space $V$, a semifree
finite pointed $CW$ $\pi$-complex $Y$ and an integral $\pi$-equivariant
$S$-duality map $\sigma:V^{\infty} \to X \wedge_\pi Y$.\\
\hfill\qed
\end{proposition}

\begin{remark} {\rm The theory of \cite{ranicki2}
deals with the `integral' stable $\pi$-equivariant homotopy groups
$$\{X;Y\}_{0,\pi}~=~\mathop{\varinjlim}\limits_V\,
[V^{\infty} \wedge X,V^{\infty}  \wedge Y]_{\pi}$$
with the direct limit running over all the finite-dimensional
inner product spaces $V$. The forgetful map for $\pi=\ZZ_2$
$$\{X;Y\}_{0,\ZZ_2} \to \{X;Y\}_{\ZZ_2}$$
is in general neither injective nor surjective, but by Adams \cite{adams3}
(cf. Proposition \ref{Z2cover} (ii)) it is an isomorphism for
finite $CW$ $\ZZ_2$-complexes $X,Y$ with $X$ semifree, and in this
case the $\ZZ_2$-equivariant $S$-duality theories of Wirthm\"uller \cite{wirth} (cf. Proposition \ref{Z2S-dual})
and \cite{ranicki2} coincide.}\\
\hfill\qed
\end{remark}

\section{The $\pi$-equivariant constructions}\label{equivconstr}

Let $X,Y$ be pointed $\pi$-spaces, and let $V$
be an inner product space with trivial $\pi$-action.
The geometric Hopf invariant (\ref{hopf}) of a $\pi$-equivariant
map $F:V^{\infty}\wedge X \to V^{\infty} \wedge Y$ is a $\ZZ_2$-equivariant map
$$h_V(F)~:~\Sigma S(LV)^+ \wedge V^{\infty} \wedge X \to LV^{\infty} \wedge
V^{\infty} \wedge Y \wedge Y$$
which is also $\pi$-equivariant, with
$$\pi \times Y \wedge Y \to Y \wedge Y ~;~(g,(y_1,y_2)) \mapsto (gy_1,gy_2)~.$$
The `$\pi$-equivariant geometric Hopf invariant of $F$' is the induced map
of the $\pi$-quotients
$$h_V(F)/\pi~:~\Sigma S(LV)^+ \wedge V^{\infty} \wedge X/\pi
\to LV^{\infty} \wedge V^{\infty} \wedge Y \wedge_{\pi}Y$$
with
$$\begin{array}{l}
Y \wedge_{\pi} Y~=~Y \wedge Y/\{(y_1,y_2) \sim (gy_1,gy_2)\,\vert\,
g \in \pi,\,y_1,y_2 \in Y\}~,\\[1ex]
T~:~Y \wedge_{\pi} Y \to Y \wedge_{\pi} Y~;~[y_1,y_2] \mapsto [y_2,y_1]~.
\end{array}$$

\begin{proposition}~ \label{sequence2}
{\rm (i)} For any pointed space $X$, any pointed $\pi$-space $Y$ and
any inner product space $V$
there is defined a long exact sequence of abelian groups/pointed sets
$$\begin{array}{l}
\dots \xymatrix{\ar[r]&} [\Sigma X,Y\wedge_{\pi}Y]_{\ZZ_2}
\xymatrix{\ar[r]^-{\di{s_{LV}^*}} &}
[\Sigma S(LV)^+  \wedge X,Y\wedge_{\pi}Y]_{\ZZ_2}
\xymatrix{\ar[r]^-{\di{\alpha^*_{LV}}}&}\\
\hskip150pt
[LV^{\infty} \wedge X,Y\wedge_{\pi}Y]_{\ZZ_2} \xymatrix{\ar[r]^-{\di{0_{LV}^*}}&}
[X,Y\wedge_{\pi}Y]_{\ZZ_2}
\end{array}$$
{\rm (ii)} For any pointed space $X$, any pointed $\pi$-space $Y$
and any inner product spaces
$V,W$ there is defined a split short exact sequence of abelian groups
$$\begin{array}{l}
0 \to [\Sigma S(LV)^+ \wedge W^{\infty} \wedge X,LV^{\infty} \wedge W^{\infty}
\wedge (Y\wedge_{\pi}Y)]_{\ZZ_2}
\xymatrix{\ar[r]^-{\di{\alpha^*_{LV}}}&}\\[1ex]
[LV^{\infty} \wedge W^{\infty} \wedge X,LV^{\infty} \wedge W^{\infty}
\wedge (Y\wedge_{\pi}Y)]_{\ZZ_2}\xymatrix{\ar[r]^-{\di{\rho}}&}
[W^{\infty} \wedge X,W^{\infty} \wedge Y/\pi] \to 0
\end{array}$$
with $\rho$ defined by the fixed points of the $\ZZ_2$-action.
The map $\rho$ is split by
$$\begin{array}{l}
\sigma~:~[W^{\infty} \wedge X,W^{\infty} \wedge Y/\pi] \to
[LV^{\infty} \wedge W^{\infty} \wedge X,LV^{\infty} \wedge W^{\infty}
\wedge (Y\wedge_{\pi}Y)]_{\ZZ_2}~;\\[1ex]
F \mapsto \sigma(F)=1_{LV^{\infty}}\wedge \Delta F~~
(\sigma(F)(v,w,x)=(v,u,y)~\hbox{\it if}~F(w,x)=(u,y))~,
\end{array}$$
with
$$\Delta~:~Y/\pi~=~(Y\wedge_{\pi}Y)^{\ZZ_2} \emb Y\wedge_{\pi}Y~;~[y] \mapsto [y,y]~.$$
For any $\ZZ_2$-equivariant map
$$G~:~LV^{\infty} \wedge W^{\infty} \wedge X \to LV^{\infty} \wedge W^{\infty}
\wedge (Y\wedge_{\pi}Y)$$
the map of fixed points
$$F~=~\rho(G)~:~W^{\infty} \wedge X \to W^{\infty} \wedge Y/\pi$$
is such that $(1 \wedge \Delta)\sigma(F)$ and $G$ agree on
$0^+ \wedge W^{\infty} \wedge X$,
with the relative difference $\ZZ_2$-equivariant map
$$\delta((1 \wedge \Delta)\sigma(F),G)~:~\Sigma S(LV)^+ \wedge W^{\infty} \wedge X \to
LV^{\infty} \wedge W^{\infty} \wedge (Y\wedge_{\pi}Y)$$
such that
$$G-(1 \wedge \Delta)\sigma(F)~=~\alpha^*_{LV}\delta((1 \wedge \Delta)\sigma(F),G)
\in {\rm im}(\alpha^*_{LV})~=~{\rm ker}(\rho)~.$$
\end{proposition}
\noindent{\it Proof\  .} This is just a special case of
Proposition \ref{sequence}, with $Y$ replaced by $Y \wedge_{\pi}Y$,
and $Y^{\ZZ_2}$ replaced by $(Y \wedge_{\pi}Y)^{\ZZ_2}=Y/\pi$.\\
\hfill\qed

\begin{proposition}~ \label{stablesequence4}
{\rm (i)} For any inner product spaces $U,V$, any pointed spaces $X$ and any pointed $\pi$-space $Y$
there is defined a commutative braid of exact sequences of stable
homotopy groups
$$\xymatrix@C-35pt{
A_1 \ar[dr] \ar@/^2pc/[rr]^-{\di{\alpha_{LV}}} && \{X;S(LV)^+ \wedge Y\wedge_{\pi}Y\}_{\ZZ_2}
\ar[dr]^-{\di{s_{LV}}}\ar@/^2pc/[rr] &&
\{S(LU)^+ \wedge X;Y\wedge_{\pi}Y\}_{\ZZ_2}\ar[dr]&\\
&~~~~~~~~~A_2~~~~~~~~~ \ar[ur] \ar[dr] &&
\{X;Y\wedge_{\pi}Y\} \ar[ur]^-{\di{s^*_{LU}}} \ar[dr]^-{\di{0_{LV}}}&&A_3\\
A_4\ar[ur] \ar@/_2pc/[rr]^-{\di{\alpha^*_{LU}}}
&&\{LU^{\infty} \wedge X;Y\wedge_{\pi}Y\}_{\ZZ_2}
\ar@/_2pc/[rr]\ar[ur]^-{\di{0^*_{LU}}} &&
\{X;V^{\infty}\wedge Y\wedge_{\pi}Y\}_{\ZZ_2}\ar[ur]&}$$

\bigskip

\noindent with
$$\begin{array}{l}
A_1~=~\{\Sigma X;LV^{\infty} \wedge Y\wedge_{\pi}Y\}_{\ZZ_2}~,~
A_2~=~\{\Sigma S(LU\oplus LV)^+\wedge X;LV^{\infty} \wedge Y\wedge_{\pi}Y\}_{\ZZ_2}~,\\[1ex]
A_3~=~\{S(LU\oplus LV)^+ \wedge X;V^{\infty} \wedge Y\wedge_{\pi}Y\}_{\ZZ_2}~,~
A_4~=~\{\Sigma S(LU)^+ \wedge X;Y\wedge_{\pi}Y\}_{\ZZ_2}~.
\end{array}$$
{\rm (ii)} For any inner product space $V \subseteq \R(\infty)$
$$\{X;LV^{\infty} \wedge Y \wedge_{\pi} Y\}_{\ZZ_2}~=~
\{X;Y/\pi\} \oplus \{X;S(L\R(\infty))/S(LV) \wedge Y \wedge_{\pi} Y\}_{\ZZ_2} ~,$$
with a split short exact sequence
$$\begin{array}{l}
0 \to \{X;S(L\R(\infty))/S(LV) \wedge Y \wedge_{\pi} Y\}_{\ZZ_2}
\xymatrix{\ar[r]^-{\di{\delta}}&} \{X;LV^{\infty} \wedge Y \wedge_{\pi}Y\}_{\ZZ_2}\\[1ex]
\hskip150pt \xymatrix{\ar[r]^-{\di{\rho}}&} \{X;Y/\pi\} \to 0
\end{array}$$
where $\delta$ is induced by the $\ZZ_2$-equivariant connecting map
$$\delta~:~
S(L\R(\infty))/S(LV)~=~\varinjlim\limits_k S(L\R^k)/S(LV) \to s S(LV)~=~LV^{\infty}~,$$
$\rho$ is defined by the fixed points of the $\ZZ_2$-action and
$\rho$ is split by
$$\sigma~:~\{X;Y/\pi\} \to \{X;LV^{\infty} \wedge
Y\wedge_{\pi} Y\}_{\ZZ_2}~;~F \mapsto (0 \wedge \Delta_{Y/\pi})F~.$$
In particular, for $V=\R(\infty)$ the fixed point map is an isomorphism
$$\rho~:~\{X;L\R(\infty)^{\infty}\wedge Y \wedge_{\pi} Y\}_{\ZZ_2}
\xymatrix{\ar[r]^-{\di{\cong}}&} \{X;Y/\pi\}~,$$
and the case $V=\{0\}$ gives
$$\{X;Y\wedge_{\pi}Y\}_{\ZZ_2}~=~\{X;Y/\pi\}\oplus
\{X;S(L\R(\infty))^+\wedge Y\wedge_\pi Y\}_{\ZZ_2}~.$$
{\rm (iii)} For any inner product space $U$ there is defined a long exact sequence
$$\begin{array}{l}
\xymatrix{
\dots \ar[r]&\{\Sigma S(LU \oplus L\R(\infty))^+\wedge
X;L\R(\infty)^{\infty} \wedge Y \wedge_{\pi} Y\}_{\ZZ_2}&}\\
\hskip25pt \xymatrix{\ar[r] &\{LU^{\infty} \wedge X;Y\wedge_{\pi} Y\}_{\ZZ_2}
\ar[r]^-{\di{\rho}} & \{X;Y/\pi\}&}\\
\hskip25pt \xymatrix{\ar[r] &\{S(LU \oplus L\R(\infty))^+\wedge
X;L\R(\infty)^{\infty} \wedge Y \wedge_{\pi} Y\}_{\ZZ_2} \ar[r] &\dots}
\end{array}
$$
with $\rho$ defined by the fixed points of the $\ZZ_2$-action, and
$$\begin{array}{l}
\{S(LU \oplus L\R(\infty))^+\wedge
X;L\R(\infty)^{\infty} \wedge Y \wedge_{\pi} Y\}_{\ZZ_2}\\[1ex]
\hspace*{50pt}=~ \{LU^\infty \wedge X;\Sigma S(LU \oplus L\R(\infty))^+
\wedge Y \wedge_{\pi} Y\}_{\ZZ_2}~.
\end{array}$$
{\rm (iv)} For any pointed $\pi$-space $X$ the morphism
$$\begin{array}{l}
0_{LV}~:~\{X/\pi;X \wedge_{\pi} X\}_{\ZZ_2}\\[1ex]
 \to \{X/\pi;LV^{\infty} \wedge X \wedge_{\pi} X\}_{\ZZ_2}~=~
 \{V^{\infty} \wedge X/\pi;V^{\infty} \wedge LV^{\infty} \wedge X \wedge_{\pi} X\}_{\ZZ_2}
\end{array}$$
 sends $\Delta_{X/\pi}$ to $0_{LV} \wedge \Delta_{X/\pi}=(\kappa^{-1}_V \wedge 1)\Delta_{V^{\infty} \wedge {X/\pi}}$.
\end{proposition}
\begin{proof} This is a $\ZZ_2$-equivariant version of the braid of
Proposition \ref{stablesequence0}; for $\pi=\{1\}$ this is just
Proposition \ref{stablesequence1}.
\hfill\qed\end{proof}

The definitions of the various constructions of \S\ref{geohopf} translate verbatim into their $\pi$-equivariant analogues, with matching properties. We shall only state the definitions here.

\begin{definition}~ \label{pisym1}
{\rm  The {\it $\pi$-equivariant symmetric construction} $\phi_V(X)$ is defined for
a pointed $\pi$-space $X$ and an inner product space $V$ to be
the $\pi\times\ZZ_2$-equivariant map
\index{$\pi$-equivariant!symmetric construction, $\phi_V(X)$}
$$\phi_V(X)~=~s_{LV}\wedge \Delta_X~:~S(LV)^+ \wedge X \to X \wedge  X~;~
(v,[x]) \mapsto [x,x]~.$$
\hfill\qed}
\end{definition}

Passing to the $\pi$-quotients there is defined a $\ZZ_2$-equivariant map
$$\phi_V(X)~=~s_{LV}\wedge \Delta_X~:~S(LV)^+ \wedge X/\pi \to X \wedge_{\pi} X~.$$

\begin{definition}~ \label{pihopf1}
{\rm  The {\it $\pi$-equivariant geometric Hopf invariant} of a
$\pi$-equivariant map $F:V^{\infty} \wedge X \to
V^{\infty} \wedge Y$ is the $\ZZ_2$-equivariant map given by the relative
difference of the $\ZZ_2$-equivariant maps
$$p~=~(1 \wedge \Delta_Y)(1 \wedge F)~,~
q~=~(\kappa^{-1}_V\wedge 1) (F \wedge F)(\kappa_V \wedge \Delta_X)$$
with \index{$\pi$-equivariant!geometric Hopf invariant, $h_V(F)$}
$$\begin{array}{ll}
h_V(F)~=~\delta(p,q)~:&\Sigma S(LV)^+ \wedge V^{\infty} \wedge X \to
LV^{\infty} \wedge V^{\infty} \wedge Y\wedge Y~;\\[1ex]
&(t,u,v,x) \mapsto \begin{cases}
p([ 1-2t,u],v,x)&\hbox{if $0 \leqslant  t \leqslant  1/2$}\\[1ex]
q([ 2t-1,u],v,x)&\hbox{if $1/2 \leqslant  t
\leqslant  1$}
\end{cases}\\[3ex]
&(t \in I,u \in S(LV),v \in V,x \in X)~.
\end{array}$$
\hfill\qed
}\end{definition}

\begin{definition}~ \label{pihopf9}
{\rm The {\it $\pi$-equivariant stable geometric Hopf invariant} of a $\pi$-equivariant map $F:V^{\infty} \wedge X \to V^{\infty} \wedge Y$ is
the stable relative difference of $p$ and $q$ (\ref{stablereldif3}),
the $\pi \times \ZZ_2$-equivariant map
\index{$\pi$-equivariant!geometric Hopf invariant, stable, $h'_V(F)$}
$$h'_V(F)~=~\delta'(p,q)~:~LV^{\infty} \wedge V^{\infty} \wedge X \to
LV^{\infty} \wedge V^{\infty} \wedge S(LV)^+ \wedge Y \wedge Y$$
which we shall regard as a stable $\ZZ_2$-equivariant map
$$h'_V(F)~:~X/\pi \sto S(LV)^+ \wedge Y \wedge_{\pi} Y~.$$}
\hfill$\qed$
\end{definition}

\begin{definition}~ \label{piquad1}
{\rm  The {\it $\pi$-equivariant quadratic construction} on a $\pi$-equivariant map
$F:V^{\infty} \wedge X \to V^{\infty} \wedge Y$ is the stable $\pi$-equivariant map
\index{$\pi$-equivariant!quadratic construction, $\psi_V(F)$}
$$\psi_V(F)~:~X \sto S(LV)^+\wedge_{\ZZ_2} (Y\wedge Y)$$
given by the image of the $\ZZ_2$-equivariant stable geometric Hopf invariant (\ref{pihopf9})
$$h'_V(F)~=~\delta'(p,q)~:~ LV^{\infty} \wedge V^{\infty} \wedge X \to
LV^{\infty} \wedge V^{\infty} \wedge S(LV)^+ \wedge Y \wedge Y$$
under the isomorphism given by Proposition \ref{Adams}
$$\{X;S(LV)^+ \wedge_{\ZZ_2} (Y \wedge Y)\}_{\pi}
\iso \{X;S(LV)^+ \wedge Y \wedge Y\}_{\pi \times \ZZ_2}~.$$
\hfill$\qed$}
\end{definition}

Passing to the $\pi$-quotients gives
$$\psi_V ~:~[V^{\infty} \wedge X,V^{\infty} \wedge Y]_{\pi}
\to \{X/\pi;S(LV)^+ \wedge_{\ZZ_2}(Y \wedge_{\pi} Y)\}~\cong~
\{X;S(LV)^+\wedge (Y \wedge_{\pi}Y)\}_{\ZZ_2}~.$$

\begin{definition}~ \label{piultraquadratic}
{\rm The {\it $\pi$-equivariant ultraquadratic construction} on a
$\pi$-equivariant map $F:\Sigma X \to \Sigma Y$
is the quadratic construction for the special case $V=\R$
\index{ultraquadratic construction, $\widehat{\psi}(F)$}
$$\widehat{\psi}(F)~=~\psi_{\R}(F)~:~X/\pi \sto
S(LV)^+\wedge_{\ZZ_2} (Y \wedge_\pi Y)~=~Y \wedge_\pi Y~,$$
identifying $S(LV)=S^0=\{1,-1\}$ with $\ZZ_2$ acting by permutation.\hfill$\qed$}
\end{definition}

\begin{definition}~ \label{pisym5}
{\rm  The {\it $\pi$-equivariant spectral Hopf invariant map} of a $\pi$-equivariant
map $F:X \to V^{\infty} \wedge Y$ is the $\pi \times \ZZ_2$-equivariant map
given by the relative difference (\ref{reldif1})
\index{$\pi$-equivariant!spectral Hopf invariant, $sh_V(F)$}
$$\begin{array}{l}sh_V(F)~=~
\delta\big((G \wedge G)\delta \phi_V(V^{\infty}\wedge Y)(1 \wedge F),
\phi_V(\Cc(F))(1\wedge GF)\big):\\[1ex]
\hskip150pt
\Sigma S(LV)^+\wedge X~\to \Cc(F)  \wedge  \Cc(F)\end{array}$$
with $G:V^{\infty} \wedge Y \to \Cc(F)$ the inclusion in the mapping cone and
$$1 \wedge GF~:~CS(LV)^+ \wedge X~=~S(LV)^+ \wedge CX\to S(LV)^+ \wedge \Cc(F)$$
the $\pi$-equivariant null-homotopy of $1\wedge GF:S(LV)^+ \wedge X \to S(LV)^+ \wedge \Cc(F)$
determined by the inclusion
$GF:CX \to\Cc(F)=(V^{\infty} \wedge Y)\cup_FCX$.\hfill\qed
}\end{definition}

Here is the $\pi$-equivariant version of the Double Point Theorem \ref{doubleHopf}.

 Let $e=(g,f):V \times X \emb V \times Y$ be an embedding of a map $f:X \to Y$. Given a regular
cover $p:\widetilde{Y} \to Y$ with group of covering translations $\pi$ let
$$\widetilde{X}~=~f^*\widetilde{Y}~=~
\{(x,\widetilde{y}) \in X \times \widetilde{Y}\,\vert\,
f(x) =p(\widetilde{y}) \in Y\}$$
be the pullback cover of $X$, so that there are defined
$\pi$-equivariant lifts of $f,e$
$$\begin{array}{l}
\widetilde{f}~:~\widetilde{X} \to \widetilde{Y}~;~(x,\widetilde{y}) \mapsto \widetilde{y}~,\\[1ex]
\widetilde{e}~=~(\widetilde{g},\widetilde{f})~:~
V \times \widetilde{X} \emb V \times \widetilde{Y}~;~
(v,(x,\widetilde{y})) \mapsto (g(v,x),\widetilde{y})
\end{array}$$
and hence a $\pi$-equivariant Umkehr map
$$\widetilde{F}~:~V^{\infty} \wedge \widetilde{Y}^{\infty}
\to V^{\infty} \wedge \widetilde{X}^{\infty}~.$$

\begin{proposition}
The $\pi$-equivariant geometric Hopf invariant
$h_V(\widetilde{F})$ is given up
to natural $\pi \times \ZZ_2$-equivariant homotopy by the composite
$$\begin{array}{l}
h_V(\widetilde{F})/\pi~=~((1 \wedge \widetilde{A})h_V(\widetilde{F})_Y)/\pi~:~\Sigma S(LV)^+ \wedge V^{\infty} \wedge Y^{\infty} \\[1ex]
\hskip50pt \to LV^{\infty} \wedge V^{\infty} \wedge (X\times_{Y}X)^{\infty} \to
LV^{\infty} \wedge V^{\infty} \wedge (\widetilde{X}\times_\pi \widetilde{X})^{\infty}
\end{array}
$$
with
$$\begin{array}{l}
h_V(\widetilde{F})_Y~=~(1 \wedge i)\widetilde{F}_1 +(1 \wedge \Delta_{\widetilde{X}})\delta(\widetilde{F}_2,\widetilde{F}_3)~:\\[1ex]
\Sigma S(LV)^+ \wedge V^{\infty} \wedge \widetilde{Y}^{\infty} \to\\[1ex]
\hskip25pt
LV^{\infty} \wedge V^{\infty} \wedge (\widetilde{X}\times_{\widetilde{Y}}\widetilde{X})^{\infty}~=~
LV^{\infty}\wedge V^{\infty} \wedge (D_2(\widetilde{f})^{\infty}\vee \widetilde{X}^{\infty})~,\\[1ex]
i~=~{\rm inclusion}~:~D_2(\widetilde{f}) \to \widetilde{X}\times_{\widetilde{Y}}\widetilde{X} ~;~(x_1,x_2) \mapsto (x_1,x_2)~,\\[1ex]
\Delta_{\widetilde{X}}~=~{\rm diagonal}~:~\widetilde{X} \to \widetilde{X} \times_{\widetilde{Y}}\widetilde{X}~;~x \mapsto (x,x)~,\\[1ex]
A~=~{\rm assembly}~=~\Delta_{\widetilde{X}} \sqcup i~:~\widetilde{X} \times_{\widetilde{Y}}\widetilde{X}~=~\widetilde{X} \sqcup D_2(\widetilde{f}) \to \widetilde{X} \times \widetilde{X}~,\\[1ex]
\widetilde{F}_1,\widetilde{F}_2,\widetilde{F}_3~\hbox{\rm as in}~\ref{Umkehrterm}.
\end{array}$$
\hfill\qed
\end{proposition}

\chapter{Surgery obstruction theory}
\label{surgeryobstruction}

In this final chapter we apply the geometric Hopf invariant to surgery obstruction theory.

\S\ref{poincare} reviews geometric Poincar\'e complexes, the Spivak normal structure, normal maps and geometric Umkehr maps.
An $n$-dimensional normal map $(f,b):M \to X$ determines a $\pi_1(X)$-equivariant geometric Umkehr map $F:\Sigma^{\infty}\widetilde{X}^+ \to \Sigma^{\infty}\widetilde{M}^+$ inducing the algebraic Umkehr
$\ZZ[\pi_1(X)]$-module chain map $f^!:C(\widetilde{X}) \to C(\widetilde{M})$, with $\widetilde{X}$ the universal cover of $X$ and
$\widetilde{M}=f^*\widetilde{X}$ the pullback cover of $M$.
In \S\ref{simple} (ignoring $\pi_1(X)$) and \S\ref{non} we shall show
that the $\pi_1(X)$-equivariant geometric Hopf invariant
$h(F)=\psi_F:H_n(X) \to Q_n(C(\widetilde{M}))$ determines the surgery obstruction $\sigma_*(f,b) \in L_n(\ZZ[\pi_1(X)])$, simplifying the indirect identification in \cite[Prop. 7.1]{ranicki2}.

In \S\ref{codim2surgery} we use the geometric Hopf invariant to construct various quadratic structures (such as the Seifert form) arising in codimension 2 surgery.

Finally, in \S\ref{total} we apply the geometric Hopf invariant to the total surgery obstruction $s(X) \in \SS_n(X)$ of Ranicki \cite{ranicki0,ranicki3,ranicki4, ranicki5}.

\section{Geometric Poincar\'e complexes and Umkehr maps}\label{poincare}

By definition, an $n$-dimensional {\it geometric Poincar\'e complex $X$}\index{geometric Poincar\'e complex} in the sense of Wall \cite{wallpoincare} is a finite $CW$ complex with a fundamental class\footnote{using twisted coefficients in the nonorientable case} $[X] \in H_n(X)$
such that the $\ZZ[\pi_1(X)]$-module chain map
$$[X] \cap - ~:~ C(\widetilde{X})^{n-*}~=~{\rm Hom}_{\ZZ[\pi_1(X)]}(C(\widetilde{X}),\ZZ[\pi_1(X)])_{n-*} \to C(\widetilde{X})$$
is a chain equivalence, with $\widetilde{X}$ the universal cover of $X$.
In particular, a compact $n$-dimensional manifold is an $n$-dimensional geometric Poincar\'e complex, as is any finite CW complex homotopy equivalent to a manifold\footnote{a differentiable manifold is triangulable and is a $CW$ complex, whereas a topological manifold need not be triangulable and only has the homotopy type of a $CW$ complex.}.

The Browder-Novikov-Sullivan-Wall surgery theory deals with the two fundamental questions:
\begin{itemize}
\item[(i)]~ is an $n$-dimensional geometric Poincar\'e complex homotopy equivalent to a manifold?
\item[(ii)]~ is a homotopy equivalence of $n$-dimensional manifolds homotopic to a diffeomorphism?
\end{itemize}
Note that (ii) is the rel $\partial$ version of (i), since the mapping cylinder
of a homotopy equivalence of manifolds is a geometric Poincar\'e cobordism between manifolds. The theory works best for $n \geqslant 5$, and there is also a version for topological manifolds and homeomorphisms.

A finite $CW$ complex $X$ is an $n$-dimensional geometric Poincar\'e complex if and only if for any regular neighbourhood $W$ of an embedding $X \subset S^{n+k}$ ($k$ large)
$$\hbox{\rm homotopy fibre}(\partial W \subset W)~\simeq~S^{k-1}$$
(Ranicki \cite[Prop. 3.11]{ranicki2}). The pair
$$(\nu_X:X \to BG(k),\rho_X:S^{n+k} \to T(\nu_X))$$
defines the {\it Spivak normal structure}\index{Spivak normal structure, $(\nu_X,\rho_X)$} of $X$, with $\nu_X$ the $(k-1)$-spherical fibration
$$\nu_X~:~S^{k-1} \to \partial W \to W~\simeq~X$$
and $\rho_M$ the Pontrjagin-Thom map
$$\rho_M~:~S^{n+k} \to S^{n+k}/{\rm cl.}(S^{n+k}\backslash W)~=~W/\partial W~=~T(\nu_X)~,$$
such that the composite of the Hurewicz map and the Thom isomorphism
$$\pi_{n+k}(T(\nu_X)) \to \widetilde{H}_{n+k}(T(\nu_X))~\cong~H_n(X)$$
sends $\rho_M$ to the fundamental class $[X] \in H_n(X)$ (using twisted
coefficients in the nonoriented case).

An $n$-dimensional {\it normal map}\index{normal map, $(f,b)$} $(f,b):M \to X$ is a degree 1 map $f:M \to X$ from an $n$-dime\-nsional manifold $M$ to an $n$-dimensional geometric Poincar\'e complex $X$, together with a bundle map
$b:\nu_M \to \eta$ over $f$ from the stable normal bundle $\nu_M$ of $M$ to a stable bundle $\eta$ over $X$ in the Spivak normal class.
The surgery obstruction $\sigma_*(f,b) \in L_n(\ZZ[\pi_1(X)])$ of Wall \cite{wall2} is such that $\sigma_*(f,b)=0$ if (and for $n \geqslant 5$ only if) $(f,b)$ is normal bordant to a homotopy equivalence. The original construction of $\sigma_*(f,b)$ used preliminary surgeries below the middle dimension to make $(f,b):M \to X$ such that $f_*:\pi_m(M) \to \pi_m(X)$
is an isomorphism for $m<[(n-1)/2]$, and then read off the surgery obstruction form the middle dimensional data. The surgery obstruction
was expressed in Ranicki \cite{ranicki1,ranicki2} directly from $(f,b)$,
as the cobordism class of an $n$-dimensional quadratic Poincar\'e complex
$(C,\psi)$ over $\ZZ[\pi_1(X)]$. The homology of $C$ is
$$H_*(C)~=~K_*(M)~=~{\rm ker}(\widetilde{f}_*:H_*(\widetilde{M}) \to H_*(\widetilde{X}))$$
with $\widetilde{X}$ = universal cover of $X$, $\widetilde{M}=f^*\widetilde{X}$ the pullback cover of $M$.
The quadratic structure $\psi \in Q_n(C)$ was obtained using the `quadratic construction' on a $\pi_1(X)$-equivariant geometric Umkehr map
$F:\Sigma^{\infty}\widetilde{X}^+ \to \Sigma^{\infty}\widetilde{M}^+$. However, the identification (\cite[Prop. 7.1]{ranicki2})
$$\sigma_*(f,b)~=~(C,\psi) \in L_n(\ZZ[\pi_1(X)])$$
was somewhat indirect. At the time, it was not obvious how to directly
express Wall's self-intersection form $\mu(S^m \imm M^{2m}) \in
Q_{(-)^m}(\ZZ[\pi_1(M)])$  in terms of the quadratic construction.
This will be remedied in \S\ref{non} below

For any $n$-dimensional geometric Poincar\'e complex $X$ with Spivak normal structure $(\nu_X,\rho_X)$ the composite
$$\xymatrix{\sigma_X~=~\Delta \rho_X:S^{n+k} \ar[r]^-{\di{\rho_X}} & T(\nu_X)
\ar[r]^-{\di{\Delta}} & X^+ \wedge T(\nu_X)}$$
is the Atiyah-Spivak-Wall $S$-duality map. (See \S\ref{Sduality} above for a treatment of $S$-duality). For any pointed space $Z$ there are defined $S$-duality isomorphisms
$$\begin{array}{l}
\sigma_X~:~\{T(\nu_X);Z\} \xymatrix{\ar[r]^-{\di{\cong}}&} \{S^{n+k};X^+ \wedge Z \}~;~
F \mapsto F\sigma_X~,\\[1ex]
\sigma_X~:~\{X^+;Z\} \xymatrix{\ar[r]^-{\di{\cong}}&} \{S^{n+k};Z \wedge T(\nu_X) \}~;~
G \mapsto G\sigma_X~.
\end{array}$$
For any $n$-dimensional geometric Poincar\'e complexes
$M,X$ there is thus defined an $S$-duality isomorphism
$$\{T(\nu_M);T(\nu_X)\}~\cong~\{X^+;M^+\}~.$$
As in Ranicki \cite{ranicki2} given a normal map $(f,b):M \to X$
 call any $S$-dual of
$T(b):T(\nu_M) \to T(\eta)=T(\nu_X)$ a {\it geometric Umkehr map}
for $(f,b)$\index{Umkehr!geometric}
$$F~=~T(b)^*~:~\Sigma^{\infty} X^+ \to \Sigma^{\infty} M^+$$
inducing the Umkehr chain map
$$f^{\,!}~:~C(X)~\simeq~C(X)^{n-*} \xymatrix@C+10pt{\ar[r]^-{\di{f^*}}&}
C(M)^{n-*}~\simeq~C(M)~,$$
such that
$$\begin{array}{l}
(1 \wedge f)F \simeq 1~:~\Sigma^{\infty}X^+ \to \Sigma^{\infty}  X^+~,\\[1ex]
ff^{\,!} \simeq 1~:~C(X) \to C(X)~.\end{array}$$
The homology and cohomology groups split as
$$\begin{array}{l}
H_*(M)~=~H_*(X) \oplus K_*(M)~,\\[1ex]
H^*(M)~=~H^*(X) \oplus K^*(M)
\end{array}
$$
with
$$\begin{array}{l}
K_m(M)~=~{\rm ker}(f_*:H_m(M) \to H_m(X))~,\\[1ex]
K^m(M)~=~{\rm ker}(f^!:H^m(M) \to H^m(X))~.\end{array}$$
The Poincar\'e duality isomorphisms $H^*(M)\cong H_{n-*}(M)$
restrict to Poincar\'e duality isomorphisms
$$K^*(M)~\cong~K_{n-*}(M)~.$$
The intersection pairing
$$\lambda~:~K_m(M) \times K_{n-m}(M) \to \ZZ$$
is such that
$$\begin{array}{l}
\lambda(x,y)~=~(-)^{m(n-m)}\lambda(y,x) \in \ZZ ~~(x \in K_m(M),y\in K_{m-n}(M))~,\\[1ex]
\lambda([M]\cap a,[M] \cap b)~=~\langle a\cup b,[M]\rangle \in \ZZ~~
(a \in K^{n-m}(M),b \in K^m(M))~.\end{array}$$

\begin{remark} Ranicki \cite[\S3]{ranicki2} developed a $\pi$-equivariant $S$-duality theory, such that for a normal map $(f,b):M \to X$ there is an  isomorphism
$$\{T(\nu_{\widetilde{M}});T(\nu_{\widetilde{X}})\}_{\pi}~
\cong~\{\widetilde{X}^+,\widetilde{M}^+\}_{\pi}$$
with $\pi=\pi_1(X)$. An $n$-dimensional normal map $(f,b):M \to X$ induces a $\pi$-equivariant map $T(\widetilde{b}):T(\nu_{\widetilde{M}})
\to T(\nu_{\widetilde{X}})$ with a $\pi$-equivariant $S$-dual geometric Umkehr map $F:\Sigma^{\infty} \widetilde{X}^+ \to \Sigma^{\infty}\widetilde{M}^+$.
The algebraic theory of surgery of Ranicki \cite{ranicki2} (outlined in \S\ref{qgroups} above) obtained the surgery obstruction $\sigma_*(f,b) \in L_n(\ZZ[\pi_1(X)])$ directly from the chain level quadratic construction of a $\pi_1(X)$-equivariant geometric Umkehr map
$F:V^\infty\wedge \widetilde{X}^+ \to V^{\infty}\wedge \widetilde{M}^+$, with $\widetilde{X}$ the universal cover of $X$ and $\widetilde{M}=f^*\widetilde{X}$  the pullback cover of $M$.
We deal with the even-dimensional simply-connected case in  section
\ref{simple} and the general non-simply-connected case in section \ref{non}.  By Proposition \ref{quadquad} the geometric Hopf invariant of $F:V^{\infty} \wedge \widetilde{X}^+ \to V^{\infty} \wedge \widetilde{M}^+$ induces the quadratic construction on the chain level of \cite{ranicki2,ranicki3}
$$h_V(F)~=~\psi_F~:~H_n(X) \to Q_n(C(\widetilde{M}))~=~H_n(S(\infty)\times_{\ZZ_2}(\widetilde{M} \times_{\pi_1(X)}\widetilde{M}))~,$$
which determines a quadratic refinement $\mu$ for the intersection form $\lambda$, allowing the quadratic Poincar\'e structure to be obtained directly from  the underlying homotopy theory. See \S\ref{non} for details.
\hfill\qed
\end{remark}

\begin{remark} For a normal map $(f,b):M \to X$
it is possible to construct a geometric Umkehr map $F:V^{\infty} \wedge \widetilde{X}^+
\to V^{\infty} \wedge \widetilde{M}^+$ geometrically, as on Ranicki \cite[p.37]{ranicki3}. Replace $X$ by an $(n+k)$-dimensional manifold with boundary $(W,\partial W)$ homotopy equivalent to
$(X \times D^k,X \times S^{k-1})$, for $k \geqslant {\rm max}(3,5-n,n+1)$,
applying the $\pi$-$\pi$ Theorem of Wall \cite[3.3]{wall2}
\footnote{which applies because the Spivak normal fibration of $X$ has
a vector bundle reduction} and the Whitney embedding theorem. It is then
possible to approximate $(f,b)$ by a framed embedding
$e:V \times M \emb W\backslash \partial W$
with $V$ a $k$-dimensional inner product space/ The adjunction Umkehr map of a lift to a $\pi_1(X)$-equivariant embedding $\widetilde{e}:V \times \widetilde{M} \emb \widetilde{W}\backslash \widetilde{\partial W}$
is a geometric Umkehr map for $(f,b)$
$$F~:~\widetilde{W}/\widetilde{\partial W}~\simeq~V^{\infty} \wedge \widetilde{X}^+ \to\widetilde{W}/{\rm cl.}(\widetilde{W} \backslash \widetilde{e}(V\times \widetilde{M}))~=~V^{\infty} \wedge \widetilde{M}^+~.$$
\hfill\qed
\end{remark}

\section{The geometric Hopf invariant and the simply-connected surgery obstruction}\label{simple}

The simply-connected surgery obstruction groups are
$$L_n(\ZZ)~=~\begin{cases}
\ZZ&{\rm if}~n \equiv 0(\bmod\,4)\\
0&{\rm if}~n \equiv 1(\bmod\,4)\\
\ZZ_2&{\rm if}~n \equiv 2(\bmod\,4)\\
0&{\rm if}~n \equiv 3(\bmod\,4)
\end{cases}$$
so only the even-dimensional case $n=2m$ need be considered.
The surgery obstruction of a $2m$-dimensional normal map $(f,b):M \to X$
is
$$\sigma_*(f,b)~=~
\begin{cases}
{\rm signature}(K_m(M),\lambda,\mu)/8\\
{\rm Arf}(K_m(M;\ZZ_2),\lambda,\mu)
\end{cases}
\in L_{2m}(\ZZ)~=~\begin{cases}
\ZZ&{\rm if}~m \equiv 0(\bmod 2)\\
\ZZ_2&{\rm if}~m \equiv 1(\bmod 2)
\end{cases}$$
with $(\lambda,\mu)$  the (intersection, self-intersection) quadratic form on the middle-dimensional homology kernel
$$\begin{cases}
K_m(M)~=~{\rm ker}(f_*:H_m(M) \to H_m(X))\\[1ex]
K_m(M;\ZZ_2)~=~{\rm ker}(f_*:H_m(M;\ZZ_2) \to H_m(X;\ZZ_2))
\end{cases}$$
(cf. Proposition \ref{mu} above). We  construct $(\lambda,\mu)$ using the geometric Hopf invariant.

Let $M^{2m}$ be a $2m$-dimensional manifold with an embedding $M \subset \R^{2m+k}$, and suppose given a map $f:M \to X$ to a space $X$ with a $k$-plane bundle $\eta:X \to BO(k)$, and a bundle map $b:\nu_{M \subset \R^{2m+k}}\to \eta$.
The bundle map $b$ determines a bundle isomorphism
$$\tau_M \oplus \eta~\cong~\epsilon^{2m+k}~.$$
The normal bundle of an immersion $e:S^m \imm M$ is an $m$-plane bundle
$\nu_e:S^m \to BO(m)$ such that
$$\nu_e\oplus \tau_{S^m}~=~e^*\tau_M~:~S^m \to BO(2m)~,$$
so that
$$\nu_e \oplus \tau_{S^m} \oplus (gf)^*\eta~=~f^*(\tau_M \oplus \eta)~
\cong~\epsilon^{2m+k}~.$$
A null-homotopy $fe\simeq *:S^m \to X$ thus determines a stable isomorphism
$$\delta\nu_e~:~\nu_e \oplus \epsilon^{m+k}~\cong~\epsilon^{2m+k}~,$$
as classified by an element
$$(\delta\nu_e,\nu_e) \in \pi_m(V_{2m+k,m+k})~=~Q_{(-1)^m}(\ZZ)~.$$

\begin{proposition}~\label{embed9}
Let $(f,b):M \to X$ be a $2m$-dimensional normal map, with geometric
Umkehr map $F:\Sigma^jX^+ \to \Sigma^jM^+$.\\
{\rm (i)} The geometric Hopf invariant map
$$h_{\R^j}(F)~:~X^+ \sto S(L\R^j)^+\wedge_{\ZZ_2}
(M^+\wedge M^+)~=~
(S(L\R^j)\times_{\ZZ_2}(M \times M))^+$$
determines a quadratic refinement of the intersection pairing $\lambda:K_m(M) \times K_m(M) \to \ZZ$
$$\begin{array}{l}
\mu~:~K_m(M) \to Q_{(-1)^m}(\ZZ)~=~\ZZ/\{1+(-1)^{m+1}\}~;\\[1ex]
x \mapsto -(x \otimes x)h_{\R^j}(F)_*[M]~=~
\begin{cases} \lambda(x,x)/2&{\it if}~m \equiv 0 (\bmod\,2)\\
\langle Sq^{m+1}_{xF}(\iota),[X]\rangle&{\it if}~m \equiv 1(\bmod\,2)
\end{cases}
\end{array}$$
such that
\begin{itemize}
\item[\rm (a)]~ $\lambda(x,x)=(1+(-1)^m)\mu(x) \in \ZZ$,
\item[\rm (b)]~ $\mu(ax)=a^2\mu(x)\in Q_{(-1)^m}(\ZZ)~ (a \in \ZZ),$
\item[\rm (c)]~ $\mu(x+y)=\mu(x)+\mu(y)+\lambda(x,y)\in Q_{(-1)^m}(\ZZ)$,
\item[\rm (d)]~ if $x \in K_m(M)$ is represented by an immersion
$e:S^m \imm M$ with a null-homotopy $fe \simeq *:S^m \to X$ then
$$\begin{array}{l}
\lambda(x,x)~=~(1+(-1)^m)\mu(e)+\chi(\nu_e) \in \ZZ~,\\[1ex]
\mu(x)~=~H(\phi)~=~\mu(e)+(\nu_e,\delta\nu_e) \in Q_{(-1)^m}(\ZZ)~{\it with}\\[1ex]
H(\phi)=h_{\R^k}(\phi)[M] \in \widetilde{H}_{2m}
((S^{k-1})^+ \wedge_{\ZZ_2} (S^m \wedge S^m))=Q_{(-1)^m}(\ZZ)\\[1ex]
\hbox{\it the Hopf invariant of the stable map}\\[1ex]
\phi~:~\Sigma^kM^+ \xymatrix{\ar[r]^-{\di{E}}&} \Sigma^kT(\nu_e)
\xymatrix@C+15pt{\ar[r]^-{\di{T(\delta\nu_e)}}&} \Sigma^kT(\epsilon^m)
\xymatrix{\ar[r]&} S^{m+k}\\[1ex]
\hbox{\it where $E$ is the Umkehr stable map of $e$.}
\end{array}$$
\end{itemize}

\noindent {\rm (ii)} If $f$ is $m$-connected then $K_m(M)=\pi_{m+1}(f)$
and every element $x \in K_m(M)$ is represented by an immersion $e:S^m
\imm M$ with a null-homotopy $fe\simeq *:S^m \to X$.  There are two
special cases:
\begin{itemize}
\item[\rm (a)]~ For all $m$ it is possible to choose a
representative with $\nu_e=\epsilon^m$ and $(\nu_e,\delta\nu_e)=0$.
\item[\rm (b)]~ If $m \geqslant 3$ and $\pi_1(M)=\{1\}$ it is possible to choose a
representative with $e$ an embedding, so that $\mu(e)=0$.
\end{itemize}
 \noindent {\rm (iii)} If $f$ is $m$-connected then
$\mu(x)=0$ if (and for $m \geqslant 3$, $\pi_1(M)=\{1\}$ only if)
it is possible to kill $x \in K_m(M)$ by surgery on $(f,b)$, i.e.
to represent $x$ by a framed embedding $e:S^m \emb M$
with a null-homotopy $fe \simeq *:S^m \to X$.
{\rm (iv)} The simply-connected surgery obstruction of Kervaire and Milnor \cite{kervmiln} and  Browder \cite{browder2}
$$\begin{array}{l}
\sigma_*(f,b)~=~\begin{cases}
{\rm signature}(K_m(M),\lambda)/8\\[1ex]
{\rm Arf}(K_m(M;\ZZ_2),\lambda,\mu)
\end{cases} \\[1ex]
\hskip100pt \in L_{2m}(\ZZ)~=~
\begin{cases}\hbox{$\ZZ$}&\hbox{if $m \equiv 0(\bmod\,2)$} \\[1ex]
\hbox{$\ZZ_2$}&\hbox{if $m \equiv 1(\bmod\,2)$}
\end{cases}
\end{array}$$
is such that $\sigma_*(f,b)=0$ if {\rm (}and for
$\pi_1(X)=\{1\}$, $m \geqslant 3$ only if\,{\rm )} $(f,b)$ is normal bordant
to a homotopy equivalence.\\
\hfill\qed
\end{proposition}

 \begin{remark} \label{embed10} {\rm The quadratic refinement $\mu$ was
constructed by Browder \cite{browder2} using functional Steenrod squares, by
Wall \cite{wall2} using immersions (assuming $(f,b)$ is $m$-connected),
and by Ranicki \cite{ranicki2} using the quadratic construction $\psi_G$.  \hfill\qed}
\end{remark}

\section{The geometric Hopf invariant and codimension 2 surgery}
\label{codim2surgery}

We shall now consider the geometric Hopf invariant for maps
$F:V^{\infty} \wedge X \to V^{\infty} \wedge Y$ in the 1-dimensional case $V=\R$, which has applications to codimension 2 surgery, such as knots.

The homotopy exact sequence of the fibration
$$S^m \to BO(m) \to BO(m+1)$$
is mapped by the $J$-homomorphism to the $EHP$ sequence of $S^m$
$$\xymatrix@R+10pt{
\pi_{m+1}(BO(m+1)) \ar[r] \ar[d]^-{\di{J}} &
\pi_m(S^m) \ar[r] \ar[d]^-{\di{\cong}} &
\pi_m(BO(m)) \ar[r] \ar[d]^-{\di{J}} &
\pi_m(BO(m+1)) \ar[d]^-{\di{J}} \\
\pi_{2m+1}(S^{m+1}) \ar[r]^-{\di{H}} &
\hbox{$\ZZ$} \ar[r]^-{\di{P}}&
\pi_{2m-1}(S^m) \ar[r]^-{\di{E}} &
\pi_{2m}(S^{m+1})}$$
where
$$HJ~=~\chi~:~\pi_{m+1}(BO(m+1)) \to \ZZ$$
sends an $(m+1)$-plane bundle over $S^{m+1}$ to its Euler number (only
determined up to sign in general). The image of the Hopf invariant is
$${\rm im}(H:\pi_{2m+1}(S^{m+1}) \to \ZZ)~
=~\begin{cases}0&\hbox{if $m$ is even}\\[1ex]
\hbox{$\ZZ$}&\hbox{if $m=1,3,7$}\\[1ex]
\hbox{$2\ZZ$}&\hbox{otherwise}~~(\hbox{\rm Adams \cite{adams1}})~.
\end{cases}$$
The morphism
$$\pi_{m+1}(BO(m+1),BO(m))~=~\ZZ \to \pi_m(BO(m))~;~1 \mapsto \tau_{S^m}$$
is $\begin{cases}{\rm injective}\\[1ex] 0\end{cases}$ for
$\begin{cases}m \neq 1,3,7\\[1ex] m=1,3,7\end{cases}$ with the generator
$$(\delta \tau_{S^m},\tau_{S^m})~=~1 \in \pi_{m+1}(BO(m+1),BO(m))~=~\ZZ$$
represented by the
$\begin{cases}\hbox{\rm unique}\\[1ex]
\hbox{\rm exotic}\end{cases}$
stable trivialization
$\delta\tau_{S^m}:\tau_{S^m} \oplus \epsilon \cong \epsilon^{m+1}$.

\begin{definition}~ \label{embed11} {\rm  (Ranicki \cite[p.818]{ranicki3})
An $n$-dimensional normal map $(f,b):M \to X$ is {\it ultranormal}
if it is obtained from an $(n+1)$-dimensional manifold $W$ and a
$\ZZ[\pi_1(X)]$-homology equivalence $h:W \to X \times S^1$ by
codimension 1 transversality
$$(f,b)~=~h\vert~:~N~=~h^{-1}(X \times \{*\}) \to X~.\eqno{\hbox{\qed}}$$
}\end{definition}

\begin{example} \label{embed12} {\rm
(\cite[p.827]{ranicki3}) Let $k$ be a $(2m-1)$-knot, i.e. a (locally flat) embedding
$k:S^{2m-1} \emb S^{2m+1}$. The knot complement
$$W~=~{\rm cl.}(S^{2m+1} \backslash (k(S^{2m-1}) \times D^2))$$
is a $(2m+1)$-dimensional manifold with boundary $\partial W=S^{2m-1} \times S^1$.
Let $N^{2m} \subset S^{2m+1}$ be a Seifert surface $k$, with $\partial N=k(S^m)$.
The inclusion $f:N \emb D^{2m+2}$
defines an ultranormal map
$$(f,b)~=~h\vert~:~(N,\partial N)~=~h^{-1}((D^{2m+2},k(S^{2m-1}))\times \{*\})
\to (D^{2m+2},k(S^{2m-1}))$$
for a homology equivalence
$$h~:~(W,\partial W) \to (D^{2m+2},k(S^{2m-1})) \times S^1$$
with both $f$ and $h$ identities on boundaries.\hfill\qed}
\end{example}

\begin{remark} \label{embed13} {\rm  Every $m$-connected $2m$-dimensional
normal map can be realized as an ultranormal map
(\cite[7.8]{ranicki3}).\hfill\qed}
\end{remark}

An ultranormal map $(f,b):M \to X$ has an Umkehr stable map involving only a
single suspension
$$F~:~\Sigma X^{\infty}~=~(X\times \R)^{\infty}
\xymatrix@C+10pt{\ar[r]^-{\di{h^{-1}}}&}
\overline{W}^{\infty} \xymatrix{\ar[r]&} \Sigma M^+$$
with $\overline{W}=h^*(X \times \R)$ the infinite cyclic cover of $W$,
and $\overline{W}^{\infty} \to \Sigma M^+$ the Pontrjagin-Thom map
of an embedding $N \times \R \emb \overline{W}$.
We shall now relate the geometric Hopf map of $F$ (\ref{hopf4})
$$h_{\R}(F)~:~X^{\infty} \sto M^+ \wedge M^+$$
in the case ${\rm dim}(M)=2m$ to the Seifert form on
$$K_m(M)~=~{\rm ker}(f_*:H_m(M) \to H_m(X))$$
which was originally defined using linking numbers.
In the first instance, we consider the (self-) intersection
properties of $(m,2m,1)$-dimensional embedd\-ing-immersion pairs
$(d \times e:M^m \emb \R \times N^{2m},e:M \imm N)$.

\begin{definition}~ \label{embed14} {\rm
The {\it degree} of an $m$-plane bundle $\xi:M \to BO(m)$
and a stable isomorphism $\delta\xi:\xi \oplus \epsilon \cong \epsilon^{m+1}$
over an $m$-dimensional manifold $M$ is
$$(\xi,\delta\xi)^{\ZZ}~=~{\rm degree}(\gamma:M \to S^m) \in \ZZ$$
with $\gamma:M \to S^m$ a classifying map for $(\xi,\delta\xi)$ (\ref{link10}).\hfill\qed
}\end{definition}

\begin{proposition}~\label{embed15}
{\rm (i)} The degree of $(\xi,\delta\xi)$ determines the Euler number of $\xi$ by
$$\chi(\xi)~=~(1+(-1)^m) (\xi,\delta\xi)^{\ZZ} \in \ZZ~.$$
{\rm (ii)} For any $m$-dimensional manifold $M$ the degree defines a morphism
$$[M,S^m] \to \ZZ~;~\gamma \mapsto {\rm degree}(\gamma:M \to S^m)~.$$
For connected $M$ this an isomorphism (by the Hopf-Whitney theorem) and the degree classifies pairs $(\xi,\delta\xi)$ as in \ref{embed14}. \\
{\rm (iii)} The degree is the linking number of the submanifolds
$$\begin{array}{l}
i_0~:~M \emb E(\epsilon^{m+1})~=~M \times \R^{m+1}~;~
x \mapsto (x,0)~,\\[1ex]
i_1~:~M \emb E(\epsilon^{m+1})~=~M \times \R^{m+1}~;~
x \mapsto E(\delta\xi)(x,0,1)~,\end{array}$$
which is the degree of the composite
$$\gamma~:~M \xymatrix{\ar[r]^{\di{i_1}}&}
(M \times \R^{m+1})\backslash i_0(M)~=~M \times (\R^{m+1}\backslash \{0\})
\xymatrix{\ar[r]&} \R^{m+1}\backslash \{0\}~\simeq~S^m~,$$
or equivalently the algebraic
number of elements in the 0-dimensional linking manifold
$L(M,\xi,\delta\xi)=\gamma^{-1}(*)$ (\ref{link11}).\hfill\qed
\end{proposition}

\begin{proposition}~ \label{embed16} Let $(d\times e:M^m
\emb \R \times N^{2m},e:M \imm N)$ be an
$(m,2m,1)$-dimensional embedding-immersion pair with $d\times e$
framed, so that $\nu_e:M \to BO(m)$ has a 1-stable trivialization
$$\delta\nu_e~:~\nu_e\oplus\epsilon~\cong~\epsilon^{m+1}~.$$
Let $E:\Sigma N^+ \to \Sigma T(\nu_e)$ be an Umkehr stable map for $e$,
so that the homotopy class of the Pontrjagin-Thom map
$$\phi~:~\Sigma N^+ \xymatrix{\ar[r]^-{\di{E}}&}
\Sigma T(\nu_e)\xymatrix@C+20pt{\ar[r]^-{\di{T(\delta\nu_e)}}&}
T(\epsilon^{m+1}) \to S^{m+1}$$
is the framed cobordism class of
$$\phi^{-1}(*)~=~M \subset \R \times N~.$$
{\rm (i)} The Hopf invariant map $h_{\R}(\phi):\Sigma^2N^+ \to S^{2m+2}$
is the Pontrjagin-Thom map for the framed 0-dimensional submanifold
$$h_{\R}(\phi)^{-1}(*)~=~L(M,N) \cup L(M,\nu_e,\delta\nu_e)
\subset \R \times \R \times N~,$$
with degree
$$H(\phi)~=~{\rm degree}(h_{\R}(\phi))~=~\mu^{\ZZ}(e) +
(\delta \nu_e,\nu_e)^{\ZZ} \in \ZZ~.$$
{\rm (ii)} If $e:M \imm N$ is framed, and $\delta\nu_e$
is the stabilization of the framing $\nu_e \cong \epsilon^m$,
then $L(M,\nu_e,\delta \nu_e)=\emptyset$ and
$$(\nu_e,\delta\nu_e)^{\ZZ}~=~0~~,~~H(\phi)~=~\mu^{\ZZ}(e) \in \ZZ~.$$
{\rm (iii)} If $e:M \imm N$ is an embedding then
$L(M,N)=\emptyset$ and
$$\mu^{\ZZ}(e)~=~0~~,~~H(\phi)~=~(\nu_e,\delta \nu_e)^{\ZZ} \in \ZZ~.$$
\end{proposition}
\noindent{\it Proof\  .} (i) This is the special case $n=2m$ of
\ref{link13}.\\
(ii)+(iii) Immediate from (i).\hfill\qed

\begin{example}\label{embed17}
{\rm (i) For every immersion $e:S^m \imm S^{2m}$ there exists a
map $d:S^m \to \R$ such that $d \times e:S^m \to \R \times
S^{2m}$ is a framed embedding, with the stable trivialization
$$\delta\nu_e~:~\nu_e \oplus \epsilon~\cong~\epsilon^{m+1}$$
such that the Pontrjagin-Thom map
$$\phi~:~\Sigma (S^{2m})^{\infty} \xymatrix{\ar[r]^-{\di{E}}&}
\Sigma T(\nu_e)\xymatrix@C+20pt{\ar[r]^-{\di{T(\delta\nu_e)}}&}
T(\epsilon^{m+1}) \xymatrix{\ar[r]&} S^{m+1}$$
is null-homotopic,
or equivalently such that $d\times e:S^m \emb \R
\times S^{2m}$ is framed null-cobordant (by $D^{m+1} \subset \R
\times S^{2m}$ for $m \geqslant 2$, and by a Seifert surface $N^2
\subset \R \times S^2$ for $m=1$). For this preferred choice of
framing \ref{embed16} (i) gives
$$H(\phi)~=~\mu^{\ZZ}(e) + (\delta \nu_e,\nu_e)^{\ZZ}~=~0 \in \ZZ~.$$
(ii) The Whitney immersion $f:S^m \imm S^{2m}$ has
$$(\delta\nu_e,\nu_e)^{\ZZ}~=~-(\delta\tau_{S^m},\tau_{S^m})^{\ZZ}~=~-1 \in
\pi_{m+1}(BO(m+1),BO(m))~=~\ZZ$$
with
$$\mu^{\ZZ}(e)~=~1~,~
H(\phi)~=~\mu^{\ZZ}(e) + (\delta \nu_e,\nu_e)^{\ZZ}~=~0 \in \ZZ~.$$
In particular, $f:S^1 \imm S^2$ is just the figure 8 immersion.\\
(iii) The integral geometric self-intersection of a framed immersion
$f:M^m \imm S^{2m}$ is identified by \ref{embed16} (ii) with the Hopf invariant of
the Pontrjagin-Thom map $\phi:\Sigma(S^{2m})^{\infty} \to S^{m+1}$
$$\mu^{\ZZ}(e)~=~H(\phi) \in
{\rm im}(H:\pi_{2m+1}(S^{m+1})\to \ZZ)~=~
\begin{cases}0&\hbox{if $m$ is even}\\[1ex]
\hbox{$\ZZ$}&\hbox{if $m=1,3,7$}\\[1ex]
\hbox{$2\ZZ$}&\hbox{otherwise~.}
\end{cases}$$
Thus the number of self-intersections is even, except in the Hopf
invariant 1 dimensions $m=1,3,7$.\hfill\qed}
\end{example}

Let $M^{2m} \subset S^{2m+1}$ be a codimension 1 framed submanifold, and
let $i:\R \times M \emb S^{2m+1}$ be the inclusion
of a tubular neighbourhood. Use the disjoint embeddings
$$\begin{array}{l}i_0~:~M \emb S^{2m+1}~;~x \mapsto i(0,x)~,\\[1ex]
i_1~:~M \emb S^{2m+1}~;~x \mapsto i(1,x)\end{array}$$
to define a pairing
$$s~:~M \times M \to S^{2m}~;~
(x,y) \mapsto {i_0(x) - i_1(y)\over  \Vert i_0(x) - i_1(y) \Vert}~.$$
The induced product in homology
$$\sigma~:~H_m(M) \times H_m(M) \to H_{2m}(S^{2m})~=~\ZZ~;~
(x,y) \mapsto s_*(x \otimes y)$$
is such that for any embeddings $e,e':S^m \emb M$
$$\begin{array}{l}\sigma(e_*[S^m],e'_*[S^m])~
=~\hbox{\rm linking number}
\big(i_0e \sqcup i_1e':S^m \sqcup S^m \emb S^{2m+1} \big)\\[1ex]
=~\hbox{\rm degree}(S^m \xymatrix{\ar[r]^-{\di{i_0e}}&} S^{2m+1}
\backslash i_1e'(S^m) \simeq S^m)  \in \ZZ~.\end{array}$$
For any embedding $e:S^m \emb M$ (assumed unknotted in $S^{2m+1}$
if $m=1$) the composite $i_0e:S^m \emb S^{2m+1}$ is isotopic
to the standard embedding $S^m \emb S^{2m+1}$, so that the
normal bundle $\nu_e:S^m \to BO(m)$ is equipped with a stable trivialization
$\delta\nu_e:\nu_e \oplus \epsilon \cong \epsilon^{m+1}$.
Define the disjoint embeddings
$$\begin{array}{l}
j_0~:~S^m \emb E(\nu_e \oplus \epsilon)~=~E(\nu_e) \times \R~;~
x \mapsto (x,0,0)~,\\[1ex]
j_1~:~S^m \emb E(\nu_e \oplus \epsilon)~=~
E(\nu_e) \times \R~;~x \mapsto (x,0,1)~.\end{array}$$
It follows from the commutative diagram
$$\xymatrix@R+10pt{
S^m \ar[r]^-{\di{j_1}} \ar[d]^-{\di{i_1e}}&
E(\nu_e\oplus \epsilon)\backslash j_0(S^m) \cong
E(\epsilon^{m+1})\backslash E(\delta\nu_e)^{-1}j_0(S^m) \ar[d]\\
S^{2m+1} \backslash M \ar[r] & S^{2m+1}\backslash i_0e(S^m)
}$$
that
$$\begin{array}{l}\sigma(e_*[S^m],e_*[S^m])~=~
\hbox{\rm linking number}(i_0e(S^m) \cup i_1e(S^m) \subset S^{2m+1})\\[1ex]
=~(\delta\nu_e,\nu_e)^{\ZZ} \in \pi_{m+1}(BO(m+1),BO(m))~=~\ZZ~.\end{array}$$
\indent
The above construction applies also to the embedding of a Seifert surface
$M^{2m} \emb S^{2m+1}$ for a $(2m-1)$-knot
$k:S^{2m-1} \emb S^{2m+1}$, with
$$\partial M~=~k(S^{2m-1}) \emb S^{2m+1}~.$$
If $x_1,x_2,\dots,x_{\ell} \in H_m(M)$ is a basis for the f.g. free
$\ZZ$-module $H_m(M)/{\rm torsion}$ then $(\sigma(x_i,x_j))$
is the Seifert matrix of $k$ with respect to $M$.

Let $M^{2m}$ be a $2m$-dimensional manifold with an embedding
$i:\R \times M\emb S^{2m+1}$, as before. Let
$$(d \times e:S^m \emb \R \times M~,~e:S^m \imm M)$$
be an $(m,2m,1)$-dimensional embedding-immersion pair, with Umkehr stable map
$E:\Sigma M^+ \to \Sigma T(\nu_e)$ and integral self-intersection
number
$$\mu^{\ZZ}(e)~=~\psi_{\R}(E)_*[M] \in \ZZ~.$$
The composite $i(d \times e):S^m \emb S^{2m+1}$ is an
embedding with normal bundle
$$\nu_{i(d \times e)}~=~\nu_e \oplus \epsilon~:~S^m \to BO(m+1)~.$$
An isotopy of $i(d\times e)$ to the standard embedding
$S^m \emb S^{2m+1}$ determines a stable trivialization
$$\delta\nu_e~:~\nu_e \oplus \epsilon ~\cong~ \epsilon^{m+1}$$
which is classified by the Hopf invariant of the homotopy equivalence
$$T(\delta\nu_e)~:~\Sigma T(\nu_e)~\simeq~S^{2m+1} \vee S^{m+1}~,$$
that is
$$(\delta\nu_e,\nu_e)^{\ZZ}~=~\psi_{\R}(T(\delta\nu_e))[T(\nu_e)] \in
\pi_{m+1}(BO(m+1),BO(m))~=~\ZZ~.$$
(For $m\neq 1,3,7$ this is independent of $d:S^m \to \R$, and
is just the integer $c \in \ZZ$ such that
$\nu_e=c\tau_{S^m} \in {\rm ker}(\pi_m(BO(m)) \to \pi_m(BO(m+1)))~.)$
Let $F:S^{2m+1} \to \Sigma M^+$ be the Pontrjagin-Thom map for
$i:\R \times M \emb S^{2m+1}$. (If $M$ is a Seifert surface
of a knot $k:S^{2m-1} \emb S^{2m+1}$ then $F$ is the Umkehr stable map
for the ultranormal normal map
$$(f,b)~=~{\rm inclusion}~:~(M,\partial M) \to (D^{2m+2},k(S^{2m-1}))$$
of 6.12). The composite
$$T(\delta\nu_e)EF~:~S^{2m+1} \xymatrix{\ar[r]^{\di{F}}&}
\Sigma M^+ \xymatrix{\ar[r]^{\di{E}}&}
\Sigma T(\nu_e)\xymatrix{\ar[r]^{\di{T(\delta\nu_e)}}&}S^{2m+1} \vee S^{m+1}$$
is the Pontrjagin-Thom map of the standard embedding $S^m \emb S^{2m+1}$,
so that
$$T(\delta\nu_e)EF~\simeq~{\rm inclusion}~:~S^{2m+1} \to S^{2m+1} \vee S^{m+1}$$
and by the composition formula \ref{quad3} (vi)
$$\begin{array}{l}\psi_{\R}(T(\delta\nu_e)FG)[S^{2m}]\\[1ex]
=~(F\otimes F)_*\psi_{\R}(F)_*[S^{2m}]+\psi_{\R}(E)_*[M]+\psi_{\R}(T(\delta\nu_e))[T(\nu_e)]\\[1ex]
=~0 \in \ZZ~.\end{array}$$
The evaluation of the ultraquadratic construction on $F$ (\ref{ultraquadratic})
$$\psi_{\R}(F)~:~S^{2m} \sto M^+ \wedge M^+$$
on $[S^{2m}]$ is a homology class
$$h'_{\R}(F)_*[S^{2m}]\in\widetilde{H}_{2m}(M^+ \wedge M^+)~
=~H_{2m}(M \times M)$$
such that
$$\sigma~:~H_m(M) \times H_m(M) \to \ZZ~;~(x,y) \mapsto
-\langle x^* \times y^*,h'_{\R}(F)_*[S^{2m}]\rangle$$
with $x^*,y^* \in H^m(M)$ the Poincar\'e duals of $x,y \in H_m(M)$.
Thus
$$\begin{array}{ll}
\sigma(e_*[S^m],e_*[S^m])&
=~-(F \otimes F)_*h'_{\R}(F)_*[S^{2m}]\\[1ex]
&=~h'_{\R}(F)_*[M]+h'_{\R}(T(\delta\nu_e))[T(\nu_e)]\\[1ex]
&=~\mu^{\ZZ}(e)+(\delta\nu_e,\nu_e)^{\ZZ} \in \ZZ\end{array}$$
In particular, if $x=e_*[S^m] \in H_m(M)$ for an embedding
$e:S^m \emb M$ then
$$\mu^{\ZZ}(e)~=~0~~,~~\sigma(x,x)~=~(\delta\nu_e,\nu_e)^{\ZZ} \in \ZZ~.$$
On the other hand, if $x=e_*[S^m] \in H_m(M)$ for a framed immersion
$e:S^m \imm M$ then
$$(\delta\nu_e,\nu_e)^{\ZZ}~=~0~~,~~\sigma(x,x)~=~\mu^{\ZZ}(e) \in \ZZ~.$$
Mow suppose that $M$ is $(m-1)$-connected. Every $x \in H_m(M)$ is
represented by a framed immersion $e:S^m \imm M$. If $m \geqslant
3$ then every $x \in H_m(M)$ is represented by an embedding
$e:S^m \emb M$ (which will not in general be framed).
An element $x \in H_m(M)$ is such that $\sigma(x,x)=0$ if (and for
$m \geqslant 3$ only if) it can be represented by a framed
embedding $e:S^m \emb M$, in which case $x$ can be
killed by ambient surgery to obtain an ambient cobordant
submanifold
$$(M \emb S^{2m+1})  \mapsto (M' \emb S^{2m+1})$$
with
$$M'~=~{\rm cl.}(M \backslash (e(S^m) \times D^m)) \cup D^{m+1} \times S^{m-1}$$

\begin{remark} \label{embed18}
{\rm See Levine \cite{levine} for the classification in the case $m
\geqslant 3$ of  $(2m-1)$-knots $k:S^{2m-1} \emb
S^{2m+1}$ which are simple, i.e. which admit an $(m-1)$-connected
Seifert surface $M^{2m} \emb S^{2m+1}$. The isotopy
classes of such knots are in one-one correspondence with the
$S$-equivalence classes of Seifert matrices
$\sigma=(\sigma_{ij})_{1 \leqslant i,j \leqslant \ell}$ with
$\sigma_{ij} \in \ZZ$ and $(\sigma_{ij}+(-1)^m\sigma_{ji})$
invertible. The knot $k$ corresponding to $\sigma$ has an
$(m-1)$-connected Seifert surface $M^{2m} \emb
S^{2m+1}$ with a handle presentation involving $\ell$ $m$-handles
$$M~=~D^{2m} \cup \bigcup_{\ell}(D^m \times D^m)~.$$
The attaching maps are embeddings
$\phi_i:S^{m-1} \times D^m \emb S^{2m-1}$ with
$$\begin{array}{l}\hbox{\rm linking number}(
\phi_i(S^{m-1} \times \{0\}) \cup \phi_j(S^{m-1} \times \{0\})
\subset S^{2m-1})\\[1ex] \hskip100pt =~\sigma_{ij}+(-1)^m\sigma_{ji}
\in \ZZ~~(1 \leqslant  i,j \leqslant  \ell,~ i \neq j)~.
\end{array}$$
The basis elements
$$x_i~=~(0,\dots,0,1,0,\dots,0) \in H_m(M)~=~\ZZ^{\ell}~~(1 \leqslant  i \leqslant  \ell)$$
are represented by $(m,2m,1)$-dimensional embedding-immersion pairs
$$(e_i\times f_i:S^m \emb \R \times M^{2m}~,~f_i:S^m \imm N)$$
with $e_i \times f_i$ framed and $f_i$ an embedding, such that
$$(\delta\nu_{f_i},\nu_{f_i})~=~\sigma_{ii}(\delta\tau_{S^m},\tau_{S^m}) \in
\pi_{m+1}(BO(m+1),BO(m))~=~\ZZ~.$$ \hfill\qed} \end{remark}

\begin{proposition}~ \label{sym10}
{\it
Let $N^n \subset M^{n+1}$ be a framed codimension 1 submanifold with
complement $P$, so that
$$M~=~N \times I \cup_{i_0,i_1} P$$
and there is defined a cofibration sequence
$$P^+ \xymatrix{\ar[r]&} M^+
\xymatrix{\ar[r]^{\di{F}}&} \Sigma N^+
\xymatrix{\ar[r]^{\di{G}}&} \Sigma P^+$$
with
$$\begin{array}{l}
F~=~projection~:~M^+ \to M/P ~=~N\times I/\partial (N \times I)~=~\Sigma N^+~,\\[1ex]
G~=~\Sigma i_0 - \Sigma i_1~:~\Sigma N^+ \to \Sigma P^+\end{array}$$
where $i_0,i_1:N \to P$ are the inclusions of the two boundary components.
The spectral Hopf invariant of $F$, the Hopf invariant of $G$ and $i_0,i_1$
are related by a stable homotopy commutative diagram}
$$\xymatrix@C+25pt@R+25pt{
\Sigma M^+ \ar[r]^-{\di{sh_{\R}(F)}} \ar[d]_-{\di{\Sigma F}}
& \Sigma P^+ \wedge \Sigma P^+ \\
\Sigma^2 N^+ \ar[ur]^-{\di{-h_{\R}(G)}} \ar[r]^-{\di{\Sigma^2\Delta_N}} &
\Sigma N^+ \wedge \Sigma N^+ \ar[u]_-{\di{G \wedge \Sigma i_1}}}$$
\end{proposition}
\noindent{\it Proof\  .} By \ref{sym8} we have the stable homotopy identity
$$sh_{\R}(F)+h_{\R}(G)\Sigma F~=~0~:~\Sigma M^+ \to \Sigma P^+ \wedge \Sigma P^+~.$$
By the sum formula \ref{hopf7} (viii)
$$h_{\R}(G)~=~h_{\R}(\Sigma i_0)+h_{\R}(-\Sigma i_1)+
(\Sigma i_0 \wedge -\Sigma i_1)\Sigma^2 \Delta_N~:~
\Sigma^2N^+ \to \Sigma P^+ \wedge \Sigma
P^+~.$$ By the suspension formula \ref{hopf7} (iii)
$$h_{\R}(\Sigma i_0)~=~0~:~\Sigma^2N^+ \to \Sigma P^+ \wedge \Sigma P^+~.$$
By Example \ref{hopf5} the Hopf invariant map of
$$-1~:~S^1 \to S^1~;~t \mapsto 1-t$$
is given by
$$h_{\R}(-1:S^1 \to S^1)~=~1~:~S^2 \to S^2~.$$
By the product formula \ref{hopf7} (vi)
$$\begin{array}{ll}
h_{\R}(-\Sigma i_1)&=~h_{\R}(-1 \wedge i_1 : S^1 \wedge N^+ \to S^1 \wedge P^+)\\[1ex]
&=~(\Sigma i_1 \wedge \Sigma i_1)\Sigma^2 \Delta_N~:~\Sigma^2N^+ \to
\Sigma P^+ \wedge \Sigma P^+~.\end{array}$$
Substitution in the expression for $h_{\R}(G)$ gives
$$\begin{array}{ll}
h_{\R}(G)&=~
((\Sigma i_0-\Sigma i_1)\wedge -\Sigma i_1)\Sigma^2 \Delta_N\\[1ex]
&=~-(G\wedge \Sigma i_1)\Sigma^2M^+:~\Sigma^2N^+ \to
\Sigma P^+ \wedge \Sigma P^+~.\end{array}$$
\hfill\qed

\begin{remark} \label{sym11}
{\rm Let $k:S^{n-1} \subset S^{n+1}$ be
an $(n-1)$-knot, and let $M^n \subset S^{n+1}$ be a Seifert
surface, so that $\partial M=k(S^{n-1})$. The knot complement is a
codimension 0 submanifold
$$X~=~{\rm cl.}(S^{n+1} \backslash (k(S^{n-1}) \times D^2)) \subset S^{n+1}$$
which can be cut along $M \subset X$ to obtain a codimension 0 submanifold
$$P~=~{\rm cl.}(X \backslash (M \times I)) \subset X$$
such that
$$\begin{array}{l}\partial P~=~M \cup_{\partial M}M~~,~~
X~=~(M \times I)\cup_{i_0,i_1} P~,\\[1ex]
S^{n+1}~=~(M \times I \cup \partial M \times D^2)\cup_{i_0,i_1} P\end{array}$$
with $i_0,i_1:M \to P$ the two inclusions.
As in \ref{sym10} there is defined a homotopy cofibration sequence
$$P^+ \xymatrix{\ar[r]&} S^{n+1} \xymatrix{\ar[r]^{\di{F}}&}
\Sigma (M/\partial M)
\xymatrix{\ar[r]^{\di{G}}&} \Sigma P^+$$
such that up to homotopy
$$\Sigma i_0 -\Sigma i_1~:~\Sigma M^+ \xymatrix{\ar[r]&} \Sigma (M/\partial M)
\xymatrix{\ar[r]^{\di{G}}&} \Sigma P^+~,$$
and up to stable homotopy
$$sh_{\R}(F)~=~(G \wedge \Sigma i_1)(\Sigma^2 \Delta_{N/\partial N})
(\Sigma F)~:~S^{n+2} \to \Sigma P^+\wedge \Sigma P^+~.$$
The expression of the spectral Hopf invariant
$sh_{\R}(F)$ of the unstable normal invariant $F$ of $M^n \subset
S^{n+1}$ in terms of the homotopy-theoretic Seifert form $i_1$ is
a generalization of Theorem 3.1 of Richter \cite{richter}. \hfill\qed}
\end{remark}

\section{The geometric Hopf invariant and the non-simply-connected surgery obstruction}\label{non}

The surgery obstruction group $L_n(A)$ was defined by Wall [30]
for any ring with involution $A$ to be the Witt group of
$(-1)^k$-quadratic forms over $A$ for $n=2k$, and the Witt group
of automorphisms of $(-1)^k$-quadratic forms for $n=2k+1$.
In the applications to non-simply-connected surgery obstruction theory
$A=\ZZ[\pi]$ is a group ring, with the involution
$$A \to A~;~\sum\limits_{g \in \pi} n_g g \mapsto \sum\limits_{g \in \pi}
n_g w(g)g^{-1}~(n_g \in \ZZ)$$
for an orientation morphism $w:\pi \to \ZZ_2=\{\pm 1\}$ ($w=1$ in the oriented case).  The surgery obstruction $\sigma_*(f,b) \in L_n(\ZZ[\pi_1(X)])$ of an
$n$-dimensional normal map $(f,b):M \to X$ was constructed by
first making $(f,b)$ $[n/2]$-connected by surgery below the middle
dimension, and then using the middle-dimensional quadratic
$\pi_1(X)$-valued self-intersection data on the kernel homology
$\ZZ[\pi_1(X)]$-modules
    $$K_i(M)~=~{\rm ker}(\widetilde{f}_*:H_i(\widetilde{M}) \to H_i(\widetilde{X}))$$
with $\widetilde{f}:\widetilde{M} \to \widetilde{X}$ a
$\pi_1(X)$-equivariant lift of $f$ to the universal cover
$\widetilde{X}$ of $X$ and the pullback cover
$\widetilde{M}=f^*\widetilde{X}$ of $M$. The algebraic $L$-group
$L_n(A)$ was interpreted in Ranicki \cite{ranicki2} as the cobordism group of
$n$-dimensional quadratic Poincar\'e complexes over $A$. The
surgery obstruction $\sigma_*(f,b)$ was represented in \cite{ranicki2} by the
cobordism class of an $n$-dimensional quadratic Poincar\'e complex
$(C,\psi)$ over $\ZZ[\pi_1(X)]$, with $H_*(C)=K_*(M)$
the kernel $\ZZ[\pi_1(X)]$-modules, and $\psi$
obtained from a $\pi_1(X)$-equivariant stable Umkehr map
$F:\Sigma^{\infty}\widetilde{X}^+ \to
\Sigma^{\infty}\widetilde{M}^+$ $S$-dual to the induced map of
$\pi_1(X)$-equivariant Thom spaces $T(\widetilde{b}):
T(\nu_{\widetilde{M}}) \to T(\nu_{\widetilde{X}})$ by a chain level quadratic
construction capturing the self-intersection data, namely the
`$\pi_1(X)$-equivariant quadratic construction' $\psi$.  The chain level
$\psi$ generalized the functional Steenrod squares used by Browder
\cite{browder2} in the homotopy-theoretic construction of the
simply-connected surgery obstruction. In \cite{ranicki2} $\psi$ was
generalized to $\pi_1(X)$-equivariant spectral quadratic and ultraquadratic
constructions. In particular:
\begin{itemize}
\item[(a)]~The $\pi_1(X)$-equivariant homotopy-theoretic expression for the
double points of an immersion $f:M^m \imm N^{2m}$
$$\mu(f) \in \ZZ[\pi_1(N)]/\{x - (-)^m\overline{x}\,\vert\,x \in \ZZ[\pi_1(N)]\}$$
of Wall \cite[\S5]{wall2}.

\item[(b)]~The $\pi_1(X)$-equivariant homotopy-theoretic expression for the
quadratic structure $$\psi_b \in Q_n(\Cc(f^!))$$
on the $\ZZ[\pi_1(X)]$-module chain complex kernel
of an $n$-dimensional normal map $(f,b):M \to X$ of Ranicki \cite{ranicki2}.
\end{itemize}

However, the proofs in \cite{ranicki2} were somewhat indirect: Proposition \ref{quadquad} gives a direct proof in the simply-connected case, and
Appendix \ref{appendix1} gives the non-simply-connected case.

\section{The geometric Hopf invariant and the total surgery obstruction}\label{total}

The Browder-Novikov-Sullivan-Wall surgery obstruction theory studies the `$n$-dimensional smooth structure set' $\SS^O(X)$ of a topological space $X$, and its topological analogue
$\SS^{TOP}(X)$. By definition, $\SS^O(X)$  is the set of equivalence classes of pairs
$$\hbox{( smooth $n$-dimensional manifold $M$ , homotopy equivalence $h:M \to X$ )}$$
with
$$\hbox{$(M,h) \simeq (M',h')$ if $h'^{-1}h:M \to M'$ is homotopic to a diffeomorphism.}$$
The structure set $\SS^O(X)$ is non-empty if and only if $X$ is homotopy equivalent to an $n$-dimensional smooth manifold.
The original 1960's surgery theory provided a two-stage obstruction in the case $n \geqslant 5$ for deciding if $\SS^O(X)$ is non-empty, and in this case obtained the `surgery exact sequence' of pointed sets
$$\dots \to L_{n+1}(\ZZ[\pi_1(X)]) \to \SS^O(X) \to [X,G/O] \to L_n(\ZZ[\pi_1(X)])$$
with $G/O$ the homotopy fibre of the forgetful map $BO \to BG$.
Let $X$ be an $n$-dimensional geometric Poincar\'e complex with
Spivak normal structure
$(\nu_X:X \to BG(k),\rho_X:S^{n+k} \to T(\nu_X))$ ($k$ large). The primary obstruction to $\SS^O(X)$ being non-empty is the homotopy class of the composite
$$\xymatrix{X \ar[r]^-{\nu_X} & BG(k) \ar[r] & B(G/O)~,}$$
which is zero if and only if $\nu_X$ admits a vector bundle reduction
$\eta:X \to BO(k)$ ($k$ large). If this obstruction vanishes choose a reduction, and apply Pontrjagin-Thom transversality to $\rho_X:S^{n+k} \to T(\eta)=T(\nu_X)$ to obtain a normal map from a smooth $n$-dimensional  manifold
$$(f,b)~=\rho_X \vert~:~M^n~=~\rho_X^{-1}(X) \to X$$
with a non-simply-connected surgery obstruction $\sigma_*(f,b) \in L_n(\ZZ[\pi_1(X)])$.
For $n \geqslant 5$ $\sigma_*(f,b)=0$ if and only if $(f,b)$ is normal bordant to a homotopy equivalence. A different choice of reduction $\eta$ will give a different surgery obstruction:  $X$ is homotopy equivalent to a smooth manifold if and only if there exists a reduction $\eta$ with zero surgery obstruction. See Ranicki \cite{ranicki6} for a general introduction to smooth surgery theory.

The Browder-Novikov-Sullivan-Wall theory was originally developed for smooth manifolds, but since the 1970 breakthrough of Kirby and Siebenmann \cite{kirbysiebenmann} there is also a version for topological manifolds.  Technically, a topological manifold may not be a geometric Poincar\'e complex, because it may not be a finite $CW$ complex,
but it is homotopy equivalent to one.  The topological structure set $\SS^{TOP}(X)$ is defined by analogy with $\SS^O(X)$, but using topological manifolds, which has better algebraic properties. Again, for $n \geqslant 5$ there is a two-stage obstruction theory to the existence, with an exact sequence of pointed sets
$$\dots \to L_{n+1}(\ZZ[\pi_1(X)]) \to \SS^{TOP}(X) \to [X,G/TOP] \to L_n(\ZZ[\pi_1(X)])$$
with $G/TOP$ the homotopy fibre of the forgetful map $BTOP \to BG$, such that
$$\pi_*(G/TOP)~=~L_*(\ZZ)~.$$
However, there is an essential difference between the smooth and topological theories. The `total surgery obstruction' of Ranicki \cite{ranicki0,ranicki3,ranicki4} unites the two obstructions in the topological category, and if non-empty $\SS^{TOP}(X)$ has the structure of an abelian group\footnote{There is no smooth analogue of the total surgery obstruction, and $\SS^O(X)$ does not have the structure of an abelian group, by Crowley \cite{crowley}.}.

An $n$-dimensional geometric normal complex in the sense of Quinn \cite{quinn}
$$(\nu_X:X \to BG(k),\rho_X:S^{n+k} \to T(\nu_X))~(k~{\rm large})$$
is a space $X$ together with  a $(k-1)$-spherical fibration $\nu_X$ and a map $\rho_X$. (We shall only be considering spaces $X$ which are finite complexes, or at worst have the homotopy type of a finite complex). The
composite
$$\xymatrix@C+20pt{
\pi_{n+k}(T(\nu_X))  \ar[r]^-{\displaystyle{\rm Hurewicz}} &
\widetilde{H}_{n+k}(T(\nu_X)) \ar[r]^-{\displaystyle{\rm Thom}}_-{\displaystyle{\cong}} & H_n(X)}$$
sends $\rho_X$ to the fundamental class $[X] \in H_n(X)$, using twisted coefficients in the non-oriented case.
An $n$-dimensional geometric Poincar\'e complex $X$ in the sense of Wall \cite{wallpoincare} is essentially the same as an
$n$-dimensional geometric normal complex $(X,\nu_X,\rho_X)$
such that the $\ZZ[\pi_1(X)]$-module chain map
$$[X] \cap - ~:~ C(\widetilde{X})^{n-*}~=~
{\rm Hom}_{\ZZ[\pi_1(X)]}(C(\widetilde{X}),\ZZ[\pi_1(X)])_{n-*} \to C(\widetilde{X})$$
is a chain equivalence, with $\widetilde{X}$ the universal cover of $X$.

An $n$-dimensional topological manifold $M$ is homotopy equivalent to an $n$-dimensional geometric Poincar\'e complex (but may not actually be a $CW$ complex), with $\nu_M$ the
sphere bundle of the topological normal block bundle
$\widetilde{\nu}_M:M \to B\widetilde{TOP}(k)$ of an embedding
$M \subset S^{n+k}$, and $\rho_M:S^{n+k} \to T(\nu_M)=T(\widetilde{\nu}_M)$ topologically transverse at
the zero section $M \subset T(\widetilde{\nu}_M)$ with
$$\rho_M\vert~=~{\rm id.}~:~(\rho_M)^{-1}(M)=M \to M~.$$
\indent
The primary obstruction is the
homotopy class $t(\nu_X) \in [X,B(G/TOP)]$ of $\nu_X$, such that
$t(\nu_X)=0$ if and only if $\nu_X$ admits a topological reduction
$\eta:X \to BTOP(k)$. If the primary obstruction is zero,
use a choice of $\eta$ to make $\rho_X:S^{n+k} \to T(\nu_X)=T(\eta)$ topologically transverse at the zero
section $X \subset T(\nu_X)=T(\widetilde{\nu}_X)$ (by Kirby and Siebenmann \cite{kirbysiebenmann}) obtaining a normal map
$$(f,b)~=~\rho_X\vert~:~M^n~=~(\rho_X)^{-1}(X) \to X$$
with $M$ an $n$-dimensional topological manifold.
The surgery obstruction $\sigma_*(f,b) \in L_n(\ZZ[\pi_1(X)])$  is the secondary obstruction (which depends on the choice of $\widetilde{\nu}_X$) such that $\sigma_*(f,b)=0$ if (and for $n \geqslant 5$ only if) $(f,b)$ is normal bordant to a homotopy equivalence.

The total surgery obstruction of an $n$-dimensional geometric Poincar\'e
complex $X$ is the invariant $s(X)\in \SS_n(X)$ introduced in Ranicki \cite{ranicki0} such that $s(X)=0 \in \SS_n(X)$ if (and for $n \geqslant 5$
only if) $X$ is homotopy equivalent to an $n$-dimensional topological
manifold. The invariant takes value in the relative group $\SS_n(X)$ of the assembly map
$A$ from the generalized $\LL_{\bullet}(\ZZ)$-homology groups to the surgery obstruction groups
$$\xymatrix@C-10pt{
\dots \ar[r] & H_n(X;\LL_{\bullet}(\ZZ)) \ar[r]^-{A} &
L_n(\ZZ[\pi_1(X)]) \ar[r] &\SS_n(X) \ar[r]^-{\partial} & H_{n-1}(X;
\LL_{\bullet}(\ZZ)) \ar[r] & \dots}$$
with $\LL_{\bullet}(\ZZ)$ the 1-connective quadratic $\LL$-spectrum of $\ZZ$ such that for $n \geqslant 1$
$$\pi_n(\LL_{\bullet}(\ZZ))~=~L_n(\ZZ)~=~
\begin{cases}
\ZZ&{\rm if}~n \equiv 0(\bmod\,4)\\
0&{\rm if}~n \equiv 1(\bmod\,4)\\
\ZZ_2&{\rm if}~n \equiv 2(\bmod\,4)\\
0&{\rm if}~n \equiv 3(\bmod\,4)~.
\end{cases}$$
The image
$$[s(X)]=t(X) \in H_{n-1}(X;\LL_{\bullet}(\ZZ))$$
is such that
$t(X)=0$ if and only if $\nu_X$ admits a topological reduction $\eta$, i.e. it is the primary obstruction to $X$ being homotopy equivalent to a topological manifold. If $t(X)=0$ then
choosing a reduction $\eta$ topological transversality gives a normal map
$(f,b):M \to X$ (as above) and the surgery obstruction $\sigma_*(f,b) \in L_n(\ZZ[\pi_1(X)])$ has image
$$\begin{array}{l}
[\sigma_*(f,b)]~=~s(X) \\[1ex]
\in {\rm im}(L_n(\ZZ[\pi_1(X)]) \to \SS_n(X))~=~
{\rm ker}(\SS_n(X) \to H_{n-1}(X;\LL_{\bullet}(\ZZ)))~.
\end{array}
$$
If $X$ is an $n$-dimensional topological manifold for $n \geqslant 5$ there is defined an isomorphism of exact sequences
$$\xymatrix@C-10pt{
\dots \ar[r] &L_{n+1}(\ZZ[\pi_1(X)]) \ar[r] \ar@{=}[d]
& \SS^{TOP}(X) \ar[r] \ar[d]^{\cong} & [X,G/TOP] \ar[d]^{\cong} \ar[r] &L_n(\ZZ[\pi_1(X)])\ar@{=}[d]\ar[r] &\dots\\
\dots \ar[r] & L_{n+1}(\ZZ[\pi_1(X)]) \ar[r]
& \SS_{n+1}(X) \ar[r]  & H_n(X;\LL_{\bullet}(\ZZ))  \ar[r]^-A &
L_n(\ZZ[\pi_1(X)]) \ar[r] & \dots}$$
The mapping cylinder of a homotopy equivalence $h:M \to X$ of $n$-dimensional topological manifolds is an $(n+1)$-dimensional geometric Poincar\'e cobordism $(W;M,X)$ with manifold boundary. The rel $\partial$ total surgery obstruction $s_{\partial}(W)\in \SS_{n+1}(X)$ defines the isomorphism
$$\SS^{TOP}(X) \to \SS_{n+1}(X)~;~(M,h) \mapsto s_{\partial}(W)~.$$
\indent The spectrum $\LL_{\bullet}(R)$ of Kan $\Delta$-sets was defined
in Ranicki \cite[p. 297]{ranicki0} using $n$-ads of quadratic Poincar\'e complexes over a ring with involution $R$ and the assembly map
$$A~:~X_+ \wedge \LL_{\bullet}(\ZZ) \to \LL_{\bullet}(\ZZ[\pi_1(X)])$$
was defined using bisimplicial methods, and the homotopy fibration sequence
of spectra
$$\xymatrix{
\LL_{\bullet}(\ZZ) \ar[r]^-{1+T} & \LL^{\bullet}(\ZZ) \ar[r]^-{J} & \widehat{\LL}^{\bullet}(\ZZ) \ar[r]^-{\partial} & \Sigma \LL_{\bullet}(\ZZ)}$$
with $\LL^{\bullet}(\ZZ)$ the symmetric $\LL$-spectrum such that
$$\pi_n(\LL^{\bullet}(\ZZ))~=~L^n(\ZZ)~=~
\begin{cases}
\ZZ&{\rm if}~n \equiv 0(\bmod\, 4)\\
\ZZ_2&{\rm if}~n \equiv 1(\bmod\, 4)\\
0&{\rm if}~n \equiv 2(\bmod\, 4)\\
0&{\rm if}~n \equiv 3(\bmod\, 4)
\end{cases}$$
and $\widehat{\LL}^{\bullet}(\ZZ)$ the hyperquadratic $\LL$-spectrum such that
$$\pi_n(\widehat{\LL}^{\bullet}(\ZZ))~=~\widehat{L}^n(\ZZ)~=~
\begin{cases}
\ZZ&{\rm if}~n=0\\
\ZZ_8&{\rm if}~n \equiv 0(\bmod\, 4)\\
\ZZ_2&{\rm if}~n \equiv 1(\bmod\, 4)\\
0&{\rm if}~n \equiv 2(\bmod\, 4)\\
\ZZ_2&{\rm if}~n \equiv 3(\bmod\, 4)
\end{cases}$$
A $(k-1)$-spherical fibration $\nu:X \to BG(k)$ has a canonical
$\widehat{\LL}^\bullet(\ZZ)$-cohomology Thom class $\widehat{U}(\nu) \in \dot H^k(T(\nu);\widehat{\LL}^{\bullet}(\ZZ))$ (which lifts to a $\LL^\bullet(\ZZ)$-cohomology Thom class $U(\nu) \in \dot H^k(T(\nu);\LL^{\bullet}(\ZZ))$ if and only if $\nu$ has a topological block bundle
reduction $\widetilde{\nu}:X \to B\widetilde{TOP}(k)$). For an $n$-dimensional geometric Poincar\'e complex $X$ with Spivak normal structure
$(\nu_X:X \to BG(k),\rho_X:S^{n+k} \to T(\nu_X))$ the $S$-dual
of $\widehat{U}(\nu_X)$ is a fundamental $\widehat{\LL}^{\bullet}(\ZZ)$-homology class $[X]_{\widehat{\LL}^{\bullet}(\ZZ)} \in
H_n(X;\widehat{\LL}^\bullet(\ZZ))$. The total surgery obstruction
$s(X) \in \SS_n(X)$ is the cobordism class of $\partial [X]_{\widehat{\LL}^{\bullet}(\ZZ)} \in H_{n-1}(X;\widehat{\LL}_\bullet(\ZZ))$
with assembly
$$A(\partial [X]_{\widehat{\LL}^{\bullet}(\ZZ)})~=~(C,\psi)
~=~0 \in L_{n-1}(\ZZ[\pi_1(X)])$$
with $(={\mathcal C}([X] \cap -:C(\widetilde{X})^{n-*} \to C(\widetilde{X}))_{*+1}$ a contractible
$\ZZ[\pi_1(X)]$-module chain complex. (Incidentally, the recent paper of K\"uhl, Macko and  Mole \cite{kmm} has given a detailed exposition of the total surgery obstruction from the point of view of \cite{ranicki0}.) A priori, it was not clear that $C$ supports an $(n-1)$-dimensional quadratic Poincar\'e
structure $\psi$, but this was remedied in Ranicki \cite[\S7.3]{ranicki3}: for any $n$-dimensional normal complex $(X,\nu_X:X \to BG(k),\rho_X:S^{n+k} \to T(\nu_X))$ the $\pi_1(X)$-equivariant $S$-dual of
$$\xymatrix{S^{n+k} \ar[r]^-{\rho_X} &T(\nu_X) ~=~T(\nu_{\widetilde{X}})/\pi_1(X)
\ar[r]^-{\Delta} & \widetilde{X}^+ \wedge_{\pi_1(X)} T(\nu_{\widetilde{X}})}$$
is a semistable $\pi_1(X)$-equivariant map
$F:T(\nu_{\widetilde{X}})^* \to \Sigma^k \widetilde{X}^+$
inducing the $\ZZ[\pi_1(X)]$-module chain map
$$f~=~ [X] \cap - ~:~\dot C(T(\nu_{\widetilde{X}})^*)_{*+k}~=~C(\widetilde{X})^{n-*}  \to C(\widetilde{X})~.$$
The $\pi_1(X)$-equivariant spectral quadratic construction
$$\psi_F~:~\dot H^k(T(\nu_X))~=~\ZZ \to Q_n({\mathcal C}(f))$$
factors through the suspension $S:Q_{n-1}({\mathcal C}(f)_{*+1}) \to Q_n({\mathcal C}(f))$ (\cite[Prop. 7.4.1. (iv)]{ranicki3}).
The construction $\psi_F$ is induced by the $\pi_1(X)$-equivariant
spectral Hopf invariant $\psi_F$ (\S\ref{spectral}, Definition \ref{pisym5}).

\begin{remark}
See Ranicki \cite[\S17]{ranicki4}, \cite[\S9]{ranicki5} for a combinatorial approach to the total surgery obstruction $s(X) \in \SS_n(X)$ of an $n$-dimensional geometric Poincar\'e simplicial complex $X$ using the $X$-local category of $(\ZZ,X)$-modules of Ranicki and Weiss \cite{ranickiweiss}.
\hfill\qed
\end{remark}

\renewcommand{\thechapter}{\Alph{chapter}}%

\setcounter{chapter}{0}

\appendix

\chapter{The homotopy Umkehr map}\label{appendix1}

Given a normal map $(f,b):M \to X$ we construct,
using the
method of Crabb and James \cite[p.261]{crabbjames} and Crabb
and Ranicki \cite{crabbranicki},
a fibrewise stable map
$$
F:X\sqcup_XX \sto \MM \sqcup_X X
$$
over $X$,
where $\MM$ is the path space of pairs
$(x,\alpha )$, $x\in M$, $\alpha : [0,1]\to X$, with
$\alpha (0)=f(x)$, fibred over $X$ by $(x,\alpha )\mapsto \alpha (1)$.

\section{Fibrewise homotopy theory}

Let $B$ be a reasonable space, such as a finite $CW$ complex.  A
{\it fibrewise space} over $B$ is a space $X$ together with a map
$p_X:X \to B$.  A {\it fibrewise map} $f:X \to Y$ is a map such
that the diagram\index{fibrewise!space}\index{fibrewise!map}
    $$\xymatrix{X \ar[rr]^-{\di{f}} \ar[dr]_-{\di{p_X}} && Y
    \ar[dl]^-{\di{p_Y}} \\ &B& }$$
commutes. For fibrewise pointed spaces $X,Y$ let $[X,Y]_B$ denote
the set of (fibrewise) homotopy classes of fibrewise maps
$f:X \to Y$.

The {\it product} of fibrewise spaces $X,Y$ is the fibrewise space
\index{fibrewise!product, $X\times_BY$}
$$\begin{array}{ll}
X\times_BY&~=~\{(x,y) \in X \times Y\,\vert\, p_X(x)=p_Y(y) \in B\}\\[1ex]
&=~\bigcup\limits_{b \in B} (p_X)^{-1}(b) \times (p_Y)^{-1}(b)\\[1ex]
&=~(p_X \times p_Y)^{-1}(\Delta_B) \subset X\times Y
\end{array}$$
with $\Delta_B \subset B \times B$ the diagonal subspace.

Let $q_B:\widetilde{B} \to B$ be the projection of
a regular cover with group of covering translations $\pi$.
For a fibrewise space $X$ over $B$ the pullback cover of $X$
$$\begin{array}{ll}
\widetilde{X}&=~(p_X)^*(\widetilde{B})~=~X \times_B \widetilde{B} \\[1ex]
&=~\{(x,\widetilde{b}) \in X \times \widetilde{B}\,\vert\,p_X(x)=q_B(\widetilde{b}) \in B\}
\end{array}$$
has projection
$$\widetilde{X} \to X~;~
(x,\widetilde{b}) \mapsto x~(p_X(x)=q_B(\widetilde{b})\in B)$$
with a $\pi$-equivariant lift of $p_X$
$$\widetilde{p}_X~:~\widetilde{X} \to \widetilde{B}~;~(x,\widetilde{b}) \mapsto
\widetilde{b}~.$$
For a fibrewise product $X \times_BY$ the pullback cover of $X\times_BY \to B$
$$\begin{array}{ll}
\widetilde{X\times_BY}&=~\widetilde{X}\times_{\widetilde{B}}\widetilde{Y}\\[1ex]
&=~\{(x,y,\widetilde{b}) \in X \times Y \times \widetilde{B}\,\vert\,
p_X(x)=p_Y(y)=q_B(\widetilde{b}) \in B\}
\end{array}$$
can be regarded as the $\pi$-equivariant subspace
$\widetilde{X}\times_{\widetilde{B}}\widetilde{Y} \subset  \widetilde{X} \times \widetilde{Y}$
consisting of all the $((x,\widetilde{b}),(y,\widetilde{c})) \in \widetilde{X} \times
\widetilde{Y}$ with $\widetilde{b}=\widetilde{c} \in\widetilde{B}$.
 Passing to the quotients by the $\pi$-actions we obtain the injective
{\it assembly map}\index{assembly map, $A$}
$$A~:~X \times_B Y \to \widetilde{X} \times_{\pi} \widetilde{Y}~;~(x,y)
\mapsto [(x,\widetilde{b}),(y,\widetilde{b})]$$
with $\widetilde{b} \in q_B^{-1}(p_X(x))=q_B^{-1}(p_Y(y)) \subseteq \widetilde{B}$.
In particular, for $X=Y=B$ the assembly is just the diagonal map
$$A~=~\Delta~:~B \times_B B~=~B \to \widetilde{B} \times_{\pi} \widetilde{B}~;~
b \mapsto [\widetilde{b},\widetilde{b}]~.$$
\indent The {\it coproduct} of fibrewise spaces $X,Y$ is the disjoint union\index{fibrewise!product, $X\sqcup_BY$}
    $$X\sqcup_BY=~X \sqcup Y$$
with the topology characterized by the property that fibrewise
maps $X \sqcup Y \to Z$ are in one-one correspondence with the
fibrewise maps $X \to Z$, $Y \to Z$. The {\it quotient} fibrewise
space of a closed subspace $A \subseteq X$ is
$$X/_BA~=~(X \sqcup_B B)/\sim$$
with $\sim$ the equivalence relation generated by $x \sim b$ for $x \in A$, $b
\in B$ with $p_X(x)=b$.

A fibrewise space $X$ is {\it pointed} if the projection $p_X:X
\to B$ is equipped with a section $s_X:B \to X$. The {\it smash
product} of pointed fibrewise spaces $X,Y$ is the pointed
fibrewise space\index{fibrewise!pointed space}\index{fibrewise!smash product, $X\wedge_BY$}
    $$X\wedge_BY~=~X \times_B Y/_B((X \times_B s_Y(B))
                   \cup (s_X(B) \times_B Y))~.$$
In particular, for any spaces $P,Q$ over $B$
$$(P \sqcup_B B) \wedge_B (Q \sqcup_BB)~=~(P \times_B Q) \sqcup_B B~.$$
A fibrewise map $f:X \to Y$ over $B$ induces a $\pi$-equivariant map
$$\widetilde{f}~:~\widetilde{X} \to \widetilde{Y}~;~(x,\widetilde{b}) \to (f(x)),
\widetilde{b})~.$$
Given a covering projection $q_B:\widetilde{B} \to B$ with group of covering
projections $\pi$ and a pointed fibrewise space $X$ over $B$ define the pointed
$\pi$-space  $\widetilde{X}/\widetilde{B}$, with the action of $\pi$ free away
from the base point. Define the {\it assembly map} from the fibrewise homotopy set to the $\pi$-equivariant homotopy set
$$A~:~[X,Y]_B \to [\widetilde{X}/\widetilde{B},\widetilde{Y}/\widetilde{B}]^{\pi}~;~f \mapsto \widetilde{f}~.$$
\indent
The {\it suspension} of $X$ is the pointed fibrewise space\index{fibrewise!suspension, $\Sigma_BX$}
    $$\Sigma_BX~=~(B \times S^1) \wedge_B X~.$$
\indent The (disk,\,sphere)-bundle of a (finite-dimensional real)
vector bundle $\xi$ over a space $M$ over $B$ is a pair of spaces
$(D(\xi),S(\xi))$ over $B$. The {\it fibrewise Thom space} of $\xi$ is the
fibrewise space over $B$ is\index{fibrewise!Thom space, $T_B(\xi)$}
$$T_B(\xi)~=~D(\xi)/_BS(\xi)~.$$
(The fibrewise Thom space is denoted by $M^{\xi}$ in \cite{crabbjames}
and \cite{crabbranicki}).
For each $x \in B$ there is a base point $*_x \in T_B(\xi)$, so that
$B \subset T_B(\xi)$, with
$$T_B(\xi)/B~=~D(\xi)/S(\xi)~=~T(\xi)$$
the ordinary Thom space of $\xi$.
The effect on the fibrewise Thom space of adding the trivial line
bundle $\epsilon$ to $\xi$ is given by
    $$T_B(\xi \oplus \epsilon)~=~\Sigma_BT_B(\xi)~.$$
In particular $$T_B(0)~=~M\sqcup_BB~,~T_B(\epsilon^k)~=~\Sigma^k_B(M\sqcup_BB)~~(k \geqslant 0)~.$$
If $q_B:\widetilde{B} \to B$ is a covering projection with group of covering translations $\pi$ the pullback $\widetilde{\xi}=q_B^*(\xi)$ is  a vector bundle over $\widetilde{B}$ such that
$$T_B(\widetilde{\xi})~=~T(\widetilde{\xi})/\widetilde{B}$$
is a semifree $\pi$-space with quotient
$$T_B(\widetilde{\xi})/\pi~=~T_B(\xi)~.$$
\indent A {\it fibrewise stable map} $F:X \sto Y$ of pointed
fibrewise spaces is an equivalence class of fibrewise maps
$$F~:~\Sigma^k_BX  \to \Sigma^k_BY~~(k \geqslant 0)~.$$
The abelian group of fibrewise stable maps $X \sto Y$ is written
$$\omega_B^0\{X;Y\}~=~\lim\limits_{k\to\infty} [\Sigma^k_BX;\Sigma^k_BY]_B~.$$

The fibrewise geometric Hopf map can be constructed by a routine
extension of the method in section \ref{hopf}  as a function
$$
\omega_B^0\{X;Y\} \to \omega_B^0\{X;S(\infty)^+\wedge_{\ZZ_2}(Y \wedge_BY)\}~;~
F \mapsto h(F)~$$
(See \cite[pp. 168--169, 306--308]{crabbjames} for details.)

\section{The homotopy Pontryagin-Thom construction}\label{A-emb}

Suppose given an embedding $f:M \emb N\backslash \partial N$ of a closed
$m$-dimensional manifold $M$ in the interior of an $n$-dimensional manifold
$N$ with boundary $\partial N$ (which may be empty).

We write the path space fibration of the embedding $f$ as
$$
f'~:~\MM~=~\{(x,\alpha)\in M \times N^{[0,1]} \,\vert\,
\alpha(0)=f(x)\}~ \to N~;~(x,\alpha) \mapsto \alpha(1)~.
$$
The projection
$$
\pi~:~\MM \to M~;~(x,\alpha) \mapsto x
$$
is a homotopy equivalence such that $f\pi \simeq f'$.

Let $\nu_f:M \to BO(n-m)$ be the normal bundle, and let
$D(\nu_f) \emb N$ be a tubular neighbourhood. The standard
Pontrjagin-Thom construction gives a map
$$
F~:~N/\partial N \to T(\nu_f)~;~y \mapsto
\begin{cases} y&{\rm if}~y \in D(\nu_f)\\
*&{\rm if}~y \in N \backslash D(\nu_f)\end{cases}
$$
which we shall now extend to a {\it  homotopy Pontryagin-Thom map}over $N$
$$
F_N~:~N\cup_{\partial N}N \to T_N(\pi^*\nu_f)
$$
constructed as follows. For $y\in D(\nu_f)$ in the fibre of
$D(\nu_f)$ over $x\in M$, let $\alpha : [0,1]\to N$ be the linear
path from $f(x)$ to $y$ in the fibre $D(\nu_f(x))$ (in $N$), and
let $z \in T_N(\pi^*\nu_f)$ be the point over $y$ given by $y$ in
the fibre of $\pi^*\nu_f$ over $(x,\alpha )\in \MM$.  Define a map
$$N \to T_N(\pi^*\nu_f)~;~y \mapsto
\begin{cases} z&{\rm if}~y \in D(\nu_f)\\
*&{\rm if}~y \in N \backslash D(\nu_f)~.\end{cases}$$ Note that
$\alpha (1)=y$, so that this is fibrewise over $N$. Extend this
map on $N$ to $N\cup_{\partial N}N$ by mapping the other
(basepoint) copy of $N$ to the basepoint in each fibre of
$T_N(\pi^*\nu_f)$, to obtain the homotopy Pontryagin-Thom map
$F_N:N\cup_{\partial N}N\to T_N(\pi^*\nu_f)$.

Passing to the quotient by $N$
recovers the classical Pontryagin-Thom map
$$
F~:~N/\partial N~=~(N\cup_{\partial N}N)/N \to
T_N(\pi^*\nu_f)/N~=~T(\pi^*\nu_f)~\simeq~T(\nu_f)~.
$$

Given a regular cover $q_N:\widetilde{N}\to N$ with group of
covering translations $\Gamma$ (e.g. the universal cover for
connected $N$, with $\Gamma =\pi_1(N)$) it is possible to lift $f$
to a $\Gamma$-equivariant embedding
$$
\widetilde{f}~:~\widetilde{M}~=~f^*\widetilde{N}~=~\{(x,y) \in M
\times \widetilde{N}\vert
f(x)=q_N(y) \in N\} \emb \widetilde{N}~;~(x,y) \mapsto y
$$
and the
standard Pontryagin-Thom construction gives a $\Gamma$-equivariant
Umkehr map
$$\widetilde{F}~:~\widetilde{N}/\widetilde{\partial N} \to T(\nu_{\widetilde{f}})~=~T(q_M^*\nu_f)$$
with $q_M$ the projection
$$q_M~:~\widetilde{M} \to M~;~(x,y) \mapsto x~.$$
The assembly of the homotopy Pontryagin-Thom map $F_N$ determines this
$\Gamma$-equivariant Umkehr map $\widetilde{F}$
$$
A~:~[N \cup_{\partial N}N,T_N(\pi^*\nu_f)]_N \to
[\widetilde{N}/\widetilde{\partial N},T(q_M^*\nu_f)]_{\Gamma}~;~
F_N \mapsto \widetilde{F}~.
$$
More precisely, lift the fibrewise map $F_N:N\cup_{\partial N}N
\to T_N(\pi^*\nu_f)$ from $N$ to $\widetilde{N}$ to get a
$\Gamma$-equivariant map over $\widetilde{N}$
$$\widetilde{F}~:\widetilde{N}\cup_{\widetilde{\partial N}}\widetilde{N}
\to T_{\widetilde{N}}(\widetilde{\pi}^*\nu_f)~,$$
where $\widetilde{\pi}$ is obtained by lifting $\pi$~:
$$\begin{array}{l}
\widetilde{\pi}~ :~ \widetilde{\MM}~=~ \{ (x,\alpha ,y) \in M \times
N^{[0,1]} \times \widetilde{N}
\mid\alpha (0)=f(x),\alpha (1)=q_N(y)\}\\[1ex]
\hskip150pt \to M~;~(x,\alpha,y) \mapsto x~.
\end{array}$$
For $(x,\alpha ,y)\in \widetilde{\MM}$, we have a unique path
$\widetilde{\alpha} : [0,1]\to \widetilde{N}$ lifting $\alpha$
with $\widetilde{\alpha}(1)=y$.  Then $\widetilde{\alpha}(0)\in
\widetilde{N}$ projects to $f(x)\in N$.  Hence
$(x,\widetilde{\alpha}(0))\in \widetilde{M}$.  This allows us to
define a map
$$r~:~\widetilde{\MM} \to \widetilde{M}~;~(x,\alpha ,y)\mapsto
(x,\widetilde{\alpha} (0))$$ such that $\widetilde{\pi}= q_M\circ
r$. We obtain the $\Gamma$-equivariant map $\widetilde{g}$ by
collapsing the basepoints $\widetilde{N}$
$$\widetilde{g}~
:~\widetilde{N}/\widetilde{\partial N}~=~(\widetilde{N}\cup_{\widetilde{\partial N}}\widetilde{N}) /(\widetilde{N} \times \{*\})
 \to T_{\widetilde{N}}(\widetilde{\pi}^*\nu_f)/(\widetilde{N} \times \{*\}) ~=~T(q_M^*\nu_f)~.
$$

In particular, if $g=F_N:N\cup_{\partial N}N \to T_N(\pi^*\nu_f)$ then
$\widetilde{g}=\widetilde{F}:\widetilde{N}/\widetilde{\partial N}
\to T(q^*_M\nu_f)$.
If $f:M \emb N$ is a framed embedding, with $\nu_f=\epsilon^{n-m}$
the trivial $(n-m)$-plane bundle, the homotopy Pontryagin-Thom
construction is a map
$$
F_N~:~N \cup_{\partial N}N\to T_N(\pi^*\nu_f)~=~\Sigma_N^{n-m}(\MM \sqcup_N N)$$
and the $\Gamma$-equivariant Umkehr is a map
$$
\widetilde{F}~:~\widetilde{N}/\widetilde{\partial N} \to T(q^*_M\nu_f)~=~\Sigma^{n-m}(\widetilde{M}^+)~.
$$

\section{The homotopy Umkehr map of an immersion}

More generally, suppose given an immersion $f:M \imm N\backslash\partial N$ of
a closed $m$-dimensional manifold in the interior of an $n$-dimensional
manifold with boundary, with normal bundle $\nu_f:M \to BO(n-m)$.

We again write the path space fibration of the map $f$ as the space
$$
f'~:~\MM~=~\{(x,\alpha)\in M \times N^{[0,1]} \,\vert\,
\alpha(0)=f(x)\}~ \to N~;~(x,\alpha) \mapsto \alpha(1)
$$
over $N$, with the homotopy equivalence $\pi :\MM \to M$.

In this section we shall associate to $f$ a fibrewise stable homotopy
class
$$
F_N \in \omega^0_N\{ N\cup_{\partial N} N;\, T_N(\pi^*\nu_f)\}
$$
that we call the {\it homotopy Umkehr of the immersion} $f$.

For
sufficiently large $j \geqslant 0$ there exists a map $e:M \to
\R^j$ such that
$$g~=~e \times f~:~M \to \R^j \times N~;~x \mapsto (e(x),f(x))$$
is an embedding which is regular homotopic to the composite
    $$M \imm N \emb \R^j\times N$$
of $f$ and the embedding
    $$N \emb \R^j \times N~;~x \mapsto (0,x)$$
with trivial normal bundle $\epsilon^j:N \to BO(j)$. The normal
bundle of $g$ is
$$\nu_g~=~\nu_f\oplus \epsilon^j~:~M \to BO(n-m+j)~.$$

The homotopy Pontryagin-Thom map of $g$  is a map
$$
G_{\R^j\times N}~:~\R^j \times (N \cup_{\partial N}N) \to
T_{\R^j \times N}((\pi')^*(\nu_{e \times f}))$$ with $\pi':\MM' \to M$
over $\R^j\times N$ with
the projection of
$$
\MM'~=~\{(x,\alpha) \in M \times (\R^j \times N)^{[0,1]}\,\vert\,
\alpha(0)=(e(x),f(x))\}~.
$$
It is possible to regard $G_{\R^j\times N}$ as a fibrewise stable map over $N$
$$
F_N~=~G_{\R^j \times N}~:~ \Sigma^j_N(N\cup_{\partial N}N) \to
\Sigma^j_NT_N(\pi^*\nu_f)~,
$$
representing an element
$$
F_N \in \omega_N^0\{N \cup_{\partial N}N;T_N(\pi^*\nu_f)\}~.
$$
The stable homotopy class of the homotopy Umkehr map $F_N$ depends only
on the regular homotopy class of the immersion $f$.

For a regular cover $q_N:\widetilde{N} \to N$ with group of
covering translations $\Gamma$, $g$ lifts to a
$\Gamma$-equivariant embedding $\widetilde{g}:\widetilde{M} \emb
\R^j\times \widetilde{N}$, and there is a $\Gamma$-equivariant
Umkehr map
$$
\widetilde{F}~:~\Sigma^j(\widetilde{N}/\widetilde{\partial N}) \to
T(\nu_{\widetilde{g}})
~=~\Sigma^jT(\nu_{\widetilde{f}})~=~ \Sigma^jT(q^*_M\nu_f)~.
$$
It is determined by $F_N$:
$$
\omega_N^0\{N \cup_{\partial N}N;T_N(\pi^*\nu_f)\} \to
\{\widetilde{N}^+;T(\nu_{\widetilde{f}})\}^{\Gamma}~;~F_N \mapsto
\widetilde{F}~.
$$

\section{The homotopy Umkehr map of a normal map}

Suppose given a Browder-Novikov-Sullivan-Wall normal map $(f,b):M \to X$ from an
$m$-dimensional manifold $M$ to an $m$-dimensional geometric
Poincar\'e complex $X$, with $b:\nu_M \to \eta$ a map over $f$ from a stable normal bundle
$\nu_M$ of $M$ to a vector bundle $\eta$ over $X$.

In this section, we write $\MM \to X$ for the path
fibration of the map $f$
$$
f'~:~\MM~=~\{(x,\alpha)\in M \times N^{[0,1]} \,\vert\,
\alpha(0)=f(x)\}~ \to X~;~(x,\alpha) \mapsto \alpha(1),
$$
with the homotopy equivalence $\pi : \MM \to M$, $(x,\alpha )\mapsto x$.
We shall construct a fibrewise stable homotopy class over $X$,
$$
F_X : \omega^0_X\{ X \sqcup_X X ;\, \MM\sqcup_X X\},
$$
that we call the {\it homotopy Umkehr of the normal map}.

As in \S\ref{surgeryobstruction}, by the $\pi$-$\pi$ theorem of Wall \cite[Chapter 3]{wall2},
for sufficiently large $j \geqslant 0$ $(X
\times D^j,X \times S^{j-1})$ is homotopy equivalent to an
$(m+j)$-dimensional manifold with boundary $(N,\partial N)$ and
$(f,b)$ can be approximated by a framed embedding $M \emb N
\backslash \partial N$. The homotopy Pontryagin-Thom map of the embedding
$$
N\cup_{\partial N}N \to \Sigma^j_N(\MM \sqcup_N N)
$$
determines the homotopy Umkehr map
$$
F_X~:~ \Sigma^j_X(X \sqcup_X X) \to \Sigma^j_X(\MM \sqcup_X X)
$$
representing an element $F_X \in \omega_X^0\{X \sqcup_X X; \MM \sqcup_XX\}$.
We shall show below that this class is independent of the choice of the
manifold $N$.

Now the fibrewise Hopf invariant is a function
$$
\begin{array}{l}
h(F_X)\in \omega_X^0(X \sqcup_X X;\,
S(\infty)^+\wedge_{\ZZ_2}((\MM \sqcup_XX)\wedge_X(M \sqcup_XX))\}\\[1ex]
\hskip75pt =~
\omega_X^0(X\sqcup_XX;\,
(S(\infty)\times_{\ZZ_2}(\MM \times_X \MM ))\sqcup_XX)\}~.
\end{array}
$$

Notice that the fibre product $\MM \times_X \MM$ can be described as
the space of triples $(x,y,\alpha )$, where $x,\, y\in M$ and $\alpha
: [-1,1]\to X$ is a path from $x=\alpha (-1)$ to $y=\alpha (1)$,
projecting to $\alpha (0)\in X$.

If $\widetilde{X}$ is the universal cover of $X$ and
$\widetilde{M}=f^*\widetilde{X}$ is the pullback cover of $M$ then
$F_X$ induces the $\pi_1(X)$-equivariant $S$-dual
$\widetilde{F}:\Sigma^j(\widetilde{X}^+) \to \Sigma^j(\widetilde{M}^+)$ of
Ranicki \cite{ranicki2}, and the fibrewise Hopf invariant induces the $\pi_1(X)$-equivariant Hopf invariant (inducing the quadratic construction on the chain level)
$$
h(\widetilde{F}) \in \{X^+,(S(\infty) \times_{\ZZ_2}
(\widetilde{M} \times_{\pi_1(X)}\widetilde{M}))^+\}~.
$$

\indent Finally, here is the proof that $F_X$ is well-defined, i.e. independent
of the choice of $N$.
It is convenient to use the notation of Crabb and Ranicki \cite{crabbranicki},
which is reprinted in Appendix \ref{appendix3} and to which we now refer.
In particular, we need to consider manifolds without boundary that are not
necessarily compact and use fibrewise stable homotopy with compact supports.

Let $M$ be a closed manifold and $N$ a manifold with empty boundary.
Consider a map $f: M\to N$. The classical Umkehr is a stable map
$$
f^!:N^+ \to M^{f^*\tau N -\tau M},
$$
from the one-point compactification of $N$.
It is constructed by choosing a smooth embedding $e: M\into V$ for
some Euclidean space $V$ and deforming $f$ to be a smooth map to get
a smooth embedding $(e,f): M\into V\times N$ with normal bundle $\nu$.
The Pontrjagin-Thom construction gives a map
$$
(V\times N)^+ =V^+\wedge N^+ \to M^\nu,
$$
which represents $f^!$.

The homotopy Umkehr is a fibrewise stable map over $N$ with compact supports
in
$$
{}_c\omega^0_N\{ N\times S^0;\, \CC_N^{\pi^*(f^*\tau N-\tau M)}\},
$$
where $\CC = \{ (x,y,\alpha ) \st  x\in M,\, y\in N,\, \alpha : [0,1]
\to N,\, \alpha (0)=f(x),\, \alpha (1)=y\}$ projecting $(x,y,\alpha)$
to $y\in N$,
and $\pi :\CC \to M$ maps $(x,y,\alpha )$ to $x$.

There is a duality isomorphism
$$
{}_c\omega^0_N\{ N\times S^0;\, \CC_N^{\pi^*(f^*\tau N-\tau M)}\}
\to
\omega^0\{ S^0;\, \CC^{\pi^*(f^*\tau N-\tau M)-\tau N}\}\,.
$$
(See  Crabb \cite[Prop 4.1]{crabbcoincidence}.)
Since $\pi :\CC \to M$ is a homotopy equivalence, this group is just
$$
\omega^0\{ S^0;\, M^{-\tau M}\} =\omega^0\{ M_+;\, S^0\}= \omega^0(M).
$$
The fibrewise Umkehr corresponds to the element $1\in\omega^0(M)$.

We can see this by rewriting the definition. Assume that $M$ is embedded
as a submanifold of $N$.
First we look at the inclusion of the open subspace  $i:E\nu \into N$
(the total space of $\nu$).
By construction, the homotopy Umkehr lies in the image of
$$
i_!:
{}_c\omega^0_{E\nu}\{ E\nu \times S^0;\, \CC_{E\nu}^{\pi^*(f^*\tau N-\tau M)}\}
\to
{}_c\omega^0_N\{ N\times S^0;\, \CC_N^{\pi^*(f^*\tau N-\tau M)}\}
$$
given by null extension over the complement.
This allows us to reduce to the case in which $N=E\nu$, which we now assume.

We have a fibrewise map $E\nu \to \CC_{E\nu}$ over $E\nu$, given by
straight lines:
$$
y\in \nu_x  \mapsto (x,y, \alpha ),
$$
where $\alpha (t)=ty\in\nu_x$.

The homotopy Umkehr is the image under
$$
{}_c\omega^0_{E\nu}\{ E\nu\times S^0;\,  (E\nu)_{E\nu}^{p^*\nu}\}
\to
{}_c\omega^0_{E\nu}\{ E\nu\times S^0;\, \CC_{E\nu}^{p^*\nu}\} ,
$$
where $p:E\nu \to M$ is the projection, of the class in
$$
{}_c\omega^0_{E\nu}\{ E\nu\times S^0;\,  (E\nu)_{E\nu}^{p^*\nu}\}
=
{}_c\omega^0_{E\nu} \{ E\nu \times S^0;\, D(p^*\nu )/_{E\nu} S(p^*\nu)\}
$$
given by the projection $E\nu_x \to D(\nu_x)/S(\nu_x)$ ($x\in M$).
This gives the suspension of $1\in \omega^0(M)$.

Now suppose that $h: N \to N'$ is a proper map and that
$a$ is a stable vector bundle isomorphism from $\tau N$ to $h^*\tau N'$.
We want to relate the homotopy Umkehr of $f$ and the homotopy Umkehr of
$f'=h\circ f : M\to N'$.
There are homomorphisms
$$
{}_c\omega^0_N\{ N\times S^0;\, \CC_N^{\pi^*(f^*\tau N-\tau M)}\}\to
{}_c\omega^0_N\{ N\times S^0;\, \DD_N^{\pi^*(f^*\tau N-\tau M)}\},
$$
and
$$
(h,a)^*:
{}_c\omega^0_{N'}\{ N'\times S^0;\,
(\CC')_{N'}^{\pi^*((f')^*\tau N'-\tau M)}\}\to
{}_c\omega^0_N\{ N\times S^0;\, \DD_N^{\pi^*(f^*\tau N-\tau M)}\},
$$
where $\DD$ is the pullback $h^*\CC'$, that is,
$\DD =\{ (x,y,\alpha' ) \st x\in M,\, y\in N,\, \alpha':[0,1]\to N',\,
\alpha'(0)=h(f(x)),\, \alpha'(1)=h(y)\}$,
and $\CC \to\DD$ maps $(x,y,\alpha )$ to $(x,y,h\circ \alpha )$.

We claim that the homotopy Umkehr of $f$ and $f'$ have the same image in
$$
{}_c\omega^0_N\{ N\times S^0;\, \DD_N^{\pi^*(f^*\tau N-\tau M)}\}.
$$

Choose a smooth embedding $e: M\into V$, with $V\not=0$. Deform $f$ to
a smooth map $M\to N$.
Then deform $h$ by a homotopy constant outside a compact subspace of $N$
to a map that is smooth in an open neighbourhood of $f(M)$.
Since $V$ is non-zero, $a$ gives a vector bundle isomorphism
$V\oplus \tau N \to V\oplus h^*\tau N'$, which we can arrange by deformation
to be smooth in a neighbourhood of $(e\times f)(M)\subseteq V\times N$.

Let us now change notation, writing $N$ instead of $V\times N$, $f$ instead
of $e\times f$ and
$N'$ instead of $V\times N'$. So we have  a smooth embedding
$f: M\into N$, $h:N\to N'$ is smooth on a neighbourhood of $f(M)$,
and $f'=h\circ f: M\into N'$ is a smooth embedding. And we have a vector
bundle isomorphism $a : \tau N\to h^*\tau N'$ which is smooth in
a neighbourhood of $f(M)$.

Choose a riemannian metric on $N'$ and pull it back using $a$ to get
a riemannian metric on a neighbourhood of $f(M)$ in $N$.
The restriction of $a$ to $f(M)$
also sets up an isomorphism between the normal bundle $\nu$ of
$f$ and the normal bundle $\nu'$ of $f'$.

We can take small tubular neighbourhoods $D(\nu )\into N$ and $D(\nu')\into
N'$ of $M$.
Then we can deform $h$ through a homotopy that is constant outside $D(\nu )$
to a map that restricts on the smaller disc bundle
$D_{1/2}(\nu )$ to the diffeomorphism
$D_{1/2}(\nu ) \to D_{1/2}(\nu')$ given by
$a| : \nu \to \nu'$.
This has arranged that $h$ is a diffeomorphism from a neighbourhood
of $f(M)$ to a neighbourhood of $f'(M)$.

Now we can just compare the constructions, proving that $F_X$ is indeed independent of $N$.
For a normal map $(f,b): M\to X$ we thus get a well-defined
homotopy Umkehr map
$$
F_X\in \omega^0_X\{ X\times S^0;\, \MM_{+X}\}.
$$

\section{A fibrewise spectral Hopf invariant}
We continue to use the notation from Appendix \ref{appendix3}.

Suppose now that a closed manifold $M$ is embedded as a submanifold of
a closed manifold $N$ with normal bundle $\nu$.
Choose again a tubular neighbourhood $D(\nu )\into M$ of $M$
and write $P=N- B(\nu )$ for the complement of the open disc bundle $B(\nu )$.
Thus $P$ is a compact manifold with boundary $\partial P= S(\nu )$.

As in section \ref{A-emb}, we write $\MM \to N$ for the space of
continuous paths $\alpha : [0,1]\to N$ such that $\alpha (0)\in M$, fibred
over $N$ by projection to the end-point $\alpha (1)$, and $\pi :\MM
\to M$ for the homotopy equivalence taking $\alpha$ to $\alpha (0)$.
Let $\PP \to N$ be the space of paths $\alpha : [0,1]\to N$ such
that $\alpha (0)\in P$, mapping to $\alpha (1)\in N$.
We write $\pi : \PP\to P$ also for the homotopy equivalence
$\alpha\mapsto \alpha (0)$.

There is the following explicit description of the fibrewise
homotopy cofibre of the homotopy Pontryagin-Thom map
$$
N\times S^0 \to \PP_N^{\pi^*\nu}
$$
(due, essentially, to Klein and Williams \cite{kleinwilliams1}).

\begin{proposition}
The homotopy Pontryagin-Thom map fits into a fibrewise homotopy
cofibration sequence
$$
\PP_{+N} \to N\times S^0 \to  \MM_N^{\pi^*\nu} \to
\Sigma_N (\PP_{+N})\to\cdots
$$
over $N$.
\qed
\end{proposition}
A proof can be found in \cite[Lemma 6.1]{crabbcoincidence}.

Now let us specialize to the case in which
the normal bundle $\nu$ is trivial, say $\nu =M\times V$,
so that the homotopy Pontryagin-Thom map has the form
$$
N\times S^0 \to (N\times V^+)\wedge_N \MM_{+N},
$$
to which we may apply fibrewise the construction of the spectral Hopf
invariant, Definition \ref{sym5}, to produce a fibrewise stable
$\Zz_2$-equivariant map
over $N$ in the equivalent forms
$$
N\times \Sigma S(LV)_+ \to
(\PP\times_N \PP)_N^{\Rr\oplus L}\, ,
$$
or
$$
N\times (LV)^+ \to
(N\times (\Rr\oplus L)^+ )\wedge_N
(S(LV) \times (\PP\times_N \PP))_{+N}\, ,
$$
or
$$
N\times S^0 \to
(S(LV) \times (\PP\times_N \PP))_{+N}^{\Rr\oplus L-LV} \, .
$$
\smallskip
\par\noindent
This determines a non-equivariant fibrewise stable map as an element of
$$
\omega^0_N\{ N\times S^0;\, (S(LV)\times_{\Zz_2} (\PP\times_N
\PP))_N^{\Rr\oplus H -HV}\},
$$
where $H$ is the Hopf line bundle $S(LV)\times_{\Zz_2} L$ over
the real projective space $S(LV)/\Zz_2$.
\begin{remark}
Suppose that $V=\Rr$. Then this group reduces to
$$
\omega_N^0\{ N\times S^0;\, \Sigma_N (\PP\times_N\PP )_{+N}\}
\, ,
$$
which maps, by collapsing fibrewise basepoints, to
$$
\omega^0\{ N_+;\, \Sigma (\PP\times_N\PP)_+\}\, .
$$

Composing with the map $\pi\times\pi :
\PP\times_N\PP\to P\times P$
(which is not, in general, a homotopy equivalence),
we obtain the stable map
$$
N_+ \to\Sigma (P\times P)_+
$$
that appears in Proposition \ref{sym10}.
\end{remark}

\newpage

\chapter{Notes on $\ZZ_2$-bordism}\label{appendix2}
\section{Introduction}
In this Appendix we give a brief account of $\ZZ_2$-equivariant bordism,
expressing geometrically defined equivariant
bordism groups as equivariant stable homotopy groups of suitably defined spaces.
As the material may be of independent interest, the presentation
is written so as to be self-contained.

Equivariant stable homotopy will be written as $\omega_*^{\ZZ_2}$ ($*\in\Zz$).
We shall use the local coefficient notation
$\omega^{\ZZ_2}_*(X;\, \xi )$ for the reduced equivariant stable homotopy
group of the Thom space $\Thom (\xi , X)$ of a real $\ZZ_2$-vector bundle $\xi$ over $X$ and,
more generally,
if $\eta$ is another $\ZZ_2$-vector bundle over the same base,
$\omega^{\ZZ_2}_*(X;\,\xi -\eta )$ for the stable homotopy of
the Thom space of the virtual vector bundle $\xi -\eta$.
Corresponding notation is used for non-equivariant stable homotopy groups.

The same symbol will often be used for a finite-dimensional
$\ZZ_2$-module $V$ and the trivial $G$-vector bundle $X\times V\to X$
over a $\ZZ_2$-space $X$.
We write $V^\cpt$ for the one-point compactification of $V$ with
basepoint at $\infty$ and $X\bpt$ for the pointed space obtained
by adjoining a disjoint basepoint to $X$.

The standard non-trivial $1$-dimensional representation
$\Rr$ of $\ZZ_2$ with the involution $-1$ is denoted by $L$.
\section{Framed bordism}
Let $X$ be a (metrisable) $\ZZ_2$-ANR (Absolute Neighbourhood Retract),
and let $\xi$ and $\eta$ be
real $\ZZ_2$-vector bundles over $X$.
We restrict the topology of $X$ for the sake of precision,
but this is not really necessary as we are always dealing with the system of
maps from a compact $\ZZ_2$-ENR (Euclidean Neighbourhood Retract) to $X$.
The fibre dimensions of $\xi$ and $\eta$ are allowed to vary over
the components of $X$.

Up to homotopy we may assume that any finite-dimensional real $\ZZ_2$-vector
bundle is equipped with an equivariant positive-definite inner product.
We write
$\O_X(\xi ,\eta )\to X$ for the Stiefel bundle with fibre at $x\in X$ the
Stiefel manifold $\O (\xi_x,\eta_x)$ of linear isometries $\xi_x\into\eta_x$.
Topologically, the space $\O_X (\xi ,\eta )$ is a $\ZZ_2$-ANR.

\smallskip

The notion of a {\it restricted} stable isomorphism between two vector
bundles plays a crucial r\^ole in the theory.
\begin{definition}
Let $\zeta_0$ and $\zeta_1$ be real $\ZZ_2$-vector bundles over a
compact $\ZZ_2$-ENR $Y$.  We define a {\it restricted stable isomorphism}
$\zeta_0\iso\zeta_1$ to be an equivalence class of $\ZZ_2$-vector bundle
isomorphisms $\zeta_0\oplus\Rr^i\iso\zeta_1\oplus\Rr^i$,
the equivalence being generated by homotopy and stabilization with respect
to $i$.
It will sometimes be convenient to think of a restricted stable isomorphism
as given by an isomorphism $\zeta_0\oplus V_0\iso  \zeta_1\oplus V_0$
for some finite-dimensional $\Rr$-vector space $V_0$ on which $\ZZ_2$ acts
trivially.
\end{definition}
We recall that a {\it stable isomorphism} from $\zeta_0$ to $\zeta_1$
is an equivalence class of vector bundle isomorphisms $\zeta_0\oplus V\iso
\zeta_1\oplus V$, where $V$ is any finite dimensional real $\ZZ_2$-module.
A restricted stable isomorphism thus determines a stable isomorphism.
The restricted stable automorphisms of a $\ZZ_2$-vector bundle
correspond to elements of
$\KO^{-1}(Y/\ZZ_2)$, while the
stable automorphisms correspond to elements of $\KO^{-1}_{\ZZ_2}(Y)$.
(The orbit space $Y/\ZZ_2$ is a compact ENR.)
In general, different restricted stable isomorphisms may give
the same stable isomorphism, but not if the action of $\ZZ_2$ on $Y$ is free.
\begin{lemma} \label{B-resfree}
Suppose that $Y$ is a compact free $\ZZ_2$-ENR. Then every stable
isomorphism $\zeta_0\iso\zeta_1$ arises from a unique restricted
stable isomorphism.
\end{lemma}
\begin{proof}
For an isomorphism
$\zeta_0\oplus V\iso \zeta_1\oplus V$, where $V$ is a $\ZZ_2$-module,
corresponds to
an isomorphism $(\zeta_0\oplus V)/\ZZ_2\iso (\zeta_1\oplus V)/\ZZ_2$
over $Y/\ZZ_2$ and $(Y\times V)/\ZZ_2$ can be embedded as a direct
summand in a trivial bundle.
\hfill\qed\end{proof}

In setting up the basic definitions it is useful to allow a $\ZZ_2$-manifold
$M$, with tangent bundle $\tau M$, to have components of (possibly)
different dimensions. This will allow us to refer to {\it the} fixed
submanifold $M^{\ZZ_2}$
without introducing supplementary notation for the components
of each dimension.

We recall first the classical notion of non-equivariant framed bordism.
\begin{definition}
For non-equivariant vector bundles $\xi_0$ and $\eta_0$ over an ANR $X_0$,
we write $\Omega_0(X_0;\,\xi_0 ,\eta_0 )$ for the bordism group
of triples $(M_0,f,a)$ where $M_0$ is a closed manifold, $f : M_0\to X_0$ is a map
and $a: \tau M_0\oplus f^*\xi_0 \iso f^*\eta_0$ is a stable isomorphism.
\end{definition}
The non-equivariant bordism group depends only on the virtual bundle $\xi_0-\eta_0$ and the Pontrjagin-Thom construction gives an isomorphism
$$
\Omega_0(X_0;\, \xi_0 ,\eta_0 ) \Rarr{\iso}{} \omega_0(X_0;\, \xi_0 -\eta_0 ).
$$

The equivariant theory is more subtle.
\begin{definition}
We define the {\it restricted framed bordism group}
$\Omega^{\ZZ_2 ,\res}_0(X;\,\xi ,\eta )$ to be the bordism group
of triples $(M,f,a)$ consisting of
a closed $\ZZ_2$-manifold $M$ together with a $\ZZ_2$-map $f:M\to X$ and
a restricted stable isomorphism $a: \tau M\oplus f^*\xi \iso f^*\eta$,
specified by a vector bundle isomorphism
$$
\tau M\oplus f^*\xi\oplus V_0 \iso f^*\eta\oplus V_0
$$
for some vector space $V_0$ with trivial $\ZZ_2$-action.

Bordisms are required to be restricted, too, that is,
$(M,f,a)$ is null-bordant if there exists a $\ZZ_2$-manifold $W$ with boundary
$\partial W=M$ equipped with
a map $g:W\to X$ extending $f$ and a restricted stable isomorphism
$b:\tau W\oplus g^*\xi \iso \Rr\oplus g^*\eta$ extending $a$.

The negative of the class $(M,f,a)$ is given by $(M,f,a')$ where
$a' : \tau M\oplus f^*\xi\oplus V_0\oplus\Rr \iso f^*\eta\oplus V_0\oplus\Rr$
is $a\oplus (-1)$.
\end{definition}

Here we do not use the notation `$\xi -\eta$', because the restricted bordism
group $\Omega^{\ZZ_2 ,\res}_0(X;\,\xi ,\eta )$
(unlike the stable homotopy group $\omega^{\ZZ_2}_0(X;\, \xi -\eta )$)
does not, in general, depend only on the virtual $\ZZ_2$-vector bundle.

\smallskip

The next construction is the key that allows us to relate equivariant bordism
to stable homotopy.
\begin{definition}
Let $F_X^{\res}(\xi ,\eta )\to X$ be the infinite Stiefel bundle
$$
\dirlim_{k\geqslant 0}
\dirlim_{i\geqslant 0}
\O_X (\xi\oplus \Rr^i, \eta\oplus \Rr^i\oplus\Rr^k)
$$
formed as the direct limit over the injective maps
$$
\O_X (\xi\oplus \Rr^i, \eta\oplus \Rr^i\oplus\Rr^k)  \to
\O_X (\xi\oplus \Rr^{i+1}, \eta\oplus \Rr^{i+1}\oplus\Rr^k)
$$
given by the direct sum with the identity $1:\Rr \to \Rr$ in the decomposition
$\Rr^{i+1}=\Rr^i\oplus\Rr$
and the injective maps
$$
\O_X (\xi\oplus \Rr^i, \eta\oplus \Rr^i\oplus\Rr^k)  \to
\O_X (\xi\oplus \Rr^i, \eta\oplus \Rr^i\oplus\Rr^{k+1})
$$
given by including $\Rr^k$ in $\Rr^{k+1}=\Rr^k\oplus\Rr$.

The group $\ZZ_2$ acts on the space $F^\res_X(\xi ,\eta )$ in the obvious way.
\end{definition}

To be precise, we topologize $F^\res_X(\xi ,\eta )$ as the union of
the ANR subspaces $\O_X (\xi\oplus\Rr^i,\eta\oplus\Rr^i\oplus\Rr^k)$
with the weak topology. Any map from a compact ENR to $F^\res_X( \xi ,\eta )$
factors through the inclusion of one of these subspaces, so that
in practice we can work with
$\O_X (\xi\oplus\Rr^i,\eta\oplus\Rr^i\oplus\Rr^k)$
for sufficiently large $i$ and $k$.
\begin{lemma} \label{B-nonequF}
The projection $F^\res_X(\xi ,\eta ) \to X$ is a non-equivariant
homotopy equivalence.
\end{lemma}
\begin{proof}
Indeed, the fibres are non-equivariantly contractible, because of the limit
over the inclusions $\Rr^k\into\Rr^{k+1}$.
\hfill\qed\end{proof}
\begin{lemma}\label{B-pullback}
Suppose that $h : X'\to X$ is a $\ZZ_2$-map of $\ZZ_2$-ANRs.
Let $\xi'$ and $\eta'$ be the pullbacks of $\xi$ and $\eta$ to $X'$.
Then $F^\res_{X'}(\xi',\eta')\to X'$ is the pullback of
$F^\res_X(\xi ,\eta )\to X$ by $h$.
\qed
\end{lemma}

The connection between
bordism and stable homotopy is made by the Pontrjagin-Thom construction.
\begin{definition}
The {\it framed Pontrjagin-Thom homomorphism}
$$
\PT :\Omega^{\ZZ_2,\res}_0(X;\,\xi ,\eta ) \to
\omega^{\ZZ_2}_0(F_X^\res (\xi,\eta );\, \xi -\eta ).
$$
is defined as follows.
Consider a map $f:M\to X$ and a restricted stable isomorphism
$a:\tau M\oplus  f^*\xi\iso f^*\eta$
given by a vector bundle isomorphism $\tau M\oplus f^*\xi\oplus\Rr^i
\iso f^*\eta\oplus \Rr^i$ for some $i$,
representing a restricted bordism class $(M,f,a)$.
By including $f^*\xi\oplus \Rr^i$ in $\tau M\oplus f^*\xi\oplus \Rr^i$ we get
a lift of $f$ to a map $\tilde f:
M\to \O_X (\xi\oplus \Rr^i,\eta\oplus \Rr^i)\to F_X^\res (\xi ,\eta )$.
The standard Pontrjagin-Thom construction applied to $\tilde f$ gives
a class
$$
\PT (M,f,a)\in\omega_0^{G}(F_X^\res (\xi ,\eta );\,\xi-\eta ).
$$

We must check that this class is $0$ if $(M,f,a)$ is the boundary of
$(W,g,b)$,  where $g: W\to X$ extends $f$ on $\partial W=M$
and $b:\tau W\oplus g^*\xi \iso \Rr\oplus g^*\eta$
is a restricted stable isomorphism extending $a$ on $M$.
From a bundle isomorphism $\tau W\oplus g^*\xi\oplus \Rr^i\iso \Rr\oplus
g^*\eta\oplus \Rr^i$ we get a lift of $g$ to a map
$\tilde g:W\to\O_X (\xi\oplus \Rr^i,\Rr\oplus\eta\oplus \Rr^i)\to
F_X^\res (\xi ,\eta )$.
Hence the framed Pontrjagin-Thom construction on an associated lift $\tilde f$
of $f$ gives $0$.
\end{definition}
This verification that the framed Pontrjagin-Thom homomorphism is well-defined
explains the definition of the bundle $F_X^\res (\xi ,\eta )$.

\smallskip

The principal result of this section can now be stated as follows.
\begin{theorem} \label{B-mainthm}
The framed Pontrjagin-Thom homomorphism
$$
\PT : \Omega^{\ZZ_2,\res}_0 (X;\, \xi ,\eta ) \to
\omega^{\ZZ_2}_0 (F_X^\res (\xi ,\eta );\, \xi -\eta )
$$
is an isomorphism.
\end{theorem}

It will be proved by analysing separately the free and fixed-point
groups in bordism and stable homotopy.
\begin{definition}
We define the {\it free bordism group}
$\Omega_0^{\ZZ_2 ,\free}(X;\, \xi ,\eta )$
to be the bordism group of closed free $\ZZ_2$-manifolds $M$ with
a $\ZZ_2$-map $ f : M\to X$ and a stable $\ZZ_2$-isomorphism
$a: \tau M\oplus f^*\xi \iso f^*\eta$.
(Bounding manifolds $W$ are also required to have a free action of
$G$.)
\end{definition}
\begin{definition}
There is a {\it free Pontrjagin-Thom homomorphism}
$$
\PT^\free :
\Omega_0^{\ZZ_2 ,\free}(X;\, \xi ,\eta )
\to \omega^{\ZZ_2}_0(E\ZZ_2 \times X ;\, \xi -\eta )
$$
given by the classical construction using the product $M\to E\ZZ_2 \times X$
of the classifying map $M\to E\ZZ_2$ of the free action with $f$.
\end{definition}
\begin{proposition}
Write $\bar X = E\ZZ_2\times_{\ZZ_2} X$ and let $\bar\xi$ and $\bar\eta$ be the
induced bundles over $\bar X$.
Then we have isomorphisms
$$
\begin{matrix}
\Omega_0( \bar{X}; \, \bar{\xi},\bar{\eta}) & \iso &
\Omega^{\ZZ_2 ,\free }_0(X;\, \xi ,\eta )\\
\noalign{\smallskip}
\Darr{\iso}{} && \Darr{}{\PT^\free} \\
\noalign{\smallskip}
\omega_0(\bar X;\, \bar\xi -\bar\eta ) & \iso &
\omega^{\ZZ_2}_0 (E\ZZ_2 \times X ;\, \xi -\eta )
\end{matrix}
$$
through which $\PT^\free$ corresponds to the classical Pontrjagin-Thom
homomorphism.
\end{proposition}
\begin{proof}
This reduces, by naturality with respect to inclusions,
to showing that, for a free $\ZZ_2$-manifold $M$,
the equivariant fundamental class $[M]\in\omega_0^{\ZZ_2}(M;\, -\tau M)$
corresponds to the non-equivariant fundamental class
$[\bar M]\in\omega_0(\bar M;\, -\tau\bar M)$
of $\bar M=M/\ZZ_2$
under the isomorphism:
$$
\omega_0(\bar M;\, -\tau\bar M) \iso \omega_0^{\ZZ_2}(M;\, -\tau M) \iso
\omega_0^{\ZZ_2}(E\ZZ_2 \times M;\, -\tau M).
$$
But the isomorphism is given by the equivariant Umkehr map
$$
\pi^!:
\omega_0^{\ZZ_2}(\bar M;\, -\tau \bar M) \to \omega_0^{\ZZ_2}(M;\,-\tau M)
$$
for the projection $\pi : M\to \bar M$,
and this takes the fundamental class
of $\bar M$ with the trivial involution to the fundamental class of $M$.
\hfill\qed\end{proof}
Now we can apply the classical theory of framed bordism
to deduce that $\PT^\free$ is an isomorphism.
\begin{corollary}
The free Pontrjagin-Thom homomorphism
$$
\PT^\free  :
\Omega_0^{\ZZ_2 ,\free}(X;\, \xi ,\eta )
\to \omega^{\ZZ_2}_0(E\ZZ_2 \times X ;\, \xi -\eta )
$$
is an isomorphism.
\qed
\end{corollary}
\begin{definition}
Using the fact, recorded in Lemma \ref{B-resfree},
that stable isomorphisms over free $\ZZ_2$-spaces correspond
to restricted stable isomorphisms we get a homomorphism
$$
\gamma : \Omega^{\ZZ_2, \free}_0(X;\, \xi, \eta )
\to \Omega^{\ZZ_2, \res}_0(X;\, \xi, \eta ).
$$
\end{definition}

\begin{lemma}\label{B-gamma}
We have a commutative diagram
$$
\begin{matrix}
 \Omega^{\ZZ_2, \free}_0(X;\, \xi, \eta )
&\Rarr{\gamma}{} & \Omega^{\ZZ_2, \res}_0(X;\, \xi, \eta )\\
\noalign{\smallskip}
\Darr{\iso}{\PT^\free} && \Darr{}{\PT} \\
\noalign{\smallskip}
\omega^{\ZZ_2}_0(E\ZZ_2 \times F^\res_X(\xi ,\eta );\, \xi -\eta )
&\Rarr{}{} &
\omega^{\ZZ_2}_0(F^\res_X(\xi ,\eta );\, \xi -\eta )\, .
\end{matrix}
$$
\end{lemma}
\begin{proof}
We use the non-equivariant homotopy equivalence
$F^\res_X(\xi ,\eta ) \to X$, Lemma \ref{B-nonequF}, to get an isomorphism
$$
\omega_0^{\ZZ_2} (E\ZZ_2 \times F^\res_X(\xi ,\eta );\, \xi -\eta )
\to\omega_0^{\ZZ_2} (E\ZZ_2 \times X;\, \xi -\eta ).
$$
The assertion is then clear from the definition of the Pontrjagin-Thom
maps.
\hfill\qed\end{proof}
\begin{lemma}
Let $\xi_{\ZZ_2}$ and $\eta_{\ZZ_2}$ be, respectively,
the orthogonal complements of $\xi^{\ZZ_2}$ in $\xi |X^{\ZZ_2}$ and
of $\eta^{\ZZ_2}$ in $\eta | X^{\ZZ_2}$.

The fixed subspace $F_X^\res (\xi ,\eta )^{\ZZ_2}$
is homotopy equivalent to
the bundle of $\ZZ_2$-monomorphisms.
$$
\O^{\ZZ_2}_{X^{\ZZ_2}}(\xi_{\ZZ_2},\eta_{\ZZ_2})
=\O_{X^{\ZZ_2}}(\xi_{\ZZ_2},\eta_{\ZZ_2})
$$
over $X^{\ZZ_2}$.
\end{lemma}
\begin{proof}
We observe that the group $\ZZ_2$ acts on $\xi_{\ZZ_2}$ and $\eta_{\ZZ_2}$
as the involution $-1$.
The fixed subspace is, by definition,
$$
\left(\dirlim_k\dirlim_i
\O_{X^{\ZZ_2}}(\xi^{\ZZ_2}\oplus\Rr^i, \eta^{\ZZ_2}\oplus\Rr^i\oplus\Rr^k)\right)
\times_{X^{\ZZ_2}}\O^{\ZZ_2}_{X^{\ZZ_2}} (\xi_{\ZZ_2},\eta_{\ZZ_2}).
$$
The result follows from Lemma \ref{B-nonequF}.
\hfill\qed\end{proof}
We recall the fundamental {\it localization exact sequence} in
$\ZZ_2$-equivariant stable homotopy for the $\ZZ_2$-space $F^\res_X(\xi ,\eta)$.
See, for example, Crabb \cite[Lemma (A.1)]{crabb}.
\begin{proposition}\label{B-loc}
There is a long exact sequence
$$
\cdots
\Rarr{\partial}{}\omega^{\ZZ_2}_0(E\ZZ_2\times F^\res_X(\xi ,\eta );\, \xi -\eta )
\Rarr{\gamma}{}
\omega^{\ZZ_2}_0(F^\res_X(\xi ,\eta );\, \xi -\eta )
$$
$$
\Rarr{\rho}{}
\omega_0(F^\res_X(\xi ,\eta )^{\ZZ_2};\, \xi^{\ZZ_2}-\eta^{\ZZ_2})
\Rarr{\partial}{} \cdots
$$
\end{proposition}
\begin{proof}
The homomorphism $\gamma$ is induced by the projection $E\ZZ_2\to *$.

The fixed-point map $\rho$
sends a stable $\ZZ_2$-map to its
restriction to the subspaces fixed by $\ZZ_2$.

The boundary homomorphism $\partial$ is the composition
$$
\omega_0(F^\res_X(\xi ,\eta )^{\ZZ_2};\, \xi^{\ZZ_2}-\eta^{\ZZ_2})
\to\omega_0^{\ZZ_2}(D(\eta_{\ZZ_2})\times_{X^{\ZZ_2}}
F^\res_X(\xi ,\eta )^{\ZZ_2};\, \xi^{\ZZ_2}-\eta^{\ZZ_2})
$$
$$
\Rarr{\partial}{}\omega^{\ZZ_2}_{-1}(S(\eta_{\ZZ_2}) \times_{X^{\ZZ_2}}
F^\res_X(\xi ,\eta )^{\ZZ_2};\, \xi^{\ZZ_2}-(\eta^{\ZZ_2}\oplus\eta_{\ZZ_2}))
$$
$$
\to \omega^{\ZZ_2}_{-1}(E\ZZ_2 \times F_X^\res (\xi ,\eta );\, \xi -\eta )
$$
of the map induced by the group homomorphism $\ZZ_2 \to 0$ to the trivial group,
the boundary homomorphism $\partial$ in the exact sequence
of the pair $(D(\eta_{\ZZ_2}),S(\eta_{\ZZ_2}))$
and the map induced by the product of the classifying map
$S(\eta_{\ZZ_2})\to E\ZZ_2$ of the free action on the sphere bundle
and the inclusion of $F^\res_X(\xi, \eta)^{\ZZ_2}$ in
$F^\res_X(\xi ,\eta )$.
\hfill\qed\end{proof}
The main theorem will be established by comparing the localization sequence
in stable homotopy
with the geometrically defined Conner-Floyd exact sequence in
equivariant bordism that we now construct.
\begin{definition}
The {\it fixed point homomorphism}
$$
\rho :
\Omega^{\ZZ_2,\res }_0(X;\, \xi ,\eta )
\to
\Omega_0(F_X^\res (\xi ,\eta)^{\ZZ_2};\,\xi^{\ZZ_2},\eta^{\ZZ_2})
$$
is defined as follows.

Given a $\ZZ_2$-manifold $M$,
a map $f: M\to X$ and a restricted stable isomorphism $a: \tau M\oplus
f^*\xi \iso f^*\eta$,
we form the fixed subspace $M^{\ZZ_2}$.
We obtain, by restricting $f$ and $a$, a map $f^H: M^H\to X^H$ and
a stable isomorphism
$\tau M^H\oplus (f^H)^*\xi^H\iso (f^H)^*\eta^H$.

The classifying map $M\to F^\res_X(\xi ,\eta )$ restricts on fixed points
to a map $M^{\ZZ_2} \to F^\res_X(\xi ,\eta )^{\ZZ_2}$.

More precisely, we see that
the orthogonal complements $\xi_{\ZZ_2}$ and $\eta_{\ZZ_2}$
of $\xi^{\ZZ_2}$ and $\eta^{\ZZ_2}$, determine
the normal bundle $\nu$ of $M^{\ZZ_2}\into M$
by an isomorphism $\nu\oplus (f^{\ZZ_2})^*\xi_{\ZZ_2}\iso (f^{\ZZ_2})^*\eta_{\ZZ_2}$.
This gives us a map $M^{\ZZ_2}\to \O_{X^{\ZZ_2}}(\xi_{\ZZ_2},\eta_{\ZZ_2})$,
that is, a map $M^{\ZZ_2} \to \O_{X^{\ZZ_2}}^{\ZZ_2}(\xi_{\ZZ_2},\eta_{\ZZ_2})$.

The manifold $M^{\ZZ_2}$ with the extra data described above represents the
image of $(M,f,a)$ under $\rho$.
\end{definition}
\begin{lemma} \label{B-rho}
There is a commutative diagram
$$
\xymatrix{
\Omega^{\ZZ_2,\res }_0(X;\, \xi ,\eta )
\ar[r]^-{\rho} \ar[d]^-{\PT} &
\Omega_0(F_X^\res (\xi ,\eta)^{\ZZ_2};\,\xi^{\ZZ_2},\eta^{\ZZ_2})
\ar[d]^-{\iso}  \\
\omega^{\ZZ_2}_0(F^\res_X(\xi ,\eta );\, \xi -\eta )
\ar[r]^-{\rho} &
\omega_0(F_X^\res (\xi ,\eta)^{\ZZ_2};\,\xi^{\ZZ_2}-\eta^{\ZZ_2})
}
$$
\end{lemma}
\begin{proof}
The right-hand map is the classical Pontrjagin-Thom homomorphism,
so an isomorphism.
Commutativity follows directly from the construction of the Pontrjagin-Thom
maps.
\hfill\qed\end{proof}
We need to extend our definition of bordism groups to a theory indexed
by the integers.
Notice, first, that there is a natural identification
$$
\Omega_0^{\ZZ_2,\res}(X;\, \xi ,\eta ) = \Omega_0^{\ZZ_2,\res}(X;\, \Rr\oplus\xi,
\Rr\oplus\eta )
$$
and, indeed,
$$
\Omega_0^{\ZZ_2,\res}(X;\, \xi ,\eta ) = \Omega_0^{\ZZ_2,\res}(X;\, V_0\oplus\xi,
V_0\oplus\eta )
$$
for any vector space $V_0$ on which $\ZZ_2$ acts trivially.
\begin{definition}
We introduce groups $\Omega_*^{\ZZ_2 ,\res}(X;\, \xi ,\eta )$ for
$*\in\Zz$ so that
$$
\Omega_0^{\ZZ_2 ,\res}(X;\, \Rr^m\oplus\xi , \Rr^n\oplus\eta )
=\Omega_{n-m}^{\ZZ_2 ,\res}(X;\, \xi , \eta )
$$
for $m,\, n\geqslant 0$.
(To be precise, the definition is made so as to be compatible with the
identification
$$
\Omega_0^{\ZZ_2,\res}(X;\, \Rr^{m+1}\oplus\xi ,\Rr^{m+1}\oplus\eta ) =
\Omega_0^{\ZZ_2,\res}(X;\, \Rr\oplus(\Rr^m\oplus \xi), \Rr\oplus(\Rr^n\oplus\eta ))
$$
in the order $\Rr^{m+1}=\Rr\oplus\Rr^m$, $\Rr^{n+1}=\Rr\oplus\Rr^n$.)
\end{definition}
In particular,
$$
\Omega_{-1}^{\ZZ_2 ,\res}(X;\, \xi , \eta )
=\Omega_0^{\ZZ_2 ,\res}(X;\, \Rr\oplus\xi , \eta ),
$$
\begin{definition} \label{B-partialdef}
The {\it boundary homomorphism}
$$
\partial :
\Omega_0(F_X^\res (\xi ,\eta)^{\ZZ_2};\,\xi^{\ZZ_2},\eta^{\ZZ_2})
\to \Omega^{\ZZ_2, \free}_{-1}(X;\, \xi ,\eta )
$$
is constructed as follows.

Consider a closed manifold $N$ with a map $g:N\to X^{\ZZ_2}$,
a stable isomorphism $\tau N\oplus g^*\xi^{\ZZ_2}\iso g^*\eta^{\ZZ_2}$,
given by a bundle isomorphism $\tau N\oplus g^*\xi^{\ZZ_2}\oplus V_0\iso
g^*\eta^{\ZZ_2}\oplus V_0$, and a
bundle $\ZZ_2$-monomorphism $g^*\xi_{\ZZ_2}\into g^*\eta_{\ZZ_2}$
with complementary bundle $\nu$ supplied by a lift of $g$ to a map
$$
N\to F^\res_X(\xi ,\eta )^{\ZZ_2} \to \O^{\ZZ_2}_{X^{\ZZ_2}}(\xi_{\ZZ_2},\eta_{\ZZ_2}).
$$

Now form the free $\ZZ_2$-manifold $S(\nu )$.
It is equipped with a $\ZZ_2$-map $h: S(\nu ) \to N
\to X^{\ZZ_2}\to X$.
We have a $\ZZ_2$-isomorphism $\Rr\oplus \tau S(\nu )
\iso \tau N\oplus \nu$.
Taking the direct sum with the isomorphism $\tau N\oplus g^*\xi^{\ZZ_2}\oplus V_0
\iso g^*\eta^H{\ZZ_2}\oplus V_0$ and $\nu\oplus g^*\xi_{\ZZ_2}\iso g^*\eta_{\ZZ_2}$ we
get a $\ZZ_2$-isomorphism
$$
\Rr\oplus \tau S(\nu )\oplus h^*\xi\oplus V_0
\iso h^*\eta\oplus V_0.
$$

The manifold $S(\nu)$ and the associated data
represents the required boundary class in $\Omega_{-1}^{\ZZ_2,\free}
(X;\, \xi ,\eta )$.
\end{definition}
\begin{lemma} \label{B-delta}
We have a commutative diagram
$$
\xymatrix{
\Omega_0(F_X^\res (\xi ,\eta)^{\ZZ_2};\,\xi^{\ZZ_2},\eta^{\ZZ_2})
\ar[d]^-{\iso}
\ar[r]^-{\partial} &
\Omega_{-1}^{\ZZ_2 ,\free}(X;\, \xi ,\eta )
 \ar[d]^-{{\PT^\free}{}} \\
\omega_0(F_X^\res (\xi ,\eta)^{\ZZ_2};\,\xi^{\ZZ_2}-\eta^{\ZZ_2}) \ar[r]^-{\partial} &
\omega_{-1}^{\ZZ_2}(E\ZZ_2 \times F^\res_X(\xi ,\eta );\, \xi -\eta )}$$
\end{lemma}
\begin{proof}
We look at an element of
$\Omega_0(F^\res_X(\xi ,\eta )^{\ZZ_2};\, \xi^{\ZZ_2}, \eta^{\ZZ_2})$
represented by a
manifold $N$ with a map $g:N\to X^H$,
a stable isomorphism $\tau N\oplus g^*\xi^{\ZZ_2}\iso g^*\eta^{\ZZ_2}$,
given by a bundle isomorphism $\tau N\oplus g^*\xi^{\ZZ_2}\oplus V_0\iso
g^*\eta^{\ZZ_2}\oplus V_0$, and a
bundle $\ZZ_2$-monomorphism $g^*\xi_{\ZZ_2}\into g^*\eta_{\ZZ_2}$
with complementary bundle $\nu$, as in the definition of
$\partial$ (Definition \ref{B-partialdef}).
By the obvious naturality under inclusions we may assume that
$X=N=X^{\ZZ_2}$, $g$ is the identity map, $\xi =0$ and
$\eta =\tau N\oplus \nu$, so that $\eta^{\ZZ_2}=\tau N$ and
$\eta_{\ZZ_2}=\nu$.

The image of the element
in $\omega_0(X^{\ZZ_2};\, \xi^{\ZZ_2}-\eta^{\ZZ_2})=\omega_0(N;\, -\tau N)$
is the fundamental class $[N]$ of $N$.

Its image in $\Omega^{\ZZ_2 ,\free}_{-1}(X;\, \xi ,\eta )
=\Omega_{-1}^{\ZZ_2 ,\free}(N;\, 0,\tau N\oplus\nu )$
is represented by the sphere bundle
$S(\nu )\to N$ with $\Rr\oplus\tau S(\nu )=\tau N\oplus\nu$.
This maps by the Pontrjagin-Thom construction to the $\ZZ_2$-fundamental class
$[S(\nu )]\in \omega^{\ZZ_2}_0(S(\nu );\, -\tau S(\nu ))
=\omega_{-1}^{\ZZ_2}(S(\nu );\, -\tau N-\nu )$
and then by the classifying map $S(\nu )\to E\ZZ_2 \times N$
to the group
$\omega^{\ZZ_2}_{-1}(E\ZZ_2 \times N;\, -\tau N -\nu )$.

From the construction of the localization exact sequence we have
a commutative square
$$
\begin{matrix}
\omega_0^{\ZZ_2}(N;\, -\tau N) & \Rarr{\partial}{} &
\omega_{-1}^{\ZZ_2}(S(\nu );\, -\tau N -\nu ) \\
\noalign{\smallskip}
\Darr{\rho}{} & & \Darr{}{} \\
\noalign{\smallskip}
\Omega_0(N;\, -\tau N) & \Rarr{\partial}{} &
\omega_{-1}^{\ZZ_2}(E\ZZ_2 \times N;\,-\tau N-\nu )
\end{matrix}
$$
The assertion now follows from the fact that the equivariant fundamental
class $[N]$ in the top line is mapped by $\partial$ to the
fundamental class $[S(\nu )]$, because $\partial$ in the exact sequence
of the pair $(D(\nu ),S(\nu ))$ is the Umkehr
homomorphism of the projection $S(\nu )\to N$.
\hfill\qed\end{proof}
\begin{proposition}
The {\it Conner-Floyd sequence}:
$$
\cdots
\Rarr{\partial}{}
\Omega^{\ZZ_2,\free}_0(X;\, \xi ,\eta )
\Rarr{\gamma}{}
\Omega^{\ZZ_2,\res }_0(X;\, \xi ,\eta )
\Rarr{\rho}{}
$$
$$
\Omega_0(F_X^\res (\xi ,\eta)^{\ZZ_2};\,\xi^{\ZZ_2},\eta^{\ZZ_2})
\Rarr{\partial}{}
\Omega^{\ZZ_2,\free}_{-1}(X;\, \xi ,\eta ) \Rarr{\gamma}{}\cdots
$$
is exact.
\end{proposition}
\begin{proof}{\it (Outline)}
We use the standard surgery arguments;
compare the proof in Crabb, Mishchenko, Morales Mel\'endez and Popelensky \cite[Theorem 4.2]{cmmp}.

It is clear that $\rho\circ\gamma =0$.
We show that $\ker\rho \subseteq\im\gamma$.
Suppose that a class represented by $f: M\to X$ and
$a :\tau M\oplus f^*\xi \iso f^*\eta$ lies in the kernel of $\rho$.
We have a manifold $W$ with $\partial W =X^{\ZZ_2}$, a map
$g : W \to X^{\ZZ_2}$ and a stable isomorphism
$b : \tau W \oplus g^*\xi^{\ZZ_2}\iso \Rr\oplus g^*\eta^{\ZZ_2}$
extending $a$.
Further, we have a vector bundle $\hat\nu$ over $W$ extending the
normal bundle $\nu$ over $M^{\ZZ_2}$ and an isomorphism
$\hat\nu \oplus g^*\xi_{\ZZ_2} \iso g^*\eta_{\ZZ_2}$.

Take an equivariant tubular neighbourhood $D(\nu )\into M$
of $M^{\ZZ_2}$ in $M$.
We construct a free $\ZZ_2$-manifold
$M'= (M-B(\nu ))\cup_{S(\nu )} S(\hat\nu )$
cobordant to $M$.
From $f$ and $g$ we get a map $f' : M'\to X$.
From $a$ and $b$ we get a restricted stable isomorphism
$a' : \tau M'\oplus (f')^*\xi \iso (f')^*\eta$.
Notice that $\Rr\oplus\tau S(\hat\nu )=\tau W\oplus \hat\nu$,
so that $\Rr\oplus \tau S(\hat\nu )\oplus g^*\xi^{\ZZ_2}\oplus g^*\xi_{\ZZ_2}
\iso\tau W \oplus\hat\nu \oplus g^*\xi^{\ZZ_2} \oplus g^*\xi_{\ZZ_2}
\iso \Rr\oplus g^*\eta^{\ZZ_2}\oplus g^*\eta_{\ZZ_2}$.
This gives a lift to $\Omega_0^{\ZZ_2 ,\free}(X;\, \xi ,\eta )$
 and shows that
$\ker\rho \subseteq \im \gamma$.

The composition $\gamma\circ\partial$ is zero, because $S(\nu )$ is the
boundary of $D(\nu)$.
To see that $\ker\gamma \subseteq \im\partial$, suppose that a manifold
$S$ represents a class in $\Omega^{\ZZ_2 ,\free}_{-1}(X;\,\xi ,\eta )$
that is the boundary of a manifold $M$. Then the fixed-point construction
applied to $M$
produces an element of $\Omega_0(F^\res_X(\xi ,\eta )^{\ZZ_2};\, \xi^{\ZZ_2},
\eta^{\ZZ_2})$ mapping to $S$ under $\partial$.
(Notice that there are no fixed points on the boundary $S$ of $M$.)

This same construction with $S=\emptyset$ shows that $\partial \circ\rho =0$.
We must check finally that $\ker\partial \subseteq \im\rho$.
Given $N$ representing a class in $\ker\partial$, so that $S(\nu)$
is the boundary of a manifold $M_0$, we form the closed
manifold $M=M_0 \cup D(\nu)$ by gluing along $\partial M_0=S(\nu )=
\partial D(\nu)$.
Then $M$ (with the associated structure) gives an element
of $\Omega_0^{\ZZ_2 ,\res}(X;\, \xi ,\eta )$ lifting the class of
$N$.
(In more detail, we have $g: N\to X^{\ZZ_2}$,
$\tau N\oplus g^*\xi^{\ZZ_2}\oplus V_0 \iso g^*\eta^{\ZZ_2}\oplus V_0$
and $\nu\oplus g^*\xi_{\ZZ_2}\iso g^*\eta_{\ZZ_2}$.
We have $h_0 : M_0 \to X$ and $\tau M_0\oplus h_0^*\xi \oplus V_0
\iso h_0^*\eta \oplus V_0$, for $V_0$ of sufficiently high dimension.
We get a map $f : M\to X$ from $h$ on $S(\nu)$ and $h_0$ on $M_0$
and a bundle isomorphism $\tau M\oplus f^*\xi \oplus V_0\iso f^*\eta\oplus
V_0$.)
\hfill\qed\end{proof}

\begin{proof} {\it (of Theorem \ref{B-mainthm})}
The theorem now follows by applying the five-lemma to the diagram
$$
\begin{matrix}
\Darr{\partial}{}&&\Darr{\partial}{} \\
\noalign{\smallskip}
\Omega^{\ZZ_2,\free}_0(X;\, \xi ,\eta ) & \Rarr{\iso}{} &
\omega^{\ZZ_2}_0 (E\ZZ_2 \times F_X^\res(\xi ,\eta );\, \xi -\eta )\\
\noalign{\smallskip}
\Darr{\gamma}{}&& \Darr{\gamma}{} \\
\noalign{\smallskip}
\Omega^{\ZZ_2,\res }_0(X;\, \xi ,\eta ) & \Rarr{\PT}{} &
\omega^{\ZZ_2}_0(F^\res_X(\xi ,\eta ); \, \xi -\eta )\\
\noalign{\smallskip}
\Darr{\rho}{}&&\Darr{\rho}{} \\
\noalign{\smallskip}
\Omega_0(F_X^\res (\xi ,\eta)^{\ZZ_2};\,\xi^{\ZZ_2},\eta^{\ZZ_2})
&\Rarr{\iso}{}& \omega_0(F_X^\res (\xi ,\eta)^{\ZZ_2};\,\xi^{\ZZ_2}-\eta^{\ZZ_2})\\
\noalign{\smallskip}
\Darr{\partial}{} &&\Darr{\partial}{} \\
\noalign{\smallskip}
\Omega^{\ZZ_2,\free}_{-1}(X;\, \xi ,\eta )&\Rarr{\iso}{} &
\omega^{\ZZ_2}_{-1} (E\ZZ_2 \times F_X^\res(\xi ,\eta );\, \xi -\eta )\\
\noalign{\smallskip}
\Darr{\gamma}{} &&\Darr{\gamma}{}
\end{matrix}
$$
relating the Conner-Floyd and localization exact sequences.
using Lemmas \ref{B-rho}, \ref{B-delta} and
\ref{B-gamma} to check commutativity of the various squares.
\hfill\qed\end{proof}
\section{Some deductions}
We begin by considering some cases in which the
Pontrjagin-Thom homomorphism
$$
\Omega^{\ZZ_2 ,\res}_0(X;\,\xi ,\eta )
\to\omega^{\ZZ_2}_0(X;\,\xi -\eta )
$$
is an isomorphism.

Our first result is a special case of a theorem
of Hauschild \cite[Satz IV.2]{hauschild1}.

\begin{proposition} \label{B-equ}
Suppose that $\xi$ is isomorphic to a subbundle of $X\times\Rr^n$
for some $n\geqslant 0$.
Then
$$
\Omega^{\ZZ_2 ,\res}_0(X;\,\xi ,\eta )
\iso\omega^{\ZZ_2}_0(X;\,\xi -\eta )\, .
$$
\end{proposition}
\begin{proof}
It is easy to see directly
from the definition of $F^\res_X(\xi ,\eta )$ as a direct limit that
the projection
$$
F^\res_X(\xi ,\eta ) \to X
$$
is a $\ZZ_2$-homotopy equivalence
with inverse given by a section of $\O_X^{\ZZ_2}(\xi ,\Rr^n)$.

(One can also see that the restriction to the fixed subspaces
$$
F^\res_X(\xi ,\eta )^{\ZZ_2}\simeq \O^{\ZZ_2}_{X^{\ZZ_2}}(\xi_{\ZZ_2},\eta_{\ZZ_2}) \to X^H
$$
is a non-equivariant homotopy equivalence, because $\xi_{\ZZ_2}=0$.)
\hfill\qed\end{proof}

If the action of $\ZZ_2$ on $\eta$ is (essentially) trivial, we may apply
the classical splitting of equivariant stable homotopy (as described,
for example, in tom Dieck \cite[Chapter II, Theorem (7.7)]{tdbook}).
\begin{proposition} \label{B-split}
Suppose that $\eta$ is isomorphic to a subbundle of $X\times\Rr^n$
for some $n\geqslant 0$.
Let $X^{\ZZ_2}_0$ be the union of those components of $X^{\ZZ_2}$ on which
$\xi^{\ZZ_2}$ is zero.

Then the bordism and stable homotopy groups split as direct sums
$$
\Omega_0^{\ZZ_2,\res}(X;\,\xi ,\eta )=
\omega_0(E\ZZ_2 \times_{\ZZ_2} X;\, \xi_\#-\eta_\#)\oplus
\omega_0(X^{\ZZ_2}_0;\, \xi^{\ZZ_2}-\eta^{\ZZ_2})
$$
and
$$
\omega_0^{\ZZ_2}(X;\,\xi -\eta)=
 \omega_0(E\ZZ_2 \times_{\ZZ_2} X;\, \xi_\#-\eta_\#)\oplus
\omega_0(X^{\ZZ_2};\, \xi^{\ZZ_2}-\eta^{\ZZ_2}),
$$
where $\xi_\#$ and $\eta_\#$ are the vector bundles over
$E\ZZ_2\times_{\ZZ_2} X$ associated with the
$\ZZ_2$-equivariant bundles $\xi$ and $\eta$ over $X$.
\end{proposition}
\begin{proof}
We have $\eta_{\ZZ_2}=0$, so that $\O^{\ZZ_2}_{X_{\ZZ_2}}(\xi_{\ZZ_2},\eta_{\ZZ_2})$ is empty over the
complement of $X^{\ZZ_2}_0$ and equal to $X^{\ZZ_2}_0$ over $X^{\ZZ_2}_0$.
\hfill\qed\end{proof}
\begin{example}
Suppose that $X=*$ is a point, $\xi =0$ and $\eta=\Rr^m$.
The general theory has established that
$$
\Omega_m^{\ZZ_2,\res}(*;\, 0,0) =
\Omega_0^{\ZZ_2,\res}(*;\, 0,\Rr^m) =
\omega_m(B\ZZ_2)\oplus \omega_m(*).
$$
It is easy to see this directly. For suppose that $M$ is an $m$-dimensional
closed $\ZZ_2$-manifold equipped with a restricted stable isomorphism
$\tau M\iso M\times\Rr^m$.
The fixed submanifold $M^{\ZZ_2}$ must also
have dimension $m$: $\tau M^{\ZZ_2}\iso M^{\ZZ_2}\times\Rr^m$. Thus, $M^{\ZZ_2}$ is
a union of components of $M$.

The manifold $M$ decomposes as a disjoint union
$$
M= M_{\free} \sqcup M^{\ZZ_2},
$$
where $M_{\free}$ is
the subspace of points with trivial isotropy group.
The free $\ZZ_2$-manifold $M_\free$ represents an element of $\omega_m(B\ZZ_2)$
and the fixed manifold $M^{\ZZ_2}$ represents an element of $\omega_m(*)$.

More generally, for any space $X$ there is a geometric
description of the splitting
$$
\Omega_m^{\ZZ_2,\res}(X;\, 0,0) =\omega_m(E\ZZ_2\times_{\ZZ_2} X)
\oplus \omega_m(X^{\ZZ_2}).
$$
\end{example}

We can also specify a range of dimensions in which the Pontrjagin-Thom
map gives an isomorphism between restricted bordism and stable homotopy.

\begin{proposition} \label{stab}
Suppose that, for each point $x\in X^{\ZZ_2}$ such that $(\xi_{\ZZ_2})_x\not=0$,
$$
\dim \eta^{\ZZ_2}_x -\dim \xi^{\ZZ_2}_x + 1 < \dim (\eta_{\ZZ_2})_x
-\dim (\xi_{\ZZ_2})_x \, .
$$
Then the Pontrjagin-Thom map
$$
\Omega^{\ZZ_2,\res}_0(X;\, \xi ,\eta )\to \omega_0^{\ZZ_2}(X;\, \xi -\eta )
$$
is an isomorphism.
\end{proposition}
\begin{proof}
We consider the commutative diagram of exact sequences
$$
\begin{matrix}
\omega_1(F_X^\res (\xi ,\eta)^{\ZZ_2};\,\xi^{\ZZ_2}-\eta^{\ZZ_2})
&\Rarr{\pi_*}{}& \omega_1(X^{\ZZ_2};\,\xi^{\ZZ_2}-\eta^{\ZZ_2})\\
\noalign{\smallskip}
\Darr{\partial}{} &&\Darr{\partial}{} \\
\noalign{\smallskip}
\omega^{\ZZ_2}_0(E\ZZ_2 \times F_X^\res(\xi ,\eta );\, \xi -\eta ) & \Rarr{\pi_*}{\iso} &
\omega^{\ZZ_2}_0 (E\ZZ_2 \times X;\, \xi -\eta )\\
\noalign{\smallskip}
\Darr{\gamma}{}&& \Darr{\gamma}{} \\
\noalign{\smallskip}
\omega^{\ZZ_2}_0(F^\res_X(\xi ,\eta); \, \xi -\eta ) & \Rarr{\pi_*}{} &
\omega^{\ZZ_2}_0(X; \, \xi -\eta )\\
\noalign{\smallskip}
\Darr{\rho}{}&&\Darr{\rho}{} \\
\noalign{\smallskip}
\omega_0(F_X^\res (\xi ,\eta)^{\ZZ_2};\,\xi^{\ZZ_2}-\eta^{\ZZ_2})
&\Rarr{\pi_*}{}& \omega_0(X^{\ZZ_2};\,\xi^{\ZZ_2}-\eta^{\ZZ_2})\\
\noalign{\smallskip}
\Darr{\partial}{} &&\Darr{\partial}{} \\
\noalign{\smallskip}
\omega^{\ZZ_2}_{-1}(E\ZZ_2 \times F_X^\res (\xi ,\eta);\, \xi -\eta )&\Rarr{\pi_*}{\iso} &
\omega^{\ZZ_2}_{-1} (E\ZZ_2 \times X;\, \xi -\eta )
\end{matrix}
$$
in which the horizontal maps are induced by the projection
$\pi : F^\res_X(\xi ,\eta )\to X$.

We can identify $F^\res_X(\xi ,\eta )^{\ZZ_2}\to X^{\ZZ_2}$ with
the Stiefel bundle
$$
\O_{X^{\ZZ_2}}(\xi_{\ZZ_2},\eta_{\ZZ_2})\to X^{\ZZ_2}.
$$
The fibre at $x$ is a single point if $(\xi_{\ZZ_2})_x=0$,
empty if $\dim (\xi_{\ZZ_2})_x > \dim (\eta_{\ZZ_2})_x$,
and satisfies
$\pi_i(\O ((\xi_{\ZZ_2})_x,(\eta_{\ZZ_2})_x)) =0$
for $i <  \dim (\eta_{\ZZ_2})_x-\dim (\xi_{\ZZ_2})_x$
if $\dim (\xi_{\ZZ_2})_x \leq \dim (\eta_{\ZZ_2})_x$.

The proof is completed by the five lemma.
(Notice that, if $\dim (\xi_{\ZZ_2})_x > \dim (\eta_{\ZZ_2})_x$,
we have $\dim \eta^{\ZZ_2}_x -\dim \xi^{\ZZ_2}_x + 1 <0$, so that all
the relevant groups are zero.)
\hfill\qed\end{proof}

The stable homotopy group $\omega_0^{\ZZ_2}(X;\, \xi -\eta )$ can
always be
represented as a direct limit of restricted bordism groups.

Suppose that $\zeta$ is a real $\ZZ_2$-vector bundle over $X$.
We define a homomorphism
$$
\Sigma_\zeta : \Omega^{\ZZ_2 ,\res}_0(X;\, \xi ,\eta )
\to \Omega^{\ZZ_2 ,\res}_0(S(\Rr\oplus\zeta );\, \xi ,\eta\oplus\zeta )
$$
by sending the class of a manifold $M$, $f: M\to X$ and a restricted stable
isomorphism $a : \tau M\oplus f^*\xi \to f^*\eta$ to the manifold
$S=S(\Rr\oplus f^*\zeta )$ (with a choice of smooth structure on $f^*\zeta$
-- different choices will give isomorphic smooth vector bundles)
with the map $S\to M \to X$ given by the
projection composed with $f$ and the restricted stable isomorphism
$$
\tau  S\oplus \Rr \oplus f^*\xi \iso \tau M\oplus f^*\zeta \oplus\Rr
\iso f^*\eta\oplus f^*\zeta\oplus\Rr .
$$
\begin{lemma}
There is a commutative diagram
$$
\begin{matrix}
\Omega^{\ZZ_2 ,\res}_0(X;\, \xi ,\eta ) &\Rarr{\Sigma_\zeta}{}&
\Omega^{\ZZ_2 ,\res}_0(S(\Rr\oplus\zeta );\, \xi ,\eta\oplus\zeta ) \\
\noalign{\smallskip}
\Darr{}{} && \Darr{}{} \\
\noalign{\smallskip}
\omega_0^{\ZZ_2}(X;\, \xi -\eta ) &\Rarr{}{\sigma_\zeta} &
\omega_0^{\ZZ_2}(X;\, \xi -\eta )\oplus \omega_0^{\ZZ_2}(X;\, \xi -\eta -\zeta ),
\end{matrix}
$$
in which the vertical maps are Pontrjagin-Thom homomorphisms
and $\sigma_\zeta$ is the inclusion of the first factor.
\end{lemma}
\begin{proof}
The homomorphism $\sigma_\zeta$ is the boundary $\partial$
in the exact sequence of the pair $(D(\Rr\oplus\zeta ),S(\Rr\oplus\zeta ))$.
The group $\omega_0^{\ZZ_2}(S(\Rr\oplus\zeta );\, \xi -\eta -\zeta )$
is split by the section $(-1,0)$ of $X\to S(\Rr \oplus \zeta )$.
\hfill\qed\end{proof}

\begin{proposition}
If $j$ is sufficiently large, there is a natural
isomorphism
$$
\Omega^{\ZZ_2 ,\res}_0(X\times S(\Rr\oplus L^j);\, \xi ,\eta\oplus L^j)
\iso
\omega^{\ZZ_2}_0(X;\,  \xi -\eta) \oplus
\Omega^{\ZZ_2, \res} _0(X;\,  \xi ,\eta \oplus L^j) \, .
$$
\end{proposition}
\begin{proof}
We take $\zeta$ to be the trivial bundle $L^j$ and
use the stability criterion Proposition \ref{stab}.

The point is that $\dim \eta_x^{\ZZ_2} =\dim \xi_x^{\ZZ_2} + 1
< \dim (\eta_{\ZZ_2})_x + j -\dim (\xi_{\ZZ_2})_x$
for $k$ sufficiently large.
\hfill\qed\end{proof}
\begin{remark}
By introducing relative bordism groups we could represent
the equivariant stable homotopy group
$\omega^{\ZZ_2}_0(X;\, \xi -\eta )$ exactly as the restricted bordism group
of the pair $X\times (S(\Rr \oplus L^j),*)$
with coefficients $(\xi ,\eta\oplus L^j)$.
\end{remark}

\begin{definition}
We define the (unrestricted) {\it framed bordism group}
$\Omega^{\ZZ_2}_0(X;\, \xi ,\eta )$ to be the bordism group
of closed $\ZZ_2$-manifolds $M$ with a $\ZZ_2$-map $f : M\to X$
and a stable isomorphism $a: \tau M\oplus f^*\xi \iso f^*\eta$,
specified by a vector bundle isomorphism $\tau M\oplus f^*\xi\oplus V
\iso f^*\eta\oplus V$ for some $\ZZ_2$-module $V$.
\end{definition}
Thus
$$
\Omega_0^{\ZZ_2}(X;\, \xi ,\eta )
=\dirlim_{V : V^{\ZZ_2}=0} \Omega_0^{\ZZ_2 ,\res}(X;\, V\oplus\xi , V\oplus\eta ).
$$
\begin{definition}
Let $F_X^{\ZZ_2}(\xi ,\eta )$ be the $\ZZ_2$-bundle
$$
\dirlim_{V : V^{\ZZ_2}=0} F^\res_X(V\oplus \xi ,V\oplus \eta ),
$$
or. more concretely,
$$
\dirlim_{j \geqslant 0} F^\res_X(L^j\oplus \xi ,L^j\oplus \eta )\, .
$$
We have an equivalence between $F^{\ZZ_2}_X(\xi ,\eta )$ and $F^{\ZZ_2}_X(V\oplus\xi,
V\oplus\eta )$ for any $\ZZ_2$-module $V$. Thus $F^{\ZZ_2}_X(\xi ,\eta )$
depends only on the virtual $\ZZ_2$-vector bundle $\xi -\eta$.
\end{definition}
Our results on the restricted bordism groups now describe $\Omega_0^{\ZZ_2}
(X;\,\xi ,\eta )$ as a stable homotopy group.
We extend the framed Pontrjagin-Thom homomorphism algebraically to
the direct limits.
\begin{theorem}
The framed Pontrjagin-Thom homomorphism
$$
\PT : \Omega_0^{\ZZ_2}(X;\, \xi ,\eta ) \to
\omega_0^{\ZZ_2}(F^{\ZZ_2}_X(\xi ,\eta );\, \xi -\eta )
$$
is an isomorphism.
\qed
\end{theorem}
\begin{lemma}
The fixed subspace $F_X^{\ZZ_2}(\xi ,\eta )^{\ZZ_2}$ is homotopy equivalent
to the bundle
$$
\dirlim_{j\geqslant 0} \O^{\ZZ_2}_{X^{\ZZ_2}}(L^j\oplus\xi_{\ZZ_2},L^j\oplus\eta_{\ZZ_2})
=\dirlim_{j\geqslant 0} \O_{X^{\ZZ_2}}(L^j\oplus\xi_{\ZZ_2},L^j\oplus\eta_{\ZZ_2}).
$$
It depends only on the virtual bundle $\xi_{\ZZ_2}-\eta_{\ZZ_2}$ over $X^{\ZZ_2}$.
\qed
\end{lemma}

\begin{example}
We look again at the special case $X=*$, $\xi =0$, $\eta =\Rr^m$.

Using Proposition \ref{B-split} we can write
$$
\Omega^{\ZZ_2}_m(*;\, 0,0) = \omega_m(B\ZZ_2 ) \oplus
\dirlim_j\, \omega_m(\O^{\ZZ_2}(L^j,L^j))
=\omega_m(B\ZZ_2) \oplus \omega_m(\O (\infty)).
$$
\end{example}
\begin{example}
It is instructive to look in detail at the $0$-dimensional case.
A closed $0$-dimensional manifold is just a finite set, and a
$0$-dimensional $\ZZ_2$-manifold $M$ is a finite $\ZZ_2$-set with tangent bundle $\tau M=M\times 0$.

A $\ZZ_2$-manifold $M$ consisting of a single point has two restricted
stable framings given by the isomorphisms $\pm 1: \tau M\oplus \Rr
\to \tau M\oplus\Rr$. The class $[M,1]$ generates the summand
$\omega_0(*)=\Zz$ in
$$
\Omega^{\ZZ_2 ,\res}_0(*;\, 0,0) = \omega_0(B\ZZ_2 )\oplus \omega_0(*)
=\Zz \oplus \Zz\, ,
$$
and $[M,-1]=-[M,1]$.
The manifold has two further stable framings represented by the isomorphism
$-1: \tau M\oplus L \to \tau M\oplus L$, which we shall call $t$
and $-1 : \tau M\oplus\Rr \oplus L \to \tau M\oplus\Rr\oplus L$,
which we call $-t$. Then $[M,-t]=-[M,t]$ generates the second
$\Zz$ summand in $\omega_0(\O (\infty ))=\Zz \oplus \Zz=\Zz [M,1]\oplus\Zz [M,t]$
in the framed bordism group
$$
\Omega^{\ZZ_2}_0(*;\, 0,0) = \omega_0(B\ZZ_2 )\oplus \omega_0(\O (\infty ))
=\Zz \oplus (\Zz\oplus\Zz )\, .
$$

A free $\ZZ_2$-manifold $N$ with two points has two $\ZZ_2$-equivariant
stable framings given by $\pm 1: \tau N\oplus\Rr \to \tau N\oplus \Rr$ and $[N,1]=-[N,-1]$ generates the summand $\omega_0(B\ZZ_2 )=\Zz$ in
$\Omega^{\ZZ_2 ,\res}_0(*;\, 0,0)$ and $\Omega^{\ZZ_2}_0(*;\, 0,0)$.

(Of course, $N$ can also be framed as the boundary of $D(L)$
by an isomorphism
$\tau N \oplus \Rr \to \tau N\oplus L$ to represent $0$ in
$\Omega^{\ZZ_2 ,\res}_0(*;\, \Rr ,L)$.)

Any $0$-dimensional $\ZZ_2$-manifold with a restricted $(0,0)$ framing
is equivalent to
a disjoint union of copies of $(M,1)$, $(M,-1)$, $(N,1)$ and $(N,-1)$.
And any $(0,0)$ framed $\ZZ_2$-manifold of dimension $0$ is equivalent
to a disjoint union of copies of these four manifolds, $(M,t)$ and $(M,-t)$.
\end{example}
We conclude this section by relating the description of $\Omega_0^{\ZZ_2} (X;\,\xi ,\eta )$
using the stable structure of the tangent bundle to the classical
description involving normal maps.
\begin{definition}
Suppose that $\eta =X\times U$ is trivial. A {\it normal map}
is prescribed by a map $f: M\to X$,
an equivariant embedding $j:M\into U\oplus V$, for some $\ZZ_2$-module $V$, and
an equivariant vector bundle isomorphism $a:\nu_M \to f^*\xi \oplus V$.
\end{definition}

Such a normal map
determines a vector bundle
isomorphism
$$
\tau M\oplus f^*\xi\oplus V \to U\oplus V
$$
and so a bordism class in $\Omega^{\ZZ_2}_0(X;\, \xi , U)$.
Its image in
$\omega_0^{\ZZ_2}(F_X^{\ZZ_2}(\xi ,U);\, \xi -U)$
under $\PT$
can be described directly in terms of the normal map.
The Pontrjagin-Thom construction applied to a tubular neighbourhood
composed with the classifying map of the vector bundle
isomorphism $a$ gives a map
$$
\begin{matrix}
(U\oplus V)^\cpt &\to &\Thom (\nu_M ,M) &\to&
\Thom (\xi\oplus V, F^{\ZZ_2}_X(\xi , U)) \\
\Uarr{=}{} &&&& \Darr{=}{}\\
U^\cpt\wedge V^\cpt &&&&
\Thom (\xi , F_X^{\ZZ_2}(\xi ,U)) \wedge V^\cpt\, ,
\end{matrix}
$$
which represents the stable class.
(Recall that $V^\cpt$ is the one-point compactification
of the vector space $V$ and  $\Thom (\nu_M,M)$ is the Thom space
of the vector bundle $\nu_M$ over $M$.)
\begin{proposition}
Two normal maps
$$
f_i : M\to X,\,
j_i: M\into U\oplus V_i,\,
a_i: \nu (j_i) \to f_i^*\xi \oplus V_i, \, (i=0,\, 1),
$$
determine the same
cobordism class in $\Omega_0^{\ZZ_2}(X;\, \xi , U)$ if and only
if they are cobordant in the sense that there is a $\ZZ_2$-manifold $W$ with
boundary $\partial W =M_0\sqcup M_1$, a $\ZZ_2$-map $g: W \to X$ restricting
to $f_i$ on $M_i$, an embedding $W\into U\oplus (V_0\oplus V_1\oplus V)$,
for some $\ZZ_2$-module $V$, and an equivariant vector bundle isomorphism
$\nu_W \to g^*\xi \oplus (V_0\oplus V_1\oplus V)$ extending, in the obvious
way, the embedding $j_i$ and isomorphism $a_i$ on $M_i$.
\qed
\end{proposition}

\section{Unoriented bordism}
In this section we shall explain how $\ZZ_2$-equivariant bordism and
related theories can be described first as restricted framed bordism
groups and then as $\ZZ_2$-stable homotopy groups.
Throughout this section $Y$ is a $\ZZ_2$-ANR.

We write the infinite Grassmannian of $m$-dimensional real subspaces of $\bigcup_{i, \, j\geqslant 0} \Rr^i\oplus L^j$ as
$$
G_m(\infty \oplus\infty L) =\dirlim_{i,j\geqslant 0} G_m(\Rr^i\oplus L^j),
$$
and let $\eta_m$ denote the canonical $m$-dimensional vector bundle.
The Grassmannian is a $\ZZ_2$-space with fixed subspace
$$
G_m(\infty \oplus \infty L )^{\ZZ_2}=
\bigsqcup_{r=0}^m G_r(\infty )
\times G_{m-r} (\infty L),
$$
and $\eta_m$ restricts on the $r$-component
to the direct sum of $\eta_m^{\ZZ_2}=\eta_r$ and $(\eta_m)_{\ZZ_2}=\eta_{m-r}$.
Here
$$
G_r(\infty )=
\dirlim_{i\geqslant 0} G_r(\Rr^i ) =B\O (r)
$$
and
$$
G_{m-r}(\infty L)=
\dirlim_{j\geqslant 0} G_{m-r}(\Rr^j\otimes L) \, .
$$

\begin{definition} \label{B-iota}
We have inclusion maps
$$
\iota_m : G_m(\infty \oplus\infty L) \to G_{m+1}(\infty \oplus\infty L)
$$
mapping an $m$-dimensional subspace $E\subseteq \Rr^i\oplus (L\otimes\Rr^j)$ to
the $(m+1)$-dimensional subspace $\Rr\oplus E$ of
$\Rr^{i+1}\oplus (L\otimes\Rr^j)$.
So $\iota_m^*(\eta_{m+1})=\Rr\oplus\eta_m$.
\end{definition}

The infinite Grassmannian $G_m(\infty\oplus\infty L)$ is a classifying space
for $\ZZ_2$-vector bundles of dimension $m$. To be precise, given
an $m$-dimensional $\ZZ_2$-bundle $\zeta$ over a compact $\ZZ_2$-ENR
$K$, there is a map
$$
g: K\to G_m(\infty\oplus \infty L)
$$
and a bundle
isomorphism $a: \zeta \to f^*\eta_m$; and any two such $(g,a)$ are
homotopic.
(A pair $(g,a)$ will be given by a  bundle monomorphism
$\zeta \into K\times V$ for some $\ZZ_2$-module $V$.)

\begin{definition}
We define the {\it unoriented $m$-dimensional bordism group of $Y$}, $\NN_m^{\ZZ_2}(Y)$,
to be the bordism group of
closed $\ZZ_2$-manifolds $M$ of dimension $m$ equipped with a $\ZZ_2$-map $f:M\to Y$.
\end{definition}
\begin{definition}
We construct a homomorphism
$$
\nn : \NN_m^{\ZZ_2}(Y) \to \dirlim_{n\geqslant m}
\Omega^{\ZZ_2,\res}_0(Y\times G_n(\infty \oplus\infty L );\, \Rr^{n-m},\eta_n),
$$
the direct limit being formed with respect to the maps $\iota_n$.

Consider a $\ZZ_2$-map $f:M\to Y$.
We have a $\ZZ_2$-map
$$
g:M\to G_m(\infty \oplus\infty L)
$$
and a bundle isomorphism $a: \tau M\to g^*\eta_m$
classifying the tangent bundle.
The pair $(f\times g,a)$ defines an element of
$\Omega^{\ZZ_2 ,\res}_0(Y\times G_m(\infty\oplus\infty L);\,0,\eta_m)$
which depends only on $f:M\to Y$.

Indeed, a homotopy from $(g,a)$ to $(g',a')$ gives a cobordism
$h:W=M\times [0,1]\to G_m(\infty\oplus\infty L)$
and an isomorphism $b:\tau W\iso \Rr\oplus h^*\eta_m$
restricting on $\partial W$ to $(g,a)\sqcup (g',a')$.
\end{definition}
The homomorphism $\nn$ is easily seen to be an isomorphism.
\begin{lemma}
The transformation
$$
\nn : \NN_m^{\ZZ_2}(Y) \to \dirlim_{n\geqslant m}
\Omega^{\ZZ_2,\res}_0(Y\times G_n(\infty \oplus\infty L );\, \Rr^{n-m},\eta_n)
$$
is an isomorphism.
\end{lemma}
\begin{proof}
A class in the restricted bordism group is
specified by an $m$-manifold $M$, a map $f:M\to Y$, a map $g_n: M\to G_n(\infty
\oplus\infty L)$ and a restricted stable isomorphism
$\tau M\oplus\Rr^{n-m}\iso g_n^*\eta_n$.
The restricted stable isomorphism can be represented by a bundle
isomorphism $\tau M\oplus \Rr^{n-m}\oplus\Rr^i\iso g_n^*\eta_n\oplus\Rr^i$
for some $i\geqslant 0$.
Mapping from $G_n(\infty\oplus \infty L)$ to $G_{n+i}(\infty\oplus\infty L)$
we get a map $g_{n+i}: M\to G_{n+i}(\infty \oplus\infty L)$ and a bundle
isomorphism $a_{n+i}: \tau M\oplus \Rr^{n-m+i}\iso g_{n+i}^*\eta_{n+i}$.
This will be homotopic to the stabilization of any classifying pair
$(g,a)$ for $M$.
\hfill\qed\end{proof}
Thus we have obtained a homotopy-theoretic description of $\NN_*^G$.
\begin{theorem}\label{B-N}
There is a natural isomorphism
$$
\NN_m^{\ZZ_2}(Y) \iso \dirlim_{n\geqslant m}
\omega_m^{\ZZ_2}(Y \times G_n(\infty \oplus\infty L );\, \Rr^n-\eta_n).
$$
\end{theorem}
\begin{proof}
The assertion now follows at once from Proposition \ref{B-equ}.
But it is useful to spell out the details.

The group
$\omega_m^{\ZZ_2}(Y \times G_n(\infty \oplus\infty L );\, \Rr^n-\eta_n)$
in the statement
is to be understood as the direct limit
$$
\dirlim_{V} \omega_m^{\ZZ_2}(Y \times G_n(V);\, \Rr^n-\eta_n)
$$
over $\ZZ_2$-modules $V$
and
$$
\omega_m^{\ZZ_2}(Y \times G_n(V);\, \Rr^n-\eta_n)
=\omega_m^{\ZZ_2}(Y \times G_n(V);\, \Rr^n\oplus\eta_n^\perp-V),
$$
where $\eta_n^\perp$ is the orthogonal complement of $\eta_n$ in
the trivial bundle $V$ over $G_n(V)$.

If we take $V$ to include a trivial summand $\Rr^n$.
$V=\Rr^n\oplus V'$, we have
$$
\omega_m^{\ZZ_2}(Y \times G_n(\Rr^n\oplus V');\, \Rr^n-\eta_n)
=\omega_m^{\ZZ_2}(Y \times G_n(\Rr^n\oplus V');\, \eta_n^\perp-V').
$$
For $n\geqslant m$, the Pontrjagin-Thom homomorphism
$$
\Omega_m^{\ZZ_2,\res}(Y \times G_n(\Rr^n \oplus V');\, \Rr^n,\eta_n)
\to
\omega_m^{\ZZ_2}(Y \times G_n(\Rr^n \oplus V');\, \Rr^n-\eta_n)
$$
is an isomorphism.
\hfill\qed\end{proof}

The classical non-equivariant unoriented bordism group
$\NN_m(Y_0)$ of an ANR $Y_0$ can be described similarly
as a framed bordism group and so as a stable homotopy group.
\begin{proposition}
There are natural isomorphisms
$$
\NN_m(Y_0) \Rarr{\iso}{}
\dirlim_{n\geqslant m} \Omega_0(Y\times G_n(\infty );\, \Rr^{n-m}, \eta_n)
\Rarr{\iso}{}
\dirlim_{n \geqslant m}\omega_m(Y_0\times G_n(\infty );\, \Rr^n-\eta_n)
$$
for any ANR $Y_0$.
\qed
\end{proposition}
\begin{definition}
Write $\NN_m^{\ZZ_2 ,\free}(Y)$ for the bordism group of
maps $f : M\to Y$ from a free $\ZZ_2$-manifold $M$.
(So $(M,f)$ is a boundary if there is a free manifold $W$
with boundary $\partial W=M$ and a $\ZZ_2$-map $g :W\to Y$
extending $f$.)
\end{definition}
The free bordism group is described classically as a non-equivariant
bordism group.
\begin{lemma}
There are natural equivalences
$$
\NN_m^{\ZZ_2 ,\free}(Y) \Rarr{\iso}{}
\NN_m^{\ZZ_2}(E\ZZ_2 \times Y) \Rarr{\iso}{}
\NN_m (E\ZZ_2 \times_{\ZZ_2} Y)
$$
for any $\ZZ_2$-ANR $Y$.
\qed
\end{lemma}
We can now recognize the Conner-Floyd exact sequence in
unoriented $\ZZ_2$-bordism as a translation of the localization
exact sequence in $\ZZ_2$-equivariant stable homotopy.
\begin{proposition}
The classical Conner-Floyd \cite{cf} sequence
$$
\cdots\to \NN_m^{\ZZ_2 ,\free} (Y)
\to
\NN_m^{\ZZ_2}(Y)
\to
\bigoplus_{s=0}^m \NN_{m-s}(Y^{\ZZ_2}\times
G_s(\infty )))
\to\cdots
$$
corresponds to the direct limit of the stable homotopy localization sequences
$$
\cdots\to
\omega_m^{\ZZ_2}(E\ZZ_2\times Y\times G_n(\infty\oplus\infty L));
\, \Rr^n -\eta_n)\Rarr{\gamma}{}
$$
$$
\omega_m^{\ZZ_2}(Y\times G_n(\infty\oplus\infty L));
\, \Rr^n -\eta_n)\Rarr{\rho}{}
\omega_m(Y^{\ZZ_2}\times G_n(\infty\oplus\infty L)^{\ZZ_2});
\, \Rr^n -\eta_n^{\ZZ_2}) \to \cdots \, .
$$
\end{proposition}
\begin{proof}
The identification on fixed points is made by the isomorphism
$$\begin{array}{l}
\dirlim_n \,\,
\omega_m(Y^{\ZZ_2}\times G_n(\infty\oplus\infty L)^{\ZZ_2};
\, \Rr^n -\eta_n^{\ZZ_2}) \\
\noalign{\medskip}
=\dirlim_n \bigoplus_{0\leq s\leq n}
\omega_{m-s} (Y^{\ZZ_2}\times G_{n-s}(\infty )\times
G_s(\infty L);
\Rr^{n-s}-\eta_{n-s}) \\
\noalign{\medskip}
=\bigoplus_{s=0}^m \NN_{m-s}(Y^{\ZZ_2}\times G_s(\infty L)).
\end{array}
$$
We can identify the free term
$$
\omega_m^{\ZZ_2}(E\ZZ_2\times Y\times G_n(\infty);
\, \Rr^n -\eta_n)
\Rarr{\iso}{}
\omega_m^{\ZZ_2}(E\ZZ_2\times Y\times G_n(\infty\oplus\infty L);
\, \Rr^n -\eta_n)
$$
with
$$
\omega_m((E\ZZ_2\times_{\ZZ_2} Y)\times G_n(\infty);
\, \Rr^n -\eta_n)
$$
and so with $\NN_m(E\ZZ_2\times_{\ZZ_2} Y)$.
\hfill\qed\end{proof}

The geometric equivariant bordism group $\NN_m^{\ZZ_2}(Y)$
has been expressed as the direct limit
$$
\dirlim_n\dirlim_{V'}
\omega_m^{\ZZ_2}(Y \times G_n(\Rr^n\oplus V');\, \Rr^n-\eta_n),
$$
where $V'$ is a $\ZZ_2$-module and
$G_n(\Rr^n\oplus V')$ is included in $G_n(\Rr^n\oplus U'\oplus V')$
through the inclusion $\Rr^n\oplus V' \into \Rr^n\oplus U'\oplus V'$
and $G_n(\Rr^n\oplus V')$ is included in $G_{n+1}(\Rr^{n+1}\oplus V')$
by the direct sum with $\Rr \subseteq \Rr\oplus\Rr^n=\Rr^{n+1}$.

\begin{definition} \label{B-homotopical}
{\it The homotopical equivariant bordism group} $N_m^{\ZZ_2}(Y)$ is
defined as the direct limit
$$
N^{\ZZ_2}_m(Y)=
\dirlim_{V,V'}
\omega_m^{\ZZ_2}(Y \times G_n(V\oplus V');\,
V-\eta_n),
$$
where $V$ and $V'$ are $\ZZ_2$-modules of dimension $n$ and $n'$
respectively and $G_n(V\oplus V')$ is included in
$G_{k+n}(U\oplus V\oplus U'\oplus V')$, for
$\ZZ_2$-modules $U$ and $U'$ of dimension $k$ and $k'$, by the direct sum with
$U$.
There is a natural transformation from the geometric to the homotopical
theory:
$$ \NN^{\ZZ_2}_m(Y) \to N^{\ZZ_2}_m(Y)\, .
$$
\end{definition}
In the non-equivariant theory we have, for an ANR $Y_0$,
$$
\NN_m(Y_0) =\dirlim_{V,V'} \omega_m (Y_0\times G_n(V_0\oplus V_0');\,
V_0-\eta_n),
$$
where limit runs over vector spaces $V_0$ and $V_0'$ of dimension $n$
and $n'$; it, therefore, makes sense to write
$$
\NN_m(Y_0) = N_m(Y_0)\, .
$$

\begin{proposition}
For the homotopical theory $N_*^{\ZZ_2}$ there is a Conner-Floyd long
exact sequence
$$
\cdots\to N_m(E\ZZ_2\times_{\ZZ_2}Y)
\to
N_m^{\ZZ_2}(Y)
\to
\bigoplus_{s=-\infty}^m N_{m-s}(Y^{\ZZ_2}\times B\O (\infty ))
\to\cdots
$$
\end{proposition}
\begin{proof}
This is the direct limit of
the $\omega_*^{\ZZ_2}$ localization sequences for the spaces
$Y\times G_{n'}(V\oplus V')$ with coefficients $\eta_{n'}-V'$.
The limit over all $\ZZ_2$-modules $V$ produces the summands in all
dimensions $m-s \geqslant 0$.

For further discussion of the sequence we refer to Sinha \cite{sinha}.
\hfill\qed\end{proof}
\begin{remark}
We have defined $N_m(Y_0)$ as the direct limit
$$
\dirlim_{V_0,V'_0} \omega_m(Y_0\times G_n(V_0\oplus V'_0);\, V_0-\eta_n)
$$
over finite dimensional $\Rr$-vector spaces $V_0$ and $V'_0$ of dimension
$n$ and $n'$ respectively.
In the direct system
the Grassmannian $G_n(V_0\oplus V'_0)$ is included in
$G_{k+n}(U_0\oplus V_0\oplus U'_0\oplus V'_0)$,
where $U_0$ and $U'_0$ are vector
spaces of dimension $k$ and $k'$, by taking the direct sum of
an $n$-dimensional subspace of $V_0 \oplus V'_0$ with $U_0$ to
get a $(k+n)$-dimensional subspace of $U_0\oplus V_0\oplus V'_0\subseteq
U_0\oplus V_0\oplus U'_0\oplus V'_0$.

In the bordism interpretation, this description classifies the
stable tangent bundle. For the more conventional interpretation
involving the stable normal bundle and the eponymous Thom spectrum,
we can identify these groups with
$$
\dirlim_{V_0,V'} \omega_m(Y_0\times G_{n'}(V_0\oplus V'_0);\, \eta_{n'}-V'_0) .
$$
The orthogonal complement identifies $G_n(V_0\oplus V'_0)$ with
$G_{n'}(V_0\oplus V'_0)$ and then $\eta_n\oplus\eta_{n'} = V_0\oplus V'_0$.
The corresponding inclusion of $G_{n'}(V_0\oplus V'_0)$ as a subspace of
$G_{k'+n'}(U_0\oplus V_0\oplus U'_0\oplus V'_0)$ is given by the direct sum
with $U'_0$.
\end{remark}

There is a similar complementary interpretation of the equivariant groups
$N^{\ZZ_2}_m(Y)$ in the notation of Definition \ref{B-homotopical}.
\begin{proposition} \label{B-normal}
The homotopical bordism group $N^{\ZZ_2}_m(Y)$ can be expressed as
the direct limit
$$
N^{\ZZ_2}_m(Y)=
\dirlim_{V,V'}
\omega_m^{\ZZ_2}(Y \times G_{n'}(V\oplus V');\, \eta_{n'}-V'),
$$
over pairs of $\ZZ_2$-modules $V$ and $V'$ of dimension $n$ and $n'$,
where the inclusion of $G_{n'}(V\oplus V')$ into
$G_{k'+n'}(U\oplus V\oplus U'\oplus V')$ for $\ZZ_2$-modules
$U$ and $U'$ of dimension $k$ and $k'$ is given by the direct sum with $U'$.
\qed
\end{proposition}
\begin{remark}
The standard definition of $N_m(Y_0)$ in terms of Thom spectra
is as the direct limit
$$
\dirlim_{V_0,V'_0}
[(\Rr^m\oplus V'_0)^\cpt;\, (Y_0)\bpt\wedge
\Thom(\eta_{n'},G_{n'}(V_0\oplus V'_0))]
$$
over vector spaces $V_0$ and $V'_0$.
But stabilization
$$
[(\Rr^m\oplus V'_0)^\cpt;\, (Y_0)\bpt\wedge
\Thom (\eta_{n'},G_{n'}(V_0\oplus V'_0))]
\to \omega_m(Y_0\times G_{n'}(V_0\oplus V'_0); \, \eta_{n'}-V'_0)
$$
is an isomorphism if $n'> m+1$.
\end{remark}
One has a similar description of the $\ZZ_2$-equivariant theory,
for what are rather formal reasons.
\begin{proposition}
There is an isomorphism
$$
\alpha :\dirlim_{V,V'}
[(\Rr^m\oplus V')^\cpt ;\, Y\bpt\wedge
\Thom(\eta_{n'}, G_{n'}(V\oplus V'))]_{\ZZ_2}
\to N^{\ZZ_2}_m(Y).
$$
\end{proposition}
\begin{proof}
Using Proposition \ref{B-normal} and
the definition of an equivariant stable map, we can write
$$
N^{\ZZ_2}_m(Y) =
\dirlim_{U', V, V'}
[(\Rr^m\oplus U'\oplus V')^\cpt ;\, Y\bpt\wedge
\Thom ( U'\oplus \eta_{n'}, G_{n'}(V\oplus V'))]_{\ZZ_2}
$$
as a direct limit over $\ZZ_2$-modules $U'$, $V$ and $V'$.
The stabilization map $\alpha$ is defined by including the $U'=0$ terms.
In the opposite direction, we get a map
$$
\beta :N^{\ZZ_2}_m(Y) \to
\dirlim_{U, V,V'}
[(\Rr^m\oplus U'\oplus V')^\cpt ;\, Y\bpt\wedge
\Thom ( \eta_{k'+n'}, G_{k'+n'}(V\oplus U'\oplus V'))]_{\ZZ_2},
$$
where $k'=\dim U'$ and $G_{n'}(V\oplus V')$ is embedded
into $G_{k'+n'}(V\oplus U'\oplus V')$ by taking the direct sum with $U'$
so that $\eta_{k'+n'}$ restricts to $U'\oplus\eta_{n'}$.

The composition $\beta\circ\alpha$ is clearly the identity.
To understand $\alpha\circ\beta$, we observe that the map
$$
\begin{matrix}
[(\Rr^m\oplus U'\oplus V')^\cpt ;\, Y\bpt \wedge
\Thom(U'\oplus \eta_{n'}, G_{n'}(V\oplus V'))]_{\ZZ_2} \\
\Darr{}{} \\
[(\Rr^m\oplus U'\oplus V')^\cpt ;\, Y\bpt\wedge
\Thom ( \eta_{k'+n'}, G_{k'+n'}(V\oplus U'\oplus V'))]_{\ZZ_2}
\end{matrix}
$$
at the finite level lifts the homomorphism
$$
\begin{matrix}
\omega_m^{\ZZ_2}(Y \times G_{n'}(V\oplus V');\, \eta_{n'}- V') \\
\Darr{}{} \\
\omega_m^{\ZZ_2}(Y \times G_{k'+n'}(V\oplus (U'\oplus V'));\, \eta_{k'+n'}-
(U'\oplus V'))
\end{matrix}
$$
in the direct system describing $N^{\ZZ_2}_m(Y)$.
Hence $\alpha\circ\beta$ is the identity.
\hfill\qed\end{proof}

A refinement of the bordism groups in which restrictions are placed
on the fixed-point sets was considered by Kosniowski \cite{ck1} and Waner \cite{waner}.
Suppose that $\AA \subseteq \RO (\ZZ_2)$ is a finite set of virtual
representations of dimension $0$.
Thus $\AA$ is a finite set of multiples of $[L]-[\Rr]$.
\begin{definition}
We define $\NN_m^{\ZZ_2, \AA} (Y)$
to be the set of bordism classes of $\ZZ_2$-maps $f : M\to Y$
from an $m$-dimensional $\ZZ_2$-manifold $M$ such that, for each $x\in M^{\ZZ_2}$,
$[\tau_xM]-m \in \RO (\ZZ_2 )$ lies in $\AA$.
\end{definition}

\begin{definition}
For a $\ZZ_2$-module $V$ and integer $n\geqslant 1$, let
$G_n^\AA (V)$ be the subset of $G_n(V)$ consisting of those
$n$-dimensional subspaces $E\subseteq V$ such that
either (i) $E\notin G_n(V)^{\ZZ_2}$ or
(ii) $E\in G_n(V)^{\ZZ_2}$ and $[E]-n\in \AA \subseteq\RO(\ZZ_2 )$.
\end{definition}

\begin{lemma}
The subset $G_n^\AA (V)$ is an open $\ZZ_2$-subspace of $G_n(V)$
and so is a $\ZZ_2$-ENR.
\qed
\end{lemma}
\begin{lemma}\label{B-AA}
Suppose that $M$ is a closed $\ZZ_2$-manifold such that,
for each $x\in M^{\ZZ_2}$,
$[\tau_xM]-m \in \RO (\ZZ_2 )$ lies in $\AA$.
Then there is a $\ZZ_2$-module $V$ and a $\ZZ_2$-monomorphism
$\tau M \into M\times V$ such that for each $x\in M$ the $m$-dimensional
subspace $\tau_x M\subseteq V$ lies in $G_m^\AA (V)$.
\end{lemma}
\begin{proof}
For each $x\in M^{\ZZ_2}$, let $V_x=\tau_xM$ and let $i_x: \tau_x M\to V_x$
be the identity map; for a point $x$ in the complement of $M^{\ZZ_2}$,
let $V_x=\tau_xM\oplus \tau_y M$, where $\{ x,y\}$ is the $\ZZ_2$-orbit
of $x$ and let $i_x :\tau_x M\into V_x$ be the inclusion.
Each inclusion $i_x$ extends to a $\ZZ_2$-map $\tau M \to M\times V_x$,
which will
be a monomorphism on some neighbourhood of the orbit of $x$ and give
a map of some smaller neighbourhood $U_x$ to $G_m^\AA (V_x)$
by the openness of $G_m^\AA (V_x)$ in $G_n(V_x)$.
Choose a finite subset $Q\subseteq X$ so that the  open sets
$U_x,\, x\in Q$, cover $M$, and
take the sum of the linear maps $\tau M \to M\times V_x$, $x\in Q$, to get
a map $\tau M \into M\times V$, where $V=\bigoplus_{x\in Q} V_x$.
\hfill\qed\end{proof}
\begin{proposition}
There is a natural isomorphism
$$
\nn :\NN_m^{\ZZ_2 , \AA} (Y) \to
\dirlim_{n\geqslant m} \Omega_0^{\ZZ_2 ,\res}(Y\times
G_n^\AA (\infty \oplus \infty L);
\, \Rr^{n-m},\eta_n).
$$
\end{proposition}
\begin{proof}
Consider a representative $f: M\to Y$ of a class in
$\NN_m^{\ZZ_2,\AA } (Y)$.
Choose a $\ZZ_2$-monomorphism $\tau M\into M\times V$ as in Lemma \ref{B-AA} above.
Its classifying map determines an element of the group
$\Omega^{\ZZ_2,\res}_0(Y\times G_m^\AA (V); 0,\eta_m)$ and so of
$\Omega^{\ZZ_2,\res}_0(Y\times G_m^\AA (\infty \oplus \infty L);
0,\eta_m)$.
This construction produces a well-defined homomorphism
$$
\NN_m^{\ZZ_2,\AA} (Y) \to
\Omega_0^{\ZZ_2,\res}(Y\times G_n^\AA  (\infty \oplus \infty L);
\, 0,\eta_m)
$$
and so to the direct limit.

Conversely, given a representative of consisting of maps $f:M\to Y$ and
$g: M\to G_n^\AA (V)$ and a bundle isomorphism
$\tau M\oplus \Rr^{n-m}\oplus \Rr^i\iso g^*\eta_n\oplus \Rr^i$, of
a class in $\Omega_0^{\ZZ_2,\res}(Y\times G_n^\AA (V);
\, \Rr^{n-m}, \eta_m)$,
the component $f: M\to Y$ gives an element of $\NN_m^{\ZZ_2, \AA} (Y)$.
\hfill\qed\end{proof}
\begin{corollary}
There is a natural isomorphism
$$
\NN_m^{\ZZ_2 , \AA} (Y) \iso \dirlim_{n\geqslant m}
\omega_m^{\ZZ_2} (Y\times G_n^\AA (\infty
\oplus \infty L); \Rr^n - \eta_n).
$$
\end{corollary}
\begin{proof}
This follows at once from Proposition \ref{B-equ}.
\hfill\qed\end{proof}

The naive version of oriented $G$-bordism,
which is appropriate for many geometric applications, can be treated
in the same way.
(More sophisticated versions may be found in Costenoble and Waner \cite{costewaner}.)
\begin{definition}
An {\it orientation} of
an $m$-dimensional $\ZZ_2$-manifold $M$ is a $\ZZ_2$-equivariant trivialization $S(\Lambda^m\tau M) \to M\times S(\Rr)$ of
the orientation bundle of $M$.
In other words, it is an orientation in the non-equivariant sense
that is fixed under the action of the group $\ZZ_2$.
\end{definition}
\begin{definition}
We define the {\it oriented $m$-dimensional bordism group of $Y$}, $\OO_m^{\ZZ_2}(Y)$,
to be the bordism group of oriented
closed $\ZZ_2$-manifolds $M$ of dimension $m$ equipped with a $\ZZ_2$-map $f:M\to Y$.
\end{definition}
For a $\ZZ_2$-module $V$, we write $G^\ori_m(V)$ for the
Grassmann manifold of oriented $m$-dimensional subspaces of $V$. The double cover $G^\ori_m(V) \to G_m(V)$ is the
sphere bundle $S(\Lambda^m\eta_m)$ over $G_m(V)$.
\begin{proposition}
There are natural equivalences
$$
\OO^{\ZZ_2}_m(Y) \Rarr{\iso}{\nn}
 \dirlim_{n\geqslant m}
\Omega^{\ZZ_2,\res}_0(Y\times G_n^\ori(\infty \oplus\infty L );\, \Rr^{n-m},\eta_n)
$$
$$
\qquad
\Rarr{\iso}{}
\dirlim_{n\geqslant m}
\omega^{\ZZ_2}_m(Y\times G_n^\ori (\infty \oplus\infty L );\, \Rr^n-\eta_n),
$$
for any $\ZZ_2$-ANR $Y$.
\qed
\end{proposition}

\chapter{{\bf The geometric Hopf invariant and double points (2010)}}
\label{appendix3}

This, with a few small differences, including the correction of a misprint in Proposition 3.1, is the text of our joint paper
{\it The geometric Hopf invariant  and double points}, Journal of Fixed Point Theory and Applications 7, 325--350 (2010). It is included here by kind permission of Springer Verlag.

\newpage

{
%
\newtheorem{thm}{Theorem}[section]
\newtheorem{cor}[thm]{Corollary}
\newtheorem{lem}[thm]{Lemma}
\newtheorem{prop}[thm]{Proposition}
\newtheorem{defn}[thm]{Definition}
\newtheorem{rem}[thm]{Remark}
\newtheorem{ex}[thm]{Example}
\long\def\Thm#1{\begin{thm} #1 \end{thm}}
\long\def\Cor#1{\begin{cor} #1 \end{cor}}
\long\def\Lem#1{\begin{lem} #1 \end{lem}}
\long\def\Prop#1{\begin{prop} #1 \end{prop}}
\long\def\Def#1{\begin{defn} \rm #1 \end{defn}}
\long\def\Rem#1{\begin{rem} \rm #1 \end{rem}}
\long\def\Ex#1{\begin{ex} \rm #1 \end{ex}}
\def\Sect{\section}
\def\imm{\looparrowright}
\def\Rarr#1#2{\xrightarrow[#2]{#1}}
\long\def\Ref#1#2#3#4#5#6{
\bibitem{#1}
{\rm #2,}
\textit{#3.}
{\rm #4}
\textbf{#5}
{\rm #6.}
}
\long\def\Refb#1#2#3#4{
\bibitem{#1}
{\rm #2,}
\textit{#3.}
#4.
}
%
\def\ZZ{{\Zz /2}}
\def\CC{{\mathcal C}}
\def\CCC{\mathfrak{L}(f)}
\def\Dp{\mathfrak{D}(f)}
\def\Dpone{\mathfrak{D}(f')}
\def\barDp{\overline{\mathfrak{D}}(f)}
\def\barDpone{\overline{\mathfrak{D}}(f')}
\def\hDp{\hbox{\rm h-${\mathfrak{D}(f)}$}}
\def\hI{\hbox{\rm h-${\mathfrak{I}(f^{(1)},f^{(2)})}$}}
\def\barhDp{\hbox{\rm h-$\overline{\mathfrak{D}}(f)$}}
\def\barD{Z}
\def\CD{\hbox{\rm h-$\widehat{\mathfrak{D}}(f)$}}
\def\CDa{\hbox{\rm h-$\widehat{\mathfrak{D}}(f^{(1)})$}}
\def\CDb{\hbox{\rm h-$\widehat{\mathfrak{D}}(f^{(2)})$}}
\def\CDone{\hbox{\rm h-$\widehat{\mathfrak{D}}(f')$}}
\def\CDg{\hbox{\rm h-$\widehat{\mathfrak{D}}(g)$}}
\def\EE{{\mathcal{D}}}
\def\O{{\rm O}}
\def\Qq{{\mathbb Q}}
\def\Zz{{\mathbb Z}}
\def\Nn{{\mathbb N}}
\def\Rr{{\mathbb R}}
\def\Cc{{\mathbb C}}
\def\Tt{{\mathbb T}}
\def\SU{{\rm SU}}
\def\Aut{{\rm Aut}}
\def\Hom{{\rm Hom}}
\def\im{\hbox{\rm im}}
\def\into{\hookrightarrow}
\def\imm{\looparrowright} 
\def\iso{\cong}
\def\cl#1{\overline{#1}} 
\def\union{\sqcup} 
\def\comp{\mathbin{\mathchoice
{\circ}
{{\scriptstyle\circ}}
{{\scriptscriptstyle\circ}}
{{\scriptscriptstyle\circ}}
}}
\def\cpt{+} 
\def\st{\mid} 
\def\Null{{\rm Null}} 
\def\hNull{\hbox{\rm h-Null}} 
\def\hInd{\hbox{\rm h-$\gamma$}} 
\def\map{{\rm map}}
\def\geq{\geqslant} \def\leq{\leqslant} 
%
%
%

\bigskip

\bigskip

\Sect*{Abstract}

The geometric Hopf invariant of a stable map $F$
is a stable $\ZZ$-equivariant map $h(F)$ such that the stable
$\ZZ$-equivariant homotopy class of $h(F)$ is the primary obstruction
to $F$ being homotopic to an unstable map. In this paper we express
the geometric Hopf invariant of the Umkehr map $F$ of an immersion
$f:M^m \imm N^n$ in terms of the double point set of $f$.
We interpret
the Smale-Hirsch-Haefliger regular homotopy classification of immersions
$f$ in the metastable dimension range $3m<2n-1$ (when a generic $f$ has no
triple points) in terms of the geometric Hopf invariant.



%

\Sect*{Introduction}

The original Hopf invariant $H(F) \in \Zz$ of a map $F:S^3 \to S^2$
was interpreted by
Steenrod as the evaluation of the cup product in the mapping cone $C(F)$.
The mod 2 Hopf invariant $H_2(F) \in \ZZ$ of a map $F:S^j \to S^k$ was then
defined using the functional Steenrod squares of $f$.
The geometric Hopf invariant of a stable map $F:\Sigma^{\infty}X \to
\Sigma^{\infty}Y$ is the stable $\ZZ$-equivariant map
$$
h(F)~=~(F\wedge F)\Delta_X-\Delta_YF~:~\Sigma^\infty X\to
\Sigma^\infty(Y\wedge Y)
$$
measuring the failure of $F$ to preserve the diagonal maps of $X$ and
$Y$, with $\ZZ$ acting by the transposition $T:(y_1,y_2) \mapsto (y_2,y_1)$
on $Y \wedge Y$. Thus $h(F)$ is a homotopy-theoretic generalization of the
functional Steenrod squares. The stable homotopy
class of $h(F)$ is the primary obstruction to $F$ being homotopic to an
unstable map.

Given an immersion $f:M^m \imm N^n$ with normal bundle $\nu (f)$ we
express the geometric Hopf invariant $h(F)$ of the Umkehr map
$F:\Sigma^{\infty}N^+\to \Sigma^{\infty}M^{\nu(f)}$ in terms of the
double point set of $f$.  There are many antecedents for this
expression! The stable homotopy class of $h(F)$ depends only on the
regular homotopy class of $f$.  If $f$ is regular homotopic to an
embedding then $h(F)$ is stably null-homotopic.  We interpret the
Smale-Hirsch-Haefliger regular homotopy classification of immersions
$f$ in the metastable dimension range $3m<2n-1$ (when a generic $f$ has
no triple points) in terms of the geometric Hopf invariant.

In ``Geometric Hopf invariant and surgery theory''\footnote{this volume!} we  provide a considerably more detailed exposition
of the geometric Hopf invariant $h(F)$ and its applications to
manifolds.  This will include the $\pi_1(N)$-equivariant geometric Hopf
invariant $\widetilde{h}(F)$ needed for a homotopy-theoretic treatment
of the double point invariant $\mu(f)$ of Wall \cite[\S5]{wall2} for a
generic immersion $f:M^m \imm N^{2m}$ which plays such an important
r\^ole in non-simply-connected surgery theory, with $M=S^m$.  When both
$M$ and $N$ are connected and oriented and $f$ induces the trivial map
$\pi_1(M)\to \pi_1(N)$, $\mu (f)$ is an element of the group
$$
\Zz[\pi_1(N)]/\langle g-(-)^mg^{-1}\,\vert\, g \in \pi_1(N)\rangle
$$
and $\mu(f)=0$ if (and, for $m > 2$, only if)  $f$ is regular
homotopic to an embedding, by the Whitney trick for removing double points.
$\widetilde{h}(F)$ induces the quadratic construction $\psi_F$ of \cite{ranicki2} on the chain level.

The present paper is set out as follows.  Section 1 describes briefly
the construction of the geometric Hopf invariant and its fibrewise
generalization.  The double point theorem is stated and proved in
Section 2.  In Section 3, building on work of Dax \cite{dax}, Hatcher
and Quinn \cite{hatcherquinn}, Salomonsen \cite{salomonsen} and Li, Liu and
Zhang \cite{liliuzhang}, we relate the geometric Hopf invariant in a stable
range to Haefliger's obstruction to the existence of a regular homotopy
from an immersion to an embedding.  The papers of Boardman and Steer
\cite{bs} and Koschorke and Sanderson \cite{ks1}
are also relevant. The variation of the geometric Hopf
invariant of an immersion under a (not necessarily regular) homotopy is
computed in Section 4 in terms of the Smale-Hirsch-Haefliger
classification.  In Section 5 we use Whitney's figure-of-eight
immersion \cite{whitney3} to construct, in a metastable range,
immersions close to a given embedding.  Prerequisites for that section,
on the differential-topological classification of vector bundle
monomorphisms, are given in an Appendix.

We shall write the one-point compactification of a locally compact
Hausdorff topological space $X$ as $X^\cpt$.
A subscript `$+$' will be used for the adjunction of a disjoint
basepoint to a space. If $X$ is compact $X^\cpt=X_+$.
For a Euclidean vector bundle $\xi$ over a general space $X$,
we write $D(\xi )$ for the closed
unit disc bundle, $S(\xi )$ for the sphere bundle and $B(\xi )$ for
open unit disc bundle. The Thom space of $\xi$ is written as
$X^\xi$.  To simplify notation, we shall sometimes write
$Y^\xi$, rather than $Y^{p^*\xi}$, for the Thom space of the
pullback $p^*\xi$ by a map $p:Y\to X$, if the map $p$ is clear from
the context. Similarly, we sometimes write $V$, instead of $X\times V$,
for the trivial vector bundle over $X$ with fibre the vector space $V$.

Methods of fibrewise homotopy theory will be used extensively.
Fibrewise constructions, such as the one-point compactification,
the Thom space, or the smash product, over a base $B$ will be indicated
by attaching a subscript `$B$' to the relevant symbol.
We follow  the notation for (fibrewise) stable homotopy adopted in
\cite{crabbcoincidence}.
Consider fibrewise pointed spaces $X\to B$ and $Y\to B$ over an
ENR base $B$.
If $B$ is compact and $A$ is a closed sub-ENR, we write
$$
\omega^0_B\{ X;\, Y\} \text{\quad and\quad}
\omega^0_{(B,A)}\{ X;\, Y\},
$$
respectively,
for the abelian group of stable fibrewise maps $X\to Y$ over $B$
and the relative group
defined in terms of homotopy classes of maps that are zero over the
subspace $A$.
(See, for example, \cite[Part II, Section 3]{crabbjames}.)
We also need to consider fibrewise maps with compact supports.
When $B$ is not necessarily compact we write
$$
{}_c\omega^0_B\{ X;\, Y\}
$$
for the group of fibrewise stable maps
that are zero outside a compact subspace of $B$.
The $\omega^0$-theories are extended,
using the fibrewise suspension $\Sigma_B$ over $B$,
to $\omega^i$-cohomology theories indexed by $i\in\Zz$.
When $Y\to B$ is a trivial bundle $B\times S^i\to B$,
there are natural identifications of the
fibrewise groups with the reduced stable cohomotopy of an
appropriate pointed space:
$$
\omega^0_B\{ X;\, B\times S^i\} =\tilde\omega^i(X/B),
\text{\quad and\quad}
\omega^0_{(B,A)}\{ X;\, B\times S^i\} =\tilde\omega^i(X/(X_A\cup B)),
$$
where $X_A\to A$ denotes the restriction of $X\to B$.

We shall also need $\ZZ$-equivariant stable homotopy theory, which we
indicate by a prefix as ${}^\ZZ\omega^*$. So, for example, the
equivariant stable cohomotopy of a point, ${}^\ZZ\omega^0(*)$, is
the direct limit over all finite-dimensional $\Rr$-vector spaces $V$ and $W$
of the homotopy classes of pointed $\ZZ$-maps
$(V\oplus LW)^\cpt\to (V\oplus LW)^\cpt$.
Here, and throughout the paper, we write $L$ for  the non-trivial
$1$-dimensional representation $\Rr$ of $\ZZ$ with
the involution $-1$, and,
for a finite-dimensional real vector space $W$, abbreviate the tensor product
$L\otimes W$ to $LW$.

\smallskip

We thank Mark Grant for pointing out the relevance of Hatcher and Quinn \cite{hatcherquinn}.
\Sect{A review of the geometric Hopf invariant}
Let $X$ and $Y$ be pointed topological spaces, and let $V$ be a finite
dimensional Euclidean space. It is convenient to assume that $X$ is
a compact ENR and that $Y$ is an ANR.

Let the generator $T \in \ZZ$ act on $X \wedge X$ by transposition
$$
T~:~X \wedge X \to X \wedge X~;~(x,y)\mapsto (y,x)\, .
$$
The diagonal map
$$
\Delta_X~:~X\to X \wedge X~;~x \mapsto (x,x)
$$
is $\ZZ$-equivariant. The diagonal map $\Delta_{V^\cpt }$
extends to a $\ZZ$-equivariant homeomorphism
$$
\kappa_V~:~LV^{\cpt} \wedge V^{\cpt} \to V^{\cpt} \wedge V^{\cpt}~;~
(u,v) \mapsto (u+v,-u+v)\, .
$$
The $\ZZ$-action $u \mapsto -u$ on $LV$ has fixed point $\{0\}$; the
$\ZZ$-action on the unit sphere
$$
S(LV)~=~\{u \in V\,\vert\, \Vert u \Vert = 1\}
$$
is free, with $\Vert~\Vert$ any inner product on $V$.

The geometric Hopf invariant of a map
$F:V^\cpt \wedge X \to V^\cpt \wedge Y$
measures the difference $(F\wedge F)\Delta_X - \Delta_YF$, given that
$(F\wedge F)\Delta_{V^\cpt \wedge X}=\Delta_{V^\cpt \wedge Y}F$.
The diagram of $\ZZ$-equivariant maps
$$
\xymatrix@R+20pt@C+40pt{
LV^\cpt \wedge V^\cpt \wedge X \ar[r]^-{\displaystyle{1 \wedge \Delta_X}}
\ar[d]_-{\displaystyle{1 \wedge F}} &
LV^\cpt \wedge V^\cpt \wedge X \wedge X
\ar[d]^-{\displaystyle{G}}\\
LV^\cpt \wedge V^\cpt \wedge Y
\ar[r]^-{\displaystyle{1 \wedge \Delta_Y}}
&LV^\cpt \wedge V^\cpt \wedge Y \wedge Y}
$$
does not commute in general, with
$$
\begin{array}{l}
G~=~(\kappa^{-1}_V \wedge 1)(F \wedge F)(\kappa_V \wedge 1)~:\\[1ex]
\hskip50pt
LV^\cpt \wedge V^\cpt \wedge X \wedge X
\to LV^\cpt \wedge V^\cpt \wedge Y \wedge Y~;\\[1ex]
\hskip75pt (u,v,x_1,x_2) \mapsto ((w_1-w_2)/2,(w_1+w_2)/2,y_1,y_2)\\[1ex]
\hskip75pt (F(u+v,x_1)=(w_1,y_1),F(-u+v,x_2)=(w_2,y_2))\, .
\end{array}
$$
However, the $\ZZ$-equivariant maps defined by
$$
\begin{array}{l}
p~=~G(1\wedge \Delta_X)~:~LV^\cpt \wedge V^\cpt \wedge X \to
LV^\cpt \wedge V^\cpt \wedge Y \wedge Y~;\\[1ex]
\hskip75pt (u,v,x) \mapsto ((w_1-w_2)/2,(w_1+w_2)/2,y_1,y_2)\\[1ex]
\hskip75pt (F(u+v,x)=(w_1,y_1),F(-u+v,x)=(w_2,y_2))\, ,\\[1ex]
q~=~(1 \wedge \Delta_Y)(1 \wedge F)~:~
LV^\cpt \wedge V^\cpt \wedge X \to
LV^\cpt \wedge V^\cpt \wedge Y \wedge Y~;\\[1ex]
\hskip75pt (u,v,x) \mapsto (u,w,y,y)~~(F(v,x)=(w,y))
\end{array}
$$
agree on $0^+ \wedge V^\cpt \wedge X =V^\cpt \wedge X \subseteq
LV^\cpt \wedge V^\cpt \wedge X$, with
$$
\begin{array}{l}
p\vert~=~q\vert~=~(\kappa^{-1}_V \wedge 1)\Delta_{V^\cpt \wedge Y}
F~=~(\kappa^{-1}_V \wedge 1)(F \wedge F)(\kappa_V \wedge 1)
\Delta_{V^\cpt \wedge X}~:\\[1ex]
\hskip100pt V^\cpt \wedge X
\to LV^\cpt \wedge V^\cpt \wedge Y\wedge Y\, .
\end{array}
$$
We adopt the following terminology: for $t \in [0,1]$, $u \in S(V)$ let
$$[t,u]~=~\dfrac{tu}{1-t} \in V^\cpt~.$$
\Def{(\cite[pp.~306--308]{crabbjames})
The {\it geometric Hopf invariant} of a pointed
map $F:V^\cpt \wedge X \to
V^\cpt \wedge Y$ is the $\ZZ$-equivariant map given by the relative
difference of the $\ZZ$-equivariant maps $p,q$
\index{geometric Hopf invariant!$h_V(F)$}
$$
\begin{array}{ll}
h_V(F)~=~\delta(p,q)~:&\Sigma S(LV)^+ \wedge V^\cpt \wedge X \to
LV^\cpt \wedge V^\cpt \wedge Y\wedge Y~;\\[1ex]
&(t,u,v,x) \mapsto \begin{cases}
q([ 1-2t,u],v,x)&\hbox{if $0 \leq  t \leq  1/2$}\\[1ex]
p([ 2t-1,u],v,x)&\hbox{if $1/2 \leq  t
\leq 1$}
\end{cases}\\[3ex]
&(t \in [0,1],\, u \in S(LV),\, v \in V,\, x \in X)\, .
\end{array}
$$
\qed
}

We are primarily interested in the stable $\ZZ$-equivariant class
of $h_V(F)$, which we write simply as
$$
h_V(F) \in {}^\ZZ\omega^0\{ \Sigma S(LV)_+\wedge X;\,
LV^\cpt\wedge (Y\wedge Y)\}\, .
$$
Using duality, we can rewrite this stable homotopy group in different ways and
thus obtain two other versions of the geometric Hopf invariant as follows.
Smashing $h_V(F)$ with the identity on the sphere $S(LV)_+$
and composing with the duality map
$$
LV^\cpt \to \Sigma S(LV)_+ \to
S(LV)_+\wedge \Sigma S(LV)_+
$$
we get a map
$$
V^\cpt \wedge LV^\cpt \wedge X \to
V^\cpt \wedge LV^\cpt \wedge S(LV)_+ \wedge (Y\wedge Y)
$$
and a second version of the geometric Hopf invariant as an element
$$
h'_V(F)\in{}^\ZZ\omega^0\{ X;\, S(LV)_+\wedge (Y\wedge Y)\}\, .
$$
\Rem{\label{remark1}
The $\ZZ$-equivariant cofibration sequence
$$
S(LV)_+ \to S^0=\{0\}^\cpt \to LV^\cpt
$$
induces an exact sequence of stable $\ZZ$-equivariant homotopy groups
$$
{}^\ZZ\omega^0\{ X;\, S(LV)_+\wedge (Y\wedge Y)\} \to
{}^\ZZ\omega^0\{ X;\, Y\wedge Y\} \to
{}^\ZZ\omega^0\{X ;\, LV^\cpt\wedge Y\wedge Y\}\, .
$$
The stable class $h'_V(F)$ has image
$(F \wedge F)\Delta_X - \Delta_YF$ in
${}^\ZZ\omega^0\{ X;\, Y\wedge Y\}$.
\smallskip

We may also use the Adams isomorphism \cite[Theorem 5.3]{adams3}
$$
{}^\ZZ\omega^0\{ X;\, S(LV)_+\wedge (Y\wedge Y)\}~\cong~
\omega^0\{ X;\, S(LV)_+\wedge_{\ZZ} (Y\wedge Y)\}
$$
to regard $h'_V(F)$ as a non-equivariant class
$$
h''_V(F)\in\omega^0\{ X;\, S(LV)_+\wedge_{\ZZ} (Y\wedge Y)\}\, .
$$
\qed}
\Rem{Instead of applying the Adams isomorphism,
we can use fibrewise techniques.
The unstable map $h_V(F)$ lifts to
an equivariant fibrewise pointed map over $S(LV)$:
$$
S(LV)\times (\Sigma V^\cpt \wedge X) \to
S(LV)\times (V^\cpt\wedge LV^\cpt \wedge  (Y\wedge Y))\, ,
$$
and we may divide by the free $\ZZ$-action to get
a fibrewise pointed map over the real projective space $P(V)=S(LV)/\ZZ$:
$$
P(V)\times (\Sigma V^\cpt \wedge X) \to
(V\oplus (H\otimes V))^\cpt_{P(V)}\wedge_{P(V)}
(S(LV)\times_{\ZZ} (Y\wedge Y))\, ,
$$
where $H=S(LV)\times_{\ZZ} L$ is the Hopf line bundle over $P(V)$
and $\{-\}^+_{P(V)}$ denotes fibrewise one-point
compactification over $P(V)$.
This fibrewise map represents a stable class in the group
$$
\omega^{-1}_{P(V)}\{ P(V)\times X;\, (H\otimes V)^\cpt_{P(V)}
\wedge_{P(V)} (S(LV)\times_{\ZZ}(Y\wedge Y))\}\, ,
$$
which may be identified by fibrewise Poincar\'e-Atiyah duality
(since $\Rr\oplus\tau P(V)= H\otimes V$) with
$$
\omega^0\{ X;\, S(LV)_+\wedge_{\ZZ} (Y\wedge Y)\}\, .
$$
(For the duality theorem, see, for example, \cite[Part II, Section 12]{crabbjames}.)
\qed}
\Rem{By \cite[pp.~60--61]{crabb} the limit over all finite-dimensional
$V$ of the exact sequences in Remark \ref{remark1} is
a split exact sequence
$$
\omega^0\{ X;\, (E\ZZ )_+\wedge_{\ZZ} (Y\wedge Y)\} \Rarr{\delta}{}
{}^\ZZ\omega^0\{ X;\, Y\wedge Y\} \Rarr{\rho}{}\omega^0\{X ;Y\}
$$
with the fixed point map $\rho$ split by
$$
\omega^0\{X ;Y\} \to {}^\ZZ\omega^0\{ X;\,
Y\wedge Y\}~;~F \mapsto \Delta_YF\, .
$$
The stable $\ZZ$-equivariant homotopy class
$$
h(F)~=~(F\wedge F)\Delta_X-\Delta_YF \in {}^\ZZ\omega^0\{ X;\, Y\wedge Y\}
$$
is the image under $\delta$ of
$$
h'(F)\,=\,\varinjlim\limits_V \,  h'_V(F) \in
\omega^0\{ X;\, (E\ZZ )_+\wedge_{\ZZ} (Y\wedge Y)\}\, .
$$
\qed}

We shall need some fundamental properties of the geometric Hopf invariant.
First, it is an obstruction to desuspension.
The Hopf invariant of a suspension is zero, but there is a more
precise result for the geometric Hopf invariant.
\Lem{\label{suspension}
Suppose that $F$ is the suspension $1_{V^\cpt}\wedge F_0$
of a map $F_0 : X\to Y$. Then $h_V(F) = 1_{V^\cpt}\wedge n_V \wedge
(\Delta_Y\circ F_0)$, where $n_V : \Sigma S(LV)_+\to LV^\cpt$ is
an explicit $\ZZ$-equivariantly null-homotopic map.
\qed}
The second property, giving a formula for the Hopf invariant of a composition,
is suggested by the stable form described in Remark \ref{remark1}.
\Prop{\label{composition}{\rm (Composition formula).}
Let $X,\, Y$ and $Z$ be pointed spaces, and
suppose that $F:V^\cpt\wedge X\to V^\cpt\wedge Y$ and $G: V^\cpt\wedge Y\to V^\cpt\wedge Z$
are pointed maps. Then
$$
h_V(G\comp F) = h_V(G)[F] + [G\wedge G]h_V(F)\in
{}^\ZZ\omega^0\{ \Sigma S(LV)_+\wedge X;\, LV^\cpt\wedge (Z\wedge Z)\},
$$
where
$$
[F]\in {}^\ZZ \omega^0\{ X;\, Y\} \text{\quad and \quad}
[G\wedge G]\in {}^\ZZ\omega^0\{ Y\wedge Y;\, Z\wedge Z\}
$$
are the stable classes determined by $F$ (with the trivial action of
$\ZZ$) and $G\wedge G$.
\qed}
Lastly, the Hopf invariant satisfies the following sum formula.
\Prop{\label{sum}{\rm (Sum formula).}
Let $F_+,\, F_-$ be maps $V^\cpt\wedge X\to V^\cpt\wedge Y$.
Suppose that $v\in S(V)$. Choose a tubular neighbourhood of
$\{ v,\, -v\}$ in $V$ and let $\nabla : V^\cpt\to V^\cpt\wedge V^\cpt$ be the
associated Pontryagin-Thom map.
Then the Hopf invariant of the sum $F=(F_+\wedge F_-)\circ (\nabla\wedge 1_X)$
is given by
$$
h_V(F)=h_V(F_+)+h_V(F_-) + \iota [(F_+\wedge F_-)\circ\Delta_X]
$$
where the induction homomorphism $\iota$ is the composition of the isomorphism
$$
\omega^0\{ X;\, Y\wedge Y\}
\iso
{}^\ZZ\omega^0\{ \Sigma S(Lv)_+\wedge X; \, (Lv)^\cpt\wedge (Y\wedge Y)\}
$$
and the map
$$
{}^\ZZ\omega^0\{ \Sigma S(Lv)_+\wedge X;
\, (Lv)^\cpt\wedge (Y\wedge Y)\}
\to {}^\ZZ\omega^0\{ \Sigma S(LV)_+\wedge X;\, LV^\cpt\wedge (Y\wedge Y)\}
$$
induced by the inclusion of the $1$-dimensional subspace $\Rr v\into V$.
\qed}
The explicit construction of the geometric Hopf invariant is readily
extended to the fibrewise theory.
Suppose now that $X\to B$ and $Y\to B$ are fibrewise pointed spaces
over an ENR $B$.  (We shall assume that both are locally fibre homotopy
trivial and that the fibres have the homotopy type of CW complexes,
finite complexes in the case of $X$.)
Consider a fibrewise pointed map $F: (B\times V^\cpt)\wedge_B X \to
(B\times V^\cpt)\wedge_B Y$. If $B$ is compact, we have a fibrewise
geometric Hopf invariant
$$
h_V(F) \in  {}^\ZZ\omega^0_B\{ (B\times \Sigma S(LV)_+)\wedge_B X;\,
(B\times LV^\cpt)\wedge_B (Y\wedge_B Y)\}\, ,
$$
and corresponding variants $h'_V(F)$ and $h''_V(F)$.
See, for example, \cite[Part II, Section 14]{crabbjames}.
This fibrewise Hopf invariant is an obstruction to fibrewise desuspension.
Indeed, suppose that the restriction of $F$ to a closed sub-ENR $A\subseteq
B$ is the fibrewise suspension of a map $X_A \to Y_A$ over $A$.
Then Lemma \ref{suspension} shows how to define a relative fibrewise
Hopf invariant in
$$
{}^\ZZ\omega^0_{(B,A)}\{ (B\times \Sigma S(LV)_+)\wedge_B X;\,
(B\times LV^\cpt)\wedge_B (Y\wedge_B Y)\}\, ,
$$
which lifts $h_V(F)$. (One uses the fact that the inclusion of $A$ in $B$
is a cofibration.)
When $B$ is not (necessarily) compact and $F$ is a fibrewise suspension
outside some compact subspace of $B$, the same method gives
a fibrewise Hopf invariant with compact supports:
$$
h_V(F) \in
{}_{\phantom{Z/}c}^\ZZ\omega^0_B\{ (B\times \Sigma S(LV)_+)\wedge_B X;\,
(B\times LV^\cpt)\wedge_B (Y\wedge_B Y)\}\, .
$$

\Sect{The double point theorem} Let $f: M\imm N$ be a (smooth)
immersion of a closed manifold $M$ in a connected manifold $N$
(without boundary) of dimension $n$, with normal bundle $\nu (f)$,
usually abbreviated to $\nu$. We do not require $M$ to be
connected, nor that all the components should have the same
dimension; the maximum dimension of a component is denoted by $m$.
Let $e: M\to V$ be a smooth map to a finite-dimensional Euclidean
space $V$ such that $e(x)\not=e(y)$ whenever $f(x)=f(y)$ for
$x\not=y$. This gives a (smooth) embedding $(e,f): M \into V\times
N$ with normal bundle $V\oplus\nu$. (To be precise, there is a
short exact sequence $0\to V\to \nu (e,f)\to \nu \to 0$ which is
split by a choice of metrics.)

We have an associated Pontryagin-Thom map (defined up to homotopy)
$$
F : V^\cpt\wedge N^\cpt \to V^\cpt\wedge M^\nu .
$$
Its geometric Hopf invariant is a stable $\ZZ$-homotopy class
$$
h_V(F)\in{}^\ZZ\omega^0\{ \Sigma S(LV)_\cpt \wedge N^\cpt  ;\,
 LV^\cpt \wedge (M^\nu\wedge M^\nu )\}\, ,
$$
where $\ZZ$ interchanges the factors of $M^\nu\wedge M^\nu$.

Suppose that the immersion $f$ is self-transverse and that
there are no $k$-tuple points for $k>2$.
(This is the case for a generic immersion if $3m<2n$.)
The {\it double point set}
$$
\Dp = \{ (x,y)\in M\times M-\Delta (M) \st f(x)=f(y)\}
$$
is then a smooth $\ZZ$-submanifold of $M\times M$
(of constant dimension $2m-n$ if $M$ is connected), on which
$\ZZ$ acts freely, and its normal bundle may be identified
with the pullback $j^*\tau N$ of
the tangent bundle of $N$ by the map $j:\Dp \to N$ mapping $(x,y)$ to
$f(x)=f(y)\in N$.
We also have a $\ZZ$-map $d :\Dp \to LV-\{ 0\}$ given
by $d(x,y)=e(x)-e(y)$, and thus an embedding
$$
(d,j) :  \Dp \into (LV-\{ 0\})\times N\, ,
$$
with normal bundle $LV\oplus k^*(\nu \times \nu )$,
where $k: \Dp\to M\times M$ is the inclusion.
The Pontryagin-Thom construction gives a $\ZZ$-map
$$
(LV-\{ 0\})^\cpt\wedge N^\cpt \to
LV^\cpt\wedge \Dp^{k^*(\nu\times\nu )}.
$$
Composing with the map induced by $k$, we get a $\ZZ$-homotopy
class
$$
\phi :\Sigma S(LV)_+ \wedge N^\cpt \to LV^\cpt\wedge
(M^\nu\wedge M^\nu )\, .
$$
\Thm{\label{dpt}
{\rm (The double point theorem).}
The geometric Hopf invariant $h_V(F)$ of the Pontryagin-Thom map
$F$ is equal to the $\ZZ$-equivariant
stable homotopy class of the map $\phi$ determined, as described above,
by the double point manifold $\Dp$.
\qed}
We may also consider the second version of the Hopf invariant
$$
h'_V(F) \in{}^\ZZ\omega^0\{ N^\cpt ;\,
S(LV)_+ \wedge (M^\nu\wedge M^\nu )\}\, .
$$
This, too, may be described directly in terms of the double points.
The embedding $\Dp \into LV\times N$ with normal
bundle $LV\oplus k^*(\nu\times\nu )$
provides a Pontryagin-Thom map
$$
LV^\cpt\wedge N^\cpt \to
LV^\cpt\wedge \Dp^{k^*(\nu\times\nu )}
$$
which we compose with the map $\Dp \to S(LV)\times (M\times M)$
given by $e$ and the inclusion $k$ to get
$$
\phi' : LV^\cpt\wedge N^\cpt \to LV^\cpt\wedge S(LV)_+
\wedge (M^\nu\wedge M^\nu )\, .
$$
\Cor{The second version $h'_V(F)$ of the geometric Hopf invariant
is represented by the $\ZZ$-map $\phi'$ defined in the text.
\qed}
We also have the non-equivariant stable Hopf invariant
$$
h''_V(F) \in\omega^0\{ N^\cpt ;\,
(S(LV) \times_{\ZZ} (M\times M))^{\nu\times\nu}\}\, .
$$
The free $\ZZ$-manifold $\Dp$ is a double cover of the set
$\barDp\subseteq N$ of double points of the immersion $f$.
We have an induced map
$$
\barDp =\Dp /\ZZ \to S(LV)\times_{\ZZ} (M\times M)
$$
and an embedding of $\barDp$ in $N$ with normal bundle the pullback of
$\nu\times\nu$. The Pontryagin-Thom construction gives a map
$$
\phi'' :  N^\cpt \to \barDp^{\nu\times\nu} \to
(S(LV) \times_{\ZZ} (M\times M))^{\nu\times\nu}.
$$
\Cor{The stable geometric Hopf invariant $h''_V(F)$ is equal to the
stabilization $[\phi'']$ of the class determined by the double point
manifold $\barDp\subseteq N$.
\qed}

These results will be obtained as consequences of a more precise
fibrewise theorem, which we describe next.

Let $\CC\to N$ be the space of  pairs $(x,\alpha )$
where $x\in M$ and
$\alpha :[0,1]\to N$  is a continuous path
such that $\alpha (0)=f(x)$, fibred over $N$ by projection to
the other endpoint $\alpha (1)$.
We have a homotopy equivalence $\pi :\CC \to M$ given by $\pi (x,\alpha )=x$.

The homotopy Pontryagin-Thom map defined by
the embedding $(e,f) : M\into V\times N$
as described in \cite[Section 6]{crabbcoincidence} (and in \cite{crabbjames})
is a pointed map
$$
\tilde F :N\times V^\cpt \to (N\times V^\cpt)\wedge_N\CC_N^{\pi^*\nu}
$$
with compact supports over $N$.
(To be exact,
the space $\CC$ in \cite{crabbcoincidence} is the space of pairs $(x,\beta )$,
where $x\in M$ and $\beta : [0,1] \to V\times N$ is a path starting at
$\beta (0)= (e(x),f(x))$. We omit here the redundant component in the
contractible space $V$.)
We may then form the fibrewise geometric Hopf invariant
$$
h_V(\tilde F) \in{}_{\phantom{Z/}c}^\ZZ\omega^0_N\{ N\times \Sigma S(LV)_+;\,
(N\times LV^\cpt)\wedge_N
(\CC_N^{\pi^*\nu} \wedge_N \CC_N^{\pi^*\nu})\}\, ,
$$
again as a stable $\ZZ$-equivariant map with compact supports over $N$.

Now the fibre product $\CC\times_N\CC$ of pairs $((x,\alpha ),(y,\beta ))$
such that $\alpha (1)=\beta (1)$ may be identified, by splicing
$\alpha$ to the reversed path $\beta$, with
the space $\CD$ of triples $(x,y,\gamma )$ with $x,\, y\in M$ and $\gamma :
[-1,1]\to N$ a continuous path from $\gamma (-1)=f(x)$ to
$\gamma (1)=f(y)$ projecting to $\gamma (0)\in N$.
(Thus $\gamma (t)$ is $\alpha (1+t)$ if $-1\leq t\leq 0$,
$\beta (1-t)$ if $0\leq t\leq 1$.)
It has an action of $\ZZ$ in which the involution interchanges
$x$ and $y$ and reverses the path $\gamma$.
The double point set $\Dp$ is included in $\CD$ as the space of constant paths.

Let $\EE$ be the space of pairs $((x,y),\alpha )$, where
$(x,y)\in \Dp$ and $\alpha :[0,1]\to N$ is a path
such that $\alpha (0)=f(x)=f(y)$.
This corresponds to the subspace of points $(x,y,\gamma )\in \CD$
with $\gamma (-t)=\gamma (t)$.
The fibrewise Pontryagin-Thom construction
on $\Dp\into (LV-\{ 0\})\times N$
gives a fibrewise map
$$
N\times \Sigma S(LV)_+ \to
(N\times LV^\cpt)\wedge_N \EE_N^{\nu\times\nu},
$$
which we compose with the inclusion
$$
\EE_N^{\nu\times\nu}\into \CD_N^{\nu\times\nu},
$$
to get an equivariant fibrewise map
$$
\tilde \phi : N\times \Sigma S(LV)_+ \to
(N\times LV^\cpt)\wedge_N \CD_N^{\nu\times\nu}.
$$

\Thm{\label{hdpt}
{\rm (Homotopy double point theorem).}
The fibrewise geometric Hopf invariant
$$
h_V(\tilde F) \in
{}_{\phantom{Z/}c}^\ZZ\omega^0_N\{ N\times \Sigma S(LV)_+;\,
(N\times LV^\cpt)\wedge_N
\CD_N^{\nu\times\nu}\}
$$
of the homotopy Pontryagin-Thom map $\tilde F$
is equal to the fibrewise stable class of the map $\tilde \phi$
determined by the double points of $f$.
\qed}

The dual version is a class
$$
h'_V(\tilde F) \in
{}_{\phantom{Z/}c}^\ZZ\omega^0_N \{ N\times S^0;\, (N\times S(LV)_+) \wedge_N
\CD_N^{\nu\times\nu}\}\, .
$$
There is also a non-equivariant form.
The stable Hopf invariant $h''_V(\tilde F)$ lies in
$$
{}_c\omega^0_N\{ N\times S^0;\,
(S(LV)\times_{\ZZ} \CD )_N^{\nu\times \nu}\}\, ,
$$
and this group can be identified with
$$
\tilde\omega_0((S(LV)\times_{\ZZ} \CD )^{\nu\times\nu -\tau N})
$$
by fibrewise Poincar\'e-Atiyah duality. (For a general treatment of
fibrewise duality see, for example, \cite[Part II, Section 12]{crabbjames}.
The duality theorem required
here is stated in \cite{crabbcoincidence} as Proposition 4.1.)

The Pontryagin-Thom construction applied to the
double point manifold $\barDp =\Dp /(\ZZ)$ equipped with the map
$\barDp\to S(LV)\times_{\ZZ} \CD$ given
by the inclusion $\Dp\to\CD$ and the map $\Dp\to S(LV)$ given by $e$
(via $d$) produces a stable homotopy class
$$
\tilde \phi'' :
S^0 \to (S(LV)\times_\ZZ \CD )^{\nu\times\nu -\tau N} .
$$
\Cor{We have
$$
h''_V(\tilde F)=\tilde \phi'' \in
\tilde\omega_0((S(LV)\times_{\ZZ} \CD )^{\nu\times\nu -\tau N})\, .
$$
\qed}
Before turning to the proof of Theorem \ref{hdpt},
we explain how the fibrewise Hopf invariant $h_V(\tilde F)$ determines
the simpler invariant $h_V(F)$.
\Lem{The Hopf invariant $h_V(F)$ is the image of $h_V(\tilde F)$ under
the composition:
$$
\begin{array}{ll}
&{}^\ZZ_{\phantom{Z/}c}\omega^0\{ N\times \Sigma S(LV)_+;\,
(N\times LV^\cpt)\wedge_N \CD^{\nu\times\nu}_N\} \\[1ex]
&\qquad\to
{}^\ZZ\omega^0\{ N^\cpt\wedge \Sigma S(LV)_+;\,
LV^\cpt\wedge \CD^{\nu\times\nu}\}\\[1ex]
&\qquad\to
{}^\ZZ\omega^0\{ N^\cpt\wedge \Sigma S(LV)_+;\,
LV^\cpt\wedge (M^\nu\wedge M^\nu )\}
\end{array}
$$
of the homomorphism defined
by collapsing the basepoint sections over $N$ to a point
and that induced by
the projection $\pi\times\pi :\CD = \CC\times_N\CC\to M\times M$.
}
\begin{proof}
This is easily seen from the explicit construction of the geometric
Hopf invariant.
\qed
\end{proof}

In a similar manner, $h'_V(F)$ and $h''_V(F)$ are the images
of the refined invariants $h'_V(\tilde F)$ and $h''_V(\tilde F)$
under homomorphisms defined by collapsing basepoint sections over $N$
and projecting from $\CD$ to $M\times M$.
\Rem{This construction also provides the $\pi$-equivariant Hopf invariant
considered, in greater detail, in this volume.
Suppose that $\Gamma$ is a discrete group and that $q:\widetilde N\to N$
is a principal $\Gamma$-bundle (for example, a universal covering
space with $\Gamma$ the fundamental group of $N$).
We let $\widetilde M=f^*\widetilde N$ be the induced bundle over $M$:
thus
$$
\widetilde M=\{ (x,z)\in M\times\widetilde N\st f(x)=q(z)\}\, .
$$
We may define a map $\CD \to (\tilde M\times \tilde M)/\Gamma$ as follows.
Given $(x,y,\gamma )\in\CD$ (so $x,\, y\in M$, $\gamma : [-1,1]\to N$,
$\gamma (-1)=f(x)$, $\gamma (1)=f(y)$),
lift $\gamma$ to a path $\tilde \gamma$ in $\widetilde N$, determined up
to multiplication by an element of $\Gamma$. The $\Gamma$-orbit of
$((x,\tilde\gamma (-1)),(y,\tilde\gamma (1)))\in \widetilde M\times\widetilde M$
is independent of the choice of the lift.
The Hopf invariant $h''_V(\tilde F)$ thus gives us an element of
$$
\omega^0\{ N^\cpt ;\,
(S(LV)\times_\ZZ ((\widetilde M\times\widetilde M)/\Gamma ))
^{\nu\times\nu }\}\, ,
$$
where $\nu\times\nu$ is lifted from $M\times M$ to
$(\widetilde M\times\widetilde M)/\Gamma$ by the obvious projection.
\qed}

In view of the relations between the various Hopf invariants, it
will suffice to prove Theorem \ref{hdpt} in order to establish Theorem \ref{dpt}
and their sundry corollaries.

\medskip

\noindent
{\it Proof of Theorem \ref{hdpt}.}
Writing $\bar \nu$ for the normal bundle of $\barDp$ in $N$, choose
a tubular neighbourhood $D(\bar\nu )\into N$ of $\barDp$.
For $(x,y)\in \Dp$, the fibre of $\bar\nu$ at $f(x)=f(y)$ is
$\nu_x\oplus\nu_y$.
We identify $f^{-1}(\barDp)$ with $\Dp$ by projecting to the first factor.
Then the inverse image of the tubular neighbourhood $D(\bar\nu )$
is a tubular neighbourhood $D(\nu')$ of $\Dp$ in $M$,
where $\nu'_x=\nu_y$. The normal bundle $\nu$ restricted to $D(\nu ')$
is then identified with the restriction of $\nu$ to $\Dp$.

We may assume that $e$ vanishes outside the tubular neighbourhood $D(\nu')$.
The homotopy Pontryagin-Thom map is then a $V$-fold suspension
outside $D(\bar\nu )$. By Proposition \ref{suspension}, the fibrewise
Hopf invariant is canonically null-homotopic outside  $D(\bar\nu )$.
The Hopf invariant $h_V(\tilde F)$ is thus represented by a fibrewise
map which is zero outside the tubular neighbourhood.
In this way we localize $h_V(\tilde F)$ to
an element of
$$
{}^\ZZ\omega^0_{(D(\bar\nu ),S(\bar\nu ))}\{ N\times \Sigma S(LV)_+;
\, (N\times LV^\cpt)\wedge_N \CD^{\nu\times\nu}_N\}
$$
constructed from the immersion data in a neighbourhood of the
double points.

The local data consists simply of the double cover $\Dp\to \barDp$
and the vector bundle $\nu$ over $\Dp$. The bundle $\nu'$ is the
pullback of $\nu$ under the covering involution and $\bar\nu$
over $\barDp$ is the push-forward of $\nu$.
The `local $M$' is the total space of $\nu'$ over $\Dp$, and the `local
$N$' is the total space of $\bar\nu$ over $\barDp$. The immersion
$f$ is given by the projection $\Dp\to\barDp$.
We also need the map $e$ and, without loss of generality, we may assume
that it is determined by a $\ZZ$-map $\Dp\to S(LV)$ (extended radially
on $D(\nu')$ to taper to zero).

In the fibre over $\{ x,y\}\in\barDp$, the manifold $M$ is the disjoint
union
$$
(\{0\}\times\nu_y)\sqcup (\nu_x\times\{ 0\}) \imm\nu_x\oplus\nu_y
$$
immersed in the fibre of $N$ by projection, with the double point at $(0,0)$.
Write $e(x)=v$, $e(y)=-v$. The Hopf invariant is calculated by
the sum formula of Proposition \ref{sum}.
In the notation used there, we take $X= \nu_x^\cpt\times\nu_y^\cpt$,
$Y=\nu_x^\cpt\vee\nu_y^\cpt$,
and the maps $F_+$ and $F_-$ are suspensions of
the compositions of the projection to $\nu_x^\cpt$
or $\nu_y^\cpt$, respectively, and the inclusion of the wedge summand.
The Hopf invariants of the suspensions $F_+$ and $F_-$ vanish
and the Hopf invariant of the sum is determined by the product
$X\to Y\times Y$, so by the projection $\nu_x^\cpt\times\nu_y^\cpt\to
\nu_x^\cpt\wedge \nu_y^\cpt$.

The same computation performed fibrewise over $\barDp$ gives
the localized fibrewise Hopf invariant
as the image of the element in
$$
\omega^0_{(D(\bar \nu ),S(\bar\nu ))}\{ D(\bar\nu )\times S^0;\,
\bar\nu^\cpt_{D(\bar\nu )}\}
$$
given by the inclusion of $D(\bar\nu )$ in $\bar\nu$.
Unravelling the definitions, one sees that this produces $[\phi ]$.
\qed

\smallskip

In the remainder of the paper we shall be concerned with the behaviour
of the geometric Hopf invariant $h_V(\tilde F)$ as one deforms the
map $f: M\to N$.
Suppose, to begin, that we have
smooth homotopies $f_t : M\to N$ and $e_t:M\to V$
such that each $f_t$ is an immersion
and such that $(e_t,f_t) : M\to V\times N$ is an
embedding for each $t\in [0,1]$.
We write $f=f_0$, $f'=f_1$, and $e=e_0$, $e'=e_1$.
Then we have, up to homotopy, an isomorphism
$\nu'=\nu (f_1)\to \nu (f_0)=\nu$
between the normal bundles of the immersions
and a fibre homotopy equivalence $\CDone \to \CD$ over $N$.
In the standard terminology, the homotopy $f_t$ is a {\it regular
homotopy} from $f$ to $f'$ and the homotopy $(e_t,f_t)$ is
an {\it isotopy} from $(e,f)$ to $(e',f')$.
Let $F'$ be the map obtained from the homotopy Pontryagin-Thom
construction on $(e',f')$.
\Prop{{\rm (Regular homotopy/isotopy invariance).}
The fibrewise Hopf invariants $h_V(\tilde F')$ and $h_V(\tilde F)$
correspond under the induced isomorphism:
$$
\begin{array}{ll}
&{}_{\phantom{Z/}c}^\ZZ\omega^0_N\{ N\times \Sigma S(LV)_+;\,
(N\times LV^\cpt)\wedge_N
\CDone_N^{\nu'\times\nu'}\} \\[1ex]
&\to
{}_{\phantom{Z/}c}^\ZZ\omega^0_N\{ N\times \Sigma S(LV)_+;\,
(N\times LV^\cpt)\wedge_N
\CD_N^{\nu \times\nu }\} \, .
\end{array}
$$
}
\begin{proof}
This follows from the homotopy invariance of the geometric Hopf invariant,
which in turn follows easily from the explicit construction described
in the Introduction.
\qed
\end{proof}

\Sect{Immersions and embeddings}
The vanishing of the stable Hopf invariant $h_V(\tilde F)$ is a
necessary condition for the existence of a regular homotopy $f_t$
and isotopy $(e_t,f_t)$ such that $f'=f_1$ is an embedding.
In the opposite direction, we record first:
\Prop{Suppose that $3m<2(n-1)$ and that $h_V(\tilde F)=0$. Then
$\tilde F$ is
homotopic (through maps with compact support over $N$) to the
$V$-fold suspension
of a map determined by a section $N\to\CC^{\pi^*\nu}_N$.
}
\begin{proof}
This is a consequence of
the fibrewise EHP-sequence (see \cite[Part II, Proposition 14.39]{crabbjames}),
applicable because $0\leq 3(n-m-1)-(n+0)$.
\qed\end{proof}
To proceed further, we shall assume
that the dimension of $V$ is large: $\dim\, V >2m$.
This guarantees that $M$ can be embedded in $V$, and
we may take the map $e:M\to V$ to be an embedding.
In a metastable range, the fibrewise Hopf invariant of $\tilde F$ is
the precise obstruction to the existence of a regular homotopy
from the immersion $f:M\imm N$ to an embedding.
\Thm{\label{HHQ}
{\rm (Haefliger \cite{haefliger}, Hatcher-Quinn \cite{hatcherquinn}).}
Suppose that $3m<2(n-1)$ and $\dim V>2m$. Then
$$
h''_V(\tilde F) \in \tilde\omega_0((S(LV)\times_\ZZ \CD)^{\nu\times\nu
-\tau N})
= \tilde\omega_0((E\ZZ \times_\ZZ \CD)^{\nu\times\nu
-\tau N})
$$
vanishes if and only if $f$ is homotopic through immersions to
an embedding of $M$ into $N$.
\qed}

We define $\hDp$ to be the subspace of $\CD$ consisting of the
triples $(x,y,\gamma )$ such that $x\not=y$
and call $\hDp$ the
space of {\it homotopy double points} of the immersion $f$.
The Hopf bundle $E\ZZ\times_\ZZ L$ over $B\ZZ$ will be denoted by $H$.
\Lem{\label{free}
Let $\CCC$ be the complement of
$\hDp$ in $\CD$. Then one has a homotopy
cofibration sequence:
$$
\begin{array}{ll}
&(E\ZZ\times_\ZZ \hDp )^{\nu\times\nu -\tau N}
\to  (E\ZZ\times_\ZZ \CD )^{\nu\times\nu -\tau N} \\ [1ex]
&\qquad\qquad
\to  (E\ZZ\times_\ZZ \CCC )^{H\otimes\tau M+\nu\times\nu -\tau N},
\end{array}
$$
in which the first map is given by the inclusion of $\hDp$ in $\CD$
and the second by the homotopy Pontryagin-Thom construction on the
diagonal submanifold $M$ in $M\times M$.
}

As a space over $M$, $\CCC$ is the pullback by $f : M\to N$ of the free loop
space of $N$, $\map (S^1,N)\to N$, fibred over $N$ by evaluation at
$1\in S^1\, (\subseteq\Cc )$.
\begin{proof}
The vector bundle $L\otimes\tau M$, corresponding to $H\otimes \tau M$
over $B\ZZ$, is the normal bundle of the diagonal
inclusion of $M$ in $M\times M$.
Choose a $\ZZ$-equivariant tubular neighbourhood
$D(L\otimes\tau M)\into M\times M$ of the diagonal
in $M\times M$.
The inclusion of the subspace of $\hDp$ consisting of the triples
$(x,y,\gamma )$ such that $(x,y)\notin B(L\otimes \tau M)$
into $\hDp$ is a homotopy equivalence.

The argument used in the proof of Lemma 6.1 in \cite{crabbcoincidence} then shows
that we have a homotopy cofibre sequence.
\qed\end{proof}
\Cor{\label{isomorphism}
The inclusion induces an isomorphism
$$
\tilde\omega_0((E\ZZ\times_\ZZ \hDp )^{\nu\times\nu -\tau N})
\to
\tilde\omega_0((E\ZZ\times_\ZZ \CD )^{\nu\times\nu -\tau N})
$$
provided that $m<n-1$.
}
\begin{proof}
This follows at once from the long exact sequence of the cofibration.
\qed\end{proof}

The involution on $\hDp$ is free, and hence, writing
$\barhDp$ for the quotient $\hDp /(\ZZ )$, we have an isomorphism
$$
\tilde\omega_0((E\ZZ\times_\ZZ \hDp )^{\nu\times\nu -\tau N})
\to
\tilde\omega_0 (\barhDp^{\nu \times\nu -\tau N})\, .
$$

Until now, we have thought of $\CD$, which arose as the
fibre product $\CC\times_N\CC$, as a space over $N$. It also
fibres over $M\times M$ and the fibrewise space $\CD \to M\times M$
is $\ZZ$-equivariantly locally fibre homotopy trivial.
(Compare \cite[Definition 2.3]{crabbcoincidence}.)
In the same way, $\hDp$ is fibred over the complement $M\times M-\Delta (M)$
and $\barhDp$ is fibred over $(M\times M-\Delta (M))/\ZZ$.

\medskip

\noindent
{\it Proof of Theorem \ref{HHQ}.}
We shall use results and terminology from \cite[Section 7]{crabbcoincidence},
which derive from work of Koschorke \cite{koschorke2} and
Klein and Williams \cite{kleinwilliams1}.

Suppose, first, that $M$ is connected.
Put $\tilde B =M\times M - B(L\otimes \tau M))$; it is a manifold
with a free $\ZZ$-action.
Let $B$ be the manifold $\tilde B/\ZZ$
with boundary $\partial B$ the projective bundle $P(\tau M)$.
Let $E$ be the bundle $(\tilde B \times (N\times N))/\ZZ$ over $B$
and $Z\subseteq E$ the diagonal sub-bundle $B\times N=(\tilde B \times N)/\ZZ$.
The fibrewise normal bundle of the inclusion
$Z\into E$ is $H\otimes \tau N$,
where $H$ is the line bundle over $B$ associated to the double
cover.
The $\ZZ$-equivariant square $f\times f: \tilde B\to N\times N$
defines a section $s$ of $E\to B$, which, if the tubular
neighbourhood is chosen to be sufficiently small, has the property that
$s(x)\notin Z_x$ for $x\in\partial B$.
Together, the fibre bundle $E\to B$, the sub-bundle $Z\to B$ and the section $s$
constitute what is called in \cite{crabbcoincidence} an intersection problem.
In the language used there, $s$ is nowhere null on the boundary
$\partial B$ and the homotopy null-set
 $\hNull (s)$, fibred over $B$,
is easily identified with the restriction of
$\barhDp\to (M\times M-\Delta (M))/\ZZ$.
The inclusion $\hNull (s)\into \barhDp$ is, as we have already
noted in the proof of Lemma \ref{free}, a homotopy equivalence.

If $f$ is an embedding, then the section $s$ is nowhere null.
Now the homotopy Euler index
(\cite[Definition 7.3]{crabbcoincidence})
$$
\hInd (s; \partial B) \in
\tilde\omega_0 (\hNull (s)^{H\otimes\tau N - \tau B)})
$$
is an obstruction to deforming $s$, through a homotopy constant on
$\partial B$, to a section that is nowhere null.
By Proposition 7.4 of \cite{crabbcoincidence}
it is the precise obstruction if $\dim B < 2(\dim N -1)$.
We conclude that if $\hInd (s;\, \partial B)=0$ and $m < n -1$,
then
$$
f\times f : M\times M - B(L\otimes\tau M) \to N\times N
$$
is $\ZZ$-equivariantly homotopic,
by a homotopy that is constant on the boundary $S(L\otimes\tau M)$,
to a map into $N\times N - \Delta (N)$.
According to Haefliger \cite[Th\'eor\`eme 2]{haefliger},
if further $3m<2(n-1)$,
then this is a sufficient condition
for $f:M\imm N$ to be homotopic through immersions to an embedding
of $M$ into $N$.

Thus far, the
argument is taken from \cite[Theorem A.4 and Corollary A.5]{kleinwilliams1}.
We now relate the index $\hInd (s;\partial B)$ to the geometric Hopf invariant.
\Lem{The homotopy Euler index
$$
\hInd (s;\partial B)\in
\tilde\omega_0 (\hNull (s)^{\nu \times\nu -\tau N})
$$
corresponds to the fibrewise Hopf invariant
$$
h_V''(\tilde F)\in\tilde\omega_0((E\ZZ\otimes_{\ZZ} \CD)^{\nu\times\nu
-\tau N})\, .
$$
}
\begin{proof}
The correspondence is made via the homotopy equivalence
$\hNull (s)\into \barhDp$ and the isomorphism from Corollary \ref{isomorphism}.

We can assume that $f$ satisfies the conditions for the homotopy double point
theorem.
Then $h_V''(\tilde F)$ is represented by the double point manifold
$\barDp$.
By Proposition 7.8 of \cite{crabbcoincidence}, the homotopy Euler index
is represented by the null-set $\Null (s)$, which is exactly
$\barDp$.
\qed\end{proof}
Hence the vanishing of $h''_V(\tilde F)$ implies, in the metastable range,
that $f$ is regularly homotopic to an embedding.

This has dealt with the case in which $M$ is connected.
Now suppose that $M$ is a disjoint union $M^{(1)}\sqcup M^{(2)}$
and write $f^{(i)}$ for the restriction of $f$ to $M^{(i)}$.
Then $\CD$ splits equivariantly over
$$
M\times M = (M^{(1)}\times M^{(1)})
\sqcup (M^{(2)}\times M^{(2)}) \sqcup
(M^{(1)}\times M^{(2)}\sqcup M^{(2)}\times M^{(1)})
$$
as a disjoint union
$$
\CDa\sqcup \CDb\sqcup (\ZZ \times \hI )\, ,
$$
where $\hI$ is the space of triples $(x,y,\gamma )$ with
$(x,y)\in M^{(1)}\times M^{(2)}$ and $\gamma (-1)=f^{(1)}(x)$,
$\gamma (1)=f^{(2)}(y)$.
The fibrewise Hopf invariant decomposes, according to
Proposition \ref{sum}, as a sum of three terms.
The first two are the fibrewise Hopf invariants of $f^{(1)}$ and
$f^{(2)}$; the third, $(12)$-component,
is a more elementary product obstruction.

Arguing by induction, we may suppose that both $f^{(1)}$ and $f^{(2)}$
are embeddings, intersecting transversely,
and that $M^{(2)}$ is connected. We consider a new
intersection problem with $B=M^{(2)}$, $E=B\times N$ and
$Z=B\times f(M^{(1)})$. Let $s : B\to E$ be the section
$s(x)=(x,f^{(2)}(x))$.
This time $\hNull (s)$ is $\hI$,
and we may identify the homotopy Euler index $\hInd (s)$ in
the same way with the $(12)$-component of the fibrewise Hopf
invariant, both being represented by the manifold
$f(M^{(1)})\cap f(M^{(2)})$.
By Proposition 7.4 of \cite{crabbcoincidence} again, the vanishing of the
homotopy Euler index implies that $f^{(2)}$ is homotopic to a map
into $N-f(M^{(1)})$, because $\dim M^{(2)} < 2(n -\dim M^{(1)}-1)$.
Finally, we may apply \cite[Theorem 1.1]{hatcherquinn}
to deduce that $f^{(2)}$ is isotopic to an embedding of
$M^{(2)}$ into the complement of $f(M^{(1)})$.
This inductive step is enough to conclude the proof of Theorem \ref{HHQ}.
\qed

\Rem{\label{fun}
Suppose that $M$ (as well as $N$) is connected. Choose basepoints $*\in M$
and $*=f(*)\in N$. Then we can include the loop-space
$\Omega N$ in $\CD$ by mapping
a loop $\gamma : [-1,1]\to N$, with $\gamma (-1)=*=\gamma (1)$,
to $(*,*,\gamma )\in\CD$, and the set of path components of $\CD$ is
in this way identified with the set of double cosets
$$
f_*\pi_1(M)\backslash\pi_1(N)/f_*\pi_1(M) =\pi_0(\CD )
$$
with the $\ZZ$-action given by the group-theoretic inverse. Thus
we may identify
the set of path components of $S(LV)\times_\ZZ \CD$, if $\dim V>1$, with
the orbit space of the involution.
\qed}
\Ex{Suppose that $n=2m$.
Then $\tilde\omega_0((E\ZZ\times_\ZZ \CD )^{\nu\times\nu -\tau N})$
is a direct sum of groups indexed by $\pi_0(\CD )/\ZZ$,
each component being isomorphic to $\Zz$ or $\Zz /2\Zz$.
When $M$ is connected we may label the summands as in Remark \ref{fun}
by equivalence classes of group elements $g\in\pi_1(N)$.
Let $w_M : \pi_1(M)\to \{ \pm 1\}$ and $w_N : \pi_1(N)\to \{\pm 1\}$
be the orientation maps (corresponding to $w_1M$ and $w_1N$).
The $g$-summand is isomorphic to $\Zz$ if and only if
for all $a,\, b\in\pi_1(M)$:
(o) if $f_*(a)g=gf_*(b)$, then $w_M(ab)=w_N(f_*(a))$ ($=w_N(f_*(b))$),
and
(i) if $f_*(a)g=g^{-1}f_*(b)$, then $(-1)^mw_M(ab)=w_N(f_*(a)g)$.
In particular, if $M$ is orientable,
$N$ is oriented and $f_*\pi_1(M)$ is trivial,
then the obstruction group is
$\Zz [\pi_1(N)]/\langle g-(-1)^mg^{-1}\st g\in\pi_1(N)\rangle$
and the Hopf invariant $h''_V(\tilde F)$
is Wall's invariant $\mu (f)$ in the form mentioned in the Introduction.
\qed}
\Sect{Homotopic immersions}
We next investigate immersions homotopic, as maps, to the given immersion
$f$.
Consider smooth homotopies $f_t : M\to N$ and $e_t : M\to V$
such that $f_0$ and $f_1$ are immersions and each map
$(e_t,f_t):M\to V\times N$, for $0\leq t\leq 1$, is an embedding with
normal bundle $\mu_t$.
We again write $f=f_0$, $f'=f_1$, and $\nu =\nu (f)$, $\nu'=\nu (f')$.
The homotopies determine, up to a homotopy,
a bundle isomorphism $a:V\oplus \nu'=\mu_1\to\mu_0 =V\oplus \nu$
and a $\ZZ$-equivariant fibre homotopy equivalence $\CDone \to \CD$ over
$M\times M$.

There is an associated class
$$
\theta (e_t,f_t)\in\tilde\omega_0((P(V)\times M)^{H\otimes\nu -\tau M})
$$
which vanishes if each $f_t$ is an immersion.
It is constructed as follows.
Let $v_0 : V\to V\oplus \nu$ be the inclusion and
let $v_1: V\to V\oplus \nu'=\mu_1 \to\mu_0=V\oplus\nu$
be the composition of the inclusion with the isomorphism $a$.
Now we have a stable cohomotopy difference class
$$
\delta (v_0,v_1)\in
\tilde\omega^{-1}((M\times P(V))^{-H\otimes (V\oplus\nu )})\, ,
$$
which is the metastable obstruction to deforming $v_0$ to $v_1$
through vector bundle monomorphisms;
see Section 6 for details.
We define $\theta (e_t,f_t)$ to be the dual class in stable homotopy.
(Recall that $\tau P(V)\oplus\Rr =H\otimes V$.)

Using the inclusion $M\into\CD$ : $x\mapsto (x,x,\gamma )$,
where $\gamma$ is the constant loop at $f(x)$, we get a map
$$
i:
\tilde\omega_0((P(V)\times M)^{H\otimes \nu -\tau M})
\to
\tilde\omega_0((S(LV)\times_{\ZZ} \CD )^{\nu\times\nu-\tau N})\, .
$$
\Rem{The map $M\to\CD$ picks out
(under the assumption that $M$ is connected)
a component $\CD_0$ of $\CD$,
which is preserved by $\ZZ$.
Hence $i$ maps into the summand
$$
\tilde\omega_0((S(LV)\times_\ZZ \CD_0)^{\nu\times\nu -\tau N})\, .
$$
\qed}

The link between the difference class $\theta$ and the Hopf invariant
is forged by a generalization of the classical description of
the Hopf invariant on the image of the $J$-homomorphism.
\Prop{\label{ImJ}
Let $a : V\oplus \nu' \to V\oplus \nu$ be a vector bundle isomorphism
over $M$. The fibrewise one-point compactification of $a$ is a map
of sphere bundles
$$
A: (M\times V^\cpt)\wedge_M (\nu')^\cpt_M \to (M\times V)^\cpt\wedge \nu^\cpt_M
$$
over $M$. Then the fibrewise Hopf invariant
$$
h_V(A) \in {}^\ZZ\omega^0_M\{ (M\times \Sigma S(LV)_+)\wedge_M (\nu')^\cpt_M
;\, (M\times LV^\cpt)\wedge_M (\nu^\cpt_M \wedge_M \nu^\cpt_M) \}
$$
of $A$
coincides up to sign,
under the identifications described below, with the difference class
$$
\delta (v_0,v_1) \in
\tilde\omega^{-1}((M\times P(V))^{-H\otimes(V\oplus \nu )}))
$$
of the monomorphisms $v_0,\, v_1 : V\to V\oplus\nu$ given, respectively,
by the inclusion of the first factor and the composition of the inclusion
$V\to V\oplus\nu'$ with $a$.
}
\begin{proof}
The fibrewise smash product $\nu^\cpt_M\wedge_M\nu^\cpt_M
=(\nu\oplus\nu )_M^\cpt$ with the action of $\ZZ$ which
interchanges the factors is equivariantly homeomorphic
(by the construction $\kappa_V$ in Section 1) to
$(\nu\oplus L\nu )^\cpt_M = \nu^\cpt_M\wedge (L\nu )^\cpt_M$. The
isomorphism $a: V\oplus \nu'\to V\oplus\nu$ gives a stable fibre
homotopy equivalence $(\nu')^\cpt_M\to\nu^\cpt_M$. Taken together,
these equivalences allow us to think of $h_V(A)$ as an element of
the group
$$
{}^\ZZ\omega^0_M\{ M\times \Sigma S(LV)_+;\,
(LV\oplus L\nu )^\cpt_M\}\, ,
$$
which is then identified with
$$
{}^\ZZ\tilde\omega^{-1}((M\times S(LV))^{-(LV\oplus L\nu)})
=
\tilde\omega^{-1}((M\times P(V))^{-H\otimes(V\oplus \nu )}) \, .
$$
Both $h_V(A)$ and $\delta (v_0,v_1)$ are defined by difference constructions.
The proof that they coincide follows from
a direct comparison of the definitions.
\qed\end{proof}
\Thm{\label{Var}
{\rm (Homotopy/isotopy variation).}
The Hopf invariant
$$
h''_V(\tilde F')\in \tilde\omega_0((S(LV)\times_\ZZ \CDone)
^{\nu'\times\nu' -\tau N})
$$
corresponds to
$$
h''_V(\tilde F) + i\theta (e_t,f_t) \in \tilde\omega_0((S(LV)\times_\ZZ \CD)
^{\nu \times\nu  -\tau N})\, .
$$
}
\begin{proof}
Recollect that $\tilde F$ is a fibrewise map with compact supports
over $N$:
$$
N\times V^\cpt \to (N\times V^\cpt)\wedge_N \CC_N^{\pi^*\nu}.
$$
Up to homotopy, the map $\tilde F'$ associated with $(e',f')$ is
the composition of $\tilde F$ with the map
$$
(N\times V^\cpt)\wedge_N \CC_N^{\pi^*\nu} =\CC_N^{\pi^*(V\oplus\nu )} \to
\CC_N^{\pi^*(V\oplus \nu')}=(N\times V^\cpt)\wedge_N \CC_N^{\pi^*\nu'}
$$
obtained by lifting the bundle isomorphism $a^{-1} : V\oplus \nu
\to V\oplus \nu'$ over $M$ via $\pi : \CC\to M$.
In the notation of Proposition \ref{ImJ} we have $\tilde F\simeq
\pi^*A\comp\tilde F'$.

The fibrewise version of the composition formula
(Proposition \ref{composition}) expresses the
difference $\alpha =h_V(\tilde F)-h_V(\tilde F')$
in terms of $A$.
Using Proposition \ref{ImJ} and the explicit form of the geometric Hopf
invariant one sees that $\alpha$ lies in the image of the diagonal
map
$$
\begin{array}{ll}
\Delta_* :
&{}^\ZZ_{\phantom{Z/}c}\omega_N^0\{ N\times \Sigma S(LV)_+;\,
(N\times LV^\cpt)\wedge_N \CC^{\pi^*(\nu \oplus L\nu )}_N\} \\
&\qquad
\to {}^\ZZ_{\phantom{Z/}c}\omega_N^0\{ N\times \Sigma S(LV)_+;\, (N\times LV^\cpt)\wedge_N
(\CC^{\pi^*\nu}_N\wedge_N\CC^{\pi^*\nu )}_N)\} \, .
\end{array}
$$
This may be rewritten in dual form as:
$$
\tilde\omega_0((P(V)\times \CC)^{\pi^*(\nu \oplus H\otimes\nu -\tau N)})
\to \tilde\omega_0((S(LV)\times_\ZZ \CD) ^{\nu \times\nu  -\tau N})\, .
$$
But $\pi :\CC \to M$ is a homotopy equivalence.
Hence $\Delta_*$ is just another manifestation of the map $i$ in
the statement of the theorem.

Now the tubular neighbourhood of $M$ in $V\times N$ gives a proper
map $B(\nu ) \to N$. To calculate $\alpha$, which is concentrated on
the diagonal, we can thus lift from $N$ to $B(\nu )$.
Here we have a fibrewise problem over $M$, and the identification
of $\alpha$ is achieved by Proposition \ref{ImJ}.
\qed\end{proof}
\Rem{In the stable range $\dim V >2m$, where the maps $e_t$ are
redundant, we can use the methods of
the previous section to give an alternative proof of Theorem \ref{Var}.
\qed}
\Thm{{\rm (Hirsch \cite{hirsch}).}
Suppose that $3m<2n-1$ and $\dim V >2m$.
Then two immersions $f$ and $f'$ are homotopic through immersions
if and only if the associated difference class
$\theta$ in $\tilde\omega_0((P(V)\times M)^{H\otimes\nu -\tau M})$
is zero.
}
\begin{proof}
The derivative of the immersion $f$ gives a vector bundle monomorphism
$df : \tau M \to f^*\tau N$ over $M$.
According to Hirsch \cite{hirsch}, for $m<n$
immersions $f': M\to N$
together with a homotopy $f_t$ from $f=f_0$ to $f'=f_1$
are classified by homotopy
classes of vector bundle monomorphisms $\tau M \to f^*\tau N$.
In the metastable range $m+1 < 2(n-m)$, that is, $3m<2n-1$,
immersions with a homotopy to $f$ are thus classified by
$$
\tilde\omega^{-1}(P(\tau M)^{-H\otimes f^*\tau N})
=
\tilde\omega_0(P(\tau M)^{H\otimes ( f^*\tau N-\tau M)-\tau M})
=
\tilde\omega_0(P(\tau M)^{H\otimes \nu-\tau M})\, .
$$
Assuming that $\dim V>2m$ we may fix an embedding $e : M\into V$.
The derivative of $e$ includes
$\tau M$ in the trivial bundle $M\times V$ and gives an isomorphism
$$
\tilde\omega_0(P(\tau M)^{H\otimes \nu -\tau M}) \to
\tilde\omega_0((M\times P(V))^{H\otimes \nu -\tau M})\, .
$$
\vskip-1.4\baselineskip
\qed\end{proof}

\Sect{Immersions close to an embedding}
Consider the special case of a closed manifold $M$ of
(constant) dimension $m$
and a real vector bundle $\nu$ of dimension $n-m$ over $M$.
Working in the metastable range $3m<2n-1$,
we take $N$ to be the total space of $\nu$ and $f: M\to N$ to be the
embedding given by the zero section of the vector bundle.
As $e:M\to V$ we may take the constant map $0$.

Let $v_0 : V\into V\oplus\nu$ be the inclusion of the first factor.
Suppose that $v_1: V\into V\oplus \nu$ is another inclusion,
which we may assume to be isometric.
Thus $v_0$ and $v_1$ give sections of the bundle $\O(V,V\oplus \nu )$
whose fibre at $x\in M$
is the Stiefel manifold of orthogonal linear maps
$v:V\to V\oplus\nu_x$.
Let $X_1(V,\nu )$ be the sub-bundle with fibre consisting of those
linear maps $v$ such that $v+i_x$ has kernel of dimension $1$,
where $i_x:V\to V\oplus\nu_x$ is the inclusion of the
first factor.
In the range of dimensions that we are considering,
for a generic smooth section $v_1$ of $\O(V,V\oplus\nu )\to M$
the kernel of $(v_1)_x+i_x$ is nowhere of dimension $>1$,
$v_1$ is transverse to the sub-bundle $X_1(V,\nu )$ and
$v_1^{-1}(X_1(V,\nu ))$ is a submanifold
$\barD$ of $M$ of dimension $2m-n$, equipped with a map $\barD\to P(V)$
classifying the $1$-dimensional kernel.
The normal bundle of $\barD$ in $M$ is identified with
$H\otimes\nu$ and $\barD$ represents the element
$\delta\in\tilde\omega_0((M\times P(V))^{H\otimes\nu -\tau M})$
dual to $\delta (v_0,v_1)$.
More details are provided in an appendix (Section 6).

We want to construct an immersion $f'$ close to the zero section $f$
with double point set $\barD$. More precisely, we shall construct
homotopies $f_t$ and $e_t$, with $f_0=f$, $f_1=f'$ and each $(e_t,f_t)$
an embedding, such that $\barDpone =\barD$ and $\theta (e_t,f_t)=\delta$.

Choose an open tubular neighbourhood $\Omega =H\otimes\nu\into M$ of $\barD$.
Since $\dim \nu =n-m >2m-n=\dim \barD$, we can split
off a trivial line from the restriction of $\nu$ to $\barD$ as:
$\nu |\barD =\Rr\oplus\zeta$.

Whitney gave in \cite{whitney3}, for a Euclidean space $U$,
an explicit `punctured figure-of-eight' immersion:
$$
w: \Rr\oplus U\to (\Rr\oplus U)\times (\Rr\oplus U)
$$
with double points at $(\pm 1,0)$.
In slightly modified form it may be written as
$$
w(s,y) = ((1-\lambda (s,y))s,y,-\lambda (s,y),\lambda (s,y) sy)\, ,
$$
where $\lambda (s,y)=\psi (s^2+\| y\|^2)$
and $\psi : [0,\infty )\to \Rr$ is
a smooth, non-negative, monotonic decreasing function,
with $\psi (1)=1$, $\psi'(1)=-1/2$ and $\psi (r)=0$ for $r\geq 2$.
(In \cite{whitney3}, $\psi (r)=2/(1+r)$.)
The derivative at the double point $(\pm 1,0)$ is
$$
{\textstyle
\dfrac{\partial w}{\partial s} = (1,0,\pm 1,0),\qquad
\dfrac{\partial w}{\partial y} = (0,1,0,\pm 1).
}
$$
Writing $w(s,y)$ in the form
$(s,y,0,0) +\lambda (s,y)(-s,0,-1,sy)$, we see that
$w(s,y) =(s,y,0,0)$ for $\|(s,y)\|$ large.
The immersion $w$ has $\ZZ$-symmetry as an equivariant map
$$
w: L\oplus LU\to (L\oplus LU)\times (\Rr\oplus U)\, .
$$
The two double points are distinguished by the $\ZZ$-map
$$
c : L\oplus LU \to L
$$
given by $c(s,y)= \lambda (s,y)s$.

Now Whitney's construction, applied on the fibres of $\zeta$,
gives an immersion
$$
\begin{array}{ll}
f': &\Omega = H\otimes\nu =H\oplus (H\otimes\zeta ) \\
&\qquad \to
(H\oplus (H\otimes\zeta ))\times (\Rr\oplus\zeta) =
(H\otimes\nu )\times\nu  = \nu\, |\, \Omega \subseteq N
\end{array}
$$
of the open subset $\Omega$ of $M$ into the total space of the restriction
of $\nu$ to $\Omega$.
We extend $f'$ to the whole of $M$ to coincide with the zero section
$f$ outside a compact subspace of $\Omega$.
Its double point set $\Dpone$ is the double-cover $S(H|\barD )$ of $\barD$ in
$\Omega$.
The map $c$ composed with the inclusion of $H$ into the trivial bundle
$\Omega\times V$ and the projection to $V$ gives a map
$$
e' : \Omega = H\oplus (H\otimes\zeta ) \to V\, ,
$$
which is zero outside a compact subset of $\Omega$ and can be extended
by $0$ to a map $e': M\to V$ which distinguishes the double point
pairs.
The required homotopies $e_t$ and $f_t$, $t\in [0,1]$, joining
$e$ to $e'$ and $f$ to $f'$ are
defined by replacing $\lambda$ in the definition of $e'$ and $f'$
by $t\lambda$.
One checks that $(e_t,f_t)$ is an embedding for all $t$, but that
$f_t$ fails to be an immersion when $t=1/\psi (0)$.
\Thm{Suppose that $3m<2n-1$.
Then Whitney's construction described above produces, for any given element
$\delta\in\tilde\omega_0((M\times P(V))^{H\otimes\nu -\tau M})$, a
homotopy $(e_t,f_t)$ with $\theta (e_t,f_t)=\delta$.
}
\begin{proof}
By construction, the double point manifold $\Dpone =S(H|\barD)$
represents $\delta$.
The assertion thus follows from
the Double Point Theorem \ref{hdpt} for $(e',f')$
in conjunction with the Homotopy Variation Theorem \ref{Var}.
\qed\end{proof}
\Rem{The same construction may be used to modify a general immersion
$f:M\imm N$, and map $e:M\to V$, in the complement of $\Dp$.
We can insert $\barD$ in the complement, because $2(2m-n)<m$.
\qed}
\Ex{Whitney's construction gives the classical immersions of the
sphere $S^m\imm S^{2m}$ close to the equatorial inclusion
and $S^m\imm S^m\times S^m$ close to the diagonal.
\qed}

The early work of Smale \cite{smale1,smale2} has been followed
by a vast literature on the homotopy-theoretic properties
of immersions, including \cite{dax}, \cite{haefligerhirsch}, \cite{hatcherquinn},
\cite{hirsch}, \cite{kleinwilliams2}, \cite{liliuzhang} and \cite{salomonsen}.
\Sect{Appendix: Monomorphisms of vector bundles}
Let $\xi$ and $\eta$ be smooth real vector bundles, of dimension
$n$ and $r$ respectively, over a closed
$m$-manifold $M$. We shall describe the differential-topological
classification of homotopy classes of
vector bundle monomorphisms $\eta\to\xi$ in the
metastable range $m +1 < 2(n - r)$.

Suppose that $v_0,\, v_1: \eta \into\xi$ are two vector bundle
monomorphisms.
Doing homotopy theory,
we may assume that $\xi$ and $\eta$ have positive-definite inner
products and that the monomorphisms are isometric embeddings.
Then $v_0$ and $v_1$ are sections of
the bundle $\O (\eta ,\xi )$ whose fibre at $x\in M$ is
the Stiefel manifold of orthogonal linear maps $v: \eta_x\to\xi_x$.
Topological obstruction theory gives a difference class
$$
\delta (v_0 ,v_1) \in \tilde\omega^{-1}(P(\eta )^{-H\otimes \xi})\, ,
$$
where $P(\eta )$ is the projective bundle of $\eta$ and $H$ is
the Hopf line bundle.
This arises as follows. A section of $\O (\eta ,\xi )$
determines a nowhere zero
section of $H\otimes \xi$ over $P(\eta )$:
over $\ell\in P(\eta_x)$ (where $\ell\subseteq \eta_x$ is a line)
a monomorphism $v:\eta_x\to \xi_x$ gives an embedding of
$\ell$ in $\xi_x$ and so a non-zero
vector in $\ell^*\otimes \xi_x$, which is the fibre of
$H\otimes \xi$.
Then $\delta (v_0, v_1)$ is defined as the difference class
$\delta (s_0, s_1)$ of the two nowhere zero sections $s_0$
and $s_1$ of the vector bundle $\xi$ over $P(\eta )$ constructed in this
way from $v_0$ and $v_1$.
We may assume that $s_0$ and $s_1$ are sections of the sphere bundle
$S(H\otimes \xi )$.
Write $s_t = (1-t)s_0 +ts_1$ for $0\leq t\leq 1$. Then $\delta (s_0,s_1)$
is represented explicitly by the map, over $P(\eta )$,
$$
\bar s: ([0,1],\partial [0,1]) \times P(\eta )
\to (D(H\otimes \xi ),S(H\otimes \xi ))
$$
given by the homotopy $s_t$.
In the metastable range, the vector bundle monomorphisms $v_0$ and $v_1$
are homotopic if and only if $\delta (v_0,v_1)=0$.
(See, for example, \cite{crabb,crabbjames,crabbcoincidence}.)

Thus far, the theory is topological.
We now use Poincar\'e duality for the manifold $P(\eta )$ to identify
$\tilde\omega^{-1}(P(\eta )^{-H\otimes \xi} )$
with $\tilde\omega_0(P(\eta )^{H\otimes (\xi -\eta ) -\tau M})$.
(Up to homotopy, the stable tangent bundle is given by
an isomorphism $\Rr\oplus\tau P(\eta )\iso (H\otimes\eta )\oplus\tau M$.)
Assuming that the monomorphisms $v_0$ and $v_1$ are smooth we shall
represent the dual obstruction class by
a submanifold $Z$ of $M$ together with a map $Z\to P(\eta )$ and appropriate
normal bundle information.
The monomorphism $v_0$ will play a special r\^ole in the description;
to emphasize this,
we write $i=v_0$ for the preferred embedding $i:\eta\into\xi$
and write $\nu$ for the orthogonal complement of $i(\eta )$ in $\xi$.
Let $X_k(\eta ,\nu )$, for $k\geq 1$,
be the sub-bundle of $\O (\eta ,\xi )=\O (\eta ,\eta\oplus\nu )$
with fibre consisting of those linear maps
$v$ such that $i_x+v$ has kernel of dimension $k$.
By Lemma \ref{Cayley} below, we may assume, if $m+1<2(n-r)$, that
$v_1$ never meets $X_k(\eta ,\nu )$ for
$k>1$ and is transverse to the sub-bundle
$X_1(\eta ,\nu )$. The inverse image $v_1^{-1}(X_1(\eta ,\nu ))$
is, therefore, a submanifold
$Z$ of $M$ of dimension $m+r-n$, equipped with a section $Z\to P(\eta )$
classifying the $1$-dimensional kernel. The normal bundle of
$Z$ in $M$ is identified, by Lemma \ref{Cayley}, with $H\otimes\nu$.
\Prop{\label{obstruction}
The submanifold $Z$ described above, with the line bundle $H$
classified by the section of $P(\eta )$ over $Z$ and the isomorphism
between the normal bundle and $H\otimes\nu$, represents the dual
of $\delta (v_0,v_1)$ in  $\tilde\omega_0(P(\eta )^{H\otimes \nu -\tau M})$.
}
\begin{proof}
Consider the sections $s_0$ and $s_1$ of $S(H\otimes\xi )$ over $P(\eta )$
associated with $v_0$ and $v_1$ as in the definition of
$\delta (v_0,v_1)$. The section $\bar s$ of $D(H\otimes\xi )$ over
$[0,1]\times P(\eta )$ given by the homotopy $s_t =(1-t)s_0+ts_1$
is transverse to the zero section and its zero-set is precisely
$\{ \frac{1}{2}\}\times Z$.
The normal bundle is $\Rr \oplus \tau_MP(\eta )\oplus (H\otimes\nu )$,
where $\tau_MP(\eta )$ is the tangent bundle along the fibres of
$P(\eta )\to M$, and this is identified with $(H\otimes\eta )\oplus
(H\otimes\nu )=H\otimes\xi$.
Hence, $Z$ with the normal bundle data represents the stable homotopy class
dual to $\delta (s_0,s_1)=\delta (v_0,v_1)$.
(This is the classical representation of the dual Euler
class of a vector bundle by the zero-set of a generic smooth section.)
\qed\end{proof}
\Rem{A more symmetric treatment may be given by looking
at sections of the fibre product $\O (\eta,\xi )\times_M\O (\eta ,\xi )$
and the sub-bundles with fibre consisting of the pairs
$(u,v)$ such that $\dim\ker (u+v)=k$.
}

The properties of $v_1$ required in the proof of Proposition \ref{obstruction}
follow from  the next lemma, in which
the Lie algebra of the orthogonal group $\O (V)$ of a Euclidean
vector space $V$ is written as ${\mathfrak o}(V)$.
\Lem{\label{Cayley}
Let $V$ and $W$ be finite-dimensional orthogonal vector spaces.
For $0\leq k\leq\dim V$, let $X_k(V,W)$  be the subspace
of the Stiefel manifold $\O (V,V\oplus W)$ consisting of
the maps $v$ such that $i+v$, where $i$ is the inclusion of the first
summand $V\into V\oplus W$, has kernel of dimension $k$.
Then $X_k(V,W)$ is a submanifold diffeomorphic to the total space of
the vector bundle
${\mathfrak o}(\zeta^\perp )\oplus\Hom (\zeta^\perp, W)$
over the Grassmann manifold
$G_k(V)$ of $k$-planes in $V$, where $\zeta^\perp$ is the
orthogonal complement in $V$ of the
canonical $k$-dimensional vector bundle $\zeta$ over $G_k(V)$.
Its normal bundle is naturally identified with
${\mathfrak o}(\zeta )\oplus \Hom (\zeta ,W)$.
}
\begin{proof}
This can be established by using
the (generalized) Cayley transform,
which is written down explicitly in \cite[Part II, Lemma 13.13]{crabbjames}.
The restriction of the normal bundle of the embedded submanifold
to the subspace $G_k(V)$ is naturally identified with
${\mathfrak o}(\zeta )\oplus \Hom (\zeta ,W)$.
The normal bundle itself is naturally identified with
the pullback by parallel translation.
\qed\end{proof}

In particular, the submanifold $X_k(\eta ,\nu )$ considered above
has codimension $(n-r)k+k(k-1)/2$ in $\O (\eta
,\eta\oplus\nu )$. So, if $k\geq 2$, the codimension is at least
$2(n-r)+1$. The condition $m<2(n-r)+1$ ensures that a generic
section of $\O (\eta ,\eta\oplus\nu )$ is transverse to $X_1(\eta
,\nu )$ and disjoint from $X_k(\eta ,\nu )$ for $k>1$. \Rem{In
\cite{koschorke1} Koschorke gave an intermediate representation of the
dual of $\delta (v_0,v_1)$ by a submanifold of $(0,1)\times M$.
Consider the section $\bar v$ of $\Hom (\eta ,\xi )$ over
$[0,1]\times M$ given by the homotopy $v_t=(1-t)v_0 +tv_1$ Let
$Y_k(\eta ,\xi )$ be the sub-bundle of $\Hom (\eta, \xi )$ with
fibre consisting of the linear maps with kernel of dimension $k$;
it has codimension $(n-r)k+k^2$, which is $\geq 2(n-r)+4 >m+1$.
Suppose that $\bar v$  can be deformed, by a homotopy through maps
coinciding with $v_0$ and $v_1$ at the endpoints, to a smooth
section $\bar v'$, that never meets $Y_k(\eta ,\xi )$ for $k>1$
and meets $Y_1(\eta ,\xi )$ transversely. This is always possible
if $m<2(n-r)+3$. The inverse image of $Y_1(\eta ,\xi )$ is then a
submanifold $Z$ of $(0,1)\times M$ of dimension $m+r-n$ equipped
with a map $Z\to P(\eta )$ given by the $1$-dimensional kernel.
This data, too, represents the dual of $\delta (v_0,v_1)$. \qed}

}


\providecommand{\bysame}{\leavevmode\hbox
 to3em{\hrulefill}\ }

\printindex


\end{document}